\documentclass[11pt]{article}
\usepackage{amssymb,amsmath, amsthm, amsopn, amsfonts,bm,mathtools}
\usepackage[english]{babel}
\usepackage[all]{xy}
\usepackage{fig4tex}
\usepackage{color}\definecolor{rd}{rgb}{1,0.3,0.35}
\usepackage{graphicx}
\usepackage{stmaryrd}
\usepackage{mathrsfs}
\usepackage{hyperref}
\usepackage{caption}
\usepackage{fancyvrb}
\usepackage{fvextra}

\DeclareMathOperator*{\card}{card}
\DeclareMathOperator{\Vect}{Span}

\usepackage[hmargin=2.5cm,vmargin=3cm,centering]{geometry}
\newtheorem{thm}{Theorem}[section]
\newtheorem{cor}[thm]{Corollary}
\newtheorem{lem}[thm]{Lemma}
\newtheorem{definition}[thm]{Definition}
\newtheorem{prop}[thm]{Proposition}
\newtheorem{remark}[thm]{Remark}
\newtheorem{hyp}[thm]{Hypothesis}
\theoremstyle{definition}

\def\ds{\displaystyle}

\def\H{{\mathcal{H}}}

\def\N{{\bf N}}
\newcommand{\field}[1]{\mathbb{#1}}
\newcommand{\rz}{\field{R}}
\newcommand{\cz}{\field{C}}
\newcommand{\nz}{\field{N}}
\newcommand{\zz}{\field{Z}}

\newcommand{\sz}{\field{S}}
\newcommand{\tz}{\field{T}}
\newcommand{\kz}{\field{K}}

\def\11{{\rm 1~\hspace{-1.2ex}l} }

\def\pa{\partial}
\def\ccup{\mathop{\cup}}

\def\supp{{\rm supp~}}
\def\Real{{\mathrm{Re~}}}

\def\Id{\mathrm{Id}}
\def\Ran{\mathrm{Ran} }
\def\supp{\mathrm{supp}~}
\begin{document}
\title{Bar codes of persistent cohomology and Arrhenius law for $p$-forms}
\author{
D.~Le~Peutrec\thanks{Dorian.LePeutrec@math.u-psud.fr,
Universit{\'e} Paris-Saclay, CNRS, Laboratoire de math{\'e}matiques d'Orsay,
F-91405 Orsay, France}\\
F.~Nier\thanks{nier@math.univ-paris13.fr, Universit{\'e} Sorbonne Paris Nord,
LAGA, UMR-CNRS 7539, 99 av.~J.B.~Cl{\'e}ment, F-93430
Villetaneuse, France.}\\
C.~Viterbo\thanks{Claude.Viterbo@ens.fr,
DMA, UMR-CNRS 8553, Ecole normale sup{\'e}rieure/PSL University, 45 rue d'Ulm, 75230 Paris Cedex, France.}
}
\maketitle
\abstract{This article shows that counting or computing the small eigenvalues of the Witten Laplacian in the semi-classical limit can be done without assuming
that the potential is a Morse function  as the authors did in \cite{LNV}. In connection with persistent cohomology, we prove that the rescaled logarithms of these small  eigenvalues are asymptotically determined by the 
lengths of the bar code of the function $f$\,. In particular, this
proves that these quantities are stable in the $C^0$ topology on the
space of functions. Additionally, our analysis provides a general
method for computing the subexponential corrections in a large number
of cases.\\


\tableofcontents{}

\section{Introduction}
\label{sec:intro}

\subsection{Motivations}
\label{sec:motiv}

Since its discovery in the late nineteenth century, Arrhenius law (see \cite{Arr}) is
one of the most robust laws of chemistry or physics. Actually, its range
of applications has increased over decades and is now also commonly used  in
biology or social sciences  as an empirical law
whose parameters be can figured out rather easily, even when the
microscopic or individual mechanisms
 are not well understood. Its early
interpretations were done within the thermodynamical or statistical
physics framework. They are now formulated in the modern
and general language of stochastic processes, more specifically of 
the Brownian motion of a particle
evolving  in a gradient field. At low temperature $h>0$ in some
dimensionless scaling, the lifetime $\tau_{\alpha,h}$
of the state $\alpha$
is exponentially large with
\begin{equation}
  \label{eq:lifetime}
\log\tau_{\alpha,h}\ {\sim}\ \frac{\ell_{\alpha}}{h}\,,
\end{equation}
where $\ell_{\alpha}$ is the energy variation between a local minimum
and the lowest saddle point that we need to cross to reach a state of  lower energy. Practically
and as an illustration of the robustness of Arrhenius law,  it is neither necessary to know  the energy landscape nor
the configuration space: in the end only the $\ell_{\alpha}$'s are
important and they are determined experimentally, e.g. in chemistry
kinetics.  
A general justification of \eqref{eq:lifetime}  was proposed by Freidlin and Wentzell in
\cite{VeFr1,VeFr2} relying on large deviation arguments (see also
\cite{FrWe} and \cite{Ber} for a wider overview).\\
In an energy landscape described by the function $2f:M\to \rz$\,, those lifetimes are 
generically
the inverses of eigenvalues of   the operator $-h\Delta+2\nabla
f\cdot\nabla$ in $L^{2}(M,e^{-\frac{2f}{h}}~dx)$\,, where
$e^{-\frac{2f}{h}}~dx$ is the associated invariant measure (it exists e.g.
when $M$ is a compact Riemannian manifold without boundary). 
After a conjugation by $e^{\frac{f}{h}}$ and a multiplication by $h$
(corresponding to a change of time scale), 
it becomes  the operator
$$
\Delta_{f,h}^{(0)}=-h^{2}\Delta +|\nabla f(x)|^{2}-h(\Delta
f)(x)=d_{f,h}^{*}d_{f,h}\quad\text{acting~in}~ L^{2}(M,dx)\,,
$$
where $d_{f,h}=e^{-\frac{f}{h}}(hd)e^{\frac{f}{h}}$ is the Witten
differential and $d_{f,h}^{*}$ its adjoint. 
This operator acts on
general differential forms as the Witten Laplacian, a
deformation of the Hodge Laplacian:
$$
\Delta_{f,h}=(d_{f,h}+d_{f,h}^{*})^{2}=\mathop{\oplus}_{p=0}^{\dim M}\Delta_{f,h}^{(p)}\,,
$$
where the direct sum separates the degrees. When $f$ is a Morse
function, Witten in \cite{Wit} (see also \cite{CFKS})
proved that as $h$ goes to zero,  the eigenvalues of $\Delta_{f,h}$
are divided into two groups, given in our scaling as one bounded from below by $C_{f}h$ for some
$C_{f}>0$\,, and one being of the order~$o(h)$\,. 
The  small (here $o(h)$) eigenvalues of $\Delta_{f,h}^{(p)}$ correspond to critical
points of index $p$: this is intuitively to be expected, since the
eigenfunctions should concentrate in the region where $\vert \nabla f
\vert$ is small, that is near the critical points of $f$\,.  This argument
provided an analytical proof of Morse inequalities, in the line of
several results relating topological quantities and spectral analysis,
one of the earliest being the Atiyah-Patodi-Singer proof of the index
theorem (see~\cite{APS}). 

 In
\cite{HeSj4}\,, Helffer and Sj{\"o}strand gave a rigorous proof of Witten's claims and proved that those small
eigenvalues were actually exponentially small, without specifying their
size. This was later extended to Morse-Bott functions by Bismut and Helffer-Sj\"ostrand (see \cite{Bis} and \cite{HeSj6}).  After this, many applications of Witten Laplacians or more
general Witten deformations 
were used to study various 
global  topological invariants of manifolds or fibre
bundles by counting the small eigenvalues of such operators (see
e.g. \cite{BiZh,Zha,ChLi}).\\

When $f$ is a Morse
function, the Arrhenius law in degree $0$ says that the $o(h)$ eigenvalues of
$\Delta_{f,h}^{(0)}$ satisfy
\begin{equation}
\label{eq:Witas}
\log\lambda_{\alpha,h}^{(0)}\ \sim\ -2\frac{f(y_{\alpha})-f(x_{\alpha})}{h}\quad
\text{as}~h\to 0^{+}\,,
\end{equation}
where $x_{\alpha}$ is a local minimum and $y_{\alpha}$ is an
associated saddle point.
Already around 1935, Eyring and Kramers (see \cite{Eyr,Kra}), motivated by the theory of the
activated complex in chemistry, 
 proposed a more accurate version 
which reads here
\begin{equation}
  \label{eq:EyrKra}
\lambda_{\alpha,h}^{(0)}\sim
\frac{h}{\pi}C_{\alpha}e^{-2\frac{f(y_{\alpha})-f(x_{\alpha})}{h}}\quad\text{as}~h\to 0^{+}\,,
\end{equation}
where the constant $C_{\alpha}$ depends on the Hessians at
the non degenerate critical points $x_{\alpha}$ and $y_{\alpha}$\,,
$x_{\alpha}$ is a local minimum (here a critical point of index $0$), and
$y_{\alpha}$ a saddle point (here a critical point of index $1$)\,.\\
The first mathematical proof of the Eyring-Kramers formula was performed in
degree~$0$ in \cite{BEGK,BGK} by using potential theoretic and
capacity arguments, and in \cite{HKN} by
improving Helffer-Sj{\"o}strand's semiclassical analysis for
$\Delta_{f,h}^{(0)}$ (see also e.g.  the prior works \cite{HKS,Micl} for results less precise
than \eqref{eq:EyrKra} but more precise than \eqref{eq:Witas}). 
These results were proved under the
assumption that $f$ is a Morse
function  with simple local
minima and simple saddle
points (a Morse function has simple critical values or critical points
if every critical value is the image of  a single critical point),
and with distinct lengths : the real numbers
$\ell_{\alpha}=2(f(y_{\alpha})-f(x_{\alpha}))$ are all distinct.
The pairing between
local minima $x_{\alpha}$ and saddle points $y_{\alpha}$ (critical
points with index $1$) was done by
extending the intuitive picture of basins of attraction, more precisely by
considering the connected components of sublevel sets of $f$\,. 
Note that this differs from the 
instantonic picture, associated with curves which are   intersections of
stable and unstable manifolds of $-\nabla f$\,, which is in some
sense local  and would lead to a complicated analysis of cancellations while
computing precisely the $\lambda_{\alpha,h}^{(0)}$'s.  This
pairing relies on global topological considerations which are robust
with respect to the $C^0$ perturbations of the energy profile $2f$\,. 
By making use of the min-max principle,  
it is  actually not
difficult to start from the analysis
done in \cite{HKN} for Morse functions and to recover
\eqref{eq:Witas}  and the results of 
\cite{VeFr1,VeFr2,HKS,Micl} in cases where the
local minima are degenerate.
\\
The situation is completely different for general differential forms of degree $p$\,. In
\cite{LNV}, we proved an Eyring-Kramers law (and therefore an
Arrhenius law) by assuming again that the function $f$
was a generic Morse function with simple critical
values and such that the difference between critical values were all
distinct.  
 Here the problem is to understand which critical values $f(x_{\alpha})$ and $f(y_{\alpha})$ are paired in order to compute the exponential factors. This pairing is obtained topologically by  using a refinement of Barannikov's
presentation of Morse theory. This can be restated in modern terms with the bar code of $f$\,, denoted
$B_{f}=([a_{\alpha}^{*},b_{\alpha}^{*+1}[)_{\alpha \in A^{*}}$\,, associated with the 
Morse function $f$ on $M$\,, with the notation $a_{\alpha}^{(p)}=f(x_{\alpha})$ and
$b_{\alpha}^{(p+1)}=f(y_{\alpha})$\,, where the critical point
$x_{\alpha}$ has index  $p$ and
$y_{\alpha}$ has index $p+1$\,. 
 Later, it was noticed in
\cite{UsZh,PoSh} that those bar codes were nothing but the bar codes of persistent homology, developed since the beginning of the
21st century (see \cite{EdHa} for a historical review).
An important
feature of the Barannikov complex, and hence  of persistent homology, is the stability result which says
in the latter framework
$$
d_{bot}(\mathcal{B}_{f},\mathcal{B}_{g})\leq \|f-g\|_{\mathcal{C}^{0}}\,,
$$
 where the bottleneck distance $d_{bot}$ estimates the variations of
 the lengths of the bars.\\

But the bar code of a function is defined for any continuous function, except the bars are now infinitely many, with the property that  for any $ \varepsilon_0 >0$\,, only finitely many are greater than $ \varepsilon_0$\,.  
It is then natural to state the following conjecture.\\

\noindent
\textbf{Main Conjecture~:} \textit{Consider a $\mathcal{C}^{\infty}$ (or even Lipschitz) function $f$ on a
compact manifold $M$ with  bar code $\mathcal{B}_{f}$\,. 
We denote by $A^{(p)}(\ell)$ the set of bars in $\mathcal{B}_{f}$ of the type $[a_\alpha^{(p)},b_\alpha^{(p+1)}[$
with $b_\alpha^{(p+1)}-a_\alpha^{(p)}>\ell$\,,  and 
 $A_c^{(p-1)}(\ell)$  the set of bars in $\mathcal{B}_{f}$ of the type $[a_\alpha^{(p-1)},b_\alpha^{(p)}[$ with  $b_\alpha^{(p+1)}-a_\alpha^{(p)}>\ell$  and $b_\alpha^{(p)}<+\infty$\,. 
Then, there exists  $\varepsilon_0 >0$ such that, for every
$ \varepsilon \in]0,\varepsilon_0]$\,,
$\Delta_{f,h}^{(p)}$  admits $\sharp\big(A^{(p)}(\ell) \cup  A_c^{(p-1)}(\ell)\big)$
eigenvalues $\lambda_{\alpha,h}^{(p)}$
smaller than  $e^{-2\frac{\ell+ \varepsilon}{h}}$ (with multiplicity), where $\alpha\in A^{(p)}(\ell) \cup  A_c^{(p-1)}(\ell)$\,.
They can moreover be labelled such that
\begin{displaymath} 
\forall\,\alpha\in A^{(p)}(\ell) \cup  A_c^{(p-1)}(\ell)\,,\ \ \ \log \lambda_{\alpha,h}^{(p)}\sim -2\frac{b_{\alpha}^{(p+1)}-a_{\alpha}^{(p)}}{h} \qquad \text{as}~h \to 0^{+}\,.
\end{displaymath} 
}\\

The goal of this paper is to prove this conjecture  under the assumption that $f$ has a
finite number of critical values.
\\

Note that we do not assume in the Main Conjecture (as well as in our theorems) that $f$ is Morse. 
One important consequence of the Main Conjecture (and hence of our main
theorems) is that the decay rate of the eigenvalues is continuous in
$f$ for the $C^0$ topology. This is not the case for subexponential
factors, since they usually depend on the eigenvalues of the Hessian
of $f$ at the critical points. 

In the case $p=0$ of functions, the Eyring-Kramers law \eqref{eq:EyrKra} has been extended in the form
$\lambda_{\alpha,h}^{(0)}\sim
C_{\alpha}(f)h^{\nu_{\alpha}(f)}e^{-2\frac{f(y_{\alpha})-f(x_{\alpha})}{h}}$
when $f$ is not a Morse function or
when $f$ is a Morse function with multiple critical values
(i.e. the preimage of a critical value may contain several critical points), the
latter appearing in  practical situations with
natural symmetries. We refer for example to
\cite{BeGe,BeDu,Mic,DLLN2, LeNe1, LeNe2}, whence  it appears that  the exponent
$\nu_{\alpha}(f)$\,, or the constant $C_{\alpha}(f)$ in the
subexponential factor, may be discontinuous when 
a general function $f$ is approximated by a sequence of generic Morse
functions. On the other hand, it will follow from our results that the
$\ell_{\alpha}=2(f(y_{\alpha})-f(x_{\alpha}))$ are stable. Understanding how the eigenvalues
$\lambda_{\alpha,h}(f)$ or the lifetimes $\tau_{\alpha,h}(f)$
depend on $f$ is also important for applications to acceleration of
stochastic algorithms  (see \cite{LeNi, DLLN1,DLLN2, LeNe1, LeNe2} and
references therein). This leads to the 
\\

\noindent
\textbf{Main Question :} \textit{Is there a way to analyze how the
subexponential factor of Eyring-Kramers law for $p$-forms varies when $f$ is
changed~? In particular, does it explain the observed
discontinuities~?}\\

Again, the answer is yes. Our presentation of Arrhenius law for
$p$-forms
provides a very general result. The method actually completely  separates
the determination of the exponential scales
$e^{-\frac{\ell_{\alpha}}{h}}$\,, related with global algebraic topological
objects, from the determination of the subexponential factors, which
rely on some local analysis. Many applications with various
discontinuous effects will be presented at the end of this text.
Actually, the discontinuities w.r.t. the energy landscape $f$
 of the leading term for the
subexponential factor $C_{\alpha}(f)h^{\nu_{\alpha}(f)}$ are easily
understood on the simple example of the Laplace integrals 
\begin{eqnarray*}
&I(\delta,h)=
  \int_{\rz}e^{-\frac{x^{4}/4-\delta x^{2}/2+1_{\rz^{+}}(\delta)\delta^{2}/4}{h}}~dx\,,
&\\
\text{which satisfy} & I(\delta,h)\stackrel{h\to 0}{\sim}C_{\delta}h^{1/2}&\text{when}~\delta\neq 0\,,\\
\text{and}& I(\delta,h)\stackrel{h\to 0}{\sim}Ch^{1/4}&\text{when}~\delta=0\,.
\end{eqnarray*}

\subsection{General assumptions and notations}
\label{sec:genass}
\textbf{The manifold $M$:}
The Riemannian manifold $(M,g)$ is assumed compact without boundary
 with $\dim_{\rz} M=d$ and non necessarily
oriented. Some non compact manifold will be considered 
 in
Subsection~\ref{sec:moregenmfld}. 
In the non-orientable case, the Hodge star operator, $\star$\,, sends
$\Lambda T^{*}M=\mathop{\oplus}_{p=0}^{d}\Lambda^{p}T^{*}M$ to $\Lambda T^{*}M\otimes_{M}\mathrm{or}_{M}$\,, where $\mathrm{or}_M$ is the orientation (line) bundle, which is of course locally trivial.  When
$N\subset M$ is a regular hypersurface admitting a global unit normal
(or conormal) vector the orientation twist $\mathrm{or}_N$ is the restriction of ${\mathrm{or}_M}$\,. \\
In local coordinates the metric will be written $g=g_{ij}(x)dx^{i}dx^{j}$
with $g^{-1}=g^{ij}(x)\frac{\partial}{\partial
 x^{i}}\frac{\partial}{\partial x^{j}}$
and the musical isomorphisms ${}^{\sharp}:T^{*}M\to TM$ and
${}^{\flat}:TM\to T^{*}M$ are given by
$$
(\omega_{i}dx^{i})^{\sharp}=g^{ij}\omega_{j}\frac{\partial}{\partial
  x^{i}}\quad\text{and}\quad (X^{i}\frac{\partial}{\partial x^{i}})^{\flat}=g_{ij}X^{j}dx^{i}\,.
$$
The
differential $d$ acts on ${\cal C}^{\infty}(M; \Lambda
T^{*}M\otimes_{M}\cz)$ or ${\cal D}'(M; \Lambda T^{*}M\otimes_{M}\cz)$ and augments the degree of forms by $1$\,.
The codifferential $d^{*}=(-1)^{\text{deg}}\star^{-1}d\star$ acts on ${\cal C}^{\infty}(M;\Lambda
T^{*}M\otimes_{M}\cz)$ and ${\cal D}'(M;\Lambda
T^{*}M \otimes_{M}\cz)$ and decreases the degree
by $1$\,.
In the sequel and unless otherwise specified, we always consider
complex valued differential forms and the tensorization
 by $\cz$ will
be omitted in the notation. The duality bracket
$\langle~\,,\,~\rangle$ between ${\cal D}'(M;\Lambda^{p}T^{*}M\otimes \mathrm{or}_{M})$ and
${\cal C}^{\infty}(M;\Lambda^{d-p}T^{*}M)$
(where ${\cal D}'$ and ${\cal C}^{\infty}$ can be interchanged)
is  assumed $\cz$-antilinear
on the left-hand side and $\cz$-linear on the right-hand side. Stokes's formula  then implies that $d^{*}$
is the formal adjoint of $d$ according to 
$$
0=\int_{M}d(\overline{\omega}\wedge \star \eta)=\int_{M}d\overline{\omega}\wedge (\star
\eta)+(-1)^{\text{deg}\,\omega}\overline{\omega}\wedge d(\star\eta)=
\langle d\omega\,,\,\eta\rangle-\langle\omega\,,\, d^{*}\eta\rangle
$$
for $\omega,\eta\in {\cal C}^{\infty}(M;\Lambda^{p-1}T^{*}M)$\,.\\

\noindent\textbf{Functional spaces:} The $L^{2}$-norm of sections of $\Lambda
T^{*}M$ is the one given by the metric $g$ and we recall
$$
\int_{M}\langle \omega\,,\, \eta\rangle_{\Lambda
  T_{q}^{*}M}~d\text{vol}_{g}(q)=\int_{M}\overline{\omega}\wedge \star \eta\,.
$$
We use the notation $W^{s,p}$ for the Sobolev space with $s$
derivatives in $L^{p}$\,. In particular, $W^{s,2}$ corresponds to the
standard Hilbertian Sobolev spaces while $W^{1,\infty}$ will be used
for the set of Lipschitz functions.
For an open domain $\Omega\subset M$ and for $s\in\rz$\,,  the notation
$W^{s,2}(\overline{\Omega};\Lambda T^{*}M)$ denotes the set of
restrictions to $\Omega$ of $W^{s,2}$-sections in $M$\,, and when there is
no ambiguity or necessity, we shall use the short version
$W^{s,2}(\overline{\Omega})$\,. The same definition holds for ${\cal
  C}^{\infty}(\overline{\Omega};\Lambda T^{*}M)$\,.
We recall that when $\Omega$ is a regular domain, that is when
$\partial \Omega$ is a ${\cal C}^{\infty}$ hypersurface,
$W^{s,2}(\overline{\Omega};\Lambda T^{*}M)$ coincides with  $W^{s,2}(\Omega;\Lambda T^{*}M)$ 
by interpolation and duality from the special cases of  $s\in \nz$ (see
e.g. \cite{ChPi}). In such a case, the trace theorem holds from
$W^{s,2}(\Omega;\Lambda T^{*}M)$ to $W^{s-1/2,2}(\Omega; \Lambda
T^{*}\partial \Omega)$ for $s>\frac{1}{2}$\,.
 The local regularity theory is not affected
when sections of $\Lambda T^{*}M\otimes \mathrm{or}_{M}$ and $\Lambda
T^{*}M\otimes \mathrm{or}_{\partial \Omega}$ are considered and we
shall use  the short notation $W^{s,2}(\overline{\Omega})$ or
$W^{s,2}(\partial \Omega)$ indifferently 
for sections of the trivial and orientation line bundles, unless we need to distinguish
the global behaviour. Other  functional
spaces will be introduced later in our analysis.\\

\noindent\textbf{Witten differential and Witten Laplacian:}
The Witten differential and the Witten Laplacian are deformations of
the differential $d$ and the Hodge Laplacian $dd^{*}+d^{*}d$
associated with a real valued function $f$ and a positive parameter
$h>0$  in the asymptotics $h\to 0$\,.
\begin{definition}
\label{de:fab}
Let $f$ be a real valued function on $M$\,.
  For $a\in \overline{\rz}= \rz \cup \{-\infty, + \infty\}$\,, we use the notations
  \begin{eqnarray*}
    &&f^{a}=\left\{x\in M\,, f(x)<a\right\}\quad,\quad f^{\leq
  a}=\left\{x\in M\,, f(x)\leq a\right\}\,,\\
&&
f_{a}=\left\{x\in M\,, f(x)>a\right\}\quad,\quad f_{\geq
   a}=\left\{x\in M\,, f(x)\geq a\right\}\,,
\end{eqnarray*}
with all the combinations like $f_{a}^{b}=\left\{x\in M\,,\quad a<f(x)<b\right\}$\,.
\end{definition}
\noindent 
Although weaker regularity assumptions for the function $f$ will be discussed later, the
following simple hypothesis will be convenient for us.
\begin{hyp}
\label{hyp:mainf}
  The function $f$  on $(M,g)$ is assumed to be  Lipschitz 
  with a finite number $N$ of values
  $c_{1},\ldots,c_{N_{f}}$ such that:
  \begin{itemize}
  \item $f\in {\cal C}^{\infty}(M\setminus
    f^{-1}(\left\{c_{1},\ldots,c_{N_{f}}\right\});\rz)$
\item $\forall x\in M\setminus
  f^{-1}(\left\{c_{1},\ldots,c_{N_{f}}\right\})\,,\quad \left|\nabla f(x)\right|\neq 0$\,.
  \end{itemize}
\end{hyp}
\noindent When $f\in {\cal C}^{\infty}(M;\rz)$\,, the above assumption simply says
that $f$ has a finite number $\leq N_{f}$ of critical values\,. For a Lipschitz function, we count also ``fake'' critical values allowing singularities
of $f$ at those values. 
We nevertheless call $c_{1},\ldots, c_{N_{f}}$ the ``critical values''
of $f$ and use the notation  
$$
M_{reg}=\left\{x\in (M\setminus\mathrm{suppsing}\; f)\,, \nabla
  f(x)\neq 0\right\}\subset 
M\setminus f^{-1}(\left\{c_{1},\ldots, c_{N_{f}}\right\})\,.
$$
\noindent When $M$ is a real analytic manifold,
Hypothesis~\ref{hyp:mainf} may be replaced by the following
 simpler natural assumption.
\begin{hyp}
\label{hyp:realana} On the real analytic  compact Riemannian
manifold $M$\,, $f$ is a Lipschitz subanalytic function.
\end{hyp}
Actually, the proof of the main result, Theorem~\ref{th:induc}, will
hold under Hypothesis~\ref{hyp:mainf} or under some milder assumptions  which
are more technical and will  appear as consequences of
Hypothesis~\ref{hyp:realana} in Subsection~\ref{sec:moregenLip}.
We  also  refer to Subsection~\ref{sec:moregenLip} for more material on Lipschitz subanalytic functions.\\

\noindent Under Hypothesis~\ref{hyp:mainf} or more generally for a
Lipschitz function $f$ and for
$h>0$\,, the differential operators $d_{f,h}$\,, $d_{f,h}^{*}$ and
$\Delta_{f,h}$ are defined by:
\begin{eqnarray}
\label{eq:dfh}
  &&d_{f,h}=e^{-\frac{f}{h}}(hd)e^{\frac{f}{h}}=hd
     +df\wedge\quad,\quad d_{f,h}\circ d_{f,h}=0\,,\\
\label{eq:dfhstar}
&&
   d_{f,h}^{*}=e^{\frac{f}{h}}(hd^{*})e^{-\frac{f}{h}}=hd^{*}+\mathbf{i}_{\nabla
   f}\quad,\quad d_{f,h}^{*}\circ d_{f,h}^{*}=0\,,\\
\label{eq:Deltafh}
&&
   \Delta_{f,h}=(d_{f,h}+d_{f,h}^{*})^{2}=d_{f,h}^{*}d_{f,h}+d_{f,h}\circ
   d_{f,h}^{*}=h^{2}\Delta_{0,1}+|\nabla f(x)|^{2}+h({\cal L}_{\nabla
   f}+{\cal L}_{\nabla f}^{*})\,.
\end{eqnarray}
The above identities make sense when considering $d_{f,h}$ and $d_{f,h}^{*}$ as operators from
 $W^{1,2}(M)$ to $L^{2}(M)$ or from 
$L^{2}(M)$ to $W^{-1,2}(M)$\,, and  for the compositions of two of them
and for $\Delta_{f,h}$\,,
as operators from $W^{1,2}(M)$ to $W^{-1,2}(M)$\,.  We shall be more precise on requirements
for  domains once we add the boundary conditions.\\

\noindent\textbf{Convention for closed operators and quadratic
  forms:}\\
We shall
consider various closed realizations in $L^{2}$ spaces
of the above differential
operators $d_{f,h}$\,, $d_{f,h}^{*}$\,, and $\Delta_{f,h}$\,, which will be
denoted $ d_{f,\bullet,h}\,,\, d_{f,\bullet,h}^{*},$ and
$\Delta_{f,\bullet,h}$\,, where the subscript $\bullet$ will specify
the realization. When $A$ is
a closed operator in a Hilbert space (resp. when $Q$ is a closed
quadratic form), writing $Au$ (resp. $Q(u)$ or $Q(u,v)$ for the
associated sesquilinear form) means that $u$ belongs to the domain of
$A$ (resp. $u$ or $u,v$ belong to the domain of $Q$). For example
$d_{f,\bullet,h}\omega=\alpha\in L^{2}$ means in particular $\omega\in
D(d_{f,\bullet,h})$\,, possibly imposing  boundary conditions.\\

\noindent\textbf{Comparing exponential scales:}
\begin{definition}
  For two functions $F,G:]0,h_{0}[\to \cz$\,, one says 
  \begin{itemize}
  \item $F(h)=\tilde{O}(G(h))$ if: 
$$
\forall \varepsilon>0\,,\exists h_{\varepsilon}, C_{\varepsilon}>0\,, \forall h\in
]0,h_{\varepsilon}[ \,, \quad |F(h)|\leq C_{\varepsilon}|G(h)|e^{\frac{\varepsilon}{h}}\,;
$$
\item $F(h)=\tilde{o}(G(h))$ if:
$$
\exists \varepsilon,h_{\varepsilon},C_{\varepsilon}>0\,,\forall
h\in]0,h_{\varepsilon}[\,, 
\quad |F(h)|\leq C_{\varepsilon}|G(h)|e^{-\frac{\varepsilon}{h}}\,;
$$
\item $F(h)\stackrel{\log}{\sim}G(h)$ if:
$$
|F(h)|=\tilde{O}(|G(h)|)\quad\text{and}\quad |G(h)|=\tilde{O}(|F(h)|)\,.
$$
  \end{itemize}
When $|F|,|G|>0$\,, the above three conditions can be written
respectively 
\begin{eqnarray*}
  &&\limsup_{h\to
      0}h\log\left(\frac{|F(h)|}{|G(h)|}\right)\leq 0\,,\\
&&\limsup_{h\to
  0}h\log\left(\frac{|F(h)|}{|G(h)|}\right)<0\\
&&
\lim_{h\to
  0}h\log\left(\frac{|F(h)|}{|G(h)|}\right)=0\,.
\end{eqnarray*}
In the two first definitions, the constant $C_{\varepsilon}$ can be
fixed to $1$ by changing $h_{\varepsilon}$ (and $\varepsilon$  in
the second definition).\\
When $F:X\times ]0,h[\to \cz$\,, the statements
``$F(x,h)=\tilde{\mathcal{O}}(G(h))$ (or $F(x,h)=\tilde{o}(G(h))$) 
(locally) uniformly''  are used when the
above definitions make sense for the corresponding suprema $\sup_{x}F(x,h)$\,.
\end{definition}

\noindent\textbf{Bar code:}\\
Although a more precise definition and construction will be recalled
especially in Appendix~\ref{app:perscohom}, we can start with a short
definition.
\begin{definition}
  Under Hypothesis~\ref{hyp:mainf}, a (persistence cohomology)
 bar code associated with $f$ is a
  finite family ${\cal B}=([a_{\alpha},b_{\alpha}[)_{\alpha\in A}$
with $-\infty<a_{\alpha}<b_{\alpha}\leq +\infty$\,, $a_{\alpha}\in
\left\{c_{1},\ldots, c_{N_{f}}\right\}$\,,\, $b_{\alpha}\in
\left\{c_{2},\ldots,c_{N_{f}},+\infty\right\}$\,, with the following
properties:
\begin{itemize}
\item 
it is graded according
 to $A=\sqcup_{p=0}^{d} A^{(p)}$\,, $[a_{\alpha},
b_{\alpha}[=[a_{\alpha}^{(p)},b_{\alpha}^{(p+1)}[$ when $\alpha\in
A^{(p)}$\,;
\item  for any pair $a,b$\,, $a<b$\,, $a,b\not
\in\left\{c_{1},\ldots c_{N_{f}}\right\}$\,, there exists a basis of
the relative homology vector space $H^{p}(f^{b},f^{a})$
indexed by the bars  of degree $p$ 
 with a unique endpoint lying in $]a,b[$\,. In particular, the
 relative Betti number is given by:
$$
\beta^{p}(f^{b},f^{a})=\dim H^{p}(f^{b},f^{a})=
\sharp\left\{\alpha\in A^{(p)}\,, \sharp
  \left\{a_{\alpha}^{(p)},b_{\alpha}^{(p+1)}\right\}\cap ]a,b[=1\right\}\,.
$$
\end{itemize}
\end{definition}
For a general Lipschitz function, such a finite bar code is well defined under
the following assumption (see Subsection~\ref{sec:moregenLipgen} and  Appendix~\ref{app:perscohom}).
\begin{hyp}
\label{hyp:Lipbar}
The function $f:M\to \rz$ is a Lipschitz function and there exists a
finite number of values $c_{1}< c_{2}\ldots<c_{N_{f}}$ such that for
any $a\in \rz \setminus\left\{c_{1},\ldots, c_{N_{f}}\right\}$\,, the
following property holds along $f^{-1}(\left\{a\right\})$:\\
For any $x_{0}\in f^{-1}(\left\{a\right\})$\,, there exists a
neighborhood $U_{x_{0}}$ of $x_{0}$ in $M$\,, a local coordinate
system $x=(x^{1},x')\in \rz\times \rz^{d-1}$\,, and a constant
$C_{x_{0}}$ such that
$$
\forall x=(x^{1},x'), y=(y^{1},x')\in U_{x_{0}}\,,\quad
\frac{1}{C_{x_{0}}}|x^{1}-y^{1}|\leq |f(x^{1},x')-f(y^{1},x')|\,.
$$
\end{hyp}
This notion of bar code, and especially the
identification of two bar codes, after possibly adding empty
intervals, is better understood after associating 
with a bar code ${\cal
  B}_{A}=([a_{\alpha},b_{\alpha}[)_{\alpha\in A}$ the constructible
sheaf $\oplus_{\alpha\in A}\kz_{[a_{\alpha},b_{\alpha}[}$ of
$\kz$-vector spaces,
 on $\rz$\,. Then, a persistence bar code 
associated with a function $f$ satisfying  
Hypothesis~\ref{hyp:mainf} is essentially unique and 
then denoted ${\cal B}(f)$\,.\\
After possibly adding empty bars such that $a_{\alpha}=b_{\alpha}$  or
$c_{\beta}=d_{\beta}$\,, two
different bar codes ${\cal
  B}_{A}=\left([a_{\alpha},b_{\alpha}[\right)_{\alpha\in A}$ and ${\cal
  B}_{B}=\left([c_{\beta},d_{\beta}[\right)_{\beta\in B}$ can be assumed
with the same cardinality, $\sharp A=\sharp B$\,. The bottleneck
distance is then defined by
$$
d_{bot}({\cal B}_{A},{\cal B}_{B})=\inf_{j:A\stackrel{\mathrm{bij}}{\to}B}\max_{\alpha\in
A}\max (|a_{\alpha}-c_{j(\alpha)}|, |b_{\alpha}-d_{j(\alpha)}|)\,,
$$
with the convention $|(+\infty)-(+\infty)|=0$\,.\\
The stability theorem for persistent (co)homology (see e.g. \cite{CEH,KaSc}) says
that for two functions $f_{1},f_{2}$ which satisfy
Hypothesis~\ref{hyp:mainf} or Hypothesis~\ref{hyp:Lipbar},
$$
d_{bot}({\cal B}(f_{2}),{\cal B}(f_{1}))\leq \|f_{2}-f_{1}\|_{{\cal C}^{0}}\,.
$$

\subsection{Simple results}
\label{sec:illres}
The method presented in this text leads to several results and can
actually be extended to other cases. Essentially, we show that
 the usual
generic assumption that the function $f$ is a Morse function can be
replaced by a very general one, after replacing the algebraic
topological information in terms of Morse  indices by the ones given
by the persistent cohomology bar code associated with $f$\,.
The following simple statements
illustrate what can be obtained.
\begin{thm}
\label{th:mainsimple}
Assume that $f$ satisfies Hypothesis~\ref{hyp:mainf} and let
$\Delta_{f,M,h}$ be the self-adjoint Witten Laplacian defined with
$D(\Delta_{f,M,h})=\left\{\omega\in W^{1,2}(M)\,, \Delta_{f,h}\omega\in
  L^{2}(M)\right\}$ and $\Delta_{f,M,h}\omega=\Delta_{f,h}\omega$
according to \eqref{eq:Deltafh}, and
$\Delta_{f,M,h}=\mathop{\oplus}_{0\leq p\leq
  d}^{\perp}\Delta_{f,M,h}^{(p)}$\,. Let ${\cal B}(f)$ be a
persistent cohomology bar code associated with $f$\,. Then, there is a
bijection between $A^{(p)}\sqcup \left\{\alpha\in A^{(p-1)}\,,
  b_{\alpha}^{(p)}\neq +\infty\right\}$ and the $\tilde{o}(1)$
eigenvalues counted with multiplicities of $\Delta_{f,M,h}^{(p)}$\,.
Precisely, there exists $\varepsilon_{0}>0$ small enough such that for
all $\varepsilon\in ]0,\varepsilon_{0}[$\,, there exists
$h_{\varepsilon}>0$ such that the $\tilde{O}(e^{-\frac{\varepsilon}{h}})$-eigenvalues of $\Delta_{f,h}^{(p)}$
counted with multiplicity for 
$h\in ]0,h_{\varepsilon}[$ are given by
$\lambda_{\alpha}^{(p)}(h)$\,, $\alpha\in A^{(p)}$ or ($\alpha\in
A^{(p-1)}$ and $b_{\alpha}^{(p)}\neq +\infty$), with
\begin{eqnarray*}
\text{either} 
   ~b_{\alpha}^{(p+1)}=+\infty\,,&& \text{and then}\qquad\lambda_{\alpha}^{(p)}(h)=0\,,\\
  \text{or}~b_{\alpha}^{*+1}<+\infty\,,&&\text{and then}\qquad
\lim_{h\to 0}-h\log \lambda_{\alpha}^{(p)}(h)=2(b_{\alpha}^{*+1}-a_{\alpha}^{*})\,.
\end{eqnarray*}
Obviously, the multiplicity of the $0$-eigenvalue of
$\Delta_{f,M,h}^{(p)}$\,, the dimension of its kernel, equals the
$p^{th}$ Betti number of $M$\,,
$\sharp\left\{\alpha\in A^{(p)}\,,
  b_{\alpha}^{(p+1)}=+\infty\right\}=\beta^{(p)}(M)$\,, and does not
depend on the function $f$\,.
\end{thm}
To summarize the situation, the logarithms  of  the exponentially
small eigenvalues of $\Delta_{f,M,h}^{(p)}$ are  given
by the lengths of the bars $b_{\alpha}^{*+1}-a_{\alpha}^{*}$ of which
one endpoint in $\rz$ is of degree $p$\,, the eigenvalues associated
with infinite lengths being identically $0$ for $h$ small enough.
A direct application of the stability results of persistent cohomology
then gives the variations of the exponentially small spectrum when the
function $f$ is perturbed.
\begin{cor}
\label{cor:mainsimple}
Assume that $f$ satisfies Hypothesis~\ref{hyp:mainf}, let ${\cal
  B}(f)=\left([a_{\alpha},b_{\alpha}[\right)_{\alpha\in A}$\,,
$a_{\alpha}<b_{\alpha}$\,, $A=\sqcup_{0\leq p\leq d}A^{(p)}$\,, 
  be a persistent bar code associated with
$f$\,, and set 
$\ell_{min}=\min\left\{b_{\alpha}-a_{\alpha}, \alpha \in A\right\}$\,.
For any other function $g$ which satisfy Hypothesis~\ref{hyp:mainf}
with $\|g-f\|_{{\cal C}^{0}}<\frac{\ell_{min}}{4}$\,, the 
$\tilde{O}(e^{-\frac{\ell_{min}}{h}})$ eigenvalues of
$\Delta_{g,M,h}^{(p)}$ can be labelled with multiplicities
$$
\lambda_{\alpha}(g,h),\ \  \alpha\in A^{(p)}\text{or}\; \alpha\in
A^{(p-1)}\,, b_{\alpha}^{(p)}\neq +\infty\,,
$$
with
\begin{eqnarray*}
  &&\lambda_{\alpha}(g,h)=0\quad \text{if}\quad b_{\alpha}^{(p+1)}=+\infty\\
\text{or}&&
2(b_{\alpha}-a_{\alpha})+4\|g-f\|_{\mathcal{C}^{0}}\geq
\lim_{h\to 0}-h\log(\lambda_{\alpha}(g,h))\geq
            2(b_{\alpha}-a_{\alpha})-4\|g-f\|_{{\cal C}^{0}}>\ell_{min}\,.
\end{eqnarray*}
\end{cor}
One rapidly realizes that we make no normal form assumption for
$f$ near the ``critical values'' 
$c_{1},\ldots c_{N_{f}}$ of $f$\,. Even if we work with  ${\cal
  C}^{\infty}$ functions,  any closed set $K$
 of $M$ can be the global
minimum of $f\in {\cal C}^{\infty}(M)$ by taking a non negative
 ${\cal C}^{\infty}$ function vanishing only on $K$ after Whitney's
 extension theorem. Having a finite number of critical values
 restricts the possible sets $K$ which still make a very large class.
Hence, no algebraic behaviour with respect to $h$ of the leading terms
of the subexponential factors 
can be expected as it is the case when $f$ is
assumed to be a Morse function. Theorem~\ref{th:mainsimple} simply says that exponentially
small eigenvalues and their exponential scales are given by the
algebraic topology without specifying a possible subexponential 
factor.
Among other results, we will obtain similar things for
$\Delta_{f,f^{-1}([a,b]),h}$\,, $-\infty\leq a<b\leq +\infty$\,,
$a,b\not\in\left\{c_{1},\ldots c_{N_{f}}\right\}$\,, when considering the proper
boundary conditions on $f^{-1}(\left\{a\right\})$ (Dirichlet type) and
$f^{-1}(\left\{b\right\})$ (Neumann type).
Actually, this leads
us to the presentation of our strategy which passes through local
problems on $\rz=f(M)$ and a recurrence argument on the number $N$ of
``critical values'' lying in $[a,b]$\,.

\subsection{Strategy and outline of the article}
\label{sec:strat}

Proving a result like Theorem~\ref{th:mainsimple}, even in this
simplified form, is a rather long process which is clearly split into
various steps.
\begin{itemize}
\item A general presentation of bar codes in persistent (co)homology
  as well as properties of Hodge Laplacians on weakly regular domains
  are recalled in Appendix~\ref{app:perscohom} and in
  Appendix~\ref{sec:abstHodge}.
\item In Section~\ref{sec:expdec}, we set up the functional analysis framework,
  the relevant boundary conditions for Witten Laplacians, the
  corresponding integration by parts formulas, as well as weighted
  integration techniques {\it \`a la} Agmon, in order to obtain exponential
  decay estimates. We especially consider self-adjoint realizations of
  Witten Laplacians $\Delta_{f,h}$ in the domain $f^{-1}([a,b])$ when
  $a<b$ are not critical values, always with Dirichlet boundary
  conditions along  $f^{-1}(\left\{a\right\})$\,, the lowest level
  set, and Neumann boundary conditions along $f^{-1}(\left\{b\right\})$\,, the
  highest level set. Those self-adjoint realizations will be denoted by
  $\Delta_{f,f^{-1}([a,b],h)}$\,, and possibly
  $\Delta_{f,f^{-1}([a,b]),h}^{(p)}$ when specifying the degree. 
Remember the intuitive picture for functions:
  Dirichlet (resp. Neumann) boundary conditions are associated with a potential
  $-\infty$ (resp. $+\infty$). Such boundary conditions are actually
  the natural ones in order to 
  avoid boundary layer phenomena along the boundaries in the
  spectral analysis.
For further applications, this analysis is done in
  a weak regularity framework, and the long series of works by
  Mitrea~{\it et~al.} were instrumental in setting up the proper
  framework. The end of this section gathers repeatedly used technical
  lemmas, deduced from the exponential decay and weighted resolvent
 estimates for boundary
  Witten Laplacians.
\item Once the geometrical issues in the weak regularity case are
  solved, the rest of the analysis becomes essentially one dimensional
  on $\rz\supset f(M)$\,, as suggested by the bar code structure. The first step consists
  in understanding what happens when there is a single critical value
  in the energy interval 
  $[a,b]$\,, $[a,b]\cap
  \left\{c_{1},\ldots,c_{N_{f}}\right\}=\left\{\tilde{c}_{1}\right\}$\,. In
  this specific case,
  the bar code of $f$ has no bar compactly included in
  $]a,b[$\,. Accordingly, $\Delta_{f,f^{-1}([a,b]),h}$ should not have any
    non zero exponentially small eigenvalue. This is the main result of
 Section~\ref{sec:exp0},
 formulated in Proposition~\ref{pr:exp0}, which states
  that all the $\tilde{o}(1)$-eigenvalues of
  $\Delta_{f,f^{-1}([a,b])h}^{(p)}$ are equal to~$0$\,. After
  preliminary notations related with variations of the min-max
  principle, the core of the proof is done in
  Subsection~\ref{sec:exp0}, and  follows in some sense Carleman's
  general scheme for
  uniqueness  results of PDE, along the energy interval $[a,b]\subset\rz$\,. Resolvent estimates and other corollaries
  are listed afterwards. Section~\ref{sec:expdec} and
  Proposition~\ref{pr:exp0} provide in particular the number of
  $\tilde{o}(1)$-eigenvalues of $\Delta_{f,f^{-1}([a,b]),h}^{(p)}$ 
counted with multiplicities in this setting: it equals the relative Betti number
$$\beta(f^{b},f^{a};\rz)=\dim \ker(\Delta_{0,f^{-1}([a,b]),1}^{(p)})=\dim\ker(\Delta_{f,f^{-1}([a,b]),h}^{(p)})\,.$$
\item Only in Section~\ref{sec:rough} really starts the relationship
  between the bar code  $\mathcal{B}_{f}$ of $f$ and the spectral
  properties of $\Delta_{f,f^{-1}([a,b]),h}$\,. It contains an
  enumeration of the  non zero $\tilde{o}(1)$-eigenvalues of
  $\Delta_{f,f^{-1}([a,b]),h}$ in terms of bars compactly embedded in
  $]a,b[$\,, while the dimension
  $\dim\ker(\Delta_{f,f^{-1}([a,b]),h}^{(p)})=\beta^{(p)}(f^{b},f^{a};\rz)$
  is also expressed in terms of $\mathcal{B}_{f}$\,.
This section ends with Proposition~\ref{pr:roughmino} which proves the
 rough lower bound $e^{-2\frac{b-a}{h}}$ 
for the non zero $\tilde{o}(1)$-eigenvalues of
$\Delta_{f,f^{-1}([a,b]),h}$ (see Proposition~\ref{pr:roughmino}).
\item An important step elucidated in \cite{HKN}, and used in many
  forthcoming articles, consisted in the trivial observation that
  the eigenvalues of $\Delta_{f,f^{-1}([a,b]),h}$\,, restricted to some
  spectral compact segment\,,  are the square of the singular values
of  the restricted
  differential $d_{f,f^{-1}([a,b]),h}$\,. Singular values are much
  more flexible spectral quantities than eigenvalues. One of their
  advantage is that, in many situations, the approximation errors appear as
  relative ones for all the singular values, a property which is not
  fulfilled by eigenvalues.  We gather several functional analysis
  preliminary results in Section~\ref{sec:singval}, which elaborates in
  a functional abstract setting how various matricial error estimates
  propagate nicely to singular values estimates. 
\item The core of the proof of Theorem~\ref{th:mainsimple} is done in
  Section~\ref{sec:accanN}. It is a rather sophisticated proof by
  induction on the number $N$ of critical values contained in the
  energy interval $[a,b]$\,,
  $\left\{c_{1},\ldots,c_{N_{f}}\right\}\cap
  [a,b]=\left\{\tilde{c}_{1},\ldots,\tilde{c}_{N}\right\}$\,. This
recurrence  is  initiated by Section~\ref{sec:exp0} for the
  case $N=1$\,. Although it contains several steps, the induction from
  $N$ to $N+1$ mimics in some way the proof of  Mayer-Vietoris'
  Theorem.   The main result of this section is Theorem~\ref{th:induc},
  which can be considered as the central result of this text, while
  Proposition~\ref{pr:exp0} proves the simplest non trivial particular
  case. This induction contains many intermediate results, 
 which lead in particular in  
 Section~\ref{sec:coroll} to 
 Theorems~\ref{th:specres} and~\ref{th:stab1}, which  generalize
 respectively Theorem~\ref{th:mainsimple} and
  Corollary~\ref{cor:mainsimple} to the boundary Witten Laplacian
  $\Delta_{f,f^{-1}([a,b]),h}$\,.
  
\item Section~\ref{sec:regass} is devoted to various
  generalizations of Theorem~\ref{th:induc} and of its spectral
  corollaries. The first one concerns results for some domains which
  are not bounded by level sets, e.g.  for (non necessarily) small deformations
  $N_{t}$ and $N_{n}$
  of the level sets $f^{-1}(\left\{a\right\})$ and
  $f^{-1}(\left\{b\right\})$ for which the conditions $\partial_{n}
  f\big|_{N_{t}}<0$ and $\partial_{n}f\big|_{N_{n}}>0$ are still valid,
  and for which all the conclusions of Theorem~\ref{th:induc} and of its
  corollaries hold true. The second generalization is about noncompact manifolds like $\rz^{d}$\,, for which the results still hold provided we make
   some assumptions on  $M$ and
  $f$ at infinity.
The most technical one concerns the extension to a general subanalytic
Lipschitz function on a real analytic manifold (see
Hypothesis~\ref{hyp:realana}). 
Even when $f$ is a subanalytic real Lipschitz
function, it is possible to define
self-adjoint realizations $\Delta_{f,f^{-1}([a,b]),h}$\,, critical
values and finite bar codes, but there is an extra difficulty to establish Agmon's type
estimates to accurately control the exponential decay estimates. This problem
is solved in Subsections~\ref{sec:moregenLipStrat} and
\ref{sec:moregenLipAgm} by modelling a collection of solutions to
Hamilton-Jacobi equations associated with some natural stratification of the
subanalytic graph of $f$ in $M\times \rz$\,. 

\item Finally, Section~\ref{sec:applications} answers precisely our
Main  Question in various explicit cases. We return to our results of
\cite{LNV}, where Morse functions with simple critical values (one
critical point for every critical value) were considered. It was too
rapidly conjectured in \cite{LNV} that some topological constant
$\kappa^{2}$ appearing in the
subexponential factor was equal to $1$\,. It is true in the case of oriented
surfaces (see \cite{Lep-2D}), but examples  are now provided with a constant
$\kappa^{2}$
 equal to any
$n^{2}$\,, $n\in \nz^{*}$\,, the first example with
$\kappa^{2}=4$ arising in the case of a Morse function on $\rz
P^{2}$\,. Additionally, in the case of Morse functions with multiple
critical values, the constant $\kappa$ has to be replaced by an
``incidence matrix'', $\bm\kappa$\,,
related with the bar code. Various examples, including non Morse functions, show that
the accurate computation of the prefactors now results from two well
separated analyses: one for the global topology of the sublevel sets
relying on the bar code, and
one for the local asymptotic expansions of Laplace integrals.
\end{itemize}

\section{Boundary Witten Laplacians}
\label{sec:expdec}
In this section we specify the domain of various operators
involved in our analysis and review the basic exponential decay
estimates. The general assumptions and notations have been set in
Subsection~\ref{sec:genass} and in particular the function $f$
satisfies Hypothesis~\ref{hyp:mainf} or Hypothesis~\ref{hyp:Lipbar}. We shall give the 
definition of Dirichlet and Neumann boundary conditions for Witten
Laplacians on strongly Lipschitz domains $\overline{\Omega}$\,. Most
of the time in the
sequel, these domains will be level set domains
$\overline{\Omega}=f^{-1}([a,b])$ with $a,b\not\in
\left\{c_{1},\ldots,c_{N_{f}}\right\}$\,. The required 
Agmon's type or exponential
decay estimates will be proved under
Hypothesis~\ref{hyp:mainf}. We are unable to prove these estimates in the general setting of Hypothesis ~\ref{hyp:Lipbar} but will prove them 
for subanalytic Lipschitz functions (see Subsection~\ref{sec:moregenLip}).

\subsection{Tangential and normal traces}
\label{sec:domains}

\subsubsection{Smooth case}
\label{sec:domainssmooth}

\begin{definition}
  \label{de:tannorm} Let $N\subset M$ is a ${\cal
    C}^{\infty}$-hypersurface of $M$\,, $n$ a unit normal
  vector and $n^{\flat}$  the associated
  covector, defined locally.
  When $\omega\in W^{s,2}(M;\Lambda T^{*}M)$\,, $s>\frac{1}{2}$\,, the
  tangential and normal traces denoted
  $\mathbf{t}_{N}\omega$ and $\mathbf{n}_{N}\omega$ are defined by
$$
\mathbf{t}_{N}\omega=\mathbf{i}_{n}(n^{\flat}\wedge \omega)\big|_{N}\quad
\text{and}\quad
\mathbf{n}_{N}=
n^{\flat}\wedge(\mathbf{i}_{n}\omega)
\big|_{N}\,.
$$
\end{definition}
Before we extend this definition to more singular forms,  let us make
explicit this definition in coordinate systems
(see e.g. \cite{Sch}):
\begin{itemize}
\item When $n$ is a  normalized  normal vector to
 $N$\,, any vector field in $X=T_{N}M$ can be
 decomposed into $X=X_{T}\oplus X_{n}n$\,. The traces
 $\mathbf{t}_{N}\omega$ and $\mathbf{n}_{N}\omega$ are then equal to 
$$
\mathbf{t}_{N}\omega(X_{1},\ldots X_{p})=\omega\big|_{N}(X_{1,T},\ldots,X_{p,T})
\quad \text{and}\quad \mathbf{n}_{N}\omega=\omega\big|_{N}-\mathbf{t}_{N}\omega\,.
$$
\item 
With local coordinates $(x^{1},\ldots,x^{d})=(x',x^{d})\in\rz^{d}$ in a neighborhood  $U_{x_{0}}^{M}$
in $M$ of $x_{0}\in N$\,,
 such that $N\cap U_{x_{0}}^{M}=\left\{(x',x^{d})\in
   U_{0}^{\rz^{d}}\,, x^{d}=0\right\}$\,,\, 
$
g=\sum_{ij<d}g_{i,j}(x',x^{d})dx^{i}dx^{j}+(dx^{d})^{2}$\,, 
$n=\frac{\partial}{\partial x^{d}}$ and $n^{\flat}=dx^{d}$\,, and when
a differential form is written
\begin{eqnarray*}
  &&\omega=\sum_{\sharp
  I'=p,\; d\not\in I'}\omega_{I'}(x',x^{d})dx^{I'}+
\sum_{\sharp J'=p-1,\; d\not\in J'}\omega_{J'}(x',x^{d})dx^{J'}\wedge
dx^{d}\,,
\\
\text{with}&&dx^{I}=dx^{i_{1}}\wedge \cdots\wedge dx^{i_{\sharp
              I}}\quad,\quad  i_{1}<\cdots<i_{\sharp I}\quad,\quad
I=\left\{i_{1},\ldots,i_{\sharp I}\right\}\,,
\end{eqnarray*}
the
tangential and normal traces are given by
$$
\mathbf{t}_{N}\omega=\sum_{\sharp
  I'=p,\; d\not\in I'}\omega_{I'}(x',0)dx^{I'}\quad\text{and}\quad
\mathbf{n}_{f=s}\omega=\sum_{\sharp J'=p-1,\; d\not\in J'}\omega_{J'}(x',0)dx^{J'}\wedge dx^{d}\,.
$$
\item From those formulas one gets at once $\star
\mathbf{t}_{N}=\mathbf{n}_{N}\star$\,, where $\star$ denotes the
Hodge $\star$ operator on $(M,g)$\,. The possible
orientation twist $\mathrm{or}_{M}$ is locally trivial so that the orientability of
$M$ is not required. 
\item  
When restricted to the  tangent space to $N$\,,
$\mathbf{t}_{N}\omega$ coincides with  $j_{N}^{*}\omega$ where
$j_{N}:N\to M$ is the natural imbedding. Therefore
$\mathbf{t}_{N}d=d\mathbf{t}_{N}$ and therefore
$\mathbf{n}_{N}d^{*}=d^{*}\mathbf{n}_{N}$\,. Note also that
$\mathbf{t}_{N}$ and $\mathbf{n}_{N}$ commute with multiplications
by functions.
\end{itemize}

\subsubsection{Lipschitz domains}
\label{sec:domainsLip}

The typical case which will be considered is when $N=\partial \Omega$
is the boundary of a Lipschitz domain of $M$ (strongly Lipschitz according to the terminology of \cite{GMM}).
This means that $\Omega$ is locally the hypograph of a Lipschitz
function in a proper coordinate system. For the notations, $\Omega$ is
an open domain in $M$ and its closed version is
$\overline{\Omega}=\Omega\sqcup N$ with $N=\partial \Omega$\,.
 Precisely we consider
the following situation.
\begin{hyp}
  \label{hyp:domain}
The domain  $\overline{\Omega}=\Omega\sqcup N\subset M$ is
a  
Lipschitz domain  with $N=N_{t}\sqcup N_{n}$ made of two disjoint
closed hypersurfaces.
\end{hyp}

\noindent
When $\Omega$ is  a regular domain, with $\mathcal{C}^{\infty}$
boundaries $N_{t}$ and $N_{n}$\,, the unit normal
vector field $n$ to
$N=\partial \Omega$ is globally defined so that the hypersurface measure
$d\sigma$\,, the orientation twist $\mathrm{or}_{N}$ and the Hodge
$\star$ operation on $N=\partial \Omega$ are deduced from
$d\mathrm{Vol}_{g}$ and $\mathrm{or}_{M}$ and the Hodge $\star$  on
$M$\,. In the general case when the domain $\Omega$ has only the
assumed Lipschitz regularity, the same things hold except that the
normal vector is defined $d\sigma$-almost everywhere along
$N_{t}\sqcup N_{n}$\,, $d\sigma$ being the
$\mathcal{H}^{d-1}$-Hausdorff measure.\\

For two forms $u,v\in W^{1,2}(\Omega, \Lambda T^{*}M)$\,,
the Green formula yields
\begin{align}
\langle du, v\rangle_{ L^2 ( \Omega )}- \langle u, d^*v\rangle_{ L^2 (
  \Omega )}&=\int_{\Omega}d(\overline{u}\wedge \star
             v)=\int_{N}\mathbf{t}_{N}(\overline{u}\wedge \star v)
\nonumber
\\
\label{eq:Greent}
&=\int_{N} \langle  u,\mathbf{i}_{n} v\rangle_{\Lambda T^*_{\sigma}M}
  d\sigma=\int_{N} \langle \mathbf{t}_{N} u,\mathbf{i}_{n}
  v\rangle_{\Lambda T^*_{\sigma}M} d\sigma,\\
\label{eq:Greenn}
&= \int_{N} \langle  n^{\flat} \wedge u,v\rangle_{\Lambda
  T^*_{\sigma}M} d\sigma= \int_{N} \langle  n^{\flat} \wedge u,
  \mathbf{n}_{N} v\rangle_{ \Lambda T^*_{\sigma}M} d\sigma
\end{align}
while the decomposition  $N=N_{t}\sqcup N_{n}$ into two
disjoint closed hypersurfaces clearly implies
\begin{eqnarray*}
  &&\left(\mathbf{t}_{N_{t}}u=0\right)\Leftrightarrow\left(\mathrm{supp}~n^{\flat}\wedge
     u\subset N_{n}\right)\\
\text{and}&&
\left(\mathbf{n}_{N_{n}}v=0\right)\Leftrightarrow
\left(\mathrm{supp}~\mathbf{i}_{n}v\subset N_{t}\right)\,.
\end{eqnarray*}
Moreover according to \cite{JMM},
when $\omega\in L^{2}(\Omega;\Lambda T^{*}M)$ and
$d\omega\in L^{2}(\Omega;\Lambda T^{*}M)$\,,
the above Green formulas provide the duality needed
 to define $n^{\flat}\wedge 
\omega\big|_{N}\in  W^{-\frac12,2}(N;\Lambda T^{*}M)$
  by
\begin{equation}
\label{eq.traceT}
\forall g\in W^{\frac{1}{2},2}(N)\,,\quad
\langle  n^{\flat}\wedge \omega\,,\,g\rangle_{W^{-\frac{1}{2},2}(N),W^{\frac{1}{2},2}(N)}
=\langle d\omega\,,\,
G\rangle_{L^{2}(\Omega)}-\langle \omega\,,\, d^{*}G\rangle_{L^{2}(\Omega)}\,,
\end{equation} 
where $G$ is any form in $W^{1,2}(\Omega;\Lambda T^{*}M)$ such that $G\big|_{N}=g\in
W^{\frac12,2}(N;\Lambda T^{*}M)$\,.
Similarly, 
when $\omega$ and
$d^{*}\omega$ belong to $L^{2}(\Omega;\Lambda T^{*}M)$\,,
one can define $\mathbf{i}_{n}\omega\big|_{N}\in  W^{-\frac12,2}(N;\Lambda T^{*}M)$
by
\begin{equation}
\label{eq.traceN}
\forall g\in W^{\frac{1}{2},2}(N)\,,\quad
\langle  \mathbf{i}_{n}\omega\,,\,g\rangle_{W^{-\frac{1}{2},2}(N),W^{\frac{1}{2},2}(N)}
=\langle \omega\,,\,
d G\rangle_{L^{2}(\Omega)}-\langle d^{*}\omega\,,\,  G\rangle_{L^{2}(\Omega)}\,.
\end{equation}
In particular, when $O$ is an open subset of $N$ and
when the trace 
$n^{\flat}\wedge \omega\big|_{O}$
defined in the sense of \eqref{eq.traceT}
(resp. of \eqref{eq.traceN})
belongs to $L^{2}(O;\Lambda T^{*}M)$\,,
the
  tangential (resp. normal) trace
  $\mathbf{t}_{O}\omega$ (resp. $\mathbf{n}_{O}\omega$)  is  well defined
  on $O$ 
  by the standard formula  from Definition~\ref{de:tannorm}:
$$
\mathbf{t}_{O}\omega=\mathbf{i}_{n}(n^{\flat}\wedge \omega)\big|_{O}\quad
\text{\big(\,resp.}\quad
\mathbf{n}_{O}=
n^{\flat}\wedge(\mathbf{i}_{n}\omega)
\big|_{O}\,\big)\,.
$$
We may thus make sense of
the boundary 
condition $\mathbf{t}_{N_{t}}\omega=0$
(resp. $\mathbf{n}_{N_{n}}\omega =0$),
which is equivalent to 
$\mathrm{supp}~n^{\flat}\wedge \omega\big|_{N}\subset N_{n}$
(resp. $\mathrm{supp}~\mathbf{i}_{n}\omega\big|_{N} \subset N_{t}$),
 for any $\omega\in L^{2}(\Omega;\Lambda T^{*}M)$ such that
$d\omega\in L^{2}(\Omega;\Lambda T^{*}M)$ (resp. $d^{*}\omega\in
L^{2}(\Omega;\Lambda T^{*}M)$).\\

\noindent
According to \cite[Proposition~3.1]{JMM},  ${\cal C}^{\infty}_{0}(\Omega\cup N_{n};\Lambda T^{*}M)$ (resp. ${\cal
  C}^{\infty}_{0}(\Omega\cup N_{t};\Lambda T^{*}M)$) is dense in 
\begin{eqnarray}
\label{eq.calT}
  && {\cal T}=\left\{\omega\in L^{2}(\Omega;\Lambda T^{*}M)\,,
     d\omega\in L^{2}(\Omega;\Lambda T^{*}M)\,, \mathbf{t}_{N_{t}}\omega=0\right\}\\
     \label{eq.calN}
\text{\big(\,resp. in}&&
{\cal N}=\left\{\omega\in L^{2}(\Omega;\Lambda T^{*}M)\,,
     d^{*}\omega\in L^{2}(\Omega;\Lambda T^{*}M)\,, \mathbf{n}_{N_{n}}\omega=0\right\}\,\big)
\end{eqnarray}
endowed with the norm $\|\omega\|_{L^{2}(\Omega)}+\|d\omega\|_{L^{2}(\Omega)}$
(resp. $\|\omega\|_{L^{2}(\Omega)}+\|d^{*}\omega\|_{L^{2}(\Omega)}$).
\noindent
Theorem~3.4 of \cite{JMM} also says that when
 $u,v\in L^{2}(\Omega)$ with $du\in  L^{2}(\Omega;\Lambda^{p+1} T^{*}M)$\,, $d^{*}v\in L^{2}(\Omega;\Lambda^{p} T^{*}M)$\,,
 and
 $$
 \supp \mathbf{i}_{n} v\subset \Gamma \quad \text{or}\quad \supp(n^{\flat} \wedge u)\subset\Gamma
 $$
 with $\Gamma=N_{t}$ or $\Gamma=N_{n}$\,,
the following Green formulas
\begin{equation}
\label{eq.gen-Green}
\begin{aligned}
\langle du, v\rangle_{ L^2 ( \Omega )}-\langle u, d^*v\rangle_{
  L^2 (  \Omega )}  &=  \int_{\Gamma}\langle
 n^{\flat} \wedge u,n^{\flat} \wedge(\mathbf{i}_{n} v)
\rangle_{T^*_{\sigma}\Omega}d\sigma\\
&  = \int_{\Gamma}\langle
\mathbf{i}_{n} (n^{\flat} \wedge u),\mathbf{i}_{n} v
\rangle_{T^*_{\sigma}\Omega} d \sigma
\end{aligned}
\end{equation}
make sense with a r.h.s. interpreted in general in a weak form specified in 
\cite[Proposition~3.3]{JMM}. Notice that under Hypothesis~\ref{hyp:domain},
the geometric assumptions concerned with $\Gamma$ in \cite{JMM}
are trivially satisfied without any locally mixed boundary conditions. Additionally, when $ \mathbf{i}_{n} v$ and $n^{\flat} \wedge u$
belong 
to $L^{2}(\Gamma)$\,, 
the r.h.s. of \eqref{eq.gen-Green} are standard  integrals along the
boundary.
\begin{definition}
\label{de:defW}
Let $\Omega$ be a Lipschitz domain of $M$ with
$\overline{\Omega}=\Omega\sqcup N$\,, $N=N_{t}\sqcup N_{n}$ like above, and let
$\mathcal{T},\mathcal{N}$ be the spaces defined in \eqref{eq.calT}\eqref{eq.calN}.\\
The space
$$
W(\Omega;\Lambda T^{*}M)=\left\{\omega\in L^{2}(\Omega;\Lambda
  T^{*}M); d\omega \in L^{2}(\Omega;\Lambda T^{*}M); d^{*}\omega\in
  L^{2}(\Omega;\Lambda T^{*}M)\right\}
$$
is endowed with its natural Hilbert space norm given by
\begin{equation}
\label{eq.W}
\|\omega\|_{W(\Omega)}^{2}\ :=\ \|\omega\|^{2}_{L^{2}(\Omega)}+\|d\omega\|^{2}_{L^{2}(\Omega)}+\|d^{*}\omega\|^{2}_{L^{2}(\Omega)}\,.
\end{equation}
The closed subspace $\mathcal{T}\cap \mathcal{N}$ of $W(\Omega;\Lambda
T^{*}M)$ will be denoted $W_{\partial}(\Omega;\Lambda T^{*}M)$ and the
restriction of the $W(\Omega;\Lambda T^{*}M)$-norm $\|~\|_{W_{\partial}(\Omega)}$\,.
\end{definition}

\begin{remark}
\label{re.Lip}
\begin{itemize}
\item[i)]
By interior elliptic regularity, note that
$$ W_{\partial}(\Omega; \Lambda T^{*}M)\subset 
W(\Omega;\Lambda T^{*}M)\subset
 W_{\text{loc}}^{1,2}(\Omega;\Lambda
  T^{*}M)$$ with continuous embeddings.
However it is known that $W(\Omega;\Lambda T^{*}M)$\,, and even 
$W_{\partial}(\Omega;\Lambda T^{*}M)$ if we add boundary conditions,
differs from 
$W^{1,2}(\Omega;\Lambda  T^{*}M)$ 
for a general Lipschitz domain (see e.g. \cite{MTV}\cite{MMMT}). 
An easy
counter example is
$u=r^{\frac{\pi}{\theta_{0}}-1}\cos(\frac{\pi}{\theta_{0}}\theta)dr-r^{\frac{\pi}{\theta_{0}}-1}\sin(\frac{\pi}{\theta_{0}}\theta)d\theta$
 in the sector $0<\theta<\theta_{0}$ of $\rz^{2}$
near $r=0$\,. It satisfies $\mathbf{n}u=0$\,, $du=0$ and $d^{*}u\in
L^{2}$ near $r=0$ while $u\not\in W^{1,2}$ near $r=0$ when
$\theta_{0}>\pi$\,.
\item[ii)] The space
$W(\Omega;\Lambda T^{*}M)$ and its subspace
$W_{\partial}(\Omega;\Lambda T^{*}M)$
are Lipschitz-module:
for any
$\varphi\in W^{1,\infty}(\Omega;\mathbb R)$ and 
 $\omega\in  W(\Omega;\Lambda T^{*}M)$\,,
 $\varphi \omega$ belongs to $W(\Omega;\Lambda T^{*}M)$ and
the mapping 
 $\omega\in  W(\Omega;\Lambda T^{*}M)\mapsto \varphi \omega \in W(\Omega;\Lambda T^{*}M) $
 is continuous. Moreover, for any  bounded sequence $(\varphi_{n})_{n\in\N}$ of $W^{1,\infty}(\Omega;\mathbb R)$
 such that $\varphi_{n}\to \varphi $ a.e. and $ d\varphi_{n}\to d\varphi $ a.e.,
 the convergence $\varphi_{n}\omega \to \varphi\omega$
 holds for the $W(\Omega;\Lambda T^{*}M)$-norm  for every $\omega\in
 W(\Omega;\Lambda T^{*}M)$\,.
\item[iii)] In our case it is proven in \cite{MMMT} and it is extended
  in \cite{JMM} that $W_{\partial}(\Omega;\Lambda T^{*}M)$ is embedded
  in $W^{1/2,2}(\Lambda T^{*}M)$\,. Again the exponent $\frac{1}{2}$
  cannot be improved for a general strongly Lipschitz domain
  $\Omega$\,.
\item[iv)] For a different approach on regularity issues for Lipschitz
  domains and relying on a generalization of Bogovski\u{\i} and Poincar{\'e}
  type integrals, we refer to \cite{CoMcI}, \cite{Mit} and \cite{MiMo}.
 \end{itemize}
\end{remark}

\begin{prop}
\label{pr:exttn0}
Let $W_{\partial}(\Omega;\Lambda T^{*}M)$ be the space of Definition~\ref{de:defW}.\\
Every $\omega\in W_{\partial}(\Omega;\Lambda T^{*}M) $ belongs to $W^{\frac12,2}(\Omega;\Lambda T^{*}M)$
and has, in the sense of \eqref{eq.traceT} and \eqref{eq.traceN},  tangential and normal traces $\mathbf{t}_{N}\omega$ and $\mathbf{n}_{N}\omega$
which actually belong to $L^{2}(N;\Lambda T^{*}M)$\,. Moreover,  there exists
 $C>0$ such that
$$
\forall \omega\in W_{\partial}(\Omega;\Lambda T^{*}M)\,,\quad
\|\omega\|_{W^{\frac12,2}(\Omega)}^{2}+
\|\omega|_{N}\|_{L^{2}(N)}^{2}
\ \leq \ 
C \|\omega\|^{2}_{W_{\partial}(\Omega)}\,,
$$
where   $\omega|_{N}:=\mathbf{t}_{N}\omega+ \mathbf{n}_{N}\omega \in L^{2}(N;\Lambda T^{*}M)$
is the total trace of $\omega$\,.\\
Finally, in the case where  $\overline{\Omega}$ is a smooth domain, 
Gaffney's inequality holds:
$$
W_{\partial}(\Omega;\Lambda T^{*}M)\ =\ \left\{\omega\in W^{1,2}(\Omega;\Lambda
  T^{*}M),\, \mathbf{t}_{N_{t}}\omega=0\,,\, \mathbf{n}_{N_{n}}\omega=0\right\}
$$
and  there
exists $C\geq 1$ such that
$$
\forall \omega\in W_{\partial}(\Omega;\Lambda T^{*}M)\,,\quad
C^{-1}\|\omega\|_{W^{1,2}(\Omega)}^{2}\ \leq\ 
 \|\omega\|_{W_{\partial}(\Omega)}
 \ \leq\ 
C\|\omega\|_{W^{1,2}(\Omega)}^{2}\,.
$$
\end{prop}
\begin{proof}
The first part of the statement is an immediate consequence of the analysis 
led in \cite{JMM} (see e.g. Theorem~1.1 there),
but our setting is actually simpler since no
locally mixed boundary conditions appear.\\
For Gaffney's inequality when the domain $\overline \Omega$ is smooth, 
consider first 
$$\omega\in W'(\Omega;\Lambda T^{*}M)\ :=\ 
\left\{u\in 
W^{1,2}(\Omega;\Lambda
  T^{*}M),\, \mathbf{t}_{N_{t}}u=0\,,\, \mathbf{n}_{N_{n}}u=0\right\}$$
and  a  function $\chi\in {\cal
  C}^{\infty}_{0}(\Omega\cup N_{t};[0,1])$ such that $\chi\equiv 1$ in
a neighborhood of $N_{t}$\,, 
 and decompose  $\omega$ as $\omega=\chi \omega +(1-\chi)\omega=\omega_{1}+\omega_{2}$\,. 
 For any
differential operator $L$ of order $\leq1$\,, note then the relation $\|L\omega_{j}\|_{L^{2}}\leq
C_{\chi,L,j}\|\omega\|_{L^{2}}+\|L\omega\|_{L^{2}}$\,, $j=1,2$\,. Now,
$\omega_{1}=\chi \omega \in W^{1,2}(\Omega;\Lambda T^{*}M)$ satisfies  $\mathbf{t}_{\partial
  \Omega}\omega_{1}=0$ and $\omega_{2}=(1-\chi) \omega \in W^{1,2}(\Omega;\Lambda T^{*}M)$
satisfies $\mathbf{n}_{\partial \Omega}\omega_{2}=0$\,. 
Gaffney's inequality for
Dirichlet boundary conditions then  says
$$
\|\omega_{1}\|_{W^{1,2}}^{2}\leq C_{1}\left[\|\omega_{1}\|_{L^{2}}^{2}+\|d\omega_{1}\|_{L^{2}}^{2}+\|d^{*}\omega_{1}\|_{L^{2}}^{2}\right]
$$
for some $C_{1}$ independent of $\omega_{1}$\,,
while Gaffney's inequality  for Neumann boundary conditions
  says 
$$
\|\omega_{2}\|_{W^{1,2}}^{2}\leq C_{2}\left[\|\omega_{2}\|_{L^{2}}^{2}+\|d\omega_{2}\|_{L^{2}}^{2}+\|d^{*}\omega_{2}\|_{L^{2}}^{2}\right]
$$
for some $C_{2}$ independent of $\omega_{2}$ (these two different boundary conditions have been treated separately
in \cite{Sch}). Adding the above two inequalities then leads 
to 
$$
\forall \omega\in W'(\Omega;\Lambda T^{*}M)\,,\ \ \ 
\|\omega\|_{W^{1,2}}^{2}\leq C\left[\|\omega\|_{L^{2}}^{2}+\|d\omega\|_{L^{2}}^{2}+\|d^{*}\omega\|_{L^{2}}^{2}\right]\,.
$$
In order to achieve the proof of 
Proposition~\ref{pr:exttn0}, it then suffices to show that 
$W'(\Omega;\Lambda T^{*}M)$ equals $W_{\partial}(\Omega;\Lambda
T^{*}M)$\,. We can forget the boundary conditions. With a regular
boundary, a simple local reflexion after identifying the domain with
a half space, leads to the problem on $\rz^{d}$ with a Lipschitz
riemannian metric,  asking if a compactly supported form in $\omega\in L^{2}_{comp}(\rz^{d})$
such $d\omega\in L^{2}(\rz^{d})$ and $d^{*}\omega \in
L^{2}(\rz^{d})$ belongs to $H^{1}_{comp}(\rz^{d})$\,. It is a
straightforward application of Lax-Milgram's theorem.
\end{proof}

\subsection{Witten's deformation}
\label{sec:wittdeff}
The function $f$ is assumed to be a Lipschitz
function and   the domain
$\overline{\Omega}$ satisfies  Hypothesis~\ref{hyp:domain}. 
Improved regularity results  are stated when $f$ and
$\overline{\Omega}$ are more regular.
\begin{definition}
\label{de:domain} 
Assume $f\in W^{1,\infty}(M;\rz)$\,, $h>0$\,,  and
Hypothesis~\ref{hyp:domain} for $\overline{\Omega}=\Omega\sqcup
N=\Omega\sqcup N_{t}\sqcup N_{n}$\,. The operators
$d_{f,\overline{\Omega},h}$ and $d^{*}_{f,\overline{\Omega},h}$ are
defined by
\begin{eqnarray*}
 && D(d_{f,\overline{\Omega},h})\ :=\ \left\{\omega\in L^{2}(\Omega;\Lambda
  T^{*}M)\,,\, d_{f,h}\omega\in L^{2}(\Omega;\Lambda T^{*}M)\,,
    \mathbf{t}_{N_{t}}\omega=0\right\}\ =\ {\cal T}
\\
\text{and}&&D(d_{f,\overline{\Omega},h}^{*})\ :=\ \left\{\omega\in L^{2}(\Omega;\Lambda
  T^{*}M)\,,\, d_{f,h}^{*}\omega\in L^{2}(\Omega;\Lambda T^{*}M)\,,
    \mathbf{n}_{N_{n}}\omega=0\right\}
    \ =\ {\cal N}\,,
\end{eqnarray*}
where ${\cal T}$ and ${\cal N}$ are the spaces defined in \eqref{eq.calT}
and \eqref{eq.calN},
and we recall that
$$
d_{f,h}=e^{-\frac{f}{h}}(hd)e^{\frac{f}{h}}=hd
     +df\wedge\quad\text{and}\quad
d_{f,h}^{*}=e^{\frac{f}{h}}(hd^{*})e^{-\frac{f}{h}}=hd^{*}+\mathbf{i}_{\nabla
   f}
$$
according to 
\eqref{eq:dfh} and \eqref{eq:dfhstar}.
\end{definition}
\noindent A particular case that we will study extensively is when
$\overline{\Omega}=f^{-1}([a,b])$\,, $N_{t}=f^{-1}(\left\{a\right\})$\,,
$N_{n}=f^{-1}(\left\{b\right\})$\,,
 and $a<b$ do not belong  to $\left\{c_{1},\ldots,
  c_{N_{f}}\right\}$ under Hypothesis~\ref{hyp:mainf}
 (in this case  $\Omega=f_{a}^{b}$ according to
Definition~\ref{de:fab}). With such an
$f$-dependent domain, it will be useful to consider
$d_{0,f^{-1}([a,b]),h}$ and $d_{f,f^{-1}([a,b]),h}$\,.

\begin{prop}
\label{pr:Dom-d}
In the framework of Definition~\ref{de:domain}, the operator
$d_{f,\overline{\Omega},h}$ (resp. $d^{*}_{f,\overline{\Omega},h}$) is
densely defined, closed, and $\mathrm{Ran}~
d_{f,\overline{\Omega},h}\subset \ker d_{f,\overline{\Omega},h}$
(resp. $\mathrm{Ran}~d_{f,\overline{\Omega},h}^{*}\subset \ker
d_{f,\overline{\Omega},h}^{*}$)\,. Its adjoint is
$d_{f,\overline{\Omega},h}^{*}$
(resp. $d_{f,\overline{\Omega},h}$). The subspace ${\cal
  C}^{\infty}_{0}(\Omega\cup N_{n};\Lambda T^{*}M)$ (resp. ${\cal
  C}^{\infty}_{0}(\Omega\cup N_{t};\Lambda T^{*}M)$) is dense in
$D(d_{f,\overline{\Omega},h})$ (resp. $D(d_{f,\overline{\Omega},h}^{*})$).
Finally, the identity
$$
D(d_{f,\overline{\Omega},h})\cap
D(d_{f,\overline{\Omega},h}^{*})\ =\ W_{\partial}(\Omega;\Lambda T^{*}M)\,,
$$
holds true when $W_{\partial}(\Omega;\Lambda T^{*}M)$ is the space of Definition~\ref{de:defW}.
\end{prop}

\begin{proof}
The operators $d_{f,\overline{\Omega},h}$ and
$d_{f,\overline{\Omega},h}^{*}$ having respective domains
${\cal T}$ and ${\cal N}$\,, with $\mathcal{T}\cap
\mathcal{N}=W_{\partial}(\Omega;\Lambda T^{*}M)$ by Definition~\ref{de:defW},
they are clearly densely defined, and they
are  bounded perturbations of $hd_{0,\overline{\Omega},1}$
and $hd^{*}_{0,\overline{\Omega},1}$ owing to $d_{f,h}=hd+df\wedge$
and $d_{f,h}^{*}=hd_{0,\overline{\Omega},1}+\mathbf{i}_{\nabla f}$\,.
The operators $d_{0,\overline{\Omega},1}$ and
$d_{0,\overline{\Omega},1}^{*}$ are moreover closed  with the density
properties, 
according to the presentation around \eqref{eq.traceT}--\eqref{eq.calN}.\\
As bounded perturbations,  the
  adjoint of $d_{f,\overline{\Omega},h}$ equals
  $d_{f,\overline{\Omega},h}^{*}$ because the adjoint of
  $d_{0,\overline{\Omega},1}$ is $d_{0,\overline{\Omega},1}^{*}$
while the adjoint of the bounded perturbation
  $df\wedge$ is $\mathbf{i}_{\nabla f}$\,. Actually $\omega$ belongs
  to the domain of the adjoint of 
  $d_{0,\overline{\Omega},1}$  iff
$$
\exists C>0\,,\ \forall u\in {\cal C}^{\infty}_{0}(\Omega\cup N_{n};\Lambda
T^{*}M)\,,\quad
|\langle du\,,\, \omega\rangle|\leq C\|u\|_{L^{2}}\,.
$$
Taking any $u\in {\cal C}^{\infty}_{0}(\Omega; \Lambda T^{*}M)$ implies
$d^{*}\omega\in L^{2}(\Omega;\Lambda T^{*}M)$ and therefore
$\mathbf{i}_{n}\omega\big|_{N}$ is well defined in
$W^{-1/2,2}(N;\Lambda T^{*}M)$\,. Using afterwards Green's formula
\eqref{eq.traceN} with a 
general $u\in {\cal
  C}^{\infty}_{0}(\Omega\cup N_{n};\Lambda T^{*}M)$ 
leads to $\mathbf{i}_{n}\omega\big|_{N_{n}}=0$\,. 
Thus the domain of
the  adjoint of
$d_{0,\overline{\Omega},1}$  is included in
$D(d^{*}_{0,\overline{\Omega},1})$\,, which is enough to
conclude.\\
It remains to check  $\mathrm{Ran}~d_{f,\overline{\Omega},h}\subset
\ker d_{f,\overline{\Omega},h} $ and
$\mathrm{Ran}~d_{f,\overline{\Omega},h}^{*}\subset \ker d_{f,\overline{\Omega},h}^{*}$\,.
The identities  \eqref{eq:dfh} and \eqref{eq:dfhstar} already
say that $d_{f,h}d_{f,\overline{\Omega},h}\omega=0$ in ${\cal
  D}^{'}(\Omega\,, \Lambda T^{*}M)$
(resp. $d_{f,h}^{*}d_{f,\overline{\Omega},h}^{*}\omega=0$) when
$\omega\in D(d_{f,\overline{\Omega},h})$ (resp. $\omega\in D(
d_{f,\overline{\Omega},h}^{*})$)\,. We can conclude that
$d_{f,\overline{\Omega},h}\omega\in \ker d_{f,\overline{\Omega},h}$
(resp. $d_{f,\overline{\Omega},h}^{*}\omega\in \ker
d^{*}_{f,\overline{\Omega},h}$)
if $\mathbf{t}_{N_{t}}d_{f,h}\omega=0$
(resp. $\mathbf{n}_{N_{n}}d_{f,h}^{*}\omega=0$) or more precisely,
with the weak formulation of the trace 
defined in \eqref{eq.traceT} (resp. in \eqref{eq.traceN}),
if
$\mathrm{supp}~n^{\flat}\wedge (d_{f,h}\omega)\big|_{N}\subset N_{n}$
(resp.
$\mathrm{supp}~\mathbf{i}_{n}(d_{f,h}^{*})\omega\big|_{N}\subset
N_{t}$).
For $\omega\in {\cal C}^{\infty}_{0}(\Omega\cup
N_{n};\Lambda T^{*}M)$ (resp. $\omega\in {\cal
  C}^{\infty}_{0}(\Omega\cup N_{t};\Lambda T^{*}M)$) 
the weakly defined trace $n^{\flat}\wedge (d_{f,h}\omega)\big|_{N_{t}}$ 
(resp. $\mathbf{i}_{n}(d_{f,h}^{*}\omega)\big|_{N_{n}}$) obviously vanishes because $N_{t}\cap \mathrm{supp}~d_{f,h}\omega
=\emptyset$ (resp. $N_{n}\cap \mathrm{supp}~d_{f,h}^{*}\omega=\emptyset$)\,. By the density of ${\cal C}^{\infty}_{0}(\Omega\cup
N_{n};\Lambda T^{*}M)$ (resp. ${\cal C}^{\infty}(\Omega\cup
N_{t};\Lambda T^{*}M)$) in $D(d_{f,\overline{\Omega},h})$
(resp. $D(d_{f,\overline{\Omega},h}^{*})$), we deduce
\begin{eqnarray*}
  && \forall \omega\in D(d_{f,\overline{\Omega},h})\,, \quad n^{\flat }\wedge
     (d_{f,\overline{\Omega},h}\omega)\big|_{N_{t}}=0\quad\text{in}~W^{-1/2,2}(N_{t})
\\
\big(\, \text{resp.}&&\forall
  \omega\in D(d_{f,\overline{\Omega},h}^{*})\,,\quad\mathbf{i}_{n}d_{f,\overline{\Omega},h}^{*}\omega\big|_{N_{n}}=0\quad\text{in}~W^{-1/2,2}(N_{n})\, \big)\,.
\end{eqnarray*}
This ends the proof.
\end{proof}
We now apply  results of the abstract Hodge theory reviewed in
Appendix~\ref{sec:abstHodge} to  our specific framework.
\begin{prop}
\label{pr:domain}
Assume Hypothesis~\ref{hyp:domain} for 
$\overline{\Omega}=\Omega\sqcup N_{t}\sqcup N_{n}$\,,  $f\in
W^{1,\infty}(\Omega;\rz)$ and let
$W_{\partial}(\Omega;\Lambda T^{*}M)$ be the space of Definition~\ref{de:defW}\,.
\begin{enumerate}
\item The operator $d_{f,\overline{\Omega},h}+d_{f,\overline{\Omega},h}^{*}$ with domain
$$
D(d_{f,\overline{\Omega},h})\cap
D(d_{f,\overline{\Omega},h}^{*})\ =\ W_{\partial}(\Omega;\Lambda T^{*}M)
$$
is self-adjoint and has a
  compact resolvent.
\item The operator $\Delta_{f,\overline{\Omega},h}:=
d_{f,\overline{\Omega},h}d^{*}_{f,\overline{\Omega},h}+d^{*}_{f,\overline{\Omega},h}d_{f,\overline{\Omega},h}$ with domain
$$
D(\Delta_{f,\overline{\Omega},h})  = \{u\in D(d_{f,\overline{\Omega},h})\cap D(d_{f,\overline{\Omega},h}^{*})\ \text{s.t.}\  d_{f,h}
u\in D(d^{*}_{f,\overline{\Omega},h})
\ \text{and}\ d^{*}_{f,h}
u\in D(d_{f,\overline{\Omega},h})\}$$
is a self-adjoint operator
  with a compact resolvent. It is  the Friedrichs
  extension associated with  the (closed) quadratic form
  $Q_{f,\overline{\Omega},h}(\omega)=\|d_{f,h}\omega\|_{L^{2}}^{2}+\|d_{f,h}^{*}\omega\|_{L^{2}}^{2}$
  with domain $D(d_{f,\overline{\Omega},h})\cap
  D(d_{f,\overline{\Omega},h}^{*})$\,.
\item The ranges of $d_{f,\overline{\Omega},h}$ and $d_{f,\overline{\Omega},h}^{*}$ are
  closed and the following Hodge decompositions hold in $L^{2}$:
$$
L^{2}(\Omega;\Lambda T^{*}M)=
\underbrace{\Ran( d_{f,\overline{\Omega},h})
 \mathop{\oplus}^{\perp}\ker(\Delta_{f,\overline{\Omega},h})
}_{\ker(d_{f,\overline{\Omega},h})} 
\hspace{-1.8cm}\overbrace{\phantom{\ker(\Delta_{f,\overline{\Omega},h})}
\mathop{\oplus}^{\perp}\Ran(
    d_{f,\overline{\Omega},h}^{*})}^{\ker{(d_{f,\overline{\Omega},h}^{*})}}
$$
\item For any $z\in \cz\setminus \sigma(\Delta_{f,\overline{\Omega},h})$\,, one has
for any compactly supported  and bounded
  measurable function  $\chi$ on $\rz$
  and for any $\omega\in D(\mathbf{d})$\,,
where $\mathbf{d}=d_{f,\overline{\Omega},h}$ or  $\mathbf{d}=d_{f,\overline{\Omega},h}^{*}$\,,
  \begin{align*}
\mathbf{d}(z-\Delta_{f,\overline{\Omega},h})^{-1}\omega
=(z-\Delta_{f,\overline{\Omega},h})^{-1}\mathbf{d}\omega
\quad\text{and}\quad
\mathbf{d}\circ \chi(\Delta_{f,\overline{\Omega},h})\omega
=\chi(\Delta_{f,\overline{\Omega},h})\circ \mathbf{d}\omega\,.
  \end{align*}
\item
When $\overline \Omega$ is smooth and $f\in {\cal
    C}^{2}(\overline{\Omega};\rz)$\,, the domain of
>   $\Delta_{f,\overline{\Omega},h}$ equals
$$D(\Delta_{f,\overline{\Omega},h})=\left\{\omega\in W^{2,2}(\Omega;\Lambda
  T^{*}M)\,,
                           \begin{array}[c]{ll}
                              \mathbf{t}_{N_{t}}\omega=0\,,&
                             \mathbf{n}_{N_{n}}\omega=0\,,\\
   		\mathbf{t}_{N_{t}}d_{f,h}^{*}\omega=0\,,&
                          \mathbf{n}_{N_{n}}d_{f,h}\omega=0
                           \end{array}
 \right\} \,. $$
  \end{enumerate}
\end{prop}
\begin{proof}
The identification of $D(d_{f,\overline{\Omega},h})\cap
D(d_{f,\overline{\Omega},h}^{*})$ is done in
Proposition~\ref{pr:Dom-d}. The statements 1), 2), 3) are then
straightforward applications of  Proposition~\ref{pr.abs-hodge} in Appendix~\ref{sec:abstHodge}.
The first identity of the statement 4) 
is an application of the general relation 
\eqref{eq.Res-commut} in Appendix~\ref{sec:abstHodge}. 
The second identity then
comes from
the functional calculus for self-adjoint operators.\\ 
Finally, for 5), it suffices to notice that
$\Delta_{f,h}=-h^{2}\Delta_{0,1}+|\nabla f|^{2}+h({\cal L}_{\nabla
  f}+{\cal L}_{\nabla f}^{*})$ and that
$\Delta_{f,\overline{\Omega},h}$ is a bounded perturbation of
$h^{2}\Delta_{0,\overline{\Omega},1}$ when
$f\in {\cal
    C}^{2}(\overline{\Omega};\rz)$\,.
But  the elliptic analysis made in
\cite{Sch}\cite{MMMT} (see also \cite{LNV} for the combination of Dirichlet on
$N_{t}$ and Neumann on $N_{n}$ boundary conditions)  ensures that the
domain of $\Delta_{0,\overline{\Omega},1}=dd^{*}+d^{*}d$ is
$$
D(\Delta_{0,\overline{\Omega},1})=
\left\{\omega\in W^{2,2}(\Omega;\Lambda
  T^{*}M)\,,
                           \begin{array}[c]{ll}
                              \mathbf{t}_{N_{t}}\omega=0\,,&
                             \mathbf{n}_{N_{n}}\omega=0\,,\\
   		\mathbf{t}_{N_{t}}d_{f,h}^{*}\omega=0\,,&
                          \mathbf{n}_{N_{n}}d_{f,h}\omega=0
                           \end{array}
 \right\} \,.
$$
\end{proof} 
\begin{remark}
\label{re:rem1}
  Let us complete the statements of Propositions~\ref{pr:Dom-d} and~\ref{pr:domain} with
  some remarks when $f$ satisfies Hypothesis~\ref{hyp:mainf} or Hypothesis~\ref{hyp:Lipbar}.
  \begin{itemize}
  \item The domain $D(d_{f,\overline{\Omega},h})$ does not contain
 any  other regularity
  assumption
    than $\omega\in L^{2}(\Omega)$\,, $d_{f,h}\omega\in L^{2}(\Omega)$\,, and
    does not  contain any
    condition on $N_{n}$\,. 
 In particular, when $a'\leq a<b$ do not belong to
 $\left\{c_{1},\ldots, c_{N_{f}}\right\}$ according to
 Hypothesis~\ref{hyp:mainf}\,, the domain $\underline{f_{a}^{b}}$
 (resp. 
 $\underline{f_{a'}^{b}}$)  equals to $f^{-1}([a,b])$ and satisfies
 Hypothesis~\ref{hyp:domain} with $N_{t}=f^{-1}(\left\{a\right\})$
 (resp. $N_{t}=f^{-1}(\left\{a'\right\})$) and
 $N_{n}=f^{-1}(\left\{b\right\})$\,. This a consequence of implicit
 functions theorem which is the classical $\mathcal{C}^{1}$-version under
 Hypothesis~\ref{hyp:mainf} and
 still holds in a Lipschitz version
 under the more general Hypothesis~\ref{hyp:Lipbar} (see Subsection~\ref{sec:moregenLipgen}).
The density of 
$\mathcal C^{\infty}_{0}(f_{a}^{b}\cup f^{-1}(\left\{b\right\});\Lambda
T^{*}M)$ in $D(d_{f,f^{-1}([a,b]),h})$ provides the following
extension result:
\begin{equation}
  \label{eq:exten0}
  \forall \omega\in D(d_{f,f^{-1}([a,b]),h})\,,\quad \tilde{\omega}\in
  D(d_{f,f^{-1}([a',b]),h})\,,\ \ \text{where}\quad\tilde{\omega}\big|_{f_{a}^{b}}=\omega\ \ 
  \text{and}\ \  \tilde{\omega}\big|_{f_{a'}^{a}}\equiv 0\,.
\end{equation}
\item Hodge decomposition in Proposition~\ref{pr:domain}-3) says that
$$
\ker(\Delta_{f,\overline{\Omega},h})\simeq\ker(d_{f,\overline{\Omega},h})/\Ran(d_{f,\overline{\Omega},h})\simeq\ker(d_{0,\overline{\Omega},1})/\Ran(d_{0,\overline{\Omega},1})\,.
$$
From
the usual Hodge theory on the manifold with boundary 
$\overline{\Omega}$\,,
the dimension of $\ker(\Delta_{f,\overline{\Omega},h}^{(p)})$ is thus the relative
Betti number $\dim H^{p}(\overline{\Omega}\,,\, N_{t}
)$ and is independent of $h>0$\,.  In particular,  when
$\Omega=f_{a}^{b}$ and $a<b$ are not in $\left\{c_{1},\ldots,c_{N_{f}}\right\}$\,,  it is 
$$
\dim \ker(\Delta_{f,f^{-1}([a,b]),h})=
\dim H^{p}(f^{b},f^{a})=:\beta^{(p)}(f^{b},f^{a})\,.
$$
If moreover $[c,d]\subset [a,b]$ and  $([a,b]\setminus]c,d[)\cap\left\{c_{1},\ldots,c_{N_{f}}\right\}=\emptyset$\,,
 then for every $a'\in [a,c]$ and $b'\in[d,b]$\,, the dimensions
 $\dim H^{p}(f^{b},f^{a})$ and $\dim H^{p}(f^{b'},f^{a'})$ are equal
and then
\begin{equation}
\label{eq.const-kernel}
\dim \ker(\Delta_{f,f^{-1}([a,b]),h})=\dim \ker(\Delta_{f,f^{-1}([a',b']),h})
\,.
\end{equation}
\item When $s\geq 0$\,, the commutation of $d_{f,\overline{\Omega},h}$
  with $1_{[0,s]}(\Delta_{f,\overline{\Omega},h})$ ensures that
 the restricted differential  $\delta_{[0,s]}=1_{[0,s]}(\Delta_{f,\overline{\Omega},h})d_{f,\overline{\Omega},h}$
  defines a finite dimensional 
complex with Betti numbers $\dim H^{p}(\overline{\Omega}, N_{t})$:
\begin{equation}
  \label{eq.redcompl}
\xymatrix{
0\ar[r]&
F_{[0,s]}^{(0)}
\ldots
F_{[0,s]}^{(p-1)}
\ar[r]^-{\delta_{[0,s]}^{(p-1)}}\ar@<1ex>[l]&
F_{[0,s]}^{(p)}
\ar[r]^-{\delta_{[0,s]}^{(p)}}\ar@<1ex>[l]^-{\delta^{(p-1)*}_{[0,s]}}&
F_{[0,s]}^{(p+1)}
\ldots
F_{[0,s]}^{(d)}\ar[r]\ar@<1ex>[l]^-{\delta^{(p)*}_{[0,s]}}&0
\ar@<1ex>[l]
}
\end{equation}
where $F^{(p)}_{[0,s]}=\Ran
1_{[0,s]}(\Delta_{f,\overline{\Omega},h}^{(p)})$\,.  This will
be studied more carefully  when 
$\Omega=f_{a}^{b}$\,, with the notations $F_{[0,s],[a,b],h}$
and $\delta_{[0,s],[a,b],h}$ in order to handle various intervals
$[a,b]$\,. 
  \end{itemize}
\end{remark}

\subsection{Agmon's type estimates}
\label{sec:agmontype}

We review a series of 
exponential decay estimates which  are
 adapted
from \cite{DiSj}\cite{HeSj2}, and \cite{LNV} for Witten Laplacians with
boundary conditions. Those are standard when the function $f$ satisfy
Hypothesis~\ref{hyp:mainf} but only a part of them can be proved when
$f$ is a general Lipschitz function which satisfies
Hypothesis~\ref{hyp:Lipbar}. 

\subsubsection{Weighted integration by parts formulas}
\label{sec:consAgm}

We present here weighted integration by parts formulas with low
regularity assumptions.
These formulas will be used in the sequel, after optimizing the weights,
in order to prove different exponential decay estimates.
Under Hypothesis~\ref{hyp:mainf}, the regular case, this will lead to the usual
Agmon estimates presented in the next section.
A variation of these arguments will be developed in
Section~\ref{sec:moregenLip} 
under Hypothesis~\ref{hyp:realana} (subanalytic case) and will require the low regularity results listed below.

\begin{lem}
\label{le:Agmon}
Assume Hypothesis~\ref{hyp:domain} for $\overline{\Omega}=\Omega\sqcup
N_{t}\sqcup N_{n}$\,.
Let $f,\varphi\in W^{1,\infty}(M;\rz)$\,,
$\Delta_{f,\overline{\Omega},h}$ 
be the self-adjoint operator defined in
Proposition~\ref{pr:domain},
and $\sum_{j=1}^{J}\chi_{j}^{2}=1$ be a smooth partition of unity in $\overline\Omega$\,.
For any $\omega \in D(Q_{f,\overline{\Omega},h})=W_{\partial}(\Omega;\Lambda T^{*}M)$ (see \eqref{eq.W} and the lines below),
with the notation 
$$
\tilde{\omega}=e^{\frac{\varphi}{h}}\omega,
$$
the following identities hold true:
\begin{align}
\label{eq:reQ1}
&\Real Q_{f,\overline{\Omega},h}(\omega\,,\,
  e^{\frac{2\varphi}{h}}\omega)
=\|d_{f,\overline{\Omega},h}\tilde{\omega}\|^{2}+\|d_{f,\overline{\Omega},h}^{*}\tilde{\omega}\|^{2}-\langle\tilde{\omega}\,,\,|\nabla
\varphi|^{2}\tilde{\omega}\rangle\,,\\
\label{eq:reQ2}
\text{and}\quad&\Real Q_{f,\overline{\Omega},h}(\omega\,,\,
  e^{\frac{2\varphi}{h}}\omega)
=\sum_{j=1}^{J} \Real Q_{f,\overline{\Omega},h}(\chi_{j}\omega\,,\,
  e^{\frac{2\varphi}{h}}\chi_{j}\omega)-h^{2}\sum_{j=1}^{J}\||\nabla \chi_{j}|\tilde{\omega}\|^{2}\,.
\end{align}
Moreover, when in addition
$f\in\mathcal{C}^{2}(M)$\,, the identity \eqref{eq:reQ1}
writes also
\begin{align}
 \nonumber
  \Real Q_{f,\overline{\Omega},h}(\omega, e^{\frac{2\varphi}{h}}\omega)&=
h^{2}\|d\tilde{\omega}\|_{L^{2}}^{2}+h^{2}\|d^{*}\tilde{\omega}\|_{L^{2}}^{2}\\
\nonumber
&\hspace{1cm}+
\langle\tilde{\omega}\,,\,
 (|\nabla f|^{2}-|\nabla\varphi|^{2}+h\mathcal{L}_{\nabla f}+h
\mathcal{L}_{\nabla f}^{*})\tilde{\omega}\rangle\\
 \label{eq.ippphi}
&\hspace{1cm}+h\left(\int_{N_{n}}-\int_{N_{t}}\right)\langle\tilde \omega\,,\,\tilde\omega\rangle_{\Lambda
  T_{\sigma}^{*}\Omega}
\left(\frac{\partial f}{\partial n}\right)(\sigma)~d\sigma\,.
\end{align}
Lastly, when $f\in W^{1,\infty}(M;\rz)$ and $\varphi\in\mathcal{C}^{2}(M)$\,,
the above quantity can be written
  \begin{align}
\nonumber
\Real Q_{f,\overline{\Omega},h}(\omega\,,\,
  e^{\frac{2\varphi}{h}}\omega)
&=
Q_{f-\varphi,\overline{\Omega},h}(\tilde{\omega}\,,\,
  \tilde{\omega})
\\
\nonumber
&\hspace{1cm}
+ \langle \tilde{\omega}
 \,,\,  (2\nabla f.\nabla\varphi-2|\nabla
\varphi|^{2}+h\mathcal{L}_{\nabla
  \varphi}+h\mathcal{L}_{\nabla \varphi}^{*})\tilde{\omega}
  \rangle 
\\
\label{eq:reQ3}
&\hspace{1cm}
+h\left(\int_{N_{n}}-\int_{N_{t}}\right)\langle\tilde \omega\,,\,\tilde\omega\rangle_{\Lambda
  T_{\sigma}^{*}\Omega}
\left(\frac{\partial \varphi}{\partial n}\right)(\sigma)~d\sigma\,.
  \end{align}
\end{lem}
\begin{proof}
We recall that according to Remark~\ref{re.Lip}, $W_{\partial}(\Omega;\Lambda T^{*}M)$
is a Lipschitz-module.
For the first statement \eqref{eq:reQ1}, simply write
\begin{align*}
\Real
Q_{f,\overline{\Omega},h}(\omega\,,\, e^{\frac{2\varphi}{h}}\omega)&=
\Real
Q_{f,\overline{\Omega},h}(e^{-\frac{\varphi}{h}}\tilde{\omega}\,,\,
                                                                 e^{\frac{\varphi}{h}}\tilde{\omega})
\\
&=
\Real
\langle (d_{f,h}-d\varphi\wedge)\tilde{\omega}\,,\,
(d_{f,h}+d\varphi\wedge)\omega\rangle
\\
&\hspace{3cm}
+
\Real
\langle (d_{f,h}^{*}+\mathbf{i}_{\nabla \varphi})\tilde{\omega}\,,\, 
(d_{f,h}^{*}-\mathbf{i}_{\nabla \varphi})\tilde{\omega}\rangle
\\
&=
\|d_{f,h}\omega\|^{2}+\|d_{f,h}^{*}\omega\|^{2}
-\langle
  d\varphi\wedge\tilde{\omega}\,,\,d\varphi\wedge\tilde{\omega}\rangle
 -\langle
   \mathbf{i}_{\nabla \varphi}\tilde{\omega}\,,\,
   \mathbf{i}_{\nabla\varphi}\tilde{\omega}
\rangle
\\
&=
\|d_{f,h}\omega\|^{2}+\|d_{f,h}^{*}\omega\|^{2}
-\langle
 \tilde{\omega}\,,\,\underbrace{(\mathbf{i}_{\nabla\varphi}(d{\varphi}\wedge)+(d{\varphi}\wedge)\mathbf{i}_{\nabla
  \varphi})}_{=|\nabla \varphi|^{2}}\tilde{\omega}\rangle\,.
\end{align*}
For \eqref{eq:reQ2}, we start from \eqref{eq:reQ1} after noticing that
$\chi_{j}\tilde{\omega}\in W_{\partial}(\Omega;\Lambda T^{*}M)$ when $\omega\in
W_{\partial}(\Omega;\Lambda T^{*}M)$\,. We compute
\begin{align*}
\|d_{f,h}\chi_{j}\tilde{\omega}\|^{2}+\|d_{f,h}^{*}\chi_{j}\tilde{\omega}\|^{2}
&=
\|\chi_{j}d_{f,h}\tilde{\omega}\|^{2}+\|\chi_{j}d_{f,h}^{*}\tilde{\omega}\|^{2}\\
&\hspace{1cm}+
2\Real\langle \chi_{j}
  d_{f,h}\tilde{\omega}\,,\,(hd\chi_{j}\wedge)\tilde{\omega}\rangle
-2\Real\langle\chi_{j}d_{f,h}^{*}\tilde{\omega}\,,\,
  h\mathbf{i}_{\nabla \chi_{j}}\tilde{\omega}\rangle\\
&\hspace{1cm}
+h^{2}\left[\langle d\chi_{j}\wedge
\tilde{\omega}\,,\, d\chi_{j}\wedge \tilde{\omega} \rangle
+\langle \mathbf{i}_{\nabla \chi_{j}}\tilde{\omega}\,,\,
  \mathbf{i}_{\nabla \chi_{j}}\tilde{\omega}\rangle
\right]
\\
&=
\|\chi_{j}d_{f,h}\tilde{\omega}\|^{2}+\|\chi_{j}d_{f,h}^{*}\tilde{\omega}\|^{2}\\
&\hspace{1cm}+
\Real\langle 
  d_{f,h}\tilde{\omega}\,,\,(hd\chi_{j}^{2}\wedge)\tilde{\omega}\rangle
-\Real\langle d_{f,h}^{*}\tilde{\omega}\,,\,
  h\mathbf{i}_{\nabla \chi_{j}^{2}}\tilde{\omega}\rangle\\
&\hspace{1cm}
+h^{2}\langle
\tilde{\omega}\,,\, \underbrace{(\mathbf{i}_{\nabla \chi_{j}}(d\chi_{j}\wedge )+
(d\chi_{j}\wedge)\mathbf{i}_{\nabla \chi_{j}})}_{=|\nabla\chi_{j}|^{2}}\tilde{\omega} \rangle\,.
\end{align*}
Summing w.r.t $j\in \left\{1,\ldots,J\right\}$ leads to
$$
Q_{f,\overline{\Omega},h}(\omega\,,\,
e^{\frac{2\varphi}{h}}\omega)-\sum_{j=1}^{J}Q_{f,\overline{\Omega},h}(\chi_{j}\omega\,,\,e^{\frac{2\varphi}{h}}\chi_{j}\omega)=-h^{2}\sum_{j=1}^{J}\||\nabla
\chi_{j}|\tilde{\omega}\|^{2}\,.
$$
Let us now assume that $f\in \mathcal C^{2}(M)$\,.
According to \eqref{eq:reQ1}, 
the identity
$$\Real Q_{f,\overline{\Omega},h}(\omega\,,\,
  e^{\frac{2\varphi}{h}}\omega)
= Q_{f,\overline{\Omega},h}(\tilde\omega\,,\,
  \tilde\omega)-\langle\tilde{\omega}\,,\,|\nabla
\varphi|^{2}\tilde{\omega}\rangle
$$
holds true
and it  suffices to prove the formula
\eqref{eq.ippphi} when $\varphi=0$\,. To this end, 
one first writes
for $\omega\in  D\left ( Q_{f,\overline{\Omega},h}\right)$\,,
\begin{align*}
\left\| d_{f,h}\omega\right\|^{2}_{L^{2}}+
\left\|d_{f,h}^{*}\omega\right\|^{2}_{L^{2}}&=
h^{2}\left\| d\omega\right\|^{2}_{L^{2}}+
h^{2}\left\| d^{*}\omega\right\|^{2}_{L^{2}}
+\left\| d f\wedge \omega\right\|^{2}_{L^{2}}
\\
&\quad
+ \left\| \mathbf i_{\nabla f} \omega\right\|^{2}_{L^{2}} + \  h\big(\langle d f\wedge \omega,
d\omega \rangle_{L^{2}} +  \langle d \omega,
d f\wedge \omega \rangle \\
&\quad 
+\langle d^{*} \omega,
\mathbf i_{\nabla f} \omega \rangle_{L^{2}} +  \langle \mathbf i_{\nabla f} \omega,
d^{*}\omega \rangle \big)\\
&=
h^{2}\left\| d\omega\right\|^{2}_{L^{2}}+
h^{2}\left\| d^{*}\omega\right\|^{2}_{L^{2}}+
\left\|\left|\nabla f\right|\omega\right\|^{2}_{L^{2}}\\
 &\quad+h\langle  \omega, (\mathcal L_{\nabla f}+  \mathcal L^{*}_{\nabla f})\omega \rangle_{L^{2}}
 +   h(\langle d f\wedge \omega,
d\omega \rangle_{L^{2}} \\
&\quad -
\langle d^{*}(d f\wedge \omega),
\omega \rangle_{L^{2}}
-\langle d \mathbf i_{\nabla f} \omega , \omega \rangle_{L^{2}}
 +   \langle \mathbf i_{\nabla f} \omega,
d^{*}\omega \rangle_{L^{2}} )\,,
\end{align*}
where the last equality holds thanks to the relations
$(df \wedge)^*=\mathbf i_{\nabla f}$\,,
\begin{align*}
&\mathcal L_{\nabla f}=d\circ \mathbf i_{\nabla f}+\mathbf i_{\nabla f}\circ d
\quad\text{and}\quad
 \mathcal L^{*}_{\nabla f}=(df\wedge)\circ d^{*}+ d^{*}\circ (df\wedge)\,.
\end{align*}
The relation \eqref{eq.ippphi} follows using in addition
the generalized Green formula \eqref{eq.gen-Green} which gives here, since $\omega\in  D\left ( Q_{f,\overline{\Omega},h}\right)$
and hence admits a total trace on $N$\,,
and  $df\wedge \omega$\,, $\mathbf i_{\nabla f} \omega\in \{v\in L^{2}, dv \in L^{2}, d^{*}v\in L^{2}\}$:
\begin{align*}
\label{eq.trace-norm}
\langle d f\wedge \omega,
d\omega \rangle_{L^{2}} -
\langle d^{*}(d f\wedge \omega),
\omega \rangle_{L^{2}} &= \int_{N_n}\langle {n}^{\flat}\wedge \omega ,  {n}^{\flat}\wedge \mathbf i_{n} (df\wedge \omega)\rangle_{ T^*_{\sigma} \Omega}d\sigma\\
&=\int_{N_n}\langle  \omega ,\mathbf i_{n}  ({n}^{\flat}\wedge \mathbf i_{n} (df\wedge \omega))\rangle_{ T^{*}_{\sigma}\Omega}d\sigma\\
&=\int_{N_n}\langle  \omega , \mathbf i_{n} (df\wedge \omega)\rangle_{ T^{*}_{\sigma}\Omega}d\sigma\\
&=\int_{N_n}(\pa_{n}f\, \langle  \omega , \omega\rangle_{ T^{*}_{\sigma}\Omega}- \langle  \omega , df\wedge \underbrace{\mathbf i_{n}\omega}_{=0}\rangle _{ T^{*}_{\sigma}\Omega}   )d\sigma\\
&=\int_{N_n}\pa_{n}f\, \langle  \omega , \omega\rangle_{ T^{*}_{\sigma}\Omega}d\sigma
\end{align*}
as well as
\begin{align*}
   \langle \mathbf i_{\nabla f} \omega,
d^{*}\omega \rangle_{L^{2}}
-\langle d \mathbf i_{\nabla f} \omega , \omega \rangle_{L^{2}}
& = 
-\int_{N_{t}}\langle {n}^{\flat}\wedge \mathbf i_{\nabla f} \omega ,  {n}^{\flat}\wedge \mathbf i_{n} \omega \rangle_{ T^*_{\sigma}\Omega}
\\
 &=-\int_{N_t}\pa_{n}f\, \langle  \omega , \omega\rangle_{ T^{*}_{\sigma}\Omega}d\sigma
d\sigma\,.
\end{align*}
Lastly, let us prove the relation  \eqref{eq:reQ3}.
By direct expansion with  $f$
and $\varphi$ Lipschitz continuous
and 
$$
d_{f-\varphi,h}=d_{f,h} -(d\varphi\wedge) = hd + (df\wedge) -
(d\varphi\wedge)
\quad\text{and}
\quad 
d^{*}_{f-\varphi,h}=d^{*}_{f,h} -{\bf i}_{\nabla \varphi}= hd^{*} + {\bf i}_{\nabla f} -
{\bf i}_{\nabla \varphi}\,,
$$
we obtain 
\begin{align*}
Q_{f-\varphi,\overline{\Omega},h}(\tilde \omega,\tilde \omega)&=
Q_{f,\overline{\Omega},h}(\tilde \omega,\tilde \omega)\\
&\quad- 2\Real\big(\langle d f\wedge \tilde \omega,
d \varphi\wedge \tilde \omega \rangle +  \langle \mathbf i_{\nabla f} \tilde\omega,
\mathbf i_{\nabla \varphi} \tilde\omega \rangle  \big)\\
&\quad
- 2h\Real\big(\langle d \tilde\omega,
d \varphi\wedge \tilde\omega \rangle
+ \langle d^{*} \tilde\omega,
\mathbf i_{\nabla \varphi}\tilde\omega \rangle\big)\\
&\quad +\left\| d \varphi\wedge \tilde\omega\right\|^{2}+\left\| \mathbf i_{\nabla \varphi} \tilde\omega\right\|^{2}\,.
\end{align*}
By adding this relation for the pairs $(f,\varphi)$ and $(0,-\varphi)$\,, we obtain
\begin{align*}
Q_{f-\varphi,\overline{\Omega},h}(\tilde \omega,\tilde \omega)+Q_{\varphi,\overline{\Omega},h}(\tilde \omega,\tilde \omega)&=
Q_{f,\overline{\Omega},h}(\tilde \omega,\tilde \omega)+Q_{0,\overline{\Omega},h}(\tilde \omega,\tilde \omega)\\
&\quad- 2\underbrace{\Real\big(\langle d f\wedge \tilde \omega,
d \varphi\wedge \tilde \omega \rangle +  \langle \mathbf i_{\nabla f} \tilde\omega,
\mathbf i_{\nabla \varphi} \tilde\omega \rangle  \big)}_{=\langle \tilde \omega,(\nabla f\cdot \nabla \varphi) \tilde\omega\rangle}\\
&\quad
+0\\
&\quad +2\left\| |\nabla\varphi| \tilde\omega\right\|^{2}\,.
\end{align*}
Finally, using the relation \eqref{eq:reQ1} gives
\begin{align*}
\Real Q_{f,\overline{\Omega},h}( \omega,e^{\frac{2\varphi}{h}}\omega)&=
Q_{f,\overline{\Omega},h}(\tilde \omega,\tilde \omega) - \left\| |\nabla\varphi| \tilde\omega\right\|^{2} \\
&=  Q_{f-\varphi,\overline{\Omega},h}(\tilde \omega,\tilde \omega)
+2 \langle \tilde \omega,(\nabla f\cdot \nabla \varphi -|\nabla \varphi|^{2}) \tilde\omega\rangle\\
&\quad + Q_{\varphi,\overline{\Omega},h}(\tilde \omega,\tilde \omega)-
Q_{0,\overline{\Omega},h}(\tilde \omega,\tilde \omega) - \left\| |\nabla\varphi| \tilde\omega\right\|^{2}.
\end{align*}
When in addition  $\varphi \in \mathcal{C}^{2}(M)$\,, 
using \eqref{eq:reQ1} and \eqref{eq.ippphi} with $f=\varphi$ leads to the relation
\eqref{eq:reQ3}.
\end{proof}

\begin{remark}
Alternatively, one could first prove the relation \eqref{eq:reQ3}
for $f,\varphi\in \mathcal{C}^{2}(M)$\,, and then approximate
a general $f\in W^{1,\infty}(M)$ by a sequence in $\mathcal C^{2}(M)$
as in Remark~\ref{re.Lip}.
\end{remark}

\subsubsection{Exponential decay estimates}
\label{sec:expdecayest}

Under Hypothesis~\ref{hyp:mainf}, 
these estimates rely on the integration by parts
formula \eqref{eq.ippphi}
of Lemma~\ref{le:Agmon}. They will be replaced 
by a new hypothesis for more general Lipschitz function $f$\,, which
will be ultimately  verified when $f$ is Lipschitz subanalytic  in
Subsection~\ref{sec:moregenLip}.

\begin{definition}
Assume Hypothesis~\ref{hyp:mainf} for $f$ and remember
$$
M_{reg}=\left\{x\in (M\setminus\mathrm{suppsing}\; f)\,, \nabla
  f(x)\neq 0\right\}\subset 
M\setminus f^{-1}(\left\{c_{1},\ldots, c_{N_{f}}\right\})\,.
$$
  The Agmon distance $d_{Ag}$ on $M$ associated with $f\in
  \mathcal{C}^{\infty}(M)$ is the geodesic pseudodistance associated
  with the degenerate metric $1_{M_{reg}}|\nabla f|^{2}g$\,, namely
$$
d_{Ag}(x,y)=\inf_{
\scriptsize
  \begin{array}[c]{l}
 \gamma\in {\cal C}^{1}([0,1];M)\,,\\  \gamma(0)=x\,,\gamma(1)=y
  \end{array}
}\int_{0}^{1}1_{M_{reg}}(\gamma(t))|\nabla f(\gamma(t))|
|\gamma'(t)|~dt\,.
$$
\end{definition}
\noindent Because $f\in W^{1,\infty}(M)\cap {\cal C}^{\infty}(M_{reg})$\,, we know
$d_{Ag}(x,y)\leq \|\nabla f\|_{L^{\infty}}d_{g}(x,y)$ where $d_{g}$ is
the geodesic distance and $d_{Ag}$ is a Lipschitz function of
$(x,y)\in M\times M$\,.
Moreover when $x,y$ belong to the same connected component of
$M\setminus f^{-1}(\left\{c_{1},\ldots,c_{N_{f}}\right\})$
 any  ${\cal C}^{1}$ curve $\gamma$ staying in this connected
 component satisfies
$$
\int_{0}^{1}|\nabla f(\gamma(t))||\gamma'|~dt\geq
|\int_{0}^{1}\nabla f(\gamma(t)).\gamma'(t)~dt|=|f(y)-f(x)|\,.
$$
For a general $\gamma\in {\cal C}^{1}([0,1];M)$ such that
$\gamma(0)=x$ and $\gamma(1)=y$\,,
$\left\{f(\gamma(t))\,, t\in [0,1]\right\}$ is a compact
 interval. Therefore,
bounding from below 
the integral $\int_{0}^{1}\ldots dt$ 
by a sum of integrals on intervals $]t_{k},t'_{k}[$\,,
 where  $f(\gamma(t))\not\in\left\{c_{1},\ldots,c_{N_{f}},\max(f\circ
   \gamma), \min (f\circ \gamma)\right\}
$\,, leads to
$$ 
\int_{0}^{1}1_{M_{reg}}(\gamma(t))\left|\nabla
  f(\gamma(t))\right||\gamma'(t)|~dt\geq \max_{t\in
  [0,1]}f(\gamma(t))-\min_{t\in[0,1]}f(\gamma(t))\geq |f(y)-f(x)|\,.
$$
We obtain
\begin{equation}
\label{eq:dAgmmin}
\forall x,y\in M\,, \quad \|\nabla f\|_{L^{\infty}}d_{g}(x,y)\geq d_{Ag}(x,y)\geq |f(y)-f(x)|\,.
\end{equation}
When $f$ is a ${\cal C}^{\infty}$ Morse function,
more details about the more general broken
 geodesic curves, which do not hold
anymore with our general assumption and which we do not
need, are given in \cite{HeSj4}.
\begin{prop}
\label{pr:Agmon}
Assume Hypothesis~\ref{hyp:mainf} for $f$ and Hypothesis~\ref{hyp:domain} for
$\overline{\Omega}=\Omega\sqcup N_{t}\sqcup N_{n}$ with
\begin{equation}
\label{eq:signNtn}
\partial \Omega=N_{t}\sqcup N_{n}\subset M_{reg}\quad,\quad
\frac{\partial f}{\partial n}\big|_{N_{t}}<0\quad,\quad
\frac{\partial f}{\partial n}\big|_{N_{n}}>0\,.
\end{equation}
Let $\Delta_{f,\overline{\Omega},h}$ be the self-adjoint operator
defined in Proposition~\ref{pr:domain} and let $U$ denote the compact
subset of $\Omega$\,, $U=(M\setminus M_{reg})\cap \Omega$\,.
All families $(\lambda_{h})_{h>0}$ in $\cz$\,,
$(r_{h})_{h>0}$ in $L^{2}(\Omega)$ and $(\omega_{h})_{h>0}$ in
$D(\Delta_{f,\overline{\Omega},h})\subset W_{\partial}(\Omega;\Lambda T^{*}M)$
such that
$$
(\Delta_{f,\overline\Omega,h}-\lambda_{h})\omega_{h}=r_{h}\quad,\quad
\supp r_{h}\subset K\quad,\quad \lim_{h\to 0}\lambda_{h}=0\,,
$$
where $K$ is a fixed compact subset of $\overline{\Omega}$\,,
satisfy  the estimate (see \eqref{eq.W} and the lines below)
$$
\|e^{\frac{d_{Ag}(\cdot,U\cup
  K)}{h}}\omega_{h}\|_{W_{\partial}(\Omega)}=\tilde{O}(1)
\times(\|r_{h}\|_{L^{2}(\Omega)}+t_{U}\|\omega_{h}\|_{L^{2}(\Omega)})
\,,
$$
where $t_{U}=1$ if $U\neq \emptyset$ and $t_{U}=0$ if $U=\emptyset$\,.
\end{prop}

\begin{proof}
For $\varepsilon\in]0,1[$\,, one  introduces $K_{\varepsilon}=\left\{y\in \overline{\Omega}\,, d_{Ag}(y,U\cup
  K)\leq \varepsilon\right\}$
  and
  $\chi_{1}=\chi_{1,\varepsilon},\chi_{2}=\chi_{2,\varepsilon}\in C^{\infty}(\overline{\Omega},[0,1])$
  such that $\chi_{1}\equiv 0$ when $U=\emptyset$ and $\chi_{1}=1$ near $U$ else, $\supp\chi_{1}\subset K_{\varepsilon}\cap\Omega  $\,, and
  $\chi^{2}_{1}+\chi_2^{2}\equiv 1$\,.
Let us also introduce
  $\varphi_{\varepsilon}:
x\mapsto (1-\varepsilon)d_{Ag}(x,K_{\varepsilon}) \in W^{1,\infty}(\Omega)$\,, so that
$\varphi_{\varepsilon}$ satisfies  $|\nabla \varphi_{\varepsilon}|\leq
(1-\varepsilon)|\nabla f| $ almost everywhere
in $\Omega$\,. Setting
$\tilde{\omega}_{h}:=e^{\frac{\varphi_{\varepsilon}}{h}}\omega_{h}$
and applying  \eqref{eq:reQ2}
with  $\varphi_{\varepsilon}= 0$ on $K_{\varepsilon} $\,, $\supp
\chi_{1},\supp r_{h}\subset K_{\varepsilon}$\,, we obtain
\begin{align*}
\langle  r_{h}, \omega_{h}\rangle_{L^{2}}+\lambda_{h}\| \tilde{\omega}_{h}\|^{2}_{L^{2}}
&=\Real Q_{f,\overline{\Omega},h}(\omega_{h},e^{2\frac{\varphi_{\varepsilon}}{h}}\omega_{h})\\
&\geq \Real Q_{f,\overline{\Omega},h}(\chi_{2}\omega_{h},\chi_{2}e^{2\frac{\varphi_{\varepsilon}}{h}}\omega_{h})
+Q_{f,\overline{\Omega},h}(\chi_{1}\omega_{h},\chi_{1}\omega_{h}) - c_{\varepsilon}h^{2}\| \tilde{\omega}_{h}\|^{2}_{L^{2}}\,.
\end{align*}
Then, applying \eqref{eq.ippphi} of Lemma~\ref{le:Agmon} with a
$\mathcal{C}^{2}$-extension to $M$ of $f\big|_{\supp \chi_{2}}$\,, 
with $|\nabla f|^{2}\geq C_{\varepsilon}$ on $\supp \chi_{2}$ and  the sign condition~\eqref{eq:signNtn} leads to 
\begin{align}
\nonumber
\|\omega_{h}\|_{L^{2}}\|r_{h}\|_{L^{2}}
&\geq Q_{f,\overline{\Omega},h}(\chi_{1}\omega_{h},\chi_{1}\omega_{h})+
h^{2}\big(\|d\chi_{2}\tilde{\omega}_{h}\|^{2}_{L^{2}}
+\|d^{*}\chi_{2}\tilde{\omega}_{h}\|^{2}_{L^{2}}\big)\\
\nonumber
&\qquad\qquad\qquad\qquad\quad  \ +(C_{\varepsilon} - \lambda_{h}-c_{\varepsilon}h^{2})\|\chi_{2}\tilde{\omega}_{h}\|^{2}_{L^{2}}
- t_{U}(\lambda_{h}+c_{\varepsilon}h^{2})\|\chi_{1}\omega_{h}\|^{2}_{L^{2}}
\\
\label{eq.comp-1}
&\geq Q_{f,\overline{\Omega},h}(\chi_{1}\omega_{h},\chi_{1}\omega_{h})+
C'_{\varepsilon}\,h^{2}\|\chi_{2}\tilde{\omega}_{h}\|^{2}_{W}-t_{U}\|\omega_{h}\|^{2}_{L^{2}}\,,
\end{align}
where we recall from Definition~\ref{de:defW} that $\|\omega\|_{W}=\|\omega\|_{L^{2}}+\|d\omega\|_{L^{2}}+\|d^{*}\omega\|_{L^{2}}$\,.\\
Since $Q_{f,\overline{\Omega},h}(\chi_{1}\omega_{h},\chi_{1}\omega_{h})\geq 0$ and $\|\omega_{h}\|_{L^{2}}\leq
C\|r_{h}\|_{L^{2}}+t_{U}\|\omega_{h}\|_{L^{2}}$
(this is obvious when $U\neq\emptyset$ and apply
\eqref{eq.ippphi} of Lemma~\ref{le:Agmon} with $\varphi=0$ else), we obtain
the estimate
\begin{equation}
\label{eq:W12tilde}
\|\chi_{2}\tilde{\omega}_{h}\|_{W_{\partial}(\Omega)}\leq
\frac{C''_{\varepsilon}}{h}
\big(\|r_{h}\|_{L^{2}}+ t_{U}\|\omega_{h}\|_{L^{2}}\big)\,.
\end{equation}
This ends the proof when $U=\emptyset$\,.\\
When  $U\neq\emptyset$\,, 
the relations \eqref{eq.comp-1} 
and
\begin{align*}
Q_{f,\overline{\Omega},h}(\chi_{1}\omega_{h},\chi_{1}\omega_{h})
&=\|(hd +df\wedge)\chi_{1}\omega_{h} \|_{L^2}^{2}+
\|(hd^{*} +\mathbf i_{\nabla f})\chi_{1}\omega_{h} \|_{L^2}^{2}\\
&\geq \frac{h^{2}}2(\|d\chi_{1}\omega_{h} \|_{L^2}^{2}+
\|d^{*}\chi_{1}\omega_{h} \|_{L^2}^{2})
-C\|\chi_{1}\omega_{h} \|_{L^2}^{2}
\end{align*} 
lead, since $\varphi_{\varepsilon}= 0$ on $\supp\chi_{1}$\,, to
\begin{equation}
\label{eq:W12tilde'}
\|\chi_{1}\tilde{\omega}_{h}\|_{W_{\partial}(\Omega)}\leq
\frac{C'}{h}
\big(\|r_{h}\|_{L^{2}}+ \|\omega_{h}\|_{L^{2}}\big)\,.
\end{equation}
The statement of  Proposition~\ref{pr:Agmon}  then follows from
\eqref{eq:W12tilde} and \eqref{eq:W12tilde'}, by using again the IMS
localization formula \eqref{eq:reQ2} with now $\varphi=f=0$ but $\omega$
replaced by $\tilde{\omega}$\,.
\end{proof}
  Following \cite{HeSj2}\cite{DiSj} we extend the
definition of $\tilde{O}$ to the kernels of 
bounded operators from $L^{2}$ to 
$W$\,, which appears to be more natural than
$W^{1,2}$ in our setting (see indeed Definition~\ref{de:defW} and
Proposition~\ref{pr:exttn0}). For more flexibility, 
boundary conditions do not appear 
in the following definition and the full
space $W(\Omega;\Lambda T^{*}M)$ of Definition~\ref{de:defW} is used.
\begin{definition}
\label{de:Otkernel}
Let  the domain $\overline{\Omega}$ satisfy
Hypothesis~\ref{hyp:domain}. Let the operator
$A_{h}$ act continuously from $L^{2}(\Omega;\Lambda T^{*}M)$ to 
$W(\Omega;\Lambda T^{*}M)$ 
and let $\Phi\in
\mathcal{C}^{0}(\overline{\Omega}\times\overline{\Omega};\rz)$\,.
We say that the kernel $A_{h}(x,y)$ of $A_{h}$ is
$\tilde{O}(e^{-\frac{\Phi(x,y)}{h}})$
if, for all $x_{0},y_{0}\in \Omega$ and $\varepsilon>0$\,,
there exist neighborhoods $U_{\varepsilon},V_{\varepsilon}$ in $M$ of $y_{0}$ and
$x_{0}$ and  constants $h_{\varepsilon}$ such that
\begin{multline*}
\forall h\in ]0,h_{\varepsilon}[\,,\,
\forall \chi\in
\mathcal{C}^{\infty}_{0}(V_{\varepsilon})\,,\,\exists C_{\chi,\varepsilon}>0\,,\,\forall u\in L^{2}(\Omega)\ \text{s.t.}\  \supp u \subset
U_{\varepsilon}\,,\ \\
\|\chi A_{h}u\|_{W(\Omega)}\leq C_{\chi,\varepsilon}e^{-\frac{\Phi(x_{0},y_{0})-\varepsilon}{h}}\|u\|_{L^{2}}\,.
\end{multline*}
For a finite family $\big(\Phi_{k}\big)_{k\in\{1,\dots,K\}}$ in 
$\mathcal{C}^{0}(\overline{\Omega}\times\overline{\Omega};\rz)$\,, 
 the kernel $A_{h}(x,y)$ of $A_{h}$ is
said to be $\tilde{O}(\sum_{k=1}^{K}e^{-\frac{\Phi_{k}(x,y)}{h}})$
when it is $\tilde{O}(e^{-\frac{\min_{1\leq
      k\leq  K}\Phi_{k}(x,y)}{h}})$\,.
\end{definition}
When $A_{h}(x,y)=\tilde{O}(e^{-\frac{\Phi(x,y)}{h}})$ and
$B_{h}(x,y)=\tilde{O}(e^{-\frac{\Psi(x,y)}{h}})$ and $D_{h}$ is a
differential operator of order $\leq 1$ which vanishes in a fixed
(independent of $h$)
neighborhood of $\partial \Omega$ 
(remember $W(\Omega;\Lambda T^{*}M)\subset
W^{1,2}_{loc}(\Omega;\Lambda T^{*}M)$), with 
$\|D_{h}\|_{\mathcal{L}(W^{1,2};L^{2})}=\tilde{O}(1)$\,, then
$(A_{h}D_{h}B_{h})(x,y)=\tilde{O}(e^{-\frac{\Theta(x,y)}{h}})$ with
$\Theta(x,y)=\min_{z\in \overline{\Omega}}\Phi(x,z)+\Psi(z,y)$\,.\\
If $A_{h}(x,y)=\tilde{O}(e^{-\frac{\Phi(x,y)}{h}})$ and
$\psi\in \mathcal{C}^{0}(\overline{\Omega})$\,, $\varphi\in
W^{1,\infty}(\overline{\Omega})$ satisfy
$\varphi(x)\leq \Phi(x,y)-\psi(y)$ for all $y\in \overline{\Omega}$\,,
then 
$\sup_{u\in
  L^{2}(\Omega)}\frac{\|e^{\frac{\varphi}{h}}A_{h}u\|_{W^{1,2}}}{\|e^{\frac{\psi}{h}}u\|_{L^{2}}}=\tilde{O}(1)$\,.\\
An easy application concerns the
 case when the gradient of $f$ does
not vanish in $\overline{\Omega}\subset M_{reg}$\,, under
Hypothesis~\ref{hyp:mainf}. 
\begin{prop}
\label{pr:Agmon1}
Assume Hypothesis~\ref{hyp:mainf} for $f$\,, and
Hypothesis~\ref{hyp:domain} for $\overline{\Omega}=\Omega\sqcup
N_{t}\sqcup N_{n}$ with now
\begin{equation}
\label{eq:signNtnfull}
\overline{\Omega}\subset M_{reg}\quad,\quad
\frac{\partial f}{\partial n}\big|_{N_{t}}<0\quad,\quad
\frac{\partial f}{\partial n}\big|_{N_{n}}>0\,,
\end{equation}
where we recall that $f\in \mathcal{C}^{\infty}(M_{reg})$ has a non
vanishing gradient.
The self-adjoint operator  $\Delta_{f,\overline{\Omega},h}$  defined
in Proposition~\ref{pr:domain} 
is bounded from below
by $c_{\Omega,f,h_{1}}>0$\,.
when $h\in ]0,h_{1}[$ with $h_{1}>0$  small
enough. 
If $\lim_{h\to 0}\rho(h)=0^{+}$\,, then the resolvent
$(\Delta_{f,\overline{\Omega},h}-z)^{-1}$\,, $|z|\leq \rho(h)$\,, 
well defined for $h\in ]0,h_{0}[$\,, $h_{0}>0$ small enough, satisfies
$$(\Delta_{f,\overline{\Omega},h}-z)^{-1}(x,y)=\tilde{O}(e^{-\frac{d_{Ag}(x,y)}{h}})\leq\tilde{O}(e^{-\frac{|f(x)-f(y)|}{h}})\,,$$
according to Definition~\ref{de:Otkernel} and uniformly with respect
to $z$\,, $|z|\leq \rho(h)$\,.
\end{prop}
\begin{proof}
The lower bound and the definition of the
 resolvent is deduced from 
\eqref{eq.ippphi} in  
Lemma~\ref{le:Agmon} applied with $\varphi\equiv 0$\,, $|\nabla
  f(x)|\geq c>0$ for all $x\in \overline{\Omega}$ and  where the
  condition~\eqref{eq:signNtnfull} ensures the positivity 
of the boundary
  terms. The estimate of the kernel is then a straightforward
  consequence of Proposition~\ref{pr:Agmon} 
with here $U=\emptyset$\,.
\end{proof}
We cannot prove Proposition~\ref{pr:Agmon} and
Proposition~\ref{pr:Agmon1} for a general Lipschitz function even under
Hypothesis~\ref{hyp:Lipbar}. We replace it by an assumption which is
proved to be fulfilled by subanalytic Lipschitz functions 
in Subsection~\ref{sec:moregenLip}
\begin{hyp}
\label{hyp:AgmonLip}
For a Lipschitz function which satisfy Hypothesis~\ref{hyp:Lipbar}
with the ``critical values'' $c_{1}<\ldots<c_{N_{f}}$\,, we assume that
Proposition~\ref{pr:Agmon} and Proposition~\ref{pr:Agmon1} 
hold true after replacing
$M_{reg}$ by $M\setminus f^{-1}\left(\left\{c_{1},\ldots,c_{N_{f}}\right\}\right)$\,,  $d_{Ag}(x,y)$
 by the pseudodistance
$|f(x)-f(y)|$\,,
and by restricting to the case
$\overline{\Omega}=f^{-1}([a,b])$\,, $a<b$\,, $a,b\not\in
\left\{c_{1},\ldots,c_{N_{f}}\right\}$\,. 
\end{hyp}
\subsubsection{Adjusting boundary conditions}
\label{sec:extlemma}

Another consequence of Agmon estimates 
is the following lemma
which will be used to correct boundary conditions
 and to extend solutions to $d_{f,h}\omega=0$ to a wider domain with suitably small errors.
Under Hypothesis~\ref{hyp:mainf}, it is stated in the more general framework of
Proposition~\ref{pr:Agmon1} with $\overline{\Omega}\subset M_{reg}$\,,
although it will be applied essentially when
$\overline{\Omega}=f^{-1}([a,b])$ with $[a,b]\cap
\left\{c_{1},\ldots,c_{N_{f}}\right\}=\emptyset$\,. 
For a more general Lipschitz function we work directly in the
framework of Hypothesis~\ref{hyp:AgmonLip}.
\begin{lem}
\label{le:exten1}
Assume Hypothesis~\ref{hyp:mainf} for $f$ and
Hypothesis~\ref{hyp:domain} for $\overline{\Omega}=\Omega\sqcup
N_{t}\sqcup N_{n}$ with $\overline{\Omega}\subset M_{reg}$ and the
sign conditions $\frac{\partial f}{\partial n}\big|_{N_{t}}<0$\,,
$\frac{\partial f}{\partial n}\big|_{N_{n}}>0$\,.
Consider the operator
$\Delta_{f,\overline{\Omega},h}$ of Proposition~\ref{pr:domain}.
There exists $c>0$ and $h_{0}>0$ determined by $f$ and
$\overline{\Omega}$ and for any pair of cut-off functions
 $\chi,\tilde{\chi}\in \mathcal{C}^{\infty}(\overline{\Omega};[0,1])$
which satisfies $d\chi,d\tilde{\chi}\in \mathcal{C}^{\infty}_{0}(\Omega)$\,,
with $\tilde{\chi}\equiv 1$ in a neighborhood of $\supp d\chi$\,,
a constant $C_{\chi,\chi'}>0$ such that the following holds.\\
When $\omega\in
W(\Omega; \Lambda T^{*}M)$\,,
the
forms
\begin{eqnarray*}
 \eta_{1}=d_{f,\overline{\Omega},h}^{*}(\Delta_{f,\overline{\Omega},h})^{-1}((hd\chi)\wedge
     \omega)\quad\text{and}
     \quad
   \eta_{2}=d_{f,\overline{\Omega},h}(\Delta_{f,\overline{\Omega},h})^{-1}(h\mathbf{i}_{\nabla \chi}\omega)
\end{eqnarray*}
both belong to 
$$
D(\Delta_{f,\overline{\Omega},h})\subset
W_{\partial}(\Omega;\Lambda T^{*}M)\subset W(\Omega;\Lambda T^{*}M)
$$
 and satisfy the following inequality with  convention $d_{Ag}(\supp d\tilde{\chi}, \supp
d\chi)=+\infty$ when $\tilde{\chi}$ is the constant function $1$:
\begin{eqnarray*}
&&
\|\eta_{1}\|_{L^{2}}\leq \frac1{\sqrt{c}} \|(hd\chi)\wedge
     \omega\|_{L^{2}}
\quad\text{and}\quad
\|\eta_{2}\|_{L^{2}}\leq \frac1{\sqrt{c}} \|(h\mathbf{i}_{\nabla\chi})
     \omega\|_{L^{2}}\,,
\end{eqnarray*}
\begin{eqnarray*}
&&
     \|d_{f,h}(\chi\omega-\tilde{\chi}\eta_{1})\|_{L^{2}}\leq 
\frac{1}{\sqrt{c}}\|(hd\chi)\wedge
     d_{f,h}\omega\|_{L^{2}}+\|\chi d_{f,h}\omega\|_{L^{2}}
\\
&&\hspace{7cm}+\tilde{O}(e^{-\frac{d_{Ag}(\supp
     d\tilde{\chi},\supp d\chi)}{h}})\|(hd\chi)\wedge
     \omega\|_{L^{2}}\,,
\end{eqnarray*}
\begin{eqnarray*}
&&
\|d_{f,h}^{*}(\tilde{\chi}\eta_{1})\|_{L^{2}}
\leq
\tilde{O}(e^{-\frac{d_{Ag}(\supp
     d\tilde{\chi},\supp d\chi)}{h}})\|(hd\chi)\wedge
     \omega\|_{L^{2}}\,,
\end{eqnarray*}
\begin{eqnarray*}
&&\|d_{f,h}^{*}(\chi\omega -\tilde{\chi}\eta_{2})\|_{L^{2}}
\leq \frac{1}{\sqrt{c}}\|h\mathbf{i}_{\nabla
   \chi}d_{f,h}^{*}\omega\|_{L^{2}}
+\|\chi
   d_{f,h}^{*}\omega\|_{L^{2}}\\
&&\hspace{7cm}
+\tilde{O}(e^{-\frac{d_{Ag}(\supp
     d\tilde{\chi},\supp d\chi)}{h}})\|h\mathbf{i}_{\nabla
   \chi}\omega\|_{L^{2}}\,,
\end{eqnarray*}
\begin{eqnarray*}
\|d_{f,h}(\tilde{\chi}\eta_{2})\|_{L^{2}}\leq 
\tilde{O}(e^{-\frac{d_{Ag}(\supp
     d\tilde{\chi},\supp d\chi)}{h}})\|h\mathbf{i}_{\nabla \chi}
     \omega\|_{L^{2}}\,,
\end{eqnarray*}
\begin{eqnarray*}
&&
 \left(\begin{array}[c]{c}
     \|d_{f,h}(\chi \omega
     -\tilde{\chi}(\eta_{1}+\eta_{2}))\|_{L^{2}}\\
+\\
\|d_{f,h}^{*}(\chi \omega-\tilde{\chi}(\eta_{1}+\eta_{2}))\|_{L^{2}}
   \end{array}
\right)
\leq
C_{\chi,\tilde{\chi}}\left[\|d_{f,h}\omega\|_{L^2}+\|d_{f,h}^{*}\omega\|_{L^{2}}\right]
\\
&&\hspace{7cm}+ \tilde{O}(e^{-\frac{d_{Ag}(\supp
     d\tilde{\chi},\supp d\chi)}{h}})\|\omega\|_{L^{2}(\supp d\chi)}
\,.
\end{eqnarray*}
When $f$ is a Lipschitz function which satisfies
Hypothesis~\ref{hyp:Lipbar} and Hypothesis~\ref{hyp:AgmonLip}
the results are the same when $\overline{\Omega}=f^{-1}([a,b])$\,,
$c_{n}<a<b<c_{n+1}$\,,  and $d_{Ag}(K,K')$
is replaced by $\inf_{x\in K\,, y\in K'}|f(x)-f(y)|$\,.
\end{lem}
\begin{remark}
  Note that $\omega$ is not assumed to belong to the domain of
 $d_{f,\overline{\Omega},h}$\,, $d_{f,\overline{\Omega},h}^{*}$ or $\Delta_{f,\overline{\Omega},h}$ (no boundary conditions) and the same holds in general for
  $\chi\omega$\,. Accordingly, we used the notations $d_{f,h}$ and
  $d_{f,h}^{*}$ for the differential operators.
In some applications $\chi$ will be chosen such that
  $\chi \omega$ and therefore  $\chi
  \omega-\tilde{\chi}(\eta_{1}+\eta_{2})$
 belong to  one of these domains. Example given, if $\chi \omega\in D(\Delta_{f,\overline{\Omega},h})$\,, 
the last inequality then
  provides a good estimate of $Q_{f,\overline{\Omega},h}(\chi
  \omega-\tilde{\chi}(\eta_{1}+\eta_{2}))$ when $\supp \chi$ and
  $\supp \tilde{\chi}$ are well chosen.
\end{remark}
\begin{proof} 
Proposition~\ref{pr:Agmon1} under Hypothesis~\ref{hyp:mainf}, or
Hypothesis~\ref{hyp:AgmonLip} with
 $\Omega=f^{-1}([a,b])$ in the more general case, ensures
$\Delta_{f,\overline{\Omega},h}\geq c>0$ for $h\in]0,h_{0}[$\,.
 When   $\Delta_{f,\overline{\Omega}h}u=v\in
L^{2}(\Omega)$\,, it implies first $\|u\|\leq \frac{1}{c}\|v\|$\,.
We apply \eqref{eq:reQ1} with
$\varphi=0$:
$$
\|d_{f,\overline{\Omega},h}u\|^{2}+\|d_{f,\overline{\Omega},h}^{*}u\|^{2}=\Real\langle
u\,,\, \Delta_{f,\overline{\Omega},h}u\rangle\leq \|u\|\|v\|\leq \frac{1}{c}\|v\|^{2}\,.
$$
This proves the two first inequalities for $\|\eta_{1}\|_{L^{2}}$ and $\|\eta_{2}\|_{L^{2}}$\,.\\
Moreover, the equality
$$
d_{f,h}(\chi \omega)=\chi(d_{f,h}\omega)+(hd\chi)\wedge \omega
$$
implies
\begin{equation}
  \label{eq:dfhchidfh}
0=d_{f,h}[\chi (d_{f,h}\omega)]+d_{f,h}\left[(hd\chi)\wedge
  \omega\right]=(hd\chi)\wedge(d_{f,h}\omega)+d_{f,h}\left[(hd\chi)\wedge
  \omega\right]\,.
\end{equation}
Our assumptions ensure $(hd\chi)\wedge \omega\in
D(d_{f,\overline{\Omega},h})$ and $\eta_{1}\in
D(\Delta_{f,\overline{\Omega},h}) \subset D(d_{f,\overline{\Omega},h})$\,.
By using
$\Delta_{f,\overline{\Omega},h}=d_{f,\overline{\Omega},h}d_{f,\overline{\Omega},h}^{*}+d_{f,\overline{\Omega},h}^{*}d_{f,\overline{\Omega},h}$
and the commutation relation stated in Proposition~\ref{pr:domain}-4),
compute: 
\begin{eqnarray*}
  d_{f,\overline{\Omega},h}\eta_{1}&=&d_{f,\overline{\Omega},h}d_{f,\overline{\Omega},h}^{*}(\Delta_{f,\overline{\Omega},h})^{-1}(hd\chi\wedge\omega)\\
&=& (hd\chi)\wedge
    \omega-d_{f,\overline{\Omega},h}^{*}d_{f,\overline{\Omega},h}(\Delta_{f,\overline{\Omega},h})^{-1}((hd\chi)\wedge
    \omega)\\
&=& (hd\chi)\wedge \omega
    -d_{f,\overline{\Omega},h}^{*}(\Delta_{f,\overline{\Omega},h})^{-1}(d_{f,\overline{\Omega},h}[(hd\chi)\wedge
    \omega])\\
&\stackrel{\eqref{eq:dfhchidfh}}{=}&
(hd\chi)\wedge \omega
    +d_{f,\overline{\Omega},h}^{*}(\Delta_{f,\overline{\Omega},h})^{-1}((hd\chi)\wedge
    d_{f,h}\omega){\color{red}\,.}
\end{eqnarray*}
With
$d_{f,h}(\tilde{\chi}\eta_{1})=\tilde{\chi}(d_{f,h}\eta_{1})+(hd\tilde{\chi})\wedge
\eta_{1}$ and $\tilde{\chi}d\chi\equiv d\chi$\,, this implies:
$$
d_{f,h}(\tilde{\chi}\eta_{1})=(hd\chi)\wedge \omega
+\tilde{\chi}d_{f,\overline{\Omega},h}^{*}(\Delta_{f,\overline{\Omega},h})^{-1}((hd\chi)\wedge
d_{f,h}\omega)
+(hd\tilde{\chi})\wedge d_{f,\overline{\Omega},h}^{*}(\Delta_{f,\overline{\Omega},h})^{-1}((hd\chi)\wedge\omega)\,.
$$
We have proved
\begin{multline*}
d_{f,h}(\chi\omega -\tilde{\chi}\eta_{1})=
\chi(d_{f,h}\omega)
-\underbrace{\tilde{\chi}d_{f,\overline{\Omega},h}^{*}(\Delta_{f,\overline{\Omega},h})^{-1}((hd\chi)\wedge
d_{f,h}\omega)}_{(I)}
\\-\underbrace{(hd\tilde{\chi})\wedge
  d_{f,\overline{\Omega},h}^{*}(\Delta_{f,\overline{\Omega},h})^{-1}((hd\chi)\wedge\omega)
}_{(II)}\,.
\end{multline*}
Since
$\|d_{f,\overline{\Omega},h}^{*}(\Delta_{f,\overline{\Omega},h})^{-1}\|\leq
\frac 1{\sqrt{c}}$ for $h$ small enough, it follows
$$
\|(I)\|_{L^{2}}\leq \frac{1}{\sqrt{c}}\|(hd\chi)\wedge d_{f,h}\omega\|_{L^{2}}\,.
$$
For the last term, Proposition~\ref{pr:Agmon1} under
Hypothesis~\ref{hyp:mainf} says
$$
\|(II)\|_{L^{2}}=\tilde{O}(e^{-\frac{d_{Ag}(\supp d\tilde{\chi},\supp
    d\chi)}{h}})\|(hd\chi)\wedge \omega\|_{L^{2}}\,,
$$
while Hypothesis \ref{hyp:AgmonLip} with
$\overline{\Omega}=f^{-1}([a,b])$ in the more general case gives
$$
\|(II)\|_{L^{2}}=\tilde{O}(e^{-\frac{\min_{x\in \supp
      d\tilde{\chi},y\in \supp
    d\chi}|f(x)-f(y)|}{h}})\|(hd\chi)\wedge \omega\|_{L^{2}}\,,
$$
Meanwhile the identities
$d_{f,h}^{*}(\tilde{\chi}\eta_{1})=
\tilde{\chi}d_{f,h}^{*}\eta_{1}
+h\mathbf{i}_{\nabla\tilde{\chi}}\eta_{1}$  
and
$d_{f,h}^{*}\eta_{1}=0$ lead to
$$
d_{f,h}^{*}(\tilde{\chi}\eta_{1})=h\mathbf{i}_{\nabla\tilde{\chi}}\eta_{1}=
h\mathbf{i}_{\nabla\tilde{\chi}}d_{f,\overline{\Omega},h}^{*}(\Delta_{f,\overline{\Omega},h})^{-1}((hd\chi)\wedge
\omega)\,.
$$
which yields the fourth inequality.\\
Working with $\eta_{2}$ is completely symmetric by exchanging the role
of $d_{f,h}$ and $d_{f,h}^{*}$\,, after starting with
\begin{eqnarray*}
  &&
     d_{f,h}^{*}(\chi\omega)=\chi(d_{f,h}^{*}\omega)+h\mathbf{i}_{\nabla
     f}\omega\\
\text{and}&&
0=d_{f,h}^{*}(d_{f,h}^{*}\chi \omega)=
h\mathbf{i}_{\nabla_{\chi}}(d_{f,h}^{*}\omega)+d_{f,h}^{*}\left[
h\mathbf{i}_{\nabla
             \chi}d_{f,h}^{*}\omega\right]\,.
\end{eqnarray*}
The last inequality is obtained by summation.
\end{proof}

\subsubsection{Resolvent estimates}
\label{sec:resest}
From this paragraph and until the end of Section~6, the analysis
becomes essentially one dimensional along $\rz\supset f(M)$\,.
Accordingly we now work specifically with
$\overline{\Omega}=f^{-1}([a,b])$\,, $N_{t}=f^{-1}(a)$\,,
$N_{n}=f^{-1}(b)$\,, $a,b\not\in \left\{c_{1},\ldots,
  c_{N_{f}}\right\}$
or possibly $\overline{\Omega}=\sqcup_{n=1}^{N}f^{-1}([a_{n},b_{n}])$\,,
$a_{n},b_{n}\not\in\left\{c_{1},\ldots,c_{N_{f}}\right\}$\,, under
Hypothesis~\ref{hyp:mainf} for $f$\,, or by assuming
Hypothesis~\ref{hyp:Lipbar} and Hypothesis~\ref{hyp:AgmonLip} for a
more general Lipschitz function $f$\,.\\
 Also the upper bounds
$\tilde{O}(e^{-\frac{d_{Ag}(K,K')}{h}})$ in
Proposition~\ref{pr:Agmon}\,, Proposition~\ref{pr:Agmon1} and
Lemma~\ref{le:exten1} are replaced by their weaker form
$\tilde{O}(e^{-\frac{\inf_{x\in K, y\in K'}|f(x)-f(y)|}{h}})$ which is
the one given in Hypothesis~\ref{hyp:AgmonLip}.\\
We present here resolvent kernel estimates when $[a,b]$ contains one or a
fixed number $N$ of ``critical values''  of $f$\,. It assumes some
spectral localization, in   \eqref{eq:spass1} and \eqref{eq:spassN},
 which is not yet proved. It will be done in the
next sections with increasing complexity and precision: first for $N=1$ in
Section~\ref{sec:locpbs} and then
for a general $N$ in Section~\ref{sec:rough}, 
followed by the accurate
version for $N\geq 1$  in
 Section~\ref{sec:accanN}. It is also presented in
a more general form where actually the $N$ critical values may be
replaced by $N$ clusters of critical values for further applications.\\
Let us first consider the case when $[a,b]$ contains one cluster of
``critical values''.
\begin{prop}
\label{pr:1well}
Assume Hypothesis~\ref{hyp:mainf}, or more generally
Hypothesis~\ref{hyp:Lipbar} and Hypothesis~\ref{hyp:AgmonLip}, for $f$ and let $a<c<b$ and  $\varepsilon_{0}\in ]0,\min(b-c,c-a)[$ 
be such that
$$
[a,b]\cap \left\{c_{1},\ldots,c_{N_{f}}\right\}\subset ]c-\frac{\varepsilon_{0}}{16},c+\frac{\varepsilon_{0}}{16}[\,.
$$
Assume also that $\Delta_{f,f^{-1}([a,b]),h}$\,, the self-adjoint
operator in $f^{-1}([a,b])\subset M$ given in
Proposition~\ref{pr:domain}
 with
$N_{t}=\{f=a\}$ and  $N_{n}=\{f=b\}$ satisfies:
\begin{equation}
\label{eq:spass1}
\exists h_{0}>0\,,\ \forall h\in ]0,h_{0}[\,,\quad
\sigma(\Delta_{f,f^{-1}([a,b]),h})\cap [0,e^{-\frac{\varepsilon_{0}}{h}}]
\subset[0,e^{-\frac{4\varepsilon_{0}}{h}}]\,.
\end{equation}
Then the estimate
\begin{equation*}
(\Delta_{f,f^{-1}([a,b]),h}-z)^{-1}(x,y)
=\tilde{O}(e^{-\frac{|f(x)-f(y)|}{h}+\frac{3\varepsilon_{0}}{h}})
\end{equation*}
holds, according to Definition~\ref{de:Otkernel}, uniformly with respect
to $z$\,, $|z|=e^{-\frac{2\varepsilon_{0}}{h}}$\,.
\end{prop}
\begin{proof}
We prove Proposition~\ref{pr:1well} by adapting the analysis made in
\cite[pp.~57--58]{DiSj}.
Let us  consider the self-adjoint realizations
$\Delta_{f,f^{-1}([a,c-\frac{\varepsilon_{0}}{16}]),h}$ 
and $\Delta_{f,f^{-1}([c+\frac{\varepsilon_{0}}{16},b]),h}$ for which
Proposition~\ref{pr:Agmon1}
says
\begin{equation}
\label{eq.res-outside}
(\Delta_{f,f^{-1}([a,c-\frac{\varepsilon_{0}}{16}]),h}-z)^{-1}(x,y)
\ \ \text{and}\ \ 
(\Delta_{f,f^{-1}([c+\frac{\varepsilon_{0}}{16},b]),h}-z)^{-1}(x,y)
\ \ \text{are}\ \
\tilde{O}(e^{-\frac{|f(x)-f(y)|}{h}})\,,
\end{equation}
 uniformly with respect to 
 $z\in \cz,|z|=e^{-2\frac{\varepsilon_{0}}{h}}$\,.
Let moreover $\theta$ and $\hat\theta$ be two cut-off functions such that
$\theta\in C^{\infty}_{0}(f^{-1}(]c-\frac{\varepsilon_{0}}{8},c+\frac{\varepsilon_{0}}{8}[);[0,1])$\,, 
$\theta\equiv 1$ around $f^{-1}([c-\frac{\varepsilon_{0}}{16},c+\frac{\varepsilon_{0}}{16}])$\,, 
and $\hat\theta\in C^{\infty}_{0}(f^{-1}(]c-\frac{3\varepsilon_{0}}{16},c+\frac{3\varepsilon_{0}}{16}[);[0,1])$\,, 
$\hat\theta\equiv 1$ around $f^{-1}([c-\frac{\varepsilon_{0}}{8},c+\frac{\varepsilon_{0}}{8}])$\,.
Let us also define $\theta_{-},\hat\theta_{-} \in C^{\infty}(f^{-1}(]-\infty,c[);[0,1])$
and  $\theta_{+},\hat\theta_{+} \in C^{\infty}(f^{-1}(]c,+\infty[);[0,1])$
such that
$$
\theta_{-}+\theta+\theta_{+}=1\quad\text{and}\quad \hat\theta_{-}+\hat\theta+\hat\theta_{+}=1\,.
$$
\figinit{0.9cm}
\figpt 0:(-6,-1)
\figpt 1000: (8,-1)
\figpt 1001:(-7,-1)
\figpt 1002:(7.5,-1)
\figpt  1:(4,-1)
\figpt 3: (3,-1)
\figpt 11:(2.8,-1)
\figpt 111:(2.6,-1)
\figpt 12:(2.5,-0.5)
\figpt 13:(2.4,0.5)
\figpt 14:(2.35,1)
\figpt 114:(2.1,1)
\figpt 115:(2,1)
\figpt 15:(-6,1)
\figpt 5:(0,-1)

\figpt 2115:(-2,1)
\figpt 2114:(-2.1,1)
\figpt 2014:(-2.35,1)
\figpt 2013:(-2.4,0.5)
\figpt 2012:(-2.5,-0.5)
\figpt 2111:(-2.6,-1)
\figpt 2011:(-2.8,-1)
\figpt 2003:(-3,-1)
\figpt 2001:(-4,-1)

\figpt  10001:(4,1)
\figpt 10003: (3,1)
\figpt 10011:(2.8,1)
\figpt 10111:(2.6,1)
\figpt 10012:(2.5,0.5)
\figpt 10013:(2.4,-0.5)
\figpt 10014:(2.35,-1)
\figpt 10114:(2.1,-1)
\figpt 10115:(2,-1)
\figpt 11001:(7.5,1)

\figpt 12115:(-2,-1)
\figpt 12114:(-2.1,-1)
\figpt 12014:(-2.35,-1)
\figpt 12013:(-2.4,-0.5)
\figpt 12012:(-2.5,0.5)
\figpt 12111:(-2.6,1)
\figpt 12011:(-2.8,1)
\figpt 12003:(-3,1)
\figpt 12001:(-4,1)
\figpt 12000:(-6,1)

\figpt 203:(-2,-1)
\figpt 211:(-1.8,-1)
\figpt 311:(-1.6,-1)
\figpt 212:(-1.5,-0.5)
\figpt 213:(-1.4,0.5)
\figpt 214:(-1.36,0.8)
\figpt 314:(-1.1,0.8)
\figpt 315:(-1,0.8)
\figpt 205:(-1,-1)
\figpt 2205:(1,-1)
\figpt 2315:(1,0.8)
\figpt 2314:(1.1,0.8)
\figpt 2214:(1.36,0.8)
\figpt 2213:(1.4,0.5)
\figpt 2212:(1.5,-0.5)
\figpt 2311:(1.6,-1)
\figpt 2211:(1.8,-1)
\figpt 2203:(2,-1)
\figpt 10203:(-2,0.8)
\figpt 10211:(-1.8,0.8)
\figpt 10311:(-1.6,0.8)
\figpt 10212:(-1.5,0.3)
\figpt 10213:(-1.4,-0.7)
\figpt 10214:(-1.36,-1)
\figpt 10314:(-1.1,-1)
\figpt 10315:(-1,-1)
\figpt 10205:(-1,0.8)
\figpt 12200:(-6,0.8)

\figpt 12205:(1,0.8)
\figpt 12315:(1,-1)
\figpt 12314:(1.1,-1)
\figpt 12214:(1.36,-1)
\figpt 12213:(1.4,-0.7)
\figpt 12212:(1.5,0.3)
\figpt 12311:(1.6,0.8)
\figpt 12211:(1.8,0.8)
\figpt 12203:(2,0.8)
\figpt 12002:(7.5,0.8)

\figdrawbegin{}
\figset(width=1.3)
\figdrawline[1000,1001]
\figset(width=3)
\figdrawBezier 1[1,1,3,11]
\figdrawBezier 1[11,111,12,13]
\figdrawBezier 1[13,14,114,115]
\figdrawBezier 1[2001,2001,2003,2011]
\figdrawBezier 1[2011,2111,2012,2013]
\figdrawBezier 1[2013,2014,2114,2115]
\figdrawline[115,2115]
\figdrawline[11,1002]
\figdrawline[2011,0]
\figset(width=1)
\figdrawBezier 1[10001,10001,10003,10011]
\figdrawBezier 1[10011,10111,10012,10013]
\figdrawBezier 1[10013,10014,10114,10115]
\figdrawline[10001,11001]
\figdrawBezier 1[12001,12001,12003,12011]
\figdrawBezier 1[12011,12111,12012,12013]
\figdrawBezier 1[12013,12014,12114,12115]
\figdrawline[12001,12000]
\figset(width=3)
\figset(dash=3)
\figdrawline[0,211]
\figdrawBezier 1[211,311,212,213]
\figdrawBezier 1[213,214,314,315]
\figdrawline[315,2315]
\figdrawBezier 1[2211,2311,2212,2213]
\figdrawBezier 1[2213,2214,2314,2315]
\figdrawline[2211,1002]
\figset(width=1.)
\figdrawBezier 1[10211,10311,10212,10213]
\figdrawBezier 1[10213,10214,10314,10315]
\figdrawline[10211,12200]
\figdrawBezier 1[12211,12311,12212,12213]
\figdrawBezier 1[12213,12214,12314,12315]
\figdrawline[12211,12002]
\figdrawend{}
\figvisu{\figBoxB}{\textbf{Figure 1:} Positions of the cut-off
functions, $\theta_{-}, \theta, \theta_{+}, \hat{\theta}_{-},
\hat{\theta}, \hat{\theta}_{+}$\,.}{
\figset write(mark=$+$)
\figwrites 0:$a$(0.3)
\figwrites 3:$c+\frac{3\varepsilon_{0}}{16}$(0.7)
\figwrites 2003:$c-\frac{3\varepsilon_{0}}{16}$(0.7)
\figwrites 5:$c$(0.3)
\figwrites 205:$c-\frac{\varepsilon_{0}}{16}$(0.7)
\figwrites 2205:$c+\frac{\varepsilon_{0}}{16}$(0.7)
\figwrites 203:$c-\frac{\varepsilon_{0}}{8}$(0.2)
\figwrites 2203:$c+\frac{\varepsilon_{0}}{8}$(0.2)
\figwrites 1002:$b$(0.3)
\figset write(mark={})
\figwriten 115:$\hat{\theta}$(0.3)
\figwritese 214:$\theta$(0.3)
\figwritesw 12002:$\theta_{+}$(0.3)
\figwritese 12200:$\theta_{-}$(0.3)
\figwritene 12000:$\hat{\theta}_{-}$(0.3)
\figwritenw 11001:$\hat{\theta}_{+}$(0.3)
}
\centerline{\box\figBoxB}

\medskip
The support conditions imply the following resolvent identity:
\begin{align*}
(\Delta_{f,f^{-1}([a,b]),h}-z)^{-1}&=
(\Delta_{f,f^{-1}([a,b]),h}-z)^{-1}\hat\theta
+\theta_{-}(\Delta_{f,f^{-1}([a,c-\frac{\varepsilon_{0}}{16}]),h}-z)^{-1} \hat\theta_{-}\\
&-(\Delta_{f,f^{-1}([a,b]),h}-z)^{-1} \hat\theta[\Delta_{f,h},\theta_{-}]
(\Delta_{f,f^{-1}([a,c-\frac{\varepsilon_{0}}{16}]),h}-z)^{-1} \hat\theta_{-}\\
&+\theta_{+}(\Delta_{f,f^{-1}([c+\frac{\varepsilon_{0}}{16},b]),h}-z)^{-1} \hat\theta_{+}\\
&-(\Delta_{f,f^{-1}([a,b]),h}-z)^{-1} \hat\theta[\Delta_{f,h},\theta_{+}]
(\Delta_{f,f^{-1}([c+\frac{\varepsilon_{0}}{16},b]),h}-z)^{-1} \hat\theta_{+}\,.
\end{align*}
Since moreover
$\|(\Delta_{f,f^{-1}([a,b]),h}-z)^{-1}\|_{\mathcal L(L^{2},L^{2})}\leq 2
e^{2\frac{\varepsilon_{0}}{h}}$
for $|z|=e^{-\frac{2\varepsilon_{0}}{h}}$\,, because the hypothesis
 ensures
 $\mathrm{dist}_{\cz}(z,
\sigma(\Delta_{f,f^{-1}([a,b]),h}))\leq
\frac{e^{-\frac{2\varepsilon_{0}}{h}}}{2}$ for
 $h>0$ small enough,
 applying Proposition~\ref{pr:Agmon} to 
$$
(\Delta_{f,f^{-1}([a,b]),h}-z)\omega_{h}=r_{h}=\hat{\theta}\hat{r}_{h}
$$
with $\supp \hat\theta\subset f^{-1}(]c-3\frac{\varepsilon_{0}}{16},c+3\frac{\varepsilon_{0}}{16}[)$
first yields
$$
[(\Delta_{f,f^{-1}([a,b]),h}-z)^{-1}\hat\theta](x,y)\ =\ 
\tilde{O}(e^{-\frac{|f(x)-c|-3\varepsilon_{0}/16}{h}+2\frac{\varepsilon_{0}}{h}})
$$
and then
\begin{eqnarray*}
[(\Delta_{f,f^{-1}([a,b]),h}-z)^{-1}\hat\theta](x,y)&=& 
\tilde{O}(e^{-\frac{|f(x)-c|-3\varepsilon_{0}/16}{h}+2\frac{\varepsilon_{0}}{h}})\,
\tilde{O}(e^{-\frac{|f(y)-c|-3\varepsilon_{0}/16}{h}})
\\
&=&\ \tilde{O}(e^{-\frac{|f(x)-f(y)|}{h}+3\frac{\varepsilon_{0}}{h}}) \,.
\end{eqnarray*}
By using \eqref{eq.res-outside},
$\|[\Delta_{f,h},\theta_{\pm}]\|_{\mathcal{L}(W^{1,2};L^{2})}=\tilde{O}(1)$\,,
$[\Delta_{f,h},\theta_{\pm}]$  vanishing in a neighborhood of
$f^{-1}(\left\{a,b\right\})$\,, and the latter  estimate   for all the left factors concerned  in the above
resolvent identity, we obtain
\begin{eqnarray*}
(\Delta_{f,f^{-1}([a,b]),h}-z)^{-1}(x,y)&=&
\tilde{O}(e^{-\frac{|f(x)-f(y)|}{h}+\frac{3\varepsilon_{0}}{h}})+
\tilde{O}(e^{-\frac{|f(x)-f(y)|}{h}})\\
&=&\tilde{O}(e^{-\frac{|f(x)-f(y)|}{h}+\frac{3\varepsilon_{0}}{h}})\,.
\end{eqnarray*}
\end{proof}

\begin{prop}
\label{pr:Nwell}
Assume Hypothesis~\ref{hyp:mainf}, or more generally
Hypothesis~\ref{hyp:Lipbar} and Hypothesis~\ref{hyp:AgmonLip} for
$f$\,. Let $a<b$ belong to
$\overline{\rz}\setminus\left\{c_{1},\ldots, c_{N_{f}}\right\}$ and
let $\overline{\Omega}=f^{-1}([a,b])$ with
$N_{t}=f^{-1}(\left\{a\right\})$\,,
$N_{n}=f^{-1}(\left\{b\right\})$\,. 
Assume that there exist
$a=\tilde{c}_{0}<\tilde{c}_{1}<\ldots<\tilde{c}_{N}<\tilde{c}_{N+1}=b$
and $\varepsilon_{0}\in]0, \frac{\min_{1\leq n\leq
    N+1}(\tilde{c}_{n}-\tilde{c}_{n-1})}{16}[ $ such that
$$
]a,b[\cap \left\{c_{1},\ldots,c_{N_{f}}\right\}\subset 
\sqcup_{n=1}^{N}]\tilde{c}_{n}-\frac{\varepsilon_{0}}{16},\tilde{c}_{n}+\frac{\varepsilon_{0}}{16}[
\,.
$$
The operator $\Delta_{f,f^{-1}([a,b]),h}$ is the self-adjoint
realization of the Witten Laplacian given in
Proposition~\ref{pr:domain} and accordingly $
     \Delta_{n}=\Delta_{f,f^{-1}([\tilde{c}_{n-1}+(1-\delta_{n,1})\varepsilon_{0},\tilde{c}_{n+1}-(1-\delta_{n,N})\varepsilon_{0}]),h}$
     is defined for $1\leq n\leq N$
     where $\delta_{m,n}$ is the Kronecker symbol. We  assume
     \begin{equation}
\label{eq:spassN}
\forall n \in \left\{1,\ldots, N\right\}\,,\quad \sigma(\Delta_{n})\cap
[0,e^{-\frac{\varepsilon_{0}}{h}}]\subset [0,e^{-\frac{4\varepsilon_{0}}{h}}]\,.
\end{equation}
Then every $z\in \cz$ such that $|z|=e^{-\frac{2\varepsilon_{0}}{h}}$
belongs to the resolvent set of $\Delta_{f,f^{-1}([a,b]),h}$ provided that
$h\in ]0,h_{0}[$ with $h_{0}>0$ small enough.
Moreover, there exists  a constant
$N_{0}\in\N^{*}$\,, determined by $b-a$ and $\min_{2\leq n\leq
  N}\tilde{c}_{n}-\tilde{c}_{n-1}$\,, such that
\begin{align*}
(\Delta_{f,f^{-1}([a,b]),h}-z)^{-1}(x,y)
&=\tilde{O}(e^{-\frac{|f(x)-f(y)|}{h}+3N_{0}\frac{\varepsilon_{0}}{h}})
\end{align*}
holds, according to Definition~\ref{de:Otkernel}
 uniformly with respect to $z$\,, $|z|=e^{-\frac{2\varepsilon_{0}}{h}}$\,.
\end{prop}

\begin{proof}
We prove Proposition~\ref{pr:Nwell} by adapting the analysis made in
\cite[pp.~58--59]{DiSj}. Call $\eta_{0}=\min_{2\leq n
  \leq N}\frac{\tilde{c}_{n}-\tilde{c}_{n-1}}{2}$ and take $\varepsilon_{0}\in
]0,\min_{1\leq n\leq N+1}\frac{\tilde{c}_{n}-\tilde{c}_{n-1}}{16}[$\,,
$\varepsilon_{0}\leq \frac{\eta_{0}}{8}$ as stated.\\
For $n\in\{1,\dots,N\}$\,, let us 
introduce $\theta_{n}\in C^{\infty}_{0}(f^{-1}(]\tilde{c}_{n}-\frac{\varepsilon_{0}}{8},\tilde{c}_{n}+\frac{\varepsilon_{0}}{8}[);[0,1])$
such that $\theta_{n}\equiv 1$ in a neighborhood of 
 $f^{-1}([\tilde{c}_{n}-\frac{\varepsilon_{0}}{16},\tilde{c}_{n}+\frac{\varepsilon_{0}}{16}])$\,,
and
$$
\chi_{n}\ :=\ \big(1-\sum_{m\neq n} \theta_{m}\big)\big|_{f^{-1}([\tilde{c}_{n-1},\tilde{c}_{n+1}])}
\ =\ 
\big(1- \theta_{n-1}-\theta_{n+1}\big)\big|_{f^{-1}([\tilde{c}_{n-1},\tilde{c}_{n+1}])}
\,.
$$
Here, we use the convention $\theta_{-1}=\theta_{N+1}=0$\,.
We also need another partition of unity
$1=\sum_{n=1}^{N}\tilde{\chi}_{n}$\,, $0\leq \tilde{\chi}_{n}\leq
1$\,, such that
\begin{eqnarray*}
  && \tilde{\chi}_{n}\equiv
     1~\text{on}~f^{-1}([\tilde{c}_{n}-\eta_{0}/2,\tilde{c}_{n}+\eta_{0}/2])\quad
     \text{for}~1\leq n\leq N\,,\\
&&\tilde{\chi}_{n}\in
\mathcal{C}^{\infty}_{0}(f^{-1}(]\tilde{c}_{n-1}+\eta_{0}/2,\tilde{c}_{n+1}-\eta_{0}/2[))\quad\text{for}~2\leq
n\leq N-1\,\\
\text{and}&&
\tilde{\chi}_{1}\equiv 0~\text{on}~f^{-1}([\tilde{c}_{2}-\eta_{0}/2,b])\quad
\tilde{\chi}_{N}\equiv 0~\text{on}~f^{-1}([a,\tilde{c}_{N-1}+\eta_{0}/2])\,.
\end{eqnarray*}
Note in particular that our conditions, $\varepsilon_{0}\leq
\frac{\eta_{0}}{8}$ and  $\supp \theta_{n}\subset
f^{-1}([\tilde{c}_{n}-\frac{\varepsilon_{0}}{8},\tilde{c}_{n}+\frac{\varepsilon_{0}}{8}])$\,,
 ensure
$\chi_{n}\equiv 1$ on $\supp\tilde{\chi}_{n}$\,.\\
We now set for every 
$z\in \cz,|z|=e^{-2\frac{\varepsilon_{0}}{h}}$:
\begin{equation}
\label{eq.R0}
R_{0}(z)\ :=\ \sum_{n=1}^{N} \chi_{n}(\Delta_{n}-z)^{-1}\tilde \chi_{n}\,,
\end{equation}
where we recall
$\Delta_{n}=\Delta_{f,f^{-1}([\tilde{c}_{n-1}+(1-\delta_{n,1})\varepsilon_{0},\tilde{c}_{n+1}-\delta_{N,n}\varepsilon_{0}]),h}$\,.
Because the boundary conditions are satisfied\,, a simple computation shows
\begin{equation}
\label{eq.Delta-R0}
(\Delta_{f,f^{-1}([a,b]),h}-z)R_{0}\ =\ I-K\,,
\end{equation}
with
\begin{equation}
\label{eq.K}
K\ =\ \sum_{n=1}^{N}\sum_{m\in\{n-1,n+1\}} [\Delta_{f,h},\theta_{m}]\big|_{f^{-1}([\tilde{c}_{n-1},\tilde{c}_{n+1}])}\,(\Delta_{n}-z)^{-1}\tilde \chi_{n}\,.
\end{equation}
Moreover Proposition~\ref{pr:1well} applied to every $\Delta_{n}$ and 
\eqref{eq.K} combined with the support conditions of $\theta_{m},
\tilde{\chi}_{n}$ imply
$$
\|K\|_{\mathcal L(L^{2},L^{2})}\ =\ \tilde{O}(e^{-\frac{C}{h}+\frac{3\varepsilon_{0}}{h}})\,,
$$
where
$$
C:=\min_{\scriptsize
  \begin{array}[c]{l}
n\in\{1,\dots,N\}\\
m\in\{n-1,n+1\}
  \end{array}
}\big(\min_{
  \scriptsize\begin{array}[c]{l}
    y\in \supp \tilde \chi_{n} \\
x\in \supp \theta_{m}
  \end{array}
}|f(x)-f(y)|\big)\geq \frac{\eta_{0}}{2}-\frac{\varepsilon_{0}}{8}\,,
$$
and $\varepsilon_{0}\leq \frac{\eta_{0}}{8}$\,, yields
$$
\|K\|_{\mathcal{L}(L^{2};L^{2})}=\tilde{O}(e^{-\frac{\eta_{0}/2-25\varepsilon_{0}/8}{h}})
=\tilde{O}(e^{-\frac{7\eta_{0}}{64 h}})\,.
$$
 For $h>0$ small
enough, $I-K:L^{2}\to L^{2}$ in then invertible and the resolvent set
of $\Delta_{f,f^{-1}([a,b]),h}$ contains $\left\{z\in \cz\,, |z|=
  e^{-\frac{2\varepsilon_{0}}{h}}\right\}$\,.\\
Let us now consider the exponential decay
estimate.  Write first
\begin{equation}
\label{eq.res-R0}
(\Delta_{f,f^{-1}([a,b]),h}-z)^{-1}\ =\ R_{0}(z)\sum_{\ell\in\N}K^{\ell}
\ =\ R_{0}(z)\sum_{\ell=0}^{N_{0}-1}K^{\ell}+ R_{0}(z) K_{N_{0}}\,,
\end{equation}
and choose $N_{0}\in\N^{*}$ such that
$N_{0}\times\frac{7\eta_{0}}{64}\geq (b-a)$ and
\begin{equation}
\label{eq.K0+1}
\|K_{N_{0}}\|_{\mathcal L(L^{2},L^{2})}\ =\
\|\sum_{\ell\geq N_0}K^{\ell}\|_{\mathcal L(L^{2},L^{2})}
\ =\ \tilde{O}(e^{-\frac{b-a}{h}})=\tilde{O}(e^{-\frac{\max_{x,y\in f^{-1}([a,b])}|f(x)-f(y)|}{h}})\,.
\end{equation}
By referring again to  Proposition~\ref{pr:1well}
and from the definition \eqref{eq.R0} or $R_{0}(z)$\,, we know:
\begin{equation}
\label{eq.R0xy}
R_{0}(z)(x,y)\ =\ \tilde{O}(e^{-\frac{|f(x)-f(y)|}{h}+3\frac{\varepsilon_{0}}{h}})\,.
\end{equation}
The relation \eqref{eq.R0xy} together with \eqref{eq.K0+1}
implies that
\begin{equation}
\label{eq.R0-K0+}
(R_{0}\circ K_{N_{0}})(x,y)\ =\ \tilde{O}(e^{-\frac{\min_{z\in M}|f(x)-f(z)|+b-a}{h}+3\frac{\varepsilon_{0}}{h}})
\ =\ \tilde{O}(e^{-\frac{|f(x)-f(y)|}{h}+3\frac{\varepsilon_{0}}{h}})\,.
\end{equation}
Moreover, the relation \eqref{eq.R0xy} together with
$$
K(x,y)\ =\ \tilde{O}(e^{-\frac{|f(x)-f(y)|}{h}+3\frac{\varepsilon_{0}}{h}})\,,
$$
which follows as well from Proposition~\ref{pr:1well}, implies that for every
$\ell\in\N$\,, one has:
\begin{equation}
\label{eq.R0-Kell}
(R_{0}(z)\circ K^{\ell})(x,y)\ =\ \tilde{O}(e^{-\frac{|f(x)-f(y)|}{h}+3(\ell+1)\frac{\varepsilon_{0}}{h}})\,.
\end{equation}
One finally deduces from \eqref{eq.res-R0}
and from \eqref{eq.R0-K0+}, \eqref{eq.R0-Kell} that the estimate
$$
(\Delta_{f,f^{-1}([a,b]),h}-z)^{-1}(x,y)\ =\ \tilde{O}(e^{-\frac{|f(x)-f(y)|}{h}+3N_{0}\frac{\varepsilon_{0}}{h}})\,,
$$
holds uniformly with respect to  $z\in
\cz,|z|=e^{-2\frac{\varepsilon_{0}}{h}}$\,.
This concludes the proof of Proposition~\ref{pr:Nwell}.
\end{proof}

\section{Local problems}
\label{sec:locpbs}

In this section we shall use Agmon type estimates to study carefully
the case when there is a unique ``critical value'' of $f$ in $]a,b[$\,,
$-\infty\leq a<b\leq +\infty$\,.
\begin{hyp}
\label{hyp:1vc}
The function $f$ is assumed to satisfy Hypothesis~\ref{hyp:mainf}, or
Hypothesis~\ref{hyp:Lipbar} and Hypothesis~\ref{hyp:AgmonLip},  and
the values $a,b$\,, $-\infty\leq a<b\leq +\infty$\,, are chosen such
that 
$$
[a,b]\cap \left\{c_{1},\ldots, c_{N_{f}}\right\}=]a,b[\cap \left\{c_{1},\ldots,c_{N_{f}}\right\}=\left\{\tilde{c}_{1}\right\}\,.
$$
The domain is $\overline{\Omega}=f^{-1}([a,b])$\,, with
$N_{t}=f^{-1}(\left\{a\right\})$ and $N_{n}=f^{-1}(\left\{b\right\})$\,,
and the operator $\Delta_{f,f^{-1}([a,b]),h}$ is the one defined in Proposition~\ref{pr:domain}.
\end{hyp}
 With this assumption
all the exponential decay estimates of Section~\ref{sec:agmontype} can
be used with the pseudodistance $|f(x)-f(y)|$\,.  
The main result of this section says that, in this
framework, the only possible exponentially small eigenvalue of
$\Delta_{f,f^{-1}([a,b]),h}$ is $0$\,.
\begin{prop}
\label{pr:exp0}
Under Hypothesis~\ref{hyp:1vc}, the spectrum of the operator $\Delta_{f,f^{-1}([a,b]),h}$  satisfies
$$
\forall \varepsilon>0\,, \exists h_{\varepsilon}>0\,, \forall h\in
]0,h_{\varepsilon}[\,, \sigma(\Delta_{f,f^{-1}([a,b]),h})\cap [0,e^{-\frac{\varepsilon}{h}}]\subset\left\{0\right\}\,.
$$
\end{prop}
Proposition~\ref{pr:exp0} will be proved in several steps. 
Consequences e.g. for resolvent
estimates will be given afterwards.
\subsection{Useful quantities and notations}
\label{sec:usqunot}
Let us first recall the following notion of distance between (spectral)
subspaces which is convenient for spectral analysis (see e.g. \cite[pp.~59--61]{DiSj}).
\begin{definition}
\label{de:dEF}
For $E,F$  two closed subspaces of a Hilbert space $\mathcal{H}$\,, the
non symmetric distance $\vec{d}(E,F)$ is defined as 
$$
\vec{d}(E,F)=\sup_{x\in E,\|x\|=1}d_{\mathcal H}(x,F)=\|\Pi_{E}-\Pi_{F}\Pi_{E}\|=\|\Pi_{E}-\Pi_{E}\Pi_{F}\|\,,
$$
where $\Pi_{E},\Pi_{F}$ are the orthogonal projection on $E$\,,$F$\,. 
\end{definition}
This distance satisfies:
\begin{itemize}
\item $\vec{d}(E,F)=0$ iff $E\subset F$\,;
\item $\vec{d}(E,G)\leq \vec{d}(E,F)+\vec{d}(F,G)$\,;
\item $\vec{d}(E,F)<1$ if and only if $\Pi_{F}\big|_{E}:E\to F$
  is one-to-one with a continuous left-inverse, and $\Pi_{E}\big|_{F}:F\to E$ is onto in this case;
\item  $\big(\,\vec{d}(E,F)<1$ and  $\vec{d}(F,E)<1\,\big)$ if and only if
  $\Pi_{F}\big|_{E}:E\to F$ and $\Pi_{E}\big|_F:F\to E$ are bijections
  with continuous inverses. In this case, the equality
  $\vec{d}(E,F)=\vec{d}(F,E)$ holds true\,;
\item if we know a priori $\dim E=\dim F<+\infty$ then
$$
(\vec{d}(E,F))<1)\Leftrightarrow
\left(\vec{d}(E,F)<1\quad\text{and}\quad
  \vec{d}(F,E)<1\right)\Leftrightarrow (\vec{d}(F,E)<1)\,.
$$
\end{itemize}
We will use a variation of the min-max principle associated with the
quantities $\gamma(\alpha,[a,b],h)$ and $\Gamma(\alpha,[a,b],h)$
defined below.
Remember that $Q^{(p)}_{f,f^{-1}([a,b]),h}$ is the quadratic form
associated with $\Delta_{f,f^{-1}([a,b]),h}$
(see the second item of Proposition~\ref{pr:domain}).
\begin{definition}
\label{de:gam}
For $p\in\left\{0,\ldots d\right\}$\,, $s\geq 0$\,, 
let $F^{(p)}_{[0,s],[a,b],h}$ denote the range of the spectral
projection  $1_{[0,s]}(\Delta_{f,f^{-1}([a,b]),h}^{(p)})$\,,  with in particular $F^{(p)}_{\left\{0\right\},f^{-1}([a,b]),h}=\ker(\Delta_{f,f^{-1}([a,b]),h}^{(p)})$\,.\\
For $\alpha>0$\,, the quantities
$\gamma^{(p)}(\alpha,[a,b],h)$ and $\Gamma^{(p)}(\alpha,[a,b],h)$ are
defined by
\begin{eqnarray}
\nonumber
\gamma^{(p)}(\alpha,[a,b],h)
&=&
    \vec{d}(F^{(p)}_{[0,e^{-\frac{\alpha}{h}}],[a,b],h},F^{(p)}_{\left\{0\right\},[a,b],h})
=
 \vec{d}(F^{(p)}_{[0,e^{-\frac{\alpha}{h}}],[a,b],h},\ker(\Delta_{f,f^{-1}([a,b]),h}^{(p)}))\\
\label{eq:gam}
&=&
\sup_{\omega_{h}\in F^{(p)}_{[0,e^{-\frac{\alpha}{h}}],[a,b],h}\setminus\{0\}}
\frac{\textrm{dist}_{L^{2}}(\omega_{h},\ker{\Delta_{f,f^{-1}([a,b]),h}^{(p)}})}{\|\omega_{h}\|_{L^{2}}}\,,\\
\label{eq:Gam}
\Gamma^{(p)}(\alpha,[a,b],h)&=&
\sup_{\|\omega_{h}\|_{L^{2}}=1\ :\ 
     Q_{f,f^{-1}([a,b]),h}^{(p)}(\omega_{h})\leq e^{-\frac{\alpha}{h}}
}
\textrm{dist}_{L^{2}}(\omega_{h},\ker(\Delta_{f,f^{-1}([a,b]),h}^{(p)}))\,.
\end{eqnarray}
\end{definition}
Those quantities satisfy simple properties:
\begin{itemize}
\item The quantities $\gamma^{(p)}(\alpha,[a,b],h)$
and 
  $\Gamma^{(p)}(\alpha,[a,b],h)$ are  decreasing w.r.t $\alpha$ and, since 
$$
F^{(p)}_{[0,e^{-\frac{\alpha}{h}}],[a,b],h}\subset\{\omega\in D(Q_{f,f^{-1}([a,b]),h})\ \text{s.t.}\ 
  Q_{f,f^{-1}([a,b]),h}^{(p)}(\omega_{h})\leq
  e^{-\frac{\alpha}{h}}\}\,,
$$
  they satisfy
$$
0\ \leq\ \gamma^{(p)}(\alpha,[a,b],h)\ \leq\ 
\Gamma^{(p)}(\alpha,[a,b],h)\,.
$$
It says in particular:
$$
\left(\lim_{h\to
  0}\Gamma^{(p)}(\alpha,[a,b],h)=0\right)\Rightarrow \left(\lim_{h\to
    0}\gamma^{(p)}(\alpha,[a,b],h)=0\right)\,.
$$

\item Since $\Delta_{f,f^{-1}([a,b]),h}$ is self-adjoint, the spectral
  theorem implies:
$$\gamma^{(p)}(\alpha,[a,b],h)=0 \quad \text{iff}\quad 
\sigma(\Delta_{f,f^{-1}([a,b]),h}^{(p)})\cap [0,e^{-\frac{\alpha}{h}}]\subset\left\{0\right\}
$$
and
$$ \gamma^{(p)}(\alpha,[a,b],h)=1 \quad \text{else}.$$
In particular, it provides the expression
$$
\gamma^{(p)}(\alpha,[a,b],h)
=
\sup_{\|\omega_{h}\|=1\ :\ 
    \left\{\begin{array}[c]{l}
     \Delta_{f,f^{-1}([a,b]),h}^{(p)}\omega_{h}=\lambda_{h}\omega_{h}\\
     \lambda_{h}\leq e^{-\frac{\alpha}{h}} 
    \end{array}\right.}
\textrm{dist}_{L^{2}}(\omega_{h},\ker{\Delta_{f,f^{-1}([a,b]),h}^{(p)}})
$$
and the convergence $\lim_{h\to 0}\gamma^{(p)}(\alpha,[a,b],h)=0$
  means precisely that:
\begin{equation}
\label{eq.gamma-equiv}
\exists h_{\alpha}>0\,, \forall h\in
]0,h_{\alpha}[\,, \sigma(\Delta_{f,f^{-1}([a,b],h)}^{(p)})\cap [0,e^{-\frac{\alpha}{h}}]\subset\left\{0\right\}\,.
\end{equation}
\item  The spectral theorem also implies
$$\Gamma^{(p)}(\alpha,[a,b],h)=1 \quad \text{iff}\quad 
\sigma(\Delta_{f,f^{-1}([a,b],h)}^{(p)})\cap ]0,e^{-\frac{\alpha}{h}}]\neq\emptyset$$
and
\begin{equation}
\label{eq.Gamma-equiv} 
 \big(\Gamma^{(p)}(\alpha,[a,b],h)\big)^{2}\in[0,\frac{e^{-\frac{\alpha}{h}}}{\min\big(\sigma(\Delta_{f,f^{-1}([a,b],h)}^{(p)})\setminus\{0\}\big)}]\subset
 [0,1[\quad  \text{else}\,.
\end{equation}
 Actually, $\sigma(\Delta_{f,f^{-1}([a,b],h)}^{(p)})\cap
 ]0,e^{-\frac{\alpha}{h}}]\neq\emptyset$ implies $\Gamma^{(p)}(\alpha,[a,b],h)\geq
 \gamma^{(p)}(\alpha,[a,b],h)\geq 1$ and obviously
 $\Gamma^{(p)}(\alpha,[a,b],h)=1$\,.\\
Reciprocally when $\sigma(\Delta_{f,f^{-1}([a,b]),h})\cap
]0,e^{-\frac{\alpha}{h}}]=\emptyset$ and
 for any $\omega_{h}$
which satisfies the inequality
 $Q_{f,f^{-1}([a,b]),h}^{(p)}(\omega_{h})\leq e^{-\frac{\alpha}{h}}\|\omega_{h}\|_{L^{2}}^{2}$\,, the spectral decomposition
$$
\omega_{h}=1_{\{0\}}(\Delta_{f,f^{-1}([a,b]),h}^{(p)})\omega_{h}
+1_{[\min(\,\sigma(\Delta_{f,f^{-1}([a,b],h)}^{(p)})\setminus\{0\}\,),+\infty[}(\Delta_{f,f^{-1}([a,b]),h}^{(p)})\omega_{h}
$$
leads to 
\begin{align*}
\textrm{dist}^{2}_{L^{2}}(\omega_{h},\ker(\Delta_{f,f^{-1}([a,b]),h}^{(p)}))&\ =\ \|1_{[\min(\,\sigma(\Delta_{f,f^{-1}([a,b],h)}^{(p)})\setminus\{0\}\,),+\infty[}(\Delta_{f,f^{-1}([a,b]),h}^{(p)})\omega_{h}\|_{L^{2}}^{2}
\\
&\ \leq\  \frac{e^{-\frac{\alpha}{h}}}{\min\big(\sigma(\Delta_{f,f^{-1}([a,b],h)}^{(p)})\setminus\{0\}\big)}\|\omega_{h}\|_{L^{2}}^{2}\,.
\end{align*}

\item We deduce  from \eqref{eq.gamma-equiv} and \eqref{eq.Gamma-equiv}  that
$$
\left(\lim_{h\to 0}\gamma^{(p)}(\alpha',[a,b],h)=0
\right)\Rightarrow \left(
\forall\alpha>\alpha'\,, \Gamma^{(p)}(\alpha,[a,b],h)\leq e^{-\frac{\alpha-\alpha'}{2h}}
\underset{h\to 0}{\longrightarrow}0\right)\,.
$$
Up to an arbitrary small change of the positive parameter $\alpha$\,, working
with $\gamma^{(p)}$ or $\Gamma^{(p)}$ is then essentially 
 equivalent.
\end{itemize}

\subsection{Exponentially small eigenvalues are zero}
\label{sec:exp0}
This section is devoted to the proof of Proposition~\ref{pr:exp0}. First of all, we can
assume $\tilde{c}_{1}=0$ if $f$ is replaced by $f-\tilde{c}_{1}$\,.
 The proof will be
done in three steps  connected  by the remarks on
$\gamma^{(p)}$ and $\Gamma^{(p)}$ from the previous subsection.\medskip

\noindent\textbf{Step 1:} Assume $[a,b]=[-\varepsilon,\varepsilon]$ with $\varepsilon>0$ (and $\tilde{c}_{1}=0$). We prove
here that
$$\forall \alpha'=2\varepsilon+c>2\varepsilon\,,\ \ \lim_{h\to 0}\gamma^{(p)}(\alpha',[-\varepsilon,\varepsilon],h)=0\,,
$$
where, owing to the
monotonicity of $\gamma^{(p)}(\alpha,[-\varepsilon,\varepsilon],h)$ 
w.r.t $\alpha$\,, we can focus on $c\in ]0,\varepsilon[$\,.\\
According to \eqref{eq.gamma-equiv}, it amounts to show there exists
$h_{c}>0$ such that
$$\left(\lambda_{h}\in\sigma(\Delta_{f,f^{-1}([a,b],h)}^{(p)})\cap
  [0,e^{-\frac{2\varepsilon+c}{h}}]\right)\Rightarrow
\left(\forall h\in ]0,h_{c}[\,,\quad \lambda_{h}=0\right)\,.$$
Take
then $\omega_{h}\in
D(\Delta_{f,f^{-1}([-\varepsilon,\varepsilon]),h}^{(p)})$
satisfying
$$
\|\omega_{h}\|_{L^{2}}=1\ \  \text{and} \ \ \Delta_{f,f^{-1}([-\varepsilon,\varepsilon]),h}^{(p)}\omega_{h}=\lambda_{h}\omega_{h}
\ \ \text{with} \ \ 0\leq\lambda_{h}\leq
e^{-\frac{2\varepsilon+c}{h}}$$
(the result is obvious for the $h$'s for which the existence of
$\omega_{h}$ fails).
 The exponential decay estimates of
Proposition~\ref{pr:Agmon} (or Hypothesis~\ref{hyp:AgmonLip} for a
general Lipschitz function) applied with
$N_{t}=f^{-1}(\left\{-\varepsilon\right\})$ and
$N_{n}=f^{-1}(\left\{\varepsilon\right\})$\,,
$K=\emptyset$\,, 
$U=f^{-1}(\left\{0\right\})$\,, $d_{Ag}(x,U)\geq |f(x)|$\,,
 and $r_{h}=0$ writes:
 \begin{equation}
   \label{eq:expdecay}
\int_{f^{-1}([-\varepsilon,\varepsilon])}e^{\frac{2|f(x)|}{h}}|\omega_{h}(x)|^{2}~dx\leq
\|e^{\frac{|f|}{h}}\omega_{h}\|^{2}_{W(f^{(-1)}([-\varepsilon,\varepsilon]))}=\tilde{O}(1)\,.
\end{equation}
Hence the mass of the probability measure with density
$|\omega_{h}|^{2}(x)$ concentrates on
$U=f^{-1}(\left\{0\right\})$ as $h\to 0$\,.
We deduce the a priori estimate
$$
\forall \delta\in ]0,\varepsilon[\,, \exists h_{\delta}>0\,, \forall
h\in]0,h_{\delta}[\,, \quad \|e^{\frac{f}{h}}1_{f_{-\delta}^{\varepsilon}}(x)\omega_{h}\|_{L^{2}}\geq \frac{e^{-\frac{\delta}{h}}}{2}\,.
$$
Once the parameter $c\in ]0,\varepsilon[$ is fixed, introduce
$s_{1}=\frac{c}{4}$ and $s_{2}\in
(\frac{c}{4},\frac{c}{2})$ and take $\chi\in
\mathcal{C}^{\infty}(M;[0,1])$ such that $\chi\equiv 0$ near
$\overline{f^{-s_{2}}}$ which contains 
a neighborhood of $\overline{f^{-\frac{c}{2}}}$ and
$\chi\equiv 1$ near
$\overline{f_{-s_{1}}}=\overline{f_{-\frac{c}{4}}}$\,.
\vspace{2cm}
\figinit{0.8cm}
\figpt1000:(-8,-1)
\figpt1001:(8,-1)
\figpt1002:(0,-1)
\figpt 0:(-7,-1)
\figpt 1:(7,-1)
\figpt 2:(-6,-1)
\figpt 3:(-3,-1)
\figpt 4:(-1.5,-1)
\figpt 5:(6,-1)
\figdrawbegin{}
\figset(width=1.3)
\figdrawline[1000,1001]
\figdrawend

\figvisu{\figBoxA}{\textbf{Figure 2:} Positions in the interval $[-\varepsilon,\varepsilon]$\,.}{
\figset write(mark=$+$)
\figwrites 1002:$0$(0.3)
\figwrites 0: $-\varepsilon$(0.3)
\figwrites 1:$\varepsilon$(0.3)
\figwrites 2:$-c$(0.3)
\figwrites 3:$-\frac{c}{2}$(0.3)
\figwrites 4:$-\frac{c}{4}$(0.3)
}
\centerline{\box\figBoxA}

\medskip
\figinit{0.9cm}
\figpt 0:(6,-1)
\figpt 1000: (-7,-1)
\figpt 1001:(7,-1)
\figpt  1:(-6,-1)
\figpt 2: (0,-1)
\figpt 3: (-4,-1)
\figpt 11:(-3.5,-1)
\figpt 111:(-3,-1)
\figpt 12:(-2.5,-0.5)
\figpt 13:(-2,0.5)
\figpt 14:(-1.75,1)
\figpt 114:(-1,1)
\figpt 115:(-0.5,1)
\figpt 15:(7,1)
\figdrawbegin{}
\figset(width=1.3)
\figdrawline[1000,1001]
\figset(width=3)
\figdrawBezier 1[1000,1,3,11]
\figdrawBezier 1[11,111,12,13]
\figdrawBezier 1[13,14,114,115]
\figdrawline[115,15]
\figdrawend{}
\figvisu{\figBoxB}{\textbf{Figure 3:} Cut-off function $\chi$ in $[-\varepsilon,\varepsilon]$\,.}{
\figset write(mark=$+$)
\figwrites 0:$0$(0.3)
\figwrites 1:$-\frac{c}{2}$(0.3)
\figwrites 2:$-\frac{c}{4}=-s_{1}$(0.3)
\figwrites 3:$-s_{2}$(0.3)
\figset write(mark={})
\figwriten 115:$\chi$(0.3)
}
\centerline{\box\figBoxB}

\medskip
Since
$$
d(\chi e^{\frac{f}{h}}\omega_{h})=
\chi d(e^{\frac{f}{h}}\omega_{h})+
d\chi\wedge
(e^{\frac{f}{h}}\omega_{h})\,,
$$
we deduce
\begin{equation}
\label{eq:somnorm}
\|d(\chi e^{\frac{f}{h}}\omega_{h})\|_{L^{2}}^{2}
\leq 2\|\chi d(e^{\frac{f}{h}}\omega_{h})\|_{L^{2}}^{2}+2\|d\chi\wedge
(e^{\frac{f}{h}}\omega_{h})\|_{L^{2}}^{2}\,.
\end{equation}
The estimate
$$
Q_{f,f^{-1}([-\varepsilon,\varepsilon]),h}(\omega_{h})=\|e^{-\frac{f}{h}}(hd)e^{\frac{f}{h}}\omega_{h}\|_{L^{2}}^{2}+
\|e^{\frac{f}{h}}(hd^{*})e^{-\frac{f}{h}}\omega_{h}\|_{L^{2}}^{2}\leq e^{-\frac{2\varepsilon+c}{h}}
$$
with $f\leq \varepsilon$ then implies that the first term in the r.h.s. 
of \eqref{eq:somnorm} is of order $\tilde{O}(e^{-\frac{c}{h}})$\,.
Meanwhile $\supp(d\chi)\subset \overline{f^{-\frac{c}{4}}}$
and the exponential decay estimate \eqref{eq:expdecay}
 imply that the
second term in the r.h.s. of
\eqref{eq:somnorm} is of order $\tilde{O}(e^{-2\frac{2c}{4h}})$\,. 
Adding the boundary conditions
$\mathbf{n}_{f=\varepsilon}\omega_{h}=0$ and $\mathbf{n}_{f=\varepsilon}d_{f,h}\omega_{h}=0$\,, i.e.
$\mathbf{n}_{f=\varepsilon}(e^{\frac{f}{h}}\omega_{h})=0$ and
$\mathbf{n}_{f=\varepsilon}d(e^{\frac{f}{h}}\omega_{h})=0$\,, 
we have thus proved that
$$
\left\{
  \begin{array}[c]{l}
    \chi e^{\frac{f}{h}}\omega_{h}\in D(\Delta_{0,f^{-1}([-
s_{2}, \varepsilon],1)})\,,\\
    \|d_{0,f^{-1}([-s_{2},\varepsilon],1)}(\chi
    e^{\frac{f}{h}}\omega)\|_{L^{2}}^{2}=\tilde{O}(e^{-\frac{c}{h}})\,,\\
\lim_{h\to 0} h\log \|\chi e^{\frac{f}{h}}\omega_{h}\|_{L^{2}}=0\,.
  \end{array}
\right.
$$ 
Set $u_{h}=\frac{\chi e^{\frac{f}{h}}\omega_{h}}{\|\chi
  e^{\frac{f}{h}}\omega_{h}\|_{L^{2}}}$ so that $\|u_{h}\|_{L^{2}}=1$\,,
$u_{h}\in D(\Delta_{0,f^{-1}([-s_{2}, \varepsilon]),1})$
and $\|du_{h}\|_{L^{2}}^{2}=\tilde{O}(e^{-\frac{c}{h}})$\,.\\
By using the Hodge decomposition (see Proposition~\ref{pr:domain}) and $\sigma
(\Delta_{0,f^{-1}([-s_{2},\varepsilon]),1})\setminus\{0\}\subset
[\mu_{1},+\infty)\subset \rz^{+*}$\,, with 
$\mu_{1}$ fixed by $\varepsilon>0$ and $s_{2}>0$\,,
we obtain the decomposition of $u_{h}$:
$$
u_{h}
\ =\ \Pi_{\ker d_{0,f^{-1}([-s_{2},\varepsilon]),1}}u_{h}+  d^{*}_{0,f^{-1}([-s_{2},\varepsilon]),1}u_{2,h}\,,
$$
where $ d^{*}_{0,f^{-1}([-s_{2},\varepsilon]),1}u_{2,h}
$ in $\big(\ker \Delta_{0,f^{-1}([-s_{2},\varepsilon]),1})^{\perp}=\Ran~1_{\{[\mu_{1},+\infty)\}}(\Delta_{0,f^{-1}([-s_{2}, \varepsilon]),1}^{(p)})$\,.
Writing shortly $\mathbf{d}=d_{0,f^{-1}([-s_{2},\varepsilon],1)}$
and $\mathbf{d}^{*}=d_{0,f^{-1}([-s_{2},\varepsilon],1)}^{*}$\,, it follows that
$$
\tilde{O}(e^{-\frac{c}{h}})=\|\mathbf{d} u_{h}\|_{L^{2}}^{2}=
\|\mathbf{d} \mathbf{d}^{*} u_{2,h}\|_{L^{2}}^{2}
= 
Q_{0,f^{-1}([-s_{2}, \varepsilon]),1}^{(p)}(  \mathbf{d}^{*} u_{2,h})
\geq \mu_{1}\|\mathbf{d}^{*} u_{2,h}\|_{L^{2}}^{2}
\,.
$$
We deduce $\textrm{dist}_{L^{2}}(u_{h},\ker
d_{0,f^{-1}([-s_{2},\varepsilon],1)})=
\tilde{O}(e^{-\frac{c}{2h}})$ and then
the existence of a form
$\eta_{h}\in \ker(d_{0,f^{-1}([-s_{2},\varepsilon],1)})$
such that 
$$
\|\chi e^{\frac{f}{h}}\omega_{h}-\eta_{h}\|_{L^{2}(f^{-1}([-s_{2},\varepsilon])}=\tilde{O}(e^{-\frac{c}{2h}})\,.
$$
By the first item of Remark~\ref{re:rem1}, the extension
$\tilde{\eta}_{h}$ of $\eta_{h}$  by $0$ in $f^{-s_{2}}_{-\varepsilon}$
belongs to $\ker(d_{0,f^{-1}([-\varepsilon,\varepsilon]),1})$ with
$\supp\tilde{\eta}_{h}\subset
\overline{f_{-s_{2}}^{\varepsilon}}$ and 
$\|\chi
e^{\frac{f}{h}}\omega_{h}-\tilde{\eta}_{h}\|_{L^{2}}=\tilde{O}(e^{-\frac{c}{2h}})$\,.\\
After multiplying by
$e^{-\frac{f}{h}}=\mathcal{O}(e^{\frac{s_{2}}{h}})$
in $f_{-s_{2}}^{\varepsilon}$\,, we obtain
$$
\left\{
\begin{array}[c]{l}
     \|\chi\omega_{h}-e^{-\frac{f}{h}}\tilde{\eta}\|_{L^{2}}=
\tilde{O}(e^{-\frac{c}{2h}+\frac{s_{2}}{h}})\quad,\quad
\frac{c}{2}>s_{2}\,,\\
 e^{-\frac{f}{h}}
\tilde{\eta}_{h}\in
\ker  (d_{f,f^{-1}([-\varepsilon,\varepsilon]),h})\,.
\end{array}
\right.
$$
We conclude with 
$\|\chi
\omega_{h}-\omega_{h}\|_{L^{2}}=\tilde{O}(e^{-\frac{c}{4h}})$
(since $\supp(1-\chi)\subset\overline{f^{-s_{1}}}
=\overline{f^{-\frac{c}{4}}}$) that 
$$
\textrm{dist}(\omega_{h},\ker(d_{f,f^{-1}([-\varepsilon,\varepsilon]),h}))=O(
e^{-\frac{c'}{h}})\quad \text{for some $c'>0$}.
$$
The duality consists in replacing $f$ by $-f$ (which does not change
$[-\varepsilon,\varepsilon]$), the differential form $\omega_{h}\in
W(f^{-1}(]-\varepsilon,\varepsilon[);\Lambda^{p} T^{*}M)$ by
$\star \omega_{h}\in
W(f^{-1}(]-\varepsilon,\varepsilon[);\Lambda^{d-p} T^{*}M\otimes\mathrm{or}_{M})$ where the orientation twist
does not change the analysis, $\mathbf t$ by $\mathbf n$ (and
conversely), $\star$ and $\star^{-1}$\,,  and $d_{f,h}$ by $d_{-f,h}^{*}$
(and conversely).
This leads to 
$$
\textrm{dist}(\star \omega_{h}\,,\,
\ker(d_{-f,f^{-1}([-\varepsilon,\varepsilon]),h}))
\ =\ \textrm{dist}(\omega_{h}\,,\,
\ker(d_{f,f^{-1}([-\varepsilon,\varepsilon]),h}^{*}))
\ =\ O(e^{-\frac{c'}{h}})\,.
$$
Assume by contradiction that $\lambda_{h}\neq 0$\,.
Since $\omega_{h}=\lambda_{h}^{-1} \Delta_{f,f^{-1}([-\varepsilon,\varepsilon]),h}^{(p)}
\omega_{h}\in\big(\ker \Delta_{f,f^{-1}([-\varepsilon,\varepsilon]),h}^{(p)}\big)^{\perp}$\,,
 the Hodge decomposition (see Proposition~\ref{pr:domain}) leads to
$$
\omega_{h}
\ =\ \Pi_{\ker d_{f,f^{-1}([-\varepsilon,\varepsilon]),h}}\omega_{h}
+ \Pi_{\ker d^{*}_{f,f^{-1}([-\varepsilon,\varepsilon]),h}}\omega_{h}\,.
$$
The squared norm $1=\|\omega_{h}\|^{2}$ thus equals
\begin{eqnarray*} 
\textrm{dist}^{2}_{L^{2}}(\omega_{h},\ker(d_{f,f^{-1}([-\varepsilon,\varepsilon]),h}))
+
\textrm{dist}^{2}_{L^{2}}(\omega_{h},\ker(d_{f,f^{-1}([-\varepsilon,\varepsilon]),h}^{*}))
=
\tilde{O}(e^{-\frac{c'}{h}})\,,
\end{eqnarray*}
which is impossible for $0<h<h_{\varepsilon}$\,,
$h_{\varepsilon}>0$ small enough.
It follows that $\sigma(\Delta_{f,f^{-1}([a,b],h)}^{(p)})\cap [0,e^{-\frac{2\varepsilon+c}{h}}]\subset \{0\}$
for $h$ small enough,
which implies $\lim_{h\to 0}\gamma^{(p)}(\alpha',[-\varepsilon,\varepsilon],h)=0$ according to the comments
following Definition~\ref{de:gam}.
\medskip

\noindent\textbf{Step~2:} From Step~1, we know $\lim_{h\to
  0}\gamma^{(p)}(\alpha',[-\varepsilon,\varepsilon],h)=0$ for any
$\alpha'>2\varepsilon$ and the comparison of the quantities
$\gamma^{(p)}$ and $\Gamma^{(p)}$ in the previous subsection  leads to 
$$
\forall \alpha> 2\varepsilon\,,\quad \lim_{h\to 0}\Gamma^{(p)}(\alpha,[-\varepsilon,\varepsilon],h)=0\,.
$$
Working with $\Gamma^{(p)}$ brings the flexibility to use some
restriction argument from $f_{a}^{b}$ to
$f_{-\varepsilon}^{\varepsilon}$\,, which of course does not send
eigenvectors onto eigenvectors.\medskip

\noindent\textbf{Step~3:}
For the general case $a<0=\tilde{c}_{1}<b$\,, we now prove
$$\forall \alpha>0\,,\ \ \sigma(\Delta_{f,f^{-1}([a,b],h)}^{(p)})\cap [0,e^{-\frac{\alpha}{h}}]\subset \{0\}\,,$$
where, by monotonicity w.r.t $\alpha$\,, it is sufficient to consider $\alpha\leq
\min(-a, b)$\,.
Let us then assume 
that $\omega_{h}$ satisfies
$\Delta_{f,f^{-1}([a,b]),h}^{(p)}\omega_{h}=\lambda_{h}
\omega_{h}$ with $\|\omega_{h}\|_{L^{2}}=1$ and $0\leq\lambda_{h}\leq
e^{-\frac{\alpha}{h}}$\,. Take $\varepsilon\in
]0,\frac{\alpha}{4}[$ and consider
$f_{-\varepsilon}^{\varepsilon}\subset f_{a}^{b}$\,. 
We know that
$$
\|d_{f,h}\omega_{h}\|_{L^{2}(f_{-\varepsilon}^{\varepsilon})}^{2}
+\|d_{f,h}^{*}\omega_{h}\|_{L^{2}(f_{-\varepsilon}^{\varepsilon})}^{2}
\leq
\|d_{f,f^{-1}([a,b]),h}\omega_{h}\|_{L^{2}(f_{a}^{b})}^{2}+
\|d_{f,f^{-1}([a,b]),h}^{*}\omega_{h}\|_{L^{2}(f_{a}^{b})}^{2}
\leq e^{-\frac{\alpha}{h}}\,,
$$
although $\omega_{h}\big|_{f_{-\varepsilon}^{\varepsilon}}$ a priori does not
belong neither to $D(\Delta_{f,f^{-1}([-\varepsilon,\varepsilon]),h}^{(p)})$
nor to $D(Q_{f,f^{-1}([-\varepsilon,\varepsilon],h)}^{(p)})$\,.\\
We now use Lemma~\ref{le:exten1} in the two subsets
$f^{-1}([-\varepsilon,-\delta])$ and
$f^{-1}([\delta,\varepsilon])$ for some $\delta\in
]0,\frac{\varepsilon}{4}[$ which will be fixed later.\\
Consider $\overline{\Omega}=f^{-1}([-\varepsilon,-\delta]) $
(the other case is symmetric) and take the cut-off $\chi_{-},\tilde{\chi}_{-}\in
\mathcal{C}^{\infty}(f^{-1}[-\varepsilon,-\delta];[0,1])$
 with
$\supp \chi_{-}\subset f^{-1}(-]\varepsilon+\delta,-\delta])$\,,
$\chi_{-}\equiv 1$ in $f^{-1}([-\varepsilon+2\delta,-\delta])$\,, and
$\supp \tilde{\chi}_{-}\subset f^{-1}([-\varepsilon,-\delta[)$\,,
$\tilde{\chi}_{-}\equiv 1$ in $f^{-1}([-\varepsilon,-2\delta])$\,.
\vspace{1cm}
\figinit{0.8cm}
\figpt1000:(-8,-1)
\figpt1001:(9,-1)
\figpt1002:(0,-1)
\figpt 0:(-6.5,-1)
\figpt 10:(-6,-1)
\figpt 1:(8.5,-1)
\figpt 2:(-1.5,-1)
\figpt 3:(-1,-1)
\figpt 4:(1,-1)
\figdrawbegin{}
\figset(width=1.3)
\figdrawline[1000,1001]
\figdrawend

\figvisu{\figBoxA}{\textbf{Figure 4:} Positions in the interval $[a,b]$\,.}{
\figset write(mark=$+$)
\figwrites 1002:$0$(0.3)
\figwrites 0: $a$(0.3)
\figwriten 10:$-\alpha$(0.3)
\figwrites 1:$b$(0.3)
\figwriten 2:$-\frac{\alpha}{4}$(0.3)
\figwrites 3:$-\varepsilon$(0.3)
\figwrites 4:$\varepsilon$(0.3)
}
\centerline{\box\figBoxA}

\medskip
\figinit{cm}
\figpt 0:(6,-1)
\figpt 1000: (-7,-1)
\figpt 1001:(7,-1)
\figpt  1:(-6,-1)
\figpt 3: (-5,-1)
\figpt 11:(-4.8,-1)
\figpt 111:(-4.6,-1)
\figpt 12:(-4.5,-0.5)
\figpt 13:(-4.4,0.5)
\figpt 14:(-4.35,1)
\figpt 2: (-4,-1)
\figpt 114:(-4.1,1)
\figpt 115:(-4,1)
\figpt 15:(5,1)
\figpt 5:(5,-1)

\figpt 203:(5,-1)
\figpt 211:(4.8,-1)
\figpt 311:(4.6,-1)
\figpt 212:(4.5,-0.5)
\figpt 213:(4.4,0.5)
\figpt 214:(4.36,0.8)
\figpt 314:(4.1,0.8)
\figpt 315:(4,0.8)
\figpt 215:(-6,0.8)
\figpt 205:(4,-1)

\figdrawbegin{}
\figset(width=1.3)
\figdrawline[1000,1001]
\figset(width=3)
\figdrawBezier 1[1,1,3,11]
\figdrawBezier 1[11,111,12,13]
\figdrawBezier 1[13,14,114,115]
\figdrawline[115,15]
\figset(dash=3)
\figdrawline[203,211]
\figdrawBezier 1[211,311,212,213]
\figdrawBezier 1[213,214,314,315]
\figdrawline[315,215]
\figdrawend{}
\figvisu{\figBoxB}{\textbf{Figure 5:} Cut-off functions $\chi_{-}$ and
$\tilde{\chi}_{-}$ in $[-\varepsilon,0]\subset[a,b]$\,.}{
\figset write(mark=$+$)
\figwrites 0:$0$(0.3)
\figwriten 1:$-\varepsilon$(0.3)
\figwrites 3:$-\varepsilon+\delta$(0.3)
\figwrites 2:$-\varepsilon+2\delta$(0.7)
\figwrites 5:$-\delta$(0.3)
\figwrites 205:$-2\delta$(0.3)
\figset write(mark={})
\figwriten 115:$\chi_{-}$(0.3)
\figwritesw 214:$\tilde{\chi}_{-}$(0.3)
}
\centerline{\box\figBoxB}

\vspace{1cm}
The form $\eta_{1,-}$ and $\eta_{2,-}$ in  $D(\Delta_{f,f^{-1}([-\varepsilon,-\delta]),h})$ are defined by
\begin{eqnarray*}
\eta_{1,-}&=&d_{f,[-\varepsilon,-\delta],h}^{*}(\Delta_{f,f^{-1}([-\varepsilon,-\delta]),h})^{-1}((hd\chi_{-})\wedge
\omega_{h})
\\
\eta_{2,-}&=&d_{f,[-\varepsilon,-\delta],h}(\Delta_{f,f^{-1}([-\varepsilon,-\delta]),h})^{-1}(h\mathbf{i}_{\nabla\chi_{-}}
\omega_{h})\,.
\end{eqnarray*}
Lemma~\ref{le:exten1} combined with $d_{Ag}(x,y)\geq |f(x)-f(y)|$
implies
\begin{multline*}
\|d_{f,h}(\chi_{-}\omega_{h}-\tilde{\chi}_{-}(\eta_{1,-}+\eta_{2,-}
)\|_{L^{2}(f_{-\varepsilon}^{-\delta})}
+
\|d_{f,h}^{*}(\chi_{-}\omega_{h}-\tilde{\chi}_{-}(\eta_{1,-}+\eta_{2,-}
)\|_{L^{2}(f_{-\varepsilon}^{-\delta})}
\\
\leq
\tilde{O}(e^{-\frac{\varepsilon-4\delta}{h}})\|\omega_{h}\|_{L^{2}(f_{-\varepsilon+\delta}^{-\varepsilon+2\delta})}
+C_{\chi_{-}}\left[\|d_{f,h}\omega_{h}\|_{L^{2}(f_{-\varepsilon}^{-\delta})}+\|d_{f,h}^{*}\omega_{h}\|_{L^{2}(f_{-\varepsilon}^{-\delta})}\right]\,.
\end{multline*}
Because $\Delta_{f,f^{-1}([a,b],h)}\omega_{h}=\lambda_{h}\omega_{h}$
with $\|\omega_{h}\|_{L^{2}}=1$\,, the Agmon estimate of
Proposition~\ref{pr:Agmon} (or Hypothesis~\ref{hyp:AgmonLip} for a
general Lipschitz function), applied with
$N_{t}=f^{-1}(\left\{a\right\})$ and
$N_{n}=f^{-1}(\left\{b\right\})$\,,
$K=\emptyset$\,, 
$U=f^{-1}(\left\{0\right\})$\,, $d_{Ag}(x,U)\geq |f(x)|$\,,
 and $r_{h}=0$\, implies
 \begin{equation}
\label{eq:expdecay2}
\|\omega_{h}\|_{L^{2}(f_{-\varepsilon+\delta}^{-\varepsilon+2\delta})}=\tilde{O}(e^{-\frac{\varepsilon-2\delta}{h}})\,,
\end{equation}
while we know
$$
  \|d_{f,h}\omega_{h}\|_{L^{2}(f_{-\varepsilon}^{-\delta})}^{2}+\|d_{f,h}^{*}\omega_{h}\|_{L^{2}(f_{-\varepsilon}^{-\delta})}^{2}\leq
  \|d_{f,h}\omega_{h}\|_{L^{2}(f_{a}^{b})}^{2}+\|d_{f,h}^{*}\omega_{h}\|_{L^{2}(f_{a}^{b})}^{2}
\leq e^{-\frac{\alpha}{h}}\,.
$$
With $\frac\alpha4>\varepsilon>4\delta$\,, we have thus
\begin{equation}
\label{eq.app-le1}
\|d_{f,h}(\chi_{-}\omega_{h}-\tilde{\chi}_{-}(\eta_{1,-}+\eta_{2,-}
)\|_{L^{2}(f_{-\varepsilon}^{-\delta})}
+
\|d_{f,h}^{*}(\chi_{-}\omega_{h}-\tilde{\chi}_{-}(\eta_{1,-}+\eta_{2,-}
)\|_{L^{2}(f_{-\varepsilon}^{-\delta})}
=
\tilde{O}(e^{-\frac{2\varepsilon-6\delta}{h}})\,.
\end{equation}
A symmetric construction provides 
two cut-off functions $\chi_{+},\tilde{\chi}_{+}\in
\mathcal{C}^{\infty}(f^{-1}([\delta,\varepsilon]))$ such that 
$$\supp
\chi_{+}\subset
 f^{-1}([\delta,\varepsilon-\delta[)\ \ \ \text{and}\ \ \  \chi_{+}\equiv 1\ 
\ \text{in}\ \  f^{-1}([\delta, \varepsilon-2\delta])\,,$$
$$\supp
\tilde{\chi}_{+}\subset f^{-1}(]\delta,\varepsilon])
\ \ \ \text{and}\ \ \ 
\tilde{\chi}_{+}\equiv 1\ \ \text{ in} \ \ f^{-1}([2\delta,\varepsilon])\,,$$
and then two forms
$\eta_{1,+}, \eta_{2,+}\in D(\Delta_{f,f^{-1}([\delta,\varepsilon]),h})$
such that
\begin{equation}
\label{eq.app-le2}
\|d_{f,h}(\chi_{+}\omega_{h}-\tilde{\chi}_{+}(\eta_{1,+}+\eta_{2,+}
)\|_{L^{2}(f_{\delta}^{\varepsilon})}
+
\|d_{f,h}^{*}(\chi_{+}\omega_{h}-\tilde{\chi}_{+}(\eta_{1,+}+\eta_{2,+}
)\|_{L^{2}(f_{\delta}^{\varepsilon})}
=
\tilde{O}(e^{-\frac{2\varepsilon-6\delta}{h}})\,.
\end{equation}
Take  now $\chi\in
\mathcal{C}^{\infty}_{0}(f^{-1}(]-\varepsilon+\delta,\varepsilon-\delta[;[0,1])$
which equals $1$ in
$f^{-1}([-\varepsilon+2\delta,\varepsilon-2\delta])$ and coincides
 with
$\chi_{-}$ (resp. $\chi_{+}$) in
$f_{-\varepsilon+\delta}^{-\varepsilon+2\delta}$
(resp. in $f_{\varepsilon-2\delta}^{\varepsilon-\delta}$) and set
$$
v_{h}=\chi \omega_{h}-\tilde{\chi}_{-}(\eta_{1,-}+\eta_{2,-}) -\tilde{\chi}_{+}(\eta_{1,+}+\eta_{2,+})\,.
$$
\medskip
\figinit{cm}
\figpt 0:(6,-1)
\figpt 1000: (-7,-1)
\figpt 1001:(7,-1)
\figpt  1:(-6,-1)
\figpt 3: (-5,-1)
\figpt 11:(-4.8,-1)
\figpt 111:(-4.6,-1)
\figpt 12:(-4.5,-0.5)
\figpt 13:(-4.4,0.5)
\figpt 14:(-4.35,1)
\figpt 2: (-4,-1)
\figpt 114:(-4.1,1)
\figpt 115:(-4,1)
\figpt 15:(-3,1)
\figpt 5:(5,-1)

\figpt 203:(5,-1)
\figpt 211:(4.8,-1)
\figpt 311:(4.6,-1)
\figpt 212:(4.5,-0.5)
\figpt 213:(4.4,0.5)
\figpt 214:(4.35,1)
\figpt 314:(4.1,1)
\figpt 315:(4,1)
\figpt 215:(3,1)
\figpt 205:(4,-1)
\figpt 206:(0,1)

\figpt 6:(0,-1)
\figpt 7:(1,-1)
\figpt 8:(-1,-1)

\figdrawbegin{}
\figset(width=1.3)
\figdrawline[1000,1001]
\figset(width=3)
\figdrawBezier 1[1,1,3,11]
\figdrawBezier 1[11,111,12,13]
\figdrawBezier 1[13,14,114,115]
\figdrawline[115,15]
\figdrawline[203,211]
\figdrawBezier 1[211,311,212,213]
\figdrawBezier 1[213,214,314,315]
\figdrawline[315,215]
\figdrawline[5,0]
\figset(dash=3)
\figdrawline[15,215]
\figdrawend{}
\figvisu{\figBoxB}{\textbf{Figure 6:}  Cut-off function $\chi$ in $[-\varepsilon,\varepsilon]$\,.}{
\figset write(mark=$+$)
\figwriten 0:$\varepsilon$(0.3)
\figwriten 1:$-\varepsilon$(0.3)
\figwrites 3:$-\varepsilon+\delta$(0.3)
\figwrites 2:$-\varepsilon+2\delta$(0.7)
\figwrites 5:$\varepsilon-\delta$(0.3)
\figwrites 205:$\varepsilon-2\delta$(0.7)
\figwrites 6:$0$(0.3)
\figwrites 7:$\delta$(0.3)
\figwrites 8:$-\delta$(0.3)
\figset write(mark={})
\figwriten 115:$\chi_{-}$(0.3)
\figwriten 214:$\chi_{+}$(0.3)
\figwriten 206:$\chi$(0.3)
}
\centerline{\box\figBoxB}

\vspace{1cm}
This form is close to
$\omega_{h}\big|_{f_{-\varepsilon}^{\varepsilon}}$\,. In fact, write
$$
v_{h}- \omega_{h}\big|_{f_{-\varepsilon}^{\varepsilon}}=(\chi-1)\omega_{h}\big|_{f_{-\varepsilon}^{\varepsilon}}-\tilde{\chi}_{-}(\eta_{1,-}+\eta_{2,-})-\tilde{\chi}_{+}(\eta_{1,+}+\eta_{2,+})\,,
$$
where, according to Lemma~\ref{le:exten1} and
to the exponential decay estimate \eqref{eq:expdecay2}  (and its
symmetric version on $[\varepsilon-2\delta,\varepsilon-\delta]$)\,,
$$
\|\tilde{\chi}_{\pm}\eta_{i,\pm}\|_{L^{2}}=\mathcal O(\|\omega_{h}\|_{\supp{d\chi_{\pm}}})
=\tilde{O}(e^{-\frac{\varepsilon-2\delta}{h}})
\quad\text{for $i\in\{1,2\}$}
$$
and
$$
\|(\chi-1)\omega_{h}\|_{L^{2}(f_{-\varepsilon}^{\varepsilon})}
=\tilde{O}(e^{-\frac{\varepsilon-2\delta}{h}})\,,
$$
which implies
$$
\|v_{h}-\omega_{h}\|_{L^{2}(f_{-\varepsilon}^{\varepsilon})}=\tilde{O}(e^{-\frac{\varepsilon-2\delta}{h}})\,.
$$
The form $v_{h}$ also
 satisfies, for $\mathbf{d}=d_{f,h}$ or $\mathbf{d}=d_{f,h}^{*}$\,,
$$
\mathbf{d}v_{h}=
[\mathbf{d}(\chi_{-}\omega_{h}-\tilde{\chi}_{-}(\eta_{1,-}+\eta_{2,-}))]\big|_{f_{-\varepsilon}^{-\delta}}
+[\mathbf{d}\omega_{h}]\big|_{f_{-\delta}^{\delta}}
+[\mathbf{d}(\chi_{+}\omega_{h}-\tilde{\chi}_{+}(\eta_{1,+}+\eta_{2,+})]\big|_{f_{\delta}^{\varepsilon}}\,.
$$
Then, since $v_{h}$ belongs to 
$D(\Delta_{f,f^{-1}([-\varepsilon,\varepsilon],h)})$ 
by construction, it
satisfies, by \eqref{eq.app-le1} and \eqref{eq.app-le2},
$$
\|d_{f,f^{-1}([-\varepsilon,\varepsilon])h}v_{h}\|^{2}+\|d^{*}_{f,f^{-1}([-\varepsilon,\varepsilon]),h}v_{h}\|^{2}
\ =\  \tilde{O}(e^{-\frac{4\varepsilon-12\delta}{h}})\,.
$$
We finally take $\delta=\frac{\varepsilon}{12}$ for which the
r.h.s. of the above relation is $\tilde{O}(e^{-\frac{3\varepsilon}{h}})$\,, with
$3\varepsilon>2\varepsilon$\,. By Step~2, this implies 
$$
\lim_{h\to 0}
\frac{\textrm{dist}_{L^{2}}(v_{h},\ker(\Delta_{f,f^{-1}([-\varepsilon,\varepsilon]),h}))}{\|v_{h}\|_{L^{2}}}=0\,.
$$
But the Agmon estimates
of
Proposition~\ref{pr:Agmon} or Hypothesis~\ref{hyp:AgmonLip} also imply
 $$
 \|\omega_{h}\|_{L^{2}(f_{-\varepsilon}^{\varepsilon})}
 = 1 + \tilde{O}(e^{-\frac{\varepsilon}{h}})
\quad  \text{and then}\quad
\|v_{h}\|_{L^{2}(f_{-\varepsilon}^{\varepsilon})}=1+\tilde{O}(e^{-\frac{\varepsilon-2\delta}{h}})\,.
 $$
Denoting by $F\subset L^{2}(f_{a}^{b})$  the subspace
$\ker(\Delta_{f,f^{-1}([-\varepsilon,\varepsilon]),h}^{(p)})$ extended by $0$ in $f_{a}^{-\varepsilon}\sqcup
f_{\varepsilon}^{b}$\,, it then follows from the preceding analysis that
$$
\lim_{h\to 0}\textrm{dist}_{L^{2}}(\omega_{h}\,,\, F)=0\,.
$$
Since $\dim F$ is finite and does not depend on $h>0$ 
(see the second item  in  Remark~\ref{re:rem1}), 
there exists $h_{\alpha}>0$ such that 
for every 
$h\in
]0,h_{\alpha}[$\,,
$$
\dim F_{[0,e^{-\frac{\alpha}{h}}],[a,b],h}^{(p)} \leq 
\dim F=\dim
\ker(\Delta_{f,f^{-1}([-\varepsilon,\varepsilon],h)}^{(p)})=\dim \ker
(\Delta_{f,f^{-1}([a,b]),h}^{(p)})\,,
$$
 where the last equality follows from
$[-\varepsilon,\varepsilon]\subset[a,b]$
and $[a,b]\cap \left\{c_{1},\ldots,c_{N_{f}}\right\}= \{\tilde{c}_{1}=0\} $
(see~\eqref{eq.const-kernel}). This implies
that
$\sigma(\Delta_{f,f^{-1}([a,b],h)}^{(p)})\cap [0,e^{-\frac{\alpha}{h}}]\subset \{0\}$
for $h\in ]0,h_{\alpha}[$
and
this ends the proof.

\subsection{Consequences}
\label{sec:resest1}
We still work under Hypothesis~\ref{hyp:1vc}: $f$ admits a unique
``critical value''
 $\tilde{c}_{1}\in [a,b]$\,, $a<\tilde{c}_{1}<b$\,. With the information of
Proposition~\ref{pr:exp0}, the resolvent estimates of 
Subsection~\ref{sec:consAgm}
 lead easily to similar
estimates for spectrally defined operators. Finally we deduce other
properties which will be used in the induction process in terms of the
number $N$ of ``critical values''.\\

\subsubsection{Estimates for spectral operators}

For a Borel set $I\subset \rz$ we introduce the notation:
\begin{equation}
  \label{eq:defPiE}
  \Pi_{I,[a,b],h}=1_{I}(\Delta_{f,f^{-1}([a,b]),h})\,.
\end{equation}

\begin{prop}
\label{pr:proj1}
  Under Hypothesis~\ref{hyp:1vc} the spectral projection on the kernel
  $\Pi_{\left\{0\right\},[a,b],h}$ satisfies 
\begin{equation*}
\Pi_{\left\{0\right\},[a,b],h}(x,y)=\tilde{O}(e^{-\frac{|f(x)-f(y)|}{h}})
\end{equation*}
according to Definition~\ref{de:Otkernel}.
\end{prop}
\begin{proof}
  It suffices to use the formula
$$
\Pi_{\left\{0\right\},[a,b],h}=\frac{1}{2i\pi}\int_{\gamma_{h}}(z-\Delta_{f,f^{-1}[a,b],h})^{-1}~dz
$$
for the suitable contour $\gamma_{h}$ such that 
$1=\tilde{O}(\textrm{dist}(\gamma_{h},\sigma
(\Delta_{f,f^{-1}([a,b]),h})))$\,, and then to apply
Proposition~\ref{pr:1well} with $\varepsilon_{0}>0$ arbitrarily small.
 Such a contour is chosen as follows.
For $n\in\nz$\,, Proposition~\ref{pr:exp0} says  
$$
\exists h_{n}>0\,, \forall h\in ]0,h_{n}[\,,
\quad
\sigma(\Delta_{f,f^{-1}([a,b]),h})\cap [0,e^{-\frac{1}{2(n+1)h}}]=\left\{0\right\}\,,
$$
and the condition $h_{n+1}<h_{n}$ can be added\,.
Take simply $\gamma_{h}=\left\{z\in\cz\,,
  |z|=e^{-\frac{1}{(n+1)h}}\right\}$ for
 $h\in [h_{n+1},h_{n}[$\,.
\end{proof}
The final result of this paragraph
 extends the exponential decay estimates of
Proposition~\ref{pr:Agmon} (or Hypothesis~\ref{hyp:AgmonLip}), when $f$ admits a single singular value
$\tilde{c}_{1}$\,, under orthogonality conditions. It will be referred to as
the ``orthogonality lemma''.\\
Because $\Delta_{f,f^{-1}([a,b]),h}$ has a discrete spectrum, the operator
\begin{equation*}
\Delta_{f,f^{-1}([a,b]),h}\big|_{\ker(\Delta_{f,f^{-1}([a,b]),h})^{\perp}}:\ker(\Delta_{f,f^{-1}([a,b]),h})^{\perp}\to \ker(\Delta_{f,f^{-1}([a,b]),h})^{\perp}\,,
\end{equation*}
is invertible. 
We now define $(\Delta_{f,f^{-1}([a,b]),h}^{\perp})^{-1}$ by extension
by $0$ on $\ker(\Delta_{f,f^{-1}([a,b]),h})$:
\begin{equation}
  \label{eq:resorth}
  (\Delta_{f,f^{-1}([a,b]),h}^{\perp})^{-1}=\underbrace{0}_{\ker(\Delta_{f,f^{-1}([a,b]),h})}\mathop{\oplus}^{\perp}
\underbrace{(\Delta_{f,f^{-1}([a,b]),h}\big|_{\ker(\Delta_{f,f^{-1}([a,b]),h})^{\perp}})^{-1}}_{\ker(\Delta_{f,f^{-1}([a,b]),h})^{\perp}}
\end{equation}
 Thus, the equality 
$\omega_{h}=(\Delta_{f,f^{-1}([a,b],h)}^{\perp})^{-1}r_{h}$ 
simply means that  $\omega_{h}$ is the
unique solution in $\ker(\Delta_{f,f^{-1}([a,b]),h})^{\perp}\cap
D(\Delta_{f,f^{-1}([a,b]),h})$ to 
$$
\Delta_{f,f^{-1}([a,b],h)}\omega_{h}=(1-\Pi_{\left\{0\right\},[a,b],h})r_{h}\,.
$$

\begin{lem}
  \label{le:orthlem}
Under Hypothesis~\ref{hyp:1vc},  the operator defined by
\eqref{eq:resorth}
satisfies
$$
(\Delta_{f,f^{-1}([a,b]),h}^{\perp})^{-1}(x,y)=\tilde{O}(e^{-\frac{|f(x)-f(y)|}{h}})
$$
in the sense of Definition~\ref{de:Otkernel}.
\end{lem}
\begin{proof}
With
$A=\Delta_{f,f^{-1}([a,b]),h}$  and
$\Pi_{\left\{0\right\},[a,b],h}=1_{\left\{0\right\}}(A)$ write simply
$\omega_{h}=(\Delta_{f,f^{-1}([a,b]),h}^{\perp})^{-1}r_{h}$ as
\begin{align*}
\omega_{h}=(1-\Pi_{\left\{0\right\},[a,b],h})\omega_{h}&=-\frac{1}{2i\pi}\int_{\gamma_{h}}\frac{A}{z(z-A)}\omega_{h}~dz
\\
&=-\frac{1}{2i\pi}\int_{\gamma_{h}}\frac{1}{z(z-A)}(1-\Pi_{\left\{0\right\},[a,b],h})r_{h}~dz\,,
\end{align*}
where $\gamma_{h}$ is the contour introduced in the proof of
Proposition~\ref{pr:proj1}. To conclude, it then suffices to combine the resolvent
estimates of Proposition~\ref{pr:1well} with $\varepsilon_{0}>0$ arbitrarily small,
as used in the proof of Proposition~\ref{pr:proj1}, 
 and the result of Proposition~\ref{pr:proj1}.
\end{proof}
\subsubsection{Changing the interval $[a,b]$}
\label{sec:changinb}
For further applications, it is useful to specify the effect of
changing $b$ in $f_{a}^{b}$\,. Rough estimates after a
change of $a$ and $b$ are followed by more accurate estimates
 after a
change of $b$ only.\\
Remember that we work under Hypothesis~\ref{hyp:1vc} which contains
Hypothesis~\ref{hyp:mainf} or for a more general Lipschitz function
Hypothesis~\ref{hyp:Lipbar} and 
Hypothesis~\ref{hyp:AgmonLip}. 
\begin{prop}
\label{pr:aba1b1}
Assume Hypothesis~\ref{hyp:1vc} and $a<a'<\tilde{c}_{1}<b'<b$\,. The kernels
$F_{\left\{0\right\},[\alpha,\beta],h}=\ker(\Delta_{f,f^{-1}([\alpha,\beta]),h})=\Ran~\Pi_{\left\{0\right\},[\alpha,\beta],h}$\,,
$\alpha\in \left\{a,a'\right\}$\,, $\beta\in \left\{b,b'\right\}$
satisfy
$$
\vec{d}(F_{\left\{0\right\},[a',b'],h},F_{\left\{0\right\},[a,b],h})=\vec{d}(F_{\left\{0\right\},[a,b],h},F_{\left\{0\right\},[a',b'],h})=\tilde{O}(e^{-\frac{\min\left\{b'-\tilde{c}_{1},\tilde{c}_{1}-a'\right\}}{h}})\,,
$$
where the second inclusion of $F_{\left\{0\right\},[a',b'],h}\subset L^{2}(f_{a'}^{b'})\subset L^{2}(f_{a}^{b})$ is implemented by
the extension by $0$ on $f_{a}^{a'}\cup f_{b'}^{b}$\,.
  \end{prop}
  \begin{proof}
    We already know $\dim F_{\left\{0\right\},[a,b],h}^{(p)}=\dim
    F_{\left\{0\right\},[a',b'],h}^{(p)}=\beta^{(p)}(f^{b},f^{a})$ for $p\in \left\{0,\ldots,d\right\}$\,. From the remarks following Definition~\ref{de:dEF}, it then suffices to prove
$$
\vec{d}(F_{\left\{0\right\},[a,b],h},F_{\left\{0\right\},[a',b'],h})=\tilde{O}(e^{-\frac{\min\left\{b'-\tilde{c}_{1},\tilde{c}_{1}-a'\right\}}{h}})\,.
$$
 For
a normalized vector
$\psi\in F_{\left\{0\right\},[a,b],h}$ the exponential decay estimate of
Proposition~\ref{pr:Agmon} (or Hypothesis~\ref{hyp:AgmonLip} for a
more general Lipschitz function $f$) with $r_{h}=0$ and $\lambda_{h}=0$ says
$$
\|e^{\frac{|f(x)-\tilde{c}_{1}|}{h}}\psi\|_{W(f_{a}^{b})}=\tilde{O}(1)\,.
$$
For any $\varepsilon>0$ small enough, take $\chi\in
\mathcal{C}^{\infty}_{0}(f_{a'+\varepsilon}^{b'-\varepsilon};[0,1])$ such that $\chi\equiv 1$
in a neighborhood of
$f^{-1}([a'+2\varepsilon,b'-2\varepsilon])$\,. The form $\chi\psi$
then belongs to  $D(\Delta_{f,f^{-1}([a',b']),h})$ with
$d_{f,h}\psi=(hd\chi)\wedge \psi$\,,
$d_{f,h}^{*}\psi=-h\mathbf{i}_{\nabla\chi}\psi$\,, and therefore
\begin{eqnarray*}
  &&\langle \chi\psi\,,\, \Delta_{f,f^{-1}([a',b']),h} (\chi\psi)\rangle
=\|d_{f,h}(\chi\psi)\|_{L^{2}}^{2}+\|d_{f,h}^{*}(\chi\psi)\|_{L^{2}}^{2}=\tilde{O}(e^{-2\frac{\min\left\{\tilde{c}_{1}-a',
     b'-\tilde{c}_{1}\right\}-2\varepsilon}{h}})\\
\text {and}&&\|\psi-\chi\psi\|_{L^{2}}^{2}=\tilde{O}(e^{-2\frac{\min\left\{\tilde{c}_{1}-a',b'-\tilde{c}_{1}\right\}-2\varepsilon}{h}})\,.
\end{eqnarray*}
Because $0$ is the only exponentially small eigenvalue of
$\Delta_{f,f^{-1}([a',b']),h}$\,, this implies
$$
\textrm{dist}_{L^{2}}(\chi\psi\,, F_{\left\{0\right\},[a',b'],h})=
\tilde{O}(e^{-\frac{\min\left\{\tilde{c}_{1}-a',b'-\tilde{c}_{1}\right\}-2\varepsilon}{h}})\,.
$$
If $F=F_{\left\{0\right\},[a',b'],h}$ is considered as a subspace of
$L^{2}(f_{a}^{b})$ after extension by $0$ on $f_{a}^{a'}\cup
f_{b'}^{b}$\,, the orthogonal projection $\Pi_{F}:L^{2}(f_{a}^{b})\to
F$ is given by $\Pi_{F}u=\Pi_{\left\{0\right\},[a',b'],h}(
u\big|_{f_{a'}^{b'}})$ again extended by $0$ on $f_{a}^{a'}\cup
f_{b'}^{b}$\,.\\
From $\|\Pi_{\left\{0\right\},[a',b'],h}\|\leq 1$ and the exponential
decay estimates for $\psi$\,, 
we deduce, by setting $E=F_{\left\{0\right\},[a,b],h}$\,,
\begin{eqnarray*}
\|(\Pi_{E}-\Pi_{F}\Pi_{E})\psi\|
&=&
\|\psi-\Pi_{\left\{0\right\},[a',b'],h}(\psi\big|_{f_{a'}^{b'}})\|_{L^{2}(f_{a}^{b})}\\
&\leq&
       \|\psi-\chi\psi\|_{L^{2}(f_{a}^{b})}+\|\chi\psi-\Pi_{\left\{0\right\},[a',b'],h}(\chi\psi)\|_{L^{2}(f_{a'}^{b'})}
+\|\chi\psi-\psi\big|_{f_{a'}^{b'}}\|_{L^{2}(f_{a'}^{b'})}\\
&\leq&\tilde{O}(e^{-\frac{\min(\tilde{c}_{1}-a',b'-\tilde{c}_{1})-2\varepsilon}{h}})\,.
\end{eqnarray*}
Since this holds for all $\psi\in E$\,, $\|\psi\|=1$\,, this proves
$\vec{d}(E,F)=\tilde{O}(e^{-\frac{\min(\tilde{c}_{1}-a',b'-\tilde{c}_{1})-2\varepsilon}{h}})$\,,
and we conclude by taking $\varepsilon>0$ arbitrarily small.
  \end{proof}
The above result implies that the mapping $A_{h}:F_{\left\{0\right\},[a,b],h}\to
F_{\left\{0\right\},[a',b'],h}\subset L^{2}(f_{a}^{b})$ defined by $A_{h}\psi=\Pi_{\left\{0\right\},[a',b'],h}(\psi\big|_{f_{a'}^{b'}})$) 
satisfies
$$
\|A_{h}^{*}A_{h}-1\|_{\mathcal{L}(F_{\left\{0\right\},[a,b],h})}
=\tilde{O}(e^{-\frac{\min
    \left\{\tilde{c}_{1}-a',b'-\tilde{c}_{1}\right\}}{h}})$$
    and then
\begin{equation*}
\|A_{h}^{*}A_{h}-1\|_{\mathcal{L}(F_{\left\{0\right\},[a,b],h})}+
\|A_{h}A_{h}^{*}-1\|_{\mathcal{L}(F_{\left\{0\right\},[a',b'],h})}
=\tilde{O}(e^{-\frac{\min
    \left\{\tilde{c}_{1}-a',b'-\tilde{c}_{1}\right\}}{h}})\,.
\end{equation*}
A more accurate version can be given when $a=a'$\,. Actually $\tilde{O}(e^{-\frac{\min
    \left\{\tilde{c}_{1}-a',b'-\tilde{c}_{1}\right\}}{h}})$ is easily replaced by
$\tilde{O}(e^{-\frac{b'-\tilde{c}_{1}}{h}})$ but additionally a small change of
$A_{h}$ allows to improve the estimates in $f_{a}^{\tilde{c}_{1}}$\,.
\begin{prop}
\label{pr:Fabab1}
Keep the same assumptions and conventions 
as in Proposition~\ref{pr:aba1b1} with now
$a=a'$\,. There exists a linear mapping
 $A_{h}:F_{\left\{0\right\},[a,b],h}\to F_{\left\{0\right\},[a,b'],h}$
such that 
$$
\|e^{\frac{b'-f(x)+b'-\tilde{c}_{1}}{h}}[\psi-A_{h}\psi]\|_{W(f_{a}^{\tilde{c}_{1}})}=\tilde{O}(1)\|\psi\|_{L^{2}}
$$
holds for all $\psi\in
F_{\left\{0\right\},[a,b],h}$
 and
\begin{equation}
\label{eq:almunitFaba1b1}
\|A_{h}^{*}A_{h}-1\|_{\mathcal{L}(F_{\left\{0\right\},[a,b],h})}+
\|A_{h}A_{h}^{*}-1\|_{\mathcal{L}(F_{\left\{0\right\},[a',b'],h})}
=\tilde{O}(e^{-\frac{b'-\tilde{c}_{1}}{h}})\,.
\end{equation} 
\end{prop}
\begin{proof}
  The proof is modelled on Lemma~\ref{le:exten1}.\\
Let $\varepsilon\in ]0,\frac{b'-\tilde{c}_{1}}{4}[$\,, and let $\chi,\tilde{\chi}\in
\mathcal{C}^{\infty}(f^{-1}([a,b']);[0,1])$ satisfy 
\begin{eqnarray*}
  &&\chi\equiv 1~\text{in}~
f_{a}^{b'-2\varepsilon}\quad,\quad \chi\equiv 0
~\text{in}~f_{b'-\varepsilon}^{b'}\,,\\
&& \tilde{\chi}\equiv
   0~\text{in}~f_{a}^{\tilde{c}_{1}+\varepsilon}\quad,\quad\tilde{\chi}\equiv 1~\text{in}~f_{\tilde{c}_{1}+2\varepsilon}^{b'}\,.
\end{eqnarray*}
A form $\psi\in F_{\left\{0\right\},[a,b],h}=\ker(\Delta_{f,f^{-1}([a,b]),h})$\,,
$\|\psi\|_{L^{2}}=1$\,,  satisfies
$d_{f,h}\psi=0$ and $d_{f,h}^{*}\psi=0$ in $f_{a}^{b'}$ but has not to
belong to $D(\Delta_{f,f^{-1}([a,b']),h})$\,. 
We introduce
$$
\tilde{\psi}_{\varepsilon}=\chi\psi-\tilde{\chi}(\eta_{1}+\eta_{2})\,,
$$
where
\begin{align*}
   \eta_{1}&=d_{f,f^{-1}([\tilde{c}_{1}+\varepsilon,b'],h)}^{*}(\Delta_{f,f^{-1}([\tilde{c}_{1}+\varepsilon,b'],h)})^{-1}(hd\chi\wedge
   \psi)\\
   &=(\Delta_{f,f^{-1}([\tilde{c}_{1}+\varepsilon,b'],h)})^{-1}[d_{f,h}^{*}(hd\chi\wedge
   \psi)]
\end{align*}   
and
\begin{align*}
\eta_{2}&=-d_{f,f^{-1}([\tilde{c}_{1}+\varepsilon,b'],h)}(
\Delta_{f,f^{-1}([\tilde{c}_{1}+\varepsilon,b'],h})^{-1}(h\mathbf{i}_{\nabla
 \chi} \psi)
\\
&=-(
\Delta_{f,f^{-1}([\tilde{c}_{1}+\varepsilon,b'],h})^{-1}[d_{f,h}(h\mathbf{i}_{\nabla
 \chi} \psi)]\,.
\end{align*}   
Note that the last equality in each of the two above relations
follows from the intertwining relations of 
Proposition~\ref{pr:domain}-4). This implies in particular that
 $\eta_{1},\eta_{2}$ both belong to the domain
$D(\Delta_{f,f^{-1}([\tilde{c}_{1}+\varepsilon],b'),h})$ and hence satisfy the boundary
conditions at $\left\{f=b'\right\}$\,. Since moreover $\psi\in
D(\Delta_{f,f^{-1}([a,b]),h})$ satisfies the boundary conditions at
$\left\{f=a\right\}$\,, $\tilde{\psi}_{\varepsilon}$ then belongs to
$D(\Delta_{f,f^{-1}([a,b'],h)})$\,.\\
Besides, the exponential decay estimates on $\psi$ given by
Proposition~\ref{pr:Agmon} (or Hypothesis~\ref{hyp:AgmonLip}) imply
$$
\|\psi\|_{W(f_{b'-2\varepsilon}^{b'})}=\tilde{O}(e^{-\frac{b'-\tilde{c}_{1}-2\varepsilon}{h}})
$$
and therefore 
$$
\|d_{f,h}^{*}(h d\chi\wedge \psi)\|_{L^{2}}=
\tilde{O}(e^{-\frac{b'-\tilde{c}_{1}-2\varepsilon}{h}})\quad,\quad
\|d_{f,h}(h\mathbf{i}_{\nabla \chi}
\psi)\|_{L^{2}}=\tilde{O}(e^{-\frac{b'-\tilde{c}_{1}-2\varepsilon}{h}}).
$$
The exponential decay estimates stated in Proposition~\ref{pr:Agmon1} (or
Hypothesis~\ref{hyp:AgmonLip}) then imply
$$
\|e^{\frac{b'-f(x)+b'-\tilde{c}_{1}-4\varepsilon}{h}}\eta_{1}\|_{W(f_{\tilde{c}_{1}+\varepsilon}^{b'})}+\|e^{\frac{b'-f(x)+b'-\tilde{c}_{1}-4\varepsilon}{h}}\eta_{2}\|_{W(f_{\tilde{c}_{1}+\varepsilon}^{b'})}=\tilde{O}(1)\,.
$$ 
Set $\omega_{h}=\tilde{\psi}_{\varepsilon}-\Pi_{\left\{0\right\},[a,b'],h}\tilde{\psi}_{\varepsilon}\in
D(\Delta_{f,f^{-1}([a,b']),h})\cap
\ker(\Delta_{f,f^{-1}([a,b']),h})^{\perp}$ and compute
\begin{eqnarray*}
  &&d_{f,f^{-1}([a,b']),h}\omega_{h}=d_{f,f^{-1}([a,b']),h}\tilde{\psi}_{\varepsilon}
\stackrel{d_{f,h}\psi=0}=-hd\tilde{\chi}\wedge(\eta_{1}+\eta_{2})\\
&& d_{f,f^{-1}([a,b']),h}^{*}\omega_{h}=d_{f,f^{-1}([a,b']),h}^{*}\tilde{\psi}_{\varepsilon}
\stackrel{d_{f,h}^{*}\psi=0}{=}h\mathbf{i}_{\nabla \tilde{\chi}}(\eta_{1}+\eta_{2})\\
&&\Delta_{f,f^{-1}([a,b']),h}\omega_{h}=r_{h}=(1-\Pi_{\left\{0\right\},[a,b'],h})r_{h}\\
&&
\|e^{\frac{b'-f(x)+b'-\tilde{c}_{1}-4\varepsilon}{h}}r_{h}\|_{L^{2}(f_{a}^{b'})}=\tilde{O}(1)\,.
\end{eqnarray*}
The ``orthogonality lemma'' (Lemma~\ref{le:orthlem}) with
$\omega_{h}=\tilde{\psi}_{\varepsilon}-\Pi_{\left\{0\right\},[a,b'],h}\tilde{\psi}_{\varepsilon}$
yields
$$
\|
e^{\frac{b'-f(x)+b'-\tilde{c}_{1}-4\varepsilon}{h}}[\tilde{\psi}_{\varepsilon}-\Pi_{\left\{0\right\},[a,b'],h}\tilde{\psi}_{\varepsilon}]
\|_{W(f_{a}^{b'})}=\tilde{O}(1)\,.
$$
By defining 
$A^{\varepsilon}_{h}\psi:=\Pi_{\left\{0\right\},[a,b'],h}\tilde{\psi}_{\varepsilon}\in F_{\left\{0\right\},[a,b'],h}\subset L^{2}(f_{a}^{b})$\,,
it then follows 
from the latter relation and from the relation
$\psi\equiv \tilde{\psi}_{\varepsilon}$ in $f_{a}^{\tilde{c}_{1}+\varepsilon}$
that 
$$
\|e^{\frac{b'-f(x)+b'-\tilde{c}_{1}}{h}}[\psi-A^{\varepsilon}_{h}\psi]\|_{W(f_{a}^{\tilde{c}_{1}})}=\tilde{O}(e^{\frac {4\varepsilon} h})
$$
and 
\begin{align*}
\|\psi-A^{\varepsilon}_{h}\psi\|_{L^{2}(f_{a}^{b})}
&\leq\|\psi-\tilde{\psi}_{\varepsilon}\|_{L^{2}(f_{a}^{b})}+\|\tilde{\psi}_{\varepsilon}-A^{\varepsilon}_{h}\psi\|_{L^{2}(f_{a}^{b'})}
\\
&\hspace{-1cm}\leq
\|(1-\chi)\psi\|_{L^{2}(f_{a}^{b'})}+
\|\tilde \chi(\eta_{1}+\eta_{2})\|_{L^{2}(f_{a}^{b'})}
+\|\psi\|_{L^{2}(f_{b'}^{b})}+\|\tilde{\psi}_{\varepsilon}-A^{\varepsilon}_{h}\psi\|_{L^{2}(f_{a}^{b'})}
\\
&= \tilde{O}(e^{-\frac{b'-\tilde{c}_{1}-4\varepsilon}{h}})\,.
\end{align*}
In order to conclude, it thus just remains to choose $\varepsilon$
depending on $h\in ]0,h_{0}[$ in a proper way. To do so, note that when
$\varepsilon=\frac{1}{n+1}$ with $n\in\nz$ large enough to ensure $\varepsilon\in ]0,\frac{b'-\tilde{c}_{1}}{4}[$\,, there exists $h_{n}>0$ such that 
for every $h\in]0,h_{n}[$\,,
$$
\|e^{\frac{b'-f(x)+b'-\tilde{c}_{1}}{h}}[\psi-A^{\varepsilon}_{h}\psi]\|_{W(f_{a}^{\tilde{c}_{1}})}\leq e^{\frac{5}{(n+1)h}}
\quad\text{and}
\quad
\|\psi-A^{\varepsilon}_{h}\psi\|_{L^{2}(f_{a}^{b})}
\leq e^{\frac{5}{(n+1)h}}\,e^{-\frac{b'-\tilde{c}_{1}}{h}}\,.
$$
The sequence $(h_{n})_{n\in\nz}$ can be chosen decreasing and 
it then suffices to define
$A_{h}:=A_{h}^{\frac 1{n+1}}$ when $h\in [h_{n+1},h_{n}[$\,.
\end{proof}

\subsubsection{Interactions of solutions to $d_{f,h}\omega=0$ with local
  spectral problems}
\label{sec:interdfh}
We conclude this section with a result which will be used in the
construction and analysis of global quasimodes (see Section~\ref{sec:accanN}).
It provides information about solutions to $d_{f,h}\omega=0$ in
$f^{\tilde{c}_{1}}$\,, in particular how the exponential decay can be combined
with local spectral information.
\begin{prop}
\label{pr:interdfh}
Assume Hypothesis~\ref{hyp:1vc} and $a_{0}\leq a<\tilde{c}_{1}<b'<b$\,. 
Let $\delta(h)>0$ satisfy $\lim_{h\to 0}\delta(h)=0$ and 
let the family 
$(\omega_{h})_{h\in]0,h_{0}[}$  satisfy 
$\omega_{h}\in W(f_{a}^{\tilde c_1-\delta(h)};\Lambda T^{*}M)$
and
$d_{f,h}\omega_{h}=0$  in
$f_{a_{0}}^{\tilde{c}_{1}-\delta(h)}$  with
$$
\|e^{\frac{f(x)-a_{0}}{h}}\omega_{h}\|_{W(f_{a}^{\tilde{c}_{1}-\delta(h)})}=\tilde{O}(1)\,.
$$
Take any cut-off function $\chi\in
  \mathcal{C}^{\infty}_{0}(f^{-1}([a,\tilde{c}_{1}[);[0,1])$  such that $\chi\equiv 1$
  in a neighborhood of $\left\{f=a\right\}$
 and assume that $h>0$ is small enough so that $\supp \chi \subset [a, \tilde{c}_{1}-\delta(h)[ $\,.
\begin{description}
\item[i)] The form
  $\Pi_{\left\{0\right\},[a,b],h}[d_{f,h}(\chi\omega_{h})]=\Pi_{\left\{0\right\},[a,b],h}[(hd\chi)\wedge
  \omega_{h}]$ does not depend on the choice of the cut-off function $\chi$\,.
\item[ii)] If $\Pi_{\left\{0\right\},[a,b],h}[d_{f,h}(\chi\omega_{h})]=0$\,,
  then there exists a family of similar cut-off functions
  $\chi_{h}$  such that
  $\tilde{\omega}_{h}=\chi_{h}\omega_{h}-d_{f,f^{-1}[a,b],h}^{*}(\Delta_{f,f^{-1}[a,b],h}^{\perp})^{-1}[(hd\chi_{h})\wedge\omega_{h})]$\,,
where, in the r.h.s.,  $\chi_{h}$ in the first term is extended by $1$ and 
the second term is extended by $0$ in $f_{a_{0}}^{a}$\,,
  satisfies
  \begin{eqnarray*}
&& \tilde{\omega}_{h}\equiv \omega_{h}\quad\text{in}~f_{a_{0}}^{a}\,,\\
    && d_{f,h}\tilde{\omega}_{h}=0\quad\text{in}~f_{a_{0}}^{b}\,,\\
\text{and}&&
\|e^{\frac{f(x)-a_{0}}{h}}\tilde{\omega}_{h}\|_{W(f_{a}^{b})}=\tilde{O}(1)\,.
  \end{eqnarray*}
\item[iii)] If $A_{h}:F_{\left\{0\right\},[a,b],h}\to F_{\left\{0\right\},[a,b'],h}$ is the operator
  introduced in Proposition~\ref{pr:Fabab1}, then 
 for any
 $\psi\in
F_{\left\{0\right\},[a,b],h}$\,, 
  the quantity 
$\langle d_{f,h}(\chi\omega_{h})\,,\,\psi-A_{h}\psi\rangle$ does not
depend on the choice of $\chi$ and 
$$
\forall \psi\in F_{\left\{0\right\},[a,b],h}\,,\quad
\langle d_{f,h}(\chi\omega_{h})\,,\,\psi-A_{h}\psi\rangle
=\tilde{O}(e^{-\frac{b'-a_{0}+b'-\tilde{c}_{1}}{h}})\|\psi\|_{L^{2}}\,.
$$
\end{description}
\end{prop}
\begin{proof}
\noindent\textbf{i)} Let $\chi_{1},\chi_{2}$ be two cut-off functions
like $\chi$ in our assumptions. Then
$\chi_{1}\omega_{h}-\chi_{2}\omega_{h}$ belongs to
$D(d_{f,f^{-1}[a,b],h})$ and 
$$
d_{f,f^{-1}[a,b],h}(\chi_{1}\omega_{h}-\chi_{2}\omega_{h})=d_{f,h}(\chi_{1}\omega_{h})-d_{f,h}(\chi_{2}\omega_{h})\,.
$$
We simply conclude with the commutation
$$
\Pi_{\left\{0\right\},[a,b],h}d_{f,f^{-1}([a,b]),h}=d_{f,f^{-1}([a,b]),h}\Pi_{\left\{0\right\},[a,b],h}=0\,.
$$
\noindent\textbf{ii)} When $\Pi_{\left\{0\right\},[a,b],h}[d_{f,h}(\chi
\omega_{h})]=0$\,, \textbf{i)} ensures that the latter relation is also
satisfied if we replace $\chi$ by $\chi_{\varepsilon}$  with
$\chi_{\varepsilon}\equiv 1$ in  $f_{a}^{a'-\varepsilon}$
and $\chi_{\varepsilon}=0$ in $f_{a'+\varepsilon}^{\tilde{c}_{1}}$ for 
$a'=\frac{a+\tilde{c}_{1}}{2}$ and some
$\varepsilon\in ]0,\frac{\tilde{c}_{1}-a}{2}[$\,. The a priori estimates on
$\omega_{h}$
and 
$\supp (hd\chi_{\varepsilon})\wedge \omega_{h}\subset f^{-1}([a'-\varepsilon,a'+\varepsilon])$
 imply 
$$
\|((hd\chi_{\varepsilon})\wedge \omega_{h}\|_{L^{2}(f_{a}^{b})}
=\tilde{O}(e^{-\frac{a'-a_{0}-\varepsilon}{h}})\,.
$$
The orthogonality lemma, Lemma~\ref{le:orthlem}, then implies that
 $$
\eta_{\varepsilon}=d_{f,f^{-1}([a,b])
   ,h}^{*}(\Delta_{f,f^{-1}([a,b],h)}^{\perp})^{-1}[(hd\chi_{\varepsilon})\wedge
 \omega_{h}]
$$
(is well defined and) satisfies
$$
\|e^{\frac{|f(x)-a'|-\varepsilon}{h}}\eta_{\varepsilon}\|_{L^{2}(f_{a}^{b})}
=\tilde{O}(e^{-\frac{a'-a_{0}-\varepsilon}{h}})\,.
$$
Since moreover 
$d_{f,h}(\chi_{\varepsilon}\omega_{h})=(hd\chi_{\varepsilon})\wedge
\omega_{h}=(1-\Pi_{\left\{0\right\},[a,b],h})[(hd\chi_{\varepsilon}\wedge \omega_{h})]$
belongs to  $D(d_{f,f^{-1}([a,b],h)})$\,, we can write
$$
d_{f,f^{-1}([a,b]),h}\eta_{\varepsilon}=\Delta_{f,f^{-1}([a,b]),h}(\Delta_{f,f^{-1}([a,b]),h}^{\perp})^{-1}((hd\chi_{\varepsilon})\wedge
\omega_{h})=(hd\chi_{\varepsilon})\wedge\omega_{h}\,.
$$
Using in addition $d_{f,f^{-1}([a,b]),h}^{*}\eta_{\varepsilon}=0$\,, 
we deduce
$$
\|e^{\frac{|f(x)-a'|-\varepsilon}{h}}\eta_{\varepsilon}\|_{W(f_{a}^{b})}=\tilde{O}(e^{-\frac{a'-a_{0}-\varepsilon}{h}})\,.
$$
If $\eta_{\varepsilon}$ denotes 
the extension by $0$ in $f_{a_{0}}^{a}$ of $\eta_{\varepsilon}\in
D(d_{f,f^{-1}([a,b],h)})$\,, it still belongs to
$D(d_{f,f^{-1}([a_{0},b]),h})$ and 
solves $d_{f,h}\eta_{\varepsilon}=(hd\chi_{\varepsilon})\wedge \omega_{h}$ in 
$f_{a_{0}}^{a}\cup f_{a}^{b}$\,. 
We have thus proved that
$\tilde{\omega}_{\varepsilon}:=\chi_{\varepsilon}\omega_{h}-\eta_{\varepsilon}$
satisfies
$$
d_{f,h}\tilde{\omega}_{\varepsilon}=0\ \ \text{in}\ \ f_{a_{0}}^{b}\quad\text{and}\quad\|e^{\frac{f(x)-a_{0}}{h}}\tilde{\omega}_{\varepsilon}\|_{W(f_{a}^{b})}=\tilde{O}(e^{\frac{2\varepsilon}{h}})\,.
$$ 
We then end the proof by choosing conveniently $\varepsilon$
depending on $h\in ]0,h_{0}[$ as we did
at the end of the proof of Proposition~\ref{pr:Fabab1}: when
$\varepsilon=\frac{1}{n+1}$\,, take $h_{n}>0$ such that 
$$
\forall h\in ]0,h_{n}[\,,\ \ 
\|e^{\frac{f(x)-a_{0}}{h}}\tilde{\omega}_{\varepsilon}\|_{W(f_{a}^{b})}\leq e^{\frac{3}{(n+1)h}}
$$
with
 $(h_{n})_{n\in\nz}$ decreasing, and choose
$\chi_{h}:=\chi_{\frac{1}{n+1}}$ when $h\in [h_{n+1},h_{n}[$\,.
\\
\noindent\textbf{iii)} Since
$$
\langle d_{f,h}(\chi\omega_{h})\,, \psi-A_{h}\psi\rangle
=
\langle \Pi_{\left\{0\right\},[a,b],h}[d_{f,h}(\chi\omega_{h})]\,,\, \psi\rangle
-\langle \Pi_{\left\{0\right\},[a,b'],h}[d_{f,h}(\chi\omega_{h})]\,,\, A_{h}\psi\rangle
$$
does not depend on $\chi$\,, we may
take the preceding $\chi=\chi_{\varepsilon}$\,.
Owing to Proposition~\ref{pr:Fabab1}, we deduce
\begin{align*}
|\langle d_{f,h}(\chi_{\varepsilon}\omega_{h})\,,\,
\psi-A_{h}\psi\rangle|
&\leq \|(hd\chi_{\varepsilon})\wedge
\omega_{h}\|_{L^{2}(f_{a'-\varepsilon}^{a'+\varepsilon})}\|\psi-A_{h}\psi\|_{L^{2}(f_{a'-\varepsilon}^{a'+\varepsilon})}
\\
&= \tilde{O}(e^{-\frac{a'-a_{0}-\varepsilon}{h}})\times
\tilde{O}(e^{-\frac{b'-a'-\varepsilon+b'-\tilde{c}_{1}}{h}})\,\|\psi\|_{L^{2}}\,.
\end{align*}
Since this holds for every $\varepsilon>0$ small enough, this yields the result.
\end{proof}

\section{Rough estimates for several ``critical values''}
\label{sec:rough}
In this section,  we give first estimates for 
the exponentially small eigenvalues of
$\Delta_{f,f^{-1}([a,b],h)}$\,.
We work under the following assumption which, like
Hypothesis~\ref{hyp:1vc} in Section~\ref{sec:locpbs}, gathers
Hypothesis~\ref{hyp:mainf} or (Hypothesis~\ref{hyp:Lipbar} and
Hypothesis~\ref{hyp:AgmonLip}), and specify some notations.
\begin{hyp}
\label{hyp:cN}
The function $f$ satisfies Hypothesis~\ref{hyp:mainf}, or
more generally Hypothesis~\ref{hyp:Lipbar} and Hypothesis~\ref{hyp:AgmonLip}, and we choose $\eta_{f}$
such that
$$
0< \eta_{f}<\frac{1}{2}\min_{1<n\leq N_{f}}|c_{n}-c_{n-1}|\,.
$$
In addition, $a,b$\,, $-\infty\leq
a< b\leq +\infty$\,, are not ``critical values'' of $f$:
$a,b\not\in \left\{c_{1},\ldots,c_{N_{f}}\right\}$\,.
\end{hyp}
\subsection{Bar code associated with $f$}
\label{sec:bar code}
We refer to Appendix~\ref{app:perscohom} for details and simply recall
the useful notations. We already 
mentionned in Subsection~\ref{sec:genass}
that 
Hypothesis~\ref{hyp:Lipbar} implies Hypothesis~\ref{hyp:weakreg} in
the beginning of Appendix~\ref{app:perscohom} (this is actually proved
in Subsection~\ref{sec:moregenLip}).\\ 
Under the assumption that $M$ is compact and
 $f$ has a finite number of ``critical values''
$c_{1}<\ldots <c_{N_{f}}$\,, there is a bar code 
$\mathcal{B}={\cal B}(f)=([a_{\alpha},b_{\alpha}[)_{\alpha\in A}$ where
$A$ is finite, $-\infty<a_{\alpha}<b_{\alpha}\leq +\infty$\,,
$a_{\alpha}\in \left\{c_{1},\ldots,c_{N_{f}}\right\}$\,, $b_{\alpha}\in
\left\{c_{2},\ldots,c_{N_{f}},+\infty\right\}$\,. The set $A$ is graded
according to $A=\sqcup_{p=0}^{\dim M} A^{(p)}$ so that, for $\alpha\in
A^{(p)}$\,, the grading of endpoints of the corresponding  bar is given by
 $[a_{\alpha},b_{\alpha}[=[a_{\alpha}^{(p)},b_{\alpha}^{(p+1)}[$\,. It contains all the
information about the relative cohomology groups $H(f^{b},f^{a};\rz)$
when $a<b$\,, $a,b\not\in\left\{c_{1},\ldots,c_{N_{f}}\right\}$\,.\\
More precisely here is the situation when  $a<b$ are not ``critical values''.
We forget  the bars with no end point in $]a,b[$\,, and among the
remaining ones we distinguish the ones with two endpoints in $]a,b[$:
\begin{eqnarray}
\label{eq:defAab}
 A^{*}(a,b)&=&\left\{\alpha\in A^{*},
               [a_{\alpha}^{*},b_{\alpha}^{*+1}[\cap ]a,b[
               \not\in\left\{\emptyset,]a,b[\right\} \right\}\,,\\
\label{eq:defAcab}
 A_{c}^{*}(a,b)&=&\left\{\alpha\in A^{*}(a,b),
   [a_{\alpha}^{*},b_{\alpha}^{*+1}[\cap]a,b[~\text{relatively~compact~in
   }~]a,b[\right\}\,,\\
\nonumber
&&
\alpha\in A^{*}(a,b) \Leftrightarrow
   a<a_{\alpha}^{*}<b~\text{or}~a<b_{\alpha}^{*+1}<b\,,\\
\nonumber
&&
\alpha\in A_{c}^{*}(a,b)\Leftrightarrow a<a_{\alpha}^{*}<
   b_{\alpha}^{*+1}<b\,.
\end{eqnarray}
We now partition the endpoints of the bars, multiple value being
distinguished by the index $\alpha\in A(a,b)$\,, according to
\begin{eqnarray}
\label{eq:defXab}
  {\cal X}^{*}(a,b)&=&\left\{(\alpha,a_{\alpha}^{*})\,, \alpha\in
                       A_{c}^{*}(a,b)\right\}\\
\label{eq:defYab}
{\cal Y}^{*}(a,b)&=&\left\{(\alpha,b_{\alpha}^{*}), \alpha\in
                       A_{c}^{*-1}(a,b)\right\}\\
\label{eq:defZab}
{\cal Z}^{*}(a,b)&=&\left\{(\alpha,a_{\alpha}^{*})\,, \alpha\in A^{*}(a,b)\setminus
                     A_{c}^{*}(a,b)\,, a<a_{\alpha}<b\right\}\\
\nonumber
&&\hspace{2cm}\sqcup
\left\{(\alpha,b_{\alpha}^{*})\,, \alpha\in A^{*-1}(a,b)\setminus
                     A_{c}^{*-1}(a,b), a<b_{\alpha}^{*}<b\right\}\,,\\
\label{eq:defJab}
{\cal J}^{*}(a,b)&=&{\cal X}^{*}(a,b)\sqcup {\cal Y}^{*}(a,b)\sqcup {\cal Z}^{*}(a,b)\,.
\end{eqnarray} 
Those definitions are illustrated in Figure~7: the
degrees of the bars and of the corresponding endpoints are
  indicated.
The bars in $A_{c}(a,b)$ are the ones with two endpoints in $]a,b[$
and  the critical values lying in $]a,b[$ are relabelled $\tilde{c}_{1}<\ldots<\tilde{c}_{N}$\,.
\vspace{1cm}
\figinit{0.7cm}
\figpt1000:(-8,-1)%
\figpt1001:(8,-1)%
\figpt1002:(5,-1)%
\figpt 0:(-7,-1)%
\figpt 1:(7,-1)%
\figpt 2:(-6,-1)%
\figpt 3:(-3,-1)%
\figpt 4:(3,-1)%
\figpt 5:(6,-1)
\figpt 6:(-1,-1)%
\figpt 7:(2,-1)%

\figpt 99:(-8,0)
\figpt 100:(-7,0)%
\figpt 101:(-3,0)%
\figpt 200:(-5,0,)
\figpt 102:(-3,4)%
\figpt 103:(3,4)%
\figpt 202:(0,4)
\figpt 104:(-6,0.8)%
\figpt 105:(-3,0.8)%
\figpt 204:(-4.5,0.8)
\figpt 106:(-6,1.6)%
\figpt 107:(3,1.6)%
\figpt 206:(-1.5,1.6)
\figpt 110:(-3,2.4)%
\figpt 111:(3,2.4)%
\figpt 210:(0,2.4)
\figpt 116:(3,3.2)%
\figpt 117:(7,3.2)%
\figpt 216:(5,3.2)
\figpt 118:(8,3.2)
\figpt 120:(5,4.8)
\figpt 121:(7,4.8)
\figpt 220:(6,4.8)
\figpt 119:(8,4.8)
\figpt 122:(-8,5.6)
\figpt 123:(8,5.6)
\figpt 222:(0,5.6)

\figdrawbegin{}
\figset(width=1.3)
\figdrawline[1000,6]
\figdrawline[7,1001]
\figset(dash=2)
\figdrawline[6,7]
\figset(dash=1)
\figset(width=2)
\figdrawline[100,101]
\figdrawline[102,103]
\figdrawline[104,105]
\figdrawline[106,107]
\figdrawline[110,111]
\figdrawline[116,117]
\figdrawline[120,121]
\figset(dash=4)
\figdrawline[99,100]
\figdrawline[117,118]
\figdrawline[120,119]
\figset(width=0.3)
\figdrawline[122,123]

\figdrawend

\figvisu{\figBoxA}{\textbf{Figure 7:} ${\cal X}^{*}={\cal X}^{*}(a,b)$ (lower)\,, 
${\cal Y}^{*}={\cal Y}^{*}(a,b)$ (upper) \,, 
${\cal Z}^{*}={\cal Z}^{*}(a,b)$ (lonely)}{
\figwriten 222: {killed in $]a,b[$}(0.1)
\figwrites 200:${\scriptstyle(p-1)}$(0.1)
\figwrites 202:${\scriptstyle(p)}$(0.1)
\figwrites 204:${\scriptstyle(p-1)}$(0.1)
\figwrites 206:${\scriptstyle(p)}$(0.1)
\figwrites 210:${\scriptstyle(p-1)}$(0.1)
\figwrites 216:${\scriptstyle(p)}$(0.1)
\figwrites 220:${\scriptstyle(p-1)}$(0.1)
\figwritee 101:${{\cal Z}^{(p)}}$(0.2)
\figwritew 102:${{\cal X}^{(p)}}$(0.2)
\figwritee 103:${{\cal Y}^{(p+1)}}$(0.2)
\figwritew 104:${{\cal X}^{(p-1)}}$(0.2)
\figwritee 105:${{\cal Y}^{(p)}}$(0.2)
\figwritew 106:${{\cal X}^{(p)}}$(0.2)
\figwritee 107:${{\cal Y}^{(p+1)}}$(0.2)
\figwritew 110:${{\cal X}^{(p-1)}}$(0.2)
\figwritee 111:${{\cal Y}^{(p)}}$(0.2)
\figwritew 116:${{\cal Z}^{(p)}}$(0.2)
\figwritew 120:${{\cal Z}^{(p-1)}}$(0.2)
\figset write(mark=$+$)
\figwrites 1002:$\tilde{c}_{N+1}$(0.3)
\figwrites 0: $a$(0.3)
\figwrites 1:$b$(0.3)
\figwrites 2:$\tilde{c}_{1}$(0.3)
\figwrites 3:$\tilde{c}_{2}$(0.3)
\figwrites 4:$\tilde{c}_{N}$(0.3)
\figset write(mark=$\bullet$)
\figwrites 102:$$(0.2)
\figwrites 104:$$(0)
\figwrites 106:$$(0)
\figwrites 110:$$(0)
\figwrites 116:$$(0)
\figwrites 120:$$(0)
\figwrites 122:$$(0)
}
\centerline{\box\figBoxA}

\vspace{1cm}
Then the relative Betti number are given by
\begin{equation}
  \label{eq:relBetti}
\beta^{(p)}(f^{b},f^{a})=\dim
H^{p}(f^{b},f^{a};\rz)=\dim F_{\left\{0\right\},[a,b],h}=\sharp \mathcal{Z}^{(p)}(a,b)\,,
\end{equation}
which counts the number of degree $p$ 
endpoints of the bar code
lying lonely  in $]a,b[$\,.\\
The rest of this section shows that there are exactly
$\sharp {\cal J}^{(p)}(a,b)$
 exponentially
small eigenvalues of $\Delta^{(p)}_{f,f^{-1}([a,b]),h}$\,, 
 and provides a priori estimates on the size of the non zero ones.
\subsection{Counting exponentially small eigenvalues}
\label{sec:countexpsm}

\begin{prop}
\label{pr:comptageest}
Under Hypothesis~\ref{hyp:cN} and with the notations of
Subsection~\ref{sec:bar code}, the exponentially small 
 eigenvalues of $\Delta_{f,f^{-1}([a,b]),h}$
are counted according to:
\begin{eqnarray}
  \label{eq:dimker}
  \dim\ker (\Delta_{f,f^{-1}([a,b],h)}^{(p)})&=&\sharp {\cal
                                                 Z}^{(p)}(a,b)\\
\label{eq:dimFab}
\dim F^{(p)}_{[0,\tilde{o}(1)],[a,b],h}&=&\sharp {\cal J}^{(p)}(a,b)=\sharp
  {\cal X}^{(p)}(a,b)+\sharp {\cal Y}^{(p)}(a,b)+\sharp {\cal
  Z}^{(p)}(a,b)\,,
\end{eqnarray}
where the second quality holds for $h\in ]0,h_{\varepsilon}[$ when
$\tilde{o}(1)$ is replaced by $e^{-2\frac{\eta_{f}-2\varepsilon}{h}}$ for
$\varepsilon\in ]0,\frac{\eta_{f}}{2}[$\,.
\end{prop}

Note that the right-hand side of \eqref{eq:dimFab} is nothing but the
total number of degree $p$ endpoints of the bar code lying in
$]a,b[$\,. This counting also says that the $\tilde{o}(1)$ eigenvalues
of $\Delta_{f,f^{-1}([a,b]),h}^{(p)}$ are actually $\tilde{O}(e^{-\frac{2\eta_{f}}{h}})$\,.
\begin{proof}
Equality \eqref{eq:dimker}, which was already stated in
Subsection~\ref{sec:bar code}, is proved in Appendix~\ref{app:perscohom}.
Equality \eqref{eq:dimFab} relies on exponential
 decay estimates 
and on the result in the case $\left\{c_{1},\ldots,c_{N_{f}}\right\}\cap [a,b]=\{\tilde{c}_{1}\}
\subset]a,b[$
 stated in Proposition~\ref{pr:exp0}. \\
The ``critical values'' of $f$ lying in $]a,b[$ are
 relabelled as $a<\tilde{c}_{1}<\ldots<\tilde{c}_{N}<b$\, according to
$$
]a,b[\cap \left\{c_{1},\ldots,c_{N_{f}}\right\}=[a,b]\cap
\left\{c_{1},\ldots,c_{N_{f}}\right\}=\left\{\tilde{c}_{1},\ldots, \tilde{c}_{N}\right\}\,.
$$ 
Consider the disjoint union $\overline{\Omega}$:
$$
\overline{\Omega}
=\mathop{\bigsqcup}_{j=1}^{N}f^{-1}([\tilde{c}_{j}-\eta_{f},\tilde{c}_{j}+\eta_{f}]\cap[a,b])
$$
for which the associated boundary Witten Laplacian is
\begin{equation}
  \label{eq:laplspl}
\Delta_{f,\overline{\Omega},h}=
\bigoplus_{j=1}^{N}\Delta_{f,f^{-1}([a,b]\cap
  [\tilde{c}_{j}-\eta_{f},\tilde{c}_{j}+\eta_{f}]),h}
\end{equation}
By Proposition~\ref{pr:exp0}, we know that the $\tilde{o}(1)$ eigenvalues
of $\Delta_{f,\overline{\Omega},h}$ are equal to $0$\,. For
$\varepsilon\in ]0,\eta_{f}/2[$\,, take $\chi\in
\mathcal{C}^{\infty}(\overline{\Omega};[0,1])$ such that $\chi(x)=1$
if
$\min_{1\leq j\leq N}|f(x)-\tilde{c}_{j}|\leq \eta_{f}-\varepsilon$ and
$\chi(x)=0$ if $\min_{1\leq j\leq N}|f(x)-\tilde{c}_{j}|\geq
\eta_{f}-\epsilon/2$\,.
For any $\omega\in \ker(\Delta_{f,\overline{\Omega},h}^{(p)})$\,, $\|\omega\|_{L^{2}}=1$\,, 
Proposition~\ref{pr:Agmon} (or Hypothesis~\ref{hyp:AgmonLip})  gives
\begin{equation*}
d_{f, h}(\chi
     \omega)=(hd\chi)\wedge\omega=\tilde{O}(e^{-\frac{\eta_{f}-\varepsilon}{h}})
\quad
\text{and}
\quad
d_{f,h}^{*}(\chi\omega)=h\mathbf{i}_{\nabla\chi}\omega=\tilde{O}(e^{-\frac{\eta_{f}-\varepsilon}{h}})\,.
\end{equation*}
Meanwhile our choice of $\chi$ ensures $\chi\omega\in
D(\Delta_{f,f^{-1}([a,b]),h}^{(p)})$ with now
\begin{equation}
\label{eq:numestimdf}
\|d_{f,f^{-1}([a,b]),h}(\chi\omega)\|_{L^{2}}^{2}+\|d_{f,f^{-1}([a,b]),h}^{*}(\chi\omega)\|_{L^{2}}^{2}\leq \tilde{O}(e^{-2\frac{\eta_{f}-\varepsilon}{h}})\,.
\end{equation}
Since
$\|\chi\omega-\omega\|_{L^{2}}=\tilde{O}(e^{-\frac{\eta_{f}-\varepsilon}{h}})$\,,
the spectral decomposition of $\Delta_{f,f^{-1}([a,b]),h}^{(p)}$ ensures
$$
\vec{d}(\ker(\Delta_{f,\overline{\Omega},h}^{(p)}),F^{(p)}_{[0,e^{-2\frac{\eta_{f}-2\varepsilon}{h}}],[a,b],h})=\tilde{O}(e^{-\frac{\varepsilon}{h}})
$$
and then (see indeed the lines following Definition~\ref{de:dEF})
\begin{equation}
\label{eq:ineqdim}
\dim(\ker(\Delta_{f,\overline{\Omega},h}^{(p)}))
\leq \dim
F_{[0,e^{-2\frac{\eta_{f}-2\varepsilon}{h}}],[a,b],h}^{(p)}\,,
\end{equation}
for $h\in ]0,h_{\varepsilon}[$ with $h_{\varepsilon}>0$ small
enough.\\ 
Reciprocally, when $\omega\in
F_{[0,e^{-\frac{\varepsilon}{h}}],[a,b],h}^{(p)}$\,, the exponential decay
estimates of Proposition~\ref{pr:Agmon} (or
Hypothesis~\ref{hyp:AgmonLip}) lead again to
\begin{equation*}
(hd\chi)\wedge\omega=\tilde{O}(e^{-\frac{\eta_{f}-\varepsilon}{h}})
\quad
\text{and}
\quad
h\mathbf{i}_{\nabla\chi}\omega=\tilde{O}(e^{-\frac{\eta_{f}-\varepsilon}{h}})
\end{equation*}
and then to
$$
\|d_{f,h}(\chi \omega)\|_{L^{2}}^{2}+\|d_{f,h}^{*}(\chi
\omega)\|_{L^{2}}^{2}\leq \tilde{O}(e^{-\frac{\varepsilon}{h}})
$$
 with now $\chi\omega\in
D(\Delta_{f,\overline{\Omega},h}^{(p)})$\,.
Again, with $\|\chi \omega
-\omega\|_{L^{2}}=\tilde{O}(e^{-\frac{\eta_{f}-\varepsilon}{h}})$\,,
the spectral decomposition of  $\Delta_{f,\overline{\Omega},h}^{(p)}$\,,
with
$1_{[0,e^{-\frac{\varepsilon}{2h}}]}(\Delta_{f,\overline{\Omega},h}^{(p)})=1_{\left\{0\right\}}(\Delta_{f,\overline{\Omega},h}^{(p)})$\,,
leads to
$$
\vec{d}(F_{[0,e^{-\frac{\varepsilon}{h}}],[a,b],h}, \ker(\Delta_{f,\overline{\Omega},h}^{(p)}))=\tilde{O}(e^{-\frac{\varepsilon}{4h}})
$$
and then to
$$
\dim F_{[0,e^{-\frac{\varepsilon}{h}}],[a,b],h}^{(p)}\leq
\dim \ker \Delta_{f,\overline{\Omega},h}^{(p)}\leq \dim
F_{[0,e^{-2\frac{\eta_{f}-2\varepsilon}{h}}],[a,b],h}^{(p)}\,,
$$
for $h\in]0,h_{\varepsilon}[$\,, $h_{\varepsilon}>0$ small enough, where the last inequality follows from
\eqref{eq:ineqdim}.
\\
In particular, we deduce that for every $\varepsilon>0$ small enough:
\begin{equation}
  \label{eq:numegF}
F_{[0,e^{-\frac{\varepsilon}{h}}],[a,b],h}^{(p)}=F_{[0,e^{-2\frac{\eta_{f}-2\varepsilon}{h}}],[a,b],h}^{(p)}
\end{equation}
and
$$
\dim F_{[0,e^{-2\frac{\eta_{f}-2\varepsilon}{h}}],[a,b],h}^{(p)}=\dim
\ker \Delta_{f,\overline{\Omega},h}^{(p)}\,.
$$
We conclude
with 
\begin{eqnarray*}
\dim \ker \Delta_{f,\overline{\Omega},h}^{(p)}
&=&\sum_{j=1}^{N}\beta^{(p)}(f^{\min(b,
    \tilde{c}_{j}+\eta_{f})}, f^{\max(a, \tilde{c}_{j}-\eta_{f})})
\\
&=&
\sum_{j=1}^{N}\sharp {\cal Z}^{(p)}(\max(a,\tilde{c}_{j}-\eta_{f}),\min(b,\tilde{c}_{j}+\eta_{f}))
=\sharp {\cal J}^{(p)}(a,b)\,,
\end{eqnarray*}
the total number of degree $p$ endpoints of the bar code lying in $]a,b[$\,.
\end{proof}
We have also proved the following result.
\begin{prop}
\label{pr:distkernsum}
In the framework of Proposition~\ref{pr:comptageest} and when
 $\Delta_{f,\overline{\Omega},h}$ is the operator defined in
\eqref{eq:laplspl}, the following inequality holds:
$$
\vec{d}(F_{[0,\tilde{o}(1)],h}^{(p)},\ker(\Delta_{f,\overline{\Omega},h}^{(p)}))+\vec{d}(\ker(\Delta_{f,\overline{\Omega},h}^{(p)}),F_{[0,\tilde{o}(1)],h}^{(p)})=\tilde{O}(e^{-\frac{\eta_{f}}{h}})\,.
$$
\end{prop}
\begin{proof}
 By \eqref{eq:numegF} we know that for 
 $\varepsilon>0$  small enough
$$
\vec{d}(\ker(\Delta_{f,\overline{\Omega},h}^{(p)}),F_{[0,e^{-\frac{\varepsilon}{h}}],h}^{(p)})=
\vec{d}(\ker(\Delta_{f,\overline{\Omega},h}^{(p)}),F_{[0,e^{-\frac{2\eta_{f}-2\varepsilon}{h}}],h}^{(p)})\,,
$$
while we are in cases  with $\vec{d}(A,B)=\vec{d}(B,A)<1$ by the result
of Proposition~\ref{pr:comptageest}\,.
From \eqref{eq:numestimdf} we deduce
$$
\vec{d}(\ker(\Delta_{f,\overline{\Omega},h}^{(p)}),F_{[0,e^{-\frac{\varepsilon}{h}}],h}^{(p)})=\tilde{O}(e^{-\frac{\eta_{f}-3\varepsilon/2}{h}})\,,
$$
which yields the result.
\end{proof}
The result of Proposition~\ref{pr:comptageest} can be translated in terms of singular values of $d_{f,f^{-1}([a,b]),h}$\,.
Remember that $d_{f,f^{-1}([a,b]),h}$ and $d_{f,f^{-1}([a,b]),h}^{*}$
are endomorphisms of $F_{[0,C],[a,b],h}$ such that
\begin{eqnarray*}
  &&
\Delta_{f,f^{-1}([a,b]),h}\big|_{F_{[0,C],[a,b],h}}=\delta_{[0,C],[a,b],h}\delta^{*}_{[0,C],[a,b],h}+\delta_{[0,C],[a,b],h}^{*}\delta_{[0,C],[a,b],h}\\
\text{with}&&
\delta_{[0,C],[a,b],h}=d_{f,f^{-1}([a,b]),h}\big|_{F_{[0,C],[a,b],h}}\,.
\end{eqnarray*}
\begin{prop}
\label{pr:countsing}
Under Hypothesis~\ref{hyp:cN} and  with the notations of
Subsection~\ref{sec:bar code}\,, the number of $\tilde{o}(1)$ 
non zero singular values of
  $\delta_{[0,\tilde{o}(1)],[a,b],h}=d_{f,f^{-1}([a,b]),h}\big|_{F_{[0,\tilde{o}(1)],[a,b],h}}$ is $\sharp A_{c}(a,b)$ for $h>0$ small
  enough. More precisely ``$h>0$ small enough'' means 
$h\in]0,h_{\varepsilon}[$ for some $h_{\varepsilon}>0$  when $\tilde{o}(1)$ is replaced by
$e^{-\frac{\varepsilon}{h}}$\,, $\varepsilon\in ]0,\frac{\eta_{f}}{2}[$\,.
\end{prop}
\begin{proof}
Eigenvalues and singular values are counted with multiplicities.
The non zero singular values of
$\underline{\delta}=\delta_{[0,e^{-\frac{\varepsilon}{h}}],[a,b],h}$
are the square roots of the non zero eigenvalues of $\underline{\delta}^{*}\underline{\delta}$
and coincide with the non zero singular values of $\underline{\delta}\,\underline{\delta}^{*}$\,, i.e. the square roots of the non zero eigenvalues of
$\underline{\delta}\,\underline{\delta}^{*}$\,.  
By  Hodge decomposition, the number of non zero eigenvalues
  of
  $\Delta_{f,f^{-1}([a,b]),h}\big|_{F_{[0,e^{-\frac{\varepsilon}{h}}],[a,b],h}}=\underline{\delta}\,\underline{\delta}^{*}+\underline{\delta}^{*}\underline{\delta}$
  is twice the number of non zero singular values of
  $\underline{\delta}$\,. For
  $h\in ]0,h_{\varepsilon}[$\,, 
Proposition~\ref{pr:comptageest} gives
\begin{align*}
\dim F_{[0,e^{-\frac{\varepsilon}{h}}],[a,b],h}
&=\sharp
  {\cal X}(a,b)+\sharp {\cal Y}(a,b)+\sharp {\cal
  Z}(a,b)\\
&=2 \sharp A_{c}(a,b) + 
\dim(\ker(\Delta_{f,f^{-1}([a,b]),h}))\,,
\end{align*}
which ends the proof.
\end{proof}

\subsection{Rough exponential estimates}
\label{sec:roughexp}
The upper bound on the $\tilde{o}(1)$ eigenvalues of
$\Delta_{f,f^{-1}([a,b]),h}$ contained in
Proposition~\ref{pr:comptageest} can be completed by a rough lower
bound for the non zero ones.
\begin{prop}
\label{pr:roughmino}
Assume Hypothesis~\ref{hyp:cN} and denote $a<\tilde{c}_{1}\ldots<\tilde{c}_{N}<b$ the
``critical values'' of $f$ in $]a,b[$\,. 
There exist $r(h)>0$  satisfying
$e^{-2\frac{\max\left\{b-\tilde{c}_{1},\tilde{c}_{N}-a\right\}}{h}}=\tilde{O}(r(h))$ and
 $R(h)=\tilde{O}(e^{-2\frac{\eta_{f}}{h}})$ such that the
 $\tilde{o}(1)$ non zero eigenvalues $\lambda(h)$ of
 $\Delta_{f,f^{-1}([a,b]),h}$ all belong to $[r(h),R(h)]$ for $h\in
 ]0,h_{0}[$\,, $h_{0}>0$ small enough.
\end{prop}
\begin{proof}
  The upper bound
  $R(h)=\tilde{O}(e^{-2\frac{\eta_{f}}{h}})$ is
   given by Proposition~\ref{pr:comptageest}.\\
For the lower bound, it suffices to check that if $\lambda(h)\in
\sigma(\Delta_{f,f^{-1}([a,b]),h})$ satisfies $\lambda(h)\leq
e^{-2\frac{\max\left\{b-\tilde{c}_{1},\tilde{c}_{N}-a\right\}+c}{h}}$ for some
fixed $c\in ]0,\min\left\{\tilde{c}_{1}-a,b-\tilde{c}_{N}\right\}[$\,, then there exists $h_{c}>0$ such that
$\lambda(h)=0$ for all $h\in ]0,h_{c}[$\,. The proof follows
the same arguments as those of Step~1  in Subsection~\ref{sec:exp0}.\\
Let us proceed by contradiction and assume that there 
exists a  decreasing sequence $(h_{n})_{n\in\mathbb N}$ tending 
to $0$
such that, for every $n\in\mathbb N$\,, $\Delta_{f,f^{-1}([a,b]),h_{n}}$
admits an eigenvalue $\lambda(h_{n})$ in the interval $]0,  e^{-2\frac{\max\left\{b-\tilde{c}_{1},\tilde{c}_{N}-a\right\}+c}{h_{n}}} ]$\,.
Let then, for every $n\in\mathbb N$\,, $\omega_{n}\in D(\Delta_{f,f^{-1}([a,b]),h_{n}})$ satisfy
$\|\omega_{n}\|_{L^{2}}=1$ and
$\Delta_{f,f^{-1}([a,b]),h_{n}}\omega_{n}=\lambda(h_{n})\omega_{n}$\,.
From the Agmon
estimates of Proposition~\ref{pr:Agmon} (or
Hypothesis~\ref{hyp:AgmonLip}) with $U\subset
f^{-1}(\left\{\tilde{c}_{1},\ldots, \tilde{c}_{N}\right\})$\,,  we know
that
\begin{eqnarray*}
  &&\forall \delta>0\,, \exists h_{\delta}>0\,, \forall h_{n}\in
     ]0,h_{\delta}[\,,\quad
     \|e^{\frac{f-\tilde{c}_{1}}{h_{n}}}\omega_{n}\|_{^{L^{2}}(f_{\tilde{c}_{1}-\delta}^{b})}\geq
  \frac{e^{-\frac{\delta}{h_{n}}}}{ 2}\\
\text{while}
&&
\|d_{f,h_{n}}\omega_{n}\|_{L^{2}(f_{\tilde{c}_{1}-\delta}^{b})}^{2}
+
\|d_{f,h_{n}}^{*}\omega_{n}\|_{L^{2}(f_{\tilde{c}_{1}-\delta}^{b})}^{2}\leq 
 e^{-\frac{2(b-\tilde{c}_{1})+c}{h_{n}}}\,.
\end{eqnarray*}
By setting $\tilde{\omega}_{n}=e^{\frac{f-\tilde{c}_{1}}{h_{n}}}\chi
\omega_{n}$\,, with $\chi\in
\mathcal{C}^{\infty}(f^{-1}([a,b]);[0,1])$\,, $\chi\equiv 1$ in
$f_{\tilde{c}_{1}-\frac{c}{4}}^{b}$ and $\chi\equiv 0$ 
in $f_{a}^{\tilde{c}_{1}-c'}$ with $c'\in (\frac{c}{4},\frac{c}{2})$\,, we get,
for every $n\in\mathbb N$\,, 
$$
\left\{
  \begin{array}[c]{l}
    \tilde{\omega}_{n}\in D(\Delta_{0, f^{-1}([\tilde{c}_{1}-c',b]),1})\,,\\
   \|d_{0,f^{-1}([\tilde{c}_{1}-c',b]),1}\tilde{\omega}_{n}\|_{L^{2}}^{2}=\tilde{O}(e^{-\frac{c}{h_{n}}})\,\\
\liminf_{n\to +\infty}h_{n}\log \|\tilde{\omega}_{n}\|_{L^{2}}\geq 0\,.
  \end{array}
\right.
$$

\medskip
\figinit{cm}
\figpt 0:(6,-1)
\figpt 1000: (-7,-1)
\figpt 1001:(7,-1)
\figpt 1002:(-1,-1)
\figpt 1003:(1,-1)
\figpt  1:(-6,-1)
\figpt 2: (0,-1)
\figpt 3: (-4,-1)
\figpt 11:(-3.5,-1)
\figpt 111:(-3,-1)
\figpt 12:(-2.5,-0.5)
\figpt 13:(-2,0.5)
\figpt 14:(-1.75,1)
\figpt 114:(-1,1)
\figpt 115:(-0.5,1)
\figpt 15:(-1,1)
\figpt 16:(1,1)
\figpt 17:(6.5,1)
\figptmap 0:=0 /111/0.2,0;0,1/
\figptmap 1:=1 /111/0.2,0;0,1/
\figptmap 2:=2 /111/0.2,0;0,1/
\figptmap 11:=11 /111/0.2,0;0,1/
\figptmap 111:=111 /111/0.2,0;0,1/
\figptmap 12:=12 /111/0.2,0;0,1/
\figptmap 13:=13 /111/0.2,0;0,1/
\figptmap 14:=14 /111/0.2,0;0,1/
\figptmap 114:=114 /111/0.2,0;0,1/
\figptmap 115:=115 /111/0.2,0;0,1/
\figpt 200:(-6.5,-1)
\figpt 201:(1.5,-1)
\figpt 202:(6.5,-1)
\figdrawbegin{}
\figset(width=1.3)
\figdrawline[1000,1002]
\figdrawline[1003,1001]
\figset(width=3)
\figdrawBezier 1[1000,1,3,11]
\figdrawBezier 1[11,111,12,13]
\figdrawBezier 1[13,14,114,115]
\figdrawline[115,15]
\figdrawline[16,17]
\figset(dash=4)
\figdrawline[15,16]
\figset(width=1)
\figdrawline[1002,1003]
\figdrawend{}
\figvisu{\figBoxB}{\textbf{Figure 8:} The cut-off $\chi$ in the
interval $[a,b]$\,.}{
\figset write(mark=$+$)
\figwrites 0:$c_{1}$(0.3)
\figwrites 1:$c_{1}-\frac{c}{2}$(0.3)
\figwrites 2:$c_{1}-\frac{c}{4}$(0.3)
\figwrites 200:$a$(0.3)
\figwrites 201:$c_{N}$(0.3)
\figwrites 202:$b$(0.3)
\figset write(mark={})
\figwriten 115:$\chi$(0.3)
}
\centerline{\box\figBoxB}

\medskip

\noindent
Besides, the Agmon estimates of
Proposition~\ref{pr:Agmon} (or
Hypothesis~\ref{hyp:AgmonLip}) with $U\subset
f^{-1}(\left\{\tilde{c}_{1},\ldots, \tilde{c}_{N}\right\})$  also imply
$$
\limsup_{n\to +\infty}h_{n}\log \|\tilde{\omega}_{n}\|_{L^{2}}\leq \tilde c_{N}-\tilde c_{1}\,.
$$
Hence, by extracting, we can assume that there exists $\ell\in[0, 2(\tilde c_{N}-\tilde c_{1})]$ 
such that
$$
\lim_{n\to +\infty}h_{n}\log \|\tilde{\omega}_{n}\|_{L^{2}}= \frac\ell2\,.
$$
The normalized form
$u_{n}=\frac{\tilde{\omega}_{n}}{\|\tilde{\omega}_{n}\|_{L^{2}}}$
thus
belongs to $D(\Delta_{0,f^{-1}([\tilde{c}_{1}-c',b]),1})$ and
$\|du_{n}\|_{L^{2}}^{2}=\tilde{O}(e^{-\frac{c+\ell}{h_{n}}})$\,.
By Hodge decomposition (see Step~1 in
Subsection~\ref{sec:exp0} for details), this implies that 
 $\eta_{n}$ belongs to 
$\ker(d_{0,f^{-1}([\tilde{c}_{1}-c',b]),1})$ and
$$
\|u_{n}-\eta_{n}\|_{L^{2}(f_{\tilde{c}_{1}-c'}^{b})}=\tilde{O}(e^{-\frac{c+\ell}{2h_{n}}})\,.
$$
Moreover, extending $\eta_{n}$ by $0$ in $f_{a}^{\tilde{c}_{1}-c'}$ gives
$\eta_{n}\in D(d_{0,f^{-1}([a,b]),1})$ and therefore
$e^{-\frac{f-\tilde{c}_{1}}{h_{n}}}\eta_{n}\in \ker(d_{f,f^{-1}([a,b]),h_{n}})$ with
$$
\|\chi \omega_{n}-\|\tilde{\omega}_{n}\|_{L^{2}}e^{-\frac{f-\tilde{c}_{1}}{h_{n}}}\eta_{n}\|_{L^{2}(f_{a}^{b})}=\tilde{O}(e^{-\frac{c/2-c'}{h_{n}}})\,.
$$ 
With $\|\omega_{n}-\chi
\omega_{n}\|_{L^{2}}=\tilde{O}(e^{-\frac{c}{4h_{n}}})$ and
$c''=\min\left\{c/2-c',c/4\right\}$\,, we deduce
$$
\mathrm{dist}_{L^{2}}(\omega_{h},
\ker(d_{f,f^{-1}([a,b]),h_{n}}))=\tilde{O}(e^{-\frac{c''}{h_{n}}})\underset{h\to
0}{\to} 0\,.
$$
 By duality, starting from $\|d_{f,h}^{*}\omega_{h_{n}}\|_{L^{2}}^{2}\leq
 e^{-\frac{2(\tilde{c}_{N}-a)+c}{h_{n}}}$ and extracting again, we also get
$$
\lim_{h\to 0}\mathrm{dist}_{L^{2}}(\omega_{h_{n}},\ker(d_{f,f^{-1}([a,b],h_{n})}^{*}))=0
$$
and Hodge decomposition implies $\lambda(h_{n})=0$ for
$n$ large enough (see indeed the end of Step~1 in
Subsection~\ref{sec:exp0}), which leads to a contradiction and achieves
the proof of Proposition~\ref{pr:roughmino}. 
\end{proof}
\begin{remark}
  The lower bound for the non zero eigenvalues is not optimal at this
  level. Actually,
  generalizing Step~3 of Subsection~\ref{sec:exp0} requires the
  propagation of exponential decay through ``critical values'', which  
is not true in general. This will be refined into
$e^{-2\frac{\tilde{c}_{N}-\tilde{c}_{1}}{h}}=\tilde{O}(r(h))$
at the end, when global
quasimodes for $d_{f,f^{-1}([a,b],h)}$
 will have been constructed by induction on $N$\,. Like e.g. in
 \cite{HKN,HeNi,Lep3,LNV}, we follow the strategy
 which consists in studying carefully the singular values of
 $d_{f,f^{-1}([a,b]),h}$\,, which brings more flexibility than studying
 the tricky problem of interacting wells for
 $\Delta_{f,f^{-1}([a,b]),h}$ in the spirit of
 \cite{HeSj2,HeSj3}.
\end{remark}
\begin{prop}
\label{pr:projAgm}
Assume Hypothesis~\ref{hyp:cN}, let
$a<\tilde{c}_{1},\ldots<\tilde{c}_{N}<b$ be the ``critical values'' of
$f$ in $]a,b[$ and let $R(h)$ be the function of
$h\in]0,h_{0}[$ given by Proposition~\ref{pr:roughmino} such that
$\sigma(\Delta_{f,f^{-1}([a,b]),h})\cap [0,\tilde{o}(1)]\subset[0,R(h)]$\,.
The projection
$\Pi_{[0,R(h)],[a,b],h}=1_{[0,R(h)]}(\Delta_{f,f^{-1}([a,b]),h})$
satisfies
$$
\Pi_{[0,R(h)],[a,b],h}=\tilde{O}(e^{-\frac{|f(x)-f(y)|}{h}})
$$
in the sense of Definition~\ref{de:Otkernel}.
\end{prop}
\begin{proof}
By Proposition~\ref{pr:roughmino}, we know that
$R(h)=\tilde{O}(e^{-\frac{2\eta_{f}}{h}})$\,. Set $\tilde{c}_{0}=a$ and
$\tilde{c}_{N+1}=b$ and take any $\varepsilon_{0}\in ]0,
\frac{\eta_{f}}{8}[$\,, where $\eta_{f}$ is defined in Hypothesis~\ref{hyp:cN}. Here
the first assumption of Proposition~\ref{pr:Nwell} is obviously
satisfied:
$$ 
]a,b[\cap
\left\{c_{1},\ldots,c_{N_{f}}\right\}=\left\{\tilde{c}_{1},\ldots,\tilde{c}_{N}\right\}\subset
\sqcup_{n=1}^{N}]\tilde{c}_{n}-\frac{\varepsilon_{0}}{16},\tilde{c}_{n}+\frac{\varepsilon_{0}}{16}[\,.
$$
For
$\Delta_{n}=\Delta_{f,f^{-1}([\tilde{c}_{n-1}+(1-\delta_{n,1})\varepsilon_{0},\tilde{c}_{n+1}-(1-\delta_{n,N})\varepsilon_{0}])}$\,,
$n\in\left\{1,\ldots, N\right\}$\,, we know moreover that
$$
\sigma(\Delta_{n})\cap [0,e^{-\frac{\varepsilon_{0}}{h}}]\subset \left\{0\right\}\subset[0,e^{-\frac{4\varepsilon_{0}}{h}}]\,,
$$
owing to Proposition~\ref{pr:exp0} because we are in the case
$[\tilde{c}_{n-1}+(1-\delta_{n,1})\varepsilon_{0},\tilde{c}_{n+1}-(1-\delta_{n,N})\varepsilon_{0}]\cap\left\{c_{1},\ldots,c_{N_{f}}\right\}=\left\{\tilde{c}_{n}\right\}$\,.\\
Then Proposition~\ref{pr:Nwell} says: for some $N\in\mathbb N^{*}$\,,
$$
(\Delta_{f,f^{-1}([a,b]),h}-z)^{-1}(x,y)=\tilde{O}(e^{-\frac{|f(x)-f(y)|}{h}+\frac{3N\varepsilon_{0}}{h}})
$$
uniformly w.r.t $z$\,, $|z|=e^{-\frac{2\varepsilon_{0}}{h}}$\,. But
our choice of $\varepsilon_{0}$\,, $\varepsilon_{0}>0$ and
$4\varepsilon_{0}\leq \frac{\eta_{f}}{2}$\,, and
$$
\sigma(\Delta_{f,f^{-1}([a,b]),h})\cap
[0,e^{-\frac{\varepsilon_{0}}{h}}]\subset [0,R(h)]
\subset[0,e^{-\frac{\eta_{f}}{h}}]\subset[0,e^{-\frac{4\varepsilon_{0}}{h}}]
$$
for $h\in]0,h_{0}[$\,, $h_{0}$ small enough, imply
$$
\Pi_{[0,R(h)],[a,b],h}=\frac{1}{2i\pi}\int_{|z|=e^{-\frac{2\varepsilon_{0}}{h}}}(z-\Delta_{f,f^{-1}([a,b]),h})^{-1}~dz\,.
$$
This proves
$$
\Pi_{[0,R(h)],[a,b],h}(x,y)=\tilde{O}(e^{-\frac{|f(x)-f(y)|}{h}+\frac{3N\varepsilon_{0}}{h}})\,,
$$
and we conclude by choosing $\varepsilon_{0}>0$ arbitrarily small.
\end{proof}

\section{Singular values}
\label{sec:singval}

Singular values of compact operators are much more flexible than
eigenvalues because they allow to work with two different orthonormal
bases instead of one. Ky~Fan inequalities recalled below 
provide uniform
multiplicative errors for all the singular values after perturbing
the orthonormal bases or moving the initial and final spaces. 
We recall those facts in a convenient way and complete those results by
some refined analysis of additive error terms. This is a better
rewriting of  techniques already used e.g. in \cite{HKN,HeNi,Lep3,LNV}\\
The singular values of a compact operator $B:E\mapsto F$\,, $E$ and
$F$ Hilbert spaces, are the square roots of the eigenvalues of
$B^{*}B$ (and $BB^{*}$) and they are labelled in the decreasing order
 $\mu_{1}(B)=\|B\|\geq\ldots\geq
\mu_{\ell}(B)\geq \mu_{\ell+1}(B)\ldots$ with  $\lim_{\ell\to
  \infty}\mu_{\ell}(B)=0$ after possibly completing the sequence by a
sequence of $0$'s\,. They satisfy $\mu_{\ell}(B)=\mu_{\ell}(B^{*})$\,.
 With this order, the min-max principle 
becomes a max-min principle applied to $B^{*}B$ and gives:
\begin{equation}
\label{eq:maxmin}
\mu_{\ell}(B)=\min_{\dim V=\ell-1}\max_{u\in V^{\perp}\setminus\left\{0\right\}}\frac{\|Bu\|}{\|u\|}\,.
\end{equation}
Note also that the definition also provides the existence of two
Hilbert bases $(\varphi_{j})_{j\in \mathcal{J}}$\,,
$\mathcal{J}\supset \mathcal{J}_{1}=\left\{\ell\in \nz\setminus\{0\}\,, \mu_{\ell}(B)>0\right\}$\,,
 of $E$\,, and
$(\psi_{k})_{k\in \mathcal{K}}$ of $F$\,, and a one-to-one mapping $j\in
\mathcal{J}_{1}\to k(j)\in \mathcal{K}$ such that 
\begin{eqnarray*}
B\varphi_{\ell}=\mu_{\ell}(B)\psi_{k(\ell)}\quad\text{and then}
\quad
\mu_{\ell}(B)=\|B\varphi_{\ell}\|=\langle \psi_{k(\ell)}\,,
     B\varphi_{\ell}\rangle && \text{if}~\ell\in \mathcal{J}_{1}\\
      B\varphi_{j}=0 &&\text{if}~j\in \mathcal{J}\setminus\mathcal{J}_{1}\,.
\end{eqnarray*}
When $E,F,G$ are three Hilbert spaces and
$A:E\to F$\,, $B:F\to F$\,, and $C:F\to G$\,, the singular values of $B$ also satisfy
$$
\forall \ell\in \nz\setminus\{0\}\,, \quad \mu_{\ell}(CBA)\leq \|C\|\mu_{\ell}(B)\|A\|\,.
$$
In order to handle accumulated multiplicative errors, 
 it is convenient to use the
function
\begin{equation}
\label{eq:deftau}
\tau: \mathop{\sqcup}_{n=1}^{\infty}[0,1[^{n}\to ]0,+\infty[\quad,\quad \tau(\varepsilon_{1},\ldots,\varepsilon_{n})=\prod_{k=1}^{n}\frac{1+\varepsilon_{k}}{1-\varepsilon_{k}}\,.
\end{equation}
In particular we have the implications {\bf i)} $\Rightarrow$ {\bf ii)} $\Rightarrow$ {\bf iii)} for 
\begin{eqnarray*}
 \textbf{i)}&& \max(\|
  CC^{*}-\Id_{G}\|,\|C^{*}C-\Id_{F}\|)\leq \varepsilon_{1}<1\\
&&\text{and}\quad \max(\|AA^{*}-\Id_{F}\|,\|A^{*}A-\Id_{E}\|)\leq
             \varepsilon_{2}<1\,;\\
\textbf{ii)}&& \max(\|C\|,\|C^{-1}\|) \leq
             \tau(\varepsilon_{1})^{1/2}\quad,\quad
\max(\|A\|,\|A^{-1}\|)\leq\tau(\varepsilon_{2})^{1/2}
             \,;\\
\textbf{iii)}&&\forall j\in \nz \setminus\left\{0\right\}\,, \tau(\varepsilon_{1},\varepsilon_{2})^{-1/2}\mu_{j}(B)\leq \mu_{j}(CBA)\leq \tau(\varepsilon_{1},\varepsilon_{2})^{1/2}\mu_{j}(B)\,.
\end{eqnarray*}
The first implication is a consequence of the following operator inequalities
$$
(1-\varepsilon\leq
\left|A\right|^{2}=A^{*}A\leq 1+\varepsilon)\Rightarrow 
\left(
\tau(\varepsilon)^{-1/2}\leq (1-\varepsilon)^{1/2}
\leq \left|A\right|\leq (1+  \varepsilon)^{1/2}\leq
  \tau(\varepsilon)^{1/2}\right)\,.
$$
\begin{definition}
\label{de:epsorth}
Let $\mathcal{H},\mathcal{H}'$ be two Hilbert spaces and let
$\varepsilon\in [0,1[$\,.\\
An operator $A:\mathcal{H}\to \mathcal{H}'$ will be said
$\varepsilon$-unitary  if it
satisfies the condition
$$
\max(\|A^{*}A-\Id_{\mathcal{H}}\|, \|AA^{*}-\Id_{\mathcal{H}'}\|)\leq \varepsilon\,,
$$
used in \textbf{i)} just above.\\
  A family of vectors  $(v_{j})_{j\in \mathcal{J}}$  is an
  $\varepsilon$-orthonormal basis of $\mathcal{H}$
  if
  \begin{itemize}
  \item it is total in $\mathcal{H}$\,,
    $\overline{\Vect(v_{j},j\in \mathcal{J})}=\mathcal{H}$\,,
\item 
$
\|(\langle v_{j}\,,\, v_{k}\rangle)_{j,k\in
  \mathcal{J}}-\Id_{\ell^{2}(\mathcal{J})}\|_{\mathcal{L}(\ell^{2}(\mathcal{J}))}\leq \varepsilon\,.$
\end{itemize}
Two closed subspaces $\mathcal{H}_{1},\mathcal{H}_{2}$ of
$\mathcal{H}$ provide an
$\varepsilon$-orthogonal decomposition
 of $\mathcal{H}$ if
$\mathcal{H}=\mathcal{H}_{1}\oplus \mathcal{H}_{2}$ 
and $\|\Pi_{\mathcal{H}_{1}}\Pi_{\mathcal{H}_{2}}\| \leq \varepsilon$\,.
\end{definition}
Before we review applications to singular values, notice the following
properties.
\begin{lem}
\label{le:image}
Let $\mathcal{H}$\,, $\mathcal{H}'$ be Hilbert spaces and let
$\varepsilon\in [0,1[$\,.
\begin{description}
\item[a)] For an operator $A: \mathcal{H}\to \mathcal{H}'$\,, the
  condition $\|A^{*}A-\Id_{\mathcal{H}}\|\leq \varepsilon$ is
  satisfied iff $\left|A\right|:\mathcal{H}\to \mathcal{H}$ is
  $\varepsilon$-unitary and iff $\Id_{\mathcal{H}}:
  (\mathcal{H},\langle~,~\rangle)\to (\mathcal{H},
  \langle~,|A|^{2}~\rangle)$ is $\varepsilon$-unitary.
\item[b)] An operator $A:\mathcal{H}\to \mathcal{H}'$ is
  $\varepsilon$-unitary iff
$$
\|A^{*}A-\Id_{\mathcal{H}}\|\leq \varepsilon\quad\text{and}\quad 
\overline{\textrm{Ran}~A}=\mathcal{H}'\,.
$$
\item[c)] A family $(v_{j})_{j\in \mathcal{J}}$ is an
  $\varepsilon$-orthonormal basis of $\mathcal{H}'$ iff the linear map
  $A:\ell^{2}(\mathcal{J})\to \mathcal{H}'$ given by $A((a_{j})_{j\in
    \mathcal{J}})=\sum_{j\in \mathcal{J}}a_{j}v_{j}$ is
  $\varepsilon$-unitary.
\item[d)] If the decomposition $\mathcal{H}'=\mathcal{H}_{1}\oplus
  \mathcal{H}_{2}$ is $\varepsilon$-orthogonal and
  $(\varphi_{j'})_{j'\in \mathcal{J}'}$ and $(\varphi_{j''})_{j''\in
    \mathcal{J}''}$ are orthonormal bases of $\mathcal{H}_{1}$ and
  $\mathcal{H}_{2}$ respectively, then $(\varphi_{j})_{j\in
    \mathcal{J}'\cup \mathcal{J}''}$ is an $\varepsilon$-orthonormal
  basis of $\mathcal{H}'$\,. Additionally, the identity map
    induces an $\varepsilon$-unitary map from
    $ \mathcal H'=\mathcal{H}_{1}\mathop{ \oplus}\limits^{\perp}\mathcal{H}_{2}$ to $\mathcal{H}'=\mathcal{H}_{1}\oplus \mathcal{H}_{2}$\,,
    where the first space is endowed with the
    scalar product $\langle\,,\, \rangle_{\mathcal{H}_{1}\mathop{ \oplus}\limits^\perp\mathcal{H}_{2}}$ making $(\varphi_{j})_{j\in \mathcal{J}'\cup
      \mathcal{J}''}$  orthonormal, i.e. defined by
      $$
  \forall \,u_{1}, v_{1}\,\in\, \mathcal{H}_{1}\,,\ 
  \forall \,u_{2}, v_{2}\,\in\, \mathcal{H}_{2}\,,\ \    \langle u_{1}+u_{2}\,,\,v_{1}+v_{2} \rangle_{\mathcal{H}_{1}\mathop{ \oplus}\limits^\perp\mathcal{H}_{2}}\ :=\ 
\langle u_{1}\,,\,v_{1} \rangle+  \langle u_{2}\,,\,v_{2} \rangle
  \,.
      $$
\end{description}
\end{lem}
\begin{proof}
  \textbf{a)} The first statement is a consequence of $|A|^{*}=|A|$
  and $|A|^{2}=A^{*}A$\,. The second one is deduced from $\Id^{*}=|A|^{2}$ when
 the identity operator maps $\mathcal{H}$ with the scalar product
 $\langle u,v\rangle$ to itself with the scalar product $\langle u\,,
 |A|^{2}v\rangle$\,.\\
\noindent\textbf{b)} It suffices to notice that the condition
$\|A^{*}A-\Id_{\mathcal{H}}\|\leq \varepsilon$ implies
$$
\forall u\in \mathcal{H}\,,\quad \sqrt{1-\varepsilon}\|u\|\leq
\|Au\|\leq \sqrt{1+\varepsilon}\|u\|\,.
$$
Thus $A$ is one-to-one with a closed range which has to be
$\mathcal{H}'$ by the second assumption and $A$\,, $A^{*}$\,, 
and $AA^{*}$ are invertible. Hence the spectrum of
$AA^{*}$ coincides with the spectrum of $A^{*}A$ 
by
$A^{*}(AA^{*}-\lambda\Id_{\mathcal{H}'})=(A^{*}A-\lambda\Id_{\mathcal{H}})A^{*}$
for $\lambda\in \cz$\,. The spectral theorem yields
$\|AA^{*}-\Id_{\mathcal{H}'}\|\leq \varepsilon$\,.\\
\noindent\textbf{c)} is a particular case of \textbf{b)} if we notice
that $\|A^{*}A-\Id_{\mathcal{H}}\|=\|(\langle
v_{j},v_{k}\rangle)_{j,k\in
  \mathcal{J}}-\Id_{\ell^{2}(\mathcal{J})}\|$ with
$\mathcal{H}=\ell^{2}(\mathcal{J})$\,,
 while the condition
$\overline{\textrm{Ran}~A}=\mathcal{H}'$ becomes equivalent to the totality of the family
$(v_{j})_{j\in \mathcal{J}}$\,.\\
\noindent\textbf{d)} The  family $(\varphi_{j})_{j\in
    \mathcal{J}'\cup \mathcal{J}''}$ is clearly total in $\mathcal{H}'$ and, 
    defining the map $A:\mathcal{H}\to \mathcal{H}'$ with 
    $\mathcal{H}=\ell^{2}(\mathcal{J})$ as
 in    \textbf{c)}, we get
 $$
A^{*}A- \Id_{\ell^{2}(\mathcal{J})}=\begin{pmatrix} 0 & B \\
B^{*} & 0\end{pmatrix}\quad\text{with}\quad B=(\langle \varphi_{k}\,,\, \varphi_{j}\rangle)_{j \in \mathcal J'',k\in
  \mathcal{J}'}\,.
 $$
To prove that $(\varphi_{j})_{j\in
    \mathcal{J}'\cup \mathcal{J}''}$ is an $\varepsilon$-orthonormal
  basis of $\mathcal{H}'$\,, it is then enough to prove that  
$\|B\|_{\mathcal{L}(\ell^{2}(\mathcal{J}''),\ell^{2}(\mathcal{J}') )}\leq \varepsilon$\,,
which follows from the observation that $B$ is unitarily equivalent to
$\Pi_{\mathcal{H}_{1}}|_{\mathcal{H}_{2}}:\mathcal{H}_{2}\to \mathcal{H}_{1}$\,.\\  
For the
  last statement, it suffices to note that
  the
  mapping $$u\in \big(\,\mathcal{H}_{1}\mathop{ \oplus}\limits^\perp\mathcal{H}_{2}\,,\,
 \langle\,,\, \rangle_{\mathcal{H}_{1}\mathop{ \oplus}\limits^\perp\mathcal{H}_{2}} \,\big)
   \ \longmapsto\  (\langle \varphi_{j}\,,\, u\rangle_{\mathcal{H}_{1}\mathop{ \oplus}\limits^\perp\mathcal{H}_{2}})_{j\in
    \mathcal{J}'\cup \mathcal{J}''}
 \in\ell^{2}(\mathcal{J'}\cup \mathcal{J''})$$ is unitary   and to apply \textbf{c)}.
\end{proof}

Below are consequences of those notions on singular values.

\begin{prop}
\label{pr:projepsort}
Let $E,F,G$ be three closed subspaces of  a Hilbert space 
$\mathcal{H}$ and assume
$\vec{d}(E,F)+\vec{d}(F,E)=\varepsilon_{1}<1$ and
$\vec{d}(F,G)+\vec{d}(G,F)=\varepsilon_{2}<1$\,.
Let $B:F\to F$ be a bounded operator and let $\Pi_{F}, \Pi_{G}$ 
be the
orthogonal projections on $F$ and $G$\,. The operator
$\tilde{B}=\Pi_{G}B\Pi_{F}\big|_{E}:E\to G$ is compact iff $B$ is compact  and
in this case:
$$
\forall \ell\in \nz \setminus\left\{0\right\}\,,\quad
\tau(\varepsilon_{1}^{2},\varepsilon_{2}^{ 2})^{-1/2}\,\mu_{\ell}(\tilde{B})\ 
\leq\ \mu_{\ell}(B)\ \leq\  
\mu_{\ell}(\tilde{B})\,\tau(\varepsilon_{1}^{
  2},\varepsilon_{2}^{ 2})^{1/2}\,.
$$
\end{prop}
\begin{proof}
Call $A_{FE}=\Pi_{F}\Pi_{E}+(1-\Pi_{F})(1-\Pi_{E})$\,, with
$1=\Id_{\mathcal{H}}$\,,  and compute
$$
A_{FE}^{*}A_{FE}-1=\Pi_{E}\Pi_{F}+\Pi_{F}\Pi_{E}-\Pi_{E}-\Pi_{F}\,.
$$
We deduce that for all $u\in \mathcal{H}$\,,
\begin{align*}
\langle u\,,\, (A_{FE}^{*}A_{FE}-1)u\rangle
&=2\,\Real \langle \Pi_{E}u\,,\, \Pi_{F}u\rangle -\|\Pi_{E}u\|^{2}-\|\Pi_{F}u\|^{2}
\\
&=-\|(\Pi_{E}-\Pi_{F})u\|^{2}
\\
&\geq
  -2\|(\Pi_{E}-\Pi_{F})\Pi_{E}u\|^{2}-2\|(\Pi_{E}-\Pi_{F})(1-\Pi_{E})u\|^{2}\\
&\geq
  -2\|(\Pi_{E}-\Pi_{F}\Pi_{E})u\|^{2}-2\|(\Pi_{F}-\Pi_{F}\Pi_{E})u\|^{2}\\
&\geq 
-2 \big(\,\vec{d}(E,F)^{2}+\vec{d}(F,E)^{2}\,\big)\|u\|^{2}\,.
\end{align*}
Since $0\leq \varepsilon_{1}<1$\,, we know that
$\vec{d}(E,F)=\vec{d}(F,E)=\frac{\varepsilon_{1}}{2}$
(see indeed the lines following Definition~\ref{de:dEF}) and we have thus proved
the operator inequalities
$$
0\leq (\Id_{\mathcal{H}}-A_{FE}^{*}A_{FE})\leq
\varepsilon_{1}^{2}\Id_{\mathcal{H}}\,.$$
Owing to the spectral theorem, it follows
$$
\|A_{FE}^{*}A_{FE}-\Id_{\mathcal{H}}\|\leq \varepsilon_{1}^{2}\,,
$$
and by symmetry, since $A_{FE}^{*}=A_{EF}$\,, we also get
$\|A_{FE}A_{FE}^{*}-\Id\|\leq \varepsilon_{1}^{2}$\,. The operator
$A_{FE}$ is thus $\varepsilon_{1}^{2}$-unitary, and similarly $A_{GF}$ is
$\varepsilon_{2}^{2}$-unitary.\\
Finally, $\tilde{B}=\Pi_{G}B\Pi_{F}\big|_{E}:E\to G$ is nothing but the
nonzero diagonal block of
$$A_{GF}\,B\,A_{FE}:\mathcal{H}=E\mathop{\oplus}\limits^{\perp}E^{\perp}\ \longrightarrow\ 
\mathcal{H}=G\mathop{\oplus}\limits^{\perp}G^{\perp}\,.$$ It is thus compact if and
only if $B$ is compact. Moreover, up to some additional irrelevant zeros, the
singular values of $\tilde{B}$ are the ones of $A_{GF}BA_{FE}$ and the
result follows from the general statement \textbf{i)} $\Rightarrow$ \textbf{iii)} above.
\end{proof}

\begin{prop}
\label{pr:epsorth}
Let $E,F$ be two Hilbert spaces, $B:E\to F$ be a bounded operator
and let $\varepsilon_{1},\varepsilon_{2}\in ]0,1[$\,.
\begin{description}
\item[a)] 
When $(\varphi_{j})_{j\in \mathcal{J}}$
 is an $\varepsilon_{1}$-orthonormal basis in $E$ and
$(\psi_{k})_{k\in\mathcal{K}}$ is an 
$\varepsilon_{2}$-orthonormal in basis $F$\,, let 
 $\tilde{B}:\ell^{2}(\mathcal{J})\to
\ell^{2}(\mathcal{K})$ be defined by $\tilde{B}\delta_{j}=\sum_{k\in \mathcal{K}}\langle
\psi_{k}\,, B\varphi_{j}\rangle\delta_{k}$\,. Then $\tilde{B}$ is 
compact  iff $B$ is compact, and in this case their singular
values satisfy
\begin{equation}
  \label{eq:ineqepsorth}
\forall \ell\in \nz\setminus\left\{0\right\}\,,\quad
\tau(\varepsilon_{1},\varepsilon_{2})^{-1/2}\mu_{\ell}(\tilde{B})\leq
\mu_{\ell}(B)\leq \mu_{\ell}(\tilde{B})
\tau(\varepsilon_{1},\varepsilon_{2})^{1/2}\,.
\end{equation}
\item[b)] Assume that $E=E'\oplus E''$ is an
  $\varepsilon_{1}$-orthogonal decomposition and $F=F'\oplus F''$
  is an $\varepsilon_{2}$-orthogonal decomposition such that
  $BE'\subset F'$ and $BE''\subset F''$\,, then  the relation
  \eqref{eq:ineqepsorth} holds with  $\ds\tilde{B}=\Pi_{F'}B\big|_{E'}\mathop{\oplus}^{\perp}\Pi_{F''}B\big|_{E''}:
E'\mathop{\oplus}^{\perp}E''\to F'\mathop{\oplus}^{\perp}F''$\,.
\item[c)] Assume that $B$ is compact and that $E=E'\oplus E''$ is an
  $\varepsilon_{1}$-orthogonal decomposition, and set
  $F'=\overline{BE'}$\,, $F''=(F')^{\perp}$\,. Assume moreover that
  $$\nu\ =\ \inf\Big(\,\left\{\mu_{\ell}(B\big|_{E'}), \ell\in \nz\setminus\left\{0\right\}\right\}\cap ]0,+\infty[\,\Big)
  \ \geq \ 
  \frac{\|B\big|_{E''}\|}{(1-\varepsilon_{1})^{\frac12}\varepsilon_{2}}\,.$$ 
  Then, the operator
$\ds\tilde{B}=B\big|_{E'}\mathop{\oplus}^{\perp}\Pi_{F''}B\big|_{E''}:
E'\mathop{\oplus}^{\perp}E''\to F'\mathop{\oplus}^{\perp}F''$ satisfies
$$
\forall \ell\in \nz\setminus\left\{0\right\}\,,\quad
\tau(\varepsilon_{1},\varepsilon_{2})^{-1}\mu_{\ell}(\tilde{B})
\leq \mu_{\ell}(B)\leq
\mu_{\ell}(\tilde{B})
\tau(\varepsilon_{1},\varepsilon_{2})
\,.
$$
\end{description}
\end{prop}
\begin{proof}
\textbf{a)}
This item simply follows 
from the general statement \textbf{i)} $\Rightarrow$ \textbf{iii)}
above and
from the relation 
$\tilde{B}=\Psi_{F}^{*}B\Phi_{E}$\,, where
$\Phi_{E}:\ell^{2}(\mathcal{J})\to E$ and
$\Psi_{F}:\ell^{2}(\mathcal{K})\to F$ are  defined by
$\Phi_{E}((u_{j})_{j\in \mathcal{J}})=\sum_{j\in
  \mathcal{J}}u_{j}\varphi_{j}$ 
  and 
  $\Psi_{F}((v_{k})_{k\in \mathcal{K}})=\sum_{k\in
  \mathcal{K}}v_{k}\psi_{k}$\,, and are thus respectively
  $\varepsilon_{1}$- and $\varepsilon_{2}$-unitary
  according to item {\bf c)} in Lemma~\ref{le:image}.\\
\textbf{b)} Let $(\varphi_{j})_{j\in \mathcal{J}'}$ and
$(\varphi_{j})_{j\in \mathcal{J}''}$ be two Hilbert bases of $E'$ and
$E''$\,, so that $(\varphi_{j})_{j\in \mathcal{J}'\cup \mathcal{J}''}$ is
an $\varepsilon_{1}$-orthonormal basis of $E$ 
according to item {\bf d)} in Lemma~\ref{le:image}. 
An
$\varepsilon_{2}$-orthonormal basis $(\psi_{k})_{k\in
  \mathcal{K}'\cup\mathcal{K}''}$ of $F$ is constructed  in a similar way. 
 It also follows from item {\bf d)} in Lemma~\ref{le:image}
that  the identity 
$\Id_{E}:E=E'\mathop{\oplus}\limits^{\perp}E''\to E=E'\mathop{\oplus}E''$
is $\varepsilon_{1}$-unitary and, similarly, $\Id_{F}$
is $\varepsilon_{2}$-unitary.
We conclude by  applying
the general statement \textbf{i)} $\Rightarrow$ \textbf{iii)}
above to the relation
$\tilde{B}=\Id_{F}^{*}\,B\,\Id_{E}$\,.\\
\textbf{c)} If $\|B\big|_{E''}\|=0$ there is nothing to do. Actually
this is a particular case of \textbf{b)} with
$BE''=\left\{0\right\}\subset F''$\,,\, $\varepsilon_{2}=0$ and of
course $\tau^{1/2}\leq \tau$\,. If
$\|B\big|_{E''}\|>0$\,, then there exists $\ell_{1}\in
\nz\setminus\left\{0\right\}$ such that
$\nu=\mu_{\ell_{1}}(B\big|_{E'})$ and
$\mathrm{rank}(B\big|_{E'})=\ell_{1}$\,. In particular, we can find two
Hilbert bases $(\varphi_{j})_{j\in \mathcal{J}'}$ of $E'$ and
$(\psi_{k})_{k\in \mathcal{K}}$ of $F$ such that
$\mathcal{J}'\cap \mathcal{K}\supset \left\{1,\ldots,
  \ell_{1}\right\}$ and
$$
\forall j\in \left\{1,\ldots,\ell_{1}\right\}\,,\quad B\varphi_{j}=\mu_{j}(B\big|_{E'})\psi_{j}\,.
$$
Set $F'=\Vect(\psi_{j},\,, j\in \left\{1,\ldots,
  \ell_{1}\right\})=\Ran B\big|_{E'}$ and $F''=(F')^{\perp}$\,, and
introduce the map $R:E\to E$ defined by 
\begin{eqnarray*}
  && R\big|_{E'}=0\,,\\
&& \forall u\in E'',\quad Ru=\sum_{j=1}^{\ell_{1}}\frac{\langle
   \psi_{j}\,, Bu\rangle}{\mu_{j}(B\big|_{E'})}\varphi_{j}\,.
\end{eqnarray*}
The norm of $R$ is not greater than $\varepsilon_{2}$ since
 for every $u=u'+u''\in E=E'\oplus E''$\,,
$$
\|Ru\|^{2}=\|Ru''\|^{2}=\sum_{j=1}^{\ell_{1}}\frac{\left|\langle \psi_{j}\,,\,
    Bu''\rangle\right|^{2}}{\mu^{2}_{j}(B\big|_{E'})}\leq
\frac{\|B\big|_{E''}\|^{2}}{\mu^{2}_{\ell_{1}}(B\big|_{E'})}\|u''\|^{2}\leq
(1-\varepsilon_{1})\varepsilon_{2}^{2}\|u''\|^{2}\leq \varepsilon^{2}_{2}\|u\|^{2} \,,
$$
where the last inequality follows from the last statement of Lemma~\ref{le:image}.
We deduce
\begin{eqnarray*}
  &&\|\Id_{E}-R\|\leq 1+\varepsilon_{2}\leq
\tau(\varepsilon_{2})\,,\\
&&\|(\Id_{E}-R)^{-1}\|\leq (1-\varepsilon_{2})^{-1}\leq
     \tau(\varepsilon_{2})\,,
    \end{eqnarray*}
and for every $\ell\in \nz\setminus\left\{0\right\}$\,, using the above general statement \textbf{ii)} $\Rightarrow$ \textbf{iii)},
\begin{eqnarray*}
 \tau(\varepsilon_{2})^{-1}\mu_{\ell}(B(\Id_{E}-R))\leq \mu_{\ell}(B)
\leq \mu_{\ell}(B(\Id_{E}-R)) \tau(\varepsilon_{2})\,.
\end{eqnarray*}
Moreover, the operator $B_{1}=B(1-R)$  clearly sends $E'$ into $F'$\,, and also
sends
$E''$ into $F''=(F')^{\perp}$ according to
$$
\forall u\in E'',\quad BRu=\sum_{j=1}^{\ell_{1}}\langle \psi_{j}\,,\,
Bu\rangle\psi_{j}
= \Pi_{F'}Bu\,.
$$
Since in addition $E'\oplus E''$ is a $\varepsilon_{1}$-orthogonal
decomposition of $E$ and $F=F'\ds\mathop{\oplus}^{\perp}F''$\,,
a direct application of  \textbf{b)} (with $\varepsilon_{2}=0$) says
that the singular values of $B_{1}:E\to F$ and
$\ds\tilde{B}_{1}=\Pi_{F'}B_{1}\big|_{E'}\mathop{\oplus}^{\perp}\Pi_{F''}B_{1}\big|_{E''}$
are related by
$$
\forall \ell\in \nz\setminus\left\{0\right\}\,,\quad
\tau(\varepsilon_{1})^{-1/2}\mu_{\ell}(\tilde{B}_{1})\leq \mu_{\ell}(B_{1})\leq \mu_{\ell}(\tilde{B}_{1})\tau(\varepsilon_{1})^{1/2}\,.
$$
We conclude with 
\begin{eqnarray*}
  &&
     \Pi_{F'}B_{1}\big|_{E'}=\Pi_{F'}[B\big|_{E'}-B\underbrace{R\big|_{E'}}_{=0}]=B\big|_{E'}\,;\\
&&
\Pi_{F''}B_{1}\big|_{E''}=\Pi_{F''}B\big|_{E''}-\underbrace{\Pi_{F''}BR\big|_{E''}}_{=0}=\Pi_{F''}B\big|_{E''}\,.
\end{eqnarray*}
\end{proof}

\begin{remark}
  \begin{description}
  \item[1)] 
 In the sequel, Propositions~\ref{pr:projepsort} and~\ref{pr:epsorth}
 will be used and combined
 with spaces
$E^{h},F^{h},G^{h}$\,, $E'^{h},E''^{h},F'^{h},F''^{h}$\,, operators $B^{h}$\,, $\tilde B^{h}$\,, and  bases $(\varphi^{h}_{j})_{j\in
  \mathcal{J}}$ and $(\psi^{h}_{k})_{k\in \mathcal{K}}$ which depend 
on a small
parameter $h>0$ and such that the hypotheses are satisfied with
$$
\lim_{h\to 0}\varepsilon_{1}(h)=\lim_{h\to
  0}\varepsilon_{2}(h)=0\,.
$$
More generally, note that  when
 $N$ parameters $\varepsilon_{1}(h),\dots, \varepsilon_{N}(h)$ are involved and satisfy
$0\leq \varepsilon_{n}(h)\leq \varrho(h)$ for $n\in \left\{1,\ldots,
  N\right\}$ with $\lim_{h\to 0}\varrho(h)=0$\,, then for
any $\alpha\geq 0$\,, the estimate
$$
\tau(\varepsilon_{1}(h),\ldots, \varepsilon_{N}(h))^{\alpha}=1+O(\varrho(h))
$$
holds uniformly in the sense that there exist $h_{\alpha,N,\varrho},
C_{\alpha,N}>0$ independent of
$\varepsilon_{1},\ldots,\varepsilon_{N}$ such that
$$
\forall h\in ]0,h_{\varrho,N,\alpha}[\,,\ \ 
\tau(\varepsilon_{1}(h),\ldots,\varepsilon_{N}(h))^{\alpha}-1\ \leq\ 
C_{N,\alpha}\,\varrho(h)\,. 
$$
Several applications of the previous results  in this setting will lead to estimates of the type
$$
\forall \ell\in \nz\setminus\left\{0\right\}\,,\quad 
\mu_{\ell}(B^{h})=\mu_{\ell}(\tilde{B}^{h})\big(1+O(\varrho(h))\big)\,.
$$
\item[2)]  A case is especially easy to handle: when
  $E^{h},F^{h},G^{h}$ are finite dimensional with dimension bounded by
  a common number $n_{F}$\,. In this case,
  one can use   any  norm $\|~\|_{n_{F}^{2}}$ on $\mathcal{M}_{n_{F},n_{F}}(\cz)$ in order
  to check the $O\big(\varepsilon_{1,2}(h)\big)$-orthonormality of the
  bases. The constants in the $O(\varrho(h))$-estimates are then fixed
  when $n_{F}$\,, the norm $\|~\|_{n_{F}^{2}}$ and possibly the above
  $N\in \nz$ and $\alpha\geq 0$ are fixed.\\
Additionally, we recall that in this case,
$\vec{d}(E^{h},F^{h})=\vec{d}(F^{h},E^{h})<1$ is equivalent to
$\vec{d}(E^{h},F^{h})<1$ and $\dim E^{h}=\dim F^{h}$\,.
\end{description}
\end{remark}
The following lemma will be useful in the sequel.
\begin{lem}
\label{le:factorization} Let $B^{h}: D(B^{h})
\to \mathcal{H}$\,, $ D(B^{h})\subset  \mathcal{H}$\,, be a closed unbounded operator and assume that the
closed subspaces $E^{h}, F^{h},G^{h}$ and the operator $B^{h}$ satisfy
\begin{itemize}
\item $E^{h}\subset D(B^{h})$\,;
\item  the restriction $B^{h}\big|_{E^{h}}$ is a left multiple of
  $\Pi_{F^{h}}B^{h}\big|_{E^{h}}$:
$$
\xymatrix{
{E}^{h}  \quad\ar[r]^{B^{h}}\ar[dr]_{\Pi_{F^{h}}B^{h}\big|_{E^{h}}}& \mathcal{H}\\
&F^{h}\ar[u]_{C^{h}}}\,;
$$
\item the distance between $F^{h}$ and $G^{h}$ satisfies
$$
\left[\vec{d}(F^{h},G^{h})+\vec{d}(G^{h},F^{h})\right]
\|C^{h}\|=O(\varrho(h))\quad\text{with}\quad
   \lim_{h\to 0}\varrho(h)=0\,.
$$
\end{itemize}
Then
$\Pi_{G^{h}}B^{h}\big|_{E^{h}}=(\Id_{\mathcal{H}}+O(\varrho(h)))\Pi_{F^{h}}B^{h}\big|_{E^{h}}$
and the restriction $B^{h}\big|_{E^{h}}$ is also a left multiple of
$\Pi_{G^{h}}B^{h}\big|_{E^{h}}$:
$$
\xymatrix{
{E}^{h}  \quad\ar[r]^{B^{h}}\ar[dr]_{\Pi_{G^{h}}B^{h}}& \mathcal{H}\\
&G^{h}\ar[u]_{\tilde{C}^{h}}}\,;
$$
with $\tilde{C}^{h}=C^{h}(\Id_{\mathcal{H}}+O(\varrho(h)))$\,. The
roles 
 of $F^{h}$ and $G^{h}$ are therefore symmetric.
\end{lem}
\begin{proof}
Note first that the relation
$B^{h}\big|_{E^{h}}=C^{h}\Pi_{F^{h}}B^{h}\big|_{E^{h}}$
implies 
$$
\|B^{h}\big|_{E^{h}}\|\leq \|C^{h}\|\|\Pi_{F^{h}}\|\|B^{h}\big|_{E^{h}}\|
$$
and then $\|C^{h}\|\geq1$ (except when $B^{h}\big|_{E^{h}}=0$\,, in which case
the statement of Lemma~\ref{le:factorization} is trivial).
Consider now the difference in $\mathcal{L}(E^{h};\mathcal{H})$:
$$
\Pi_{G^{h}}B^{h}\big|_{E^{h}}-\Pi_{G^{h}}
\Pi_{F^{h}}B^{h}\big|_{E^{h}}
=(\Pi_{G^{h}}-\Pi_{G^{h}}\Pi_{F^{h}})C^{h}\Pi_{F^{h}}B^{h}\big|_{E^{h}}\,.
$$
By introducing the operator
$$C_{G^{h}F^{h}}=\Pi_{G^{h}}\Pi_{F^{h}}+(1-\Pi_{G^{h}})(1-\Pi_{F^{h}})
=\Id_{\mathcal{H}}+O(\frac{ \varrho(h)^{2} }{ \|C^{h}\|^{2} })
=
\Id_{\mathcal{H}}+O(\varrho(h)^{2})
$$
like in the proof of Proposition~\ref{pr:projepsort}, we obtain
$$
\Pi_{G^{h}}B^{h}\big|_{E^{h}}=\left[C_{G^{h}F^{h}}+\left(\Pi_{G^{h}}-\Pi_{G^{h}}\Pi_{F^{h}}\right)C^{h}\right]\Pi_{F^{h}}B^{h}\big|_{E^{h}}=[\Id_{\mathcal{H}}+O(\varrho(h))]\Pi_{F^{h}}B^{h}\big|_{E^{h}}\,.
$$
We get
$\Pi_{F^{h}}B^{h}\big|_{E^{h}}=[\Id_{\mathcal{H}}+O(\varrho(h))]^{-1}\Pi_{G^{h}}B^{h}\big|_{E^{h}}$
and we take $\tilde{C}^{h}=C^{h}[\Id_{\mathcal{H}}+O(\varrho(h))]^{-1}$\,.
\end{proof}
We now consider additive error terms which arise in our applications.
\begin{prop}
\label{pr:adderr} Let $B_{1}^{h}, B_{2}^{h}:E^{h}\to F^{h}$ be two
compact operators parametrized by $h>0$\,, like possibly the Hilbert
spaces $E^{h},F^{h}$\,. 
Fix $\ell_{0}\in
\nz\setminus\left\{0\right\}$ and let $\varrho(h)>0$ satisfy
$\lim_{h\to 0}\varrho(h)=0$\,.
\begin{description}
\item[a)] When
  $\|B_{2}^{h}-B_{1}^{h}\|=O(\varrho(h))\, \max\big(\mu_{\ell_{0}}(B_{1}^{h}),
  \mu_{\ell_{0}}(B_{2}^{h})\big)$\,, the singular values are related by
  $$
\forall \ell\in \left\{1,\ldots,\ell_{0}\right\}\,,\quad
\mu_{\ell}(B_{2}^{h})=\mu_{\ell}(B_{1}^{h})\big(1+O(\varrho(h)) \big)\,.
$$
\item[b)] The two following statements are equivalent:
  \begin{eqnarray*}
    &&\min\big(\mu_{\ell_{0}+1}(B_{1}^{h}),\mu_{\ell_{0}+1}(B_{2}^{h})\big)+\|B_{2}^{h}-B_{1}^{h}\|=O\Big(\varrho(h)\max\big(\mu_{\ell_{0}}(B_{1}^{h}),\mu_{\ell_{0}}(B_{2}^{h})\big)\Big)\\
\text{and}&&
\max\big(\mu_{\ell_{0}+1}(B_{1}^{h}),\mu_{\ell_{0}+1}(B_{2}^{h})\big)+\|B_{2}^{h}-B_{1}^{h}\|=O\Big(\varrho(h)\min\big(\mu_{\ell_{0}}(B_{1}^{h}),\mu_{\ell_{0}}(B_{2}^{h})\big)\Big)\,.
  \end{eqnarray*}
\end{description}
\end{prop}
\begin{proof}
The two results are simple consequences of the max-min principle.\\
\textbf{a)} 
Assume $\|B^{h}_{2}-B^{h}_{1}\|\leq \varepsilon \mu_{\ell_{0}}(B^{h}_{1})$ with
$\varepsilon<1$\,. For 
 $\ell\in\{1,\dots,\ell_{0}\}$ and
$V\subset E^{h}$\,, $\dim V=\ell-1$\,, we write
$$
\forall u\in V^{\perp}\,,\ \ 
\|B^{h}_{1}u\|-\varepsilon\mu_{\ell_{0}}(B^{h}_{1})\|u\| \leq \|B^{h}_{2}u\|\leq \|B^{h}_{1}u\|+\varepsilon\mu_{\ell_{0}}(B^{h}_{1})\|u\|
$$
and then, using $\mu_{\ell_{0}}(B^{h}_{1})\leq
\mu_{\ell}(B^{h}_{1})$\,,
\begin{eqnarray*}
\forall u\in V^{\perp},\quad&&
\frac{\|B^{h}_{2}u\|}{\|u\|}\leq
\max_{v\in V^{\perp}\setminus
                              \left\{0\right\}}\frac{\|B^{h}_{1}v\|}{\|v\|}+\varepsilon\mu_{\ell}(B^{h}_{1})\\
&&\frac{\|B^{h}_{1}u\|}{\|u\|}-\varepsilon\mu_{\ell}(B^{h}_{1})
\leq \max_{v\in V^{\perp}\setminus\left\{0\right\}}\frac{\|B^{h}_{2}v\|}{\|v\|}\,.
\end{eqnarray*}
Therefore, for every $\ell\in\{1,\dots,\ell_{0}\}$\,, we deduce
$$
\max_{u\in V^{\perp}\setminus
                              \left\{0\right\}}\frac{\|B^{h}_{1}u\|}{\|u\|}-\varepsilon\mu_{\ell}(B^{h}_{1})\leq
\max_{u\in V^{\perp}\setminus
                              \left\{0\right\}}\frac{\|B^{h}_{2}u\|}{\|u\|}\leq
\max_{u\in V^{\perp}\setminus
                              \left\{0\right\}}\frac{\|B^{h}_{1}u\|}{\|u\|}+\varepsilon\mu_{\ell}(B^{h}_{1})
$$
for any subspace $V$ such that $\dim V=\ell-1$\,. Continuing by taking
the minimum w.r.t $V$ finally leads to 
$$
\forall \ell\in \left\{1,\ldots,\ell_{0}\right\}\,,\quad
\mu_{\ell}(B^{h}_{1})(1-\varepsilon)\ \leq\  \mu_{\ell}(B^{h}_{2})\ \leq\  (1+\varepsilon)\mu_{\ell}(B^{h}_{1})\,.
$$
The $h$-dependent assumption and the symmetry 
$B^{h}_{1}\leftrightarrow B^{h}_{2}$ 
in the above proof yield the result.\\
\textbf{b)}
First, since $\min \leq \max$\,, the second condition obviously implies the
first one. Moreover, the first condition implies 
$\|B_{2}^{h}-B_{1}^{h}\|=O\Big(\varrho(h)\max\big(\mu_{\ell_{0}}(B_{1}^{h}),\mu_{\ell_{0}}(B_{2}^{h})\big)\Big)$
and we deduce from \textbf{a)}
$\max\big(\mu_{\ell_{0}}(B_{1}^{h}),\mu_{\ell_{0}}(B_{2}^{h})\big)=O\big(\min\big(\mu_{\ell_{0}}(B_{1}^{h}),\mu_{\ell_{0}}(B_{2}^{h})\big)$\,.
We have then to show that the second condition is implied by
\begin{equation}
\label{eq.min-min}
\min\big(\mu_{\ell_{0}+1}(B_{1}^{h}),\mu_{\ell_{0}+1}(B_{2}^{h})\big)+\|B_{2}^{h}-B_{1}^{h}\|=O\Big(\varrho(h)\min\big(\mu_{\ell_{0}}(B_{1}^{h}),\mu_{\ell_{0}}(B_{2}^{h})\big)\Big)\,.
\end{equation}
But assuming this and reasoning as in the proof of \textbf{a)}
with $V\subset E^{h}$\,, $\dim V=\ell_{0}$\,, and
 using now $\|B_{2}^{h}-B_{1}^{h}\|=O\big(\varrho(h)\mu_{\ell_{0}}(B_{1}^{h})\big)$\,, leads to
$$
\max_{u\in V^{\perp}\setminus
                              \left\{0\right\}}\frac{\|B^{h}_{2}u\|}{\|u\|}\ =\ 
\max_{u\in V^{\perp}\setminus
                              \left\{0\right\}}\frac{\|B^{h}_{1}u\|}{\|u\|}+O\big(\varrho(h)\mu_{\ell_{0}}(B_{1}^{h})\big)
$$ 
and then, by taking
the minimum w.r.t $V$\,,  to 
$$
 \mu_{\ell_{0}+1}(B^{h}_{2})\ =\  \mu_{\ell_{0}+1}(B^{h}_{1}) +
O\big(\varrho(h)\mu_{\ell_{0}}(B_{1}^{h})\big)\,.
$$
It follows that 
$$\max\big(\mu_{\ell_{0}+1}(B_{1}^{h}),\mu_{\ell_{0}+1}(B_{2}^{h})\big)=
\min\big(\mu_{\ell_{0}+1}(B_{1}^{h}),\mu_{\ell_{0}+1}(B_{2}^{h})\big)+O\big(\varrho(h)\mu_{\ell_{0}}(B_{1}^{h})\big)\,.$$
Then, since $\mu_{\ell_{0}}(B_{1}^{h})=O\big(\min\big(\mu_{\ell_{0}}(B_{1}^{h}),\mu_{\ell_{0}}(B_{2}^{h})\big)$\,, 
 \eqref{eq.min-min} leads to
$$
\max\big(\mu_{\ell_{0}+1}(B_{1}^{h}),\mu_{\ell_{0}+1}(B_{2}^{h})\big)
+\|B_{2}^{h}-B_{1}^{h}\|=O\Big(\varrho(h)\min\big(\mu_{\ell_{0}}(B_{1}^{h}),\mu_{\ell_{0}}(B_{2}^{h})\big)\Big)\,,
$$
which concludes the proof.
\end{proof}

The final result of this section combines multiplicative and
additive error estimates.
\begin{prop}
\label{pr:multadd}
Let $(B^{h},D(B^{h}))$ be a densely defined closed operator in
$\mathcal{H}$\,. Let $E^{h}$\,, $F^{h}$\,, and $G^{h}$ be finite dimensional 
subspaces of $\mathcal{H}$ and let $\varrho(h)>0$ satisfy $\lim_{h\to 0}\varrho(h)=0$\,.
 Assume that both $E^{h}$ and
$F^{h}$ are contained in $D(B^{h})$  and that  the space $E^{h}$ admits the $\varrho(h)$-orthogonal
  decomposition $E^{h}={E'}^{h}\oplus
  {E''}^{h}$\,, such that:
\begin{enumerate}
\item $\Pi_{F^{h}}B^{h}=B^{h}\Pi_{F^{h}}$ on $D(B^{h})$\,;
\item $\Pi_{F^{h}}B^{h}\Pi_{F^{h}}$ has a  fixed
finite rank $\ell_{0}\in\nz$\,;
\item
  $B^{h}\big|_{{E'}^{h}}$ is a
  left multiple of
  $\Pi_{F^{h}}B^{h}\Pi_{F^{h}}\big|_{E'^{h}}=\Pi_{F^{h}}B^{h}\big|_{E'^{h}}$\,:
\begin{displaymath}
\xymatrix{
{E'}^{h}  \quad\ar[r]^{B^{h}}\ar[dr]_{\Pi_{F^{h}}B^{h}\Pi_{F^{h}}}& \mathcal{H}\\
&F^{h}\ar[u]_{C^{h}}
}
\end{displaymath}
\item with the convention $\mu_{0}(A)=+\infty$ for any compact
  operator $A$ and when 
  $\ell_{1}^{h}$ denotes the rank of $\Pi_{G^{h}}B^{h}\big|_{{E'}^{h}}$\,, the following
  inequalities are satisfied:
  \begin{eqnarray}
\label{eq:hypineq1}
  \!\!\!\!\!\!\!\!\!\!\!\!\!\!\!\!\!\!\!\!\!\!\!\!\!  &&\vec{d}(E^{h},F^{h})+\vec{d}(F^{h},E^{h})+\|C^{h}\|\left(\vec{d}(F^{h},G^{h})+\vec{d}(G^{h},F^{h})\right)
=O(\varrho(h))\,,\\
\label{eq:hypineq2}
\!\!\!\!\!\!\!\!\!\!\!\!\!\!\!\!\!\!\!\!\!\!\!\!\! &&
\|B^{h}\big|_{{E''}^{h}}\|\left[
\frac{1}{\mu_{\ell_{1}^{h}}(\Pi_{G^{h}}B^{h}\big|_{{E'}^{h}})}
+
   \frac{\|C^{h}\|(\vec{d}(F^{h},G^{h})+\vec{d}(G^{h},F^{h}))}{\max(\mu_{\ell_{0}}(\Pi_{G^{h}}B^{h}\big|_{E^{h}}),\mu_{\ell_{0}}(B^{h}\big|_{F^{h}}))}
\right]=O(\varrho(h)).
 \end{eqnarray}
\end{enumerate}
Then, the $\ell_{0}$ first singular values of
$\Pi_{F^{h}}B^{h}\Pi_{F^{h}}$ and $\Pi_{G^{h}}B^{h}\Pi_{E^{h}}$
satisfy
\begin{equation}
\label{eq:multaddl}
\forall \ell\in \left\{1,\ldots,\ell_{0}\right\}\,,\quad
\underbrace{
\mu_{\ell}(\Pi_{G^{h}}B^{h}\Pi_{E^{h}})}_{=\mu_{\ell}(\Pi_{G^{h}}B^{h}\big|_{E^{h}})}=
\underbrace{\mu_{\ell}(\Pi_{F^{h}}B^{h}\Pi_{F^{h}})}_{=\mu_{\ell}(B^{h}\big|_{F^{h}})}\big(1+O(\varrho(h))\big)\,.
\end{equation}
Moreover, the $\ell_{0}+1$-th singular value of
$\Pi_{G^{h}}B^{h}\Pi_{E^{h}}$ satisfies
\begin{equation}
\label{eq:multaddl0p1}
\frac{\mu_{\ell_{0}+1}(\Pi_{G^{h}}B^{h}\Pi_{E^{h}})}{\mu_{\ell_{0}}(\Pi_{G^{h}}B^{h}\Pi_{E^{h}})}
=\frac{\mu_{\ell_{0}+1}(\Pi_{G^{h}}B^{h}\big|_{E^{h}})}{\mu_{\ell_{0}}(\Pi_{G^{h}}B^{h}\big|_{E^{h}})}
\stackrel{h\to 0}{\sim} \frac{\mu_{\ell_{0}+1}(\Pi_{G^{h}}B^{h}\big|_{E^{h}})}{\mu_{\ell_{0}}(B^{h}\big|_{F^{h}})}=O(\varrho(h))\,.
\end{equation}
\end{prop}
\begin{proof}
 Since the statement is trivial when $\ell_{0}= 0$\,, we assume here that $\ell_{0}\geq 1$\,.
The assumptions 1. and  3. then imply
$\|C^{h}\|\geq
1$ because
$$
\|B^{h}\big|_{E'^{h}}\|\leq \|C^{h}\|\|\Pi_{F^{h}}\|\|B^{h}\big|_{E'^{h}}\|
$$
(except when $B^{h}\big|_{E'^{h}}=0$\,, in which case one chooses $C^{h}=\Pi_{F^{h}}$
so that $\|C^{h}\|=1$).
Therefore, the first estimate \eqref{eq:hypineq1} of 4. implies $\dim E^{h}=\dim F^{h}=\dim
G^{h}<\infty$ as well
as $\vec{d}(E^{h},F^{h})=\vec{d}(F^{h},E^{h})$ and $\vec{d}(F^{h},G^{h})=\vec{d}(G^{h},F^{h})$\,.\\
About dimensions, the assumptions 1. and 3. also imply
$$\mathrm{rank}(\Pi_{F^{h}}B^{h}\big|_{E'^{h}})\ =\ \mathrm{rank}(B^{h}\big|_{E'^{h}})\ =\ \ell_{1}^{h}\ \leq
\ 
\ell_{0}\,.$$ This rank $\ell_{1}^{h}$\,, 
which is not assumed to be independent of $h$\,, will 
be proved to be equal to $\text{rank}(\Pi_{G^{h}}B\big|_{E'^{h}})$\,.\medskip

\noindent
\textbf{Multiplicative estimates:} By replacing $E^{h}$ by $E'^{h}$ in
Lemma~\ref{le:factorization}, we get
$$
\Pi_{G^{h}}B^{h}\big|_{E'^{h}}=[\Id_{\mathcal{H}}+O(\varrho(h))]\Pi_{F^{h}}B^{h}\big|_{E'^{h}}
$$
and therefore
\begin{equation}
\label{eq.vs-l-h}
\forall \ell\in \left\{1,\ldots,\dim E'^{h}\right\}\,,\quad
\mu_{\ell}(\Pi_{G^{h}}B^{h}\big|_{E'^{h}})=\mu_{\ell}(\Pi_{F^{h}}B^{h}\Pi_{F^{h}}\big|_{E'^{h}})\big(1+O(\varrho(h))\big)\,.
\end{equation}
In particular,
$$\text{rank}(\Pi_{G^{h}}B^{h}\big|_{E'^{h}})
\ =\ \text{rank}(\Pi_{F^{h}}B^{h}\big|_{E'^{h}})\ =\ \ell_{1}^{h}\ =\ \text{rank}(B^{h}\big|_{E'^{h}})\,.
$$
An accurate information about the orthogonal projections
on $F'^{h}:=\Ran \,\Pi_{F^{h}}B^{h}\big|_{E'^{h}}$ and on $G'^{h}:=\Ran\,
\Pi_{G^{h}}B^{h}\big|_{E'^{h}}$ is achieved as follows. There exist
two  orthonormal systems $(\varphi_{j}^{h})_{1\leq j\leq \ell_{1}^{h}}$
in $E'^{h}$ and $(\psi_{j}^{h})_{1\leq j\leq \ell_{1}^{h}}$ in
$F'^{h}\subset F^{h}$
such that 
$$
\forall j\in \left\{1,\ldots, \ell_{1}^{h}\right\}\,,\quad 
\Pi_{F^{h}}B^{h}\varphi_{j}^{h}=\mu_{j}^{h}\psi_{j}^{h}\ , \quad\text{where}\quad \mu_{j}^{h}=\mu_{j}(\Pi_{F^{h}}B^{h}\big|_{E'^{h}})>0\,.
$$
By computing
\begin{eqnarray*}
  \psi_{j}^{h}-\frac{1}{\mu_{j}^{h}}\Pi_{G^{h}}B^{h}\varphi_{j}^{h}&=&
\frac{1}{\mu_{j}^{h}}\left[\Pi_{F^{h}}B^{h}\varphi_{j}^{h}-\Pi_{G^{h}}B^{h}\varphi_{j}^{h}\right]\\
&=&\frac{1}{\mu_{j}^{h}}(\Pi_{F^{h}}-\Pi_{G^{h}})(C^{h}\Pi_{F^{h}}B^{h}\Pi_{F^{h}}\varphi_{j}^{h})\\
&=&
(\Pi_{F^{h}}-\Pi_{F^{h}}\Pi_{G^{h}})C^{h}\psi_{j}^{h}
+(\Pi_{F^{h}}\Pi_{G^{h}}-\Pi_{G^{h}})C^{h}\psi_{j}^{h}\,,
\end{eqnarray*}
we obtain the estimates:
$$
\forall j\in \left\{1,\ldots, \ell_{1}^{h}\right\}\,,\quad
\|\psi_{j}^{h}-\frac{1}{\mu_{j}^{h}}\Pi_{G^{h}}B^{h}\varphi_{j}^{h}\|\leq
\|C^{h}\|\big(\vec{d}(F^{h},G^{h})+\vec{d}(G^{h},F^{h})\big)
\underbrace{=O(\varrho(h))}_{\eqref{eq:hypineq1}}\,.
$$
Since moreover $\text{rank}~\Pi_{G^{h}}B^{h}\big|_{E'^{h}}=\dim
G'^{h}=\ell_{1}^{h}\leq \ell_{0}$\,, it follows that 
$\left(\frac{1}{\mu_{j}^{h}}\Pi_{G^{h}}B^{h}\varphi_{j}^{h}\right)_{1\leq
j\leq \ell_{1}^{h}}$ is an
$O\Big( \|C^{h}\|\big(\vec{d}(F^{h},G^{h})+\vec{d}(G^{h},F^{h})\big) \Big)$-orthonormal 
basis of $G'^{h}$ (see Definition~\ref{de:epsorth}) and then that
$$
\|\Pi_{F'^{h}}-\Pi_{G'^{h}}\|= O\Big( \|C^{h}\|\big(\vec{d}(F^{h},G^{h})+\vec{d}(G^{h},F^{h})\big) \Big) =O(\varrho(h)).
$$
By calling $F''^{h}$ the orthogonal of $F'^{h}$ in $F^{h}$  and
$G''^{h}$ the orthogonal of $G'^{h}$ in $G^{h}$\,, the equality
\begin{eqnarray*}
\Pi_{F''^{h}}-\Pi_{G''^{h}}&=&(1-\Pi_{F'^{h}})\Pi_{F^{h}}-(1-\Pi_{G'^{h}})\Pi_{G^{h}}\\
&=&(1-\Pi_{F'^{h}})(\Pi_{F^{h}}-\Pi_{G^{h}}\Pi_{F^{h}})-(\Pi_{F'^{h}}-\Pi_{G'^{h}})\Pi_{G^{h}}\Pi_{F^{h}}\\
&&\hspace{5cm}-(1-\Pi_{G'^{h}})(\Pi_{G^{h}}-\Pi_{G^{h}}\Pi_{F^{h}})
\end{eqnarray*}
now implies (using also $\|C^{h}\|\geq
1$)
\begin{equation}
\label{eq:pisecdiff}
\|\Pi_{F''^{h}}-\Pi_{G''^{h}}\|=O\Big( \|C^{h}\|\big(\vec{d}(F^{h},G^{h})+\vec{d}(G^{h},F^{h})\big) \Big)=O(\varrho(h))\,.
\end{equation}
The separation between the $\ell_{1}^{h}$ first singular values and the
smaller ones is obtained by
applying Proposition~\ref{pr:epsorth}-\textbf{c)} to
$B=\Pi_{F^{h}}B^{h}\big|_{E^{h}}: E^{h}\to F^{h}$ and to
$B=\Pi_{G^{h}}B^{h}\big|_{E^{h}}: E^{h}\to G^{h}$ with: 
the $\varrho(h)$-orthogonal decomposition $E^{h}=E'^{h}\oplus E''^{h}$\,,
\begin{eqnarray*}
  &&
\mu_{\ell_{1}^{h}}(\Pi_{F^{h}}B^{h}\big|_{E'^{h}})=\mu_{\ell_{1}^{h}}(\Pi_{G^{h}}B^{h}\big|_{E'^{h}})\big(1+O(\varrho(h))\big)\,,\\
\text{and}&&
             \frac{\|\Pi_{F^{h}}B^{h}\big|_{E''^{h}}\|}{\mu_{\ell_{1}^{h}}(\Pi_{F^{h}}B^{h}\big|_{E'^{h}})}+
             \frac{\|\Pi_{G^{h}}B^{h}\big|_{E''^{h}}\|}{\mu_{\ell_{1}^{h}}(\Pi_{G^{h}}B^{h}\big|_{E'^{h}})}\leq
             C\frac{\|B^{h}\big|_{E''^{h}}\|}{\mu_{\ell_{1}^{h}}(\Pi_{G^{h}}B^{h}\big|_{E'^{h}})}\underbrace{=}_{\eqref{eq:hypineq2}}O(\varrho(h))\,.
\end{eqnarray*}
This implies that the singular values of $\Pi_{F^{h}}B^{h}\big|_{E^{h}}$
and of $\Pi_{G^{h}}B^{h}\big|_{E^{h}}$ satisfy
\begin{align}
\label{eq.vs-l-h'}
\forall \ell\in \left\{1,\ldots,\ell_{1}^{h}\right\}\,,\quad
 &\mu_{\ell}(\Pi_{F^{h}}B^{h}\big|_{E^{h}})
=\mu_{\ell}(\Pi_{F^{h}}B^{h}\big|_{E'^{h}})\big(1+O(\varrho(h))\big)\,,&\\
\label{eq.vs-l-h''}
& \mu_{\ell}(\Pi_{G^{h}}B^{h}\big|_{E^{h}})
=\mu_{\ell}(\Pi_{G^{h}}B^{h}\big|_{E'^{h}})\big(1+O(\varrho(h))\big)\,,&
\end{align}
and, for every $k\geq 1$\,,
\begin{align}
\label{eq.vs-l-h-3}
&  \mu_{\ell_{1}^{h}+k}(\Pi_{F^{h}}B^{h}\big|_{E^{h}})
=\mu_{k}(\Pi_{F''^{h}}B^{h}\big|_{E''^{h}})\big(1+O(\varrho(h))\big)
=O\big(\mu_{\ell_{1}^{h}}(\Pi_{F^{h}}B^{h}\big|_{E^{h}})\varrho(h) \big)\,,\\
\label{eq.vs-l-h-4}
 &\mu_{\ell_{1}^{h}+k}(\Pi_{G^{h}}B^{h}\big|_{E^{h}})
=\mu_{k}(\Pi_{G''^{h}}B^{h}\big|_{E''^{h}})\big(1+O(\varrho(h))\big)=O\big(\mu_{\ell_{1}^{h}}(\Pi_{G^{h}}B^{h}\big|_{E^{h}})\varrho(h)\big)\,.
\end{align}
 Besides, using $\vec{d}(E^{h},F^{h})+\vec{d}(F^{h},E^{h})=O(\varrho(h))$ 
and the commutation
$\Pi_{F^{h}}B^{h}\big|_{E^{h}}=\Pi_{F^{h}}B^{h}\Pi_{F^{h}}\big|_{E^{h}}$\,, 
 a direct application of Proposition~\ref{pr:projepsort} with
$B=\Pi_{F^{h}}B^{h}\big|_{F^{h}}:F^{h}\to F^{h}$ and
$\tilde{B}=\Pi_{F^{h}}B^{h}\Pi_{F^{h}}\big|_{E^{h}}=\Pi_{F^{h}}B^{h}\big|_{E^{h}}:E^{h}\to
F^{h}$ leads to:
\begin{equation}
\label{eq.vs-l-h-5}
\forall \ell\in\nz\setminus\left\{0\right\}\,,\quad
\mu_{\ell}(\Pi_{F^{h}}B^{h}\big|_{F^{h}})
=\mu_{\ell}(\Pi_{F^{h}}B^{h}\big|_{E^{h}})\big(1+O(\varrho(h)^{2})\big)\,.
\end{equation}
\textbf{Additive estimates:} 
When $\ell_{0}=\ell_{1}^{h}$\,, the statement of Proposition~\ref{pr:multadd}
follows from the equations \eqref{eq.vs-l-h} and \eqref{eq.vs-l-h'}--\eqref{eq.vs-l-h-5} and,
when $\ell_{0}>\ell_{1}^{h}$\,, these equations reduce
the problem  to the comparison
of the singular values $\mu_{k}$\,, $1\leq k\leq
\ell_{0}-\ell_{1}^{h}$\,, of the two operators
$$
\Pi_{G''^{h}}B^{h}\big|_{E''^{h}}\quad\text{and}\quad
\Pi_{F''^{h}}B^{h}\big|_{E''^{h}}\,.
$$
By \eqref{eq:pisecdiff} and \eqref{eq:hypineq2}, we know that
$$
\frac{\|B^{h}\big|_{E''^{h}}\|\|\Pi_{G''^{h}}-\Pi_{F''^{h}}\|}{\max(\mu_{\ell_{0}-\ell_{1}^{h}}(\Pi_{G''^{h}}B^{h}\big|_{E''^{h}}),\mu_{\ell_{0}-\ell_{1}^{h}}(\Pi_{F''^{h}}B^{h}\big|_{E''^{h}}))}=O(\varrho(h))\,.
$$
The first result \eqref{eq:multaddl} is thus an application of Proposition~\ref{pr:adderr}-\textbf{a)} with
$$
B_{1}^{h}=\Pi_{F''^{h}}B^{h}\big|_{E''^{h}}\quad\text{and}\quad
B_{2}^{h}=\Pi_{G''^{h}}B^{h}\big|_{E''^{h}}
$$
while replacing $\ell_{0}$ by $\ell_{0}-\ell_{1}^{h}$\,.\\
Lastly, the definition of $\ell_{0}$ in 2. implies
$$
\min(\mu_{\ell_{0}-\ell_{1}^{h}+1}(\Pi_{G''^{h}}B^{h}\big|_{E''^{h}}),\mu_{\ell_{0}-\ell_{1}^{h}+1}(\Pi_{F''^{h}}B^{h}\big|_{E''^{h}}))=
\mu_{\ell_{0}-\ell_{1}^{h}+1}(\Pi_{F''^{h}}B^{h}\big|_{E''^{h}})=0\,.
$$
The remaining statement \eqref{eq:multaddl0p1} is then given by Proposition~\ref{pr:adderr}-\textbf{b)}.
\end{proof}

\section{Accurate analysis with $N$ ``critical values''}
\label{sec:accanN}
This section is the core and the most technical part of our
analysis. It combines: i) the exponential decay estimates of eigenvectors solving
$\Delta_{f,f^{-1}([a,b]),h}\omega_{h}=\lambda_{h}\omega_{h}$\,,
$\lambda_{h}\stackrel{h\to 0}{\to}0$\,, and all
the properties of solutions to $d_{f,h}\omega_{h}=0$ stated in
Sections~\ref{sec:expdec} and~\ref{sec:locpbs}; ii) the information on local problems, that is
when $\sharp ([a,b]\cap\left\{c_{1},\ldots,c_{N_{f}}\right\})=1$\,, from
Section~\ref{sec:locpbs}; iii) 
the rough estimates when $\sharp
([a,b]\cap\left\{c_{1},\ldots,c_{N_{f}}\right\})=N$ of
Section~\ref{sec:rough}. Finally,  the recurrence analysis with respect to $N$
is modelled on linear algebra lemmas about singular values given in
Section~\ref{sec:singval}.
In the first paragraph, we review and complete previous useful notations
before stating a general result which leads easily to
Theorem~\ref{th:mainsimple},  
variations of which will be given
afterwards. It is about the construction of global quasimodes for
$\Delta_{f,f^{-1}([a,b]),h}$\,, and more precisely of a suitable basis of
 widely extended
solutions to 
$d_{f,h}\omega_{h}=0$\,,
 which, contrarily to the  eigenfunctions of
$\Delta_{f,f^{-1}([a,b]),h}$\,, provide a high flexibility when changing
the geometrical domain, in particular the values $a,b$\,. After
specifying the framework in the first paragraph, 
we check in Subsection~\ref{sec:initrec}
 the first step of the recursive
construction of such global quasimodes and presents the strategy of our
iterative method, developed in the other paragraphs.

\subsection{Assumptions, notations and main result}
\label{sec:notmainN}

We assume Hypothesis~\ref{hyp:cN}
which is: The function $f$ has a finite number of ``critical values'', $c_{1}<\ldots<c_{N_{f}}$\,,
according to Hypothesis~\ref{hyp:mainf} or
Hypothesis~\ref{hyp:Lipbar}, while  Hypothesis~\ref{hyp:AgmonLip} is
assumed for a  general Lipschitz function $f$\,, and we choose
\begin{equation}
\label{eq:etafbis}
\eta_{f}\in ]0,\frac{1}{2}\min_{1<n\leq N_{f}}|c_{n}-c_{n-1}|[\,.
\end{equation}
Moreover, the values $a,b$\,, $-\infty\leq a<b\leq +\infty$\,, are not
``critical values'' of $f$\,.\\
Like in Sections~\ref{sec:locpbs} and \ref{sec:rough}, we use the
the space $W(f_{a}^{b};\Lambda T^{*}M)$ of Definition~\ref{de:defW}.
We recall that it coincides with
$W^{1,2}(f_{a}^{b};\Lambda T^{*}M)$ under Hypothesis~\ref{hyp:mainf},
while we only know $W(f_{a}^{b};\Lambda
T^{*}M)\subset W^{1,2}_{loc}(f_{a}^{b};\Lambda
T^{*}M)$ in general (when $a,b\not\in \left\{c_{1},\ldots,c_{N_{f}}\right\}$).
According to this remark,
when $E=\sqcup_{k=1}^{K}]a_{k},b_{k}[$\,, $a_{k},b_{k}\not\in
\left\{c_{1},\ldots,c_{N_{f}}\right\}$\,, the space $W(f^{-1}(E);\Lambda
T^{*}M)$ is nothing but the direct sum
$\mathop{\oplus}_{k=1}^{K}W(f_{a_{k}}^{b_{k}};\Lambda T^{*}M)$\,,
which is included in $W^{1,2}_{loc}(f^{-1}(E);\Lambda T^{*}M)$\,.\\
The set of ``critical values'' lying in $[a,b]$ are relabelled according
to 
$$
[a,b]\cap
\left\{c_{1},\ldots,c_{N_{f}}\right\}=\left\{\tilde{c}_{1},\ldots,\tilde{c}_{N}\right\}\quad,\quad a<\tilde{c}_{1}<\ldots<\tilde{c}_{N}<b\,.
$$
The bar code associated with $f$ is still denoted by
$\mathcal{B}=\mathcal{B}(f)=([a_{\alpha},b_{\alpha}[)_{\alpha\in A}$\,.
We keep the notation $A^{*}(a,b)$\,, $A_{c}^{*}(a,b)$ given in
\eqref{eq:defAab},\eqref{eq:defAcab}, while the endpoints of the
bars with a non trivial intersection with $]a,b[$ are partitionned
into 
$$
\mathcal{J}^{*}(a,b)=\mathcal{X}^{*}(a,b)\sqcup \mathcal{Y}^{*}(a,b)\sqcup\mathcal{Z}^{*}(a,b)\,,
$$
where the definition of those sets are given in
\eqref{eq:defXab},\eqref{eq:defYab},\eqref{eq:defZab}, and
\eqref{eq:defJab}. Remember that an element $j\in \mathcal{J}^{(p)}(a,b)$ is
a pair $j=(\alpha,\tilde{c})$ with $\alpha\in A^{*}(a,b)$ and
$\tilde{c}\in \left\{\tilde{c}_{1},\ldots, \tilde{c}_{N}\right\}$\,,
$\tilde{c}=x_{\alpha}^{(p)}$\,, $y_{\alpha}^{(p)}$\,, or $z_{\alpha}^{(p)} $\,,
depending on wether:
\begin{itemize}
\item $j\in \mathcal{X}^{(p)}(a,b)$\,, which means $\alpha\in
A^{(p)}_{c}(a,b)$ and $\tilde{c}=x_{\alpha}^{(p)}$\,, 
\item $j\in \mathcal{Y}^{(p)}(a,b)$\,, which means $\alpha\in
A_{c}^{(p-1)}(a,b)$ and $\tilde{c}=y_{\alpha}^{(p)}$\,, 
\item or
$j\in \mathcal{Z}^{(p)}(a,b)$\,, which means $\alpha\in
A^{*}(a,b)\setminus A_{c}^{*}(a,b)$ and
$\tilde{c}=z^{(p)}_{\alpha}$\,. 
\end{itemize}
Below are figures which summarize the three different cases.\\

\medskip
\figinit{cm}
\figpt1000:(-7,-1)
\figpt1001:(7,-1)
\figpt 0:(-5,-1)
\figpt 1:(5,-1)
\figpt 2:(-3,-1)
\figpt 3:(2,-1)
\figpt 4:(-3,0)
\figpt 5:(2,0)

\figpt 6:(-0.5,0)
\figdrawbegin{}
\figset(width=1.3)
\figdrawline[1000,1001]
\figset(width=2)
\figdrawline[4,5]
\figset(with=1)
\figset(dash=4)
\figdrawline[2,4]
\figdrawline[3,5]
\figdrawend

\figvisu{\figBoxA}{{\textbf{Figure 9:} A bar
$[x_{\alpha}^{(p)},y_{\alpha}^{(p+1)}[=[\tilde{c}_{m},\tilde{c}_{n}[$\,, $\alpha^{(p)}\in
A_{c}^{(p)}(a,b)$\,.}}{
\figwriten 6:$\alpha^{(p)}$(0.3)
\figwritew 4:$x_{\alpha}^{(p)}$(0.2)
\figwritee 5:$y_{\alpha}^{(p+1)}$(0.2)
\figset write(mark=$+$)
\figwrites 0: $a$(0.3)
\figwrites 1:$b$(0.3)
\figwrites 2:$\tilde{c}_{m}$(0.3)
\figwrites 3:$\tilde{c}_{n}$(0.3)
\figset write(mark=$\bullet$)
\figwrites 4:$$(0)
}
\centerline{\box\figBoxA}
\begin{center}
 There are  two extreme points
 $j=(\alpha^{(p)},\tilde{c}_{m})\in
\mathcal{X}^{(p)}(a,b)$ and $j'=(\alpha^{(p)},\tilde{c}_{n})\in
 \mathcal{Y}^{(p+1)}(a,b)$\,.
\end{center}

\medskip
\figinit{cm}
\figpt1000:(-7,-1)
\figpt1001:(7,-1)
\figpt 0:(-5,-1)
\figpt 1:(5,-1)
\figpt 2:(-1,-1)
\figpt 3:(-5,0)
\figpt 4:(-1,0)
\figpt 5:(-1,0.5)
\figpt 6:(5,0.5)

\figpt 8:(-3,0)
\figpt 9:(2,0.5)
\figdrawbegin{}
\figset(width=1.3)
\figdrawline[1000,1001]
\figset(width=2)
\figdrawline[3,4]
\figdrawline[5,6]
\figset(with=1)
\figset(dash=4)
\figdrawline[2,5]
\figdrawend
\figvisu{\figBoxA}{{\textbf{Figure 10:} An extreme point
$j=(\alpha,\tilde{c})\in \mathcal{Z}^{(p)}(a,b)$\,.}}{
\figwriten 8:$\alpha^{(p)}$(0.3)
\figwriten 9:$\alpha^{(p+1)}$(0.3)
\figset write(mark=$+$)
\figwrites 0: $a$(0.3)
\figwrites 1:$b$(0.3)
\figwrites 2:$\tilde{c}=z_{\alpha}^{(p)}$(0.3)
\figset write(mark=$\bullet$)
\figwrites 5:$$(0)
}
\centerline{\box\figBoxA}
\begin{center}
After restriction to $[a,b]$\,, this represents the two possible
equivalent ways of having $j=(\alpha^{*},z_{\alpha}^{(p)})\in {\cal
Z}^{(p)}(a,b)$\,.
\end{center}

\medskip

We recall that
\begin{eqnarray*}
\delta_{[0,e^{-\varepsilon/h}],[a,b],h}^{(p)}&=&
1_{[0,e^{-\varepsilon/h}]}(\Delta^{(p+1)}_{f,f^{-1}([a,b]),h})
d^{(p)}_{f,f^{-1}([a,b]),h}
1_{[0,e^{-\varepsilon/h}]}(\Delta^{(p)}_{f,f^{-1}([a,b]),h})\\
\text{and}\quad F_{[0,e^{-\varepsilon/h}],[a,b],h}^{(p)}&=&\Ran~1_{[0,e^{-\varepsilon/h}]}(\Delta^{(p)}_{f,f^{-1}([a,b]),h})
\end{eqnarray*}
do not depend on $\varepsilon\in ]0,\varepsilon_{0}[$\,,
provided that  $h_{\varepsilon}>0$ ($h\in]0,h_{\varepsilon}[$)
 is chosen small enough when $\varepsilon$ is fixed. We then use  the notation
\begin{equation}
\label{eq:deltap}
  \delta_{[0,\tilde{o}(1)],[a,b],h}^{(p)}=\delta_{[0,e^{-\varepsilon/h}],[a,b],h}^{(p)}\quad\text{and}\quad
F_{[0,\tilde{o}(1)],[a,b],h}^{(p)}=F_{[0,e^{-\varepsilon/h}],[a,b] ,h}^{(p)}
\end{equation}
without mentioning $\varepsilon>0$\,.\\
The exponent ${}^{(p)}$ is
 forgotten when the direct sum w.r.t $p\in
\left\{0,\ldots, \dim M\right\}$ is considered.\\
The distance between vector spaces $\vec{d}(E,F)$ is the one defined
in Subsection~\ref{sec:usqunot} (see Definition~\ref{de:dEF}) and used in
Section~\ref{sec:singval}. We also recall that for $\varepsilon>0$\,, an
$\tilde{O}(e^{-\frac{\varepsilon}{h}})$-orthonormal family of vectors
$(\varphi_{\ell}^{h})_{1\leq \ell\leq L}$ in a Hilbert space
$\mathcal{H}$ is a family such that $|\langle
\varphi_{\ell}^{h},\varphi_{\ell'}^{h}\rangle-\delta_{\ell,\ell'}|=\tilde{O}(e^{-\frac{\varepsilon}{h}})$
according to Definition~\ref{de:epsorth}.\\
With the family $\mathcal{J}^{*}(a,b)$ of endpoints of bars with a
non trivial intersection with $]a,b[$\,, we will associate an
$\tilde{o}(1)$-orthonormal family of solutions to
$d_{f,h}\omega_{h}=0$ in the proper range. 
\begin{definition}
\label{de:adapted}
Under Hypothesis~\ref{hyp:cN} and for $\delta_{1}\in
]0,\frac{\eta_{f}}{8}]$\,, let
 \begin{equation}
\label{eq:defSd}
S_{\delta_{1}}:=\left\{\tilde{c}_{n}-\delta_{1},
  \tilde{c}_{n}+\delta_{1}, \quad 1\leq n\leq N\right\}\,.
\end{equation}
A family $(\varphi_{j}^{*,h})_{j\in \mathcal{J}^{*}(a,b)}$\,,
$\varphi_{j}^{*,h}=\varphi_{j}^{(p),h}$ when $j\in
\mathcal{J}^{(p)}(a,b)$\,, is called a $\delta_{1}$-family of
quasimodes if there exists $\gamma:]0,h_{0}[\to]0,+\infty[$ with
$\lim_{h\to 0}\gamma(h)=0$ such that:
  \begin{itemize}
  \item $(\varphi_{j}^{(p),h})_{j\in \mathcal{J}^{(p)}(a,b)}$ is a
    linearly independent family of $D(d_{f,f^{-1}([a,b]),h}^{(p)})$ for all $p\in
    \left\{0,\ldots, d\right\}$\,;
\item by setting $j=(\alpha,\tilde{c})$ and 
$I_{j}^{h}=[x_{\alpha}^{(p)}-\delta_{1},
  y_{\alpha}^{(p+1)}-\gamma(h)]=[\tilde{c}-\delta_{1},y_{\alpha}^{(p+1)}-\gamma(h)]$
  when $j\in \mathcal{X}^{(p)}(a,b)$\,,
and $I_{j}^{h}=[\tilde{c}-\delta_{1}, b]$ when $j\in \mathcal{Y}^{(p)}(a,b)\cup \mathcal{Z}^{(p)}(a,b)$\,:
  \begin{eqnarray}
\label{eq:suppfj}
&&
\supp \varphi_{j}^{(p),h}\subset
   f^{-1}\left((I_{j}^{h}
+[0,\gamma(h)/2])\cap
   [a,b]\right)\,,
\\
\label{eq:decfj}
&&
\|e^{\frac{|f-\tilde{c}|}{h}}\varphi_{j}^{(p),h}\|_{W(f^{-1}([a,b])\setminus S_{\delta_{1}})}
=\tilde{O}(1)
\,,\\
\label{eq:dffj}
&&
d_{f,h}\varphi_{j}^{(p),h}\equiv 0\quad
   \text{in}~f^{-1}(I_{j}^{h}\cap
   [a,b])\ \ \text{and then
   in}~f^{-1}([a,\tilde{c}-\delta_{1}]\cup (I_{j}^{h} \cap
   [a,b])).
  \end{eqnarray}
\end{itemize}
For such a family of quasimodes, we will use the notation:
$$
\mathcal{V}^{(p),h}=\Vect(\varphi_{j}^{(p),h}\,, j\in
  \mathcal{J}^{(p)}(a,b))\quad,\quad\mathcal{V}^{h}=\mathop{\oplus}_{p=0}^{d}\mathcal{V}^{(p),h}\,.
$$
\end{definition}
The idea is that the quasimode associated with the endpoint $j=(\alpha,\tilde{c})\in
\mathcal{J}^{*}(a,b)$ is supported in $[\tilde{c}-\delta_{1},b]$\,,
decays exponentially away from $\tilde{c}$\,, and solves
$d_{f,h}\varphi_{j}^{*,h}=0$ in a region essentially covered by 
the bar indexed by $\alpha$\,. Global quasimodes for $d_{f,h}$
are constructed by climbing along the values of $f$\,.  The reason why
$W$-estimates fail in a neighborhood of 
$f^{-1}(S_{\delta_{1}})$ will appear in  the construction of such a family (see in
particular Remark~\ref{re.S-delta-1} 
about the
values $\tilde{c}_{n}+\delta_{1}$).\\
The following
definition specifies how such quasimodes are truncated around
 the
upper endpoints $y_{\alpha}^{(p+1)}$ when
$j=(\alpha,x_{\alpha}^{(p)})\in \mathcal{X}^{(p)}(a,b)$\,. This
truncation operator preserves  the spaces 
$W(f^{-1}(I))$ for $I\subset [a,b]$ and
$D(d_{f,f^{-1}([a,b]),h})$ with its boundary conditions.
 \begin{definition}
\label{de:Td2}
  In the framework of Definition~\ref{de:adapted} and for
  $\delta_{2}\in ]0,\frac{\eta_{f}}{8}]$\,, let
\begin{equation}
\label{eq:defchic}
\chi_{\tilde{c}_{n},\delta_{2}}(x)=\chi\left(\frac{f(x)-\tilde{c}_{n}}{\delta_{2}}\right)
\end{equation}
for $n\in \left\{2,\ldots, N\right\}$ and a fixed $\chi\in
\mathcal{C}^{\infty}(\rz;[0,1])$ such that $\chi\equiv 1$ on
$]-\infty,-2]$ and $\textrm{supp}~\chi\subset]-\infty,-1[$\,.\\
The operator $T_{\delta_{2}}$ is defined on 
$\mathcal{V}^{h}$ by
\begin{equation}
\label{eq:defTd2}
T_{\delta_{2}}\varphi_{j}^{(p),h}=
\left\{
    \begin{array}[c]{ll}
\chi_{y_{\alpha}^{(p+1)},\delta_{2}}\varphi_{j}^{(p),h}&\text{if}~j=(\alpha,x_{\alpha}^{(p)})\in
\mathcal{X}^{(p)}(a,b)
\\
\varphi_{j}^{(p),h}&\text{if}~j\in \mathcal{Y}^{(p)}(a,b)\cup
\mathcal{Z}^{(p)}(a,b)\,.
    \end{array}
\right.
\end{equation}
\end{definition}
\begin{thm}
\label{th:induc}
Assume Hypothesis~\ref{hyp:cN} with $\eta_{f}$ given by
\eqref{eq:etafbis}.
\begin{description}
\item[a)] For any $p\in \left\{0,\ldots,\dim M\right\}$\,, the
  $\tilde{o}(1)$  non zero singular values  
of $d^{(p)}_{f,f^{-1}([a,b]),h}$\,, that is the non zero singular values of
$\delta^{(p)}_{[0,\tilde{o}(1)],[a,b],h}$\,, can be labelled by the family $\left(\mu_{j}^{h}\right)_{j\in \mathcal{X}^{(p)}(a,b)}$ (with possible multiplicities) with 
$$\mu_{j}^{h}\stackrel{log}{\sim}e^{-\frac{y_{\alpha}^{(p+1)}-x_{\alpha}^{(p)}}{h}}\,,\quad
j=(\alpha,x_{\alpha}^{(p)})\in \mathcal{X}^{(p)}(a,b)\,.
$$
\item[b)] For any $\delta_{1}\in ]0,\frac{\eta_{f}}{8}]$\,, there exists
  a $\delta_{1}$-family  $(\varphi_{j}^{*,h})_{j\in
    \mathcal{J}^{*}(a,b)}$ of quasimodes in the sense of
  Definition~\ref{de:adapted} which is
  $\tilde{O}(e^{-\frac{\delta_{1}}{h}})$-orthonormal in
  $L^{2}(f_{a}^{b})$\,.\\
The  vector space $\mathcal{V}^{h}$ spanned by those quasimodes  satisfies:
$$
\forall p\in \{0,\dots,\dim M\}\,,\ \ \ \vec{d}(\mathcal{V}^{(p),h},F_{[0,\tilde{o}(1)],[a,b],h}^{(p)})+\vec{d}(F_{[0,\tilde{o}(1)],[a,b],h}^{(p)},
\mathcal{V}^{(p),h})=\tilde{O}(e^{-\frac{\delta_{1}}{h}})\,.
$$
\item[c)] If   $T_{\delta_{2}}$ is the
  truncation operator of Definition~\ref{de:Td2} for $\delta_{2}\in
  ]0,\frac{\eta_{f}}{8}]$\,, then the map
 $d_{f,f^{-1}([a,b]),h}^{(p)}T_{\delta_{2}}:\mathcal{V}^{(p),h}\to
  L^{2}(f^{-1}([a,b]))$ is a left multiple of
  $\delta_{[0,\tilde{o}(1)],[a,b],h}^{(p)}T_{\delta_{2}}$~:
\begin{equation}
\label{eq:recleftmF}
\xymatrix@C=3cm{
\mathcal{V}^{(p),h}
\quad\ar[r]^{d_{f,f^{-1}([a,b]),h}^{(p)}T_{\delta_{2}}}\ar[dr]_{\underbrace{\delta_{[0,\tilde{o}(1)],[a,b],h}^{(p)}}_{~\eqref{eq:deltap}}T_{\delta_{2}}}&
L^{2}(f^{-1}([a,b]))\\
&F_{[0,\tilde{o}(1)],[a,b],h}^{(p)}\ar[u]_{C^{h}}
}
\end{equation}
 with $\|C^{h}\|=\tilde{O}(e^{\frac{2\delta_{2}}{h}})$\,.
\end{description}
\end{thm}
The proof will be done in several
steps, by induction on the number
of ``critical values'' $N$\,. 
 Because the graduation w.r.t $p\in \left\{0,\ldots,\dim
    M\right\}$ is associated with an obvious orthogonal decomposition
  of $F_{[0,\tilde{o}(1)], [a,b],h}$ and
  $\delta_{[0,\tilde{o}(1)],[a,b],h}$\,, and clear partitions of the sets
  of indices for bars and endpoints,
  $A(a,b)=\sqcup_{p=0}^{d}A^{(p)}(a,b)$\,, $\mathcal{J}(a,b)=\sqcup_{p=0}^{d}\mathcal{J}^{(p)}(a,b)$\,, etc., we can treat globally
  $F_{[0,\tilde{o}(1)],[a,b],h}$ and
  $\delta_{[0,\tilde{o}(1),[a,b],h}$ and forget the degree $p$\,.

\subsection{Initialisation and outline of the recurrence}
\label{sec:initrec}
\noindent\textbf{The result holds true for $N=1$:} 
According
to Proposition~\ref{pr:exp0},
we know that $\mathcal{J}(a,b)=\mathcal{Z}(a,b)$ and that the
$\tilde{o}(1)$-eigenvalues of
$\Delta_{f,f^{-1}([a,b]),h}$\,, and therefore the $\tilde{o}(1)$-singular
values of $\delta_{[0,\tilde{o}(1)],[a,b],h}$\,,  all vanish . This proves
\textbf{a)}. 
To prove \textbf{b)}, it suffices to take an orthonormal basis $(\varphi_{j}^{h})_{j\in \mathcal{J}(a,b)}$ of
$\ker(\Delta_{f,f^{-1}([\max(a,\tilde{c}_{1}-\delta_{1}),b]),h})$\,, 
extended by $0$
on
$f_{a}^{\tilde{c}_{1}-\delta_{1}}$ when $a<\tilde{c}_{1}-\delta_{1}$\,. Note that in the latter case, 
the extended family $(\varphi_{j}^{h})_{j\in \mathcal{J}(a,b)}$ is
still included in $D(d_{f,f^{-1}([a,b]),h})$\,, and actually in
$\ker(d_{f,f^{-1}([a,b]),h})$\,. 
The exponential decay estimate \eqref{eq:decfj}
comes from the exponential decay estimates on the $\varphi_{j}^{h}\in
\ker(\Delta_{f,f^{-1}[\max(a,\tilde{c}_{1}-\delta_{1}),b],h})$ given by 
Proposition~\ref{pr:Agmon} or Hypothesis~\ref{hyp:AgmonLip}
applied with $\overline{\Omega}=f^{-1}([\max(a,\tilde{c}_{1}-\delta_{1}),b])$\,,
$r_{h}=0$\,, $\lambda_{h}=0$\,,
$U=f^{-1}(\left\{\tilde{c}_{1}\right\})$ 
and $d_{Ag}(x,y)\geq |f(x)-f(y)|$\,. The distance between $\mathcal{V}^{h}$ and
$F_{[0,\tilde{o}(1)],[a,b],h}$ is also deduced from the exponential
decay estimates on the $\varphi_{j}^{h}\in
\ker(\Delta_{f,f^{-1}([\max(a,\tilde{c}_{1}-\delta_{1}),b]),h})$ as we
did  in the proofs of Propositions~\ref{pr:comptageest} and~\ref{pr:distkernsum}. 
The statement \textbf{c)} is obvious in this case because 
\begin{eqnarray*}
&&  
T_{\delta_{2}}=\Id_{\mathcal{V}^{h}}\,,\\
&&
d_{f,f^{-1}([a,b]),h}\big|_{\mathcal{V}^{h}}=0\,,\\
\text{and}&&
\delta_{[0,\tilde{o}(1)],[a,b],h}\big|_{\mathcal{V}^{h}}=\Pi_{[0,\tilde{o}(1)],[a,b],h}d_{f,f^{-1}([a,b]),h}\big|_{\mathcal{V}^{h}}=0\,.
\end{eqnarray*}

\noindent\textbf{Strategy of the proof by induction:} 
\begin{enumerate}
\item Already while checking the initial step $N=1$ or when proving
  e.g. Proposition~\ref{pr:exp0} in Subsection~\ref{sec:exp0}, it was
  convenient to work with different values of $a$ and $b$\,. From this
  point of view, the construction of $\delta_{1}$-quasimodes in the sense of
  Definition~\ref{de:adapted}, which are some specific solutions to
  $d_{f,h}\omega_{h}=0$\,, is more flexible than working with spectral
  eigenvectors of $\Delta_{f,f^{-1}([a,b]),h}$\,. Note that even though the
  extension by $0$ in $f^{a}$ of $\varphi\in
  \ker(d_{f,f^{-1}([a,b]),h})$ does not belong to
  $W(f_{a'}^{b})$ for $a'<a$\,, it belongs to
  $\ker(d_{f,f^{-1}([a',b]),h})$\,. This provides a way to extend the
  quasimodes in the area of the lower values of $f$\,. The extension to
  $f_{a}^{b'}$ with $b<b'$ will be done with a repeated use of
  Proposition~\ref{pr:interdfh}. Note for example
 that if there is no  ``critical value'' in $[b,b']$\,, a solution to
 $d_{f,f^{-1}([a,b]),h}\varphi_{h}=0$\,,  which satisfies some
 exponential decay estimates of the type
 $\|e^{\frac{f(x)}{h}}\varphi_{h}\|_{W(f_{a}^{b})}\leq \tilde{O}(C_{h})$\,, can be
 ``extended'' to a solution to
 $d_{f,f^{-1}([a,b']),h}\tilde{\varphi}_{h}=0$\,, with the same decay estimates
 in $W(f_{a}^{b'}\setminus f^{-1}(\{b-\delta\})$  for some $\delta>0$ small
 enough. To prove this, it suffices to consider $b$ as
 an artificial new ``critical value'' $\tilde{c}$ and to apply
 Proposition~\ref{pr:interdfh}-{\bf i)} with $a_{0},a,\tilde{c}_{1},b$ there
 replaced by $a, b-\delta, \tilde{c}=b, b'$\,.
Note that with this extension procedure,
$\tilde{\varphi}_{h}$  fails in general to belong to $W$
in  a neighborhood of 
$f^{-1}(\{b-\delta\})$ 
(see \eqref{eq:decfj} in this connection).
  If there is a ``critical value'' $\tilde{c}_{n}\in ]b,b'[$\,, then
 one has to study more carefully the orthogonality condition of
 Proposition~\ref{pr:interdfh}-{\bf ii)}.
\item Now Theorem~\ref{th:induc} will be assumed to be true in the case of $N$ ``critical
  values'' in $[a,b]$\,, we can deduce several consequences. The
  aforementionned  flexibility of a family of quasimodes in
  $\mathcal{V}^{h}$\,, as compared to
  a family of eigenvectors for the initial space
  $F_{[0,\tilde{o}(1)],[a,b],h}$ , can be completed 
  by replacing the arrival space 
  $F_{[0,\tilde{o}(1)],[a,b],h}$  in the diagram
  \eqref{eq:recleftmF} by a more 
flexible approximation. Moreover, the
  $\tilde{O}(e^{-\frac{\delta_{1}}{h}})$-orthogonality of the  $\delta_{1}$-family of
 quasimodes can be preserved while ensuring  true
  orthogonality properties on the images $d_{f,h}T_{\delta_{2}}\varphi_{j}^{h}$\,. 
This will be done in Subsection~\ref{sec:conseqN}. The corresponding
results will be used in the rest of the proof and for other
constructions later.
\item Let us now explain how we pass from the case of
$N$ critical values $\tilde c_{1}<\cdots <\tilde c_{N}$ to the case of 
$N+1$ critical values $\tilde c_{1}<\cdots <\tilde c_{N+1}$ in $[a,b]$\,. To do so,
introduce $a_{2}\in ]\tilde{c}_{1},\tilde{c}_{2}[$ and
$b_{1}\in]\tilde{c}_{N},\tilde{c}_{N+1}[$\,, set $a_{1}=a$\,,
$b_{2}=b$\,, and apply the result valid for $N$ ``critical values'' to
$a_{1}=a<\tilde{c}_{1}<\ldots<\tilde{c}_{N}< b_{1}$ and
to
$a_{2}<\tilde{c}_{2}\ldots <\tilde{c}_{N+1}<b_{2}=b$\,. From the
$\delta_{1}$-families of quasimodes for  the intervals 
$[a_{1},b_{1}]$ and $[a_{2},b_{2}]$\,, we can extract a
partial $\delta_{1}$-family of quasimodes for $[a,b]$ which satisfies the required properties for all bars of length
strictly smaller than
$\tilde{c}_{N+1}-\tilde{c}_{1}$\,. This construction, and all the
information coming from step $N$ in the intervals 
$[a_{1},b_{1}]$ and $[a_{2},b_{2}]$\,, is collected in
Subsection~\ref{sec:collect}. After this, in
Subsection~\ref{sec:longlength}, 
 the construction of $\delta_{1}$-quasimodes associated
with bars $j=(\alpha,x_{\alpha})\in \mathcal{X}(a,b)$ with
$x_{\alpha}=\tilde{c}_{1}$ and
$y_{\alpha}=\tilde{c}_{N+1}$ must be specified.
This leads to the definition of ``intermediate $\delta_{1}$-family
of quasimodes''
(see Definition~\ref{de:interdelta1})
 which, comparatively to Definition~\ref{de:adapted},
does not yet elucidate the interaction with the local spectral
problems 
around the ``critical value'' $\tilde{c}_{N+1}$\,.
This strategy is summarized in Figure~11 below. It is related to
Mayer-Vietoris type arguments in algebraic topology, but handling  and
propagating all the estimates on exponentially small quantities
requires some care. From this point of view, the inspiration is also
taken from the standard techniques for handling exponential decay estimates, and several \underline{up and down} inductions
on $n\in \left\{1,\ldots, N+1\right\}$ are used.\\

\medskip
\figinit{0.9cm}
\figpt1000:(-8,-1)
\figpt1001:(8,-1)
\figpt1002:(5,-1)
\figpt 0:(-7,-1)
\figpt 1:(7,-1)
\figpt 2:(-6,-1)
\figpt 3:(-3,-1)
\figpt 4:(3,-1)
\figpt 5:(6,-1)
\figpt 6:(-1,-1)
\figpt 7:(2,-1)
\figpt 8:(4,-1)
\figpt 9:(-4.5,-1)

\figpt 100:(-7,0)
\figpt 101:(-3,0)
\figpt 102:(-3,1.5)
\figpt 103:(3,1.5)
\figpt 104:(-6,0.3)
\figpt 105:(-3,0.3)
\figpt 106:(-6,0.6)
\figpt 107:(3,0.6)
\figpt 108:(-7,0.9)
\figpt 109:(3,0.9)
\figpt 110:(-3,1.2)
\figpt 111:(3,1.2)
\figpt 112:(-3,1.8)
\figpt 113:(5,1.8)
\figpt 114:(3,2.1)
\figpt 115:(5,2.1)
\figpt 116:(3,2.4)
\figpt 117:(7,2.4)
\figpt 118:(-7,2.7)
\figpt 119:(5,2.7)
\figpt 120:(5,3.0)
\figpt 121:(7,3.0)
\figpt 122:(-6,4.8)
\figpt 123:(5,4.8)
\figpt 124:(7,4.8)
\figpt 125:(-6,5.1)
\figpt 126:(5,5.1)
\figpt 127:(7,5.1)

\figpt 200:(-4.5,-0.4)
\figpt 201:(7,-0.4)
\figpt 202:(7,4.2)
\figpt 203:(-4.5,4.2)
\figpt 210:(-7,-0.8)
\figpt 211:(4,-0.8)
\figpt 212:(4,3.8)
\figpt 213:(-7,3.8)
\figpt 220:(6,4.95)
\def\Rx{2.}
\def\Ry{0.4}
\figdrawbegin{}
\figset(width=1.3)
\figdrawline[1000,6]
\figdrawline[7,1001]
\figset(dash=2)
\figdrawline[6,7]
\figset(dash=1)
\figset(width=2)
\figdrawline[100,101]
\figdrawline[102,103]
\figdrawline[104,105]
\figdrawline[106,107]
\figdrawline[108,109]
\figdrawline[110,111]
\figdrawline[112,113]
\figdrawline[114,115]
\figdrawline[116,117]
\figdrawline[118,119]
\figdrawline[120,121]
\figdrawline[122,123]
\figdrawline[125,126]
\figset(dash=4)
\figdrawline[123,124]
\figdrawline[126,127]

\figset(width=1)
\figset(dash=3)
\figdrawline[200,201,202,203,200]
\figdrawline[210,211,212,213,210]
\figdrawarcell 220;\Rx,\Ry (0,360,0)
\figdrawend

\figvisu{\figBoxA}{\textbf{Figure 11:} Positions of the bars while the
interval $[a,b]$ is covered by $[a_{1},b_{1}]\cup [a_{2},b_{2}]$\,.}{
\figwriten 220:$??$(0.6)
\figset write(mark=$+$)
\figwrites 1002:$\tilde{c}_{N+1}$(0.3)
\figwrites 0: $a=a_{1}$(0.3)
\figwrites 1:$b=b_{2}$(0.3)
\figwrites 2:$\tilde{c}_{1}$(0.3)
\figwrites 3:$\tilde{c}_{2}$(0.3)
\figwrites 4:$\tilde{c}_{N}$(0.3)
\figwrites 8:$b_{1}$(0.3)
\figwrites 9:$a_{2}$(0.3)
\figset write(mark=$\bullet$)
\figwrites 102:$$(0.2)
\figwrites 104:$$(0)
\figwrites 106:$$(0)
\figwrites 112:$$(0)
\figwrites 110:$$(0)
\figwrites 114:$$(0)
\figwrites 116:$$(0)
\figwrites 120:$$(0)
\figwrites 122:$$(0)
\figwrites 125:$$(0)
}
\centerline{\box\figBoxA}

\begin{center}
We will use the recurrence hypothesis at step $N$ first in the
interval $]a_{2},b_{2}[$ and then in the interval $]a_{1},b_{1}[$\,,
where the corresponding  proper bars (not equal to $]a_{i},b_{i}[$)
are collected in dashed rectangles. Quasimodes in $]a_{2},b_{2}[$ are
extended by $0$ in $f_{a}^{a_{2}}$\,, while the extension of
quasimodes in $]a_{1},b_{1}[$ to $f_{b_{1}}^{b}$ requires more care.
\end{center}

Once the latter ``intermediate $\delta_{1}$-family of quasimodes'' is
constructed, it is used in order to prove
Theorem~\ref{th:induc}-\noindent\textbf{a)} in Subsection~\ref{sec:lowerN+1}.
 Like in the proof of
Proposition~\ref{pr:exp0} for $N=1$\,, we have to play with different
values of $a,b$ such that
$a<\tilde{c}_{1}<\ldots<\tilde{c}_{N+1}<b$\,. 
Using
Proposition~\ref{pr:roughmino}, we deduce firstly
a lower bound
$r(h)\stackrel{\log}{\sim}e^{-\frac{\tilde{c}_{N+1}-\tilde{c}_{1}+\max(\delta_{1},\delta_{3})}{h}}$
when $a=\tilde{c}_{1}-\delta_{1}$ and $b=\tilde{c}_{N+1}+\delta_{3}$\,,
and translate it in the various variations of the operator
$\delta_{[0,\tilde{o}(1)],[a,b],h}T_{\delta_{2}}$ that we have
introduced. Secondly, we study the effect of changing $a$ and $b$ while
keeping $N+1$ ``critical values'' in $[a,b]$ as it was done in
Subsection~\ref{sec:changinb} for the case of one ``critical value''
in $[a,b]$\,. 
Thirdly, and only 
after proving Theorem~\ref{th:induc}-\textbf{a)}, we can
construct in Subsection~\ref{sec:constN+1} 
the $\delta_{1}$-family of quasimodes $(\varphi_{j}^{h})_{j\in
  \mathcal{J}(a,b)}$ at step $N+1$\,,
and check all the conditions stated in
the items \textbf{b)} and \textbf{c)} of
Theorem~\ref{th:induc}.
\end{enumerate}

\subsection{Consequences of Theorem~\ref{th:induc} at step $N$}
\label{sec:conseqN}

We assume in this section that Theorem~\ref{th:induc} holds true at step $N$\,. We
refer in particular to the Definition~\ref{de:adapted}  of
$\delta_{1}$-quasimodes $(\varphi_{j}^{h})_{j\in
  \mathcal{J}(a,b)}$ and of $\mathcal{V}^{h}=\Vect(\varphi_{j}^{h})_{j\in
\mathcal{J}(a,b)}$\,, and to the Definition~\ref{de:Td2} of the truncation
$T_{\delta_{2}}:\mathcal{V}^{h}\to D(d_{f,f^{-1}([a,b]),h})$\,, 
for $\delta_{1},\delta_{2}\in]0,\frac{\eta_{f}}{8}]$\,.\\
While keeping the initial space $\mathcal{V}^{h}$ for
$d_{f,f^{-1}([a,b]),h}T_{\delta_{2}}$\,, we replace the arrival space
$F_{[0,\tilde{o}(1)],[a,b],h}$\,, and therefore the left-multiplying
projection $\Pi_{[0,\tilde{o}(1)],[a,b],h}$\,, by a more flexible
space $G^{h}$ and a projection $\Pi_{G^{h}}$\,.
In view of Lemma~\ref{le:factorization} and of the general analysis
of singular values led in Section~\ref{sec:singval}, consider
\begin{equation}
\label{eq:defGh}
G^{h}=\ker(\Delta_{f,\overline{\Omega},h})
\quad\text{and}\quad
F^{h}=F_{[0,\tilde{o}(1)],[a,b],h}\,,
\end{equation}
where  $\Delta_{f,\overline{\Omega},h}$ is
 the operator introduced in
\eqref{eq:laplspl} with
\begin{equation}
\label{eq:defOmbar}
\overline{\Omega}=
\mathop{\bigsqcup}_{n=1}^{N}f^{-1}([\tilde{c}_{n}-\eta_{f},\tilde{c}_{n}+\eta_{f}]\cap
[a,b])\,.
\end{equation}
We recall that according to Proposition~\ref{pr:distkernsum},
\begin{eqnarray*}
&&\vec{d}(G^{h},F^{h})+\vec{d}(F^{h},G^{h})=\tilde{O}(e^{-\frac{\eta_{f}}{h}})
\\
\text{and}&&
\dim G^{h}=\dim F^{h}\,.
\end{eqnarray*}
The interest of the space $G^{h}$ is that it is defined in terms of
local spectral problems, actually kernels of local Witten 
Laplacians,
around the ``critical values'' $\tilde{c}_{1},\ldots,
\tilde{c}_{N}$\,. 
\begin{prop}
\label{pr:piGhdfhTd}
Assume that Theorem~\ref{th:induc} holds true at step $N$ and let
$G^{h}$ be defined by \eqref{eq:defGh}. The operator
$$
\Pi_{G^{h}}d_{f,f^{-1}([a,b]),h}T_{\delta_{2}}=\Pi_{G^{h}}d_{f,h}T_{\delta_{2}}:\mathcal{V}^{h}\to
L^{2}(\Omega)\subset L^{2}(f_{a}^{b})
$$
does not depend on $\delta_{2}\in
]0,\frac{\eta_{f}}{8}]$ for $h>0$ small
enough. Namely, for two different choices $\delta_{2},\delta_{2}'\in
]0,\frac{\eta_{f}}{8}]$\,, there exists $h_{\delta_{2},\delta_{2}'}>0$
such that the equality
$\Pi_{G^{h}}d_{f,h}T_{\delta_{2}}=\Pi_{G^{h}}d_{f,h}T_{\delta_{2}'}$
holds for all $h\in ]0,h_{\delta_{2},\delta_{2}'}[$\,.\\
Its singular values satisfy:
\begin{equation}
  \label{eq:muGdT}
\forall \ell \in \{1,\ldots,\dim F^{h}\}\,,\quad
   \mu_{\ell}(\Pi_{G^{h}}d_{f,f^{-1}([a,b]),h}T_{\delta_{2}}\big|_{\mathcal{V}^{h}})=\mu_{\ell}(\delta_{[0,\tilde{o}(1)],[a,b],h})(1+\tilde{O}(e^{-\frac{\delta_{1}}{h}}))\,.
 \end{equation}
Its kernel equals
\begin{equation}
\label{eq:eqker}
\ker(\Pi_{G^{h}}d_{f,h}T_{\delta_{2}})=\Vect(\varphi_{j}^{h}\,,
j\in \mathcal{Y}(a,b)\cup \mathcal{Z}(a,b))\,.
\end{equation}
In particular, when the non zero
singular values of $\delta_{[0,\tilde{o}(1)],[a,b],h}$ are labelled as
$(\mu_{j}^{h})_{j\in \mathcal{X}(a,b)}$
with $\mu_{j}^{h}\stackrel{log}{\sim}
e^{-\frac{y_{\alpha}-x_{\alpha}}{h}}$ for $j=(\alpha,x_{\alpha})$\,, the same result holds for the
$\delta_{2}$-independent operator
$\Pi_{G^{h}}d_{f,h}T_{\delta_{2}}$\,.
\end{prop}
\begin{proof}
The Definition~\ref{de:adapted} of $(\varphi_{j}^{h})_{j\in
  \mathcal{J}(a,b)}$ and
the Definition~\ref{de:Td2} of $T_{\delta_{2}}$ give
$$
d_{f,h}T_{\delta_{2}}\varphi_{j}^{h}=d_{f,h}\varphi_{j}^{h}=0\quad
\text{if}~j\in \mathcal{Y}(a,b)\cup \mathcal{Z}(a,b)\,,
$$
and
$$
d_{f,h}T_{\delta_{2}}\varphi_{j}^{h}=0\quad\text{in}\;
f^{-1}([a,y_{\alpha}-2\delta_{2}])
\cup f^{-1}([y_{\alpha}-\delta_{2},b])\quad\text{if}~j=(\alpha,x_{\alpha})\in \mathcal{X}(a,b)\,.
$$
We deduce firstly
$$
\ker(\Pi_{G^{h}}d_{f,h}T_{\delta_{2}})\supset \Vect(\varphi_{j}^{h}\,,
j\in \mathcal{Y}(a,b)\cup \mathcal{Z}(a,b))\,.
$$
In the case
$j=(\alpha,x_{\alpha})\in \mathcal{X}(a,b)$\,,
the equality
$\Pi_{G^{h}}d_{f,h}T_{\delta_{2}}=\Pi_{G^{h}}d_{f,h}T_{\delta_{2}'}$
for $h>0$ small enough is secondly a direct consequence of 
 Proposition~\ref{pr:interdfh}-\textbf{i)} applied  around the
 ``critical value'' 
 $y_{\alpha}$\,, owing to 
$\supp
d_{f,h}T_{\delta_{2}}\varphi_{j}^{h}\subset f^{-1}(]y_{\alpha}-\eta_{f},y_{\alpha}[)$ 
and to
$$
\Pi_{G^{h}}d_{f,h}T_{\delta_{2}}\varphi_{j}^{h}=
\Pi_{\{0\},[y_{\alpha}-\eta_{f},y_{\alpha}+\eta_{f}]\cap
[a,b],h}d_{f,h}T_{\delta_{2}}\varphi_{j}^{h}\,.
$$
The result \eqref{eq:muGdT} on singular values implies 
$\dim \ker(\Pi_{G^{h}}d_{f,h}T_{\delta_{2}})=\sharp
\mathcal{Y}(a,b)\cup \mathcal{Z}(a,b)$ and yields the
equality~\eqref{eq:eqker}. Let us now prove \eqref{eq:muGdT}.\\
Consider the initial vector space $E^{h}=T_{\delta_{2}}
\mathcal{V}^{h}=\Vect(T_{\delta_{2}}\varphi_{j}^{h}, j\in
\mathcal{J}(a,b))$
and the mapping
$B^{h}=d_{f,f^{-1}([a,b]),h}:E^{h}\to L^{2}(\Omega)\subset
L^{2}(f_{a}^{b})$\,. The distance to $\mathcal{V}^{h}$ is estimated by
\begin{equation}
\label{eq.T-delta-close}
\vec{d}(E^{h},\mathcal{V}^{h})+\vec{d}(\mathcal{V}^{h},E^{h})=\tilde{O}(e^{-\frac{\eta_{f}}{h}})
\leq\tilde{O}(e^{-\frac{\delta_{1}}{h}})\,.
\end{equation}
  With the factorization
\eqref{eq:recleftmF} stated in Theorem~\ref{th:induc}-{\bf c)} and
$\vec{d}(G^{h},F^{h})+\vec{d}(F^{h},G^{h})=\tilde{O}(e^{-\frac{\eta_{f}}{h}})$ with
$2\delta_{2}<\frac{\eta_{f}}{2}<\eta_{f}$\,, 
 we are
exactly in the framework of Lemma~\ref{le:factorization} with
$\varrho(h)=\tilde{O}(e^{\frac{2\delta_{2}-\eta_{f}}{h}})$\,. Therefore,
$d_{f,f^{-1}([a,b]),h}\big|_{E^{h}}$ is a left multiple of
$\Pi_{G^{h}}d_{f,f^{-1}([a,b]),h}\big|_{E^{h}}$\,,
\begin{equation}
\label{eq:factorGh}
\xymatrix{
{E}^{h}  \quad\ar[r]^{B^{h}}\ar[dr]_{\Pi_{G^{h}}B^{h}}& L^{2}(f_{a}^{b})\\
&G^{h}\ar[u]_{\tilde{C}^{h}}}
\end{equation}
with
$\tilde{C}^{h}=C^{h}(\Id_{L^{2}(f_{a}^{b})}+\tilde{O}(e^{\frac{2\delta_{2}-\eta_{f}}{h}}))$\,,
and
\begin{eqnarray*}
  &&\Pi_{G^{h}}d_{f,f^{-1}([a,b]),h}\big|_{E^{h}}=(\Id_{L^{2}(f_{a}^{b})}+\tilde{O}(e^{\frac{2\delta_{2}-\eta_{f}}{h}}))\underbrace{\Pi_{F^{h}}d_{f,f^{-1}([a,b],h)}\big|_{E^{h}}}_{=\delta_{[0,\tilde{o}(1)],[a,b],h}\Pi_{F^{h}}\big|_{E^{h}}}\,.
\end{eqnarray*}
Using additionally 
Proposition~\ref{pr:projepsort} and the relation
$\vec{d}(E^{h},F^{h})+\vec{d}(F^{h},E^{h})=\tilde{O}(e^{-\frac{\delta_{1}}{h}})$
arising from Theorem~\ref{th:induc}-{\bf b)} and \eqref{eq.T-delta-close}, this leads to
\begin{eqnarray*}
&& \forall \ell \in \{1,\ldots,\dim F^{h}\}\,,\quad
   \mu_{\ell}(\Pi_{G^{h}}d_{f,f^{-1}([a,b]),h}\big|_{E^{h}})=\mu_{\ell}(\delta_{[0,\tilde{o}(1)],[a,b],h})(1+r(h))\,,
\end{eqnarray*}
where $r(h)=\max(\tilde{O}(e^{\frac{2\delta_{2}-\eta_{f}}{h}}),
\tilde{O}(e^{-\frac{2\delta_{1}}{h}}))\leq\tilde{O}(e^{-\frac{\delta_{1}}{h}})$\,.
The comparison \eqref{eq:muGdT} for
$\Pi_{G^{h}}d_{f,h}T_{\delta_{2}}\big|_{\mathcal{V}^{h}}$ is then a
consequence of 
\begin{equation}
\label{eq.T-delta-unit}
\|T_{\delta_{2}}^{*}T_{\delta_{2}}-\Id_{\mathcal{V}^{h}}\|+\|T_{\delta_{2}}T_{\delta_{2}}^{*}-\Id_{E^{h}}\|=\tilde{O}(e^{-\frac{\eta_{f}}{h}})=\tilde{O}(e^{-\frac{\delta_{1}}{h}})\,.
\end{equation}
\end{proof}
Below are details about a useful
 block decomposition  of the operator
$\Pi_{G^{h}}d_{f,h}T_{\delta_{2}}:\mathcal{V}^{h}\to
L^{2}(\Omega)$\,.  Of course, there is the orthogonal block
decomposition with respect to the degre $p$ according to
$\Pi_{G^{h}}d_{f,h}^{(p)}T_{\delta_{2}}:\mathcal{V}^{(p),h}\to
L^{2}(\Omega;\Lambda^{p+1}T^{*}M)$\,. But we consider here a block
decomposition according to the length of the bars, which correspond to
 clusters of singular values. Again, we forget the degree $p$ here.
We need some notations.
Let 
\begin{eqnarray}
\label{eq:defXn}
&& \mathcal{X}_{n}(a,b)=\left\{j=(\alpha,x_{\alpha})\in \mathcal{X}(a,b), \;
   y_{\alpha}=\tilde{c}_{n}\right\}\,,\quad 2\leq n\leq N\,,\\
\label{eq:defXmn}&&
\mathcal{X}_{m,n}(a,b)=\left\{j=(\alpha,x_{\alpha})\in \mathcal{X}_{n}(a,b)\,,\;
   x_{\alpha}=\tilde{c}_{m}\right\}\,, \quad 1\leq m<n\leq N\,,
\end{eqnarray}
and consider the following $\tilde{O}(e^{-\frac{\delta_{1}}{h}})$-orthogonal
decompositions:
\begin{eqnarray}
\label{eq:defVmn}
&& \mathcal{V}_{m,n}^{h}=\Vect(\varphi_{j}^{h}, \; j\in
   \mathcal{X}_{m,n}(a,b))\quad\text{for}\quad 1 \leq m<n\leq N\,,\\ 
\label{eq:defVn}&&
\mathcal{V}_{n}^{h}=\mathop{\oplus}_{m=1}^{n-1}\mathcal{V}^{h}_{m,n}=\Vect(\varphi_{j}^{h}\,,\;
                   j\in \mathcal{X}_{n}(a,b))\quad\text{for}\quad 2<n\leq N\,,
\\
\label{eq:defV+}
&&\mathcal{V}_{+}^{h}=\mathop{\oplus}_{n=2}^{N}\mathcal{V}_{n}^{h}\,,\\
\label{eq:defV0}
&&\mathcal{V}_{0}^{h}=\Vect(\varphi_{j}^{h}\,,\;j\in
     \mathcal{Y}(a,b)\cup \mathcal{Z}(a,b))=
\ker(\Pi_{G^{h}}d_{f,h}T_{\delta_{2}}\big|_{\mathcal{V}^{h}})\,,
\\
\label{eq:V=V+V0}
&&\mathcal{V}^{h}=\mathcal{V}_{+}^{h}\oplus
   \mathcal{V}_{0}^{h}\,,\\
\nonumber
\text{with}&&
\Pi_{G^{h}}d_{f,h}T_{\delta_{2}}\mathcal{V}_{n}^{h}\subset
              \ker(\Delta_{f,f^{-1}([\tilde{c}_{n}-\eta_{f},\tilde{c}_{n}+\eta_{f}]\cap[a,b]),h})\,,\\
\label{eq:defGn}
\text{while}&&
G^{h}=\mathop{\oplus}^{\perp}_{n\in \left\{1,\ldots, N\right\}}
              \underbrace{\ker(\Delta_{f,f^{-1}([\tilde{c}_{n}-\eta_{f},\tilde{c}_{n}+\eta_{f}]\cap[a,b]),h})}_{=:G_{n}^{h}}\,.
\end{eqnarray}
\begin{prop}
\label{pr:Vmn}
Under the assumptions of Proposition~\ref{pr:piGhdfhTd} and with the
notations \eqref{eq:defVmn}--\eqref{eq:defGn}, the operator
$\Pi_{G^{h}}d_{f,h}T_{\delta_{2}}\big|_{\mathcal{V}_{m,n}^{h}}:\mathcal{V}_{m,n}^{h}\to
G_{n}^{h}$ 
 is, for $1\leq m<n\leq N$\,, one to one, and,
 when $\dim \mathcal{V}_{m,n}^{h}\neq 0$\,,
  its 
singular values
all satisfy
$\mu^{h}\stackrel{log}{\sim}e^{-\frac{\tilde{c}_{n}-\tilde{c}_{m}}{h}}$\,.\\
Moreover, the non zero singular values of
$\Pi_{G^{h}}d_{f,h}T_{\delta_{2}}:\mathcal{V}^{h}\to L^{2}(\Omega)$
(resp. of
$\Pi_{G_{n}^{h}}d_{f,h}T_{\delta_{2}}\big|_{\mathcal{V}_{n}^{h}}:\mathcal{V}_{n}^{h}\to
L^{2}(\Omega)$\,, where $n\in\{2,\dots,N\}$ is fixed)
are obtained by collecting all those non zero 
singular values for $1\leq
m<n\leq N$ (resp. for $1\leq m <n$)\,, with an $\tilde{O}(e^{-\frac{\delta_{1}}{h}})$ relative error.
\end{prop}
\begin{proof}
For every $1\leq m<n\leq N$ such that $\mathcal{X}_{m,n}(a,b)\neq \emptyset$\,,
the composition of  the exponential decay estimates
  on the $\varphi_{j}^{h}$ given in \eqref{eq:decfj}, $j\in
\mathcal{X}_{m,n}(a,b)$\,,  and on the elements of any
orthonormal basis $(\psi_{k}^{h})_{1\leq k\leq K_{n}}$ of
$G_{n}^{h}=\ker(\Delta_{f,f^{-1}([\tilde{c}_{n}-\eta_{f},\tilde{c}_{n}+\eta_{f}]\cap[a,b]),h})$ 
leads to
\begin{equation}
\label{eq:estimupuh}
\forall u\in \mathcal{V}_{m,n}^{h}\,,\quad
\|\Pi_{G^{h}}d_{f,h}T_{\delta_{2}}u\|=\tilde{O}(e^{-\frac{\tilde{c}_{n}-\tilde{c}_{m}}{h}})\|u\|\,.
\end{equation}
Let us now prove by reductio ad absurdum that
$$
\forall \,1\leq m<n\leq N\ \ \text{such that $\mathcal{X}_{m,n}(a,b)\neq \emptyset$}\,,\ \ \forall u\in \mathcal{V}_{m,n}^{h}\,,\quad
\|u\|=\tilde{O}(e^{\frac{\tilde{c}_{n}-\tilde{c}_{m}}{h}})
\|\Pi_{G^{h}}d_{f,h}T_{\delta_{2}}u\|\,.
$$
Let us then assume that there exist $\varepsilon_{1}>0$\,, $1\leq m_{0}<n_{0}\leq N$\,, a strictly decreasing sequence $(h_{k})_{k\in\nz}$
converging to $0$ and, for every $k\in\nz$\,,
$u_{h_{k}}\in \mathcal{V}^{h_{k}}_{m_{0},n_{0}}\setminus\left\{0\right\}$ such
that 
\begin{equation}
\label{eq:estcontruh}
\|\Pi_{G^{h_{k}}}d_{f,f^{-1}([a,b]),h_{k}}T_{\delta_{2}}u_{h_{k}}\|\leq 
e^{-\frac{\tilde{c}_{n_{0}}-\tilde{c}_{m_{0}}+\varepsilon_{1}}{h_{k}}}\|u_{h_{k}}\|\,.
\end{equation}
Without restriction, we choose
the pair $(m_{0},n_{0})$ among the pairs
for which \eqref{eq:estcontruh} holds such that
 $\lambda_{0}:= \tilde{c}_{n_{0}}-\tilde{c}_{m_{0}}$
is minimal.
Set
$$
\ell:=\sharp\left\{(m,n)\in \mathcal{X}(a,b)\,, \tilde{c}_{n}-\tilde{c}_{m}\leq \lambda_{0}\right\}\,.
$$
Theorem~\ref{th:induc}-\textbf{a)} says that the $\ell$-th singular value of
$\delta_{[0,\tilde{o}(1)],[a,b],h}$ and therefore, with \eqref{eq:muGdT}, the $\ell$-th singular value
of  $\Pi_{G^{h}}d_{f,f^{-1}[a,b],h}T_{\delta_{2}}\big|_{\mathcal{V}^{h}}$ satisfy
$$
\lim_{h\to 0}-h\log \mu_{\ell}(\Pi_{G^{h}}d_{f,f^{-1}([a,b]),h}T_{\delta_{2}}\big|_{\mathcal{V}^{h}})=
\lim_{h\to 0}-h\log\mu_{\ell}(\delta_{[0,\tilde{o}(1)],[a,b],h})=\lambda_{0}\,.
$$
By using in addition the
$\tilde{O}(e^{-\frac{\delta_{1}}{h}})$-orthogonal decomposition
$$
\mathcal{V}^{h}=\mathcal{V}_{+}^{h}\oplus \mathcal{V}_{0}^{h}\qquad\text{with}
\quad
\mathcal{V}_{0}^{h}=\ker(\Pi_{G^{h}}d_{f,f^{-1}([a,b]),h}T_{\delta_{2}}\big|_{\mathcal{V}^{h}})\,,
$$
applying Proposition~\ref{pr:epsorth}-{\bf b)} gives
$$
\lim_{h\to 0}-h\log \mu_{\ell}(\Pi_{G^{h}}d_{f,f^{-1}([a,b]),h}T_{\delta_{2}}\big|_{\mathcal{V}_{+}^{h}})=
\lim_{h\to 0}-h\log\mu_{\ell}(\delta_{[0,\tilde{o}(1)],[a,b],h})=\lambda_{0}\,.
$$
Because $\mathcal{V}_{+}^{h}$ is finite dimensional, $\dim \mathcal{V}_{+}^{h}=\sharp
\mathcal{X}(a,b)$\,, the max-min principle implies
$$
\mu_{\ell}(\Pi_{G^{h}}d_{f,h}T_{\delta_{2}}\big|_{\mathcal{V}_{+}^{h}})
=\min_{\dim W=\sharp \mathcal{X}(a,b)-\ell+1}\max_{v\in W\setminus\left\{0\right\}}\frac{\|\Pi_{G^{h}}d_{f,h}T_{\delta_{2}}v\|}{\|v\|}\,.
$$
We obtain a contradiction by considering
$$
W=\left(\mathop{\oplus}_{\tilde{c}_{n}-\tilde{c}_{m}>\lambda_{0}}\mathcal{V}_{m,n}^{h_{k}}\right)
\oplus \cz u_{h_{k}}\,.
$$
This ends the proof of the first statement.\\
By applying again Proposition~\ref{pr:epsorth}-\textbf{b)}  with now
$B=\Pi_{G^{h}}d_{f,h}T_{\delta_{2}}$ acting on $\mathcal{V}_{+}^{h}$\,,  the singular
values of $\Pi_{G^{h}}d_{f,h}T_{\delta_{2}}\big|_{\mathcal{V}_{+}^{h}}$ are
obtained, modulo
some $\tilde{O}(e^{-\frac{\delta_{1}}{h}})$ relative error, 
by collecting all the singular values of
$\Pi_{G_{n}^{h}}d_{f,h}T_{\delta_{2}}\big|_{\mathcal{V}_{n}^{h}}$\,, $n\in
\left\{2,\ldots, N\right\}$\,. Actually, $G_{n}^{h}\perp G_{n'}^{h}$ and
$\mathcal{V}_{n}^{h}$ and $\mathcal{V}_{n'}^{h}$ are
$\tilde{O}(e^{-\frac{\delta_{1}}{h}})$-orthogonal for $n\neq
n'$\,. This reduces the problem to
the computation of the singular values of
 $\Pi_{G_{n}^{h}}d_{f,h}T_{\delta_{2}}\big|_{\mathcal{V}_{n}^{h}}$\,.
For $n\in \left\{2,\ldots,N\right\}$\,, we solve it by reverse
induction  on
$m\in \left\{1,\ldots, n-1\right\}$ by considering 
$\mathop{\oplus}_{m\leq m'<n}\mathcal{V}_{m',n}^{h}$\,.
 Simply apply  Proposition~\ref{pr:epsorth}-\textbf{c)} with 
\begin{eqnarray*}
  && {E'}^{h}=\mathop{\oplus}_{m\leq
     m'<n}\mathcal{V}_{m,n}^{h}\,,\quad
 \mu_{\dim{E'}^{h}}(\Pi_{G_{n}^{h}}d_{f,h}T_{\delta_{2}}\big|_{{E'}^{h}})\stackrel{log}{\sim}e^{-\frac{\tilde{c}_{n}-\tilde{c}_{m}}{h}}\,,
\\
&&
{E''}^{h}=\mathcal{V}_{m-1,n}^{h} \,,\quad
\|\Pi_{G_{n}^{h}}d_{f,h}T_{\delta_{2}}\big|_{{E''}^{h}}\|=
\tilde{O}(e^{-\frac{\tilde{c}_{n}-\tilde{c}_{m-1}}{h}})\leq\tilde{O}(e^{-\frac{\tilde{c}_{n}-\tilde{c}_{m}+2\eta_{f}}{h}})\leq\tilde{O}(e^{-\frac{\tilde{c}_{n}-\tilde{c}_{m}+\delta_{1}}{h}})\,,
\end{eqnarray*}
by starting from the first case when $\dim {E'}^{h}\neq 0$\,. This
implies that the non zero singular values of
$\Pi_{G_{n}^{h}}d_{f,h}T_{\delta_{2}}:\mathop{\oplus}_{m-1\leq
  m'<n}\mathcal{V}_{m',n}^{h}$ are obtained, modulo
some $\tilde{O}(e^{-\frac{\delta_{1}}{h}})$ relative error,  by collecting the non zero
singular values of
$\Pi_{G_{n}^{h}}d_{f,h}T_{\delta_{2}}\big|_{\mathcal{V}_{m',n}^{h}}$
for $m-1\leq m'<n$\,.
This ends the proof of the second statement.
\end{proof}
\begin{prop}
\label{pr:diagVh}
Assume that Theorem~\ref{th:induc} holds true at step $N$ and let
$G^{h}$ be defined by \eqref{eq:defGh}.
There exists a basis $(\phi_{j}^{h})_{j\in
  \mathcal{J}(a,b)}$ of $\mathcal{V}^{h}$ such that the
$\phi_{j}^{h}$'s 
satisfy the same properties as the $\varphi_{j}^{h}$'s, that is the ones of
Definition~\ref{de:adapted} and  of Theorem~\ref{th:induc}, 
as well as the additional following one: 
\begin{eqnarray}
\label{eq:Psiorth}
 &&\Psi_{j}^{h}\perp \Psi_{j'}^{h}\quad\text{for}~j\neq j'\\
\text{where}&&
\label{eq:defPsi}
\Psi_{j}^{h}=\Pi_{G^{h}}d_{f,h}T_{\delta_{2}}\phi_{j}^{h}\,.
\end{eqnarray}
In particular,
according to Proposition~\ref{pr:Vmn},
 the singular values of
$\Pi_{G^{h}}d_{f,h}T_{\delta_{2}}:\mathcal{V}^{h}\to L^{2}(\Omega)$
are given by the numbers
$\|\Psi_{j}^{h}\|_{L^{2}}(1+\tilde{O}(e^{-\frac{\delta_{1}}{h}}))$\,,
$j\in \mathcal{J}(a,b)$\,, where
$\|\Psi_{j}^{h}\|_{L^{2}}\stackrel{log}{\sim}e^{-\frac{y_{\alpha}-x_{\alpha}}{h}}$
when $j=(\alpha,x_{\alpha})\in \mathcal{X}(a,b)$\,.
\end{prop}
\begin{proof}
We keep $\phi_{j}^{h}=\varphi_{j}^{h}$ if $j\in
\mathcal{Y}(a,b)\cup\mathcal{Z}(a,b)$\,. Because $G_{n}^{h}\perp
G_{n'}^{h}$ for $n\neq n'$ and
$\Pi_{G^{h}}d_{f,h}T_{\delta_{2}}\mathcal{V}_{n}^{h}\subset G_{n}^{h}$ for
$2\leq n\leq N$\,, it suffices to construct the family 
$(\phi_{j}^{h})_{j\in \mathcal{X}_{n}(a,b)}$ for any $n\in
\left\{2,\ldots, N\right\}$\,.
Take some fixed $n\in \left\{2,\ldots, N\right\}$\,. 
While keeping the $\tilde{O}(e^{-\frac{\delta_{1}}{h}})$-orthogonal
decomposition 
$$
\mathcal{V}_{n}^{h}=\mathop{\oplus_{1\leq m<n}}\mathcal{V}_{m,n}^{h}\,,
$$
the first result of Proposition~\ref{pr:Vmn} says that, for a fixed
pair $(m,n)$\,, the
$\tilde{O}(e^{-\frac{\delta_{1}}{h}})$-orthonormal basis
$(\varphi_{j}^{h})_{j\in \mathcal{X}_{m,n}(a,b)}$ can be replaced by
an orthonormal one $(\tilde{\varphi}_{j}^{h})_{j\in
  \mathcal{X}_{m,n}(a,b)}$ such that 
\begin{eqnarray*}
  &&
\Pi_{G_{n}^{h}}d_{f,h}T_{\delta_{2}}\tilde{\varphi}_{j}^{h}\perp 
\Pi_{G_{n}^{h}}d_{f,h}T_{\delta_{2}}\tilde{\varphi}_{j'}^{h}\quad
\text{for}~j\neq j',\  j,j'\in \mathcal{X}_{m,n}(a,b)\\
\text{and}&&
\|\Pi_{G_{n}^{h}}d_{f,h}T_{\delta_{2}}\tilde{\varphi}_{j}^{h}\|\stackrel{log}{\sim}e^{-\frac{\tilde{c}_{n}-\tilde{c}_{m}}{h}}\quad\text{for}~j\in
             \mathcal{X}_{m,n}(a,b)\,.
\end{eqnarray*}
Because the change of basis $P_{m,n}^{h}\in
\mathcal{L}(\mathcal{V}_{m,n}^{h})$ given by 
 $\tilde{\varphi}_{j}^{h}=P_{m,n}^{h}\varphi_{j}^{h}$ satisfies
$$
\|(P_{m,n}^{h})^{*}P_{m,n}^{h}-\Id_{\mathcal{V}_{m,n}}^{h}\|=\tilde{O}(e^{-\frac{\delta_{1}}{h}})\,,
$$
the new family  $(\tilde{\varphi}_{j}^{h})_{j\in
  \mathcal{X}_{m,n}(a,b)}$ keeps all the properties of the initial one
$(\varphi_{j}^{h})_{j\in \mathcal{X}_{m,n}(a,b)}$\,. In
Theorem~\ref{th:induc} at step $N$\,, nothing is changed when the
$\varphi_{j}^{h}$\,, $j\in \mathcal{X}_{m,n}(a,b)$\,, are replaced by
the $\tilde{\varphi}_{j}^{h}$\,, $j\in \mathcal{X}_{m,n}(a,b)$\,, and
this can be done for all pairs $(m,n)$ and with any
initial guess of the family $(\varphi_{j}^{h})_{j\in \mathcal{J}(a,b)}$\,.\\
Thus, it suffices to construct the family $(\phi_{j}^{h})_{j\in
  \mathcal{X}_{n}(a,b)}$ such that \eqref{eq:Psiorth} and \eqref{eq:defPsi}
hold when $j\in \mathcal{X}_{m_{1},n}(a,b)$\,, $j'\in
\mathcal{X}_{m_{2},n}(a,b)$\,, $m_{1}\neq m_{2}$\,.
Like at the end of the previous proof,
we do it by reverse induction on $m\in
\left\{1,\ldots,n-1\right\}$\,.
\begin{itemize}
\item For $m=n-1$\,, simply take $\phi_{j}^{h}=\tilde{\varphi}_{j}^{h}$
  and  set $$\mathcal{W}_{n-1,n}^{h}=\Vect(\phi_{j}^{h},h\in \mathcal{X}_{n-1,n}(a,b))=\mathcal{V}_{n-1,n}^{h}\,.$$
\item Assume that the $\phi_{j}^{h}$'s have been constructed for 
 $j\in \mathcal{X}_{m',n}(a,b)$\,, for all $m'\in
 \left\{m,\ldots,n-1\right\}$\,, with 
 $\mathcal{W}_{m',n}^{h}=\Vect(\phi_{j}^{h}, j\in
 \mathcal{X}_{m',n}(a,b))$ and the equality of the
 $\tilde{O}(e^{-\frac{\delta_{1}}{h}})$-orthogonal decompositions
$$
\mathop{\oplus}_{m\leq m'<n}\mathcal{W}_{m',n}^{h}
=\mathop{\oplus}_{m\leq m'<n}\mathcal{V}_{m',n}^{h}\,.
$$
Set, for  $j\in \mathcal{X}_{m-1, n}(a,b)$\,,
$$
\hat{\varphi}_{j}^{h}=\varphi_{j}^{h}-
\sum_{j'\in\mathop{\sqcup}_{m\leq m'<n}\mathcal{X}_{m',n}(a,b)}
\frac{\langle \Psi_{j'}^{h}\,,
  \Pi_{G_{n}^{h}}d_{f,h}T_{\delta_{2}}\varphi_{j}^{h}\rangle}{\|\Psi_{j'}^{h}\|^{2}}\phi_{j'}^{h}\,,
$$
and define
$$
\mathcal{W}_{m-1,n}^{h}:=\Vect(\hat{\varphi}_{j}^{h}\,, j\in \mathcal{X}_{m-1,n}(a,b))\,.
$$
We have clearly
$$
\Pi_{G_{n}^{h}}d_{f,h}T_{\delta_{2}}\hat{\varphi}_{j}^{h}\perp
\Vect(\Psi_{j'}^{h}, j'\in \mathcal{X}_{m',n}(a,b), m\leq m'<n)
$$
and
 $$
\mathop{\oplus}_{m-1\leq m'< n}\mathcal{V}^{h}_{m',n}
=\mathcal{W}_{m-1,n}^{h}\oplus \left(\mathop{\oplus}_{m\leq
    m'<n}\mathcal{V}_{m',n}^{h}\right)\,.
$$
All the properties of Theorem~\ref{th:induc} at step $N$ are
verified for the $\delta_{1}$-family of quasimodes given by
the $\hat{\varphi}_{j}^{h}$\,, $j\in
\mathcal{X}_{m-1,n}(a,b)$\,, and the $\phi_{j}^{h}$\,,
 $j\in \mathcal{X}_{m',n}(a,b)$\,.
The estimates on $\hat{\varphi}_{j}^{h}$\,, $j\in
\mathcal{X}_{m-1,n}(a,b)$\,, are consequences of:
\begin{eqnarray*}
  && \Psi_{j'}^{h}=\Pi_{G_{n}^{h}}d_{f,h}\phi_{j'}^{h}\in
     \ker(\Delta_{f,f^{-1}([\tilde{c}_{n}-\eta_{f},\tilde{c}_{n}+\eta_{f}]\cap[a,b]),h})\quad\text{for}~j'\in
    \mathop{\sqcup}_{m\leq m'<n}\mathcal{X}_{ m',n}(a,b)\,,\\
\text{where}&&
\|\Psi_{j'}^{h}\|_{L^{2}}\stackrel{log}{\sim}e^{-\frac{\tilde{c}_{n}-\tilde{c}_{m'}}{h}}\quad\text{when}~j'\in
              \mathcal{X}_{m',n}(a,b)\\
\text{and}&& 
\frac{\langle \Psi_{j'}^{h}\,,
  \Pi_{G_{n}^{h}}d_{f,h}T_{\delta_{2}}\varphi_{j}^{h}\rangle}{\|\Psi_{j'}^{h}\|^{2}}=
\tilde{O}(e^{\frac{\tilde{c}_{n}-\tilde{c}_{m'}}{h}})
\times
             \tilde{O}(e^{-\frac{\tilde{c}_{n}-\tilde{c}_{m-1}}{h}})=\tilde{O}(e^{-\frac{\tilde{c}_{m'}-\tilde{c}_{m-1}}{h}})\,,\\
&&
\|e^{\frac{|f-\tilde{c}_{m'}|}{h}}\phi_{j'}^{h}
\|_{W(f^{-1}([a,b])\setminus S_{\delta_{1}})}
=\tilde{O}(1)\,,\\
&&\tilde{c}_{m'}-\tilde{c}_{m-1}\geq 2\eta_{f}\geq \delta_{1}\,.
\end{eqnarray*}
Hence, the vectors $\hat{\varphi}_{j}^{h}$\,, $j\in \mathcal{X}_{m-1,n}(a,b)$\,, satisfy
$$
\|e^{\frac{|f-\tilde{c}_{m-1}|}{h}}\hat{\varphi}_{j}^{h}
\|_{W(f^{-1}([a,b])\setminus S_{\delta_{1}})}
=\tilde{O}(1)\,.
$$
Note in particular that the total space $\mathcal{V}^{h}$ is not
changed, so the statement of Theorem~\ref{th:induc}-{\bf b)} and the factorization in Theorem~\ref{th:induc}-{\bf c)} are
obviously true.\\
Once we have the $\tilde{O}(e^{-\frac{\delta_{1}}{h}})$-orthogonal
decomposition
$$
\mathcal{V}_{+}^{h}=\left(\mathop{\oplus}_{m-1\leq
    m'<n}\mathcal{W}_{m',n}^{h}\right)\oplus 
\left(\mathop{\oplus}_{1\leq
    m'<m-2}\mathcal{V}_{m',n}^{h}\right)\,,
$$
we just apply our first argument with 
$\mathcal{V}_{m-1,n}^{h}$
now replaced by
$\mathcal{W}_{m-1,n}^{h}$\,,
which permits to replace
  the
$\tilde{O}(e^{-\frac{\delta_{1}}{h}})$-orthonormal basis
$(\hat{\varphi}_{j}^{h})_{j\in \mathcal{X}_{m-1,n}(a,b)}$ of
$\mathcal{W}_{m-1,n}^{h}$ by  an orthonormal basis 
$(\tilde{\varphi}_{j}^{h})_{j\in \mathcal{X}_{m-1,n}(a,b)}$ such that
$$
\Pi_{G_{n}^{h}}d_{f,h}T_{\delta_{2}}\tilde{\varphi}_{j}^{h}\perp 
\Pi_{G_{n}^{h}}d_{f,h}T_{\delta_{2}}\tilde{\varphi}_{j'}^{h}\quad
\text{for}~j\neq j',\  j,j'\in \mathcal{X}_{m-1,n}(a,b)\,.
$$
We finally define $\phi_{j}^{h}=\tilde{\varphi}_{j}^{h}$ for $j\in \mathcal{X}_{m-1,n}(a,b)$\,.
\end{itemize}
\end{proof}

\subsection{$N\to N+1$: Collecting the information from
  step $N$}
\label{sec:collect}
We assume that Theorem~\ref{th:induc} holds at step $N$\,, i.e.
when $\sharp([a,b]\cap \left\{c_{1},\ldots,c_{N_{f}}\right\})=N$\,,
and we  consider the case
$$
[a,b]\cap\left\{c_{1},\ldots,
  c_{N_{f}}\right\}=\left\{\tilde{c}_{1},\ldots,\tilde{c}_{N+1}\right\}\,. 
$$
Define
$$
a_{1}=a\quad,\quad b_{1}=\tilde{c}_{N}+\eta_{f}\quad
\text{and}
\quad
a_{2}=\tilde{c}_{2}-\eta_{f}\quad,\quad b_{2}=b\,.
$$
We can use Theorem~\ref{th:induc} and its consequences given in
Subsection~\ref{sec:conseqN}
for
$$
[a_{1},b_{1}]\cap
\left\{c_{1},\ldots,c_{N_{f}}\right\}=\left\{\tilde{c}_{1},\ldots,
  \tilde{c}_{N}\right\}\quad\text{and}\quad
[a_{2},b_{2}]\cap \left\{c_{1},\ldots,c_{N_{f}}\right\}=\left\{\tilde{c}_{2},\ldots,
  \tilde{c}_{N+1}\right\}\,.
$$
Let us start with the interval  $[a_{2},b_{2}]$\,. Consider
$\Delta_{f,\overline{\Omega}_{2},h}$ and let $G_{2}^{h}$ and
$F_{2}^{h}$ be defined like  $G^{h}$ and $F^{h}$ in 
\eqref{eq:defOmbar} and \eqref{eq:defGh} while replacing $(a,b)$ by
$(a_{2},b_{2})$\,, with
$$
G_{2}^{h}=\mathop{\oplus}_{2\leq n\leq N+1}^{\perp}G_{2,n}^{h}
=\mathop{\oplus}_{2\leq n\leq
  N+1}\ker(\Delta_{f,f^{-1}([\tilde{c}_{n}-\eta_{f},\tilde{c}_{n}+\eta_{f}]\cap[a,b]),h})\,.
$$
For this interval $[a_{2},b_{2}]$\,, the family of quasimodes
$(\phi_{2,j}^{h})_{j\in \mathcal{J}(a_{2},b_{2})}$ is given by
Proposition~\ref{pr:diagVh} with the orthogonality condition
\eqref{eq:Psiorth},\eqref{eq:defPsi}, and we set
$$
\mathcal{W}^{h}_{m,n}(a_{2},b_{2})=\Vect(\phi_{2,j}^{h}\,, j\in
\mathcal{X}_{m,n}(a_{2},b_{2}))\,,\quad 2\leq m<n\leq N+1\,.
$$
For the interval $[a_{1},b_{1}]$\,, we use similar notations
$\Delta_{f,\overline{\Omega}_{1},h}$\,, $G_{1}^{h}$\,,
$F_{1}^{h}$
with now
$$
G_{1}^{h}=\mathop{\oplus}_{1\leq n\leq N}^{\perp}G_{1,n}^{h}
=\mathop{\oplus}_{1\leq n\leq
  N}\ker(\Delta_{f,f^{-1}([\tilde{c}_{n}-\eta_{f},\tilde{c}_{n}+\eta_{f}]\cap[a,b]),h})\,.
$$
We start with a family of
quasimodes 
\begin{equation}
\label{eq:phi0j}
 (\varphi_{0,j}^{h})_{j\in \mathcal{J}(a_{1},b_{1})} 
\end{equation}
given by Theorem~\ref{th:induc}
and
merge this family with $(\phi_{2,j}^{h})_{j\in
  \mathcal{J}(a_{2},b_{2})\cap \mathcal{J}(a_{1},b_{1})}$\,, after considering the
restrictions $\phi_{2,j}\big|_{f^{-1}([a_{2},b_{1}])}$ 
extended by $0$ in $f_{a_{1}=a}^{a_{2}}$\,,
according to the following procedure:
\begin{eqnarray*}
  && \varphi_{1,j}^{h}= \phi_{2,j}^{h}
\quad
\text{if}~j\in (\alpha,\tilde{c})\in \mathcal{J}(a_{1},b_{1})\,,\ 
      \tilde{c}\geq \tilde{c}_{2}\,,
\\
&&\varphi_{1,j}^{h}= \varphi_{0,j}^{h}
\quad
\text{if}~j\in (\alpha,\tilde{c}_{1})\in \mathcal{Z}(a_{1},b_{1})\,,
\\
&&\varphi_{1,j}^{h}= \varphi_{0,j}^{h}-
\ds\sum_{j'\in\mathcal{X}(a_{2},b_{2})\cap \mathcal{X}(a_{1},b_{1})}
\frac{\langle \Psi_{2,j'}^{h}\,,
  \Pi_{G_{2}^{h}}d_{f,h}T_{\delta_{2}}\varphi_{0,j}^{h}\rangle}{\|\Psi_{2,j'}^{h}\|^{2}}\phi_{2,j'}^{h}\quad\text{if}~j=(\alpha,\tilde{c}_{1})\in
   \mathcal{X}(a_{1},b_{1})\,,
\\
&&\text{with}\quad
\Psi^{h}_{2,j'}=\Pi_{G_{2}^{h}}d_{f,h}T_{\delta_{2}}\phi_{2,j'}^{h}=\Pi_{G_{1}^{h}}d_{f,h}T_{\delta_{2}}\varphi_{1,j'}^{h}\quad\text{for}~j'\in
   \mathcal{X}(a_{2},b_{2})\cap \mathcal{X}(a_{1},b_{1})\,,
\end{eqnarray*}
where we recall that
 $j=(\alpha,\tilde{c})\in \mathcal{X}(a_{2},b_{2})\cap
\mathcal{X}(a_{1},b_{1})$ means $\tilde{c}_{2}\leq
x_{\alpha}<y_{\alpha}\leq \tilde{c}_{N}$\,.

\begin{remark}
\label{re.gamma-id}
Assume that $\gamma_{1}(h), (\varphi_{1,j}^h)_{j\in \mathcal{J}(a_{1},b_{1}) }$ and
$\gamma_{2}(h), (\varphi_{2,j}^h)_{j\in \mathcal{J}(a_{2},b_{2}) }$ are
given 
by Theorem~\ref{th:induc} and Definition~\ref{de:adapted} at step
$N$\,, respectively in $[a_{1},b_{1}]$ 
and in $[a_{2},b_{2}]$\,.
Let us then define  
 $\gamma(h):=\max (\gamma_{1}(h),\gamma_{2}(h))$
 and, for $i\in\{1,2\}$\,,
$$
\tilde{\varphi}_{i,j}^{h}:=\begin{cases}
\varphi_{i,j}^{h} & \text{when $j\in \mathcal{Y}(a_{i},b_{i})\cup  \mathcal{Z}(a_{i},b_{i})$}\,,
\\
\chi_{y_{\alpha}^{(p+1)},\gamma(h)}\varphi_{i,j}^{h}
\quad & \text{when $j=(\alpha,x_{\alpha}^{(p)})\in \mathcal{X}^{(p)}(a_{i},b_{i})$\,, $p\in\{0,\dots,N-1\}$}\,,
\end{cases}
$$
where $\chi_{y_{\alpha}^{(p+1)},\gamma(h)}$ is defined by \eqref{eq:defchic} in 
Definition~\ref{de:Td2}.
Then, the families $(\tilde \varphi_{1,j}^h)_{j\in \mathcal{J}(a_{1},b_{1}) }$
and $(\tilde\varphi_{2,j}^h)_{j\in \mathcal{J}(a_{2},b_{2}) }$
both satisfy the properties of 
Theorem~\ref{th:induc} and Definition~\ref{de:adapted}, respectively 
 in $[a_{1},b_{1}]$ 
and in $[a_{2},b_{2}]$\,, but now with the same $\gamma(h)$\,.
Hence, we will  assume here that the properties of the families
$(\phi_{2,j}^{h})_{j\in \mathcal{J}(a_{2},b_{2})}$ and
$(\varphi_{0,j}^{h})_{j\in \mathcal{J}(a_{1},b_{1})}$
are satisfied
with the same $\gamma(h)$\,.
\end{remark}

The spaces generated by those quasimodes are denoted by
$$
\mathcal{V}^{h}(a_{1},b_{1})=\Vect(\varphi_{1,j}^{h},j\in
\mathcal{J}(a_{1},b_{1}))\quad\text{and}\quad
\mathcal{V}^{h}(a_{2},b_{2})=\Vect(\phi_{2,j}^{h},j\in
\mathcal{J}(a_{2},b_{2}))\,,
$$
and the same rule applies for $\mathcal{V}_{m,n}^{h}$\,,
$\mathcal{V}_{n}^{h}$\,, $1\leq m<n\leq N+1$\,, 
$\mathcal{V}_{+}^{h}$\,, $\mathcal{V}_{0}^{h}$ defined in
\eqref{eq:defVmn}--\eqref{eq:defV0},
while
writing $\mathcal{W}_{m,n}^{h}(a_{2},b_{2})$ instead of ${\cal
  V}_{m,n}^{h}(a_{2},b_{2})$
refers to the additional orthogonality property of
Proposition~\ref{pr:diagVh}.
\begin{prop}
\label{pr:collG1}
The family $(\varphi_{1,j}^{h})_{j\in \mathcal{J}(a_{1},b_{1})}$
satisfies all the properties of Theorem~\ref{th:induc} at step
$N$\,. 
Moreover, the family $(\phi_{1,j})_{j\in \mathcal{J}(a_{1},b_{1})}$
deduced from $(\varphi_{1,j}^{h})_{j\in \mathcal{J}(a_{1},b_{1})}$ in
Proposition~\ref{pr:diagVh} can be constructed such that
$$
\forall j\in
\mathcal{X}(a_{1},b_{1})\cap \mathcal{X}(a_{2},b_{2})\quad,\quad
\phi_{1,j}^{h}=\phi_{2,j}^{h}\,.
$$ 
\end{prop}
\begin{proof}
By construction (and Remark~\ref{re.gamma-id}), the family $(\varphi^{h}_{1,j})_{j\in
  \mathcal{J}(a_{1},b_{1})}$ is a
$\tilde{O}(e^{-\frac{\delta_{1}}{h}})$-orthonormal
$\delta_{1}$-family of quasimodes,  and
$G_{1,n}^{h}=G_{2,n}^{h}$ for $2\leq n \leq N$ and
\begin{eqnarray*}
&&\vec{d}(\Vect(\varphi_{0,j}^{h},j\in
\mathcal{J}(a_{1},b_{1})),F_{1}^{h})+\vec{d}(F_{1}^{h},\Vect(\varphi_{0,j}^{h},j\in
\mathcal{J}(a_{1},b_{1}))=\tilde{O}(e^{-\frac{\delta_{1}}{h}})\,,\\
 &&\vec{d}(\mathcal{V}^{h}(a_{2},b_{2}),F_{2}^{h})+\vec{d}(F_{2}^{h},\mathcal{V}^{h}(a_{2},b_{2}))=\tilde{O}(e^{-\frac{\delta_{1}}{h}})\,,\\
&&\vec{d}(F_{1}^{h},G_{1}^{h})+\vec{d}(G_{1}^{h},F_{1}^{h})=\tilde{O}(e^{-\frac{\eta_{f}}{h}})\leq\tilde{O}(e^{-\frac{\delta_{1}}{h}})\,,\\
\text{and} &&\vec{d}(F_{2}^{h},G_{2}^{h})+\vec{d}(G_{2}^{h},F_{2}^{h})=\tilde{O}(e^{-\frac{\eta_{f}}{h}})\leq\tilde{O}(e^{-\frac{\delta_{1}}{h}})
\end{eqnarray*}
ensure the validity of the last statement of \textbf{b)} in Theorem~\ref{th:induc}, that is
$$
\vec{d}(\mathcal{V}^{h}(a_{1},b_{1}),F_{1}^{h})+\vec{d}(F_{1}^{h},\mathcal{V}^{h}(a_{1},b_{1}))=\tilde{O}(e^{-\frac{\delta_{1}}{h}})\,.
$$
 The
exponential decay estimates on the $\varphi_{1,j}^{h}$\,,
$j=(\alpha,\tilde{c}_{1})\in \mathcal{X}(a_{1},b_{1})$\,, are actually
 obtained
like in the proof of Proposition~\ref{pr:diagVh} by noticing that
$$
\forall j\in \mathcal{X}_{1,n}(a_{1},b_{1})\,,\quad
\varphi_{1,j}^{h}= \varphi_{0,j}^{h}-
\ds\sum_{j'\in\ds\mathop{\sqcup}_{2\leq m'<n\leq N}\mathcal{X}_{m',n}(a_{2},b_{2})}
\frac{\langle \Psi_{2,j'}^{h}\,,
  \Pi_{G_{2,n}^{h}}d_{f,h}T_{\delta_{2}}\varphi_{0,j}^{h}\rangle}{\|\Psi_{2,j'}^{h}\|^{2}}\phi_{2,j'}^{h}\,,
$$
where $G_{1,n}^{h}=G_{2,n}^{h}$ for $3\leq n\leq N$ and 
$$
\Psi_{2,j'}=\Pi_{G_{2,n}^{h}}d_{f,h}T_{\delta_{2}}\phi_{2,j'}^{h}=\Pi_{G_{1,n}^{h}}d_{f,h}T_{\delta_{2}}\varphi_{1,j'}^{h}\quad\text{for}~j'\in
\mathcal{X}_{m',n}(a_{1},b_{1})\,,\ 
2\leq m'<n\leq N\,.
$$
We still have to check the factorization of
Theorem~\ref{th:induc}-\textbf{c)}, namely
$$
d_{f,f^{-1}([a_{1},b_{1}])h}T_{\delta_{2}}\big|_{\mathcal{V}^{h}(a_{1},b_{1})}=C^{h}\Pi_{F_{1}^{h}}d_{f,f^{-1}([a_{1},b_{1}]),h}T_{\delta_{2}}\big|_{\mathcal{V}^{h}(a_{1},b_{1})}\quad \text{with}~\|C^{h}\|=\tilde{O}(e^{\frac{2\delta_{2}}{h}})\,.
$$
We will do it by first considering the operator
$\Pi_{G_{1}^{h}}d_{f,h}T_{\delta_{2}}$\,.\\
From the properties of the $\varphi_{1,j}^{h}$\,, $j\in
\mathcal{J}(a_{1},b_{1})$\,, we already  know that (see indeed \eqref{eq:estimupuh})
$$
\|\Pi_{G_{1}^{h}}d_{f,h}T_{\delta_{2}}\big|_{\mathcal{V}_{m,n}^{h}(a_{1},b_{1})}\|=\tilde{O}(e^{-\frac{\tilde{c}_{n}-\tilde{c}_{m}}{h}})
\quad\text{and}\quad \mathcal{V}_{0}^{h}(a_{1},b_{1})\subset \ker(d_{f,h}T_{\delta_{2}})\,.
$$
We now check  that $\Pi_{G_{1}^{h}}d_{f,h}T_{\delta_{2}}\big|_{\mathcal{V}_{m,n}^{h}(a_{1},b_{1})}$
is one to one 
and that its singular values, which thus do not vanish, all satisfy
$\mu_{h}\stackrel{log}{\sim}e^{-\frac{\tilde{c}_{n}-\tilde{c}_{m}}{h}}$
for every $1\leq m<n\leq N$ such that $\mathcal{X}_{m,n}(a_{1},b_{1})\neq \emptyset$:
\begin{itemize}
\item Since the vectors
  $\Psi_{2,j}^{h}=\Pi_{G_{1}^{h}}d_{f,h}T_{\delta_{2}}\phi_{2,j}^{h}=
  \Pi_{G_{2}^{h}}d_{f,h}T_{\delta_{2}}\phi_{2,j}^{h}$
  are, according to Proposition~\ref{pr:diagVh} applied in $[a_{2},b_{2}]$\,, mutually orthogonal with
  $\|\Psi_{2,j}^{h}\|\stackrel{log}{\sim}e^{-\frac{\tilde{c}_{n}-\tilde{c}_{m}}{h}}$
  when $j\in \mathcal{X}_{m,n}(a_{1},b_{1})$\,, $2\leq m<n\leq N$\,,
  the result holds for $m\geq 2$\,.
\item Case $m=1$:
as in the proof of Proposition~\ref{pr:Vmn}, assume
by reductio ad absurdum  that there exist  $2\leq n \leq N$\,, a strictly decreasing sequence 
$(h_{k})_{k\in\nz}$
converging to $0$ and, for every $k\in\nz$\,,
$u_{h_{k}}\in \mathcal{V}_{1,n}^{h_{k}}(a_{1},b_{1})\setminus\left\{0\right\}$ such
that 
\begin{equation*}
\|\Pi_{G^{h_{k}}_{1}}d_{f,h_{k}}T_{\delta_{2}}u_{h_{k}}\|=\tilde{o}(e^{-\frac{\tilde{c}_{n}-\tilde{c}_{1}}{h_{k}}})\|u_{h_{k}}\|\,,
\end{equation*}
and let
$n_{0}\in \{2,\ldots N\}$ be the smallest $n$ 
such that the above holds.
Consider then
\begin{eqnarray*}
  &&
E''^{h_{k}}=(\cz u_{h_{k}})\oplus \mathcal{V}_{0}^{h_{k}}(a_{1} ,b_{1})\oplus\left(
\mathop{\oplus}_{\tilde{c}_{n}-\tilde{c}_{m}>
     \tilde{c}_{n_{0}}-\tilde{c}_{1}}\mathcal{V}_{m,n}^{h_{k}}(a_{1},b_{1})\right)\,,\\
\text{so that}
&&\dim \mathcal{V}^{h_{k}}(a_{1},b_{1})-\dim(E''^{h_{k}})=
\sharp\left(\mathop{\sqcup}_{\tilde{c}_{n}-\tilde{c}_{m}\leq \tilde{c}_{n_{0}}-\tilde{c}_{1}}
\mathcal{X}_{m,n}(a_{1},b_{1})\right)-1
=:\ell_{0}-1\,.
\end{eqnarray*}
Owing to the exponential decay estimates on the quasimodes, we obtain
  \begin{equation}
  \label{eq.norm-dfh}
     \|d_{f,h_{k}}T_{\delta_{2}}\big|_{E''^{h_{k}}}\|=\tilde{O}(e^{-\frac{\tilde{c}_{n_{0}}-\tilde{c}_{1}-2\delta_{2}}{h_{k}}})
\end{equation}
and (see  \eqref{eq:estimupuh})
$$
\|\Pi_{G_{1}^{h_{k}}}d_{f,h_{k}}T_{\delta_{2}}\big|_{E''^{h_{k}}}\|=\tilde{o}(e^{-\frac{\tilde{c}_{n_{0}}-\tilde{c}_{1}}{h_{k}}})\,.
$$
Since moreover $\|\Pi_{F_{1}^{h}}-\Pi_{F_{1}^{h}}\Pi_{G_{1}^{h}}\|=\tilde{O}(e^{-\frac{\eta_{f}}{h}})$\,,
we deduce
$\|\Pi_{F_{1}^{h_{k}}}d_{f,h_{k}}T_{\delta_{2}}\big|_{E''^{h_{k}}}\|=\tilde{o}(e^{-\frac{\tilde{c}_{n_{0}}-\tilde{c}_{1}}{h_{k}}})$
and then, applying the max-min principle
as in the proof of Proposition~\ref{pr:Vmn}
with here $W=E''^{h_{k}}$\,,
$$
\mu_{\ell_{0}}(\Pi_{F_{1}^{h_{k}}}d_{f,h_{k}}T_{\delta_{2}} \big|_{\mathcal{V}^{h_{k}}(a_{1},b_{1})})=\tilde{o}(e^{-\frac{\tilde{c}_{n_{0}}-\tilde{c}_{1}}{h_{k}}})\,.$$
Hence, 
since  $T_{\delta_{2}}$
is $\tilde{O}(e^{-\frac{\delta_{1}}{h}}) $-unitary (see \eqref{eq.T-delta-unit})
and
 $\vec{d}(F_{1}^{h},T_{\delta_{2}}\mathcal{V}^{h}(a_{1},b_{1}))+\vec{d}(T_{\delta_{2}}\mathcal{V}^{h}(a_{1},b_{1}),F_{1}^{h})=\tilde{O}(e^{-\frac{\delta_{1}}{h}})$ (see \eqref{eq.T-delta-close}), it follows from
Proposition~\ref{pr:projepsort} that
$$
\mu_{\ell_{0}}(\Pi_{[0,\tilde{o}(1)],[a_{1},b_{1}],h_{k}}d_{f,f^{-1}([a_{1},b_{1}]),h_{k}})=\tilde{o}(e^{-\frac{\tilde{c}_{n_{0}}-\tilde{c}_{1}}{h_{k}}})\;\text{with}\;
\ell_{0}=
\sharp\left(\mathop{\sqcup}_{\tilde{c}_{n}-\tilde{c}_{m}\leq \tilde{c}_{n_{0}}-\tilde{c}_{1}}
\mathcal{X}_{m,n}(a_{1},b_{1})\right),
$$
in contradiction with Theorem~\ref{th:induc}-\textbf{a)} in $[a_{1},b_{1}]$\,.
\end{itemize}
Because the spaces $\mathcal{V}_{m,n}^{h}(a_{1},b_{1})$
have mutually orthogonal images by $\Pi_{G_{1}^{h}}d_{f,h}T_{\delta_{2}}$\,, i.e.
$$
\Pi_{G_{1}^{h}}d_{f,h}T_{\delta_{2}}\mathcal{V}_{m_{1},n_{1}}^{h}(a_{1},b_{1})
\perp
\Pi_{G_{1}^{h}}d_{f,h}T_{\delta_{2}}\mathcal{V}_{m_{2},n_{2}}^{h}(a_{1},b_{1})\quad\text{for}~(m_{1},n_{1})\neq (m_{2},n_{2})\,,
$$
we can conclude like at the end of the proof of
Proposition~\ref{pr:diagVh} that there exists a basis
$(\phi_{1,j}^{h})_{j\in \mathcal{J}(a_{1},b_{1})}$
such that \eqref{eq:Psiorth} and
\eqref{eq:defPsi} hold, and
in which nothing needs
to be changed when $j\in \mathcal{X}(a_{2},b_{2})$\,.\\
It follows from the above analysis that $\Pi_{G_{1}^{h}}d_{f,h}T_{\delta_{2}}\big|_{\mathcal{V}_{+}^{h}(a_{1},b_{1})}$
is one to one, and
the factorization
$d_{f,h}T_{\delta_{2}}=\tilde{C}^{h}\Pi_{G_{1}^{h}}d_{f,h}T_{\delta_{2}}$ is then satisfied
with $\tilde{C}^{h}:G_{1}^{h}\to L^{2}(f_{a_{1}}^{b_{1}}) $ defined by
$\tilde{C}^{h} =0$ on the orthogonal complement of $\Pi_{G_{1}^{h}}d_{f,h}T_{\delta_{2}}\big(\mathcal{V}_{+}^{h}(a_{1},b_{1})\big)$
in $G_{1}^{h}$
and
$$
\forall j\in \mathcal{X}(a_{1},b_{1})\,, \quad
\tilde{C}^{h}\Psi_{1,j}^{h}=d_{f,h}T_{\delta_{2}}\phi_{1,j}^{h}\,.
$$
Moreover, the relation $\|\tilde{C}^{h}\|=\tilde{O}(e^{\frac{2\delta_{2}}{h}})$ follows  from
the orthogonality of the family $(\Psi_{1,j}^{h})_{j\in \mathcal{X}(a_{1},b_{1})}$ and from
$\|\Psi_{1,j}^{h}\|\stackrel{log}{\sim}e^{\frac{y_{\alpha}-x_{\alpha}}{h}}$
and
$\|d_{f,h}T_{\delta_{2}}\phi_{1,j}^{h}\|=\tilde{O}(e^{-\frac{y_{\alpha}-x_{\alpha}-2\delta_{2}}{h}})$
for $j=(\alpha,\tilde{c})\in \mathcal{X}(a_{1},b_{1})$ (see \eqref{eq.norm-dfh}).\\
Finally, applying the symmetric version
of 
Lemma~\ref{le:factorization}, that is exchanging $F^{h}$ and
$G^{h}$\,,
 yields the
 factorization
$d_{f,h}T_{\delta_{2}}=C^{h}\Pi_{F_{1}^{h}}d_{f,h}T_{\delta_{2}}:\mathcal{V}^{h}(a_{1},b_{1})\to
L^{2}(f_{a_{1}}^{b_{1}})$ stated in
Theorem~\ref{th:induc}-\textbf{c)}. 
\end{proof}
We have now spaces $\mathcal{W}_{m,n}^{h}(a_{1},b_{1})$\,, $1\leq
m<n\leq N$\,, and $\mathcal{W}_{m,n}^{h}(a_{2},b_{2})$\,,  $2\leq m<n\leq
N+1$\,, such that
$$
\mathcal{W}_{m,n}^{h}(a_{1},b_{1})=\mathcal{W}_{m,n}^{h}(a_{2},b_{2})\quad
\text{when}\quad2\leq m<n\leq N\,.
$$
We now work in the interval $[a,b]$ and we 
consider $\Delta_{f,\overline{\Omega},h}$\,, $G^{h}$\,, and $F^{h}$
according to \eqref{eq:defOmbar} and \eqref{eq:defGh}, after replacing
$N$ by $N+1$ and $\left\{\tilde{c}_{1},\ldots,\tilde{c}_{N}\right\}$ by
$\left\{\tilde{c}_{1},\ldots,\tilde{c}_{N+1}\right\}$\,. We set
\begin{eqnarray}
\label{eq:wmnab}
  && \hspace{-1.5cm}\mathcal{W}_{m,n}^{h}(a,b)=\left\{
     \begin{array}[c]{ll}
       \mathcal{W}_{m,n}^{h}(a_{2},b_{2})=\Vect(\phi_{2,j}^{h},
       j\in \mathcal{X}_{m,n}(a_{2},b_{2}))&\text{for}~2\leq m<n\leq
                                             N+1\,,\\
 \mathcal{W}_{1,n}^{h}(a_{1},b_{1})=\Vect(\phi_{1,j}^{h}, j\in
       \mathcal{X}_{1,n}(a_{1},b_{1}))&\text{for}~1=m <n\leq N\,,
     \end{array}
\right.\\
\label{eq:V0prim}
&&\mathcal{V}_{0}'^{h}(a,b)=\Vect(\phi^{h}_{2,j}\,, j\in
   \mathcal{Y}(a_{2},b_{2})\cup \mathcal{Z}(a_{2},b_{2}))\,,
\\
\label{eq:Vprim}
\text{and}
&& \mathcal{V}'^{h}(a,b)= \underbrace{(\mathop{\oplus}_{0<n-m\leq
   N-1}\mathcal{W}_{m,n}^{h}(a,b))}_{\mathcal{V}'^{h}_{+}}
\oplus \mathcal{V}_{0}'^{h}\,.
\end{eqnarray}
Accordingly, we introduce
\begin{eqnarray}
\label{eq:Jprim+}
  &&\mathcal{J}_{+}'(a,b)=\mathcal{X}(a_{2},b_{2})\sqcup (\mathop{\sqcup}_{2<n\leq N}
\mathcal{X}_{1,n}(a_{1},b_{1}))
=\mathop{\sqcup}_{0<n-m\leq N-1}\mathcal{X}_{m,n}(a,b)\,,\\
\label{eq:Jprim0}
&&
\mathcal{J}'_{0}(a,b)=
\mathcal{Y}(a_{2},b_{2})\sqcup \mathcal{Z}(a_{2},b_{2})
\quad\text{and}\quad \mathcal{J}'(a,b)=\mathcal{J}'_{+}(a,b)\sqcup
   \mathcal{J}'_{0}(a,b)\,,\\
   \label{eq.part-basis}
&& \varphi_{j}^{h}=\phi_{j}^{h}=
\left\{
   \begin{array}[c]{ll}
     \phi_{2,j}^{h}&\text{if}~j=(\alpha,\tilde{c})\in
                     \mathcal{J}_{+}'(a,b)\,, \tilde{c}_{2}\leq
                     \tilde{c}\,,
\\
    \phi_{2,j}^{h}&\text{if}~j\in\mathcal{J}_{0}'(a,b)\,,
\\
\phi_{1,j}^{h}&\text{if}~j=(\alpha,\tilde{c}_{1})\in \mathcal{J}'_{+}(a,b)\,.
   \end{array}
\right.
\end{eqnarray}
In  the perspective of applying
Proposition~\ref{pr:multadd}, we now consider the space $E'^{h}=T_{\delta_{2}}\mathcal{V}'^{h}$\,. 
\begin{prop}
\label{pr:Eprim}
With the notation \eqref{eq:Vprim}, consider
$E'^{h}=T_{\delta_{2}}\mathcal{V}'^{h}(a,b)$\,,
$E'^{h}_{0}=T_{\delta_{2}}\mathcal{V}'^{h}_{0}(a,b)$\,, and let $G^{h}$ be
defined by \eqref{eq:defGh} with $N$ replaced by $N+1$\,.
The operator $\Pi_{G^{h}}d_{f,f^{-1}([a,b])h}\big|_{E'^{h}}$ satisfies
\begin{eqnarray*}
  &&\mathrm{rank}\,(\Pi_{G^{h}}d_{f,h}\big|_{E'^{h}})=\sharp \mathcal{J}'_{+}(a,b)=:\ell_{1}\\
\text{and} &&\ker(\Pi_{G^{h}}d_{f,f^{-1}([a,b])h}\big|_{E'^{h}})=E'^{h}_{0}\,,
\end{eqnarray*}
and its non zero singular values can be written $(\mu_{j}^{h})_{j\in \mathcal{J}'_{+}(a,b)}$ with
$$
\mu_{j}^{h}\stackrel{log}{\sim}e^{-\frac{y_{\alpha}-x_{\alpha}}{h}}\quad
\quad\text{for every}~j=(\alpha,x_{\alpha})\in \mathcal{J}'_{+}(a,b)\,.
$$
In particular, its $\ell_{1}$-th singular value
satisfies
$$
e^{-\frac{\max(\tilde{c}_{N+1}-\tilde{c}_{2},\tilde{c}_{N}-\tilde{c}_{1})}{h}}=
\tilde{O}\left(\mu_{\ell_{1}}(\Pi_{G^{h}}d_{f,f^{-1}([a,b]),h}\big|_{E'^{h}})\right)\,.
$$
Moreover, the operator $d_{f,f^{-1}([a,b]),h}\big|_{E'^{h}}$ 
is a left multiple of $\Pi_{G^{h}}d_{f,f^{-1}([a,b]),h}\big|_{E'^{h}}$:
\begin{displaymath}
\xymatrix@C=3cm{
E'^{h}
\quad\ar[r]^{d_{f,f^{-1}([a,b]),h}}\ar[dr]_{\Pi_{G^{h}}d_{f,f^{-1}([a,b]),h}\hspace{1cm}}&
L^{2}(f^{-1}([a,b]))\\
&G^{h}\ar[u]_{\tilde{C}^{h}}
}
\end{displaymath}
with $\|\tilde{C}^{h}\|=\tilde{O}(e^{\frac{2\delta_{2}}{h}})$\,.\\
Finally, the same results hold when $G^{h}$ is replaced by $F^{h}=F_{[0,\tilde{o}(1)],[a,b],h}$\,.
\end{prop} 
\begin{proof}
 The basis $(\phi_{j}^{h})_{j\in
   \mathcal{J}'(a,b)}$ of $\mathcal{V}'^{h}(a,b)$
   defined in \eqref{eq.part-basis}
 (note that the inclusion 
$\mathcal{J}'(a,b)\subset \mathcal{J}(a,b)$ is strict in general)
has
 been constructed so that it is a partial 
$\delta_{1}$-family of
 quasimodes in the sense of Definition~\ref{de:adapted}, with 
 the additional orthogonality property
 \eqref{eq:Psiorth},\eqref{eq:defPsi}. Moreover, we know
 that (see indeed Proposition~\ref{pr:diagVh})
 \begin{eqnarray}
\label{eq:simlogphij}&&
\|\Psi_{j}^{h}\|=\|\Pi_{G^{h}}d_{f,f^{-1}([a,b]),h}T_{\delta_{2}}\phi_{j}^{h}\|_{L^{2}}\stackrel{log}{\sim}e^{-\frac{y_{\alpha}-x_{\alpha}}{h}}\ \ 
\text{when}~j=(\alpha,x_{\alpha})\in \mathcal{J}'_{+}(a,b)
\\
\label{eq:phijker}
\text{and} &&
\Psi_{j}^{h}=\Pi_{G^{h}}d_{f,h}T_{\delta_{2}}\phi_{j}^{h}=0\ \ \text{when}~j\in
   \mathcal{J}'_{0}(a,b)\,.
 \end{eqnarray}
Again, with  (see \eqref{eq.T-delta-unit})
$$
\|T_{\delta_{2}}T_{\delta_{2}}^{*}-\Id_{E'^{h}}\|+\|T_{\delta_{2}}^{*}T_{\delta_{2}}-\Id_{\mathcal{V}^{h}}\|=\tilde{O}(e^{-\frac{\eta_{f}}{h}})\,,
$$
this proves the results about the rank, the kernel, and the singular
values of $\Pi_{G^{h}}d_{f,f^{-1}([a,b]),h}\big|_{E'^{h}}$\,.\\
Moreover, reasoning with the orthogonality of the family $(\Psi_{j}^{h})_{j\in
  \mathcal{J}'(a,b)}$ and \eqref{eq:simlogphij},\eqref{eq:phijker},  like at the end of the proof of Proposition~\ref{pr:collG1},
leads to  the factorization 
$$\tilde{C}^{h}\Pi_{G^{h}}d_{f,f^{-1}([a,b]),h}T_{\delta_{2}}\big|_{\mathcal{V}'^{h}(a,b)}=d_{f,f^{-1}([a,b]),h}T_{\delta_{2}}\big|_{\mathcal{V}'^{h}(a,b)}
$$
with $\|\tilde{C}^{h}\|=\tilde{O}(e^{\frac{2\delta_{2}}{h}})$\,. We
conclude with the invertibility of
$T_{\delta_{2}}:\mathcal{V}'^{h}(a,b)\to E'^{h}$\,.\\
Finally, replacing $G^{h}$ by $F^{h}$ simply relies on Lemma~\ref{le:factorization}
used as we did around \eqref{eq:factorGh}.
\end{proof}

\subsection{$N\to N+1$: Handling the bars containing
  $[\tilde{c}_{1},\tilde{c}_{N+1}[$}
\label{sec:longlength}
We continue in the framework of the previous paragraph with
$$
[a,b]\cap\left\{c_{1},\ldots,
  c_{N_{f}}\right\}=\left\{\tilde{c}_{1},\ldots,\tilde{c}_{N+1}\right\}
$$
and 
$$
a_{1}=a\quad,\quad b_{1}=\tilde{c}_{N}+\eta_{f}\quad,\quad
a_{2}=\tilde{c}_{2}-\eta_{f}\quad,\quad b_{2}=b\,.
$$
We use the partition
$$
\mathcal{J}(a,b)=\mathcal{J}'(a,b)\sqcup \mathcal{J}''(a,b)\,,
$$
where $\mathcal{J}'(a,b)$ is defined in \eqref{eq:Jprim0}
and
\begin{eqnarray}
\label{eq:defJsec}
\mathcal{J}''(a,b)&=&\left\{j=(\alpha,\tilde{c}_{1})\in \mathcal{X}(a,b)\,,
  y_{\alpha}=\tilde{c}_{N+1}\right\}\sqcup
\left\{j=(\alpha,\tilde{c}_{1})\in \mathcal{Z}(a,b)\right\}\\
\nonumber&=&\left\{j=(\alpha,\tilde{c}_{1})\in
             \mathcal{Z}(a_{1},b_{1})\right\}=\mathcal{J}(a_{1},b_{1})\setminus
             (\mathcal{J}(a_{1},b_{1})\cap \mathcal{J}'(a,b))\,.
\end{eqnarray}
If we remember that
 $(\alpha,\tilde{c})\in \mathcal{Z}(a,b)$
can be represented by
 the bar $[\tilde{c},b[$\,, the set $\mathcal{J}''(a,b)$ actually
 collects the lower endpoints
 (which are multiple copies of $\tilde{c}_{1}$) of 
 bars containing
 $[\tilde{c}_{1},\tilde{c}_{N+1}[$\,. 
Thus, 
 the partition of $\mathcal{J}(a,b)$ and the identifications of
 $\mathcal{J}''(a,b)$ are
 clear. 
 In the preceding section, we started
with a $\delta_{1}$-family of quasimodes $(\varphi_{0,j}^{h})_{j\in
  \mathcal{J}(a_{1},b_{1})}$
 in the interval $[a_{1},b_{1}]=[a,\tilde{c}_{N}+\eta_{f}]$ 
   (see \eqref{eq:phi0j}), and only used  for 
the construction
of $E'^{h}$ in Proposition~\ref{pr:Eprim}, among the corresponding $j\in
  \mathcal{J}(a_{1},b_{1})$\,,  the indexes $j\in
\mathcal{J}(a_{1},b_{1})\cap \mathcal{J}'(a,b)$ (see \eqref{eq.part-basis}). We now use the vectors
$\varphi_{0,j}^{h}$ for $j\in \mathcal{J}''(a,b)$\,.

\begin{prop}
  \label{pr:extend}
  The vectors $\varphi_{0,j}^{h}$\,, $j\in\mathcal{J}''(a,b)$\,,
  introduced in \eqref{eq:phi0j},
  where
  $b_{1}=\tilde{c}_{N}+\eta_{f}$\,,  can
  be ``extended'' to $f^{-1}([a,b_{2}])$
    into vectors
  $\varphi_{j}^{h}\in D(d_{f,f^{-1}([a,b_{2}]),h})$ such that
  $\varphi_{j}^{h}\big|_{f_{a}^{\tilde{c}_{N}+\delta_{1}}}=
  \varphi_{0,j}^{h}\big|_{f_{a}^{\tilde{c}_{N}+\delta_{1}}}$
  and such that all the properties of   Definition~\ref{de:adapted} hold
 on
the interval $[a,b]=[a_{1},b_{2}]$ 
   with 
$I_{j}^{h}=[\tilde{c}_{1}-\delta_{1},\tilde{c}_{N+1}-\gamma''(h)]$\,, $\lim_{h\to   0}\gamma''(h)=0$\,. 
\end{prop}
\begin{proof}
 For $j\in \mathcal{J}''(a,b)$\,, $j$ has the form $j=(\alpha, \tilde c_{1})\in \mathcal{Z}(a_{1},b_{1})$ and the vector
 $\varphi_{0,j}^{h}$ then satisfies the support condition
 \eqref{eq:suppfj} (that is more precisely 
 $\supp \varphi_{0,j}^{h} \subset f^{-1}([\tilde c_{1}-\delta_{1},b_1] \cap [a_{1},b_{1}])$),
 the exponential decay estimate \eqref{eq:decfj}
 with $a,b$ replaced by $a_{1}=a, b_{1}=\tilde{c}_{N}+\eta_{f}$\,, and
 $\varphi_{0,j}^{h}\in \ker(d_{f,f^{-1}([a_{1},b_{1}]),h})$\,. For any
 $\gamma\in ]0,\eta_{f}/2[$\,, we consider the domain
 $[\tilde{c}_{N}+\delta_{1}, \tilde{c}_{N+1}-\gamma]$ and we consider
 $\tilde{c}_{N}+\eta_{f}$ as a new artificial ``critical value'', for
 which we know
 $$
\ker(\Delta_{f,f^{-1}([\tilde{c}_{N}+\delta_{1},\tilde{c}_{N+1}-\gamma]),h})=\left\{0\right\}\,.
$$
We
 then apply Proposition~\ref{pr:interdfh}-\textbf{ii)} with 
$a_{0},a,\tilde{c}_{1},b$ there replaced here by 
$\tilde{c}_{1}, \tilde{c}_{N}+\delta_{1}, 
\tilde{c}_{N}+\eta_{f}, 
\tilde{c}_{N+1}-\gamma $ and $\omega_{h}$ replaced by
$\varphi_{0,j}^{h}$\,. This provides us a new $\tilde{\omega}_{j,h}\in
D(d_{f,f^{-1}([a,\tilde{c}_{N+1}-\gamma]),h})$ which satisfies
\eqref{eq:suppfj}--\eqref{eq:dffj}, now
on $[a,\tilde{c}_{N+1}-\gamma]$
 with
$I_{j}^{h}=[\tilde c_{1}-\delta_{1},\tilde{c}_{N+1}-\gamma]$\,. 
With the cut-off
$\chi_{\tilde{c}_{N+1},\gamma}$ defined like in
Definition~\ref{de:Td2}\,, set
$$
\varphi_{j}^{\gamma,h}=\chi_{c_{N+1},\gamma}\tilde{\omega}_{j,h}
\in
D(d_{f,f^{-1}([a,b_{2}]),h})
\,.
$$
It does satisfy, on the interval $[a,b_{2}]$\,,
the conditions \eqref{eq:suppfj}--\eqref{eq:dffj}
with $I_{j}^{h}$ and $\gamma(h)$ there replaced by $[a,\tilde{c}_{N+1}-2\gamma]$
and $2\gamma$\,.
\\
For $n\in\nz$\,, take $\gamma=\frac{1}{n+1}$\,. The estimate
$B_{h}=\tilde{O}(A_{h})$ implies $B_{h}\leq e^{\frac{1}{(n+1)h}}A_{h}$
for $h\in ]0,h_{n}[$\,, and $(h_{n})_{n\in\nz}$ can be chosen to be strictly decreasing. We then adjust $\gamma''(h)=2\gamma=\frac{2}{n+1}$
for $h\in [h_{n+1},h_{n}[$ as we did at the end of the proof
of Proposition~\ref{pr:Fabab1}. This ends the proof.
 \end{proof}
\begin{remark}
\label{re.S-delta-1}
In the construction of Proposition~\ref{pr:interdfh}-\textbf{ii)}, we used the
extension by  $0$\,, here on $f_{a}^{\tilde{c}_{N}+\delta_{1}}$\,, of 
$$
d^{*}_{f,f^{-1}([\tilde{c}_{N}+\delta_{1},\tilde{c}_{N+1}-\gamma]),h}(\Delta_{f,f^{-1}([\tilde{c}_{N}+\delta_{1},\tilde{c}_{N+1}-\gamma]),h})^{-1}(hd\chi_{h}\wedge
\varphi_{0,j}^{h})\,.
$$
Because of this, the point $\tilde{c}_{N}+\delta_{1}$ must be included
in
the set $S_{\delta_{1}}$ introduced in Definition~\ref{de:adapted}.
\end{remark}

When the family $(\varphi_{j}^{h})_{j\in \mathcal{J}''(a,b)}$ is given
by Proposition~\ref{pr:extend}, the operator $T_{\delta_{2}}$ is
defined on $\Vect(\varphi_{j}^{h}, j\in \mathcal{J}''(a,b))$ by
$$
\forall j\in \mathcal{J}''(a,b), \quad T_{\delta_{2}}\varphi_{j}^{h}=\chi_{\tilde{c}_{N+1},\delta_{2}}\varphi_{j}^{h}\,,
$$ 
like in Definition~\ref{de:Td2} when $j\in \mathcal{X}(a,b)$\,.
Moreover, following the procedure of Remark~\ref{re.gamma-id},
we can assume without loss of generality that
$\gamma''(h)$\,,
given by Proposition~\ref{pr:extend},
equals $\gamma(h)$\,,
considered in Section~\ref{sec:collect} (see Remark~\ref{re.gamma-id}).
Now, the orthogonalization process  of Proposition~\ref{pr:diagVh} can
be continued by setting
\begin{equation}
\label{eq.end-ortho}
\forall j\in \mathcal{J}''(a,b),\quad
\hat{\varphi}_{j}^{h}=\varphi_{j}^{h}-
\sum_{j'\in\ds\mathop{\sqcup}_{2\leq  m'\leq N}\mathcal{X}_{m',N+1}(a,b)}
\frac{\langle\Psi_{2,j'}^{h}\,, \Pi_{G_{2,N+1}^{h}}d_{f,h}T_{\delta_{2}}\varphi_{j}^{h}\rangle}{\|\Psi_{2,j'}^{h}\|^{2}}\phi_{2,j'}^{h}\,,
\end{equation}
where $\phi_{2,j'}^{h}=\phi_{j'}^{h}$ (see \eqref{eq.part-basis})
and
$\Pi_{G_{2,N+1}^{h}}=\Pi_{G_{N+1}^{h}}$\,. Moreover, without knowing
the singular values of
$\Pi_{G_{N+1}^{h}}d_{f,h}T_{\delta_{2}}\big|_{\Vect(\hat{\varphi}_{j}^{h},
  j\in \mathcal{J}''(a,b))}$\,, we can replace the basis
$(\hat{\varphi}_{j}^{h})_{j\in J''(a,b)}$ by an orthonormal basis
$(\phi_{j}^{h})_{j\in \mathcal{J}''(a,b)}$
 such that
$\Pi_{G_{N+1}^{h}}d_{f,h}T_{\delta_{2}}\phi_{j}^{h}=\Psi_{j}^{h}$\,,
with $\Psi_{j}^{h}\perp \Psi_{j'}^{h}$ when $j\neq j'$\,, $j,j'\in
\mathcal{J}''(a,b)$\,, without changing its characteristic properties.\\
The construction of the new quasimode basis  at step $N+1$ 
is almost achieved, 
 except that the family $(\varphi_{j}^{h})_{j\in
  \mathcal{J}(a,b)}$ is not exactly a $\delta_{1}$-family of
quasimodes in the sense of Definition~\ref{de:adapted}.
In fact, we have not distinguished the endpoints of bars $j\in
\mathcal{X}_{1,N+1}(a,b)$ from the endpoints
 $j=(\alpha,\tilde{c}_{1})\in \mathcal{Z}(a,b)$ in
 \eqref{eq:defJsec}. 
For this reason, we prefer to introduce a different notation.
\begin{definition}
\label{de:interdelta1}
  The family $(\tilde{\varphi}_{j}^{h})_{j\in \mathcal{J}(a,b)}$\,, 
where
  we keep the notation $\varphi_{j}^{h}=\tilde{\varphi}_{j}^{h}$ for
  $j\in \mathcal{J}'(a,b)$\,, 
 is called
  an intermediate $\delta_{1}$-family of quasimodes if the following
  conditions are satisfied:
  \begin{enumerate}
  \item It is $\tilde{O}(e^{-\frac{\delta_{1}}{h}})$-orthonormal like
    in Theorem~\ref{th:induc} and all
  the properties of $\delta_{1}$-quasimodes in
  Definition~\ref{de:adapted} are verified, with the only difference that
  $I_{j}^{h}=[\tilde{c}_{1}-\delta_{1}, \tilde{c}_{N+1}-\gamma(h)]$
  for all $j\in \mathcal{J}''(a,b)$\,. For such a family, we set
 $\tilde{\mathcal{V}}^{h}(a,b)=\Vect(\tilde{\varphi}_{j}^{h},
  j\in \mathcal{J}(a,b))$\,, and  the operator
  $T_{\delta_{2}}:\tilde{\mathcal{V}}^{h}(a,b)\to D(d_{f,f^{-1}([a,b]),h})$
keeps the same definition
$T_{\delta_{2}}\tilde{\varphi}_{j}^{h}=T_{\delta_{2}}\varphi_{j}^{h}$
as in Definition~\ref{de:Td2} for $j\in \mathcal{J}'(a,b)$\,, while
$$
T_{\delta_{2}}\tilde{\varphi}_{j}^{h}=\chi_{\tilde{c}_{N+1},\delta_{2}}\tilde{\varphi}_{j}^{h}\quad\text{for}~j\in \mathcal{J}''(a,b)\,.
$$
\item The space $\tilde{\mathcal{V}}^{h}(a,b)$ is $\tilde{O}(e^{-\frac{\delta_{1}}{h}})$-close to $F^{h}=F_{[0,\tilde{o}(1)],[a,b],h}$:
$$
\vec{d}(\tilde{\mathcal{V}}^{h}(a,b), F_{[0,\tilde{o}(1)],[a,b],h})+\vec{d}(F_{[0,\tilde{o}(1)],[a,b],h},\tilde{\mathcal{V}}^{h}(a,b))=\tilde{O}(e^{-\frac{\delta_{1}}{h}})\,.
$$
\item When $\mathcal{V}'^{h}(a,b)=\Vect(\varphi_{j}^{h}, j\in
  \mathcal{J}'(a,b))$\,, $\mathcal{V}'^{h}_{0}(a,b)=\Vect(\varphi_{j}^{h}, j\in
  \mathcal{J}_{0}'(a,b))$\,, $\mathcal{V}'^{h}_{+}(a,b)=\Vect(\varphi_{j}^{h}, j\in
  \mathcal{J}_{+}'(a,b))$\,, all the properties of
  Proposition~\ref{pr:Eprim} hold true.
  \end{enumerate}
If $G^{h}$ is defined  like in \eqref{eq:defGh}, we use the
notation $(\tilde{\phi}_{j}^{h})_{j\in
  \mathcal{J}(a,b)}$ and  $\phi_{j}^{h}=\tilde{\phi}_{j}^{h}$ for
$j\in \mathcal{J}'(a,b)$ when the following additional orthogonality property holds:
\begin{eqnarray}
\label{eq:tPsiorth}
 &&\tilde{\Psi}_{j}^{h}\perp \tilde{\Psi}_{j'}^{h}\quad\text{for}~j\neq j'\\
\text{with}&&
\label{eq:deftPsi}
\tilde{\Psi}_{j}^{h}=\Pi_{G^{h}}d_{f,h}T_{\delta_{2}}\tilde{\phi}_{j}^{h}\,.
\end{eqnarray}
When $1\leq m<n\leq N+1$ and
$\tilde{c}_{n}-\tilde{c}_{m}<\tilde{c}_{N+1}-\tilde{c}_{1}$\,,
the corresponding spaces will be denoted
$$
\mathcal{V}_{m,n}^{h}(a,b)=\Vect(\varphi_{j}^{h}\,, j\in
     \mathcal{X}_{m,n}(a,b))\quad,\quad
\mathcal{W}_{m,n}^{h}(a,b)=\Vect(\phi_{j}^{h}\,, j\in
     \mathcal{X}_{m,n}(a,b))\,,
$$
 while
$$
 \tilde{\mathcal{V}}^{h}_{1,N+1}(a,b)=\Vect(\tilde{\varphi}_{j}^{h},
   j\in \mathcal{J}''(a,b))\quad,\quad
 \tilde{\mathcal{W}}^{h}_{1,N+1}(a,b)=\Vect(\tilde{\phi}_{j}^{h},
   j\in \mathcal{J}''(a,b))\,.
$$
\end{definition}
Our construction, and especially Proposition~\ref{pr:extend}, provides
such a family $(\tilde{\varphi}_{j}^{h})_{j\in
  \mathcal{J}(a,b)}$\,. 
 More precisely, according to
\eqref{eq.end-ortho} and the lines below,  and since $G_{n}^{h}\perp
G_{n'}^{h}$ for $1\leq n < n' \leq N+1$\,,
our construction actually provides
  a family $(\tilde{\phi}_{j}^{h})_{j\in
  \mathcal{J}(a,b)}$\,, that is satisfying in addition \eqref{eq:tPsiorth} and \eqref{eq:deftPsi}. 
Note that,
like in Proposition~\ref{pr:piGhdfhTd}, the operator 
$$
\Pi_{G^{h}}d_{f,h}T_{\delta_{2}}:\mathcal{V}^{h}(a,b)\to L^{2}(f_{a}^{b})
$$
does not depend on $\delta_{2}\in ]0,\eta_{f}[$\,.\\
In the remaining steps, we will consider various values of $a$ and $b$ and
the above properties, especially the ones involving $G^{h}$ and
$G^{h}_{N+1}=\ker(\Delta_{f,f^{-1}([\tilde{c}_{N+1}-\eta_{f},\min(b,\tilde{c}_{N+1}+\eta_{f})]),h})$\,,
which depend on $b$\,. More precisely, an intermediate
$\delta_{1}$-family of quasimodes in the sense of Definition~\ref{de:interdelta1},
and  constructed for the pair $a<b$\,,
will have to be conveniently adapted  for another pair $a'<b'$
so that it satisfies Definition~\ref{de:interdelta1} for this new pair.

\subsection{Lower bound for non zero singular values at step $N+1$}
\label{sec:lowerN+1}
This paragraph will end with the proof of
Theorem~\ref{th:induc}-\textbf{a)} at step $N+1$\,.
We are in the case 
\begin{equation}
\label{eq:N+1}
\left\{c_{1},\ldots, c_{N_{f}}\right\}\cap [a,b]=
\left\{c_{1},\ldots, c_{N_{f}}\right\}\cap ]a,b[=\left\{\tilde{c}_{1},\ldots,\tilde{c}_{N+1}\right\}\,.
\end{equation}
The notations
 $\mathcal{J}'_{+}(a,b)$\,,
$\mathcal{J}'_{0}(a,b)$\,, $\mathcal{J}'(a,b)$\,, and
$\mathcal{J}''(a,b)$ are the ones introduced  in
\eqref{eq:Jprim+}, \eqref{eq:Jprim0}, and
\eqref{eq:defJsec},
and the spaces $\mathcal{V}_{m,n}^{h}(a,b)$\,,
$\mathcal{W}_{m,n}^{h}(a,b)$\,,
$\tilde{c}_{n}-\tilde{c}_{m}<\tilde{c}_{N+1}-\tilde{c}_{1}$\,, 
$\tilde{\mathcal{V}}_{1,N+1}^{h}(a,b)$\,,
$\tilde{\mathcal{W}}_{1,N+1}^{h}(a,b)$\,,
$\mathcal{V}'^{h}_{0}(a,b)$\,, $\mathcal{V}'^{h}_{+}(a,b)$\,,
$\mathcal{V}'^{h}(a,b)$\,, are the ones of
Definition~\ref{de:interdelta1}. 
We set
$$
\ell_{0}:=\sharp \mathcal{X}(a,b)=\sharp A_{c}(a,b)=\mathrm{rank}~\delta_{[0,\tilde{o}(1)],[a,b],h}\,,
$$
where the last equality was proved in Proposition~\ref{pr:countsing}
and $\sharp A_{c}(a,b)=\sharp \mathcal{X}(a,b)$ since the number
of  bars $\alpha$ such that in $a<x_{\alpha}<y_{\alpha}<b$
 equals the number of their lower endpoints.\\
Meanwhile, we set
$$
\ell_{1}:=\ell_{0} -\sharp
\mathcal{X}_{1,N+1}(a,b)=\sharp \left\{j=(\alpha,x_{\alpha})\in
  \mathcal{X}(a,b),~
  y_{\alpha}-x_{\alpha}<\tilde{c}_{N+1}-\tilde{c}_{1}\right\}=\dim \mathcal{V}_{+}'^{h}(a,b)\,.
$$
\begin{prop}
\label{pr:roughminoG}
Consider the case $\tilde{c}_{1}-\delta_{1}\leq a<\tilde{c}_{1}$\,, $\tilde{c}_{N+1}<b\leq \tilde{c}_{N+1}+\delta_{3}$\,,
and assume $\delta_{1},\delta_{2},\delta_{3}\in
]0,\frac{\eta_{f}}{8}]$\,. Let $G^{h}$ be given by
\eqref{eq:defGh}, define $\mathcal{V}'^{h}(a,b)$\,, $\tilde{\mathcal{V}}_{1,N+1}^{h}(a,b)$\,,
and $T_{\delta_{2}}$ like in Definition~\ref{de:interdelta1},
 and consider
$$
E^{h}:=T_{\delta_{2}} \tilde{\mathcal{V}}^{h}(a,b) = T_{\delta_{2}}[\mathcal{V}'^{h}(a,b)\oplus
\tilde{\mathcal{V}}_{1,N+1}^{h}(a,b)]\,.
$$
Then, the $\ell_{0}$-th singular value of
$\Pi_{G^{h}}d_{f,h}\big|_{E^{h}}$ is bounded from below by
$$
e^{-\frac{\tilde{c}_{N+1}-\tilde{c}_{1}+\max(\delta_{1},\delta_{3})}{h}}\leq
e^{-\frac{\max(b-\tilde{c}_{1},\tilde{c}_{N+1}-a)}{h}}=\tilde{O}(\mu_{\ell_{0}}(\Pi_{G^{h}}d_{f,h}\big|_{E^{h}}))\,.
$$
\end{prop}
\begin{proof}
  With our choice $\tilde{c}_{1}-\delta_{1}\leq a<\tilde{c}_{1}$ and
  $\tilde{c}_{N+1}<b\leq\tilde{c}_{N+1}+\delta_{3}$\,, Proposition~\ref{pr:roughmino}
  says
\begin{eqnarray}
\label{eq:minomuell0}
&&
e^{-\frac{\tilde{c}_{N+1}-\tilde{c}_{1}+\max(\delta_{1},\delta_{3})}{h}}\leq
   e^{-\frac{\max(b-\tilde{c}_{1},\tilde{c}_{N+1}-a)}{h}}=\tilde{O}\left(\mu_{\ell_{0}}(\delta_{[0,\tilde{o}(1)],[a,b],h})\right)\,,\\
\text{with}&&
\nonumber
\delta_{[0,\tilde{o}(1)],[a,b],h}=\Pi_{F^{h}}d_{f,f^{-1}([a,b]),h}\big|_{F^{h}}\,,
\end{eqnarray}
where we recall $F^{h}=F_{[0,\tilde{o}(1)],[a,b],h}$\,.\\
Write 
$$
E'^{h}=T_{\delta_{2}}\mathcal{V}'^{h}(a,b)\quad\text{and}\quad 
E''^{h}=T_{\delta_{2}}\tilde{\mathcal{V}}_{1,N+1}^{h}(a,b)\,.
$$
The assumed exponential decay and the definition of
$T_{\delta_{2}}$ in Definition~\ref{de:interdelta1}  yield
$$
\vec{d}(E^{h},\tilde{\mathcal{V}}^{h}(a,b))+\vec{d}(\tilde{\mathcal{V}}^{h}(a,b),E^{h})=
\tilde{O}(e^{-\frac{\eta_{f}}{h}})
\leq\tilde{O}(e^{-\frac{\delta_{1}}{h}})
$$
and therefore
$$
\vec{d}(F^{h},E^{h})+\vec{d}(E^{h},F^{h})
=
\tilde{O}(e^{-\frac{\delta_{1}}{h}})\,.
$$
Moreover, the decomposition $E^{h}=E'^{h}\oplus E''^{h}$ is  $\tilde{O}(e^{-\frac{\delta_{1}}{h}})$-orthogonal
  and we know that
\begin{eqnarray*}
&&E^{h}\subset D(d_{f,f^{-1}([a,b],h)})\quad,\quad d_{f,f^{-1}([a,b]),h}\big|_{E^{h}}=d_{f,h}\big|_{E^{h}}\\
\text{and}&& 
\vec{d}(F^{h},G^{h})+\vec{d}(G^{h},F^{h})
=\tilde{O}(e^{-\frac{\eta_{f}}{h}})\,.
\end{eqnarray*}
In addition, Proposition~\ref{pr:Eprim}, whose properties are
ensured by the condition~3 of Definition~\ref{de:interdelta1},
 provides the factorization
\begin{eqnarray*}
  &&d_{f,f^{-1}([a,b],h)}\big|_{E'^{h}}=C^{h}\Pi_{F^{h}}d_{f,f^{-1}([a,b]),h}\big|_{E'^{h}}\\
\text{with}&&
\|C^{h}\|=\tilde{O}(e^{\frac{2\delta_{2}}{h}})\quad\text{and then}\quad
\|C^{h}\|\left[\vec{d}(F^{h},G^{h})+\vec{d}(G^{h},F^{h})\right]=\tilde{O}(e^{\frac{2\delta_{2}-\eta_{f}}{h}})
\leq\tilde{O}(e^{-\frac{\delta_{1}}{h}})\,.
\end{eqnarray*}
So, Hypotheses 1,2,3,  and the inequality
\eqref{eq:hypineq1} of Hypothesis~4 in Proposition~\ref{pr:multadd} are
satisfied with $B^{h}=d_{f,f^{-1}([a,b]),h}$ and
$\varrho(h)=\tilde{O}(e^{-\frac{\delta_{1}}{h}})$ when
$\delta_{1},\delta_{2},\delta_{3}\in ]0,\frac{\eta_{f}}{8}]$\,.
  Moreover, we know  from Proposition~\ref{pr:Eprim} that
  \begin{eqnarray*}
    &&
\mathrm{rank}(\Pi_{G^{h}}d_{f,f^{-1}([a,b]),h}\big|_{E'^{h}})=\ell_{1}=\dim
    \mathcal{V}'^{h}_{+}(a,b)=\sharp \mathcal{J}'_{+}(a,b)\\
\text{and}&&
e^{-\frac{\max(\tilde{c}_{N}-\tilde{c}_{1},\tilde{c}_{N+1}-\tilde{c}_{2})}{h}}=\tilde{O}(\mu_{\ell_{1}}(\Pi_{G^{h}}d_{f,f^{-1}([a,b],h)}\big|_{E'^{h}})\,,\\
\text{with}&&\max(\tilde{c}_{N}-\tilde{c}_{1},\tilde{c}_{N+1}-\tilde{c}_{2})\leq \tilde{c}_{N+1}-\tilde{c}_{1}-2\eta_{f}\,.
  \end{eqnarray*}
With $B_{h}=d_{f,f^{-1}([a,b]),h}$\,, the upper bound
$\|d_{f,f^{-1}([a,b],h)}\big|_{E''^{h}}\|=\tilde{O}(e^{-\frac{\tilde{c}_{N+1}-\tilde{c}_{1}-2\delta_{2}}{h}})$ (see \eqref{eq.norm-dfh}),
and
\eqref{eq:minomuell0}, 
the inequality \eqref{eq:hypineq2}
of
Hypothesis~4 is deduced from
\begin{align*}
\|B^{h}\big|_{{E''}^{h}}\|
&\left[
\frac{1}{\mu_{\ell_{1}}(\Pi_{G^{h}}B^{h}\big|_{{E'}^{h}})}
+
   \frac{\|C^{h}\|(\vec{d}(F^{h},G^{h})+\vec{d}(G^{h},F^{h}))}{\max(\mu_{\ell_{0}}(\Pi_{G^{h}}B^{h}\big|_{E^{h}}),\mu_{\ell_{0}}(B^{h}\big|_{F^{h}}))}\right]
\\
&=
\tilde{O}(e^{-\frac{\tilde{c}_{N+1}-\tilde{c}_{1}-2\delta_{2}}{h}})\times\Big[
  \tilde{O}(e^{\frac{\tilde{c}_{N+1}-\tilde{c}_{1}-2\eta_{f}}{h}})+\frac{\tilde{O}(e^{\frac{2\delta_{2}-\eta_{f}}{h}})}{\underbrace{\mu_{\ell_{0}}(B^{h}\big|_{F^{h}})}_{see~\eqref{eq:minomuell0}}}
\Big]
\\
&=\tilde{O}(e^{-\frac{\eta_{f}}{h}})+
  \tilde{O}(e^{\frac{4\delta_{2}+\max(\delta_{1},\delta_{3})-\eta_{f}}{h}})
=\tilde{O}(e^{-\frac{\delta_{1}}{h}})\,,
\end{align*}
if $\delta_{1},\delta_{2},\delta_{3}\in ]0,\frac{\eta_{f}}{8}]$\,.\\
The first result of Proposition~\ref{pr:multadd} then implies
$$
\forall \ell\in \left\{1,\ldots,\ell_{0}\right\}\,,\quad 
\mu_{\ell}(\Pi_{G^{h}}d_{f,f^{-1}([a,b]),h}\big|_{E^{h}})=\mu_{\ell}(\delta_{[0,\tilde{o}(1)],[a,b],h})(1+\tilde{O}(e^{-\delta_{1}/h}))\,,
$$
which yields in particular (see \eqref{eq:minomuell0})
$$
e^{-\frac{\max(b-\tilde{c}_{1},\tilde{c}_{N+1}-a)}{h}}=\tilde{O}(\mu_{\ell_{0}}(\Pi_{G^{h}}d_{f,f^{-1}([a,b]),h}\big|_{E^{h}}))\,.
$$
\end{proof}
In the spirit of the proof of Proposition~\ref{pr:exp0}, and in
particular  of Step~3 in Subsection~\ref{sec:exp0}, we transfer our
estimates from $[\tilde{c}_{1}-\delta_{1},\tilde{c}_{N+1}+\delta_{3}]$
to a generally wider interval $[a,b]$\,.
\begin{prop}
\label{pr:changingab}
Assume $\delta_{1},\delta_{2}\in ]0,\frac{\eta_{f}}{8}]$\,, let $a,b$ satisfy \eqref{eq:N+1},
and
let $G^{h}$ be defined by \eqref{eq:defGh}.   There exists an
intermediate $\delta_{1}$-family of quasimodes in the sense of
Definition~\ref{de:interdelta1} such that 
$$
e^{-\frac{\tilde{c}_{N+1}-\tilde{c}_{1}+\delta_{1}}{h}}
=\tilde{O}(\mu_{\ell_{0}}(\Pi_{G^{h}}d_{f,h}\big|_{E^{h}}))\quad\text{with}\quad
\ell_{0}=\sharp \mathcal{X}(a,b)\,,
$$
holds true by defining
$E^{h}=T_{\delta_{2}}\tilde{\mathcal{V}}^{h}(a,b)=T_{\delta_{2}}\Vect(\tilde{\varphi}_{j}^{h},~j\in
\mathcal{J}(a,b))$\,.
\end{prop}
\begin{proof}
Let $\delta_{1},\delta_{2}\in ]0,\frac{\eta_{f}}{8}]$\,.
When $\tilde c_{1}-\delta_{1}\leq a<\tilde c_{1}$
and $\tilde{c}_{N+1}<b\leq \tilde{c}_{N+1}+\delta_{1}$\,,
the statement of Proposition~\ref{pr:changingab}
is an immediate consequence of
Proposition~\ref{pr:roughminoG}. Moreover, when
$a<\tilde c_{1} -\delta_{1}$ and
$\tilde{c}_{N+1}<b\leq \tilde{c}_{N+1}+\delta_{1}$\,, the statement
of Proposition~\ref{pr:changingab}
simply follows after extending the quasimodes by $0$
on $f_{a}^{a'}$\,.
We thus focus on the case $b> \tilde{c}_{N+1}+\delta_{1}$\,.
Let then $\delta_{3}\in ]\delta_{1},\frac{\eta_{f}}8]$ be such that
$b':=\tilde{c}_{N+1}+\delta_{3}<b$
and set 
$a':=\max(\tilde{c}_{1}-\delta_{1},a)$\,.\\
We start from an intermediate  $\delta_{1}$-family of quasimodes
$(\tilde{\phi}_{j}^{h})_{j\in \mathcal{J}(a',b')}$\,, for the interval
$[a',b']$\,, with
the orthogonality property \eqref{eq:tPsiorth},\eqref{eq:deftPsi}.
When $a<\tilde{c}_{1}-\delta_{1}=a'$\,, these quasimodes are extended by
$0$ on $f_{a}^{a'}$\,.
We will use the spaces
$$
E^{h}(a',b')=
\underbrace{
 T_{\delta_{2}}\mathcal{V}'^{h}_{0}(a',b')
\oplus\left(\mathop{\oplus}_{1\leq n-m\leq
               N-1}T_{\delta_{2}}\mathcal{W}_{m,n}^{h}(a',b')\right)}_{E'^{h}(a',b')}
\oplus \underbrace{T_{\delta_{2}}\tilde{\mathcal{W}}^{h}_{1,N+1}(a',b')}_{E''^{h}(a',b')}
$$
and, for $(\bar{a},\bar{b})=(a',b')$ or $(\bar{a},\bar{b})=(a,b)$\,, 
$$
G^{h}(\bar{a},\bar{b})=\mathop{\oplus}_{1\leq n\leq
  N+1}^{\perp}
\underbrace{\ker(\Delta_{f,f^{-1}([\tilde{c}_{n}-\eta_{f},
    \tilde{c}_{n}+\eta_{f}]\cap[\bar{a},\bar{b}],h})}_{G_{n}^{h}(\bar{a},\bar{b})}\,.
$$
According to \eqref{eq:estimupuh} and to Propositions~\ref{pr:Eprim}
and~\ref{pr:roughminoG}, we know that
\begin{eqnarray*}
  && T_{\delta_{2}}\mathcal{V}'^{h}_{0}(a',b')\subset
     \ker(\Pi_{G^{h}(a',b')}d_{f,h}\big|_{E^{h}(a',b')})\,,\\
&&
   \|\Pi_{G^{h}(a',b')}d_{f,h}\big|_{E''^{h}(a',b')}\|=\tilde{O}(e^{-\frac{\tilde{c}_{N+1}-\tilde{c}_{1}}{h}})\,,\\
\text{and}&&e^{-\frac{\tilde{c}_{N+1}-\tilde{c}_{1}+\max(\delta_{1},\delta_{3})}{h}}=
\tilde{O}(\mu_{\ell_{0}}(\Pi_{G^{h}(a',b')}d_{f,h}\big|_{E^{h}(a',b')}))\,.
\end{eqnarray*}
Comparing the singular values of
$\Pi_{G^{h}(a',b')}d_{f,h}T_{\delta_{2}}\big|_{\mathcal{V}^{h}(a',b')}$
and of
$\Pi_{G^{h}(a',b')}d_{f,h}\big|_{E^{h}(a',b')}$
is straightforward owing to
$$
\|T_{\delta_{2}}T_{\delta_{2}}^{*}-\Id_{E^{h}(a',b')}\|+\|T_{\delta_{2}}^{*}T_{\delta_{2}}-\Id_{\mathcal{V}^{h}(a',b')}\|=\tilde{O}(e^{-\frac{\eta_{f}}{h}})\,.
$$
Meanwhile, the spaces 
$\Pi_{G^{h}(a',b')}d_{f,h}(T_{\delta_{2}}\mathcal{W}_{m,n}^{h}(a',b'))$ are mutually
orthogonal and orthogonal to
$\Pi_{G^{h}(a',b')}d_{f,h}(T_{\delta_{2}}\tilde{\mathcal{W}}_{1,N+1}(a',b'))$\,,
thanks to the orthogonality property
\eqref{eq:tPsiorth},\eqref{eq:deftPsi}. Owing
to Proposition~\ref{pr:epsorth}-\textbf{b)}, the non zero 
singular values of $\Pi_{G^{h}(a',b')}d_{f,h}\big|_{E^{h}(a',b')}$ are then
obtained by collecting the ones of
$\Pi_{G_{n}^{h}(a',b')}d_{f,h}\big|_{T_{\delta_{2}}\mathcal{W}_{m,n}^{h}(a',b')}$\,,
$1\leq n-m\leq N-1$\,, and of 
$\Pi_{G_{N+1}^{h}(a',b')}d_{f,h}\big|_{T_{\delta_{2}}\tilde{\mathcal{W}}_{1,N+1}^{h}(a',b')}$\,.\\
Moreover,
since the family $(\tilde{\phi}_{j}^{h})_{j\in \mathcal{J}'(a',b')}$ 
satisfies Definition~\ref{de:interdelta1}, and thus 
the statement of Proposition~\ref{pr:Eprim},
 the singular values of
$\Pi_{G_{n}^{h}(a',b')}d_{f,h}T_{\delta_{2}}\big|_{\mathcal{W}_{m,n}^{h}(a',b')}$
satisfy
$\mu_{h}\stackrel{log}{\sim}e^{-\frac{\tilde{c}_{n}-\tilde{c}_{m}}{h}}$
when $n-m<N-1$ (see indeed \eqref{eq:simlogphij}), while we know that
the ones of
$\Pi_{G_{N+1}^{h}(a',b')}d_{f,h}T_{\delta_{2}}\big|_{\tilde{\mathcal{W}}_{1,N+1}^{h}(a',b')}$
satisfy, for $\ell\leq \sharp \mathcal{X}_{1,N+1}(a',b')=\sharp
\mathcal{X}_{1,N+1}(a,b)$\,,
$$
e^{-\frac{\tilde{c}_{N+1}-\tilde{c}_{1}+\max(\delta_{1},\delta_{3})}{h}}
=\tilde{O}(\mu_{\ell}(\Pi_{G_{N+1}^{h}}d_{f,h}T_{\delta_{2}}\big|_{\tilde{\mathcal{W}}_{1,N+1}^{h}(a',b')}))=\tilde{O}(e^{-\frac{\tilde{c}_{N+1}-\tilde{c}_{1}}{h}})\,.
$$
Let us now construct the family $(\tilde{\varphi}_{j}^{h})_{j\in
  \mathcal{J}(a,b)}$ for the interval $[a,b]$\,.
\begin{itemize}
\item For the  $j=(\alpha,\tilde{c}_{N+1})\in \mathcal{J}'_{0}(a,b)$\,,
we take an orthonormal basis $(\tilde{\varphi}_{j}^{h})_{j=(\alpha,\tilde{c}_{N+1})
\in \mathcal{J}'_{0}(a,b)}$ of
$\ker(\Delta_{f,f^{-1}([\tilde{c}_{N+1}-\delta_{1},b]),h})$ (extended by
$0$ on $f_{a}^{\tilde{c}_{N+1}-\delta_{1}}$)\,.
\item For $j=(\alpha,\tilde{c})\in \mathcal{J}'_{0}(a,b)$ with
  $\tilde{c}<\tilde{c}_{N+1}$\,, we ``extend'' the quasimode
  $\tilde{\phi}_{j}^{h}$ as a solution to
  $d_{f,h}\tilde{\varphi}_{j}^{h}=0$ in  $[a,b]$\,,  as we did in
  Proposition~\ref{pr:extend} by referring to
  Proposition~\ref{pr:interdfh}-\textbf{ii)}, with the new artificial
  ``critical value'' $b'=\tilde{c}_{N+1}+\delta_{3}> \tilde{c}_{N+1}+\delta_{1}$\,, in the interval
  $[\tilde{c}_{N+1}+\delta_{1},b]$\,.
\item For $j\in \mathcal{X}_{m,n}(a,b)$ with $1\leq m<n\leq N$\,, we
  simply keep $\tilde{\varphi}_{j}^{h}=\tilde{\phi}_{j}^{h}$\,.
\item For the $j=(\alpha,x_{\alpha})\in \mathcal{X}(a,b)$ such that
  $y_{\alpha}=\tilde{c}_{N+1}$  and the $j=(\alpha,\tilde{c}_{1})\in \mathcal{Z}(a,b)$\,, the
  construction is detailed below after comparing,  for $m_{0}\in
  \left\{1,\ldots, N\right\}$\,, the two maps
$$
\Pi_{G_{N+1}^{h}(a',b')}d_{f,h}T_{\delta_{2}}\big|_{V_{m_{0},N+1}^{h}}\quad\text{and}\quad
\Pi_{G_{N+1}^{h}(a,b)}d_{f,h}T_{\delta_{2}}\big|_{V_{m_{0},N+1}^{h}}
=\Pi_{G_{N+1}^{h}(a',b)}d_{f,h}T_{\delta_{2}}\big|_{V_{m_{0},N+1}^{h}}\,,
$$
with
$$
V_{m_{0},N+1}^{h}=\left(\mathop{\oplus}_{\max\{2,m_{0}\}\leq m
  <N+1}\mathcal{W}_{m,N+1}^{h}(a',b')\right)
\underbrace{\oplus \tilde{\mathcal{W}}_{1,N+1}^{h}(a',b')}_{\text{if}~m_{0}=1}\,.
$$
\end{itemize}
We recall that
\begin{eqnarray*}
  &&
\dim \mathcal{W}_{m,N+1}^{h}(a',b')=\sharp
\mathcal{X}_{m,N+1}(\bar{a},\bar{b})\quad\text{when}~2\leq m<N+1\\
\text{and}&&
\dim \tilde{\mathcal{W}}_{1,N+1}^{h}(a',b')=\sharp
\mathcal{X}_{1,N+1}(\bar{a},\bar{b})\sqcup \sharp
\left\{j=(\alpha,\tilde{c}_{1})\in \mathcal{Z}(\bar{a},\bar{b})\right\}\,,
\end{eqnarray*}
where $(\bar{a},\bar{b})=(a',b')$ or $(\bar{a},\bar{b})=(a,b)$\,,
and we set, for $m_{0}\in \left\{1,\ldots N\right\}$\,,
$$
J_{m_{0},N+1}=\left(\mathop{\sqcup}_{m_{0}\leq m< N+1}
\mathcal{X}_{m,N+1}(a,b)\right)\underbrace{\sqcup \left\{j=(\alpha,\tilde{c}_{1})\in \mathcal{Z}(a,b)\right\}}_{\text{if}~m_{0}=1}\,.
$$
Since the 
$\tilde{\Psi}_{j}^{h}=\Pi_{G^{h}_{N+1}(a',b')}d_{f,h}T_{\delta_{2}}\tilde{\phi}_{j}^{h}$\,,
$j\in J_{m_{0},N+1}$\,, are mutually orthogonal and owing to 
the information on the
singular values, there exists an orthonormal basis  $(\psi_{k})_{1\leq
k\leq \dim G^{h}_{N+1}(a',b')}$ of $G_{N+1}^{h}(a',b')$ 
such that the matrix 
\begin{align*}
M^{h}&=\left(\langle \psi_{k}\,,\, \Pi_{G_{N+1}^{h}(a',b')}
  d_{f,h}T_{\delta_{2}}\tilde{\phi}_{j}^{h}\rangle\right)_{ 1\leq
       k\leq \dim G_{N+1}^{h}(a',b')\,,\; j\in
  J_{m_{0},N+1}}\\
&=
\left(\langle \psi_{k}\,,\,
  d_{f,h}T_{\delta_{2}}\tilde{\phi}_{j}^{h}\rangle\right)_{ 1\leq
  k\leq \dim G_{N+1}^{h}(a',b')\,,\; j\in
  J_{m_{0},N+1}}
\end{align*}
has the following  block diagonal structure:
\begin{itemize}
\item When $m_{0}>1$:
\begin{eqnarray*}
  &&
M^{h}=
     \begin{pmatrix}
       D^{h}\\0
     \end{pmatrix}\quad,\quad
D^{h}=\text{diag}(\lambda_{j}^{h}\,,\;
 j
\in J_{m_{0},N+1})
\,,\\
\text{where}&&\lambda_{j}^{h}\stackrel{log}{\sim}e^{-\frac{\tilde{c}_{N+1}-x_{\alpha}}{h}}~\text{for}~ j=(\alpha,x_{\alpha})
\in J_{m_{0},N+1}\,.
\end{eqnarray*}
\item When $m_{0}=1$:
\begin{eqnarray}
\nonumber
  &&
M^{h}=
\begin{pmatrix}
  D^{h}&0\\
0&R^{h}
\end{pmatrix}
\quad,\quad D^{h}=\text{diag}(\lambda_{j}^{h}\,,\;
 j=(\alpha,x_{\alpha})
\in J_{1,N+1}, x_{\alpha}\geq \tilde{c}_{2})\,,\\
\nonumber
\text{where}&&
\lambda_{j}^{h}\stackrel{log}{\sim}e^{-\frac{\tilde{c}_{N+1}-x_{\alpha}}{h}}~\text{for}~ j=(\alpha,x_{\alpha})
\in J_{1,N+1},\  x_{\alpha}\geq \tilde{c}_{2}\,,
\\
\label{eq.R-h}
\text{and}&& \|R^{h}\|=\tilde{O}(e^{-\frac{\tilde{c}_{N+1}-\tilde{c}_{1}}{h}})\,,
\end{eqnarray}
while, for $\ell_{0}'=\sharp\big(J_{1,N+1}\cap
\mathcal{X}(a,b)\big)$\,, the $\ell_{0}'$-th singular value   is bounded from below by
$$
e^{-\frac{\tilde{c}_{N+1}-\tilde{c}_{1}+\max(\delta_{1},\delta_{3})}{h}}=\tilde{O}(\mu_{\ell_{0}'}(M^{h}))\,.
$$
\end{itemize}
Proposition~\ref{pr:interdfh}-\textbf{iii)} provides an isomorphism
$A_{h}: G_{N+1}^{h}(a,b)\to G_{N+1}^{h}(a',b')=G_{N+1}^{h}(a,b')$ such that
\begin{eqnarray}
\nonumber
  &&\left\|A_{h}^{*}A_{h}-\text{Id}_{G_{N+1}^{h}(a,b)}\right\|+\|A_{h}A_{h}^{*}-\Id_{G_{N+1}^{h}(a',b')}\|=\tilde{O}(e^{-\frac{\delta_{3}}{h}})\,\\
  \nonumber
&& \forall j=(\alpha,\tilde{c})\in J_{m_{0},N+1}, \forall \psi\in
   G_{N+1}^{h}(a,b)\,,\\
   \label{eq.d-T-delta}
&&\hspace{4cm}
\langle d_{f,h}T_{\delta_{2}}\phi_{j}^{h}\,,\, \psi-A_{h}\psi\rangle=\tilde{O}(e^{-\frac{\tilde{c}_{N+1}-\tilde{c}+2\delta_{3}}{h}})\|\psi\|\,.
\end{eqnarray}
By using the $\tilde{O}(e^{-\frac{\delta_{1}}{h}})$-orthonormal basis
$(\tilde{\phi}_{j}^{h})_{j\in J_{m_{0},N+1}}$ of
$V_{m_{0},N+1}^{h}$ and the
$\tilde{O}(e^{-\frac{\delta_{3}}{h}})$-orthonormal basis
$(A_{h}^{-1}\psi_{k}^{h})_{1\leq k\leq \dim G_{N+1}^{h}(a,b)}$ of
$G_{N+1}^{h}(a,b)$\,, the singular values of the  matrix
$$
M'^{h}=\left( \langle A_{h}^{-1}\psi_{k}^{h}\,,
  d_{f,h}T_{\delta_{2}}\phi_{j}^{h}\rangle\right)_{1\leq
  k\leq \dim G_{N+1}^{h}(a,b)\,,\; j\in J_{m_{0},N+1}}
$$
coincide modulo a
$\tilde{O}(e^{-\frac{\min(\delta_{1},\delta_{3})}{h}})$-relative error
with the ones of
$\Pi_{G^{h}(a,b)}d_{f,h}T_{\delta_{2}}\big|_{V_{m_{0},N+1}^{h}}$
according to Proposition~\ref{pr:epsorth}-\textbf{a)}.
With the above inequality \eqref{eq.d-T-delta}, the $j$-th columns of $M'^{h}$ and of $M^{h}$\,,
for $j=(\alpha,x_{\alpha})\in J_{m_{0}, N+1}$\,, $x_{\alpha}\geq \tilde{c}_{2}$\,,
differ by a $\tilde{O}(\lambda_{j}^{h}\times
e^{-\frac{2\delta_{3}}{h}})$\,. When $m_{0}=1$ and
$j=(\alpha,\tilde{c}_{1})\in J_{1,N+1}$\,, the $j$-th columns
of $M'^{h}$ and of $ M^{h}$ differ by a
$\tilde{O}(e^{-\frac{\tilde{c}_{N+1}-\tilde{c}_{1}+2\delta_{3}}{h}})$ error.
Hence, we can write
$$
M'^{h}=(\text{Id}+\tilde{O}(e^{-\frac{2\delta_{3}}{h}}))M^{h}
\underbrace{+\tilde{O}(e^{-\frac{\tilde{c}_{N+1}-\tilde{c}_{1}+2\delta_{3}}{h}})}_{\text{if}~m_{0}=1}\,.
$$
When $m_{0}>1$\,, the singular values of $M'^{h}$ coincide with the
ones of $M''^{h}:=(\text{Id}+\tilde{O}(e^{-\frac{2\delta_{3}}{h}}))M^{h}$ with a
$\tilde{O}(e^{-\frac{2\delta_{3}}{h}})$-relative error.\\
When $m_{0}=1$\,, the $\ell_{0}'$-th singular value of
$M''^{h}:=(\text{Id}+\tilde{O}(e^{-\frac{2\delta_{3}}{h}}))M^{h}$ satisfies
$$
e^{-\frac{\tilde{c}_{N+1}-\tilde{c}_{1}+\max(\delta_{1},\delta_{3})}{h}}=\tilde{O}(\mu_{\ell_{0}'}(M''^{h}))\,.
$$
Hence, we get
$$
M'^{h}=M''^{h}+\tilde{O}(e^{-\frac{2\delta_{3}-\max(\delta_{1},\delta_{3})}{h}}\mu_{\ell_{0}'}(M''^{h}))\,.
$$
Since $\delta_{1}< \delta_{3}$\,, Proposition~\ref{pr:adderr}
 implies:
$$
\forall \ell\in \left\{1,\ldots,\ell_{0}'\right\}\,,\quad 
\mu_{\ell}(M'^{h})=\mu_{\ell}(M''^{h})(1+\tilde{O}(e^{-\frac{\delta_{1}}{h}}))\,.
$$
We  have thus proved  that for all $m_{0}\in
\left\{1,\ldots N\right\}$: 
\begin{multline*}
\forall \ell\in \left\{1,\ldots, \min(\sharp
  J_{m_{0},N+1},\ell_{0}')\right\}\,,\\
\mu_{\ell}(\Pi_{G_{N+1}^{h}(a,b)}d_{f,h}T_{\delta_{2}}\big|_{V_{m_{0},N+1}^{h}})=\mu_{\ell}(\Pi_{G_{N+1}^{h}(a',b')}d_{f,h}T_{\delta_{2}}\big|_{V_{m_{0},N+1}^{h}})(1+\tilde{O}(e^{-\frac{\delta_{1}}{h}}))\,.
\end{multline*}
In particular, since 
$$
\forall \delta_{3}\in ]\delta_{1},\min(\frac{\eta_{f}}{8},b-\tilde{c}_{N+1})[\,,\quad
e^{-\frac{\tilde{c}_{N+1}-\tilde{c}_{1}+\delta_{3}}{h}}=\tilde{O}(\mu_{\ell'_{0}}(\Pi_{G_{N+1}^{h}(a,b)}d_{f,h}T_{\delta_{2}}\big|_{V_{1,N+1}^{h}}))\,,
$$
and 
the right-hand side in the latter equality does not depend on $\delta_{3}$\,, we
get
\begin{equation}
\label{eq:minell0G}
e^{-\frac{\tilde{c}_{N+1}-\tilde{c}_{1}+\delta_{1}}{h}}=\tilde{O}(\mu_{\ell'_{0}}(\Pi_{G_{N+1}^{h}(a,b)}d_{f,h}T_{\delta_{2}}\big|_{V_{1,N+1}^{h}}))\,.
\end{equation}
We  now finish the presentation of our quasimodes
$(\tilde{\varphi}_{j}^{h})_{j\in J_{1,N+1}}$\,.  Like in the proof of
 Proposition~\ref{pr:diagVh}, we construct
by reverse induction from $m_{0}=N$ to $m_{0}=1$\,, starting from 
the family $(\tilde{\phi}_{j}^{h})_{j\in J_{1,N+1}}$\,,
a
basis  $(\tilde{\varphi}_{j}^{h})_{j\in J_{m_{0},N+1}}$
of $V_{m_{0},N+1}^{h}$ and 
an orthonormal  basis of $G_{N+1}^{h}(a,b)$\,,
independent of $m_{0}$\,, such that
the matrix of
$\Pi_{G_{N+1}^{h}(a,b)}d_{f,h}T_{\delta_{2}}\big|_{V_{m_{0},N+1}^{h}}$
in these bases
is diagonal (add possibly lines or columns of  zeros to make it
square). 
Since 
this  process preserves the flag
$(V_{m_{0},N+1}^{h})_{1\leq m_{0}<N+1}$\,,
the support condition and the exponential decay estimates are
valid for this new basis of $V_{1,N+1}^{h}$\,. The
$\tilde{O}(e^{-\frac{\delta_{1}}{h}})$-orthonormality of  the full new
family $(\tilde{\varphi}_{j}^{h})_{j\in \mathcal{J}(a,b)}$ and the
$\tilde{O}(e^{-\frac{\delta_{1}}{h}})$-proximity to
$F_{[0,\tilde{o}(1)],[a,b],h}$ hold
true, especially with our choice for $j=(\alpha,\tilde{c}_{N+1})\in
\mathcal{Z}(a,b)$\,. This proves the conditions~1 and 2
of Definition~\ref{de:interdelta1}. For the third condition, we notice
that the spaces $\mathcal{V}'^{h}_{+}(a,b)$ and
$\mathcal{V}_{+}'^{h}(a',b')$ are equal, like the spaces
 $G_{n}^{h}(a,b)$  and $G_{n}^{h}(a',b')$ when $2\leq n\leq N$\,, while
 $T_{\delta_{2}}$ is not changed\,.
Moreover, in the case $n=N+1$\,, 
the above orthogonalization process until
$V_{2,N+1}^{h}$ and the asymptotics of the singular values of
$\Pi_{G_{N+1}^{h}(a,b)}d_{f,h}T_{\delta_{2}}\big|_{V_{2,N+1}^{h}}$
finish the verification of the properties stated in
Proposition~\ref{pr:Eprim} for
$E'^{h}=T_{\delta_{2}}\mathcal{V}'^{h}(a,b)$ with $G^{h}=G^{h}(a,b)$\,.\\
Finally, it then follows from 
 \eqref{eq.R-h} and \eqref{eq:minell0G} that
\begin{align}
\label{eq.borne-sup}
&\mu_{\ell_{0}}(\Pi_{G^{h}(a,b)}d_{f,h}\big|_{T_{\delta_{2}}\mathcal{V}^{h}(a,b)=E^{h}})=
\tilde{O}(e^{-\frac{\tilde{c}_{N+1}-\tilde{c}_{1}}{h}})\\
 \text{and}\ \ &
e^{-\frac{\tilde{c}_{N+1}-\tilde{c}_{1}+\delta_{1}}{h}}=\tilde{O}(\mu_{\ell_{0}}(\Pi_{G_{N+1}^{h}(a,b)}d_{f,h}\big|_{T_{\delta_{2}}\mathcal{V}^{h}(a,b)=E^{h}}))\,.
\end{align}
\end{proof}
\begin{remark}
  Although we used the notation $(\tilde{\varphi}_{j}^{h})_{j\in
    \mathcal{J}(a,b)}$\,, notice that we obtain at the end of the
  proof an intermediate $\delta_{1}$-family of quasimodes $(\tilde{\phi}_{j}^{h})_{j\in
    \mathcal{J}(a,b)}$ which satisfies the orthogonality property
  \eqref{eq:tPsiorth},\eqref{eq:deftPsi} in the interval $[a,b]$\,. It
  was actually more important in the proof to put the stress on this
  property for the initial family given for the interval
  $[a',b']=[\tilde{c}_{1}-\delta_{1},\tilde{c}_{N+1}+\delta_{3}]$\,. However,
  the orthogonalization process can always be carried out afterwards.
\end{remark}
\begin{proof}[Proof of Theorem~\ref{th:induc}-{\bf a)}]
Let $a,b$ satisfy \eqref{eq:N+1} and take 
$\delta_{1},\delta_{2}\in
]0,\frac{\eta_{f}}{8}]$\,.
We reconsider the proof of Proposition~\ref{pr:roughminoG}  for the
pair $(a,b)$ with the new lower bound of Proposition~\ref{pr:changingab}:
 $$
e^{-\frac{\tilde{c}_{N+1}-\tilde{c}_{1}+\delta_{1}}{h}}
=\tilde{O}(\mu_{\ell_{0}}(\Pi_{G^{h}}d_{f,h}\big|_{E^{h}}))\quad\text{with}\quad
\ell_{0}=\sharp \mathcal{X}(a,b)\,.
$$
We then set $E^{h}=E'^{h}\oplus E''^{h}$\,,
$$
E'^{h}=T_{\delta_{2}}\mathcal{V}'^{h}(a,b)\quad,\quad
E''^{h}=T_{\delta_{2}}\tilde{\mathcal{V}}_{1,N+1}^{h}(a,b)\,,
$$
where $\mathcal{V}'^{h}(a,b)$ and
$\tilde{\mathcal{V}}_{1,N+1}^{h}(a,b)$ are associated with the
intermediate $\delta_{1}$-family of quasimodes $(\tilde{\varphi}_{j}^{h})_{j\in
  \mathcal{J}(a,b)}$ provided by Proposition~\ref{pr:changingab}.
In particular,  the verification of the inequality \eqref{eq:hypineq2}
in Proposition~\ref{pr:multadd} now becomes:
\begin{align*}
\|B^{h}\big|_{{E''}^{h}}\|
&\left[
\frac{1}{\mu_{\ell_{1}}(\Pi_{G^{h}}B^{h}\big|_{{E'}^{h}})}
+
   \frac{\|C^{h}\|(\vec{d}(F^{h},G^{h})+\vec{d}(G^{h},F^{h}))}{\max(\mu_{\ell_{0}}(\Pi_{G^{h}}B^{h}\big|_{E^{h}}),\mu_{\ell_{0}}(B^{h}\big|_{F^{h}}))}\right]
\\
&=
\tilde{O}(e^{-\frac{\tilde{c}_{N+1}-\tilde{c}_{1}-2\delta_{2}}{h}})\times\left[
\tilde{O}(e^{\frac{\tilde{c}_{N+1}-\tilde{c}_{1}-2\eta_{f}}{h}})+\tilde{O}(e^{\frac{2\delta_{2}-\eta_{f}}{h}})\times\tilde{O}(e^{\frac{\tilde{c}_{N+1}-\tilde{c}_{1}+\delta_{1}}{h}})
\right]
\\
&=\tilde{O}(e^{-\frac{\eta_{f}}{h}})+
  \tilde{O}(e^{\frac{4\delta_{2}+\delta_{1}-\eta_{f}}{h}})
=\tilde{O}(e^{-\frac{\delta_{1}}{h}})\,,
\end{align*}
with $\delta_{1},\delta_{2}\leq \frac{\eta_{f}}{8}$\,.\\
The conclusion of Proposition~\ref{pr:multadd} is then
\begin{eqnarray*}
&&
  \forall \ell\in \left\{1,\ldots,\ell_{0}\right\}\,,\quad
\mu_{\ell}(\delta_{[0,\tilde{o}(1)],[a,b],h})=\mu_{\ell}(\Pi_{G^{h}}d_{f,f^{-1}([a,b]),h}\big|_{E^{h}})(1+\tilde{O}(e^{-\frac{\delta_{1}}{h}}))\,,
\\
\text{and}&&
  \mu_{\ell_{0}+1}(\Pi_{G^{h}}d_{f,f^{-1}([a,b]),h}\big|_{E^{h}})
=\tilde{O}(e^{-\frac{\delta_{1}}{h}})\mu_{\ell_{0}}(\delta_{[0,\tilde{o}(1)],
[a,b],h})\,.
\end{eqnarray*}
In particular, we obtain
$$
e^{-\frac{\tilde{c}_{N+1}-\tilde{c}_{1}+\delta_{1}}{h}}=\tilde{O}(\mu_{\ell_{0}}(\delta_{[0,\tilde{o}(1)],[a,b],h}))
$$
and therefore, since the right-hand side of the latter equality does not depend on $\delta_{1}$\,,
$$
e^{-\frac{\tilde{c}_{N+1}-\tilde{c}_{1}}{h}}=\tilde{O}(\mu_{\ell_{0}}(\delta_{[0,\tilde{o}(1)],[a,b],h}))\,.
$$
Using in addition \eqref{eq.borne-sup} (together with Proposition~~\ref{pr:Eprim})
 leads to the statement of Theorem~\ref{th:induc}-\textbf{a)} at step $N+1$\,.
\end{proof}

We also proved 
\begin{equation}
\label{eq:majomul0p1}
\mu_{\ell_{0}+1}(\Pi_{G^{h}}d_{f,f^{-1}([a,b]),h}\big|_{E^{h}})
=\tilde{O}(e^{-\frac{\delta_{1}}{h}})\mu_{\ell_{0}}(\delta_{[0,\tilde{o}(1)],
[a,b],h})=\tilde{O}(e^{-\frac{\tilde{c}_{N+1}-\tilde{c}_{1}+\delta_{1}}{h}})\,.
\end{equation}
Moreover, according to the comments made around \eqref{eq.end-ortho},
one can choose the intermediate $\delta_{1}$-family
$(\tilde{\phi}_{j}^{h})_{j\in \mathcal{J}(a,b)}$ such that the
orthogonality property \eqref{eq:tPsiorth},\eqref{eq:deftPsi} holds, and then such that
\begin{equation}
\label{eq:eqPsij}
\|\tilde{\Psi}_{j}^{h}\|\stackrel{log}{\sim}e^{-\frac{y_{\alpha}-x_{\alpha}}{h}}\quad\text{for every}
\quad
j=(\alpha,x_{\alpha})\in\mathcal{X}(a,b)\,.
\end{equation}
\subsection{Construction of the family 
$(\varphi_{j}^{h})_{j\in\mathcal{J}(a,b)}$ at step $N+1$}
\label{sec:constN+1}

We now end the proof of Theorem~\ref{th:induc} at step $N+1$ by
finishing the construction of the $\delta_{1}$-family of quasimodes
$(\varphi_{j}^{h})_{j\in \mathcal{J}(a,b)}$\,. The  
statements \textbf{b)} and \textbf{c)}
in  Theorem~\ref{th:induc} will be easily  checked  at the end.\\
Let $a,b$ satisfy \eqref{eq:N+1}, let $G^{h}$ be defined by \eqref{eq:defGh}, and let $\delta_{1},\delta_{2}\in
]0,\frac{\eta_{f}}{8}]$\,.
We start with an intermediate $\delta_{1}$-family of quasimodes for
the interval $[a,b]$  which satisfies the orthogonality condition
\eqref{eq:tPsiorth},\eqref{eq:deftPsi} and the estimates
\eqref{eq:majomul0p1} and \eqref{eq:eqPsij}.\\
We first work in the interval $[a',b]$ with $a'=\max(a,\tilde c_{1}-\delta_{1})$\,.
Note that, since the quasimodes are all supported in $[a',b]$ 
and  
$G_{n}^{h}(a,b)=G_{n}^{h}(a',b)$  for every $2\leq n\leq
N+1$\,, the family $(\tilde{\phi}_{j}^{h})_{j\in \mathcal{J}(a,b)=\mathcal{J}(a',b)}$
is still, for the interval $[a',b]$\,,  an intermediate $\delta_{1}$-family of quasimodes
which satisfies the orthogonality condition
\eqref{eq:tPsiorth},\eqref{eq:deftPsi} and the estimates
\eqref{eq:majomul0p1} and \eqref{eq:eqPsij}.\\
The quasimodes $(\varphi_{j}^{h})_{j\in \mathcal{J}(a',b)}$ are not
changed, i.e.
$$
\varphi_{j}^{h}=\tilde{\phi}_{j}^{h}\,,
$$
when
\begin{eqnarray*}
  &&j\in \mathcal{J}_{0}'(a',b)=\mathcal{Y}(a',b)\sqcup
     \left\{j=(\alpha,\tilde{c})\in \mathcal{Z}(a',b), \tilde{c}>
     \tilde{c}_{1}\right\}\\
\text{or}&& j\in \mathcal{X}(a',b)=\mathop{\sqcup}_{1\leq m<n\leq N+1}\mathcal{X}_{m,n}(a',b)\,.
\end{eqnarray*}
We must now construct the remaining quasimodes
$\varphi_{j}^{h}$\,, $j=(\alpha,\tilde{c}_{1})\in \mathcal{Z}(a',b)$\,,
  in order to ensure 
$$
\varphi_{j}^{h}\in \ker(d_{f,f^{-1}([a',b]),h})\quad
\text{for every}~j=(\alpha,\tilde{c}_{1})\in \mathcal{Z}(a',b)\,,
$$
while we  only know for the moment that, for those $j$\,,
\eqref{eq:majomul0p1} implies
$$
\|\tilde{\Psi}_{j}^{h}\|=
\|\Pi_{G^{h}}d_{f,h}T_{\delta_{2}}\tilde{\phi}_{j}^{h}\|=\tilde{O}(e^{-\frac{\tilde{c}_{N+1}-\tilde{c}_{1}+\delta_{1}}{h}})\,.
$$
We recall that those quasimodes $\tilde{\phi}_{j}^{h}$\,, $j=(\alpha,\tilde{c}_{1})
\in \mathcal{Z}(a',b)$\,, were until now considered  in the space
$\tilde{\mathcal{W}}_{1,N+1}^{h}(a',b)$\,, together
with the quasimodes  $\tilde{\phi}_{j}^{h}$\,,
 $j\in \mathcal{X}_{1,N+1}(a',b)$\,. Let us also recall
that the rank of $\delta_{[0,\tilde{o}(1)]
  [a',b],h}$ satisfies (see Proposition~\ref{pr:countsing}):
\begin{equation}
\label{eq.rank}
\mathrm{rank}~\delta_{[0,\tilde{o}(1)],[a',b],h}=\ell_{0}=\sharp \mathcal{X}(a',b)\,.
\end{equation}
\begin{prop}
\label{pr:finalorth}
For $j=(\alpha,\tilde{c}_{1})\in \mathcal{Z}(a',b)$\,, where $a'=\max(a,\tilde{c}_{1}-\delta_{1})$\,, there exists
$(\alpha_{j,j'}^{h})_{j'\in \mathcal{X}(a',b)}$ such that
\begin{eqnarray*}
  &&\tilde{\phi}_{j}^{h}-\sum_{j'\in \mathcal{X}(a',b)}\frac{\alpha_{j,j'}^{h}}{\|\tilde{\Psi}_{j'}^{h}\|}\tilde{\phi}_{j'}^{h}
\end{eqnarray*}
belongs to $\ker(\delta_{[0,\tilde{o}(1)],[a',b],h}T_{\delta_{2}})$
with, for every $j'\in \mathcal{X}(a',b)$\,,
$$
\alpha_{j,j'}^{h}=\tilde{O}(e^{-\frac{\tilde{c}_{N+1}-\tilde{c}_{1}+\delta_{1}}{h}})\,.
$$
\end{prop}
\begin{proof}
For every $j'\in \mathcal{X}(a',b)$\,, we set
$$
 \psi_{j'}^{h}:=\frac{\tilde{\Psi}_{j'}^{h}}{\|\tilde{\Psi}_{j'}^{h}\|}\,,
$$ 
so that, when $j'=(\alpha,x_{\alpha})\in \mathcal{X}(a',b)$\,, 
$$
\Pi_{G^{h}}d_{f,h}T_{\delta_{2}}\tilde{\phi}_{j'}^{h}=\|\tilde{\Psi}_{j'}^{h}\|\psi_{j'}^{h}\quad,\quad \|\tilde{\Psi_{j'}^{h}}\|\stackrel{log}{\sim}e^{-\frac{y_{\alpha}-x_{\alpha}}{h}}
$$
and $(\psi_{j'}^{h})_{j'\in \mathcal{X}(a',b)}$ is an orthonormal system
in $G^{h}$\,.\\
By writing, for $j'\in \mathcal{X}(a',b)$\,,
\begin{align}
\nonumber
\delta_{[0,\tilde{o}(1)],[a',b],h}T_{\delta_{2}}\tilde{\phi}_{j'}^{h}&=
\Pi_{F^{h}}d_{f,h}T_{\delta_{2}}\tilde{\phi}_{j'}^{h}\\
\nonumber
&=\Pi_{G^{h}}d_{f,h}T_{\delta_{2}}\tilde{\phi}_{j'}^{h}
-(\Pi_{G^{h}}-\Pi_{F^{h}}\Pi_{G^{h}})d_{f,h}T_{\delta_{2}}\tilde{\phi}_{j'}^{h}\\
\label{eq:decomdphi}
&\qquad\qquad\qquad\quad\ \ +(\Pi_{F^{h}}-\Pi_{F^{h}}\Pi_{G^{h}})d_{f,h}T_{\delta_{2}}\tilde{\phi}_{j'}^{h}
\end{align}
with $F^{h}=F_{[0,\tilde{o}(1)],[a',b],h}$\,, $\vec{d}(F^{h},G^{h})+\vec{d}(G^{h},F^{h})=\tilde{O}(e^{-\frac{\eta_{f}}{h}})$\,,
and (see \eqref{eq.norm-dfh})
\begin{eqnarray*}
&&
\|d_{f,h}T_{\delta_{2}}\tilde{\phi}_{j'}^{h}\|=\|\tilde{\Psi}_{j'}^{h}\|\tilde{O}(e^{\frac{2\delta_{2}}{h}}),
\end{eqnarray*}
we deduce from \eqref{eq:decomdphi} that 
the family made of the
$$
\theta_{j'}^{h}=\frac{\delta_{[0,\tilde{o}(1)],[a',b],h}T_{\delta_{2}}\tilde{\phi}_{j'}^{h}}{\|\tilde{\Psi}_{j'}^{h}\|}\quad,\quad
j'\in \mathcal{X}(a',b)\,,
$$
defines an $\tilde{O}(e^{-\frac{\delta_{1}}{h}})$-orthornormal system of
$R^{h}:=\textrm{Ran}~\delta_{[0,\tilde{o}(1)],[a',b],h}$\,.
Owing to \eqref{eq.rank},  
the family $(\theta_{j'}^{h})_{j'\in \mathcal{X}(a',b)}$ is thus
an $\tilde{O}(e^{-\frac{\delta_{1}}{h}})$-orthonormal basis of $R^{h}$\,.
Denoting now by $(\hat{\theta}_{j'}^{h})_{j'\in \mathcal{X}(a',b)}$ the dual basis
of
$(\theta_{j'}^{h})_{j'\in \mathcal{X}(a',b)}$ in $R^{h}$\,, that is the unique family
satisfying
\begin{equation*}
 \forall \,j'_1, j'_{2}\, \in\,  \mathcal{X}(a',b)\,,\ \  \   \hat{\theta}_{j'_{1}}^{h}\in R_{h}\quad\text{and}\quad \langle
  \hat{\theta}_{j'_{1}}\,,\, \theta_{j'_{2}}\rangle=\delta_{j'_{1},j'_{2}}\,,
\end{equation*}
the family $(\hat{\theta}_{j'}^{h})_{j'\in \mathcal{X}(a',b)}$ is also an
$\tilde{O}(e^{-\frac{\delta_{1}}{h}})$-orthonormal basis of $R^{h}$
and
 the
orthogonal projection
on $R^{h}$ is given by 
$$
\forall u\in F^{h}\,, \quad
\Pi_{^{R^{h}}}u=\sum_{j'\in
  \mathcal{X}(a',b)}\langle\hat{\theta}_{j'}^{h}\,, u\rangle \theta_{j'}=\sum_{j'\in
  \mathcal{X}(a',b)}\frac{\langle\hat{\theta}_{j'}^{h}\,, u\rangle}{\|\tilde{\Psi}_{j'}^{h}\|}\delta_{[0,\tilde{o}(1)],[a',b],h}T_{\delta_{2}}\tilde{\phi}_{j'}^{h}\,. 
$$
For $j=(\alpha,\tilde{c}_{1})\in \mathcal{Z}(a',b)$\,, the same
decomposition as \eqref{eq:decomdphi} with now
$\|\tilde{\Psi}_{j}^{h}\|=\tilde{O}(e^{-\frac{\tilde{c}_{N+1}-\tilde{c}_{1}+\delta_{1}}{h}})$
and 
$\|d_{f,h}T_{\delta_{2}}\tilde{\phi}_{j}^{h}\|=\tilde{O}(e^{-\frac{\tilde{c}_{N+1}-\tilde{c}_{1}-2\delta_{2}}{h}})$
leads to 
$$
\|\delta_{[0,\tilde{o}(1)],[a',b],h}T_{\delta_{2}}\tilde{\phi}_{j}^{h}\|=\tilde{O}(e^{-\frac{\tilde{c}_{N+1}-\tilde{c}_{1}+\delta_{1}}{h}})\,.
$$
The statement of Proposition~\ref{pr:finalorth} follows easily by taking,
for every $j=(\alpha,\tilde{c}_{1})\in \mathcal{Z}(a',b)$
and $j'\in \mathcal{X}(a',b)$\,,
$$
\alpha_{j,j'}^{h}=\langle\hat{\theta}_{j'}^{h}\,, \delta_{[0,\tilde{o}(1)],[a',b],h}T_{\delta_{2}}\tilde{\phi}_{j}^{h}\rangle\,.
$$
\end{proof}

The following statement finishes the proof of Theorem~\ref{th:induc}. 

\begin{prop}
Assume that $a,b$ satisfy \eqref{eq:N+1}, let
$\delta_{1},\delta_{2}\in ]0,\frac{\eta_{f}}{8}]$\,, and set $a'=\max(a, \tilde{c}_{1}-\delta_{1})$\,.
  The family $(\varphi_{j}^{h})_{j\in \mathcal{J}(a,b)}$ defined by
  $$\varphi_{j}^{h}=\tilde{\phi}_{j}^{h}\quad \text{when $j\in
  \mathcal{X}(a,b)\sqcup \mathcal{Y}(a,b)\sqcup
  \left\{(\alpha,\tilde{c})\in \mathcal{Z}(a,b)\,,
    \tilde{c}>\tilde{c}_{1}\right\}$}$$  
    and
$$
\varphi_{j}^{h}=1_{f_{a'}^{b}}\times \Pi_{[0,\tilde{o}(1)],[a',b],h}T_{\delta_{2}}\Big(
\tilde{\phi}_{j}^{h}-\sum_{j'\in
\mathcal{X}(a,b)}\frac{\alpha_{j,j'}^{h}}{\|\tilde{\Psi}_{j'}^{h}\|}\tilde{\phi}_{j'}^{h}\Big)\quad\text{when}~j=(\alpha,\tilde{c}_{1})\in \mathcal{Z}(a,b)\,,
$$
where the coefficients $\alpha_{j,j'}^{h}$ are given by
Proposition~\ref{pr:finalorth}, fulfills all the conditions of
Theorem~\ref{th:induc} at step $N+1$\,.
\end{prop}
\begin{proof}
We use here the notations $a'=\max(a,\tilde{c}_{1}-\delta_{1})$ and, in
order to avoid confusions, 
\begin{eqnarray*}
  &&
\tilde{\mathcal{W}}^{h}(a,b)=\Vect(\tilde{\phi}_{j}^{h}, j\in
     \mathcal{J}(a,b))\\
\text{and}&&\tilde{\mathcal{W}}_{+}^{h}(a,b)=\Vect(\tilde{\phi}_{j}^{h}, j\in \mathcal{X}(a,b))\,,
\end{eqnarray*}
where $(\tilde{\phi}_{j}^{h})_{j\in \mathcal{J}(a,b)}$ is the
intermediate $\delta_{1}$-family of quasimodes
 we started with.\\
From the estimates
$\alpha_{j,j'}^{h}=\tilde{O}(e^{-\frac{\tilde{c}_{N+1}-\tilde{c}_{1}+\delta_{1}}{h}})$
(see Proposition~\ref{pr:finalorth})
and
$\|\tilde{\Psi}_{j'}^{h}\|\stackrel{log}{\sim}e^{-\frac{y_{\alpha}-x_{\alpha}}{h}}$
for $j=(\alpha,\tilde{c}_{1})\in \mathcal{Z}(a,b)$
and $j'=(\alpha,x_{\alpha})\in \mathcal{X}(a,b)$\,, we deduce
that
$$
\forall\,j=(\alpha,\tilde{c}_{1})\in \mathcal{Z}(a,b)\,,\ \ \ \Big\|\sum_{j'\in \mathcal{X}(a,b)}\frac{\alpha_{j,j'}^{h}}{\|\tilde{\Psi}_{j'}^{h}\|}\tilde{\phi}_{j'}^{h}\Big\|_{L^{2}}=\tilde{O}(e^{-\frac{\delta_{1}}{h}})\,.
$$
Since in addition
$\vec{d}(F_{[0,\tilde{o}(1)],[a',b],h},\tilde{\mathcal{W}}^{h}(a,b))+\vec{d}(\tilde{\mathcal{W}}^{h}(a,b),F_{[0,\tilde{o}(1)],[a',b],h})=\tilde{O}(e^{-\frac{\delta_{1}}{h}})$\,,
it follows that
$$
\|\varphi_{j}^{h}-\tilde{\phi}_{j}^{h}\|
=\tilde{O}(e^{-\frac{\delta_{1}}{h}})\quad \text{for}~j=(\alpha,\tilde{c}_{1})\in \mathcal{Z}(a,b)\,,
$$
and the family $(\varphi_{j}^{h})_{j\in \mathcal{J}(a,b)}$ is thus
$\tilde{O}(e^{-\frac{\delta_{1}}{h}})$-orthonormal.
Moreover,
the exponential decay estimates on the $\tilde{\phi}_{j'}^{h}$\,, $j'\in \mathcal{X}(a,b)$\,, lead to
$$
\forall\,j=(\alpha,\tilde{c}_{1})\in \mathcal{Z}(a,b)\,,\ \ \ \Big\|e^{\frac{|f-\tilde{c}_{1}|}{h}}\Big(\sum_{j'\in
    \mathcal{X}(a,b)}\frac{\alpha_{j,j'}^{h}}{\|\tilde{\Psi}_{j'}^{h}\|}\tilde{\phi}_{j'}^{h}\Big)
\Big\|_{W(f^{-1}([a',b]\setminus S_{\delta_{1}})}=\tilde{O}(e^{-\frac{\delta_{1}}{h}})\,.
$$
This implies,
together with Proposition~\ref{pr:projAgm},
the required exponential decay estimates on the $\varphi_{j}^{h}$\,, 
$j=(\alpha,\tilde{c}_{1})\in \mathcal{Z}(a,b)$\,.
Besides,
Proposition~\ref{pr:finalorth} gives
\begin{align*}
d_{f,h}\Pi_{[0,\tilde{o}(1)],[a',b],h}T_{\delta_{2}}\Big(\tilde{\phi}_{j}^{h}-\sum_{j'\in
  \mathcal{X}(a,b)}\frac{\alpha_{j,j'}^{h}}{\|\tilde{\Psi}_{j'}^{h}\|}\tilde{\phi}_{j}^{h}\Big)
=
\delta_{[0,\tilde{o}(1)],[a',b],h}T_{\delta_{2}}\Big(\tilde{\phi}_{j}^{h}-\sum_{j'\in
  \mathcal{X}(a,b)}\frac{\alpha_{j,j'}^{h}}{\|\tilde{\Psi}_{j'}^{h}\|}\tilde{\phi}_{j}^{h}\Big)=0\,.
\end{align*}
All those properties are preserved after extending those quasimodes by
$0$ on $f_{a}^{a'}$ when $a<a'$\,. 
Therefore, the family $(\varphi_{j}^{h})_{j\in \mathcal{J}(a,b)}$ satisfies all
the conditions of Definition~\ref{de:adapted} and is thus
a
$\tilde{O}(e^{-\frac{\delta_{1}}{h}})$-orthonormal
 $\delta_{1}$-family
of quasimodes. 
Since in addition
$$
\vec{d}(F_{[0,\tilde{o}(1)],[a',b],h},F_{[0,\tilde{o}(1)],[a,b],h})+\vec{d}(F_{[0,\tilde{o}(1)],[a,b],h},F_{[0,\tilde{o}(1)],[a',b],h})=\tilde{O}(e^{-\frac{\delta_{1}}{h}})\,,
$$
the statement \textbf{b)} of  Theorem~\ref{th:induc} is  also
satisfied.\\
It only remains to check the factorization stated in
Theorem~\ref{th:induc}-\textbf{c)}. Since
$$
d_{f,h}T_{\delta_{2}}\varphi_{j}^{h}=d_{f,h}\varphi_{j}^{h}=0\quad\text{for every}~j\not\in \mathcal{X}(a,b)\,,
$$ 
it suffices to prove the existence of $C^h$ such that
$$
\xymatrix@C=3cm{
\mathcal{V}_{+}^{h}(a,b)=\tilde{\mathcal{W}}_{+}^{h}(a,b)
\quad\ar[r]^{d_{f,f^{-1}([a,b]),h}T_{\delta_{2}}}
\ar[dr]_{\Pi_{[0,\tilde{o}(1)],[a,b],h}d_{f,h}T_{\delta_{2}}\hspace{1
    cm}}&
 L^{2}(f^{-1}([a,b]))\\
 &F_{[0,\tilde{o}(1)],[a,b],h}\ar[u]_{C^{h}}
}
$$
with $\|C^{h}\|=\tilde{O}(e^{\frac{2\delta_{2}}{h}})$\,.
Since
$\Pi_{G^{h}}d_{f,h}T_{\delta_{2}}\tilde{\phi}_{j}^{h}=\tilde{\Psi}_{j}^{h}$
with
$\|\tilde{\Psi}_{j}^{h}\|\stackrel{log}{\sim}e^{-\frac{y_{\alpha}-x_{\alpha}}{h}}$
when $j=(\alpha,x_{\alpha})\in \mathcal{X}(a,b)$ with the
orthogonality property \eqref{eq:tPsiorth},\eqref{eq:deftPsi}, 
reasoning as at the ends of the proofs of Propositions~\ref{pr:collG1}
and~\ref{pr:Eprim},
we obtain the diagram
$$
\xymatrix@C=3cm{
\mathcal{V}_{+}^{h}(a,b)=\tilde{\mathcal{W}}_{+}^{h}(a,b)
\quad\ar[r]^{d_{f,f^{-1}([a,b]),h}T_{\delta_{2}}}
\ar[dr]_{\Pi_{G^{h}}d_{f,h}T_{\delta_{2}}}&
L^{2}(f^{-1}([a,b]))\\
&G^{h}\ar[u]_{\tilde{C}^{h}}
}
$$
with $\|\tilde{C}^{h}\|=\tilde{O}(e^{\frac{2\delta_{2}}{h}})$\,.
We conclude by applying Lemma~\ref{le:factorization} with
$B^{h}=d_{f,f^{-1}([a,b]),h}T_{\delta_{2}}$\,, 
$F^{h}=F_{[0,\tilde{o}(1)],[a,b],h}$\,, and
$$
\vec{d}(F^{h},G^{h})+\vec{d}(G^{h},F^{h})=\tilde{O}(e^{-\frac{\eta_{f}}{h}})\,.
$$
\end{proof}

\section{Corollaries of Theorem~\ref{th:induc}}
\label{sec:coroll}
The statement or Theorem~\ref{th:induc} is much more flexible than its
illustrative statement, Theorem~\ref{th:mainsimple}, given in the
introduction.
Actually, even its proof, and especially the intermediate propositions of
Subsection~\ref{sec:conseqN}, have easily derived consequences which
are listed here. Subsection~\ref{sec:specres} reviews consequences on
the eigenvalues and eigenvectors of the Witten Laplacian
$\Delta_{f,f^{-1}([a,b]),h}$ when $f$ is
fixed. Subsection~\ref{sec:firststab} studies how the logarithms of
the singular values of $d_{f,f^{-1}([a,b]),h}$ vary when $f$ is
changed. It contains a generalization of
Corollary~\ref{cor:mainsimple}. 
Remember that Theorem~\ref{th:induc} is proved under
Hypothesis~\ref{hyp:cN} which gathers  Hypothesis~\ref{hyp:mainf} or
(Hypothesis~\ref{hyp:Lipbar} and Hypothesis~\ref{hyp:AgmonLip}) for
a more general Lipschitz function $f$\,. 
Hypothesis~\ref{hyp:mainf} or Hypothesis~\ref{hyp:Lipbar} ensure that
$f$ has finitely many ``critical values'' $c_{1}<\ldots<c_{N_{f}}$\,. \\
When $a,b\not\in \left\{c_{1},\ldots,c_{N_{f}}\right\}$\,,
$\Delta_{f,f^{-1}([a,b]),h}$ is the self-adjoint Witten Laplacian in
$f_{a}^{b}$\,, with Dirichlet boundary conditions on
$f^{-1}(\left\{a\right\})$ and Neumann boundary conditions on
$f^{-1}(\left\{b\right\})$\,, according to Section~\ref{sec:expdec}.\\
Finally, the bar code associated with $f$\,, under
Hypothesis~\ref{hyp:mainf} or Hypothesis~\ref{hyp:Lipbar} (see
Subsection~\ref{sec:moregenLipgen}), is
$\mathcal{B}(f)=([a_{\alpha},b_{\alpha}[)_{\alpha\in A}$\,, defined in
Subsection~\ref{sec:bar code}
and in Appendix~\ref{app:perscohom}. The restricted bar code
$\mathcal{B}(f;a,b)$\,, and the set of endpoints 
$\mathcal{J}(a,b)$\,, $\mathcal{X}(a,b)$\,, $\mathcal{Y}(a,b)$\,, $\mathcal{Z}(a,b)$\,, all
graded according to the degree $p\in \left\{0,\ldots,d\right\}$\,, are
the ones introduced in Subsection~\ref{sec:bar code}.
\subsection{Spectral results}
\label{sec:specres}
The first result generalizes Theorem~\ref{th:mainsimple}.
\begin{thm}
\label{th:specres}
Assume Hypothesis~\ref{hyp:mainf} or (Hypothesis~\ref{hyp:Lipbar} and
Hypothesis~\ref{hyp:AgmonLip}) for a more general Lipschitz function
$f$\,.  Let $a,b\not\in
\left\{c_{1},\ldots, c_{N_{f}}\right\}$ with $a<b$ and let
$\Delta_{f,f^{-1}([a,b]),h}=\mathop{\oplus}_{p=0}^{d}\Delta_{f,f^{-1}([a,b]),h}^{(p)}$ 
be defined like in
Proposition~\ref{pr:domain} with 
$N_{t}=f^{-1}(\left\{a\right\})$ and
$N_{n}=f^{-1}(\left\{b\right\})$\,.\\
The number of $\tilde{o}(1)$-eigenvalues of
$\Delta_{f,f^{-1}([a,b]),h}^{(p)}$ equals $\sharp
\mathcal{J}^{(p)}(a,b)$\,, while
$$
\dim \ker(\Delta_{f,f^{-1}([a,b]),h})=\beta^{(p)}(f^{b},f^{a})=\sharp
\mathcal{Z}^{(p)}(a,b)\,.
$$
Moreover, the non zero $\tilde{o}(1)$-eigenvalues
of
$\Delta_{f,f^{-1}([a,b]),h}^{(p)}$
counted with multiplicity can be labelled
 $\lambda_{\alpha}^{(p)}(h)$\,, $\alpha\in
A_{c}^{(p)}(a,b)\sqcup A_{c}^{(p-1)}(a,b)$\,, with
$$
\lambda_{\alpha}^{(p)}(h)\stackrel{log}{\sim}e^{-2\frac{y_{\alpha}^{*+1}-x_{\alpha}^{*}}{h}}\quad,\quad
\alpha\in A^{(p)}_{c}(a,b)\sqcup A^{(p-1)}_{c}(a,b)\,.
$$
\end{thm}
With the usual supersymmetric argument which was already recalled in
Proposition~\ref{pr:countsing}, it is a straightforward 
consequence of Theorem~\ref{th:induc}-\textbf{a)}.\\
The above result can be completed by some information on the
eigenvectors.
We start with the link between the singular values of
$\delta_{f,f^{-1}([a,b]),h}$\,, and their approximation via the
introduction of a  basis made of quasimodes, and the spectral elements of 
the operator
$\delta_{[0,\tilde{o}(1)],[a,b]),h}^{*}\delta_{[0,\tilde{o}(1)],[a,b]),h}$\,.
The spectral elements of
$$
\Pi_{[0,\tilde{o}(1)],[a,b],h}\Delta_{f,f^{-1}([a,b]),h}=\delta_{[0,\tilde{o}(1)],[a,b]),h}^{*}\delta_{[0,\tilde{o}(1)],[a,b]),h}+\delta_{[0,\tilde{o}(1)],[a,b]),h}\delta_{[0,\tilde{o}(1)],[a,b]),h}^{*}
$$
will be described afterwards by referring to Hodge decomposition and
to
duality.
\begin{prop}
\label{pr:singvdelta}
 Keep the same assumptions as in Theorem~\ref{th:specres} and define 
$\eta_{f}>0$ like in Hypothesis~\ref{hyp:cN}. Let
 $\delta_{[0,\tilde{o}(1)],[a,b],h}^{(p)}$
denote the restriction of
 $d_{f,f^{-1}([a,b]),h}$ to $F_{[0,\tilde{o}(1)],[a,b],h}^{(p)}$\,, $\delta_{[0,\tilde{o}(1)],[a,b],h}^{(p)}:F_{[0,\tilde{o}(1)],[a,b],h}^{(p)}\to
 F_{[0,\tilde{o}(1)],[a,b],h}^{(p+1)}$\,,  according
 to \eqref{eq:deltap}, and set
 \begin{eqnarray*}
   &&L^{(p)}=\left\{b_{\alpha}^{(p+1)}-a_{\alpha}^{(p)}, \alpha\in
   A_{c}^{(p)}(a,b)\right\}\,,\\
&&\delta_{f}=\min(\frac{\eta_{f}}{8}, \frac{|\ell-\ell'|}{8},
\ell\neq \ell'\in L^{(p)})>0\,.
 \end{eqnarray*}
Take the $\delta_{1}$-family of quasimodes $(\varphi_{j}^{h})_{j\in
  \mathcal{J}(a,b)}$ given by Theorem~\ref{th:induc} with
$\delta_{1}=\frac{\eta_{f}}{8}$ (and with any $\delta_{2}\in
]0,\frac{\eta_{f}}{8}]$) and define, for $\ell\in L^{(p)}$\,,
$$
\mathcal{U}_{\ell}^{(p),h}:=\Vect\left(\varphi_{j}^{h}\,,~
  j=(\alpha,x_{\alpha}^{(p)})\in \mathcal{X}^{(p)}(a,b)\,,~y_{\alpha}^{(p+1)}-x_{\alpha}^{(p)}=\ell\right)\,,
$$
and
$$
\mathcal{U}_{+\infty}^{(p),h}:=\Vect(\varphi_{j}^{h}\,,  j\in
\mathcal{Y}^{(p)}(a,b)\sqcup \mathcal{Z}^{(p)}(a,b))\,.
$$
Then, for every $\ell\in
L^{(p)}\sqcup \left\{+\infty\right\}$ and $p\in \left\{0,\ldots
  d\right\}$\,, the distance between $\mathcal{U}_{\ell}^{(p),h}$ and
> $F_{\ell}^{(p),h}$ is estimated by
$$
\vec{d}(\mathcal{U}_{\ell}^{(p),h},
F_{\ell}^{(p),h})
+\vec{d}(F_{\ell}^{(p),h},\mathcal{U}^{(p),h}_{\ell})=\tilde{O}(e^{-\frac{\delta_{f}}{h}})\,,
$$
where 
$F_{\ell}^{(p),h}\subset F_{[0,\tilde{o}(1)],[a,b],h}^{(p)}\subset
L^{2}(f_{a}^{b};\Lambda^{p} T^{*}M)$
is the spectral subspace
of 
$\delta_{[0,\tilde{o}(1)],[a,b],h}^{(p),*}\delta_{[0,\tilde{o}(1)],[a,b],h}^{(p)}$
 for the spectral range
$[e^{-2\frac{\ell+\delta_{f}}{h}},e^{-2\frac{\ell-\delta_{f}}{h}}]$\,.
\end{prop}
\begin{proof}
With our choice $\delta_{1}=\frac{\eta_{f}}{8}$\,, the basis
$(\varphi_{j}^{h})_{j\in \mathcal{J}^{(p)}(a,b)}$ is a
$\tilde{O}(e^{-\frac{\eta_{f}}{8h}})$-orthonormal family such that,
according to Theorem~\ref{th:induc}-{\bf b)} and to the definition of $T_{\delta_{2}}$
(see Definition~\ref{de:Td2}),
\begin{equation}
\label{eq.Pi-T}
\forall j\in \mathcal{J}^{(p)}(a,b),\quad \|\Pi_{[0,\tilde{o}(1)],[a,b],h}T_{\delta_{2}}\varphi_{j}^{h}-\varphi_{j}^{h}\|=\tilde{O}(e^{-\frac{\eta_{f}}{8h}})\,.
\end{equation}
For $j\in \mathcal{Y}^{(p)}(a,b)\sqcup \mathcal{Z}^{(p)}(a,b)$\,, the equality
$$
\delta_{[0,\tilde{o}(1)],[a,b],h}\Pi_{[0,\tilde{o}(1)],[a,b],h}T_{\delta_{2}}\varphi_{j}^{h}=\Pi_{[0,\tilde{o}(1)],[a,b],h}d_{f,h}\varphi_{j}^{h}=0
$$
then implies that
$(\Pi_{[0,\tilde{o}(1)],[a,b],h}T_{\delta_{2}}\varphi_{j}^{h})_{j\in
  \mathcal{Y}^{(p)}(a,b)\sqcup \mathcal{Z}^{(p)}(a,b)}$ is a
$\tilde{O}(e^{-\frac{\eta_{f}}{8h}})$-orthonormal basis of
$$
\ker(\delta_{[0,\tilde{o}(1)],[a,b],h}^{(p)})=F^{(p),h}_{+\infty}\,.
$$
This leads to the result for $\ell=+\infty$ and initializes the
decreasing induction with respect to $\ell$\,.\\
Assume now that for all $\ell>\ell_{0}$ in 
$L^{(p)}$\,, we have proved
$$
\vec{d}(\mathcal{U}^{(p),h}_{\ell},F^{(p),h}_{\ell})+\vec{d}(F^{(p),h}_{\ell},\mathcal{U}^{(p),h}_{\ell})
=\tilde{O}(e^{-\frac{\delta_{f}}{h}})\,.
$$
Let us check that it is still true for $\ell=\ell_{0}$\,. Like in
Subsection~\ref{sec:conseqN}, we introduce $G^{h}$ defined by 
\eqref{eq:defGh},\eqref{eq:defOmbar},
$G_{n}^{h}=\ker(\Delta_{f,f^{-1}([\tilde{c}_{n}-\eta_{f},\tilde{c}_{n}+\eta_{f}]\cap[a,b]),h})$
defined in \eqref{eq:defGn},
 and the spaces
$\mathcal{V}_{m,n}^{h}$ defined in \eqref{eq:defVmn} by
$$
\mathcal{V}_{m,n}^{(p),h}=\Vect(\varphi_{j}^{h}, j\in
\mathcal{X}^{(p)}_{m,n}(a,b))\,.
$$
In particular,  we have
$$
\mathcal{U}_{\ell_{0}}^{(p),h}=\mathop{\oplus}_{\tilde{c}_{n}-\tilde{c}_{m}=\ell_{0}}\mathcal{V}_{m,n}^{(p),h}\,,
$$
while
$\Pi_{G^{h,(p+1)}}d_{f,h}^{(p)}T_{\delta_{2}}(\mathcal{V}_{m,n}^{(p),h})\subset
G^{h,(p+1)}_{n}$ with $G^{h,(p+1)}_{n}\perp G^{h,(p+1)}_{n'}$ for
$n\neq n'$\,. From Proposition~\ref{pr:piGhdfhTd}, we know that the
mapping
$\Pi_{G^{h,(p+1)}}d_{f,h}^{(p)}T_{\delta_{2}}:\mathcal{U}_{\ell_{0}}^{(p),h}\to
G^{h,(p+1)}$ does not depend on $\delta_{2}\in ]0,\frac{\eta_{f}}{8}]$\,,
while Proposition~\ref{pr:Vmn} and
Proposition~\ref{pr:epsorth}-\textbf{b)} ensure that  it is one to
one with (only non zero) singular values 
all satisfying
$\mu_{h}\stackrel{log}{\sim}e^{-\frac{\ell_{0}}{h}}$\,.
Moreover, following the analysis made in the proof of Proposition~\ref{pr:piGhdfhTd},
the factorization \eqref{eq:factorGh} holds with here
$E^{h}=T_{\delta_{2}}\mathcal{U}^{(p),h}_{\ell_{0}}$\,, 
$B_{h}=d_{f,f^{-1}([a,b]),h}^{(p)}$\,, and
$\|\tilde{C}^{h}\|=\tilde{O}(e^{\frac{2\delta_{2}}{h}})$\,.
Hence, using Lemma~\ref{le:factorization} with
the relation
$$
\vec{d}(G^{h,(p+1)},F_{[0,\tilde{o}(1)],[a,b],h}^{(p+1)})+\vec{d}(F_{[0,\tilde{o}(1)],[a,b],h}^{(p+1)},G^{h,(p+1)})=\tilde{O}(e^{-\frac{\eta_{f}}{h}})
$$
leads 
  to
\begin{eqnarray*}
  &&\Pi_{G^{h,(p+1)}}d^{(p)}_{f,f^{-1}([a,b]),h}\big|_{E^{h}}=(\Id_{L^{2}(f_{a}^{b})}+\tilde{O}(e^{\frac{2\delta_{2}-\eta_{f}}{h}}))\underbrace{\Pi_{[0,\tilde{o}(1)],[a,b],h}d^{(p)}_{f,f^{-1}([a,b],h)}\big|_{E^{h}}}_{=\delta^{(p)}_{[0,\tilde{o}(1)],[a,b],h}\big|_{E^{h}}}\,.
\end{eqnarray*}
Thus, since $T_{\delta_{2}}:\mathcal{U}^{(p),h}_{\ell_{0}}\to E^{h}$
is $\tilde{O}(e^{-\frac{\eta_{f}}{h}}) $-unitary,
the operator 
$\delta_{[0,\tilde{o}(1)],[a,b],h}^{(p)}
:T_{\delta_{2}}\mathcal{U}_{\ell_{0}}^{(p),h}\to
F_{[0,\tilde{o}(1)],[a,b],h}^{(p+1)}$ is, as 
$\Pi_{G^{h,(p+1)}}d_{f,h}^{(p)}T_{\delta_{2}}:\mathcal{U}_{\ell_{0}}^{(p),h}\to
G^{h,(p+1)}$\,, 
 one to one with
singular values all logarithmically equivalent to
$e^{-\frac{\ell_{0}}{h}}$\,.
In particular, for all $j=(\alpha,x_{\alpha}^{(p)})\in
\mathcal{X}^{(p)}(a,b)$ such that
$y_{\alpha}^{(p+1)}-x_{\alpha}^{(p)}=\ell_{0}$\,, we must have
$$
\|\delta_{[0,\tilde{o}(1)],[a,b],h}^{(p)}\Pi_{[0,\tilde{o}(1)],[a,b],h}T_{\delta_{2}}\varphi_{j}^{h}\|\stackrel{log}{\sim}e^{-\frac{\ell_{0}}{h}}\,.
$$
From 
the previous estimates, the new
family of vectors $(u_{j}^{h})$ defined by
$$
u_{j}^{h}=(1-\sum_{\ell>\ell_{0}}\Pi_{F_{\ell}^{(p),h}})\Pi_{[0,\tilde{o}(1)],[a,b],h}T_{\delta_{2}}\varphi_{j}^{h}
$$
and indexed by $j=(\alpha,x_{\alpha}^{(p)})\in \mathcal{X}^{(p)}(a,b)$\,,
$y_{\alpha}^{(p+1)}-x^{(p)}_{\alpha}=\ell_{0}$
satisfies
\begin{eqnarray}
\label{eq.u-j}
  &&
     \langle u_{j}^{h}\,,
     \delta_{[0,\tilde{o}(1)],[a,b],h}^{*,(p)}\delta_{[0,\tilde{o}(1)],[a,b],h}^{(p)}u_{j}^{h}\rangle
\stackrel{log}{\sim}e^{-2\frac{\ell_{0}}{h}}\,,\\
\label{eq.u-j'}
 && u_{j}^{h}\perp 
\text{Ran}~1_{[0,e^{-2\frac{\ell_{0}+\delta_{f}}{h}}[}(\delta_{[0,\tilde{o}(1)],[a,b],h}^{*,(p)}\delta_{[0,\tilde{o}(1)],[a,b],h}^{(p)})   \,,\\
\text{and}\quad && 
\label{eq.u-j''}
 \|u_{j}^{h}-\varphi_{j}^{h}\|=\tilde{O}(e^{-\frac{\delta_{f}}{h}})
\,.
\end{eqnarray}
Note that
\eqref{eq.u-j} and \eqref{eq.u-j'}
follow easily from the definition of the family $(u_{j}^{h})$\,, while
\eqref{eq.u-j''},
which also implies the
$\tilde{O}(e^{-\frac{\delta_{f}}{h}})$-orthonormality of the family $(u_{j}^{h})$\,, follows from
   \eqref{eq.Pi-T} together with the estimate, for 
  $\ell>\ell_{0}$ and 
   $j=(\alpha,x_{\alpha}^{(p)})\in \mathcal{X}^{(p)}(a,b)$\,,
$y_{\alpha}^{(p+1)}-x^{(p)}_{\alpha}=\ell_{0}$\,,
\begin{align*}
\Pi_{F_{\ell}^{(p),h}}\varphi_{j}^{h}
&=  \big(\Pi_{F_{\ell}^{(p),h}}-\Pi_{F_{\ell}^{(p),h}}\Pi_{\mathcal U_{\ell}^{(p),h}}\big)\varphi_{j}^{h}
+ \Pi_{F_{\ell}^{(p),h}}\Pi_{\mathcal U_{\ell}^{(p),h}}\varphi_{j}^{h}\\
&= 
\tilde{O}(e^{-\frac{\delta_{f}}{h}}) +
\tilde{O}(e^{-\frac{\eta_{f}}{8h}})
\leq \tilde{O}(e^{-\frac{\delta_{f}}{h}})\,,
\end{align*}
where the last line follows from the induction hypothesis
and from the $\tilde{O}(e^{-\frac{\eta_{f}}{8h}})$-orthonormality of 
the family $(\varphi_{j}^{h})_{j\in \mathcal{J}^{(p)}(a,b)}$\,.
The relations 
\eqref{eq.u-j} and \eqref{eq.u-j''}
 imply that the vector
$$
v_{j}^{h}=1_{[0,e^{-2\frac{\ell_{0}-\delta_{f}}{h}}]}(\delta_{[0,\tilde{o}(1)],[a,b],h}^{*,(p)}\delta_{[0,\tilde{o}(1)],[a,b],h}^{(p)})u_{j}^{h}
$$
satisfies
$$
\|v_{j}^{h}-u_{j}^{h}\|=\tilde{O}(e^{-\frac{\delta_{f}}{h}})\quad\text{and thus}\quad \|v_{j}^{h}-\varphi_{j}^{h}\|=\tilde{O}(e^{-\frac{\delta_{f}}{h}})\,,
$$
while \eqref{eq.u-j'} yields
$$
v_{j}^{h}\in F_{\ell_{0}}^{(p),h}.
$$
Hence, we have proved
$\vec{d}(\mathcal{U}_{\ell_{0}}^{(p),h},F_{\ell_{0}}^{(p),h})=\tilde{O}(e^{-\frac{\delta_{f}}{h}})$
and thus, using 
$$
\dim
\mathcal{U}_{\ell_{0}}^{(p),h}=\sharp\left\{j=(\alpha,x_{\alpha}^{(p)})\in
\mathcal{X}^{(p)}(a,b)\,,
y_{\alpha}^{(p+1)}-x_{\alpha}^{(p)}=\ell_{0}\right\}=\dim F_{\ell_{0}}^{(p),h}\,,
$$
implies $\vec{d}(F_{\ell_{0}}^{(p),h},\mathcal{U}_{\ell_{0}}^{(p),h})+\vec{d}(\mathcal{U}_{\ell_{0}}^{(p),h},F_{\ell_{0}}^{(p),h})=\tilde{O}(e^{-\frac{\delta_{f}}{h}})$\,.
This ends the proof of Proposition~\ref{pr:singvdelta}.
\end{proof}

Now quasimodes have been
 constructed for $d_{f,f^{-1}([a,b]),h}$\,, the dual version can be
 given. Remember that 
$$
d_{f,h}^{*}=
(-1)^{\text{deg}}\star^{-1}e^{\frac{f}{h}}(hd) e^{-\frac{f}{h}}\star 
$$
and the construction of $\delta_{1}$-quasimodes for
$d_{f,f^{-1}([a,b],h)}^{*}$ is equivalent to the construction of 
$\delta_{1}$-quasimodes for $d_{-f,(-f)^{-1}([-b,-a]),h}$\,, where the
fiber bundle $\Lambda T^{*}M$ is replaced by $\Lambda T^{*}M\otimes
\mathrm{or}_{M}$\,. Accordingly, the degree $p$ is changed into
$d-p$\,,  the order of critical values is reversed and,
in the interval $[a,b]$\,, the role of lower and upper 
endpoints in the sets  $\mathcal{X}^{*}(a,b)$ and
$\mathcal{Y}^{*}(a,b)$ are interchanged.
\begin{definition}
\label{de:dualquasi}
Under Hypothesis~\ref{hyp:cN} and with $\delta_{1}\in
]0,\frac{\eta_{f}}{8}]$\,, a dual $\delta_{1}$-family of quasimodes
denoted by $\left(\widehat{\varphi_{j}^{*,h}}\right)_{j\in
  \mathcal{J}(a,b)}$ is defined like the family
$\left(\varphi_{j}^{*,h}\right)_{j\in \mathcal{J}(a,b)}$ in Definition~\ref{de:adapted} after
replacing:
\begin{itemize}
\item $d_{f,f^{-1}([a,b]),h}$ by $d_{f,f^{-1}([a,b]),h}^{*}$\,,
\item
  $I_{j}^{h}=[x_{\alpha}^{(p)}-\delta_{1},y_{\alpha}^{(p+1)}-\gamma(h)]$
  when $j=(\alpha,x_{\alpha}^{(p)})\in\mathcal{X}^{(p)}(a,b)$ 
by
$$
\widehat{I_{j}^{h}}=
[x_{\alpha}^{(p-1)}+\gamma(h),y_{\alpha}^{(p)}+\delta_{1}]\quad\text{when}\quad j=(\alpha,y_{\alpha}^{(p)})\in \mathcal{Y}^{(p)}(a,b)\,,
$$
\item and
  $I_{j}^{h}=[\tilde{c}-\delta_{1},b]$ when $j=(\alpha,\tilde{c})\in
  \mathcal{Y}^{(p)}(a,b)\sqcup \mathcal{Z}^{(p)}(a,b)$ by
$$
\widehat{I_{j}^{h}}=
[a,\tilde{c}+\delta_{1}]
\quad
\text{when}\quad j=(\alpha,\tilde{c})\in\mathcal{X}^{(p)}(a,b)\sqcup \mathcal{Z}^{(p)}(a,b)\,.
$$
\end{itemize}
Finally, the truncation operator $T_{\delta_{2}}$ introduced for
$\delta_{2}\in ]0,\frac{\eta_{f}}{8}]$ in Definition~\ref{de:Td2} has to
be replaced by $\widehat{T_{\delta_{2}}}$ defined by
\begin{eqnarray*}
&&\widehat{T_{\delta_{2}}}\widehat{\varphi_{j}^{(p),h}}=
\left\{
    \begin{array}[c]{ll}
\widehat{\chi_{x_{\alpha}^{(p-1)},\delta_{2}}}\widehat{\varphi_{j}^{(p),h}}&\text{if}~j=(\alpha,y_{\alpha}^{(p)})\in
\mathcal{Y}^{(p)}(a,b)
\\
\widehat{\varphi_{j}^{(p),h}}&\text{if}~j\in \mathcal{X}^{(p)}(a,b)\cup
\mathcal{Z}^{(p)}(a,b)\,,
    \end{array}
\right.
\\
\text{where}&&
\widehat{\chi_{\tilde{c},\delta_{2}}}(x)=\widehat{\chi}\left(\frac{f(x)-\tilde{c}}{\delta_{2}}\right)\,,
\end{eqnarray*}
for 
a fixed $\widehat{\chi}\in
\mathcal{C}^{\infty}(\rz;[0,1])$ such that $\widehat{\chi}\equiv 1$ on
$[2,+\infty[$ and $\mathrm{supp}~\widehat{\chi}\subset]1,+\infty[$\,.
\end{definition}
\begin{thm}
\label{th:eigensp}
Like in Theorem~\ref{th:specres}, assume Hypothesis~\ref{hyp:mainf} or (Hypothesis~\ref{hyp:Lipbar} and
Hypothesis~\ref{hyp:AgmonLip}) for a more general $f$\,, which is
equivalent to Hypothesis~\ref{hyp:cN} when the definition of
$\eta_{f}>0$ is added. Let $a,b\not\in
\left\{c_{1},\ldots, c_{N_{f}}\right\}$ and let
$\Delta_{f,f^{-1}([a,b]),h}=\mathop{\oplus}_{p=0}^{d}\Delta_{f,f^{-1}([a,b]),h}^{(p)}$ 
be defined like in
Proposition~\ref{pr:domain} with 
$N_{t}=f^{-1}(\left\{a\right\})$ and
$N_{n}=f^{-1}(\left\{b\right\})$\,. We set, like in
Proposition~\ref{pr:singvdelta},
 \begin{eqnarray*}
   &&L^{(p)}=\left\{b_{\alpha}^{(p+1)}-a_{\alpha}^{(p)}, \alpha\in
   A_{c}^{(p)}(a,b)\right\}\\
\text{and}
&&\delta_{f}=\min(\frac{\eta_{f}}{8}, \frac{|\ell-\ell'|}{8},
\ell\neq \ell'\in L^{(p)})>0\,.
 \end{eqnarray*}
 The $\delta_{1}$-family of quasimodes $(\varphi_{j}^{*,h})_{j\in
   \mathcal{J}(a,b)}$ is given
by Theorem~\ref{th:induc}
 with $\delta_{1}=\frac{\eta_{f}}{8}$\,,   
    and its dual version
 $(\widehat{\varphi_{j}^{*,h}})_{j\in \mathcal{J}(a,b)}$ by
 Definition~\ref{de:dualquasi}. For $\ell\in L^{(p)}$\,, we define lastly
 \begin{eqnarray*}
   && 
\overline{\mathcal U}_{\ell}^{(p),h}:=\mathcal{U}_{\ell}^{(p),h}\mathop{\oplus}\widehat{\mathcal{U}}_{\ell}^{(p),h}\,,\\
\text{where}&&
\mathcal{U}_{\ell}^{(p),h}=\Vect\left(\varphi_{j}^{h}\,,~
  j=(\alpha,x_{\alpha}^{(p)})\in
              \mathcal{X}^{(p)}(a,b)\,,~y_{\alpha}^{(p+1)}-x_{\alpha}^{(p)}=\ell\right)\\
\text{and}&&
\widehat{\mathcal{U}}_{\ell}^{(p),h}=\Vect\left(\widehat{\varphi_{j}^{h}}\,,~
  j=(\alpha,y_{\alpha}^{(p)})\in \mathcal{Y}^{(p)}(a,b)\,,~y_{\alpha}^{(p)}-x_{\alpha}^{(p-1)}=\ell\right)\,.
 \end{eqnarray*}
Then, for every $\ell\in L^{(p)}$\,, the space $\overline{\mathcal  U}_{\ell}^{(p),h}$ is close to
$F^{(p)}_{[e^{-2\frac{\ell+\delta_{f}}{h}},e^{2\frac{\ell-\delta_{f}}{h}}],[a,b],h}$
according to 
$$
\vec{d}\left(\overline{\mathcal U}_{\ell}^{(p),h},F^{(p)}_{[e^{-2\frac{\ell+\delta_{f}}{h}},e^{2\frac{\ell-\delta_{f}}{h}}],[a,b],h}\right)
+\vec{d}\left(F^{(p)}_{[e^{-2\frac{\ell+\delta_{f}}{h}},e^{2\frac{\ell-\delta_{f}}{h}}],[a,b],h},
\overline{\mathcal U}_{\ell}^{(p),h}\right)=\tilde{O}(e^{-\frac{\delta_{f}}{h}})\,.
$$
\end{thm}
\begin{proof}
Let us first recall  the relation
$$
\Delta_{f,f^{-1}([a,b]),h}^{(p)}\Pi_{[0,\tilde{o}(1)],[a,b],h}^{(p)}=
\underbrace{\delta_{[0,\tilde{o(1)}],[a,b],h}^{(p-1)}\delta_{[0,\tilde{o}(1)],[a,b],h}^{(p-1),*}}_{A}
+\underbrace{\delta^{(p),*}_{[0,\tilde{o}(1)],[a,b],h}\delta^{(p)}_{[0,\tilde{o}(1)],[a,b],h}}_{B}\,,
$$
where $A$ and $B$ are  self-adjoint and satisfy
$AB=BA=0$\,. We deduce from this observation and from
the Hodge decomposition that,
 for $\lambda_{h}\neq 0$\,,
$\lambda_{h}=\tilde{o}(1)$\,,
$$
\ker(\Delta_{f,f^{-1}([a,b]),h}-\lambda_{h})=
\ker(A-\lambda_{h})
\mathop{\oplus}^{\perp}
\ker(B-\lambda_{h})\,.
$$
Moreover Proposition~\ref{pr:singvdelta} says
$$
\vec{d}(\mathcal{U}_{\ell}^{(p),h},F_{\ell}^{(p),h})
+\vec{d}(F^{(p),h}_{\ell},\mathcal{U}^{(p),h}_{\ell})
=\tilde{O}(e^{-\frac{\delta_{f}}{h}})\,,
$$
where
$$
F_{\ell}^{(p)}=\mathop{\oplus}^{\perp}_{e^{-2\frac{\ell+\delta_{f}}{h}}\leq
  \lambda_{h}\leq
e^{-2\frac{\ell-\delta_{f}}{h}}} \ker(B-\lambda_{h})\,.
$$
The proximity of $\widehat{\mathcal{U}}^{(p),h}_{\ell}$ to
$\mathop{\oplus}^{\perp}_{e^{-2\frac{\ell+\delta_{f}}{h}}\leq \lambda_{h}\leq
e^{-2\frac{\ell-\delta_{f}}{h}}} \ker(A-\lambda_{h})$ is the dual version.
\end{proof}
\begin{remark}
  The last result about the eigenvectors of
  $\Delta_{f,f^{-1}([a,b]),h}$ arouses several comments.
  \begin{itemize}
  \item When there is a single bar $\alpha\in A_{c}^{(p)}(a,b)$ with length
    $\ell$\,, then $\Delta^{(p)}_{f,f^{-1}([a,b]),h}$
    (resp. $\Delta_{f,f^{-1}([a,b]),h}^{(p+1),h}$)  has one
    eigenvector associated with the eigenvalue
    $\lambda_{h}\stackrel{log}{\sim}e^{-\frac{2\ell}{h}}$ 
localized around $f^{-1}(x_{\alpha}^{(p)})$
(resp. $f^{-1}(y_{\alpha}^{(p+1)})$) and
    $\tilde{O}(e^{-\frac{\delta_{f}}{h}})$-close to the corresponding
    quasimode $\varphi_{j}^{(p),h}$
    (resp. $\widehat{\varphi_{j}^{(p+1),h}}$) with
    $j=(\alpha,x_{\alpha}^{(p)})\in \mathcal{X}^{(p)}(a,b)$
    (resp. $j=(\alpha,y_{\alpha}^{(p+1)})\in
    \mathcal{Y}^{(p+1)}(a,b)$).
\item Once we have approximated the eigenvectors associated with the non zero eigenvalues
  by the quasimodes $\varphi_{j}^{h}$ or
  $\widehat{\varphi_{j}^{h}}$\,, one can recover an approximate
  description of $\ker(\Delta_{f,f^{-1}([a,b]),h})$ by considering
  a basis of $\Vect(\varphi_{j}^{h}, j\in
  \mathcal{Y}(a,b)\sqcup \mathcal{Z}(a,b))$ whose elements are
  $\tilde{O}(e^{-\frac{\delta_{1}}{h}})$-orthogonal to all the
  $\widehat{\varphi_{j'}^{h}}$\,, $j'\in \mathcal{Y}(a,b)$\,. 
\item Actually, the description of the eigenvectors with a
  $\tilde{O}(e^{-\frac{\delta_{f}}{h}})$ error in the $L^{2}$-norm is
  much less precise than what we were able to do with the quasimodes
  $\varphi_{j}^{*,h}$\,, with a wide range control of the exponential
  decay estimates. We also know from the proof of
  Theorem~\ref{th:induc}, and this is again   illustrated in the proof of
  Proposition~\ref{pr:singvdelta}, that working with the family of
  quasimodes $(\varphi_{j}^{h})_{j\in \mathcal{J}(a,b)}$ is much more
  flexible and informative than working with the eigenvectors of
  $\Delta_{f,f^{-1}([a,b]),h}$\,. Note specifically,
  in the proof of
  Proposition~\ref{pr:singvdelta},
   the use of the
  orthogonality $G_{n}^{h}\perp G_{n'}^{h}$ for $n\neq n'$ in the
  separation of the different exponential scales associated with
  the different lengths of bars. This really relies on the fact that
  $G^{h}$ is made of kernels of separated local problems. Such an
  exact property is completely lost if we use instead the full
  spectral space $F_{[0,\tilde{o}(1)],[a,b],h}$\,.
\item From the modeling interpretation, it is interesting to note that
  the quasimodes $(\varphi_{j}^{h})_{j\in \mathcal{J}(a,b)}$
carry the same heuristic as the true eigenvectors
  for small times although they do not belong to
  $D(\Delta_{f,f^{-1}([a,b]),h})$\,. For simplicity, assume that there
  is a single bar $\alpha\in A_{c}^{(p)}(a,b)$ with length
  $\ell$\,. Then $\varphi_{j}^{h}$\,,
  $j=(\alpha,x_{\alpha}^{})\in \mathcal{X}(a,b)$\,, satisfies
\begin{align*}
\|e^{-t\Delta_{f,f^{-1}([a,b],h)}}\varphi_{j}^{h}-e^{-t\lambda_{h}}\varphi_{j}^{h}\|&=\|(e^{-t\Delta_{f,f^{-1}([a,b],h)}}-e^{-t\lambda_{h}})(\varphi_{j}^{h}-u_{h})\|
\\
&\leq 2\|\varphi_{j}^{h}-u_{h}\|=\tilde{O}(e^{-\frac{\delta_{f}}{h}})\,,
\end{align*}
where $u_{h}$ is the unitary eigenvector associated with the eigenvalue
$\lambda_{h}\stackrel{log}{\sim}e^{-2\frac{\ell}{h}}$\,.
In particular, $e^{-t\Delta_{f,f^{-1}([a,b]),h}}\varphi_{j}^{h}\sim
e^{-t\lambda_{h}}\varphi_{j}^{h}$ makes sense for times longer 
 than the lifetime $\frac{1}{\lambda_{h}}\stackrel{log}{\sim}e^{2\frac{\ell}{h}}$ of the metastable 
state $u_{h}$ as $h\to 0$\,.
  \end{itemize}
\end{remark}
\subsection{Stability theorem}
\label{sec:firststab}
The following stability theorem, of which a simple version, Corollary~\ref{cor:mainsimple}, was given in the
introduction, is a direct consequence
of  Theorem~\ref{th:specres} and of the topological stability result
$$
d_{bot}(\mathcal{B}(f), \mathcal{B}(g))\leq \|f-g\|_{\mathcal{C}^{0}}
$$
recalled in Appendix~\ref{sec:stab}.

\begin{thm}
\label{th:stab1}
In the framework of Theorem~\ref{th:specres}, namely
Hypothesis~\ref{hyp:mainf}, or (Hypothesis~\ref{hyp:Lipbar} and
Hypothesis~\ref{hyp:AgmonLip}) for a more general Lipschitz function $f$\,, and $a,b\not\in
\{c_{1},\ldots, c_{N_{f}}\}$\,, set 
$$
\ell_{min}\ :=\ \min\Big(\left\{y_{\alpha}-x_{\alpha},
  \alpha\in A_{c}(a,b)\right\}
\cup \left\{|c_{n}-b|, |c_{n}-a|, 1\leq n\leq N_{f}\right\}\Big)
\,,
$$
where $A_{c}(a,b)=A_{c}(f;a,b)$ is the set
defined in \eqref{eq:defAcab} for the function $f$\,, that is indexing the bars
of $f$ with two endpoints in $]a,b[$\,.\\
Let moreover $g$ be  any other function satisfying
Hypothesis~\ref{hyp:mainf}, or (Hypothesis~\ref{hyp:Lipbar} and
Hypothesis~\ref{hyp:AgmonLip}),
as well as
$$
\|g-f\|_{\mathcal{C}^{0}}\ <\ \frac{\ell_{min}}{4}\,,
$$ 
and such that
$a,b$ do not belong to the set $\{c'_{1},\ldots,c'_{N_{g}}\}$
made of its ``critical values''.\\
Then, the
$\tilde{O}(e^{-\frac{\ell_{min}}{h}})$ non zero eigenvalues of
$\Delta_{g,g^{-1}([a,b]),h}^{(p)}$ can be labelled
$$
\lambda_{\alpha}^{(p)}(g;h)\quad,\quad \alpha\in
A^{(p)}_{c}(a,b) \sqcup  A^{(p-1)}_{c}(a,b)\,,
$$
with, for every $\alpha\in A^{(p)}_{c}(a,b) \sqcup  A^{(p-1)}_{c}(a,b)$\,, 
$$
\ell_{min}
\ <\ 2(y_{\alpha}^{*+1}-x_{\alpha}^{*})
-4\|g-f\|_{\mathcal{C}^{0}}
\ \leq\  \lim_{h\to 0}
-h \log (\lambda_{\alpha}^{(p)}(g,h))
\ \leq\  2(y_{\alpha}^{*+1}-x_{\alpha}^{*})+4\|g-f\|_{\mathcal{C}^{0}}\,.
$$
Meanwhile, for ${\bf f}=f$ or ${\bf f}=g$\,,  the dimension
$\dim \ker (\Delta_{{\bf f},{\bf f}^{-1}([a,b]),h}^{(p)})$ equals $\beta^{(p)}({\bf f}^{b},{\bf f}^{a})$\,,
and thus
$$\dim\ker (\Delta_{f,f^{-1}([a,b]),h}^{(p)})\ =\ \ker (\Delta_{g,g^{-1}([a,b]),h}^{(p)})
\quad \text{if and only if}\quad \beta^{(p)}(f^{b},f^{a})=\beta^{(p)}(g^{b},g^{a})\,.$$
\end{thm}
\begin{proof}
After possibly adding empty bars, the bar codes associated with $f$
and $g$ can be written
 $\mathcal{B}(f)=([a_{\alpha},b_{\alpha}[)_{\alpha\in A}$ and
$\mathcal{B}(g)=([c_{\alpha},d_{\alpha}[)_{\alpha\in A}$\,, where 
$$
\max\left\{|a_{\alpha}-c_{\alpha}|, |d_{\alpha}-b_{\alpha}|\,,
  \alpha\in A, b_{\alpha}<+\infty\right\}
\leq d_{bot}(\mathcal{B}(g),\mathcal{B}(f))\leq
\|g-f\|_{\mathcal{C}^{0}}<\frac{\ell_{min}}{4}\,.
$$
The definition 
$$
\ell_{min}:=
\min(\left\{y_{\alpha}-x_{\alpha}, \alpha\in A_{c}(f;a,b)\right\}\cup
\left\{|c_{n}-a|, |c_{n}-b|, 1\leq n\leq N_{f}\right\})
$$
implies that the number of bars $\alpha\in A_{c}(g;a,b)$ such that
$y_{\alpha}-x_{\alpha}>\frac{\ell_{min}}{2}$\,, for the function
$g$\,, is in bijection with the whole set of bars
$A_{c}(f; a,b)$ for the function $f$\,, which is made by assumption of bars not smaller than
$\ell_{min}$\,. The other potential bars of $A_{c}(g;a,b)$ have  a length 
strictly smaller than~$\frac{\ell_{min}}{2}$\,.\\
Moreover,  for $\alpha\in A_{c}(a,b)$\,, the expression of $\lim_{h\to
  0}-h\log\lambda_{\alpha}^{(p)}(h)$  given
in Theorem~\ref{th:specres}, respectively applied with $g$ and  $f$\,, provides the
inequalities for the $\tilde{O}(e^{-\frac{2(\ell_{min}/2)}{h}})$ non zero 
eigenvalues
of $\Delta_{g,g^{-1}([a,b]),h}$\,.\\
Finally, the last statement of Theorem~\ref{th:stab1} is a direct consequence of
the comments made in the second item of Remark~\ref{re:rem1}.
\end{proof}

\section{Generalizations}
\label{sec:regass}
Our proofs are definitely done under Hypothesis~\ref{hyp:mainf}, while,
for a more general Lipschitz function~$f$\,, consequences of
Hypothesis~\ref{hyp:Lipbar} have not yet been checked and the
exponential decay estimates of Propositions~\ref{pr:Agmon} and~\ref{pr:Agmon1} 
have simply been replaced by
assumptions. \\
This framework was chosen in order to put the stress on the
essentially one-dimensional analysis on $\rz\supset f(M)$\,. 
Once this is well understood, it
is rather easy to adapt the analysis and the results in order to consider
more general domains, manifolds, or Lipschitz functions $f$\,. The
first generalizations will be presented for the sake of simplicity in
the framework of Hypothesis~\ref{hyp:mainf}.\\
Additionally,  we will check that Hypothesis~\ref{hyp:Lipbar} and
Hypothesis~\ref{hyp:AgmonLip} hold true under the simple assumption
that $f$ is a subanalytic  Lipschitz function (see
Hypothesis~\ref{hyp:realana}), which describes,
 in some sense, a wider class of functions than 
 Hypothesis~\ref{hyp:mainf} in a real analytic
geometry.

\subsection{More general domains}
\label{sec:moregendom}

It is not difficult to adapt all the analysis to some simple cases
when the geometrical domain $\overline{\Omega}$  differs from
$f^{-1}([a,b])$\, by tamed deformations of $\partial \Omega$\,.
\begin{prop}
\label{pr:moregendom}
Let $(M,g)$ be a compact Riemannian manifold and let $f$ satisfy
Hypothesis~\ref{hyp:mainf}.
If there exist $m_{0},n_{0}\in \left\{1,\ldots, N_{f}\right\}$
such that
$m_{0}<n_{0}$ and the boundary of the domain 
$\overline{\Omega}=\Omega \sqcup N_{t}\sqcup N_{n}$ satisfies
\begin{eqnarray*}
  &&
f(N_{t})\subset ]c_{m_{0}},c_{m_{0}+1}[\quad,\quad f(N_{n})\subset ]c_{n_{0}},c_{n_{0}+1}[\,,\\
\text{and}
&&\frac{\partial f}{\partial n}\big|_{N_{t}}<0\quad,\quad
   \frac{\partial f}{\partial n}\big|_{N_{n}}>0\,.
\end{eqnarray*}
then all  the results or Theorem~\ref{th:induc} hold true with
$\tilde{c}_{1}=c_{m_{0}+1}$\,, $\tilde{c}_{N}=c_{n_{0}}$ when
$\eta_{f}$ is chosen in the interval
\begin{eqnarray*}
&&0< \eta_{f}<\frac{1}{2}\min_{1<n\leq
   N_{f}}c_{n}-c_{n-1}\,,\\\text{and}
&& \eta_{f}<\min_{x\in N_{t}} (c_{m_{0}+1}-f(x))\quad,\quad \eta_{f}<\min_{x\in
   N_{n}}(f(x)-c_{n_{0}})\,.
\end{eqnarray*}
\end{prop}
\begin{proof}
All the proof of Theorem~\ref{th:induc} relies on the construction of
the $\delta_{1}$-family of quasimodes $(\varphi_{j}^{h})_{j\in
  \mathcal{J}(a,b)}$ when $\overline{\Omega}=f^{-1}([a,b])$\,.
We fix $a=c_{m_{0}+1}-\eta_{f}=\tilde{c}_{1}-\eta_{f}$ and
$b=c_{n_{0}}+\eta_{f}=\tilde{c}_{N}+\eta_{f}$\,.
Because the gradient lines provide a homotopy between the pairs
$(\overline{\Omega}, N_{t})$ and $(f^{-1}([a,b]),
f^{-1}\left\{a\right\})$\,, the bar code for $f$ in
$\overline{\Omega}$ relatively to $N_{t}$ can be
identified with  $\mathcal{B}(f;[a,b])$\,. Now, the quasimodes
$(\varphi_{j}^{h})_{j\in \mathcal{J}(a,b)}$ are extended by $0$ in
$f^{a}\cap \overline{\Omega}$ and, when $j\in \mathcal{Y}(a,b)\cup
\mathcal{Z}(a,b)$, they are ``extended'' in $f_{b}\cap \Omega$ as
$$
\chi \varphi_{j}^{h}-d_{f,f^{-1}([c_{n_{0}}+\delta_{1}, +\infty[\cap
  \overline{\Omega}),h}^{*}(\Delta_{f,f^{-1}([c_{n_{0}}+\delta_{1},+\infty[\cap
  \overline{\Omega}),h})^{-1}(hd\chi\wedge\varphi_{j}^{h})\,,
$$ 
like in Proposition~\ref{pr:interdfh}-\textbf{ii)}, with
$\delta_{1}\in]0,\frac{\eta_{f}}{8}]$,
 $\chi\in
\mathcal{C}^{\infty}(M;[0,1])$, $\chi\equiv 1$ in
$f^{b-\eta_{f}/2}$\,, $\chi\equiv 0$ in $f_{b-\eta_{f}/4}$, and where
Dirichlet (resp. Neumann) boundary conditions are put
on 
$f^{-1}(\left\{c_{n_{0}}+\delta_{1}\right\})$ 
(resp. on $N_{n}$), for the
domain 
$f^{-1}([c_{n_{0}}+\delta_{1},+\infty[\cap \overline{\Omega})$
\end{proof}
\begin{remark}
  Another interesting case is when the Neumann boundary conditions
  on $N=N_{n}$, where $\frac{\partial f}{\partial n}>0$, are replaced by
  Dirichlet boundary conditions. Then, generalized critical values
  corresponding to 
  critical values of $f\big|_{N}$
 appear following what is known for a Morse function~$f$ (see
 e.g. \cite{ChLi,HeNi,Lep3,LeNi,Lau}). As a topological tool,
 bar codes make sense for boundary manifolds.  But the analysis has to
 be reconsidered from the beginning,
 especially by introducing mixed Dirichlet-Neumann
 problems along the upper boundary of $\overline{\Omega}\cap f^{\leq
   t}$. We do not develop this point here
   (see however \cite{DLLN1} where such conditions are considered).
\end{remark}

\subsection{More general manifolds}
\label{sec:moregenmfld}

The following generalization aims at including the particular case
when $M=\rz^{d}$ is not compact and the gradient of $f$ dos not vanish outside
 a compact set. More specifically, we assume
\begin{hyp}
\label{hyp:noncompact}
Let $(M,g)$ be a complete Riemannian manifold and assume $f\in
\mathcal{C}^{\infty}(M;\rz)$ for the sake of simplicity. We suppose
that there exist $-\infty<a_{0}<b_{0}<+\infty$ and $\kappa>0$  such
that
\begin{itemize}
\item $K_{0}=f^{-1}([a_{0},b_{0}])$ is compact,
\item for all $x\in M\setminus K_{0}$\,, $|\nabla f(x)|\geq \kappa$,
\item $f$ has a finite number of critical values $c_{1},\ldots,
  c_{N_{f}}$ in $[a_{0},b_{0}]$ which belong to  $]a_{0},b_{0}[$.
\end{itemize}
\end{hyp}
Under this assumption, the definition of the bar code
$\mathcal{B}(f)=([a_{\alpha}^{*},b_{\alpha}^{*+1}[)_{\alpha\in A}$ is
essentially the same
as in the compact case, except that bars with $a_{\alpha}^{*}=-\infty$ and
$b_{\alpha}^{*+1}\in \rz$ are possible, according to the topology of $f^{t}$
as $t\to -\infty$\,. In such a case, $b_{\alpha}^{*+1}\in
\mathcal{Z}^{*+1}(a,b)$ for all $a,b\in
[-\infty,+\infty]\setminus\left\{c_{1},\ldots,c_{N_{f}}\right\}$ such that
$a<b_{\alpha}^{*+1}<b$\,.\\
The domain $f^{-1}([a,b])\subset M$ is actually $f^{-1}([a,b]\cap
]-\infty,+\infty[)$ when $a=-\infty$ or $b=+\infty$. Accordingly,
$\Delta_{f,f^{-1}([a,b]),h}$\,, $d_{f,f^{-1}([a,b]),h}$, and $d^{*}_{f,f^{-1}([a,b]),h}$
do not include   boundary conditions on the
infinite end in the definition of their domains.  
\begin{prop}
\label{pr:moregenmfld}
Under Hypothesis~\ref{hyp:noncompact}, all the results of
Theorem~\ref{th:induc} still hold.
\end{prop}
\begin{proof}
  The completeness of the manifold ensures that the scalar Laplacian
  is essentially self-adjoint on
  $\mathcal{C}^{\infty}_{0}(M)$. Adapting the proof of Simader's
  theorem ensures that $\Delta_{f,h}$ is essentially self-adjoint on
  $M$ and that $\Delta_{f,f^{-1}([a,b]),h}$ is essentially
  self-adjoint on the subspace of
  $\mathcal{C}^{\infty}_{0}(f^{-1}([a,b]);\Lambda T^{*}M)$ containing the
  boundary conditions, of Dirichlet type on $f^{-1}(\left\{a\right\})$ when
  $-\infty<a$ and of Neumann type on $f^{-1}(\left\{b\right\})$ when
  $b<+\infty$\,.\\
Agmon estimates and the compactness of $K_{0}=f^{-1}([a_{0},b_{0}])$
with $|\nabla f|\geq \kappa>0$ in $M\setminus K_{0}$ implies that
the solutions to
$\Delta_{f,f^{-1}([a,b]),h}\omega_{h}=\lambda_{h}\omega_{h}$ with
$\lambda_{h}\to 0$ as $h\to 0$ must satisfy
$$
\|e^{\frac{\kappa d_{g}(x, K_{0})}{h}}\omega_{h}\|_{W^{1,2}}\leq \tilde{O}(1)\,.
$$
One can then localize the analysis of exponentially small eigenvalues
to $K_{0}'=f^{-1}([a_{0}-1,b_{0}+1])$, which amounts to the case of a
compact manifold treated in Theorem~\ref{th:induc}.
\end{proof}

\subsection{More general Lipschitz functions}
\label{sec:moregenLip}
We consider more accurately the situation of a general Lipschitz
function~$f$\,, while the analysis was presented under conjectural
assumption. As a first step we recall in
Subsection~\ref{sec:moregenLipgen} how Hypothesis~\ref{hyp:Lipbar}
implies Hypothesis~\ref{hyp:weakreg} of Appendix~\ref{app:perscohom}
and therefore provides a finite bar code $\mathcal{B}_{f}$\,.\\
Once this is clarified we prove that Hypothesis~\ref{hyp:Lipbar} and
Hypothesis~\ref{hyp:AgmonLip} are
satisfied when $f$ is a subanalytic Lipschitz function, after the
suitable specification  of   the ``critical values''\,,
$c_{1}<\ldots<c_{N_{f}}$\,. It relies on the stratification of the
subanalytic graph of $f$\,, of which the properties are recalled in
Subsection~\ref{sec:moregenLipStrat}. A variation of Agmon~distance
will also be constructed after  solving the Hamilton-Jacobi
equation $|\nabla' \varphi|=|\nabla' f|$\,, where $\nabla'$ concerns only
tangential  partial
derivatives in some tubular neighoborhoods of every stratum. From this point of view, the analysis of
this Lipschitz subanalytic case, via a stratification technique, takes
some inspiration from \cite{GeNi}. Finally in  Subsubsection~\ref{sec:moregenLipAgm},
Hypothesis~\ref{hyp:AgmonLip} is checked to hold true, via some
partition of unity adapted with the stratification.

\subsubsection{Hypothesis~\ref{hyp:Lipbar} and consequences}
\label{sec:moregenLipgen}
The manifold $M$ is assumed to be compact without boundary although it
could be extended to more general cases like in
Subsection~\ref{sec:moregenmfld}.\\
Let us first define the critical
values  of a Lipschitz function $f$ or more exactly, its contrary.
\begin{definition}
\label{de:critLip}
When $f:M\to \rz$ is a Lipschitz function a value $a$ is not a critical
value  if for any $x_{0}\in f^{-1}(\left\{a\right\})$ there exists a
neighborhood $U_{x_{0}}$ 
of $x_{0}$ and a local coordinate system $x=(x^{1},x')\in
\rz\times \rz^{d-1}$ and a constant $C_{x_{0}}>0$ such that 
\begin{equation}
  \label{eq:minoLip}
\forall x=(x^{1},x'),y=(y^{1},x')\in U_{x_{0}}\,, \quad
\frac{1}{C_{x_{0}}}|x^{1}-y^{1}|\leq|f(x^{1},x')-f(y^{1},x')|\,.
\end{equation}
A critical value  $a\in f(M)\subset \rz$ is a point where the above
property fails.
\end{definition}

Since the function $f$ is continuous, the local condition condition
\eqref{eq:minoLip} can be replaced by
$$
\forall x=(x^{1},x'),y=(y^{1},x')\in U_{x_{0}}\,, \quad
\frac{1}{C_{x_{0}}}(x^{1}-y^{1})\leq f(x^{1},x')-f(y^{1},x')\quad\text{when}~x^{1}>y^{1}\,.
$$
Hypothesis~\ref{hyp:Lipbar} simply says that the Lipschitz function
$f$ has a finite number of critical values. But
the set $\left\{c_{1},\ldots,c_{N_{f}}\right\}$ of Hypothesis~\ref{hyp:Lipbar} may be strictly larger
that the set of critical values as defined above, and this a reason
why the values $c_{1},\ldots,c_{N_{f}}$ were called ``critical
values''. Actually this flexibility is especially usefull when we
consider subanalytic Lipschitz functions below. \\

\medskip
The above definition ensures that the implicit functions theorem in
the Lipschitz case can be applied locally  around $x\in f^{-1}(\left\{a,b\right\})$ with the following straightforward
consequences for the domain $f_{a}^{b}$ when $a,b$ are not ``critical
values'':
\begin{description}
\item[i)] $f_{a}^{b}$ is a strongly Lipschitz domain of $M$ according to
  the terminology of \cite{GMM}, meaning that
  it is locally the hypograph of a Lipschitz function in the proper
  coordinate system.
\item[ii)] $\overline{f^{b}_{a}}=f^{-1}([a,b])$\,.
\item[iii)] When $a=-\infty$\,, $f^{b}$ with $c<b<c'$ and no critical values
  in $]c,c'[$\,, is homotopic
  to $\Omega$ a $\mathcal{C}^{\infty}$ domain with $\partial
  \Omega\subset f_{c}^{c'}$\,.
\end{description}
The last statement can be checked by using
 finitely many local homotopies in
coordinate systems, but one could also use the
 global construction
of a
smooth transverse vector field
as proposed in \cite{Ver}-Theorem~1.12-vi).\\
The above three properties were used in our analysis. In particular
the finiteness of $N_{f}$ and \textbf{iii)} appear in
Hypothesis~\ref{hyp:weakreg} which allows the introduction of a finite
bar code $\mathcal{B}_{f}$\,. The properties \textbf{i)} and
\textbf{ii)} are used in the definition of
$\Delta_{f,f^{-1}([a,b]),h}$ according to
Proposition~\ref{pr:domain}\\

\medskip
Critical points and values can actually 
be defined in a coordinate free way, in terms of the
standard notion in non smooth analysis of Clarke's generalized
gradient and Clarke's critical points:
In $\rz^{d}$ or locally in a coordinate system in $M$\,, a Lipschitz
function admits a differential $df(x)$ almost every where by
Rademacher's theorem and the domain $\mathop{Dom}(df)$ is the set of $x$
where $df(x)$ exists. Clarke's generalized gradient at $x$ then equals
the closed convex set 
$$
\partial^{\circ} f(x)=\textrm{co}~\left\{\zeta\in\rz^{d}\,, \exists
  (x_{n})_{n\in\nz}\in \text{Dom}(df)^{\nz}\,, \lim_{n\to
    \infty}x_{n}=x\quad\text{and}\quad
  \lim_{n\to\infty}df(x_{n})=\zeta\right\}
$$
where $\textrm{co}$ denotes the convex hull\,. A Clarke critical point $x$
is a point where $0\subset \partial^{\circ}f(x)$ and a Clarke critical
value of $f$ is a value $a$ where $f^{-1}(\left\{a\right\})$ contains
a critical points. In the case of subanalytic Lipschitz functions
which will be considered more specificaly in the other paragraphs,
this definition actually coincides with the wavefront naturally
introduced in \cite{DeLe}. Staying at the local level the local
condition \eqref{eq:minoLip} for $x_{0}\in f^{-1}\left\{a\right\}$\,,
actually means that for all 
$x\in \mathop{Dom}(df)\cap U_{x_{0}}$\,, $df(x)$ lies in the intersection
of some closed salient ($\zeta\neq 0$ and
$-\zeta$ cannot both belong to it) convex cone $\mathcal{C}_{x_{0}}$ with a
shell $S_{x_{0}}=\left\{\zeta\in \rz^{d}\,, r<
|\zeta|\leq R\right\}$\,, $0<r<R<+\infty$\,. This writing is
equivalent to the fact that for all $x\in f^{-1}(\left\{a\right\})$\,,
Clarke's generalized gradient is included in the intersection of a
salient convex cone and a closed shell. This property is independent of the coordinate system and of the
 metric if we replace the differential $df$ by the gradient 
$\nabla f$\,.\\

Even in the subanalytic setting, those critical values (according to
Definition~\ref{de:critLip} or Clarke) may overestimate what the
intuition and even the final result would retain.
 Warga's example carefully analyzed in
\cite{CzRi}\,,
$$
f(x^{1},x^{2})=||x^{1}|+x^{2}|+\frac{1}{2}x^{1}\,,
$$
with the level curves in the picture below, satisfies the above
consequences \textbf{i)},\textbf{ii)} and \textbf{iii)}
 for any value $b\in \rz$ although $0$ is
a critical value of $f$\,. Note also that $0$ will be a critical value
of non well chosen
regularizations of $f$ and we refer to \cite{CzRi} for a thorough discussion 
of this point.
\begin{figure}[h]
\centering{
\includegraphics[width=7cm]{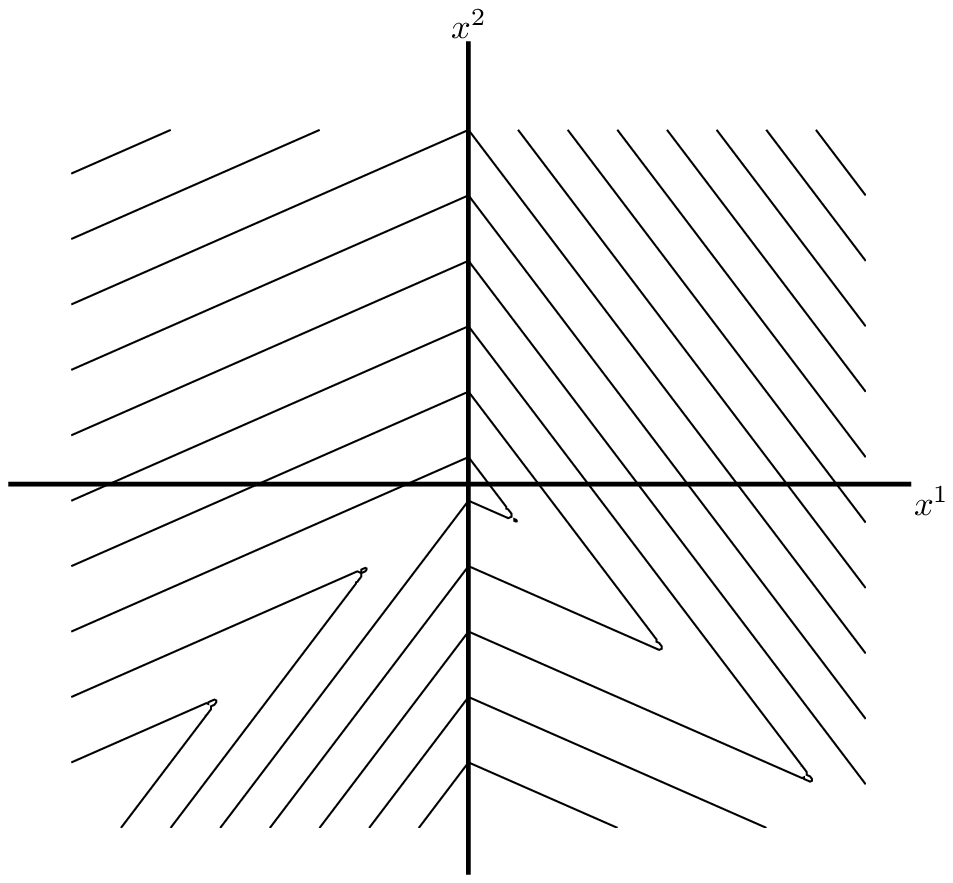}
\captionsetup{labelformat=empty,width=14cm}
\caption{\textbf{Figure 12:} Level curves of Warga's function
  $f(x^{1},x^{2})=||x^{1}|+x^{2}|+\frac{1}{2}x^{1}$\,.}
}
\end{figure}

\medskip

Actually in the subanalytic setting an even larger, but still finite,
set of values
$\left\{c_{1},\ldots, c_{N_{f}}\right\}$ will be introduced in order
to verify the second assumption, Hypothesis~\ref{hyp:AgmonLip}, used
in our analysis. 
\subsubsection{Stratification of Lipschitz subanalytic functions}
\label{sec:moregenLipStrat}
According to \cite{BDLS} a Lipschitz subanalytic function has a finite
number of  critical values and Hypothesis~\ref{hyp:Lipbar} holds true. We also
recalled in the previous paragraph that Clarke's gradient coincides with the wavefront set of
subanalytic Lipschitz functions introduced in \cite{DeLe}. However
such a notion of gradient or wavefront above a point $x\in M$\,, is
a wide closed convex set  which contains all the convex
combinations of limits of neighboring gradients without
discriminating the information which can be deduced from the
stratified structure. We specify the corresponding constructions when
$f$ is a real subanalytic Lipschitz function on a real analytic
compact Riemannian
manifold $M$ according to Hypothesis~\ref{hyp:realana}.\\
Let us first remind  the basic notions about subanalytic sets and
functions. We 
refer the reader to the founding articles \cite{Hardt}\cite{Hiro} and  to
\cite{Loja}\cite{BiMi} for a panoramic or historical presentation. 
A part but not all  of the material, presented or recalled here,
may be found in  \cite{DeLe} for the specific case of subanalytic Lipschitz
functions.\\
\medskip\\
\textbf{Review of subanalytic notions and results:}
\begin{itemize}
\item In the real analytic category, the class of subanalytic sets
is the one which contains the semianalytic sets, characterized by real
analytic equations or inequalities, and which is stable by finite set
operations (finite union, finite intersection and complement) and by proper
real analytic projections. The name ``subanalytic'' was introduced by
Hironaka and Hardt used the name ``analytic shadow'' in \cite{Hardt}
although they finally happened to describe the same class (see
\cite{Loja}).
\item Any subanalytic set $E$ of a real analytic manifold $X$ admits a stratification, that
  is a locally finite partition in real analytic connected submanifold of $X$
  called strata $E=\mathop{\sqcup}_{S\in \mathcal{S}} S$ such that
$S\cap \overline{S'}\neq \emptyset$\,, $S\neq S'$\,, implies
$S\subset \partial S'$  with $\dim S<\dim S'$\,, or equivalently because $\mathcal{S}$ is a
partition,  $S\cap \partial{S'}\neq \emptyset$\,, $S\neq S'$ implies
$S\subset \partial S'$ with $\dim S<\dim S'$\,.\\
Such a stratification can always be refined in order to satisfy
Whitney's local condition~B which reads in $\rz^{n}$ or in a
coordinate system:
$$
\left((x_{n})_{n\in \nz}\in (S')^{\nz}\quad,\quad \lim_{n\to
  \infty}x_{n}=x\in S\subset \overline{S'}\right)
\Rightarrow (T_{x}S\subset \lim_{n\to \infty}T_{x_{n}}S').
$$
When $\mathcal{C}$ is a locally finite family of subanalytic sets, the
stratification $\mathcal{S}$ can also be chosen in order to be
compatible with $\mathcal{C}$\,, which means that for all $S\in
\mathcal{S}$ and $C\in \mathcal{C}$\,, either $S\cap C=\emptyset$ or
$S\subset C$\,.
\item A subanalytic function $X\to Y$ is a function of which the graph
  is a subanalytic set of $X\times Y$\,.
\item When $f:X\to Y$ is a real analytic mapping, a stratification of
  $f$ is made of two stratifications $\mathcal{S}$ of $X$ and
  $\mathcal{F}$ of $Y$ such that 
$$
\forall S\in \mathcal{S}\,, f(S)\in \mathcal{F}\quad,\quad
\mathrm{rank}\,(f\big|_{S})=\dim f(S)\,.
$$
\item Corollary 4.4 of \cite{Hardt} assumes that $f:X\to Y$ is real
  analytic and $\mathcal{C}$ and $\mathcal{D}$ are two locally finite
  families of subanalytic sets of $X$ and $Y$ and $\Omega$ is a
  subanalytic open set such that $f\big|_{\overline{\Omega}}$ is
  proper. It then says that there exists a stratification
  $(\mathcal{S},\mathcal{F})$ of $f\big|_{\Omega}$ which is compatible
  with $\mathcal{C}$ and $\mathcal{D}$\,.
\item Famous Hironaka's
  desinguralisation theorem says that any compact subanalytic set is
  the image of a compact real analytic manifold with same dimension by
  a real analytic mapping. We will not use it specifically.
\end{itemize}
When $f:M\to \rz$ is a Lipschitz subanalytic function we consider the
two projections 
$p_{1}:M\times \rz\to M$ and $p_{2}:M\times \rz\to \rz$\,. 
From Hardt's result we know that there is a stratification of
$p_{2}:M\times \rz\to \rz$ which is compatible with 
$\mathcal{C}=\textrm{graph}\,(f)\sqcup (M\times\rz \setminus
\textrm{graph}\,(f))$ and $\mathcal{D}=\rz$\,. From this we deduce
that there is a stratification $\tilde{\mathcal{S}}$ of
$\textrm{graph}\,(f)$ and a finite number of points
$\left\{c_{1},\ldots c_{N_{f}}\right\}\in \rz$ such that all
$\tilde{S}\in \tilde{\mathcal{S}}$ satisfies
\begin{itemize}
\item either $p_{2}$ is constantly equal to some $c_{n}$ along
  $\tilde{S}$\,;
\item or there exists $n$ such that $p_{2}(\tilde{S})=]c_{n},c_{n+1}[$ and $\mathrm{rank}\,(p_{2}\big|_{\tilde{S}})=1$\,.
\end{itemize}
\begin{definition}
  \label{de:horstrata}For such a stratification of $\textrm{graph}\,(f)$\,, strata
corresponding to the first case will be called horizontal strata.
\end{definition}
Because $f$ is a Lipschitz function the projection $p_{1}:M\times
\rz\to M$ makes a diffeomorphism from $\tilde{S}$ to
$S=p_{1}(\tilde{S})$ which is a submanifold of $M$\,. The family
$\mathcal{S}=\left\{p_{1}(\tilde{S}),\tilde{S}\in \tilde{S}\right\}$
is now a stratification of $M$\,. When $\tilde{S}$ is a horizontal
stratum, then $f\big|_{S}$ is constant along
$S=p_{1}(\tilde{S})$\,. On the contrary when $\tilde{S}$ is not
horizontal $f\big|_{S}$ is a real analytic function with no critical
point on $S=p_{1}(\tilde{S})$\,.\\
Whitney's condition B also has a nice interpretation. 
It simply says
 in a local coordinate system (which allows the local 
identification of
 $T_{y}M$ with $\rz^{d}$ around any point $x\in M$)
$$
\left((x_{n})_{n\in\nz}\in (S')^{\nz},\;\lim_{n\to
    \infty}x_{n}=x\right)\Rightarrow\left(\forall T\in T_{x}S\sim \rz^{d}\,,\;\lim_{n\to\infty}(d(f\big|_{S'})_{x_{n}}[T]=d(f\big|_{S})_{x}[T]\right)\,.
$$ 
With the Riemannian structure it can be expressed in terms of
gradients. More exactly for any relatively compact open subset  $\omega_{S}$ of the stratum
$S$\,, and for $\varepsilon\in ]0,\varepsilon_{\omega_{S}}[$\,,
$\varepsilon_{\omega_{S}}>0$ small enough, the
exponential map $\exp(x,t)=\exp_{x}(t)\in M$ for $(x,t)\in TM$ is a
diffeomorphism from $\left\{(x,t)\in N\omega_{S}\,, |t|<
  \varepsilon\right\}$\,, where $N_{\omega_{S}}$ is the normal fiber
bundle over $\omega_{S}$\,, to its range
$\mathcal{T}_{\omega_{S},\varepsilon}\subset \left\{x\in M,
  d(x,\overline{\omega_{S}})<\varepsilon\right\}$\,, that we call a
tubular neighborhood of $\omega_{S}$\,. We refer the reader to
\cite{Lee} where tubular neighborhoods of closed submanifold are
introduced in this way and \cite{Lan} for further details and
generalizations with more general pseudo Riemannian structures. 
Another presentation using the embedding of $M$ in some $\rz^{N_{M}}$ is given in \cite{Hirs}.
Such a tubular neighborhood
$\mathcal{T}_{\omega_{S},\varepsilon}\subset M$ is an open subset of
the fiber bundle $\pi_{S}:N\omega_{S}\to \omega_{S}$ and is endowed
with the metric $g$ defined on $M$\,. Therefore the tangent bundle
$T_{x}\mathcal{T}_{\omega_{S},\varepsilon}=T_{x}M$ for $x\in
\mathcal{T}_{\omega_{S},\varepsilon}$\,, has an orthonomormal
decomposition $T_{x}M=T_{x}^{V}M\oplus^{\perp}T_{x}^{H}M$ where
$T_{x}^{V}M=\ker (d\pi_{S})\sim N_{\pi_{S}(x)}\omega_{S}$\,. For 
$x\in \mathcal{T}_{\omega_{S},\varepsilon}$ and $t\in
T_{x}M=T_{x}\mathcal{T}_{\omega_{S},\varepsilon}$ we define $\Pi_{S}t$
as the horizontal component of $t$ in this decomposition.
For $x\in \mathcal{T}_{\omega_{S},\varepsilon}$\,, the function
$f_{S}(x)=f(\pi_{S}x)$  is a real analytic function of $x\in
\mathcal{T}_{\omega_{S},\varepsilon}$\,. 
Because $f$ is a regular function along a stratum $S'\in \mathcal{S}$ its
gradient along $S'$ (with the metric induced by $g$) is denoted
$\nabla_{S'}f$\,.  With those notations the previous property can be
written
$$
\left((x_{n})_{n\in \nz}\in (S'\cap
  \mathcal{T}_{\omega_{S},\varepsilon})^{\nz}\,, \lim_{n\to\infty}x_{n}=x\in
  \omega_{S}\right)\Rightarrow 
\left(\lim_{n\to\infty}|\Pi_{S}\nabla_{S'}f(x_{n})-\nabla f_{S}(x_{n})|=0\right)\,.
$$ 
Let us summarize our notations:
\begin{itemize}
\item $\omega_{S}$ is a relatively compact open set of the stratum
  $S$\,.
\item $\mathcal{T}_{\omega_{S},\varepsilon}$ is a tubular neighborhood
  of $\omega_{S}$ diffeomorphic to $\left\{(x,t)\in N\omega_{S},
    |t|<\varepsilon\right\}$\,. It will be convenient to extend the
  notation to $\varepsilon=0$ with the large inequality and
  $\omega_{S}=S$\,, namely
  $\mathcal{T}_{S,0}=S$\,, which makes sense as
  $S=\limsup_{\varepsilon\to 0}\mathcal{T}_{\omega_{S,\varepsilon},
    \varepsilon}$ where $\omega_{S,\varepsilon}$ relatively compact in
  $S$ is well chosen when $\varepsilon>0$ is small\,.
\item When $S'$ is a stratum $\nabla_{S'}f$ is the gradient of $f$
  along $S'$ and for $x\in S'\cap
  \mathcal{T}_{\omega_{S},\varepsilon}$\,, $\Pi_{S}\nabla_{S'}f(x)$ is the horizontal component
  of $\nabla_{S'}f(x)$ in the orthogonal decomposition
  $T_{x}M=T_{x}^{V}\mathcal{T}_{\omega_{S},\varepsilon}\oplus^{\perp}T_{x}^{H}\mathcal{T}_{\omega_{S},\varepsilon}$\,.
\item Finally in $\mathcal{T}_{\omega_{S},\varepsilon}$\,, which is diffeomorphic
  to a subset of $N\omega_{S}$\,, one defines the
  regular function $f_{S}(x)=f(\pi_{S}x)$ where $\pi_{S}$ 
is the natural projection
  $\pi_{S}:N\omega_{S}\to \omega_{S}$\,.
\end{itemize}
\figinit{0.7cm}
\figpt 1:(-6,0)%
\figpt 2:(-4,0.4)%
\figpt 3:(-1,0.6)%
\figpt 4:(1.5,0.5)%
\figpt 5:(4,0.3)
\figpt 6:(7,0.15)%
\figpt 7:(7,0.1)%
\figpt 15:(3,4)

\figvectC 20:(-1,5)
\figvectC 21:(6,2)
\def\Vtrans{20}
\figptstra 101=1,2,3,4,5,6,7/0.5,\Vtrans/
\figptstra 201=1,2,3,4,5,6,7/-0.5,\Vtrans/
\figptstra 10=4/0.2,\Vtrans/
\figptstra 11=4/0.8,\Vtrans/
\figptstra 16=15/-0.362,\Vtrans/
\def\Vtrans{21}
\figptstra 12=10/0.8,\Vtrans/
\figptstra 13=10/-0.6,\Vtrans/

\figdrawbegin{}
\figset(width=1.3)
\figset arrowhead(ratio=0.2)
\figdrawarrow[10,15]
\figdrawarrow[10,16]
\figdrawcurve[1,1,2,3,4,5,6,6]
\figset(dash=3)
\figdrawcurve[101,101,102,103,104,105,106,106]
\figdrawcurve[201,201,202,203,204,205,206,206]
\figset(dash=5)
\figset(width=0.5)
\figdrawline[4,11]
\figdrawline[15,16]
\figdrawline[12,13]
\figdrawend

\figvisu{\figBoxA}{\textbf{Figure 13:} Picture of the projections $\pi_{S}$ and $\Pi_{S}$ when
  $x\in \mathcal{T}_{\omega_{S},\varepsilon}\cap S'$.}{
\figset write(mark=$\bullet$)
\figwritew 10:$x$(0.2)
\figset write(mark=$+$)
\figwrites 4:${\pi_{S}(x)}$(0.2)
\figset write(mark=)
\figwritene 15:${\nabla_{S'}f(x)}$(0.2)
\figwritese 16:${\Pi_{S}\nabla_{S'}f(x)}$(0.2)
\figwritene 12:${T_{x}^{H}\mathcal{T}_{\omega_{S},\varepsilon}}$(0.2)
\figwritenw 11:${T_{x}^{V}\mathcal{T}_{\omega_{S},\varepsilon}}$(0.2)
\figwrites 6:$\omega_{S}\subset S$(0.2)
\figwrites 2:${\mathcal{T}_{\omega_{S},\varepsilon}}$(1.4)
}
\centerline{\box\figBoxA}

\vspace{1cm}
With the compactness of $\overline{\omega}_{S}$ in $S$\,, Whitney 's
condition B actually implies the following uniform convergence result.
\begin{lem}
\label{le:unifcv}
Fix the  relatively compact open set $\omega_{S}$ of the stratum $S$
and let $\mathcal{T}_{\omega_{S},\varepsilon}$ denote the tubular
neighborhood defined for $\varepsilon>0$ small enough.
Then the quantities
$$
\max_{S'\in \mathcal{S}}\sup_{x\in
  \mathcal{T}_{\omega_{S},\varepsilon}\cap
  S'}|\Pi_{S}\nabla_{S'}f(x)-\nabla f_{S}(x)|
$$
and 
$$
\sup_{x\in \mathcal{T}_{\omega_{S},\varepsilon}}|\nabla
f_{S}(x)-\nabla f_{S}(\pi_{S}x)|\,,
$$
tend to $0$ as $\varepsilon\to 0^{+}$\,.
\end{lem}
\begin{proof}
  Ad absurdum if there is a sequence $(x_{n})_{n\in\nz}$ such that 
$|\Pi_{S}\nabla_{S'}f(x_{n})-\nabla f_{S}(x_{n})|\geq \eta>0$ while
$x_{n}\in \mathcal{T}_{\omega_{S},\frac{1}{n}}\cap S'$\,, then by the
compactness of $\overline{\omega_{S}}$ and the finiteness of
$\mathcal{S}$\,, we can assume that $S'$ is fixed and that
$\lim_{n\to\infty}x_{n}=x\in \overline{\omega_{S}}$\,. The lower bound
$|\Pi_{S}\nabla_{S'}f(x_{n})-\nabla f_{S}(x_{n})|\geq \eta>0$ while
$\lim_{n\to \infty}|\nabla f_{S}(x_{n})-\nabla f_{S}(x)|=0$\,, $\nabla
f_{S}(x)=\nabla_{S}f(x)$\,, then 
contradicts Whitney's condition~B.\\
Finally the last convergence is a consequence of the uniform
continuity of $\nabla f_{S}$ which can be defined on a compact neighborhood
of 
$\mathcal{T}_{\omega_{S},\varepsilon}$ for $\varepsilon\in
]0,\varepsilon_{S}[$\,, $\varepsilon_{S}>0$ small enough.
\end{proof}
\begin{prop}
\label{pr:verhypLipf}
When $f$ is a Lipschitz subanalytic function on $M$\,,
Hypothesis~\ref{hyp:Lipbar} is satified with $c_{1},\ldots,c_{N_{f}}\in
\rz$ being the values associated with horizontal strata in the
stratification of $\text{graph}(f)\subset M\times \rz$ described above.
\end{prop}
\begin{proof}
Let $x_{0}\in M\setminus
f^{-1}(\left\{c_{1},\ldots,c_{N_{f}}\right\})$\,. It belongs to a
stratum $S\in \mathcal{S}$ and we can find a relatively compact open
set $\omega_{S}\subset S$ such that $x_{0}\in \omega_{S}\subset S$\,.
The function $f_{S}$ is a real analytic-function  defined 
in the tubular open
$\mathcal{T}_{\omega_{S},\varepsilon}$ for
$\varepsilon\in]0,\varepsilon_{x_{0}}[$ with $\varepsilon_{x_{0}}>0$ small
enough. For $y\in \mathcal{T}_{\omega_{S},\varepsilon}\cap S'$\,, with
$S'\in \mathcal{S}$\,, $\dim S'=d$\,, we write
$$
\nabla f_{S}(y).\nabla f(y)=\left|\Pi_{S}\nabla f(y)\right|^{2}
-\left(\nabla f_{S}(y)-\Pi_{S}\nabla f(y)\right).\nabla f(y)
$$
and
$$
\left|\nabla f_{S}(y).\nabla f(y)-|\nabla f_{S}(x_{0})|^{2}\right|\leq 
\left||\Pi_{S}\nabla f(y)|^{2}-|\nabla f_{S}(x_{0})|^{2}\right|
+ M_{f}\left|\Pi_{S}\nabla f(y)-\nabla f_{S}(y)\right|\,.
$$ 
We know that $\left|\nabla f_{S}(x_{0})\right|>0$ because $S$ cannot
be an horizontal stratum.
By Lemma~\ref{le:unifcv}, $\varepsilon\in]0,\varepsilon_{x_{0}}[$ can be
chosen such that the right-hand side is smaller than
$\frac{1}{2}\left| \nabla f_{S}(x_{0})\right|^{2}$ for all $S'\in
\mathcal{S}$\,, such that $\dim S'=d$ and $x_{0}\in S\cap S'$\,. 
We have found a tubular neighborhood $U_{x_{0}}$ of $x_{0}$ and a coordinate system
$(x^{1},\ldots,x^{d})$ around $x_{0}$ by taking $x^{1}=f_{S}(x)$ such that
$$
\forall S'\in \mathcal{S}\,, \dim S'=d\,, \forall x\in U_{x_{0}}\cap
S'\,,\quad 
\partial_{x_{1}}f(x)\geq \frac{1}{C_{x_{0}}}\,.
$$
This neighborhood $U_{x_{0}}$ can then be reduced to
$$
U_{x_{0}}=\left\{x=(x^{1},x')=(x^{1},x^{2},\ldots,x^{d})\,,
  |x^{1}-x^{1}_{0}|<\delta\,,\, |x'|<\delta\right\}
$$
for some $\delta>0$\,.  The set
 $E=U_{x_{0}}\setminus (\cup_{\dim
  S'=d}S'\cap U_{x_{0}})$ has measure $0$ and $\nabla f(x)$ is well
defined for all $x\in U_{x_{0}}\setminus E$\,.
By Fubini's theorem the set of $x'$\,, $|x'|<\delta$\,, such that
$\left\{(x^{1},x')\,, |x^{1}-x^{1}_{0}|<\delta\right\}\cap E$ has a non zero
one dimensional measure, has Lebesgue's measure $0$ and we can write
for almost all $x'$\,, $|x'|<\delta$
$$
\forall x^{1},y^{1}\in ]x^{1}_{0}-\delta,x^{1}_{0}+\delta[\,,\quad
f(x^{1},x')-f(y^{1},x')=\int_{0}^{1}(x^{1}-y^{1})\partial_{x^{1}}f(x^{1}+t(y^{1}-x^{1}))~dt
$$
where the integrand is well defined for almost every $t\in [0,1]$ and
bounded from below by $\frac{1}{C_{x_{0}}}(x^{1}-y^{1})$ when
$x^{1}>y^{1}$\,. The continuity of $f$ then implies
$$
\forall (x^{1},x'), (y^{1},x')\in U_{x_{0}}\,,\quad
\frac{1}{C_{x_{0}}}|x^{1}-y^{1}|\leq |f(x^{1},x')-f(y^{1},x')|\,.
$$
\end{proof}
We will use open coverings of $f^{-1}([a,b])$ when
$[a,b]\sharp\left\{c_{1},\ldots, c_{N_{f}}\right\}=\emptyset$\,, made
of tubes $\mathcal{T}_{\omega_{S},\varepsilon_{S}}$ with 
$\varepsilon_{\omega_{S}}>0$\,. They will be constructed by induction on the
dimensions of the strata. They will be associated with a family of
parameters $(\varepsilon_{1},\ldots, \varepsilon_{d})$\,, with
$\varepsilon_{\omega_{S}}=\varepsilon_{\dim S}$\,. In the induction
process we authorize 
$\varepsilon_{\dim S}=0$ for $m<\dim S\leq d$\,, in which case
$\omega_{S}=S$ for every stratum $S$ of dimension $\dim S>m$\,.
\begin{definition}
\label{de:opcovm}
  Let $a<b$ belong to $\rz$ and set $\mathcal{S}_{[a,b]}=\left\{S\in
    \mathcal{S}\,, S\cap f^{-1}([a,b])\neq \emptyset\right\}$\,. A
tubular covering of $f^{-1}([a,b])$ 
contains two data, a family
  $(\varepsilon_{0},\varepsilon_{1},\ldots,\varepsilon_{d})\in [0,+\infty[^{d+1}$ and
  for every $S\in \mathcal{S}_{[a,b]}$\,, 
 a subset $\omega_{S}$ of $S$ which is either open 
 and relatively compact in $S$ if
 $\varepsilon_{\dim S}>0$ or equal to $S$ 
if $\varepsilon_{\dim S}=0$ such that for all $m\leq d$
\begin{eqnarray}
  \label{eq:tubcov}
&&f^{-1}([a,b])\cap \left(\ccup _{S\in \mathcal{S}_{[a,b]},\dim S\leq
  m}S\right)\subset
\ccup_{S\in \mathcal{S}_{[a,b]}, \dim S\leq
  m}\mathcal{T}_{\omega_{S},\varepsilon_{\dim S}}\,,\\
&&\mathcal{T}_{\omega_{S_{1}},\varepsilon_{m'}}\cap
\mathcal{T}_{\omega_{S_{2}},\varepsilon_{m'}}=\emptyset\quad
   \text{if}\quad S_{1}\neq S_{2}\,,\quad\dim S_{1}=\dim S_{2}=m'\leq m\,.
\end{eqnarray}
Such a tubular covering is said $\varepsilon$-adapted for
$\varepsilon\in ]0,1]$\,, if for any $S,S'\in \mathcal{S}_{[a,b]}$\,,
\begin{equation}
  \label{eq:tubepsadapt}
\sup_{x\in \mathcal{T}_{\omega_{S},\varepsilon_{\dim S}}\cap S'}
|\Pi_{S}\nabla_{S'}f(x)-\nabla f_{S}(x)|\leq \varepsilon\,,
\end{equation}
and 
\begin{equation}
\label{eq:tubesadapt2}
\sup_{x\in \mathcal{T}_{\omega_{S},\varepsilon_{\dim S}}}|\nabla
  f_{S}(x)-\nabla f_{S}(\pi_{S}x)|\leq \varepsilon\,.
\end{equation}
Such a covering is clearly an open covering when all the
$\varepsilon_{i}$'s are positive.\\
We will sumarize those situations by speaking of a (possibly ``an
$\varepsilon$-adapted'')(possibly ``open'') tubular covering
$(\mathcal{T}_{\omega_{S},,\varepsilon_{S}})_{S\in
  \mathcal{S}_{[a,b]}}$ associated with the parameters
$(\varepsilon_{0},\ldots,\varepsilon_{d})$\,.\\
\end{definition}
\vspace{-1cm}
\begin{figure}[h]
\centering{
\includegraphics[width=10cm]{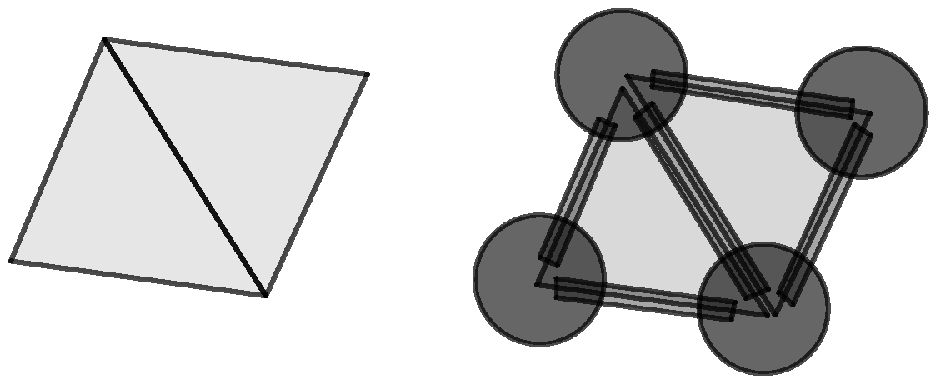}}
\captionsetup{labelformat=empty,width=14cm}
\caption{\textbf{Figure 14:} A schematic example of an open
  covering:
The stratification is on the
    left-hand side made of two triangles, the
    edges and the vertices; the open
    tubular covering with positive values for
    $\varepsilon_{0},\varepsilon_{1},\varepsilon_{2}$
is on the right-hand side. The outside of the two
    triangles is forgotten or one can compactify by identifying
    opposite external edges}
\end{figure}

A trivial example is given by
$(\varepsilon_{0},\varepsilon_{1},\ldots, \varepsilon_{d})=(0,\ldots,
0)$ and $\omega_{S}=S$ for all $S\in \mathcal{S}_{[a,b]}$\,. 
When $\mathcal{S}_{[a,b]}$ contains no stratum of dimension $m$\,, any
value $\varepsilon_{m}\geq 0$ can be used in the above definitions. 
When all the parameters $\varepsilon_{m}$\,, $0\leq m\leq d$\,, are
positive, this provides an open covering of $f^{-1}([a,b])$\,. Note
that when $\varepsilon_{m}=0$ and $\omega_{S}=S$ for $\dim S=m$\,, the condition
$x\in \mathcal{T}_{\omega_{S},\varepsilon_{m}}\cap S'$ actually
implies $x\in S=S'$ so that the condition \eqref{eq:tubepsadapt} is void for
strata $S$ of
dimension $m$\,. As a consequence, if
$(\mathcal{T}_{\omega_{S},\varepsilon_{\dim S}})_{S\in
  \mathcal{S}_{[a,b]}}$ is a (resp. an $\varepsilon$-adapted) tubular
covering of $f^{-1}([a,b])$ associated with the parameters
$(\varepsilon_{0},\ldots,\varepsilon_{d})$\,, then for any $m$
replacing $\varepsilon_{m'}$ by $0$ for $m'>m$\,, $\omega_{S}$ by $S$
if $\dim S>m$\,, $\varepsilon_{m}$ by $\varepsilon_{m}'\in
]0,\varepsilon_{m}]$ and leaving the other data, $\omega_{S}$ for
$\dim S\leq m$\,, $\varepsilon_{0},\ldots,\varepsilon_{m-1}$\,, unchanged give another
(resp. $\varepsilon$-adapted) tubular covering.\\
The following proposition implements the induction which leads to the
construction of families of $\varepsilon$-adapted open tubular 
coverings of
$f^{-1}([a,b])$\,, especially  when $[a,b]$ contains no ``critical value''.
\begin{prop}
\label{pr:tubopcov}
Assume first $f^{-1}([a,b])\cap
\left\{c_{1},\ldots,c_{N_{f}}\right\}=\emptyset$ where $c_{1},\ldots
c_{N_{f}}$ are values of $f$ associated with horizontal strata. Then
$\mathcal{S}_{[a,b]}$ contains no $0$-dimensional stratum and there
exists a (resp. an $\varepsilon$-adapted) tubular covering associated with $(\varepsilon_{0},0,\ldots,
0)$ for any $\varepsilon_{0}>0$\,.\\
Assume that there exists a (resp. an $\varepsilon$-adapted) tubular
covering associated with the parameters
$(\varepsilon_{0},\ldots,\varepsilon_{m-1},0,\ldots, 0)$ for $1\leq
m\leq d$ with $\varepsilon_{0}>0,\ldots, \varepsilon_{m-1}>0$\,, then
there exists $\varepsilon_{m}^{0}>0$  and for any $S\in
\mathcal{S}_{[a,b]}$\,, $\dim S=m$\,, a subset $\omega_{S}\subset S$ open and relatively
compact  $S$ such that for all $\varepsilon_{m}\in
]0,\varepsilon_{m}^{0}]$\,, the family
$(\mathcal{T}_{\omega_{S},\varepsilon_{\dim S}})_{S\in
  \mathcal{S}_{[a,b]}}$ associated with
$(\varepsilon_{0},\ldots,\varepsilon_{m},0,\ldots,0)$ and $\omega_{S}$
unchanged if $\dim S\leq m-1$\,, is another
(resp. $\varepsilon$-adapted) tubular covering of $f^{-1}([a,b])$\,.
\end{prop}
\begin{proof}
Because $\mathcal{S}_{[a,b]}$ contains no stratum of
  dimension $0$ a tubular covering is given by $\omega_{S}=S$ where
  all $S\in \mathcal{S}_{[a,b]}$ satisfy $\dim S\geq 1$ and any value
  of $\varepsilon_{0}>0$ makes sense.\\
Additionally every $S\in \mathcal{S}_{[a,b]}$ of
  dimension $1$ satisfies $\nabla_{S}f(x)\neq 0$ for every $x\in S\cap
  f^{-1}([a,b])$ and hence $f^{-1}([a,b])\cap S$ is a compact subset
  of $S$\,. We can choose $\omega_{S}$ open and relatively compact in
  $S$ such that $\omega_{S}\cap f^{-1}([a,b])$ is a neighborhood in
  $S$ of $S\cap f^{-1}([a,b])$\,. This is done for every $S\in \mathcal{S}_{[a,b]}$ such that
  $\dim S=1$\,. We can then choose $\varepsilon_{1}>0$ such that
  $\varepsilon_{1}<\frac{1}{2}
  d_{g}(\overline{\omega_{S_{1}}},\overline{\omega_{S_{2}}})$ for any
  $S_{1},S_{2}\in\mathcal{S}_{[a,b]}$\,, $\dim S_{1}=\dim S_{2}=1$\,, in
  order to ensure $\mathcal{T}_{\omega_{S_{1}},\varepsilon_{1}}\cap
  \mathcal{T}_{\omega_{S_{2}},\varepsilon_{1}}=\emptyset$ for
  $S_{1}\neq S_{2}$\,.\\
Assume now that the result holds for a given $m$\,,
 $1\leq m\leq d$\,. For $\dim S\leq m$\,, the set
  $\mathcal{T}_{\omega_{S},\varepsilon_{\dim S}}$  is an open set and
  $\ccup_{S\in
    \mathcal{S}_{[a,b]}, \dim S\leq
    m}\mathcal{T}_{\omega_{S},\varepsilon_{\dim S}}$ is an
  of $K_{[a,b],m}=f^{-1}([a,b])\cap (\ccup_{S\in \mathcal{S}_{[a,b]},\dim S\leq m}S)$\,. Consider the compact subset
$K_{[a,b],m+1}=f^{-1}([a,b])\cap 
(\ccup_{S\in \mathcal{S}_{[a,b]},\dim S=m+1}S)$\,. It is a a compact
set and so is $K_{[a,b],m+1}\setminus (\ccup_{S\in
  \mathcal{S}_{[a,b]},\dim S\leq
  m}\mathcal{T}_{\omega_{S},\varepsilon_{\dim S}})$ which by the
definition of the stratification $\mathcal{S}$ can be decomposed into
$\ccup_{S\in \mathcal{S}_{[a,b]},\dim S=m+1}K_{S}$ where $K_{S}$ is a
compact subset of $S$\,. We choose for $\omega_{S}$\,, $S\in
\mathcal{S}_{[a,b]}$\,, $\dim S=m+1$\,, a relatively compact
neighborhood of $K_{S}$ and then fix $\varepsilon_{m+1}>0$ small
enough such that $\mathcal{T}_{\omega_{S_{1}},\varepsilon_{m+1}}\cap
\mathcal{T}_{\omega_{S_{2}},\varepsilon_{m+1}}=\emptyset$ for any
$S_{1},S_{2}\in \mathcal{S}_{[a,b]}$\,, $\dim S_{1}=\dim S_{2}=m+1$
like in the case $m+1=1$\,.\\
Following this induction and by assuming that
$(\mathcal{T}_{\omega_{S},S})_{S\in \mathcal{S}_{[a,b]}}$ is an
$\varepsilon$-adapted tubular covering associated with
$(\varepsilon_{0},\ldots,\varepsilon_{m},0,\ldots,0)$\,,
$\varepsilon_{0}\ldots\varepsilon_{m}>0$\,,  $\varepsilon_{m+1}^{0}>0$
can be chosen such that
$$
\sup_{x\in \mathcal{T}_{\omega_{S},\varepsilon_{m+1}^{0}}\cap
  S'}|\Pi_{S}\nabla_{S'}f(x)-\nabla f_{S}(x)|\leq \varepsilon
$$
and
$$
\sup_{x\in \mathcal{T}_{\omega_{S},\varepsilon_{m+1}^{0}}}|\nabla
f_{S}(x)-\nabla f_{S}(\pi_{S}x)|\leq \varepsilon\,,
$$
for all $S\in \mathcal{S}_{[a,b]}$\,, $\dim S=m+1$\,, and all $S'\in
\mathcal{S}_{[a,b]}$\,. This still holds if $\varepsilon_{m+1}^{0}$ is
replaced by any $\varepsilon_{m+1}\in ]0,\varepsilon_{m+1}^{0}]$\,,
without changing the $\omega_{S}$\,, and
this ends the proof.
\end{proof}
\begin{definition}
  \label{de:FFFGGG}
Assume $f^{-1}([a,b])\cap
\left\{c_{1},\ldots,c_{N_{f}}\right\}=\emptyset$ and let
$(\mathcal{T}_{\omega_{S},\varepsilon_{\dim S}})_{S\in
  \mathcal{S}_{[a,b]}}$ be a tubular covering associated with the
parameters
$(\varepsilon_{0},\varepsilon_{1},\ldots,\varepsilon_{d})\in [0,+\infty[^{d+1}$\,.
The functions $\tilde{F}_{(\varepsilon_{0},\ldots,\varepsilon_{d})}$
and $F_{(\varepsilon_{0},\ldots,\varepsilon_{d})}$ are defined on
$f^{-1}([a,b])$ by
\begin{eqnarray}
\label{eq:deftF}
&& \tilde{F}_{(\varepsilon_{0},\ldots, \varepsilon_{d})}(x)=
\min_{x\in\mathcal{T}_{\omega_{S},\varepsilon_{\dim S}}\cap
   S'}|\Pi_{S}\nabla_{S'}f(x)|\,,\\
\label{eq:defF}&& F_{(\varepsilon_{0},\ldots,\varepsilon_{d})}(x)=
\min_{x\in\mathcal{T}_{\omega_{S},\varepsilon_{\dim S}}}|\nabla f_{S}(x)|\,,
\end{eqnarray}
where the minima are taken over $S,S'\in \mathcal{S}_{[a,b]}$\,.\\
On $f^{-1}([a,b])\times f^{-1}([a,b])$ the functions
$\tilde{G}_{(\varepsilon_{0},\ldots, \varepsilon_{d})}$ and
$G_{(\varepsilon_{0},\ldots,\varepsilon_{d})}$ are given by
\begin{eqnarray}
  &&
G_{(\varepsilon_{0},\ldots,\varepsilon_{d})}(x,y)=
\inf_{\tiny
     \begin{array}[c]{l}
       \gamma\in
     \mathcal{C}^{1}([0,1];f^{-1}([a,b]))\\
     \gamma(0)=x\,;\,
\gamma(1)=y
     \end{array}
  }\int_{0}^{1}F_{(\varepsilon_{0},\ldots,\varepsilon_{d})}(\gamma(t))|\gamma'(t)|~dt
\end{eqnarray}
with the same definition for $\tilde{G}_{(\varepsilon_{0},\ldots,\varepsilon_{d})}$\,.
\end{definition}
Before proving some results about those functions let us list some
simple properties:
\begin{itemize}
\item Because $\mathcal{S}_{[a,b]}$ is a finite collections of
  mesurable sets, the functions
  $\tilde{F}_{(\varepsilon_{0},\ldots,\varepsilon_{d})}$ and
  $F_{(\varepsilon_{0},\varepsilon_{d})}$ are measurable and the
  functions $\tilde{G}_{(\varepsilon_{0},\ldots,\varepsilon_{d})}$ and
  $G_{(\varepsilon_{0},\ldots,\varepsilon_{d})}$ are well defined.
\item When $\varepsilon_{1}=\ldots=\varepsilon_{d}=0$\,, the functions
  $\tilde{F}_{(0,\ldots,0)}$ and $F_{(0,\ldots,0)}$ are equal to
$$
\tilde{F}_{(0,\ldots,0)}(x)=F_{(0,\ldots,0)}(x)=\sum_{x\in S}1_{S}(x)|\nabla_{S}f(x)|\,,
$$
which is a lower semicontinuous function on $f^{-1}([a,b])$ due to
Whitney's condition~B and $|\Pi_{S}\nabla_{S'}f(x)|\leq
|\nabla_{S'}f(x)|$ for $x\in S'$ close enough to
 $S\subset \partial S'$\,.
\item Because $f$ is a Lipschitz function, the function
  $\tilde{F}_{(\varepsilon_{0},\ldots,\varepsilon_{d})}$ and
  $F_{(\varepsilon_{0},\ldots,\varepsilon_{d})}$ are uniformly bounded
  by $M_{f}=1+\|\nabla f\|_{L^{\infty}}$ when $\varepsilon\leq 1$
  because of
  $\left|\Pi_{S}\nabla_{S'}f(x)\right|\leq\left|\nabla_{S'}f(x)\right|\leq
  \left\|f\right\|_{W^{1,\infty}}$ and \eqref{eq:tubepsadapt}\,. 
Therefore the functions
  $\tilde{G}_{(\varepsilon_{0},\ldots,\varepsilon_{d})}$ and
  $G_{\varepsilon_{0},\ldots,\varepsilon_{d})}$ are $M_{f}$-Lipschitz 
pseudodistances on $f^{-1}([a,b])\times f^{-1}([a,b])$\,.
\item When $(\mathcal{T}_{\omega_{S},\varepsilon_{\dim S}})_{S\in
    \mathcal{S}_{[a,b]}}$ is an $\varepsilon$-adapted tubular covering
  of $f^{-1}([a,b])$\,, then 
$$
\sum_{x\in
  f^{-1}([a,b])}\left|\tilde{F}_{(\varepsilon_{0},\ldots,\varepsilon_{d})}(x)-F_{(\varepsilon_{0},\ldots,\varepsilon_{d})}(x)\right|\leq \varepsilon
$$
and hence
$$
\sup_{(x,y)\in
  f^{-1}([a,b])}\left|\tilde{G}_{(\varepsilon_{0},\ldots,\varepsilon_{d})}(x,y)-G_{(\varepsilon_{0},\ldots,\varepsilon_{d})}(x,y)\right|\leq
\varepsilon\times\text{diam}(f^{-1}([a,b]))\,,
$$
where $\text{diam}$ is the diameter for the metric $g$\,.
\item Let $(\mathcal{T}_{\omega_{S},\varepsilon_{\dim S}})_{S\in
    \mathcal{S}_{[a,b]}}$ be a tubular covering of $f^{-1}([a,b])$
  associated with the parameters $(\varepsilon_{0},\ldots,\varepsilon_{m},0,\ldots,0)$ with
  $\varepsilon_{0}\ldots\varepsilon_{m}>0$\,, $1\leq m\leq d$\,. For
  any $\varepsilon_{m}'\in ]0,\varepsilon_{m}]$\,, one gets another
  tubular covering of $f^{-1}([a,b])$ while keeping all the other data
  unchanged and for $\varepsilon_{m}'=0$ simply change $\omega_{S}$ into
  $S$ when $\dim S=m$\,. Then the functions
  $H_{(\varepsilon_{0},\ldots, \varepsilon_{m}',0,\ldots,0)}$\,, with
  $H=\tilde{F},F,\tilde{G},G$\,, are well defined  for any
  $\varepsilon_{m}'\in [0,\varepsilon_{m}]$ and they are decreasing
  with respect to $\varepsilon_{m}'$\,, i.e. increase as
  $\varepsilon_{m}'$ decays.
\end{itemize}
\begin{lem}
\label{le:FFFGGG}
In the framework of Definition~\ref{de:FFFGGG},
the function $F_{0,\ldots, 0}(x)=\tilde{F}_{(0,\ldots,0)}(x)$ is lower
semi-continuous bounded by $M_{f}=1+\|\nabla f\|_{L^{\infty}}$ and bounded
from below by a positive constant $m_{a,b,f}>0$\,. The function
$\tilde{G}_{(0,\ldots,0)}(x,y)=G_{(0,\ldots,0)}(x,y)$ is a
pseudodistance (fullfilling the symmetry and the triangular
inequality) which satisfies
$$
\forall x,y\in f^{-1}([a,b])\,,\quad
|f(x)-f(y)|\leq G_{(0,\ldots,0)}(x,y)\leq M_{f}d_{g}(x,y)\,,
$$
where $d_{g}$ is the geodesic distance between $x$ and $y$ in the
metric $g$\,.
\end{lem}
\begin{proof}
  We already noticed that
  $F_{(0,\ldots,0)}=\tilde{F}_{(0,\ldots,0)}$ is a lower
  semicontinuous function, bounded by $\|\nabla
  f\|_{L^{\infty}}$\,. Since $f^{-1}([a,b])$ contains no horizontal
  stratum 
$$
F_{(0,\ldots,0)}(x)=\sum_{S\in \mathcal{S}_{[a,b]}}1_{S}(x)|\nabla_{S} f(x)|
$$
does not vanish. The achieved minimum $m_{a,b,f}$ on the compact set
$f^{-1}([a,b])$ must be positive. With the estimate
$F_{(0,\ldots,0)}(x)\leq M_{f}$ for all $x\in f^{-1}([a,b])$\,, the fact
that $G_{(0,\ldots,0)}(x,y)$ defines a pseudodistance with the upper bound
$G_{(0,\ldots,0)}(x,y)\leq M_{f}d_{g}(x,y)$ is standard. 
For the lower bound because $M$-valued real analytic functions are
dense in $\mathcal{C}^{1}([0,1];M)$\,, the function $G_{(0,\ldots,0)}$ can be defined
as
$$
G_{(0,\ldots,0)}(x,y)=
\inf_{\tiny
  \begin{array}[c]{l}
\gamma\in \mathcal{C}^{\omega}([0,1];f^{-1}([a,b]))\\
\gamma(0)=x\,;\,\gamma(1)=y  
  \end{array}
}
\int_{0}^{1}F_{(0,\ldots,0)}(\gamma(t))|\gamma'(t)|~dt\,.
$$
Let $\gamma:[0,1]\to f^{-1}([a,b])\subset M$ be a real analytic
function such that
$$
G_{(0,\ldots,0)}(x,y)+\eta\geq
\int_{0}^{1}F_{(0,\ldots,0)}(\gamma(t))|\gamma'(t)|~dt\geq G_{(0,\ldots,0)}(x,y)\,.
$$
By using the recalled Hardt's result in \cite{Hardt} about the
stratification of real analytic mapping, now applied to $\gamma$ from
$[0,1]$ with the trivial stratification to $M$ with the stratification
$\mathcal{S}$\,, there exists a stratification of $[0,1]$\,, that is a
finite partition into open intervals and points $[0,1]=\sqcup_{I\in
  \mathcal{I}}I$ such that for any $I\in \mathcal{I}$ there exists
$S_{I}\in \mathcal{S}$ such that $\gamma(I)\subset S_{I}$\,. Hence
there exist $N\in\nz$\,, $0=t_{0}<t_{1}<\ldots<t_{N}=1$ and for any
$1\leq n\leq N$ a stratum  $S_{n}\in \mathcal{S}_{[a,b]}$ such that 
$\gamma(]t_{n-1},t_{n}[)\subset S_{n}$\,. We deduce
\begin{align*}
\int_{0}^{1}F_{(0,\ldots,0)}(\gamma(t))|\gamma'(t)|~dt
&=\sum_{n=1}^{N}\int_{t_{n-1}}^{t_{n}}|\nabla_{S_{n}}f(\gamma(t))||\gamma'(t)|~dt
\\
&\geq \sum_{n=1}^{N}|f(\gamma(t_{n}))-f(\gamma(t_{n-1}))|\geq |f(x)-f(y)|\,.
\end{align*}
We have proved for all $\eta>0$ the lower bound
$$
G_{(0,\ldots,0)}(x,y)+\eta\geq |f(x)-f(y)|\,,
$$
which ends the proof.
\end{proof}
\begin{prop}
\label{pr:minoGoptub}  Assume that $[a,b]\cap
\left\{c_{1},\ldots,c_{N_{f}}\right\}=\emptyset$\,. For any
$\varepsilon\in]0,1[$ there exist parameters
$(\varepsilon_{0},\ldots,\varepsilon_{d})\in ]0,+\infty[^{d+1}$
and an $\varepsilon$-adapted open tubular covering
$(\mathcal{T}_{\omega_{S},\varepsilon_{\dim S}})_{S\in
  \mathcal{S}_{[a,b]}}$ associated with the parameters
$(\varepsilon_{0},\ldots,\varepsilon_{d})$\,, such that the function
$G_{(\varepsilon_{0},\ldots,\varepsilon_{d})}$ defined in
Definition~\ref{de:FFFGGG} satisfies the uniform estimates:
\begin{equation}
  \label{eq:lowboundG}
\forall x,y\in f^{-1}([a,b])\,,\quad
|f(x)-f(y)|-\varepsilon\leq
G_{(\varepsilon_{0},\ldots,\varepsilon_{d})}(x,y)\leq 
M_{f}d_{g}(x,y)
\end{equation}
where $M_{f}=1+\|\nabla f\|_{L^{\infty}}$ and $d_{g}$ is the geodesic
distance on $(M,g)$\,.\\
For any $\varepsilon'\in]0,1[$\,, this tubular covering can be chosen, after
taking $\varepsilon>0$ small enough, such
that
\begin{equation}
  \label{eq:ineqfs}
\forall S\in \mathcal{S}_{[a,b]}\,, \nabla f(x). \nabla
f_{S}(x)-(1-\varepsilon')|\nabla f_{S}(x)|^{2}\quad \text{for~a.e.}~x\in
\mathcal{T}_{\omega_{S},\varepsilon_{\dim S}}\,,
\end{equation}
and
\begin{equation}
  \label{eq:lowbdfs}
\forall S,S'\in \mathcal{S}_{[a,b]}\,, \forall x\in
\mathcal{T}_{\omega_{S},\varepsilon_{\dim S}}\cap S'\,,\quad
|\nabla f_{S}(x)|\geq
\frac{m_{f,a,b}}{2}\quad,\quad|\Pi_{S}\nabla_{S'}f(x)|\geq \frac{m_{f,a,b}}{2}\,,
\end{equation}
where  $m_{f,a,b}>0$ was introduced in Lemma~\ref{le:FFFGGG}.
 \end{prop}
 \begin{proof}The diameter $\text{diam}(f^{-1}([a,b]))$ for the geodesic distance on
$(M,g)$ is denoted by
$$
\Delta_{a,b,f}=\text{diam}(f^{-1}([a,b]))\,.
$$
   The proof is made by induction on $m$\,, where $m$ is the maximal
   number such that $\varepsilon_{0}\ldots\varepsilon_{m}>0$\,, while
   playing with the two functions
   $\tilde{G}_{(\varepsilon_{0},\ldots,\varepsilon_{m},0,\ldots,0)}$
   and $G_{(\varepsilon_{0},\ldots,\varepsilon_{m},0,\ldots,0)}$\,.\\
More precisely we will prove that for $0\leq m\leq d$\,, there exists
$(\varepsilon_{0},\ldots,\varepsilon_{m})\in ]0,+\infty[^{m+1}$ and
an $\frac{\varepsilon}{d(2\Delta_{a,b,f}+1)}$-adapted tubular
covering $(\mathcal{T}_{\omega_{S},\varepsilon_{\dim S}})_{S\in
  \mathcal{S}([a,b])}$ associated with the parameters
$(\varepsilon_{0},\ldots,\varepsilon_{m},0,\ldots,0)$ such that
$$
|f(x)-f(y)|-\frac{m\varepsilon}{d}\leq G_{(\varepsilon_{0},\ldots,\varepsilon_{m},0,\ldots,0)}(x,y)\,.
$$
Notice that ``$\frac{\varepsilon}{d(2\Delta_{a,b,f}+1)}$-adapted''  is
stronger than ``$\varepsilon$-adapted''.\\
The statement is clearly true for $m=0$ because our assumption says that
$\mathcal{S}_{[a,b]}$ contains no $0$-dimensional stratum and 
$G_{(\varepsilon_{0},0,\ldots,0)}=\tilde{G}_{(\varepsilon_{0},0,\ldots,0)}$
  does not depend on $\varepsilon_{0}\in [0,+\infty[$\,, while the
  lower bound $G_{(0,\ldots,0)}(x,y)\geq|f(x)-f(y)|$ was proved in
  Lemma~\ref{le:FFFGGG}. Note additionally that the tubular covering
  $(\mathcal{T}_{\omega_{S},\varepsilon_{\dim S}})_{S\in
    \mathcal{S}_{[a,b]}}$\,, $\mathcal{T}_{\omega_{S},0}=S$ for $S\in
  \mathcal{S}_{[a,b]}$ is an $\frac{\varepsilon}{d(2\Delta_{a,b,f}+1)}$-adapted tubular
  covering of $f^{-1}([a,b])$\,.\\
Assume now that we have found
$(\varepsilon_{0},\ldots,\varepsilon_{m})\in ]0,+\infty[^{m+1}$ and an
$\frac{\varepsilon}{d(2\Delta_{a,b,f}+1)}$-adapted tubular covering
$(\mathcal{T}_{\omega_{S},\varepsilon_{\dim S}})_{S\in
  \mathcal{S}_{[a,b]}}$ such that 
$$
|f(x)-f(y)|-\frac{m\varepsilon}{d}\leq
G_{(\varepsilon_{0},\ldots,\varepsilon_{m},0,\ldots,0)}(x,y)\,.
$$
By Proposition~\ref{pr:tubopcov} 
$\varepsilon_{m+1}^{0}>0$ can be chosen such that for any
$\varepsilon_{m+1}\in]0,\varepsilon_{m+1}^{0}]$ there exists an $\frac{\varepsilon}{d(2\Delta_{a,b,f}+1)}$-adapted tubular
covering $(\mathcal{T}_{\omega_{S},\varepsilon_{\dim S}})_{S\in
  \mathcal{S}_{[a,b]}}$ associated with
$(\varepsilon_{1},\ldots,\varepsilon_{m+1},0,\ldots,0)$\,,
with $\omega_{S}$ independent of $\varepsilon_{m+1}>0$\,. 
For any $\varepsilon_{m+1}\in [0,\varepsilon_{m+1}^{(0)}]$ we deduce 
$$
\sup_{x,y\in^{f^{-1}([a,b])}}|\tilde{G}_{(\varepsilon_{0},\ldots,\varepsilon_{m+1},0,\ldots,0)}(x,y)-G_{(\varepsilon_{0},\ldots,\varepsilon_{m+1},0,\ldots,0)}(x,y)|
\leq \frac{\varepsilon}{d(2\Delta_{a,b,f}+1)}\times\Delta_{a,b,f}\,.
$$ 
Meanwhile we observed that
$\tilde{G}_{(\varepsilon_{0},\ldots,\varepsilon_{m},0,\ldots,0)}$ is
the monotonous limit as $\varepsilon_{m+1}\to 0^{+}$ of
$\tilde{G}_{(\varepsilon_{0},\ldots,\varepsilon_{m+1},0,\ldots,0)}$\,,
in the class of Lipschitz continuous functions on the compact set
$f^{-1}([a,b])\times f^{-1}([a,b])$\,. Dini's convergence theorem then
ensures that this convergence is uniform and we can choose
$\varepsilon_{m+1}\in ]0,\varepsilon_{m+1}^{0}]$ such that
$$
\sup_{x,y\in
  f^{-1}([a,b])}|\tilde{G}_{(\varepsilon_{0},\ldots,\varepsilon_{m+1},0,\ldots,0)}(x,y)-\tilde{G}_{(\varepsilon_{0},\ldots,\varepsilon_{m},0,\ldots,0)}(x,y)|\leq \frac{\varepsilon}{d(2\Delta_{a,b,f}+1)}\,.
$$
Gathering all those inequalities yields
\begin{align*}
  |f(x)-f(y)|-\frac{m\varepsilon}{d}&\leq
                                                         G_{(\varepsilon_{0},\ldots,\varepsilon_{m},0,\ldots,0)}(x,y)
\\
&\leq
\tilde{G}_{(\varepsilon_{0},\ldots,\varepsilon_{m},0,\ldots,0)}(x,y)+\frac{\varepsilon}{d(2\Delta_{a,b,f}+1)}\Delta_{a,b,f}\\
&
\leq
                                                         \tilde{G}_{(\varepsilon_{0},\ldots,\varepsilon_{m+1},0,\ldots,0)}(x,y)
+\frac{\varepsilon}{d(2\Delta_{a,b,f}+1)}(\Delta_{a,b,f}+1)\\
&\leq 
G_{(\varepsilon_{0},\ldots,\varepsilon_{m+1},0,\ldots,0)}(x,y)
+\frac{\varepsilon}{d(2\Delta_{a,b,f}+1)}(2\Delta_{a,b,f}+1)
\\
&\leq
G_{(\varepsilon_{0},\ldots,\varepsilon_{m+1},0,\ldots,0)}(x,y)
+\frac{\varepsilon}{d}\,.
\end{align*}
This ends the recurrence. The lower bound in \eqref{eq:lowboundG} is
finally proved when $m=d$
is reached.\\
For \eqref{eq:lowbdfs} it suffices to write
\begin{eqnarray*}
  &&\left|\nabla f_{S}(x)-\nabla f_{S}(\pi_{S}x)\right|\leq \varepsilon\quad,\quad
     |\nabla f_{S}(\pi_{S}x)|=G_{(0,\ldots,0)}(\pi_{S}x)\geq m_{f,a,b}\,,\\
&&|\Pi_{S}\nabla_{S'}f(x)-\nabla f_{S}(x)|\leq \varepsilon\,,
\end{eqnarray*}
and then to choose $\varepsilon\leq \frac{m_{f,a,b}}{4}$\,.\\
Finally with $S, S'\in \mathcal{S}_{[a,b]}$\,, $\dim S'=d$\,, and
$x\in \mathcal{T}_{\omega_{S},\varepsilon_{\dim S}}$\,, we have
$$
\nabla f(x).\Pi_{S}\nabla f(x)-(1-\frac{\varepsilon'}{2})\left|\Pi_{S}\nabla
  f(x)\right|^{2}=\frac{\varepsilon'}{2}\left|\Pi_{S}\nabla f(x)\right|^{2}\geq \frac{\varepsilon' m_{f,a,b}^{2}}{8}\,,
$$
while $\left\|\nabla f\right\|_{L^{\infty}}\leq M_{f}$ and 
$$
\left|\Pi_{S}\nabla f(x)-\nabla f_{S}(x)\right|\leq \varepsilon\,.
$$
By choosing $\varepsilon>0$ small enough we obtain for all $S,S'\in
\mathcal{S}_{[a,b]}$\,, $\dim S'=d$\,, and all $x\in S'$\,,
$$
\nabla f(x).\nabla f_{S}(x)-(1-\varepsilon')|\nabla f_{S}(x)|^{2}\geq 0\,.
$$
 \end{proof}
\subsubsection{Agmon type estimate for Lipschitz subanalytic
  potential}
\label{sec:moregenLipAgm}
 Proposition~\ref{pr:verhypLipf} says that Hypothesis~\ref{hyp:Lipbar}
 is satisfied when $f$ is a real analytic function on the compact
 Riemannian real analytic manifold $M$
 (Hypothesis~\ref{hyp:realana}), where the values
 $c_{1}<\ldots<c_{N_{f}}$ are the values associated with horizontal
 strata of $f$\,.\\
We now prove that Hypothesis~\ref{hyp:AgmonLip} is a consequence of
Hypothesis~\ref{hyp:realana} so that  Theorem~\ref{th:induc} and its
consequences in Section~\ref{sec:coroll} hold true under
Hypothesis~\ref{hyp:realana}.\\
Remember that Hypothesis~\ref{hyp:AgmonLip} gathers the results of
Proposition~\ref{pr:Agmon} and Proposition~\ref{pr:Agmon1} adapted to
a general Lipschitz function $f$\,. We will first prove the analogous
of Proposition~\ref{pr:Agmon1} in Proposition~\ref{pr:subLipAgm1} and
then deduce in Proposition~\ref{pr:subLipAgm} the analogous of Proposition~\ref{pr:Agmon}.

\begin{prop}
\label{pr:subLipAgm1}
Under Hypothesis~\ref{hyp:realana} and when 
$c_{1}<\ldots<c_{N_{f}}$ are the values associated with horizontal
strata according to Proposition~\ref{pr:verhypLipf}, choose $a<b$ such
that $[a,b] \cap \left\{c_{1},\ldots,c_{N_{f}}\right\}=\emptyset$\,.
If $\lim_{h\to
    0}\lambda_{h}=0$\,, the resolvent kernel
  $(\Delta_{f,f^{-1}([a,b]),h}-\lambda_{h})^{-1}(x,y)$ is well defined and
  satisfies 
$$
(\Delta_{f,f^{-1}([a,b]),h}-\lambda_{h})^{-1}(x,y)=\tilde{O}(e^{-\frac{|f(x)-f(y)|}{h}})\,,
$$
according to Definition~\ref{de:Otkernel}.
\end{prop}
\begin{proof}
This result relies on the stratification tools introduced in the
previous paragraph. It is proved in several steps, the first one
 being a
localization in suitable open subsets.
Let us  fix $x_{0}\in f^{-1}([a,b])$ with $f(x_{0})=t_{0}$ and we fix
the neighborhood of $x_{0}$ in $f^{-1}([a,b])$ as
$$
\mathcal{V}_{x_{0}}=f^{-1}([a,b])\cap f^{-1}(]t_{0}-\eta;t_{0}+\eta)
$$
where $\eta>0$ is a small parameter to be fixed at the end of the
analysis.\\
We want to prove that for any $\varepsilon>0$\,, any $h\in
]0,h_{\varepsilon}[$\,, $\Delta_{f,f^{-1}([a,b]),h}-\lambda_{h}$ is
invertible and that for any $r_{h}\in L^{2}(f_{a}^{b})$ such that
$\supp r_{h}\subset \mathcal{V}_{x_{0}}$\,,
$\omega_{h}=(\Delta_{f,f^{-1}([a,b]),h}-\lambda_{h})^{-1}r_{h}$ satisfies
$$
\|e^{\frac{|f(x)-f(x_{0})|}{h}}\omega_{h}\|_{W_{\partial}(f_{a}^{b})}=\tilde{O}(1)\|r_{h}\|\,.
$$
It will be convenient to call $a=t_{1}$ and $b=t_{2}$ especially when
 the arguments gather the three levels $t_{k}$\,, $k=0,1,2$\,.\\

\medskip
\noindent\textbf{i) Open covering of $f^{-1}([a,b])$:}
Because $[a,b]\cap\left\{c_{1},\ldots,c_{N_{f}}\right\}=\emptyset$\,,
for any $x\in f^{-1}([a,b])$ there exist a neighborhood
$U_{x}$ of $x$ in $M$ and a smooth function $\varphi_{x}$ on
$U_{x}$ and a constant $C_{x}>0$ such that
$$
\nabla f(y).\nabla \varphi_{x}(y)\geq \frac{1}{C_{x}}\quad
\text{and} \left|\nabla \varphi_{x}(y)\right|\leq C_{x}
\quad \text{for~a.e.~}y\in U_{x}\,. 
$$ 
Take for $\varphi_{x}(y)$ the coordinate function
$\varphi_{x}(y)=y^{1}$ given in Hypothesis~\ref{hyp:Lipbar} (see also Proposition~\ref{pr:verhypLipf}). 
By the compactness of $f^{-1}(\left\{t_{0},t_{1},t_{2}\right\})$\,,
there exists a finite family $(x_{i})_{i\in I}$ and constant
$\kappa>0$ small enough such that
$$
\nabla f.(\kappa \nabla \varphi_{x_{i}}(y))\geq 2|\kappa \nabla
\varphi_{x_{i}}(x)|^{2}\geq 2\kappa^{3}>0\quad \text{for~a.e.~}y\in U_{x_{i}}
$$
and for all $i\in I$\,.\\
Once this open covering $f^{-1}(\left\{t_{k},k=0,1,2\right\})\subset
\cup_{i\in I}U_{x_{i}}$ is fixed, we can choose the
parameter $\eta>0$ such that
$$
f^{-1}\left(\left\{t_{k},k=0,1,2\right\}+]-\frac{\eta}{2},\frac{\eta}{2}[\right)\subset
\cup_{i\in I}U_{x_{i}}\,.
$$
Again when $\eta>0$ is fixed and the stratification 
 $\mathcal{S}_{[a,b]}$ is introduced as in
 Subsection~\ref{sec:moregenLipStrat}\,, 
Proposition~\ref{pr:minoGoptub} provides us an open covering 
$$
(\mathcal{T}_{\omega_{S},\varepsilon_{\dim S}})_{S\in
  \mathcal{S}_{[a,b]}}
$$ such that the associated functions, $f_{S}$\,, $S\in
\mathcal{S}_{[a,b]}$\,, and $G_{(\varepsilon_{0},\ldots,\varepsilon_{d})}$ satisfy
\begin{eqnarray*}
&&\forall x,y\in f^{-1}([a,b])\,,\quad
|f(x)-f(y)|-\eta \leq
G_{(\varepsilon_{0},\ldots,\varepsilon_{d})}(x,y)\leq 
M_{f}d(x,y)\,,
\\
&&
\forall S\in \mathcal{S}_{[a,b]}\,, \nabla f(x). \nabla
f_{S}(x)-(1-\frac{\eta}{2})|\nabla f_{S}(x)|^{2}\leq 0 \quad \text{for~a.e.}~x\in
\mathcal{T}_{\omega_{S},\varepsilon_{\dim S}}\,,
\\
&&\forall S\in \mathcal{S}_{[a,b]}\,, \forall x\in
\mathcal{T}_{\omega_{S},\varepsilon_{\dim S}}\,,\quad
|\nabla f_{S}(x)|\geq
\frac{m_{f,a,b}}{2}\,.
\end{eqnarray*}
We now choose our open covering $f^{-1}([a,b])\subset \cup_{j\in
  J}\Omega_{j}$:
\begin{itemize}
\item $J=\mathcal{S}_{[a,b]}\cup I$\,; 
\item when $j=S\in \mathcal{S}_{[a,b]}$\,,
  $\Omega_{j}=\left\{x\in \mathcal{T}_{\omega_{S},\varepsilon_{\dim
        S}}\,, |f(x)-t_{k}|> \frac{\eta}{4}\,, k=0,1,2\right\}$
 and 
$\varphi_{j}=f_{S}$\,;
\item when $j=i\in I$\,, $\Omega_{j}=U_{x_{i}}\cap
  f^{-1}(\left\{t_{k},k=0,1,2\right\}+]-\frac{\eta}{2},\frac{\eta}{2}[)$\,,
and $\varphi_{j}=\kappa \varphi_{x_{i}}$\,.
\end{itemize}

\medskip
\noindent\textbf{ii) Choice of a global function $\varphi$:}
Once the open covering $f^{-1}([a,b])\subset \cup_{j\in J}\Omega_{j}$
is fixed  we choose
$$
\varphi(x)=(1-\eta)
\inf_{
\tiny \begin{array}[c]{l}
\gamma\in \mathcal{C}^{1}([0,1];f^{-1}([a,b])\\
  \gamma(0)=x_{0}\,,\, \gamma(1)=x
  \end{array}
}
\int_{0}^{1}1_{[a,b]\setminus \cup_{k=0}^{2}]t_{k}
    -\eta,t_{k}+\eta[}(f(\gamma(t))) F_{(\varepsilon_{0},\ldots,\varepsilon_{d})}(\gamma(t))|\gamma'(t)|~dt.
$$
Because the integrand is $0$ when $f(\gamma(t))\in
]t_{k}-\eta,t_{k}+\eta[$ the integral $\int_{0}^{1}[\ldots]dt$ can be
replaced by $\int_{T_{0}}^{T_{1}}[\ldots]dt$ where
$T_{0}=\max\left\{t\in[0,1],, f(\gamma(t))\in
  [t_{0}-\eta,t_{0}+\eta]\right\}$ and 
\begin{eqnarray*}
  && T_{1}=\min\left\{t\in [0,1]\,, f(\gamma(t))\geq
     f(b)-\eta\right\}\quad\text{if}~f(x)>f(b)-\eta\,, b=t_{2}\,,\\
&& T_{1}=\min\left\{t\in [0,1]\,, f(\gamma(t))\leq
     f(a)+\eta\right\}\quad\text{if}~f(x)<f(a)+\eta\,, a=t_{1}\,.
\end{eqnarray*}
The comparison with $G_{(\varepsilon_{0},\ldots,\varepsilon_{d})}(x,x_{0})$ then gives
$$
\frac{\varphi(x)}{1-\eta}\geq G_{(\varepsilon_{0},\ldots,\varepsilon_{d})}(x,x_{0})-2\eta \geq |f(x)-f(x_{0})|-3\eta
$$
and
\begin{equation}
\label{eq:minovphi}
\forall x\in f^{-1}([a,b])\,,\quad
\varphi(x)\geq |f(x)-f(x_{0})|-(b-a+3)\eta\,.
\end{equation}
The function $\varphi$ is a Lipschitz function of which the gradient
can be estimated almost surely in any $\Omega_{j}$\,, $j\in J$\,.
The triangle inequality for a pseudodistance implies for all $x,x'\in
f^{-1}([a,b])\cap \Omega_{j}$
\begin{align*}
\frac{|\varphi(x)-\varphi(x')|}{(1-\eta)}
&\leq 
\inf_{
\tiny \begin{array}[c]{l}
\gamma\in \mathcal{C}^{1}([0,1];f^{-1}([a,b])\\
  \gamma(0)=x\,,\,\gamma(1)=x'
  \end{array}
}
\int_{0}^{1}1_{[a,b]\setminus \cup_{k=0}^{2}]t_{k}
    -\eta,t_{k}+\eta[}(f(\gamma(t)))
  F_{(\varepsilon_{0},\ldots,\varepsilon_{d})}(\gamma(t))|\gamma'(t)|~dt
\\
&\leq
\inf_{
\tiny \begin{array}[c]{l}
\gamma\in \mathcal{C}^{1}([0,1];f^{-1}([a,b]\cap \Omega_{j})\\
  \gamma(0)=x\,,\,\gamma(1)=x'
  \end{array}
}
\int_{0}^{1}\left|\nabla \varphi_{j}(\gamma(t))\right|(\gamma(t))|\gamma'(t)|~dt\,.
\end{align*}
We used  that
$$
1_{[a,b]\setminus \cup_{k=0}^{2}]t_{k}
    -\eta,t_{k}+\eta[}(f(\gamma(t)))
  F_{(\varepsilon_{0},\ldots,\varepsilon_{d})}(\gamma(t))|\gamma'(t)|
$$
is
\begin{itemize}
\item $0$ and therefore bounded by $\left|\nabla
    \varphi_{j}(\gamma(t))\right|$ when $\gamma(t)\in
  \Omega_{j}\subset f^{-1}(\cup_{k=0}^{2}]t_{k}-\eta,t_{k}+\eta[)$
  when $j\in I$\,;
\item bounded by $\left|\nabla
    f_{S}(\gamma(t))\right||\gamma'(t)|$when $\gamma(t)\in \Omega_{j}$
  with $j=S\in \mathcal{S}_{[a,b]}$\,.
\end{itemize}
We deduce
\begin{equation}
  \label{eq:majoNvphi}
\forall j\in J\,, \quad \left|\nabla \varphi(x)\right|\leq
(1-\eta)|\nabla \varphi_{j}(x)|
\quad\text{for~a.e.}~x\in \Omega_{j}\,.
\end{equation}

\medskip
\noindent\textbf{iii) Partition of unity:}
Let $\sum_{j\in J}\chi_{j}^{2}\equiv 1$ in a neighborhood of
$f^{-1}([a,b])$ be a partition of unity with $\chi_{j}\in
\mathcal{C}^{\infty}_{0}(\Omega_{j};[0,1])$ where
$f^{-1}([a,b])\subset \cup_{j\in J}\Omega_{j}$ is the open covering
introduced in \textbf{i)}\,.  Accordingly the function $\varphi\in
W^{1,\infty}(f^{-1}[a,b])$ is the one introduced in \textbf{ii)}.
For any $\omega\in
W_{\partial}(f_{a}^{b};\Lambda T^{*}M)$\,, the relations \eqref{eq:reQ2} and
\eqref{eq:reQ1} of
Lemma~\ref{le:Agmon} give
\begin{align*}
\Real Q_{f,f^{-1}([a,b]),h}(\omega\,,\, e^{\frac{2\varphi}{h}}\omega)
&=\sum_{j\in J}\Real
Q_{f,f^{-1}([a,b]),h}(\chi_{j}\omega\,,\,
e^{\frac{2\varphi}{h}}\chi_{j}\omega)-h^{2}\||\nabla \chi_{j}|\tilde{\omega}\|^{2}\,.
\\
&=
\sum_{j\in J}
\left\|d_{f,f^{-1}([a,b]),h}\chi_{j}\tilde{\omega}\right\|^{2}
+
\left\|d_{f,f^{-1}([a,b]),h}^{*}\chi_{j}\tilde{\omega}\right\|^{2}
\\
&\hspace{1cm}-\langle
\chi_{j}\tilde{\omega}\,, |\nabla
  \varphi|^{2}\chi_{j}\tilde{\omega}\rangle
-h^{2}\||\nabla \chi_{j}|\tilde{\omega}\|^{2}\,.
\end{align*}
With \eqref{eq:majoNvphi} we deduce
\begin{align*}
\Real Q_{f,f^{-1}([a,b]),h}(\omega\,,\, e^{\frac{2\varphi}{h}}\omega)
&=
\sum_{j\in J}
\left\|d_{f,f^{-1}([a,b]),h}\chi_{j}\tilde{\omega}\right\|^{2}
+
\left\|d_{f,f^{-1}([a,b]),h}^{*}\chi_{j}\tilde{\omega}\right\|^{2}
\\
&\hspace{1cm}-(1-\eta)^{2}\langle
\chi_{j}\tilde{\omega}\,, |\nabla
  \varphi_{j}|^{2}\chi_{j}\tilde{\omega}\rangle
-h^{2}\||\nabla \chi_{j}|\tilde{\omega}\|^{2}\,.
\end{align*}
Now $\varphi_{j}$ can be extended to a $\mathcal{C}^{\infty}$ function
away from a neighborhood of $\supp \chi_{j}$ without changing the
expression and using \eqref{eq:reQ1} and \eqref{eq:reQ3} with
$\omega_{j}=e^{-(1-\eta)\frac{\varphi_{j}}{h}}\chi_{j}\tilde{\omega}\in
W_{\partial}(f_{a}^{b};\Lambda T^{*}M)$ and $\varphi$ replaced by $(1-\eta)\varphi_{j}$\,, we obtain
\begin{align*}
&\left\|d_{f,f^{-1}([a,b]),h}\chi_{j}\tilde{\omega}\right\|^{2}
+
\left\|d_{f,f^{-1}([a,b]),h}^{*}\chi_{j}\tilde{\omega}\right\|^{2}
-(1-\eta)^{2}\langle
\chi_{j}\tilde{\omega}\,, |\nabla
  \varphi_{j}|^{2}\chi_{j}\tilde{\omega}\rangle
\\
&\hspace{1cm}
=
Q_{f-(1-\eta)\varphi_{j},f^{-1}([a,b]),h}(\chi_{j}\tilde{\omega}\,,\,
  \chi_{j}\tilde{\omega})
\\
&\hspace{2cm}+
(1-\eta)\langle
(2\nabla f.\nabla \varphi_{j}-2(1-\eta)\left|\nabla \varphi_{j}\right|^{2}
+h\mathcal{L}_{\nabla\varphi_{j}}+h\mathcal{L}_{\nabla
  \varphi_{j}}^{*})\chi_{j}\tilde{\omega}\,,\, \chi_{j}\tilde{\omega}
\rangle
\\
&\hspace{2cm}
+
h(1-\eta)\left(\int_{f=b}-\int_{f=a}\right)
\langle\chi_{j}\tilde{\omega}\,,\,
  \chi_{j}\tilde{\omega}\rangle_{\Lambda
  T^{*}_{\sigma}M}
\left(\frac{\partial \varphi_{j}}{\partial n}\right)
(\sigma)~d\sigma\,.
\end{align*}
Because all the $\varphi_{j}$ are $\mathcal{C}^{\infty}$ functions
there exists $C>0$ such that
$$
\left|
\langle (\mathcal{L}_{\nabla \varphi_{j}}+\mathcal{L}_{\nabla
  \varphi_{j}}^{*})\chi_{j}\tilde{\omega}\,,\,
\chi_{j}\tilde{\omega})\rangle
\right|
\leq C\left\|\chi_{j}\tilde{\omega}\right\|^{2}\,.
$$
We have proved
\begin{align}
\label{eq:minopart1}
  \Real Q_{f,f^{-1}([a,b]),h}(\omega\,,\, e^{\frac{2\varphi}{h}}\omega)
&=
  \sum_{j\in J}Q_{f-(1-\eta)\varphi_{j}}(\chi_{j}\tilde{\omega}\,,\,
  \chi_{j}\tilde{\omega})
\\
\label{eq:minopart2}
&+
2(1-\eta)\langle
(\nabla f.\nabla \varphi_{j}-(1-\eta)\left|\nabla \varphi_{j}\right|^{2}
)\chi_{j}\tilde{\omega}\,,\, \chi_{j}\tilde{\omega}
\rangle
\\
\label{eq:minopart3}
&
+
h(1-\eta)\left(\int_{f=b}-\int_{f=a}\right)
\langle\chi_{j}\tilde{\omega}\,,\,
  \chi_{j}\tilde{\omega}\rangle_{\Lambda
  T^{*}_{\sigma}M}
\left(\frac{\partial \varphi_{j}}{\partial n}\right)
(\sigma)~d\sigma\,\\
\nonumber
&
+R_{h}(\tilde{\omega})
\end{align}
where the constant $C_{\eta}>0$ in
$$
\left|R_{h}(\tilde{\omega})\right|\leq C_{\eta}h\|\tilde{\omega}\|^{2}
$$
depends on $\eta>0$ via the construction of the open covering
$f^{-1}([a,b])\subset \cup_{j\in J}\Omega_{j}$\,, the functions
$\varphi_{j}$ and the partition of unity $\sum_{j\in
  J}\chi_{j}^{2}\equiv 1$\,.\\

\medskip
\noindent\textbf{iv) Local lower bounds:} We give a lower bound for
every individual $j\in J$ for the three terms
\eqref{eq:minopart1}\eqref{eq:minopart2} and \eqref{eq:minopart3}. The
first one \eqref{eq:minopart1} is obviously non negative according to
$$
Q_{f-(1-\eta)\varphi_{j}}(\chi_{j}\tilde{\omega}\,,\, \chi_{j}\tilde{\omega})
=\left\|d_{f-(1-\eta)\varphi_{j}}(\chi_{j}\tilde{\omega})\right\|^{2}
+
\left\|d_{f-(1-\eta)\varphi_{j}}^{*}(\chi_{j}\tilde{\omega})\right\|^{2}\geq 0
\,.
$$
For the other terms we distinguish according to $j\in I$ and
$j=S\in \mathcal{S}_{[a,b]}$\,.\\
\noindent\textbf{$\bullet~j\in I$:} In this case by recalling the choice
$\varphi_{j}=\kappa \varphi_{x_{j}}$\,, we know
$$
\nabla f.\nabla \varphi_{j}\geq 2\left|\nabla
  \varphi_{j}\right|^{2}\geq \kappa^{2}>0\quad\text{for~a.e.}~x\in \Omega_{j}\,.
$$
This implies 
$$
2(1-\eta)\left[\nabla f.\nabla \varphi_{j}-(1-\eta)\left|\nabla
    \varphi_{j}\right|^{2}\right]\geq 2(1-\eta)\left\|\nabla
  \varphi_{j}\right\|^{2}\geq (1-\eta)\kappa^{2}
\quad\text{for~a.e.}~x\in \Omega_{j}\,,
$$
where the positive constant $(1-\eta)\kappa^{2}$ is uniform w.r.t
$j\in I$\,.\\
Finally the condition $\nabla f.\nabla \varphi_{j}\geq 0$ makes sense
almost surely along the boundary $f^{-1}(\left\{a,b\right\})$ so
 that the integral terms  \eqref{eq:minopart3} are non negative.\\
\noindent\textbf{$\bullet~j=S\in \mathcal{S}_{[a,b]}$ :} Our choice of
$\Omega_{j}\subset \left\{x\in M\,, |f(x)-t_{k}|>\eta\,,
  k=0,1,2\right\}$ implies that the boundary terms
\eqref{eq:minopart3} vanish.
Finally  our choice $\varphi_{j}=f_{S}$ in \textbf{i)} implies
$$
\nabla f.\nabla\varphi_{j}-(1-\frac{\eta}{2})\left|\nabla
  \varphi_{j}\right|^{2}\geq 0\quad\text{for~a.e.}~x\in \Omega_{j}\,,
$$
We deduce
$$
2(1-\eta)\left[\nabla f.\nabla\varphi_{j}-(1-\eta)\left|\nabla
  \varphi_{j}\right|^{2}\right]
\leq 2(1-\eta)\frac{\eta}{2}\left|\nabla \varphi_{j}\right|^{2}
\geq (1-\eta)\frac{m_{f,a,b}^{2}}{4}
$$
almost every where in $\Omega_{j}$ with the positive constant
independent $(1-\eta)\frac{m_{f,a,b}^{2}}{4}$ independent of $j=S\in
\mathcal{S}_{[a,b]}$\,.\\
\noindent\textbf{v) Gathering all the lower bounds and conclusion:}\\
We
take
$\nu_{\eta}=(1-\eta)\min\left\{\kappa^{2},\frac{m_{f,a,b}^{2}}{4}\right\}$
and summing the previous lower bound w.r.t $j\in J$ leads to 
$$
\Real Q_{f,f^{-1}([a,b]),h}(\omega\,,\, e^{\frac{2\varphi}{h}}\omega)-\lambda_{h}\|e^{\frac{\varphi}{h}}\omega\|^{2}=
\geq (\nu_{\eta}-C_{\eta}h-\lambda_{h})\|\tilde{\omega}\|^{2}\geq
\frac{\nu_{\eta}}{2}\|\tilde{\omega}\|^{2}
$$
by taking $h\in ]0,h_{\eta}[$ for some small enough $h_{\eta}>0$\,.\\
Because $\Delta_{f,f^{-1}([a,b]),h}$ is self-adjoint the inequality
$$
\Real \langle e^{\frac{2\varphi}{h}}\omega\,,\,
(\Delta_{f,f^{-1}([a,b]),h)}-\lambda_{h})\omega \rangle
\geq \frac{\nu_{\eta}}{2}\|^{2}\tilde{\omega}\|\geq 
c_{\eta,h}\|\omega\|^{2}\,,\quad \tilde{\omega}=e^{\frac{\varphi}{h}}\omega\,,
$$
valid for all $\omega\in D(\Delta_{f,f^{-1}([a,b]),h})$ for some $c_{\eta,h}>0$\,,
implies that $\lambda_{h}$ belongs to the resolvent set of $\Delta_{f,f^{-1}([a,b]),h}$\,.\\
When $\omega_{h}$ solves
$(\Delta_{f,f^{-1}([a,b]),h}-\lambda_{h})\omega_{h}=r_{h}$\,, the same
inequality with $\varphi\equiv
0$ on $\supp r_{h}\subset f^{-1}(]t_{0}-\eta,t_{0}+\eta[)$\,, gives
$$
\|r_{h}\|\|\tilde{\omega}_{h}\|\geq \frac{\nu_{\eta}}{2}\|\tilde{\omega}_{h}\|^{2}\,,
$$
and $\|\tilde{\omega}_{h}\|\leq \frac{2}{\nu_{\eta}}\|r_{h}\|$\,. 
By using again \eqref{eq:reQ1} we deduce
\begin{align*}
\frac{2}{\nu_{\eta}}\|r_{h}\|^{2}\geq
\|r_{h}\|\|\tilde{\omega}_{h}\|&\geq \Real
Q_{f,f^{-1}([a,b]),h}(\omega_{h}\,,\, e^{\frac{2\varphi}{h}}\omega_{h})
-\lambda_{h}\|\tilde{\omega}_{h}\|^{2}
\\
&\geq
\|d_{f,h)}\tilde{\omega}_{h}\|^{2}+\|d_{f,h}^{*}\tilde{\omega}_{h}\|^{2}
-\left|\nabla \varphi\right|^{2}\tilde{\omega}_{h}\rangle
-\lambda_{h}\|\tilde{\omega}_{h}\|^{2}\,.
\end{align*}
And finally there exists a constant $M_{\eta}>0$ such that
$$
\frac{M_{\eta}}{h^{2}}\|r_{h}\|^{2}\geq
\|\tilde{\omega}_{h}\|^{2}+\|d\tilde{\omega}_{h}\|^{2}+\|d^{*}\tilde{\omega}_{h}\|^{2}
=\|e^{\frac{\varphi}{h}}\omega_{h}\|_{W_{\partial}(f_{a}^{b},\Lambda T^{*}M)}\,,
$$
with  $\varphi(x)\geq |f(x)-f(x_{0})|-(b-a+3)\eta$\,.\\
We conclude by taking $\eta>0$\,, on which all the construction
depends, arbitrarily small, the limit $h\to 0$ being taken for any
fixed $\eta>0$\,.
\end{proof}
\begin{remark}
\label{re:HJ}
In this proof, we did not use the global solution $\varphi$ to the
inequation $|\nabla
\varphi|^{2}-|\nabla f|^{2}\leq 0$ provided in \textbf{ii)} because such a
solution has no better regularity than the Lipschitz one. Instead we
really introduce the partition of unity in the process of obtaining
exponential decay estimates with all the functions $\varphi_{j}$ which
are regular enough and allow to use the various integration tricks of
Lemma~\ref{le:Agmon}, used  in particular in order 
to absorb the singularity of the
term $h(\mathcal{L}_{\nabla f}+\mathcal{L}_{\nabla f}^{*})$\,.
\end{remark}

\begin{prop}
\label{pr:subLipAgm}
Under Hypothesis~\ref{hyp:realana} and when 
$c_{1}<\ldots<c_{N_{f}}$ are the values associated with horizontal
strata according to Proposition~\ref{pr:verhypLipf}, choose $a<b$\,, $a,b\not\in
\left\{c_{1},\ldots,c_{N_{f}}\right\}$ and call $U$ the
compact set $f^{-1}(\left\{c_{1},\ldots,c_{N_{f}}\right\}\cap[a,b])$\,. 
All families $(\lambda_{h})_{h>0}\in \cz$\,, $(r_{h})_{h>0}\in
L^{2}(f_{a}^{b})$ and $\omega_{h}\in
D(\Delta_{f,f^{-1}([a,b]),h})\subset W_{\partial}(f_{a}^{b};\Lambda T^{*}M)$ such
that
$$
(\Delta_{f,f^{-1}([a,b]),h}-\lambda_{h})\omega_{h}=r_{h}\quad,\quad
\supp r_{h}\subset K\quad,\quad \lim_{h\to 0}\lambda_{h}=0\,,
$$
where $K$ is a fixed compact subset of $f^{-1}([a,b])$\,,
satisfy  the estimate  
$$
\|e^{\frac{\min_{y\in U\cup K}|f(.)-f(y)|}{h}}\omega_{h}\|_{W_{\partial}(f_{a}^{b})}=\tilde{O}(1)\left[\|r_{h}\|_{L^{2}(f_{a}^{b})}+t_{U}\|\omega_{h}\|_{L^{2}(f_{a}^{b})}\right]\,,
$$
where $t_{U}=1$ if $U\neq \emptyset$ and $t_{U}=0$ if $U=\emptyset$\,.
\end{prop}
\begin{proof}
  The case when $U=\emptyset$ is contained in
  Proposition~\ref{pr:subLipAgm1}. Let us consider the case when
  $U\neq \emptyset$\,. First of all, the positivity of
  $\Delta_{f,f^{-1}([a,b]),h}$ implies 
  \begin{multline*}
\left\|d_{f,h}\omega_{h}\right\|^{2}+\|d_{f,h}^{*}\omega_{h}\|^{2}+
(C-\Real \lambda_{h})\|\omega_{h}\|^{2}=\Real \langle \omega_{h}\,, (\Delta_{f,f^{-1}([a,b]),h}+C-\lambda_{h})\omega_{h}\rangle
\\
\leq \left\|r_{h}\right\|\left\|\omega_{h}\right\|+C\|\omega_{h}\|^{2}\,.
\end{multline*}
By taking $C> 2(1+\|f\|_{W^{1,\infty}})$ we obtain
$$
\|\omega_{h}\|_{W_{\partial}(f_{a}^{b})}=\tilde{O}(1)(\|r_{h}\|_{L^{2}(f_{a}^{b})}+\|\omega_{h}\|_{L^{2}(f_{a}^{b})})
$$
which provides $W^{1,2}$ estimates of $\omega_{h}$ in any compact
subset of $f_{a}^{b}=f^{-1}(]a,b[)$\,.\\
For $\varepsilon>0$ small enough, consider a cut-off function
$\chi_{\varepsilon} \in
\mathcal{C}^{\infty}(M;[0,1])$ equal to $1$ in
$K_{\varepsilon}=f^{-1}((\cup_{k=1}^{N_{f}}[c_{k}-\varepsilon,c_{k}+\varepsilon])\cap
[a,b])$ and to $0$ in the complement of $K_{2\varepsilon}$\,. The form
$\chi_{\varepsilon}\omega_{h}$ solves
$$
(\Delta_{f,h}-\lambda_{h})((1-\chi_{\varepsilon})\omega_{h})=(1-\chi_{\varepsilon})r_{h}+P_{\chi_{\varepsilon}}\omega_{h}\,,
$$ 
where $P_{\chi_{\varepsilon}}$ is a first order differential operator
with coefficients supported in $K_{2\varepsilon}\setminus
K_{\varepsilon}$ and $\chi_{\varepsilon}\omega_{h}\in
\mathop{\oplus}_{k=1}^{N_{f}-1}\Delta_{f,f^{-1}([\max(c_{k}+\varepsilon,a),
  \min(c_{k+1}-\varepsilon,b)]),h}$\,. The resolvent estimate of
Proposition~\ref{pr:subLipAgm1} applied to every $\Delta_{f,f^{-1}(([\max(c_{k}+\varepsilon,a),
  \min(c_{k+1}-\varepsilon,b)]),h}$ then implies
$$
\|e^{\frac{\min_{y\in U\cup
      K}|f(.)-f(y)|}{h}}\omega_{h}\|_{W_{\partial}(f_{a}^{b})}\leq \tilde{O}(e^{\frac{10\varepsilon}{h}})\left[\|r_{h}\|_{L^{2}(f_{a}^{b})}+\|\omega_{h}\|\right]_{L^{2}(f_{a}^{b})}\,,
$$
and then we choose $\varepsilon>0$ arbitrarily small before taking the
limit $h\to 0$\,.
\end{proof}

\section{Applications}

\label{sec:applications}

The spectral version of the stability theorem,
Corollary~\ref{cor:mainsimple} in the Introduction or
Theorem~\ref{th:stab1} for a more general version, corresponds to
what can be expected at the level of Arrhenius law identifying the
exponential scales. It is a
straightforward consequence of Theorem~\ref{th:induc}.
But the construction of global quasimodes for Theorem~\ref{th:induc}
is actually much more informative. It allows to compute the
subexponential factor, a la Eyring-Kramers, in many situations which
lead to different kind of asymptotic behaviours. As it was discussed
in the Introduction, no continuity with respect to $f$ can be expected
in the asymptotic leading term. Nevertheless some robust integral
formulation allow to follow the effect of deformations of $f$ on the
spectral quantities and to explain the emerging discontinuities. 
Contrary to Theorem~\ref{th:induc} and its consequences in
Section~\ref{sec:coroll}, we do not have a satisfactory general
formulation of this kind of refined stability property 
and we prefer to make explicit various examples, corresponding to
interesting practical cases.

\subsection{The generic Morse case}

In this subsection, we recall the results of \cite{LNV}.
Although they were presented in the oriented case, those results hold
in the more general case of non necessarily oriented compact Riemannian manifolds.
The proofs are simply adapted by paying attention to the duality arguments,
the Hodge~$\star$ operator sending  the sections of
$\Lambda^{p}T^{*}M$ to sections of $\Lambda^{p}T^{*}M\otimes \mathrm{or}_{M}$\,.
The important assumption which was made in \cite{LNV} concerns the simplicity of the critical values
of the Morse function $f$: the latter function
has  distinct critical values, which allows in particular to identify critical points
with critical values. 
In \cite{LNV}, the set $\mathcal U$ of critical points was partitioned into
 lower  $\mathcal U_{\text{L}}=\cup_{p\in\{0,\dots,d\}}\mathcal U_{\text{L}}^{(p)}$\,, upper
 $\mathcal U_{\text{U}}=\cup_{p\in\{0,\dots,d\}}\mathcal U_{\text{U}}^{(p)}$\,, and homological
$\mathcal U_{\text{H}}=\cup_{p\in\{0,\dots,d\}}\mathcal U_{\text{H}}^{(p)}$ 
  critical points.
  This partition actually coincides with the partition of bar endpoints $\mathcal J=\mathcal X\cup \mathcal Y
  \cup \mathcal Z$ in this order.
  In \cite{LNV}, we defined a boundary map $\partial_{\cal B}:\mathcal U_{\text{U}}^{(p+1)}\to \mathcal U_{\text{L}}^{(p)}$
  and $\mathcal U_{\text{U}}\cup \mathcal U_{\text{H}}\subset \ker \partial_{\cal B}$\,.
  It is exactly the dual version of the map $\mathbf{d}_{{\cal B}}$
  of Appendix~\ref{sec:triv}
  defined by
   $\mathbf{d}_{{\cal B}}: \mathcal X^{(p)}\to \mathcal Y^{(p+1)}$
   and ${\cal Y} \cup {\cal Z} \subset \ker \mathbf{d}_{{\cal B}}$\,.
   Actually, in \cite{LNV}, we started with the homological point of view before we realized
   that working directly in terms of cohomology was more natural for this analysis.
   The link with relative cohomology groups of sublevel sets of $f$\,,
   which is detailed at the end of Appendix~\ref{sec:sheaf},
   can be handled with elementary arguments under the assumptions of \cite{LNV}
   (Morse function with distinct critical values). Note that this generic Morse situation is often
   used as a simple way to introduce persistent homology (see
   e.g. \cite{EdHa}). Although it is an obvious bijection under the
   assumption that the Morse function $f$ has simple critical values,
   we use the notations, when it is  necessary, 
$\underline{x}_{\alpha}$\,, $\underline{y}_{\alpha}$ or
$\underline{z}_{\alpha}$ for the critical points associated with
values $x_{\alpha}=f(\underline{x}_{\alpha})\in \mathcal{X}$\,,
$y_{\alpha}=f(\underline{y}_{\alpha})\in \mathcal{Y}$ and
$z_{\alpha}=f(\underline{z}_{\alpha})\in \mathcal{Z}$\,. As a
comparison with the notations of Subsection~\ref{sec:bar code}, 
it is not necessary nor useful to
distinguish $x_{\alpha}=(a_{\alpha},\alpha)\in \rz\times A$ 
from the value $a_{\alpha}=f(\underline{x}_{\alpha})$\,.\\
Finally note  that the result of \cite{LNV} can be recovered while
combining Theorem~\ref{th:induc} of the present text with the final
computations of \cite{LNV}-Section~4 which rely on local WKB
approximations valid locally for any Morse function $f$\,.

Here is the main result of \cite{LNV}
with the above modified notations.

\begin{thm}
\label{th.LNV}
Assume that $f$ is a Morse function with simple critical values.
For any $p\in\{0,\dots,d\}$\,, there exists $c>0$ such that for every $h>0$ small enough,  the spectrum of $\Delta^{(p)}_{f,M,h}$
 satisfies
 $$\sigma(\Delta^{(p)}_{f,M,h})\cap [0,ch]=\sigma(\Delta_{f,h}^{(p)})\cap [0,e^{-\frac ch}]\,,$$
 and the latter set consists in $\card(\mathcal J^{(p)})$ eigenvalues counted with multiplicity. 
  For every  $h>0$ small enough,
 there exists moreover a bijection $j:\mathcal J^{(p)} \to \sigma(\Delta^{(p)}_{f,M,h})\cap [0,ch]$\,,
 where the latter set is counted with multiplicity, 
such that:
\begin{enumerate} 
\item 
For every $z_{\alpha}$ in $\mathcal Z^{(p)}$\,,
the associated eigenvalue is
$$j(z_{\alpha})=0\,.$$
\item For  every $x_{\alpha}$ in $ \mathcal X^{(p)}$\,,
$x_{\alpha}$ being the lower endpoint of the bar $[x_{\alpha},y_{\alpha}[$\,,
and hence $y_{\alpha}=\mathbf{d}_{{\cal B}} x_{\alpha}$\,,
there exists a homological  constant $\kappa_{\alpha}\in \mathbb Q^{*}$
such that
$$
j(x_{\alpha})\ =\ \kappa^{2}_{\alpha}\,\frac{h}{\pi}\,
\frac{|\lambda_{1}(\underline{y}_{\alpha})\cdots\lambda_{p+1}(\underline{y}_{\alpha})|}{
|\lambda_{1}(\underline{x}_{\alpha})\cdots\lambda_{p}(\underline{x}_{\alpha})|}
\frac{{|{\rm det\,} {\rm Hess\,} f(\underline{x}_{\alpha})|}^\frac12}
{{|{\rm det\,}{\rm Hess\,} f(\underline{y}_{\alpha})|}^\frac12}\,
e^{-2\frac{y_{\alpha}-x_{\alpha}}{h}}\big(1+\mathcal{O}(h)\big)\,,
$$
where, for any critical point $\underline{s}$
of $f$ with index $\ell$ and critical value $s=f(\underline{s})$\,, $\lambda_{1}(\underline{s}),\dots,\lambda_{\ell}(\underline{s})$
denote the negative eigenvalues of ${\rm Hess\,} f(\underline{s})$\,.
\item And $y_{\alpha}$ in $ \mathcal Y^{(p)}$\,,
$y_{\alpha}$ being the upper endpoint of the bar $[x_{\alpha},y_{\alpha}[$\,,
and hence $y_{\alpha}=\mathbf{d}_{{\cal B}} x_{\alpha}$\,,
there exists a homological  constant $\kappa_{\alpha}\in \mathbb Q^{*}$
such that

$$
j(y_{\alpha})
\ =\ \kappa^{2}_{\alpha}\,\frac{h}{\pi}\,
\frac{|\lambda_{1}(\underline{y}_{\alpha})\cdots\lambda_{p}(\underline{y}_{\alpha})|}{
|\lambda_{1}(\underline{x}_{\alpha})\cdots\lambda_{p-1}(\underline{x}_{\alpha})|}
\frac{{|{\rm det\,} {\rm Hess\,} f(\underline{x}_{\alpha})|}^\frac12}
{{|{\rm det\,} {\rm Hess\,} f(\underline{y}_{\alpha})|}^\frac12}\,
e^{-2\frac{y_{\alpha}-x_{\alpha}}{h}}\big(1+\mathcal{O}(h)\big)\,,
$$
where, for any critical point $\underline{s}$
of $f$ with index $\ell$ and critical value $s=f(\underline{s})$\,, $\lambda_{1}(\underline{s}),\dots,\lambda_{\ell}(\underline{s})$
denote the negative eigenvalues of ${\rm Hess\,} f(\underline{s})$\,.
\end{enumerate} \end{thm}

\begin{remark}
\begin{enumerate}
\item
Theorem~\ref{th.LNV} is a refinement of
Theorem~\ref{th:mainsimple} in this generic  Morse
situation. It extends Eyring-Kramers asymptotic formulas known in the
case $p=0$\,. The boundary version in $f^{-1}([a,b])$
corresponding to Theorem~\ref{th:specres} is also found in  \cite[Theorem~4.5]{LNV}.
In both papers, the general strategy consists in a recurrence with respect to the number of critical values, carried out 
by increasing the interval $[a,b]$\,.
The setting in \cite{LNV},  was simpler because: a) the critical values were assumed to be simple
while here they may be multiple or very degenerate; b) the
subexponential factors
of exponentially small quantities had explicit leading terms derived
from the WKB
approximations
(this is not possible here).
\footnote{A small confusion occurred in the construction of accurate global quasimodes 
in \cite[Section~4.2]{LNV}: a version of Proposition~\ref{pr:finalorth} is missing and can  be easily corrected.}
In this section, we will combine Theorem~\ref{th:induc} with the local computations of
 \cite{LNV}-Section~4
to provide a more general approach. 

\item In \cite{LNV}, hanks to the Morse assumption, we could compute the subexponential factors using WKB and Laplace methods..
 On the other hand,
the  exponential  factors are given by global topological quantities:
the lengths of the bar code. In the present paper we manage to compute  the logarithmic equivalents of the small eigenvalues without any knowledge of the exponential factor.  
\item The connexion between the local computation around the lower endpoint $x_{\alpha}$
and the upper one $y_{\alpha}=\mathbf{d}_{{\cal B}} x_{\alpha}$
is implemented by an application of Stokes's formula.
The boundary operator $\partial $ for chains induces a linear application  
from $H_{p+1}(f^{y_{\alpha}+\varepsilon},f^{{y_{\alpha}}-\varepsilon})$
into $H_{p}(f^{x_{\alpha}+\varepsilon},f^{{x_{\alpha}}-\varepsilon})$\,.
Under the generic Morse assumption, these spaces are $1$-dimensional
with natural bases given by the stable manifolds of $\nabla f$\,.
This actually provides the coefficients $\kappa_{\alpha}$ (see \cite[Proposition~2.12]{LNV}).
When the critical values correspond to multiple critical points, such a construction
has to be replaced by more general linear algebra (see Subsection~\ref{sec:appmoregenMorse}).
\end{enumerate}
\end{remark}

As shown in
\cite{HKN},  the homological constants $\kappa^{2}_{\alpha}$
equal $1$
when $p=0$\,, and also when $p=d$ and $M$ is oriented
by duality. 
In the case of oriented
surfaces treated in \cite{Lep-2D}, a combination of these results together
with simple duality and chain complex
arguments then implies that
these constants equal $1$ for any $p\in\{0,1,2\}$\,.
Nevertheless, contrary to this indication that it could be true in general,
which was moreover our intuition when we wrote \cite{LNV}, 
this is not the case 
as soon as $d\geq 3$ and even when $d= 2$
in the non-oriented case.
The simplest example comes from Morse theory on the projective plane.
It is more generally related to the ``open book picture'' exhibited on the front
cover of \cite{Lau-Book}.
 
 To be more specific, we shall prove the following result.
 
\begin{prop} Let $X$ be a $d$-dimensional manifold. \begin{enumerate} 
\item If $d=1,2$\,, and $X$ is orientable, then $\kappa_\alpha^2=1$\,,  
\item The coefficient $\kappa_{\alpha}^2$ may be equal to $4$ for some well chosen Morse functions on $ {\mathbb R} P^2$ and on $ {\mathbb R} P^3$\,.
\item For $d\geq 3$ and each integer $n$\,, there exists a manifold $X_n$ of dimension $d$ such that  $\kappa_\alpha^2=n^2$\,.
\item For $d\geq 4$\,,  for any integer $n$ and any closed manifold $X$\,,
  there is a function $f_n$ on $X$ such that  $\kappa_\alpha^2$ takes
  the value $n^2$\,. 
\end{enumerate} \end{prop} 
\begin{proof} 
The number $\kappa_{\alpha}$ is obtained as follows:  consider the sphere
$S^-(\underline{y}_{\alpha})$ in the unit disc bundle of the
descending manifold from $\underline{y}_{\alpha}$\,, the stable manifold
of $\nabla f$\,. It is homologous to a multiple, $\kappa_{\alpha}$ of
the descending manifold from $\underline{x}_{\alpha}$\,,
$W^-(\underline{x}_{\alpha})$\,, with $\mathbf{d}_{{\cal B}}
x_{\alpha}=y_{\alpha}$\,. But since the ascending manifold from
$\underline{x}_{\alpha}$\,, the unstable manifold of $\nabla f$\,, 
$W^+(\underline{x}_{\alpha})$ has intersection $+1$ with
$W^-(\underline{x}_{\alpha})$\,, the number
 $\kappa_{\alpha}$ is the intersection number of $S^-(\underline{y}_{\alpha})$ and $W^+(\underline{x}_{\alpha})$\,.
We work here  under the generic Morse-Smale assumption saying that all the stable and unstable manifolds
are mutually transverse, which 
ensures the finiteness of $\kappa_{\alpha}$\,,
within the construction of the Thom-Smale complex (see \cite{Lau-Book}). 
 In homological terms, if we set
 $x_{\alpha}=f(\underline{x}_{\alpha}),
 y_{\alpha}=f(\underline{y}_{\alpha})$\,, and $ \varepsilon >0$ small
 enough,
we have
the maps

$$\xymatrix{& H_*(f^{y_\alpha+\varepsilon }, f^{y_\alpha- \varepsilon } )\ar[d]^{\partial}\\H_{*-1}(f^{x_\alpha+ \varepsilon }, f^{x_\alpha- \varepsilon })\ar[r]& H_{*-1}(f^{y_\alpha-\varepsilon }, f^{x_\alpha- \varepsilon } )\ar[d]\\
&H_*(f^{y_\alpha+\varepsilon }, f^{x_\alpha- \varepsilon } )}
$$
Now since $H_*(f^{y_\alpha+\varepsilon }, f^{y_\alpha- \varepsilon }; {\mathbb Z}  )$ and $H_*(f^{x_\alpha+\varepsilon }, f^{x_\alpha- \varepsilon }; {\mathbb Z} )$  are isomorphic to $ {\mathbb Z} $\,,
the $ {\mathbb R}$-vector spaces  $H_*(f^{y_\alpha+\varepsilon }, f^{y_\alpha- \varepsilon })$ and $H_*(f^{x_\alpha+\varepsilon }, f^{x_\alpha- \varepsilon })$ have canonical generators (i.e. well defined, and not just up to a constant multiple).

But a generator on the left-hand side has its image zero in $H_*(f^{y_\alpha+\varepsilon }, f^{x_\alpha- \varepsilon } )$ by assumption. Therefore this generator has an image in $H_{*-1}(f^{y_\alpha- \varepsilon }, f^{x_\alpha- \varepsilon } )$ that lies in the image of $\partial$\,. It is thus equal to the image by $\partial$ of $\kappa(\alpha)$ times a generator. \\
Now consider the Morse function on $ {\mathbb R} P^2$ obtained by perturbing the following Morse-Bott function: 
$$[x_0,x_1,x_2] \mapsto x_2^2$$
where $[x_0,x_1,x_2]$ is the class of $(x_0,x_1,x_2) \in S^2$ by the equivalence relation $(x_0,x_1,x_2)\simeq (-x_0,-x_1,-x_2)$\,. This Morse-Bott function has a point of index $2$ at $[0,0,1]$\,, and a circle of index $0$ at $[\cos(\theta),\sin(\theta),0]$ for $\theta \in [0,\pi]$\,.
Perturbing this circle yields a pair of critical points of index $0$ and $1$\,, and the Thom-Smale complex is then 
$$\partial \underline{z}= 2\cdot \underline{y}, \partial \underline{y}=0, \partial \underline{x}=0$$
represented as
\centerline{
\xymatrix{&\underline{z}\ar[d]^{2}&\\&\underline{y}\ar[d]^{0}&&\\&\underline{x}&&
}}
The Barannikov complex (on $\mathbb Q$ or $\mathbb R$) is then 

\centerline{
\xymatrix{&&z\ar@{-}[d]^{}&\\&&y&&\\&&x&&
}}
\vskip .3cm
But necessarily  $\kappa_{\alpha}=2$\,, hence $\kappa_{\alpha}^2=4$\,. 

For $ {\mathbb R} P^3$\,, which is orientable, we have the similar function $[x_0,x_1,x_2,x_3] \mapsto x_3^2$ where $x_0^2+x_1^2, + x_2^2+x_3^2=1$ and we identify $(x_0,x_1,x_2,x_3)$ and  $(-x_0,-x_1,-x_2,-x_3)$\,. We then have a maximum $x_3=\pm1$ of index $3$\,, and an $ {\mathbb R} P^2$  Morse-Bott critical submanifold, which after perturbation yields a critical point of index $0$\,, one of index $1$ and one of index $2$\,. 

The Thom-Smale complex is then

\centerline{
\xymatrix{&\underline{t}\ar[d]^{0}&\\&\underline{z}\ar[d]^{2}&\\&\underline{y}\ar[d]^{0}&&\\&\underline{x}&&
}}
so again $\kappa (z)=2$\,. 

To obtain any squared integer, we can consider the lens space $L(n,1)$ quotient of $S^3= \{ (z_0,z_1)\in {\mathbb C} ^2 \mid \vert z_0\vert^2+\vert z_1\vert^2=1\}$ by
 $$(z_0,z_1) \simeq
(\omega z_0, \omega z_1) $$
where $\omega$ is a primitive $n$-th root of unity. The function $(z_0,z_1)\mapsto \vert z_0\vert^2$ has two critical circles a minimum and a  maximum. After  perturbation we get 

\centerline{
\xymatrix{&\underline{t}\ar[d]^{0}&\\&\underline{z}\ar[d]^{n}&\\&\underline{y}\ar[d]^{0}&&\\&\underline{x}&&
}}

and then $\vert \kappa_{\alpha} \vert =n$\,, since $H_1(L(n,1), \mathbb Z)=\mathbb Z / n \mathbb Z, H_2(L(n,1), \mathbb Z)=0 $\,.

Now assume there is some function  $f$ on the manifold $V$ with  a given bar code $\mathcal{B}_{f}$\,, and we embed  $V$ into a manifold $X$\,. Consider the function $g_ \varepsilon (x)=d(x,V)^2+ \varepsilon \rho( d(x,V)^2) f(p(x))$ where $\rho$ is nonnegative, equal to $1$ near $0$ and  vanishes outside a neighbourhood of $0$\,. then for $\varepsilon >0$ mall enough, the lower part of the bar code of $g_ \varepsilon $ coincides with the bar code of $f$\,. As a result if there is a function with some $\kappa_\alpha=n$ on $V$\,, the same holds for $X$\,. Consider the function $f$ above on $L(n,1)$\,, and normalize it so that the critical points are $0,1/3, 2/3, 1$\,. 
Consider the subset $\Lambda(n,1)=\{x \in L(n,1) \mid 1/4 \leq f(x) \leq 3/4\}$\,. This is a Lens space with two punctures, hence embeds in  ${\mathbb R}^4$ as a subset  of a compact  hypesurface $\Sigma_{n,1}$: if $\Lambda(n,1)$ is contained in $\{x \in {\mathbb R}^4 \mid \psi(x)=0\}$ and extending $\psi$ to a proper function having $0$ as a regular value, we set $\Sigma_{n,1}=\psi^{-1}(0)$\,. 
Now we can extend $f$ to a function $\tilde f$ on $\Sigma_{n,1}$ and its bar code contains $\mathcal{B}_{f}$\,. Applying the previous argument, we get a function close to $d(x,\Sigma_{n,1})^2$  containg $\mathcal{B}_{f}$ in its bar code. Since near infinity, $d(x,\Sigma_{n,1})^2$ is close to $ \vert x \vert ^2$\,, we get a function $F$ on the ball, with $F\leq c$  and $F=c$ near the boundary with arbitrary $\kappa_\alpha$\,. By embedding the ball in any $4$- manifold $M$\,, we get a function on $M$ with $\kappa_\alpha =n$\,. Again by embedding, this works on any manifold of dimension $\geq 4$\,.

\end{proof} 

More generally if for some prime $p$\,, the homology mod $p$ has rank  different from the rational  homology there must be a $y$ such that $p$ divides $\kappa_{\alpha}$\,. 

\begin{remark} \begin{enumerate} 
\item 
The converse does not hold, i.e. we may have $\kappa_{\alpha}\neq \pm 1$ while the homology has the same rank for all fields. For example if we have a Morse complex containing the following diagram 

\centerline{
\xymatrix{&\underline{z}\ar[dl]^{p}\ar[ddr]^{1}&\\\underline{y}_0\ar[ddr]^{1}&&\\ &&\underline{y}_1\ar[dl]^{-p}\\&\underline{x}&&
}}

the corresponding homology vanishes and the rational  Barannikov complex is  

\centerline{
\xymatrix{&z\ar@{-}[d]^{}&\\&y_0&&\\ &&y_1\ar@{-}[d]^{}\\&&x&
}}

but mod $p$\,, we get 
\centerline{
\xymatrix{&z\ar@{-}[dd]^{}&\\&&y_0\ar@{-}[dd]&&\\ &y_1&\\&&x&
}}

But in both cases the homology vanishes. Note however that if we look at the homology of sublevels, we can distinguish the two situations : 
if $a<f(x)< f(y_1)<c<f(y_0)$ the rank of the  homology $H^*(f^c,f^a)$ depends on the coefficient field : for $k= \mathbb Q$ we get $0$ while mod $p$\,, we get $2$\,. 
\item When several critical  values coincide, the numbers $k_\alpha$ are replaced by integral matrices. 
For example if we have the following bar code

\centerline{
\xymatrix{y_0\ar@{-}[d]^{}&y_1\ar@{-}[d]^{}&y_2\ar@{-}[d]^{}&y_3\ar@{-}[dd]^{}\\x_0&x_1&x_2&\\ &&&x_3
}}
 and if $a< x_3 <b < x_2 < c < y_2 < d$\,, we have the map 

$$\xymatrix{& H_*(f^{y_\alpha+\varepsilon }, f^{y_\alpha- \varepsilon } )\simeq  {\mathbb Z}^4\ar[d]^{\partial}\\H_{*-1}(f^{x_\alpha+ \varepsilon }, f^{x_\alpha- \varepsilon })\simeq {\mathbb Z}^3\ar[r]& H_{*-1}(f^{y_\alpha-\varepsilon }, f^{x_\alpha- \varepsilon } )\ar[d]\\
&H_*(f^{y_\alpha+\varepsilon }, f^{x_\alpha- \varepsilon } )\simeq {\mathbb Z}}
$$  hence we get a matrix $\bm\kappa \in M(4,{3, \mathbb Z})$ such that $\bm\kappa\otimes {\mathbb R} $ is surjective. We can then consider the singular values of $M$\,, and we get three numbers $\kappa_1,\kappa_2, \kappa_3$\,, however these are not the homological constants that will yield the prefactor of the eigenvalues, since we must first  compose with diagonal matrices  depending on the Hessian at each critical point involved (see Proposition \ref{pr:genalMorse})
. 

\end{enumerate} 
 \end{remark}

\subsection{Simple critical values for non Morse functions}
\label{sec:appnonMorse}

We consider here cases where
changing the function $f$ from $f_{1}$ to $f_{2}$ leads to explicit
changes of the global quasimodes $(\varphi_{j}^{(p)})_{j\in
  \mathcal{J}^{(p)}(a,b)}$ and provides accurate formulas, even for the
  subexponential factor, already known when $f_{1}$ 
is a generic Morse
function. It works especially well for functions, i.e. for $p=0$\,, and
although we are not considering Dirichlet boundary conditions at
$f^{-1}(\left\{b\right\})$ in $f_{a}^{b}$\,, 
like it is done in the  study of
quasi-stationary distributions, 
this sketches possible generalizations of the analyses made in
\cite{LeNi, DLLN1,DLLN2, LeNe1, LeNe2}.
Note however that, though obtaining precise estimates on the low spectrum
of the corresponding Witten Laplacians
with Dirichlet boundary conditions
 is an important step in the studies made  in \cite{LeNi, DLLN1,DLLN2, LeNe1, LeNe2},
these works actually focus on further issues
such as
the
exit events or
 the  concentration of  the associated quasi-stationary distributions. In particular,  
in \cite{DLLN1} are considered rare exit events,
which are actually rather related with the low spectrum
of appropriate Witten Laplacians with mixed
Dirichlet--Neumann boundary conditions.
 Simple cases when $p\neq 0$ will also be discussed afterwards.
\subsubsection{Degenerate local minima}
\label{sec:degmin}
We consider a reference function $f_{1}$ which is a generic Morse
function like in Theorem~\ref{th.LNV} with a bar code
$\mathcal{B}_{f_{1}}=([a_{1,\alpha}^{*},b_{1,\alpha}^{*+1}[)_{\alpha\in
  \mathcal{A}_{1}}$\,. In particular in degree $0$\,, there is one bar
$[a_{1,0}^{(0)},+\infty[=[x_{1,0}^{(0)},y_{1,0}^{(1)}[$ associated with the global minimum
$a_{1,0}^{(0)}$ and the sublevel set 
$\Omega_{1,0}^{(0)}=M=f_{1}^{+\infty}$\,, and there are bars
$[x_{1,k}^{(0)},y_{1,k}^{(1)}[\in \mathcal{A}_{1,c}$\,, $1\leq k\leq K_{0}$ where
$y_{1,k}^{(1)}$ is the value of saddle point and $x_{1,k}^{(0)}$ is the global
minimum value of the newly created connected component  
$\Omega_{1,k}^{(0)}$
of $f_{1}^{y_{1,k}^{(1)}}$\,, when we pass from the
sublevel set 
$f_{1}^{y_{1,k}^{(1)+0}}$ to $f_{1}^{y_{1,k}^{(1)}-0}$\,.\\
We take  $\ell_{min}^{(0)}<\min
\left\{y_{1,k}^{(1)}-x_{1,k}^{(0)}, |x_{1,k}^{(0)}-x_{1,k'}^{(0)}|\,, 0\leq k<k'\leq K_{0}\right\}$ and we
assume that the function $f_{2}$ satisfies Hypothesis~\ref{hyp:mainf}
and coincides with $f_{1}$ except around the local minima. The open
set called
$$
\omega_{k}^{(0)}=\Omega_{1,k}^{(0)}\cap f_{1}^{x_{1,k}^{(0)}+\frac{\ell_{0}}{2}}\,,
$$
is a connected open neighborhood of $x_{1,k}^{(0)}$ for all $k=0,\ldots,K_{0}$\,.
 The two functions $f_{1}$ and $f_{2}$ are compared by:
\begin{description}
\item[i)] $f_{1}\equiv f_{2}$ in a neighborhood of
  $M\setminus(\sqcup_{k=0}^{K_{0}}\omega_{k}^{(0)})$\,;
\item[ii)]$\|f_{1}-f_{2}\|_{\mathcal{C}^{0}}\leq \frac{\ell_{0}}{4}$\,.
\end{description}
Those two assumptions combined with the stability theorem
$$
d_{bot}(\mathcal{B}(f), \mathcal{B}(g))\leq \|f-g\|_{\mathcal{C}^{0}}
$$
recalled in Appendix~\ref{sec:stab}, ensure that there are exactly $K_{0}+1$ bars
$[x_{2,k}^{(0)},y_{2,k}^{(0)}[$ of degree $0$ and length larger that
$\frac{\ell_{0}}{2}$\,, where the saddle points are not changed
$y_{2,k}^{(1)}=y_{1,k}^{(1)}$ for $1\leq k\leq K_{0}$\,.
 Additionally  and especially because with our
choice of $\ell_{0}<\min\left\{|x_{1,k}^{(0)}-x_{1,k'}^{(0)}|\,, k<k'\right\}$ and \textbf{ii)},
the associated connected
component  remain unchanged
as well $\Omega_{2,k}^{(0)}=\Omega_{1,k}^{(0)}$ for  $0\leq k\leq K_{0}$\,.  We drop the index $j=1,2$ for
$\Omega_{k}^{(0)}$ and $y_{k}^{(1)}$\,. 
\begin{figure}[h]
\centering{
\includegraphics[width=10cm]{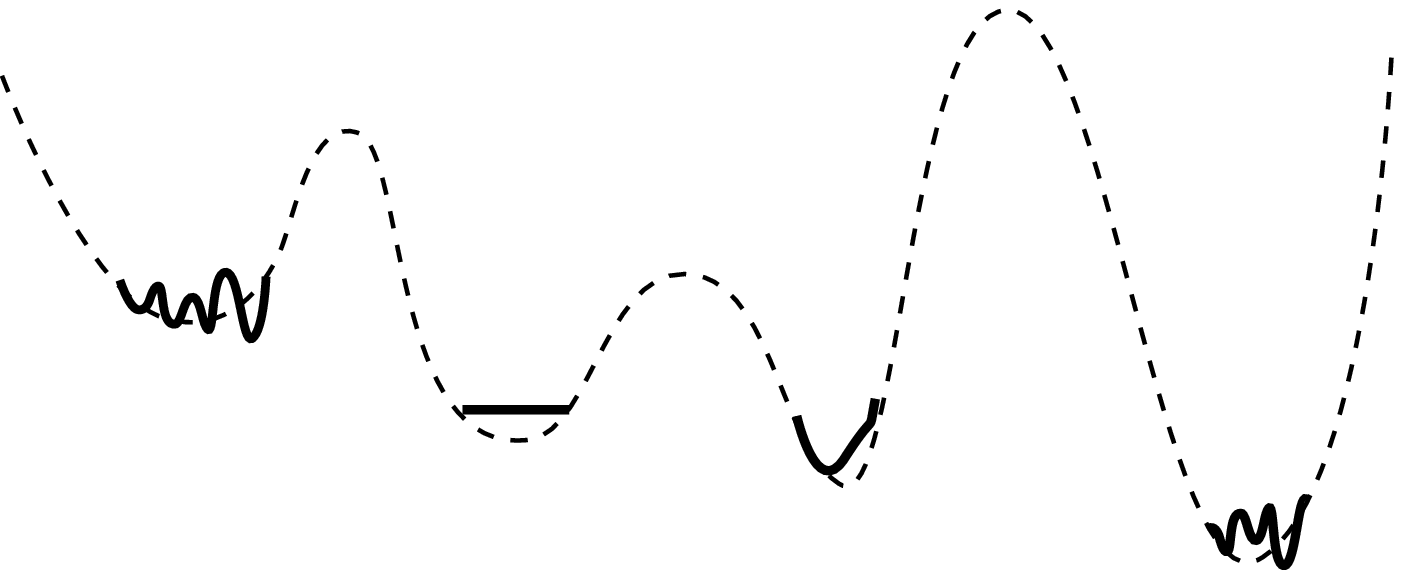}
\captionsetup{labelformat=empty,width=14cm}
\caption{\textbf{Figure 15:}
 The function $f_{1}$ is represented by
  dashed lines
and the
  modification giving $f_{2}$ by plain lines.}
}
\end{figure}
Like in the previous Subsection, we use the notation $\underline{s}$
for the point associated with the critical value $s$\,, when it is
uniquely defined.
\begin{prop}
 Under the above assumptions and in particular the comparison
 i)ii) between $f_{1}$ and $f_{2}$\,, the
 $\tilde{o}(e^{-\frac{\ell_{0}}{h}})$ eigenvalues of
 $\Delta_{f_{2},h}^{(0)}$ are given by
 \begin{equation}
\label{eq:asympdegmin}
\frac{h|\lambda_{1}(\underline{y}_{k}^{(1)})|}{\pi \left|\det
    \mathrm{Hess}~f_{1}(\underline{y}_{k}^{(1)})\right|^{1/2}}\frac{e^{-2\frac{y_{k}^{(1)}-x_{2,k}^{(0)}}{h}}}{\left(\pi
  h\right)^{-d/2}\int_{\Omega_{k}^{(0)}}e^{-2\frac{f_{2}(x)-x_{2,k}^{(0)}}{h}}~dx}\times (1+\mathcal{O}(h))
\end{equation}
as $h\to 0$ for all $k=0,1,\ldots, K_{0}$ (it is exactly $0$ for
$k=0$)\,.
\end{prop}
With this formula it then suffices to apply the Laplace method for the integral
$\int_{\Omega_{k}^{(0)}}e^{-2\frac{f(x)-x_{2,k}^{(0)}}{h}}~dx$ in order
to exhibit various asymptotic behaviours as $h\to 0$ of the
subexponential factor. We refer in particular to \cite{AGV} for the case when $f$  
is a multidimensional  polynomial  function.
\begin{proof}
When we work with functions, we are actually in the simpler framework
of \cite{HKN} for the generic Morse function $f_{1}$\,.
The problem consists in computing the square modulus of
the interaction $\langle \psi_{k}^{(1),h}\,,\,
d_{f,h}T_{\delta_{2}}\varphi_{k}^{(0),h}\rangle$ where $\psi_{k}^{(1),h}$ is a
local WKB-approximation of eigenvectors of $\Delta_{f,h}^{(1)}$ 
around the point
$\underline{y}_{k}^{(1)}$ while $\varphi_{1,k}^{(0),h}$
is a global quasimode associated with the bar
$[x_{1,k}^{(0)},y_{k}^{(1)}[$\,, solving
$d_{f,h}\varphi_{1,k}^{h}=0$ in $\Omega_{k}^{(0)}\cap
f^{y_{k}^{(1)}-\delta(h)}$ with $\lim_{h\to 0}\delta(h)= 0$\,. The truncation
$T_{\delta_{2}}$ is a smooth truncation around the level
$y_{k}^{(1)}-\delta_{2}$ with $\delta_{2}>0$ small. 
By Theorem~\ref{th:induc} and Theorem~\ref{th:specres} the same method
holds by replacing the global quasimodes $\varphi_{1,k}^{(0),h}$ by global
quasimodes $\varphi_{2,k}^{(0),h}$constructed in
Theorem~\ref{th:induc}. In details we refer more specifically to
the consequences stated in 
Subsection~\ref{sec:conseqN}. Moreover we can focus on the
bars of length larger that $\frac{\ell_{0}}{2}$ which are
$([x_{2,k}^{(0)},y_{k}^{(1)}[)_{k=0,\ldots,K_{0}}$\,. Since those
quasimodes satisfy $d_{f_{2},h}\varphi_{2,k}^{(0),h}=0$ in
$\Omega_{k}^{(0)}\cap f^{y_{k}^{(1)}-\delta(h)}$ they equal
$\sqrt{C_{k,h}}e^{-\frac{f_{2}(x)-x_{2,k}^{(0)}}{h}}$ where $C_{k,h}$ is the
normalization constant
$$
C_{k,h}=
\frac{1}{
\int_{\Omega_{k}^{(0)}\cap f^{y_{k}^{(1)}-\delta(h)}}e^{-2\frac{f_{2}(x)-x_{2,k}^{(0)}}{h}}~dx
}
=
\frac{1+\tilde{o}(1)}{
\int_{\Omega_{k}^{(0)}}e^{-2\frac{f_{2}(x)-x_{2,k}^{(0)}}{h}}~dx
}
$$
which replaces 
$$
 \frac{1}{\int_{\Omega_{k}^{(0)}}
 e^{-2\frac{f_{1}(x)-x_{1,k}^{(0)}}{h}}~dx}
=\left(\pi h\right)^{-d/2}
|\mathrm{det}~\mathrm{Hess}~f(x_{1,k}^{(0)})|^{1/2}
\times (1+\mathcal{O}(h))
\,.
$$
Finally it suffices to notice that up to the normalization constant and the
change of the length of the bar which brings another constant factor, the functions $\varphi_{1,k}^{(0),h}$
and $\varphi_{2,k}^{(0),h}$ coincide in the neigborhood of
$\underline{y}_{k}^{(1)}$ and the local computation of the interaction
is not changed.
\end{proof}
The above formula shows a good stability when $f_{1}$ is changed into
$f_{2}$ although such a stability may not appear when we make an
explicit asymptotic expansion of the Laplace integral
$\int_{\Omega_{k}^{(0)}}e^{-2\frac{f(x)-x_{2,k}^{(0)}}{h}}~dx$\,. Here
is an example in dimension $1$\,, that is for functions defined on
$\sz^{1}=\rz/(2\pi \zz)$\,. The function $f_{1}$ is assumed to have
four non degenerate critical points:
\begin{itemize}
\item at
 $\underline{x}_{1,1}^{(0)}=0$ with value $x_{1,1}^{(0)}=0$ and second derivative
 $1$\,;
\item at $\underline{x}_{0,1}^{(0)}=\pi$ with value $x_{0,1}^{(0)}=-1$\,, the global
  minimum;
\item at $\underline{y}_{1,1}^{(1)}=\frac{\pi}{2}$ with value $y_{1,1}^{(1)}=1$ and
  the second derivative equal to $-\lambda_{1}$\,;
\item at $\underline{y}_{0,1}^{(1)}=\frac{3\pi}{2}$ with value
  $y_{0,1}^{(1)}=2$\,, the global maximum.
\end{itemize}
The modified function $f_{2,\delta}$ parametrized by $\delta\in
\rz$\,, $\delta$  small, and consists in replacing
$f_{1}(x)=\frac{x^{2}}{2}+O(x^{3})$  in a small neighborhood
$[-\varepsilon,\varepsilon]$  of
$\underline{x}_{1,1}^{(0)}=0$ by
$$
f_{2,\delta}(x)=\frac{x^{4}+2\delta
  x^{2}+1_{(-\infty,0]}(\delta)\delta^{2}}{4}\,,
$$
while $f_{2,\delta}\equiv f_{1}$ outside
$[-2\varepsilon,2\varepsilon]$\,.
Formula \eqref{eq:asympdegmin} then says that the
$\tilde{o}(e^{-\frac{1}{h}})$ non zero eigenvalue of
$\Delta_{f_{2,\delta},\sz^{1},h}^{(0)}$ (for $\delta>0$ and
$\varepsilon>0$ small enough) equals
$$
\frac{h\sqrt{\lambda_{1}}e^{-\frac{2}{h}}}{\pi (\pi
  h)^{-1/2}\int_{\rz}e^{-\frac{x^{4}+2\delta x^{2}+1_{(-\infty,0]}(\delta)\delta^{2}}{2h}}~dx}
\times (1+\mathcal{O}(h))\,.
$$
It is equivalent as $h\to 0$ to 
\begin{eqnarray*}
  &\frac{h\sqrt{\lambda_{1}\delta}e^{-\frac{2}{h}}}{\pi}
\quad
&\text{when}~\delta> 0\,\\
&\frac{h^{5/4}\sqrt{\lambda_{1}}e^{-\frac{2}{h}}}{\sqrt{\pi}\int_{\rz}e^{-\frac{u^{4}}2}~du}&\quad\text{when}~\delta=0\,,\\
&\frac{h\sqrt{\lambda_{1}|\delta|}e^{-\frac{2}{h}}}{\sqrt 2\pi}
&\quad\text{when}~\delta<0\,.
\end{eqnarray*}
So the apparent discontinuity in the exponent of $h$ at $\delta=0$ is
a simple consequence of the discontinuity of the Laplace
integral. Actually the stability of persistence homology has a
stronger spectral counterpart than what is stated in
Theorem~\ref{th:stab1}: It does not concern only the exponential
scales but also allows to study the deformations of the asymptotic
subexponential factors provided that a robust formula can be proved for them.
The rest of this section explores different cases for which we are
able to prove such formulas.
\subsubsection{Variations}
\label{sec:varianonMorse}
In the previous paragraph we used a good enough knowledge of the
global quasimodes
$\varphi_{2,k}^{(0),h}=\sqrt{C_{k,h}}e^{-\frac{f_{2}(.)-x_{2,k})}{h}}$
in degree $p=0$\,, in order to get the explicit change in the asymptotic
formulas when we pass from the Morse function $f_{1}$ with simple
critical values to the function $f_{2}$ with degenerate local
minima. Such an analysis can be done in more general degree if we have
explicit enough information on the local forms of quasimodes the
global ones $\varphi_{k}^{(p),h}$ and the local ones
$\psi_{k}^{(p+1),h}$\,.  By duality this is obviously true in
dimension $1$ and we start with this example. We then consider other
possible extensions.

\paragraph{The one dimensional case with degenerate critical values}
Consider a $\mathcal{C}^{\infty}$ Morse funlction $f_{1}$ on $\rz$ such
that $|\partial_{x}f_{1}|\geq c$ for some positive constant $c$ when $x\in
\rz\setminus [-R,R]$ for $R>0$ large enough. For $-a=|a|$ and $b=|b|$
large enough the bar code $\mathcal{B}_{f_{1}}(a,b)$ does not change
when $a,b$ are changed, except for the value of the endpoints $a,b$\,,
while for such a fixed pair $(a,b)$ it can be viewed as a
restricted bar code $B_{\tilde{f}_{1}}(a,b)$ of a function
$\tilde{f}_{1}$ defined on $\sz_{1}$\,. This solves the compactness
problem in order to fit with our general framework. It can be checked
easily that in all such cases the exponentially small eigenvalues of
$\Delta_{f_{1},f_{1}^{-1}([a,b]),h}^{(p)}$ are close to the ones of
$\Delta_{f_{1},\rz,h}^{(p)}$  for $p=0,1$ and even that the endpoints of
the interval $f_{1}^{-1}([a,b])$ can be moved as long as they do not meet
the critical point without changing the final approximate 
spectral result (the same will be
 true for the function $f_{2}$). So let
us focus on $f_{1}^{-1}([a,b])$ with $-a=|a|$ and $b=|b|$ large. The
bar code is made of bars $[x_{k,1}^{(0)},y_{k,1}^{(1)}[$\,,
$k=1,\ldots, K$ with an additional bar:
\begin{itemize}
\item 
$[x_{0,1}^{(0)},b[$ if
$f_{1}\big|_{f_{1}^{-1}([a,b])}$ admits an interior global minimum at
$x_{0,1}^{(0)}=f(\underline{x}_{0,1}^{(0)})$\,,
$\underline{x}_{0,1}^{(0)}\in f_{1}^{-1}(]a,b[)$\,;
\item or $[a,y_{0,1}^{(1)}[$ if $f_{1}\big|_{f_{1}^{-1}([a,b])}$ admits an interior global maximum at
$y_{0,1}^{(1)}=f(\underline{y}_{0,1}^{(0)})$\,,
$\underline{y}_{0,1}^{(0)}\in f_{1}^{-1}(]a,b[)$\,.
\end{itemize}
\begin{figure}[h]
\centering{
\includegraphics[width=5cm]{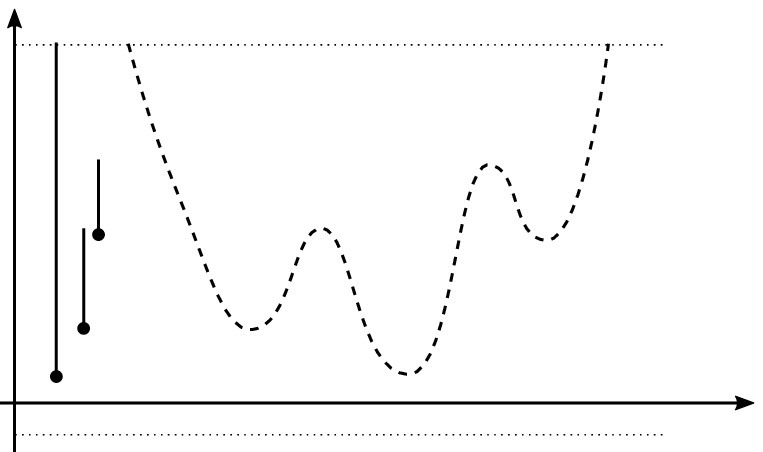}
\includegraphics[width=5cm]{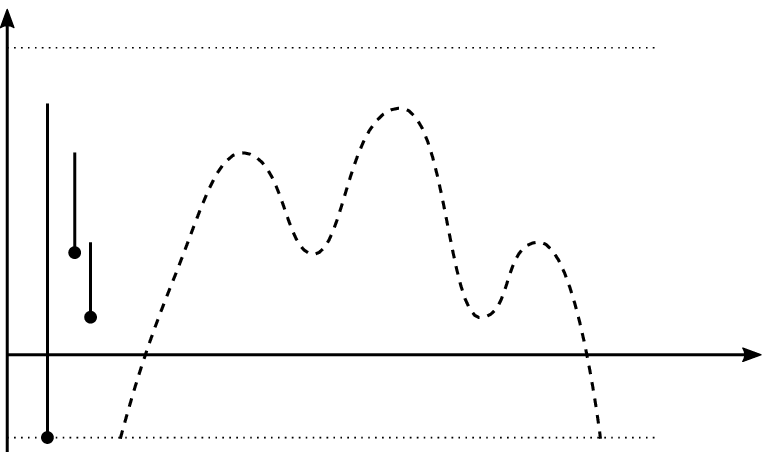}
\includegraphics[width=5cm]{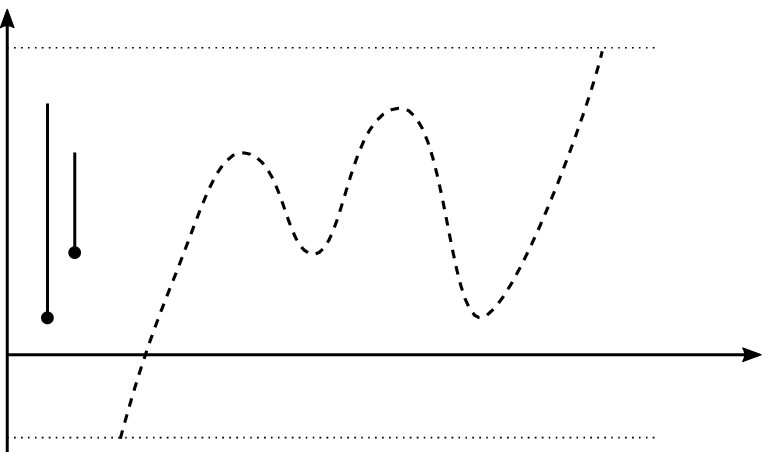}
\captionsetup{labelformat=empty, width=14cm}
\caption{\textbf{Figure 16:} Three different cases for $f_{1}$ between the level $a$ and
  $b$\,, from left-hand side to right-hand side with an interior global
  minimum, an interior
  global maximum in the interior and none of them. The bar code in
  $[a,b]$ is represented beside the $y$-axis.}
}
\end{figure}
Only in the first case, the Witten Laplacian
$\Delta_{f_{1},f_{1}^{-1}([a,b])}^{(0)}$ has a non trivial kernel $\cz
e^{-\frac{f_{1}(.)-x_{0,1}^{(0)}}{h}}$\,. Only in the second case, the
Witten Laplacian $\Delta_{f_{1},f_{1}^{-1}([a,b]),h}^{(1)}$ has a non
trivial kernel $\cz e^{\frac{f_{1}(.)-y_{0,1}^{(1)}}{h}}dx$\,. The two
cases are exclusive and a third one is when the global minimum value  of
$f_{1}\big|_{f_{1}^{-1}([a,b])}$ is $a$ and the global maximum value is
$b$\,. Depending on the cases   $f_{1}$ admits $2K+1$ or $2K$ distinct
critical values and their set  in $[a,b]$ is denoted $\mathcal{C}$\,.\\
In order to specify our modified function $f_{2}$ we first choose
$\ell_{0}<\min\left\{|c-c'|,c\neq c'\,, c,c'\in
  \mathcal{C}\right\}$\,. The connected open set $\Omega_{k,1}^{(0)}$
as the  connected component of $(f_{1})^{y_{k,1}^{(1)}}$ which
contains $\underline{x}_{k,1}^{(0)}$ for
$1\leq k\leq K$\,, with $\Omega_{0,1}^{(0)}=f_{1}^{-1}(]a,b[)$ if
the global minimum $\underline{x}_{0,1}^{(0)}\in f_{1}^{-1}(]a,b[)$
exists. By duality one defines $\Omega_{k,1}^{(1)}$ as the connected
component of $(f_{1})_{x_{k,1}^{(0)}}$ for $1\leq k\leq K$\,, with
$\Omega_{0,1}^{(1)}=f^{-1}(]a,b[)$ if the global maximum
$\underline{y}_{0,1}^{(1)}\in f_{1}^{-1}(]a,b[)$ exists.
Then the connected open sets $\omega_{k}^{(0)}$ and $\omega_{k}^{1}$
are defined by
$$
\omega_{k}^{(0)}=\Omega_{k,1}^{(0)}\cap
(f_{1})^{x_{k,1}^{(0)}+\frac{\ell_{0}}{4}}
\quad,\quad
\omega_{k,1}^{(1)}=\Omega_{k,1}^{(1)}\cap (f_{1})_{y_{k,1}^{1}-\frac{\ell_{0}}{4}}\,.
$$
The function $f_{2}$ satisfies Hypothesis~\ref{hyp:mainf} and 
\begin{itemize}
\item $f_{1}\equiv f_{2}$ in a neighborhood of
  $\rz\setminus(\sqcup_{0\leq k\leq K}(\omega_{k}^{(0)}\sqcup \omega_{k}^{(1)}))$
  where $\omega_{0}^{(0)}$ and $\omega_{0}^{(1)}$ are replaced by the
  empty set when they are not defined;
\item $\|f_{1}-f_{2}\|\leq \frac{\ell_{0}}{4}$\,.
\end{itemize}
Note in particular $f_{1}^{-1}([a,b])=f_{2}^{-1}([a,b])$\,.\\
\begin{figure}[h]
\centering{
\includegraphics[width=12cm]{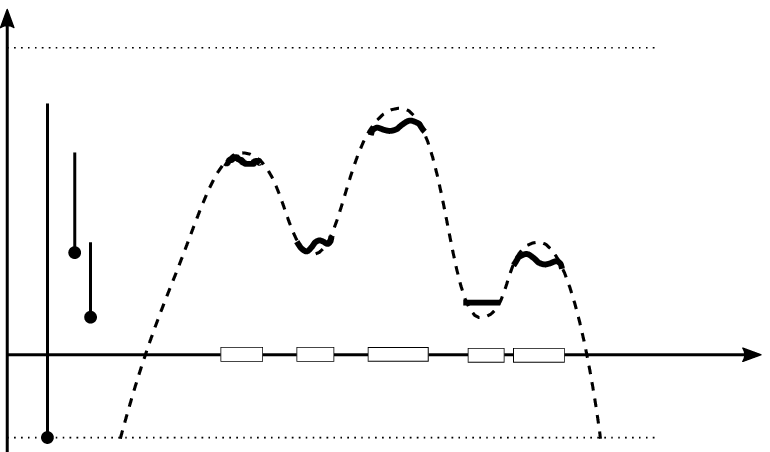}
\captionsetup{labelformat=empty, width=14cm}
\caption{\textbf{Figure 17:} The function $f_{1}$ is represented by the dashed curve, the
  open sets $\omega_{k}^{(p)}$ are materialized by the white rectangles
  along the $x$-axis and the modifications leading to $f_{2}$ by the
  plain curve.}
}
\end{figure}

Owing to the stability theorem
$$
d_{bot}(\mathcal{B}(f), \mathcal{B}(g))\leq \|f-g\|_{\mathcal{C}^{0}}
$$
the bars $[x_{k,1}^{(0)},y_{k,1}^{(1)}[$ are transformed into bars
$[x_{k,2}^{(0)},y_{k,2}^{(1)}[$ of length
$y_{k,2}^{(1)}-x_{k,2}^{(0)}>\frac{\ell_{0}}{2}$ while all the other
bars have length smaller thant $\frac{\ell_{0}}{2}$\,.
After those assumptions the spectral result take a nice simple form.
\begin{prop}
\label{pr:as1d}
For the values
$a,b$  and the function $f_{2}$ chosen like above, there are $K$ non zero $\tilde{o}(e^{-\frac{\ell_{0}}{h}})$
eigenvalues of
$\Delta_{f_{2},f_{2}^{-1}([a,b]),h}^{(0)\text{~or~}(1)}$ which are
equal to 
$$
\frac{1+\tilde{o}(1)}{(h^{-1}\int_{\omega_{k}^{(1)}}e^{2\frac{f(x)}{h}}~dx)\times
  (h^{-1}\int_{\omega_{k}^{(0)}}e^{-2\frac{f(x)}{h}}~dx)}\,,
\quad
k=1,\ldots, K\,.
$$ 
\end{prop}
\begin{proof}
By the usual supersymmetric arguments the non zero eigenvalues of
$\Delta_{f_{2},f_{2}^{-1}([a,b]),h}^{(0)}$ and
$\Delta_{f_{2},f_{2}^{-1}([a,b]),h}^{(1)}$ are the same in dimension
$1$ and we thus focus on $\Delta_{f_{2},f_{2}^{-1}([a,b]),h}^{(0)}$ or
more precisely on the non zero singular values of the restricted differential.
 We follow the general method which consist in computing the
 interaction scalar product $\langle \psi_{k}^{(1)}\,,
 d_{f,h}T_{\delta_{2}}\varphi_{k}^{(0)}\rangle$ where $\psi_{k}^{(1)}$
 is a local quasimode for $\Delta_{f_{2},h}$ in the neighborhood
 $\omega_{k}^{(1)}$ around $\underline{y}_{k,1}^{(1)}$ while
 $\varphi_{k}^{(0)}$ is a global quasimode associated with the bar
 $[x_{2,k}^{(0)},y_{2,k}^{(1)}[$ solving $d_{f,h}\varphi_{k}^{(0)}=0$
 in the connected component which contains $\omega_{k}^{(0)}$ 
of $f^{y_{k}^{(2)}-\delta(h)}$\,, with $\lim_{h\to 0}\delta(h)=0$\,. We
 work directly with the function $f_{2}$ the global quasimode
 $\varphi_{k}^{(0)}$ equals 
$$
\frac{1+\tilde{o}(1)}{\sqrt{\int_{\omega_{k}^{(0)}}e^{-2\frac{f_{2}(x)-x_{2,k}^{(0)}}{h}}dx}}e^{-\frac{f_{2}(.)-x_{2,k}^{(0)}}{h}}
\quad \text{in}~\omega_{k}^{(0)}
\,.
$$ 
By noticing that $\partial_{n}f_{2}\big|_{\partial
  \omega_{k}^{(1)}}=\partial_{n}f_{1}\big|_{\partial
  \omega_{k}^{(1)}}<0$\,, and by using Dirichlet boundary conditions
on $\partial \omega_{k}^{(1)}$ in degree $p=1$\,, we find that
$\psi_{k}^{(1)}$ can be chosen as
$$
\frac{1+\tilde{o}(1)}{\sqrt{\int_{\omega_{k}^{(1)}}e^{2\frac{f_{2}(x)-y_{2,k}^{(1)}}{h}}~dx}}e^{\frac{f_{2}(x)-y_{2,k}^{(1)}}{h}}~dx\quad \text{in}~\omega_{k}^{(1)}\,.
$$
A direct computation gives
$$
\langle \psi_{k}^{(1)}\,,\,
d_{f,h}T_{\delta_{2}}\varphi_{k}^{(0)}\rangle
=\pm \frac{he^{-\frac{y_{2,k}^{(1)}-x_{2,k}^{(0)}}{h}}(1+\tilde{o}(1))}{\sqrt{\int_{\omega_{k}^{(1)}}e^{2\frac{f_{2}(x)-y_{2,k}^{(1)}}{h}}~dx}\times
\sqrt{\int_{\omega_{k}^{(0)}}e^{-2\frac{f_{2}(x)-x_{2,k}^{(0)}}{h}}~dx}}
$$
where the factor $e^{-\frac{y_{2,k}^{(1)}-x_{2,k}^{(0)}}{h}}$ can be simplified.
\end{proof}
\begin{remark}
Note that in this proof the result on the generic Morse function
$f_{1}$ is not used. The function $f_{1}$ was introduced in order to
have a simple formulation of the assumptions fulfilled by
$f_{2}$\,. The result actually comes from a direct computation
when we know well enough the local forms of the global
($\varphi_{k}^{(0)}$) and local ($\psi_{k}^{(1)}$)
quasimodes\,. We have explicit form in dimension $1$ and the
computation is straightforward. It is not the same in  the multidimensional
case although Stokes's formula allows to perform the computation when
local approximations of local and global quasimodes are well known.
\end{remark}

\paragraph{Piecewise affine functions}
In this paragraph we make more explicit the one dimensional result
when $f$ is a continuous piecewise affine function and discuss the possible
extension to the multidimensional case.
Let $f$ be a piecewise affine function on $\rz$ such that:
\begin{itemize}
\item the derivative $f'$ does not vanish when it is defined;
\item there exists $R>0$ such that the derivative $f'$ is a 
 constant on $[R,+\infty)$ and on $(-\infty,-R]$\,;
\item the values $f(x)$ of the points $x$ where $f'$ is discontinuous
  are all distinct.
\end{itemize}

\begin{figure}[h]
\centering{
\includegraphics[width=10cm]{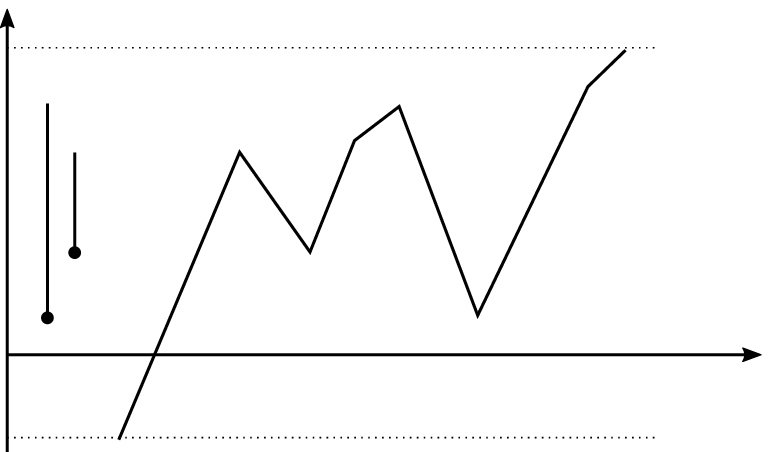}
\captionsetup{labelformat=empty, width=14cm}
\caption{\textbf{Figure 18:} A piecewise affine potential in $1D$ with distinct and some
  fake critical values.}
}
\end{figure}
Such a function $f$ can be written as a function $f_{2}$ of the
previous paragraph (simply regularize locally the discontinuous 
change of slopes in order to 
get the Morse function $f_{1}$). \\
The extension of Proposition~\ref{pr:as1d} to $a=-\infty$ and
$b=+\infty$ says that the $\tilde{o}(1)$-eigenvalues of
$\Delta_{f,\rz,h}^{(0)}$ (and by duality of $\Delta_{f,\rz,h}^{(1)}$)
are given by
$$
H[|f'(\underline{y}_{k}^{(1)}+0)|,f'(\underline{y}_{k}^{(1)}-0)]H[f''(\underline{x}_{k}^{(0)}+0),|f'(\underline{x}_{k}^{(0)}-0)|]e^{-2\frac{y_{k}^{(1)}-x_{k}^{(0)}}{h}}(1+\tilde{o}(1))\quad,\quad
k=1,\ldots, K\,,
$$  
where the finite length bars of $\mathcal{B}_{f}$ are
$[x_{k}^{(0)},y_{k}^{(1)}[$\,, $k=1,\ldots,K$\,;
$x_{k}^{(0)}$ is the local strict minimal value around the point
$\underline{x}_{k}^{(0)}$\,; $y_{k}^{(1)}$ is the local maximal value
around the point $\underline{y}_{k}^{(1)}$\,; $f'(x+0)$ and $f'(x-0)$
denote respectively the right and left derivative and
$H[s,t]=\frac{2st}{s+t}$ is the harmonic mean of $s,t>0$\,.\\
The computation when $f$ is constant on some intervals is also possible
with a subexponential factor behaving like $h$ or $h^{2}$\,, depending on the
different cases (left to the reader).\\

Now let $f$ be a piecewise affine function defined on a finite
triangulation of $\rz^{d}=\sqcup_{1\leq i\leq I}\mathcal{T}_{i}$ where
$\mathcal{T}_{i}$ is a $d$-dimensional non degenerate simplex  with endpoints
$A_{i}^{0},\ldots, A_{i}^{d}$ and where non finite simplices are
roughly taken into account by sending the first endpoint to infinity
$A_{i}^{0}=\infty$ (a more precise description is not necessary
here). We assume that $\lim_{x\to \infty}f(x)=+\infty$\,. The function
$f$ is a subanalytic function on $\rz^{d}$ of which the restriction to
any ball $B(0,R)$ can be viewed as the restriction to $B(0,R)$ of a
subanalytic function defined on $\sz^{d}$\,. This solves the
compactness problem or the questions about the topology at
infinity (alternatively we could work on the $d$-dimensional flat torus). The function $f$ has a finite number of horizontal strata
according to the terminology of Definition~\ref{de:horstrata}\,, which
contain all the critical values and the possible endpoints of the
bar code $\mathcal{B}_{f}$\,. We may consider either $\Delta_{f,\rz^{d},h}$ or by approximation
$\Delta_{f,f^{-1}([a,b]),h}$ with $-a,b>0$ large enough. According to
our analysis in Subsection~\ref{sec:moregenLip}, in particular
Proposition~\ref{pr:verhypLipf}, Proposition~\ref{pr:subLipAgm1} and Proposition~\ref{pr:subLipAgm}, the results of
Theorem~\ref{th:induc} hold in this case and we know that the
exponentially small eigenvalues of $\Delta_{f,f^{-1}([a,b]),h}^{(p)}$ satisfy
$$
\lambda_{\alpha}^{(p)}(h)\stackrel{log}{\sim}e^{-2\frac{y_{\alpha}^{*+1}-x_{\alpha}^{*}}{h}}\,,
$$
where $\alpha$ belongs to  $A_{c}^{(p)}(a,b)\sqcup
A_{c}^{(p-1)}(a,b)$\,.\\
The question is whether it is possible to give  algebraic formulas 
for the accurate asymptotic behaviour
 as this is done easily in the one dimensional
case. For such a function $f$\,, the Witten Laplacian $\Delta_{f,h}$ is a matricial Schr{\"o}dinger
operator with a singular potential. Many things are known on
scalar Schr{\"o}dinger operators with singular potentials
(see e.g. \cite{AGHKH}\cite{BGP}), but little seems to be known for
those Witten Laplacians, and especially when we think about the
algebraic topology subtleties. 
We may also start directly, instead of
$\rz^{d}$\,, on a Lipschitz manifold made of glued simplexes, with a
function $f$ which has a constant gradient along every simplex. The
functional analysis of Hodge Laplacian on Lipschitz manifold has been
considered in \cite{GMM,MMMT}. An accurate analysis of the low lying
spectrum of such Witten Laplacians would provide a large family  of
discrete and easily encoded models, from the point of view of data and
hopefully of results, which could be used as approximations of complicated
realistic situations. It would be interesting to compare with the
approach starting from purely discrete models on graphs as presented
in \cite{CdVPY}.

\paragraph{Critical submanifolds}
This case is related with degenerate Witten Laplacians studied in connection with Bott-Morse
inequalities (see e.g. \cite{Bis,HeSj6}).
We consider here simple examples where we have a critical submanifold
instead of a critical point.
We start with the mexican hat function
$f(r,\theta)=\frac{r^{4}}{4}-\frac{r^{2}}{2}+\frac{1}{4}$ in polar
coordinates $(r,\theta)$ in $\rz^{2}$ with the euclidean metric $dr^{2}+(rd\theta)^{2}$\,, which admits a non degenerate
maximum at $r=0$ with $f(0_{\rz^{2}})=\frac{1}{4}$ and a degenerate minimum at
$r=1$ with $f(1,\theta)=0$\,.\\
The bar code of the function $f$ is made of the bar $[0,+\infty[$ in
degree $0$ and the bar $[0,\frac{1}{4}[$ in degree $1$\,. We compute
the non zero exponentially small eigenvalue of
$\Delta_{f,\rz^{2},h}^{(p)}$ with $p=1$ or $2$ by computing 
the interaction scalar product $\langle \psi_{1}^{(2)}\,,\,
d_{f,h}^{(1)}T_{\delta_{2}}\varphi_{1}^{(1)}\rangle$ where
$\varphi_{1}^{1}$ is a global quasimode $1$-form 
associated with the bar
$[0,\frac{1}{4}[$  and $\psi_{1}^{(2)}$ is a local quasimode $2$-form
around $r=0$\,.\\ 
In this particular example we have explicit forms
for $\varphi_{1}^{(1)}$ and $\psi_{1}^{(2)}$:
\begin{itemize}
\item We take $\nu>0$ smaller than the truncation parameter
  $\delta_{2}$\,. Then a explicit normalized element of
  $\ker(\Delta_{f,f^{-1}([-1,\frac{1}{4}-\nu]),h}^{(1)})$ is given by
$$
\varphi_{1}^{(1)}=\frac{1}{\sqrt{\int_{f^{\frac{1}{4}-\nu}}e^{-\frac{\frac{r^{4}}{2}-r^{2}+\frac{1}{2}}{h}}r^{-2}dr (rd\theta)}}e^{-\frac{\frac{r^{4}}{4}-\frac{r^{2}}{2}+\frac{1}{4}}{h}}d\theta\,.
$$ 
\item For the local quasimode $\psi_{1}^{(2)}$ defined around $r=0$\,, we
can use either a WKB approximation, or by duality 
 the exact normalized element of
 $\ker(\Delta_{f,f^{-1}([\frac{1}{4}-\delta]),h}^{(2)})$ ($\delta>0$
 is small enough but bigger than $2\delta_{2}$) given by
equal to 
$$
\psi_{1}^{(2)}=\frac{1}{\sqrt{\int_{f_{\frac{1}{4}-\delta}}
    e^{\frac{\frac{r^{4}}{2}-r^{2}}{h}}~dr(rd\theta)}}
e^{\frac{\frac{r^{4}}{4}-\frac{r^{2}}{2}}{h}}dr\wedge
(rd\theta)\,.
$$
\end{itemize}
The scalar product 
$\langle \psi_{1}^{(2)}\,,\,
d_{f,h}(T_{\delta_{2}}\varphi_{1}^{(1)})\rangle$ is then equal to 
$$
\frac{1}{\sqrt{\int_{f_{\frac{1}{4}-\delta}}
    e^{\frac{\frac{r^{4}}{2}-r^{2}}{h}}~dr(rd\theta)}\sqrt{\int_{f^{\frac{1}{4}-\nu}}e^{-\frac{\frac{r^{4}}{2}-r^{2}+\frac{1}{2}}{h}}r^{-2}dr (rd\theta)}}
\langle dr\wedge(rd\theta)\,,\, h\chi'_{\delta_{2}}(r)dr \wedge d\theta\rangle e^{-\frac{1}{4h}}\,,
$$
where
$$
\langle dr\wedge (rd\theta)\,,\, h\chi'_{\delta_{2}}(r)dr \wedge d\theta\rangle
=\pm h\int_{r=\varrho}\frac{rd\theta}{r}=\pm 2\pi h
$$
does not depend on the value $\varrho>0$ (This is an explicit
illustration of Stokes's formula argument used in \cite{LNV} when $f$ is
a Morse function).\\
Using the asymptotics of non degenerate Laplace integrals, the
non zero exponentially small of $\Delta_{f,\rz^{2},h}^{(p)}$\,, for
$p=1,2$ equals
$$
\frac{1+O(h)}{\pi h}\times\frac{1+O(h)}{\pi (2\pi h)^{1/2}
}\times (2\pi h)^{2}e^{-\frac{1}{2h}}=\frac{2\sqrt{2}h^{1/2}+O(h^{3/2})}{\sqrt{\pi}}e^{-\frac{1}{2h}}\,.
$$
The subexponential factor $Cte\times \sqrt{\frac{h}{\pi}}$ differs from the
asymptotic behaviour $Cte\times \frac{h}{\pi}$ obtained when $f$ is a
generic Morse function. Actually it is possible to study the
transition from the Morse generic case to this degenerate case
 by taking $f_{\delta}(r,\theta)=f(r,\theta)+\delta
\gamma(r)\cos(\theta)$ where $\gamma\in \mathcal{C}^{\infty}(
]0,+\infty[;[0,1])$ equals $1$ in a neigborhood of $1$\,, and
$\delta\in\rz$ is chosen small enough. Let us illustrate this in a
larger framework. Note that the above formula is not changed if the
metric $dr^{2}+r^{2}d\theta^{2}$ is replaced by $dr^{2}+d\theta^{2}$
in a neighborhood of $r=1$\,. This will make the forthcoming analysis
simpler.\\
We consider a $\mathcal{C}^{\infty}$  function $f$ on the compact
Riemannian manifold $M$ with a finite number of critical values, which
are all non degenerate and simple except the critical value fixed to be $0$\,. We
further assume:
\begin{itemize}
\item the critical set  around the value $0$ is a closed orientable submanifold
  $M'$ of dimension $p$\,;
\item there  is a tubular neighborhood of $M'$ which is a product of
  two Riemannian manifolds $M'\times M''$ with the metric
  $g=g'\ds\mathop{\oplus}^{\perp}g''$\,; a corresponding local coordinate system
  is written $x=(x',x'')$\,;
\item in the tubular neighborhood $M'\times M''$ the function $f$ is a
  function of $x''\in M''$  and has a unique minimum $f(x_{0}'')=0$\,;
\item the bar code $\mathcal{B}_{f}$ contains a unique bar
  $[0,y_{1}^{(p+1)}[$ of degree $p$ with lower endpoint $0$ and upper
  endpoint $y_{1}^{(p+1)}<+\infty$\,; the eigenvalues of the Hessian
  at the corresponding point $\underline{y}_{1}^{(p+1)}$ are denoted
  $-\lambda_{1}(\underline{y}_{1}^{(p+1)}),\ldots,
  -\lambda_{p+1}(\underline{y}_{1}^{(p+1)})$ and $\lambda_{p+2}(\underline{y}_{1}^{(p+1)}),\ldots,
  \lambda_{d}(\underline{y}_{1}^{(p+1)})$\,;
\item a local unstable (for $-\nabla f$) closed cell around the non
  degenerate critical point $\underline{y}_{1}^{(p+1)}$ is denoted $e_{1}^{(p+1)}$ and its boundary in $M$ which
  is a $p$-dimensional sphere is denoted by $\partial
  e_{1}^{(p+1)}$\,;
\item if $\phi$ is  $\mathcal{C}^{\infty}$ Morse 
function on $M'$ with the
  maximal value $0$ and
  $\chi\in \mathcal{C}^{\infty}_{0}(M'';[0,1])$ is equal to $1$ in a
  neighborhood of $x_{0}''$ and such that $f(x'')\geq c>0$ on $\supp d\chi $\,, the function $f_{\delta}$ is defined as
  $f_{\delta}=f+\delta \chi(x'') \phi (x')$\,;
\item for the sake of simplicity we work in the energy interval
  $[a,b]$ with $a=-\varepsilon$ and $b=y_{1}^{(p+1)}+\varepsilon$
  where $\varepsilon>0$ is fixed so that the critical values of $f$ in
  $[a,b]$ are the ones contained in $[0,y_{1}^{(p+1)}]$\,.
\end{itemize}

\begin{figure}[h]
\centering{
\includegraphics[width=8cm]{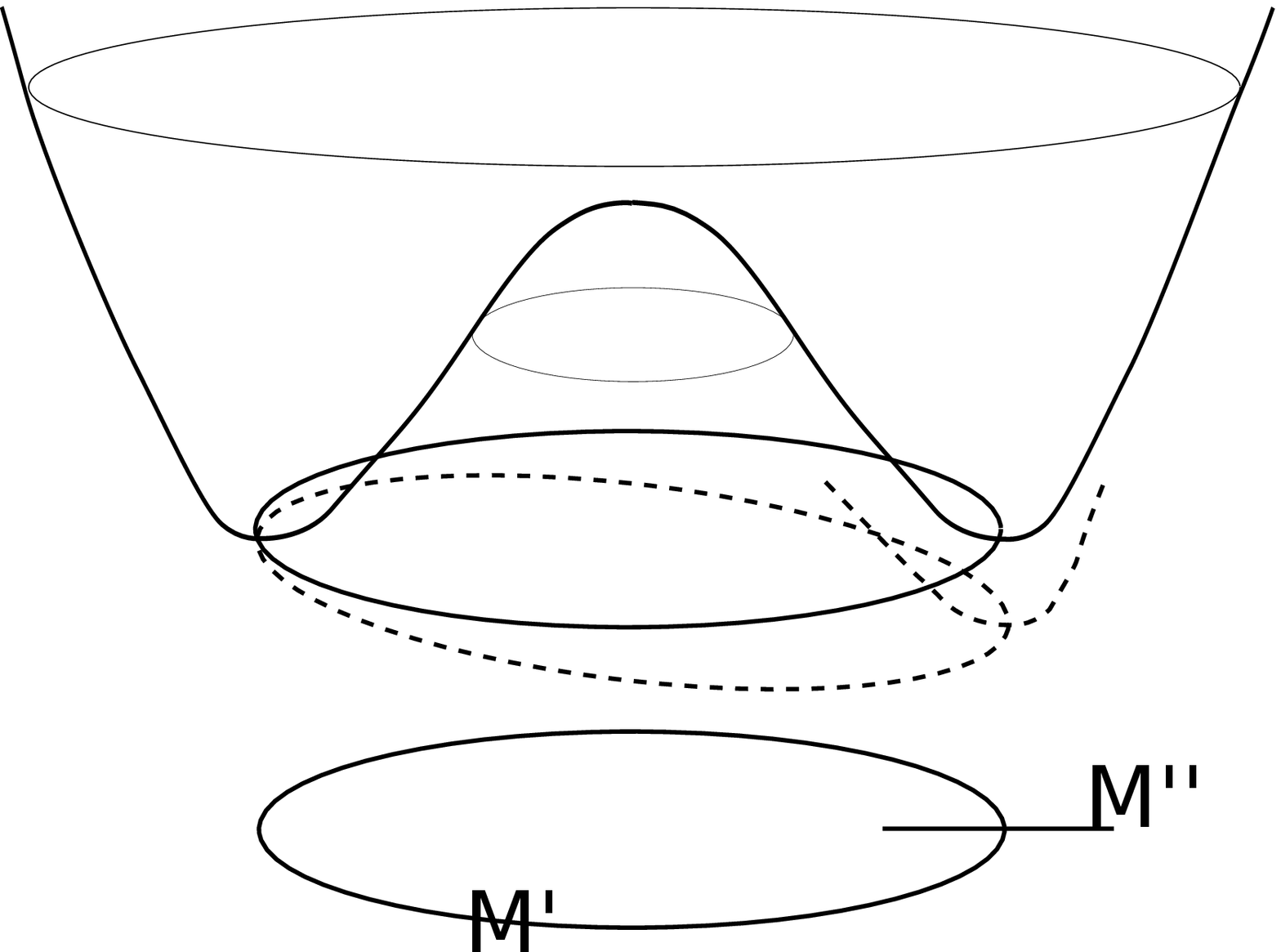}
\captionsetup{labelformat=empty, width=14cm}
\caption{\textbf{Figure 19:} Case of a critical submanifold (plain line) and its perturbation (dashed
  line): The above example is modelled on the mexican hat function
  $\frac{r^{4}}{4}-\frac{r^{2}}{2}$ with the manifold $M'=\sz^{1}$
with the metric $d\theta^{2}$
and $M''\sim \rz$ (around $r=1$) with the metric $dr^{2}$\,. The
function $\Phi(\theta)=-1-\cos(\theta)$ is a negative Morse function
with maximum value $0$ when $\theta=\pi$\,.}  
}
\end{figure}
\begin{prop}
\label{pr:critsubmfld}
Under the above assumptions, the boundary of the unstable cell
$\partial e_{1}^{(p+1)}$ is homologous to $\kappa M'$\,, for some
constant $\kappa$ and relatively to
$f^{-\varepsilon}$\,.\\
For $\delta\geq 0$ small enough, the bar code
$\mathcal{B}_{f_{\delta}}(a,b)$ admits the unique bar
$[0,y_{1}^{(p+1)}[$ of degree $p$ and length $y_{1}^{(p+1)}$\,.\\
The corresponding eigenvalue of
$\Delta_{f_{\delta},f_{\delta}^{-1}([a,b]),h}^{(p)~\text{or}~(p+1)}$
equals
$$
\frac{h}{\pi}\times 
\frac{|\lambda_{1}(\underline{y}_{1}^{(p+1)})\cdots
  \lambda_{p+1}(\underline{y}_{1}^{(p+1)})|^{1/2}}{|\lambda_{p+2}(\underline{y}_{1}^{(p+1)})\cdots
  \lambda_{d}(\underline{y}_{1}^{(p+1)})|^{1/2}}
\times \frac{(\pi h)^{-p}\left(\kappa \int_{M'}e^{\frac{2\delta
        \phi(x')}{h}}~dx'\right)^{2}}{(\pi h)^{-d/2}\int_{M'\times
    M''}e^{-2\frac{f-\delta
      \chi(x'')\phi(x')}{h}}~dx}\times e^{-\frac{2y_{1}^{(p+1)}}{h}}\times (1+\mathcal{O}(h))\,.
$$
\end{prop}
\begin{proof}
The first statement is due to the fact that the bar
$[0,y_{1}^{(p+1)}[$ of degre $p$ provides a non null linear
application from the relative homology vector space
$H_{p+1}(f^{y_{1}^{(p+1)}+\varepsilon};f^{y_{1}^{(p+1)}-\varepsilon})$\,, of which
$e_{1}^{(p+1)}$ is a representant, via the boundary map to
$H_{p}(f^{\varepsilon};f^{-\varepsilon})$\,, of which 
the cycle
$M'$ is a representant. Therefore there exists a constant
$\kappa$ such that $\partial e_{1}^{(p+1)}-\kappa M'$ is a boundary
relatively to $f^{-\varepsilon}$\,. In particular if $\omega$ is a
regular $p$-form in $\ker d_{0,f^{-1}([-\varepsilon,+\infty[),1}$ then 
\begin{equation}
\label{eq:Stokessubmfld}
\int_{\partial e_{1}^{(p+1)}}\omega=\kappa\int_{M'}\omega\,.
\end{equation}
The fact that $[0,y_{1}^{(p+1)}[$ remains the only bars of degree $p$
and length $y_{1}^{(p+1)}$ for $\delta>0$ small enough is a
consequence of the stability theorem (Note that for $\delta>0$\,,
$f_{\delta}$ is a Morse function if $x''\mapsto f(x'')$ has a non degenerate
minimum at $x_{0}''$\,.)\\
Let $\varphi_{1}^{(p)}$ be a global quasimode and $\psi_{1}^{(p+1)}$
be a local quasimode associated with the bar $[0,y_{1}^{(p+1)}[$ and
let us compute the scalar product
$$
\langle \psi_{1}^{(p+1)}\,,\, d_{f_{\delta},h}T_{\delta_{2}}\varphi_{1}^{(p)}\rangle\,.\\
$$
Because we have a non degenerate critical point at
$\underline{y}_{1}^{(p+1)}$\,, the computations of
\cite{LNV}-Section~4.3, which rely on the WKB approximation for
$\psi_{1}^{(p+1)}$ around $\underline{y}_{1}^{(p+1)}$ and 
$d_{f_{\delta},h}\varphi_{1}^{(p)}\equiv 0$ in
$f_{\delta}^{y_{1}^{(p+1)}-\delta(h)}=f^{y_{1}^{(p+1)}-\delta(h)}$\,, leads to 
\begin{multline*}
\langle \psi_{1}^{(p+1)}\,,\,
d_{f_{\delta},h}T_{\delta_{2}}\varphi_{1}^{(p)}\rangle
=\pm \left(\frac{h}{\pi}\right)^{1/2}
\times 
\frac{|\lambda_{1}(\underline{y}_{1}^{(p+1)})\cdots
  \lambda_{p+1}(\underline{y}_{1}^{(p+1)})|^{1/4}}{|\lambda_{p+2}(\underline{y}_{1}^{(p+1)})\cdots
  \lambda_{d}(\underline{y}_{1}^{(p+1)})|^{1/4}}\times(\pi
h)^{\frac{d}{4}-\frac{p}{2}}\\
\times
\int_{\partial e_{1}^{(p+1)}}e^{\frac{f_{\delta}}{h}}\varphi_{1}^{(p)}\times
e^{-\frac{y_{1}^{(p+1)}}{h}}\times (1+\mathcal{O}(h))\,.
\end{multline*}
Because $d(e^{\frac{f_{\delta}}{h}}\varphi_{1}^{(p+1)})\equiv 0$ in
$f^{y_{1}^{(p+1)}-\delta(h)}$ we may apply \eqref{eq:Stokessubmfld}
with $\omega=e^{\frac{f_{\delta}}{h}}\varphi_{1}^{(p+1)}$ and the integral
$\int_{\partial e_{1}^{(p+1)}}$ can be replaced by
$\kappa\int_{M'}$\,. Thus it suffices to know $\varphi_{1}^{(p)}$ in a
neighborhood of $M'$\,. A good approximation is given by a normalized
element of
$\ker(\Delta_{f_{\delta},f_{\delta}^{-1}([-\varepsilon,\varepsilon],h)}^{(p)})$ which is
exponentially close (in any Sobolev norm) to  the $p$-form constructed
by the separation of variables in $M'\times M''$
$$
\frac{1}{\left(\int_{M'\times M''}e^{-2\frac{f-\delta \chi(x')\phi(x')}{h}}~dx\right)^{1/2}}
e^{-\frac{f-\delta \chi(x'')\phi(x')}{h}} |\det g'(x')|^{\frac12} dx_{1}\wedge\ldots\wedge dx_{p}\,.
$$
The final result follows by taking the square.
\end{proof}
When $f(x'')$ near $x_{0}''\in M''$ and $\phi(x')$\,, $x'\in M'$\,,
are Morse functions, the above formula allows again to study the
transition between the case when $f$ is a Morse function on $M$ for
$\delta>0$ small and when $0$ is a degenerate critical value with
critical set $M'$ for $\delta=0$\,. We get the following asymptotic
behaviour for the eigenvalue of
$\Delta_{f_{\delta},f_{\delta}^{-1}([a,b]),h}^{(p)\text{~or~}(p+1)}$
associated with the bar $[0,y_{1}^{(p+1)}[$:
\begin{eqnarray*}
  &
    C_{\delta}\frac{h}{\pi}e^{-2\frac{y_{1}^{(p+1)}}{h}}(1+\mathcal{O}(h))&\text{when}~\delta>0\,,\\
& C_{0}\frac{h}{\pi}(\pi h)^{-p/2}e^{-2\frac{y_{1}^{(p+1)}}{h}}(1+\mathcal{O}(h))&\text{when}~\delta=0\,.
\end{eqnarray*}
In general degre $p$ it is possible to have a good information on the
local approximations of the global quasimodes $\varphi_{k}^{(p)}$
either when the critical value is $x_{k}^{(p)}$ is non degenerate via
a WKB approximation of when we can use some separation of
variables. Otherwise it is not clear that we could get a general
robust integral formula for the subexponential factor. Note also that we used the
fact that $y_{k}^{(p+1)}$ is a non degenerate critical value when we
reduced the computation of $\langle \psi_{1}^{(p+1)}\,,\,
d_{f,h}T_{\delta_{2}}\varphi_{1}^{(p)}\rangle$ to a integral along the
explicit cycle $\partial e_{1}^{(p+1)}$\,. Again it is not clear that
such a simple argument can be used  when $y_{k}^{(p+1)}$ is a degenerate
critical value without some other specific assumptions.
\subsection{More general Morse functions}
\label{sec:appmoregenMorse}
We consider in this paragraph a Morse function $f$ which 
may admit
multiple critical values.
For the sake of simplicity, we work in the following situation:
\begin{itemize}
\item $c<c'$\,, $c,c'\in
\left\{c_{1},\ldots, c_{N_{f}}\right\}$ are the only multiple critical
values.
\item All the critical points with critical value $c$ (resp. $c'$)\,,
  $\underline{x}_{k}^{(p)}$\,, $1\leq k\leq K$\,,
  (resp. $\underline{y}_{k'}^{(p+1)}$\,, $1\leq k'\leq K'$) have the
  index $p$ (resp. $p+1$).
\item  All the bars of $\mathcal{B}_{f}$ with the lower (resp. upper) endpoint $c$ (resp. $c'$) have a length larger or equal to $c'-c$\,. The
  numbers of such bars of length equal to $c'-c$ (the bar is a copy of
  $[c,c'[$)\,, is denoted by 
$K_{0}\leq \min(K,K')$\,.
\begin{figure}[h]
\centering{
\includegraphics[width=12cm]{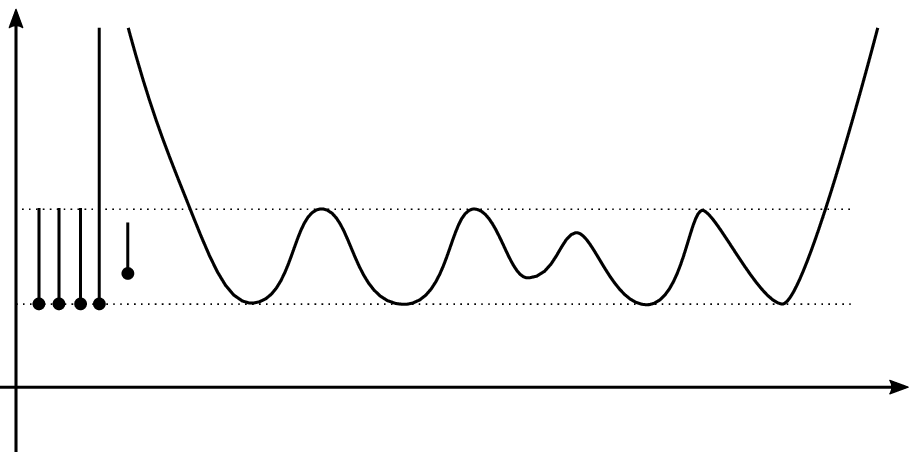}
\captionsetup{labelformat=empty,width=14cm}
\caption{\textbf{Figure 20:} A simple  example in dimension $1$ with $K=4$\,, $K_{0}=K'=3$\,.}  
}
\end{figure}
\item We will consider the energy interval $[a,b]$ such that $c$
  (resp. $c'$) is the smallest (resp. largest) critical value in
  $[a,b]$\,.
\item When $\underline{x}_{k}^{(p)}$\,, $k=1,\ldots, K$
  (resp. $\underline{y}_{k'}^{(p+1)}$\,, $k'=1,\ldots, K'$) denote the
  critical points for the value $c$ (resp. $c'$) the function
  $\chi_{k}^{(p)}\in \mathcal{C}^{\infty}(M;[0,1])$ (resp
  $\chi_{k}^{(p+1)}\in \mathcal{C}^{\infty}(M;[0,1])$) is supported in
  a neighborhood  and equals $1$ in a
  smaller neighborhood of $\underline{x}_{k}^{(p)}$
  (resp.$\underline{y}_{k'}^{(p+1)}$) for $k=1,\ldots,K$ (resp
  $k'=1,\ldots, K'$)\,. Let $t_{k}^{(p)}$\,, $k=1,\ldots, K$\,,
  (resp. $t_{k'}^{(p+1)}$\,, $k'=1,\ldots, K'$)
 be real numbers. 
For $\delta\in \rz$ small, we consider 
$$
f_{\delta}=f+\delta\left[
\sum_{k=1}^{K}t_{k}^{(p)}\chi_{k}^{(p)}+
\sum_{k'=1}^{K'}t_{k'}^{(p+1)}\chi_{k'}^{(p+1)}\right]\,.
$$
\end{itemize}
Because $f$ is a Morse function we may find $\varepsilon>0$ small
enough such that the homology vector space
$H_{p}(f^{c+\varepsilon},f^{c-\varepsilon};\rz)$ (resp. 
$H_{p+1}(f^{c'+\varepsilon},f^{c'-\varepsilon};\rz)$)
have a basis made of the descending (unstable of $-\nabla f$) 
manifolds $e_{k}^{(p)}$\,, $1\leq
k\leq K$ (resp. $e_{k'}^{(p)}$\,, $1\leq k'\leq K'$) restricted to
$f_{c-\varepsilon}$ (resp. $f_{c'-\varepsilon}$). The boundary of
$e_{k}^{(p)}$ (resp $e_{k'}^{(p+1)}$) is a $p-1$-dimensional (resp. $p$-dimensional) sphere
$\partial e_{k}^{(p)}$ (resp. $\partial e_{k'}^{(p+1)}$) lying in
$f^{-1}(\{c-\varepsilon\})$ (resp. in
$f^{-1}(\{c'-\varepsilon\})$)\,.\\
On the Witten Laplacian side,
$\ker(\Delta_{f,f^{-1}([c-\varepsilon,c+\varepsilon]),h}^{(p)})$
(resp. $\ker(\Delta_{f,f^{-1}([c'-\varepsilon,c'+\varepsilon],h)}^{(p+1)})$) may
be  approximated with a $\tilde{O}(e^{-\frac{\varepsilon}{h}})$-distance
by 
$\mathop{\oplus}_{1\leq k\leq K}^{\perp} \cz \psi_{k}^{(p)}$ 
(resp. $\mathop{\oplus}_{1\leq k'\leq K'}^{\perp} \cz\psi_{k'}^{(p+1)}$), where
 $\psi_{k}^{(p)}$ (resp. $\psi_{k'}^{(p+1)}$) is a normalized ground state of 
$\Delta_{f,k}^{(p)}$
(resp. $\Delta_{f,k'}^{(p+1)}$), 
the Witten Laplacian in degree $p$ (resp. $p+1$) with full Dirichlet boundary
conditions in $B(\underline{x}_{k}^{(p)},R\sqrt{\varepsilon})$
(resp. $B(\underline{y}_{k'}^{(p)},R\sqrt{\varepsilon})$) for $R>0$ chosen large enough\,.
We refer to \cite{Hel} and \cite{HeSj4} and we recall that for the Witten Laplacian associated with a Morse function $f$\,, the local Agmon distance
to a critical point $s$\,, $\phi$ solving $|\nabla \phi|^{2}=|\nabla f|^{2}$ and satisfying $\phi(x)\geq |f(x)-f(s)|$\,,
behaves like the square of the geodesic distance to $s$\,.
Additionally, the $L^{2}$ estimate between $\psi_{k}^{(p)}$ (resp. $\psi_{k'}^{(p+1)}$) 
and its projection onto
$\ker\Delta_{f,f^{-1}([c-\varepsilon,c+\varepsilon],h)}^{(p)}$
(resp. $\ker
\Delta_{f,f^{-1}([c'-\varepsilon,c'+\varepsilon],h)}^{(p+1)}$)
can be
completed by a $\tilde{O}(e^{-\frac{\varepsilon}{4h}})$ error estimate in any Sobolev norm
on the open set $f_{c-\frac\varepsilon2}^{c+\frac\varepsilon 2}\cap B(\underline{x}_{k}^{(p)},\frac R2\sqrt{\varepsilon})$
(resp. $f_{c'-\frac\varepsilon2}^{c'+\frac\varepsilon 2}\cap B(\underline{y}_{k'}^{(p+1)},\frac R2\sqrt{\varepsilon})$).\\
We also have WKB-approximations for all the $\psi_{k}^{(p)}$ (resp. $\psi_{k'}^{(p+1)}$)
$1\leq k\leq K$ (resp. $1\leq k'\leq K'$) in
$B(\underline{x}_{k}^{(p)},\frac R2\sqrt{\varepsilon)}$
(resp. $B(\underline{y}_{k'}^{(p+1)},\frac R2\sqrt{\varepsilon})$) which
are valid in $W^{s,2}$-norm.\\
By the construction of Theorem~\ref{th:induc} there is a
$\tilde{O}(e^{-\frac{\varepsilon}{h}})$-orthonormal family of quasimodes
$\varphi_{k}^{(p)}$\,, $1\leq k\leq p$\,, which are
approximated by the $\Pi_{\ker
  (\Delta_{f,f^{-1}([c-\varepsilon,c+\varepsilon]),h}^{(p)})}\psi_{k}^{(p)}$ and therefore
by $\psi_{k}^{(p)}$ or their WKB-approximation and which solve 
$d_{f,h}\varphi_{k}^{(p)}=0$ in $f^{-1}([c-\varepsilon,
c'-\frac{\varepsilon}{2}])$\,, vanish in $f^{c-\varepsilon}$ and
satisfy the exponential decay property.\\
At the level $c'$ the local quasimodes are
$\Pi_{\ker(\Delta_{f,f^{-1}([c'-\varepsilon,c'+\varepsilon]),h}^{(p+1)})}\psi_{k'}^{(p+1)}$
and are therefore close to $\psi_{k'}^{(p+1)}$\,.\\
For a generic choice of the coefficients $t_{k}^{(p)}$ and $t_{k'}^{(p+1)}$\,,
the perturbation $f_{\delta}$ is a  Morse function with simple
critical values  as soon as $\delta\in \rz$ is chosen small enough. Moreover
the stability theorem says that the bars with endpoints $c$ and $c' $
are simply modified by $\mathcal{O}(\delta)$ variations of the
endpoints while all the other bars are not changed owing to our choice
of $f_{\delta}$\,. We can even be more specific. The above parameter
$\varepsilon>0$\,, $R$  being fixed, $\varepsilon$ small enough,
we may take the cut-off function $\chi_{k}^{(p)}$\,, $k=1,\ldots,K$\,,
(resp. $\chi_{k'}^{(p+1)}$\,, $k'=1,\ldots,K'$) such that the equal $1$
in $B(\underline{x}_{k}^{(p)},2R\sqrt{\varepsilon})$ (resp
$B(\underline{y}_{k'}^{(p+1)},2R\sqrt{\varepsilon})$)\,. Finally
$\delta>0$ is chosen small enough such that all the critical values of
$f_{\delta}$ close to $c$ (resp $c'$) are in
$[c-\varepsilon/2,c+\varepsilon/2]$
(resp. $[c'-\varepsilon/2,c'+\varepsilon/2]$)\,. With this choice of
$f_{\delta}$\,, $(e_{k}^{(p)})_{k=1,\ldots,K}$
(resp. $(e_{k'}^{(p+1)})_{k'=1,\ldots, K'}$) defines a basis of
$H_{p}((f_{\delta})^{c+\varepsilon},(f_{\delta})^{c-\varepsilon};\rz)$
(resp. $H_{p+1}((f_{\delta})^{c'+\varepsilon};f_{\delta}^{c'-\varepsilon})$)\,. The
quasimodes $\psi_{k}^{(p)}$\,, $\psi_{k'}^{(p+1)}$\,, and their WKB-approximations are not changed because we have just changed $f$ by a
constant in $B(\underline{x}_{k}^{(p)},R\sqrt{\varepsilon})$  (resp. $B(\underline{y}_{k+1}^{(p+1)},R\sqrt{\varepsilon})$)\,.
\begin{lem}
In the above framework and for $\delta\in\rz$ small enough the
boundary map $\partial:
H_{p+1}((f_{\delta})^{c'+\varepsilon},(f_{\delta})^{c'-\varepsilon};\rz)\overset{\text{can.}}{\sim}
\mathop{\oplus}_{k'=1}^{K'}\rz e_{k'}^{(p+1)}$ induces a linear map to
$H_{p}((f_{\delta})^{c+\varepsilon},(f_{\delta})^{c-\varepsilon};\rz)\overset{\text{can.}}{\sim}\mathop{\oplus}_{k=1}^{K}\rz
e_{k}^{(p)}$ of rank $K_{0}$ which is written
$$
\partial:e_{k'}^{(p+1)}\mapsto \sum_{k=1}^{K}\bm\kappa_{k,k'}e_{k}^{(p)}\,.
$$
The matrix $\bm\kappa$ does not depend on $\delta$\,.
\end{lem}
\begin{proof}
When $\delta=0$\,,
 the boundary map sends
  $H_{p+1}(f^{c'+\varepsilon},f^{c'-\varepsilon},\rz)$
  to $H_{p}(f^{c'-\varepsilon},f^{c-\varepsilon};\rz)$ of
  which a dual basis (in cohomology) is indexed by the $K_{0}$
 bars
  $[c,c'[$\,,
  $k=1,\ldots,K_{0}$\,. It suffices to follow the bars to the lower
  endpoint to define a linear map to
  $H_{p}(f^{c+\varepsilon},f^{c-\varepsilon};\rz)$\,. 
 For a general $\delta$ small enough, $f_{\delta}$ differs from $f$
 only by a constant in each ball of radius $R\sqrt \varepsilon$
 around the critical points $x_{k}^{(p)}$\,, $y_{k'}^{(p+1)}$\,. Therefore,
 the gradient vector fields and the Morse models remain unchanged around
 these points. The homotopy 
 becomes trivial by replacing locally the level set $f^{-1}(\{c-\varepsilon\})$ (resp. $f^{-1}(\{c'-\varepsilon\})$)
 by  $f_{\delta}^{-1}(\{c-\varepsilon\})=f^{-1}(\{c-\varepsilon-\delta t_{k}^{(p)}\})$
 (resp. $f_{\delta}^{-1}(\{c'-\varepsilon\})=f^{-1}(\{c'-\varepsilon-\delta t_{k'}^{(p+1)}\})$).
 Hence,
 $(e_{k}^{(p)})_{k\in\{1,\dots,K\}}$ (resp. $(e_{k'}^{(p+1)})_{k'\in\{1,\dots,K'\}}$) appears as a canonical basis of 
 $H_{p}((f_{\delta})^{c+\varepsilon},(f_{\delta})^{c-\varepsilon};\rz)$
 (resp. $H_{p+1}((f_{\delta})^{c'+\varepsilon},(f_{\delta})^{c'-\varepsilon};\rz)$)
 in which the matrix $\bm\kappa$ of the topological linear map $\partial:H_{p+1}((f_{\delta})^{c'+\varepsilon},(f_{\delta})^{c'-\varepsilon};\rz)\to
 H_{p}((f_{\delta})^{c+\varepsilon},(f_{\delta})^{c-\varepsilon};\rz)$ remains unchanged.
\end{proof}

\begin{figure}[h]
\centering{
\includegraphics[width=12cm]{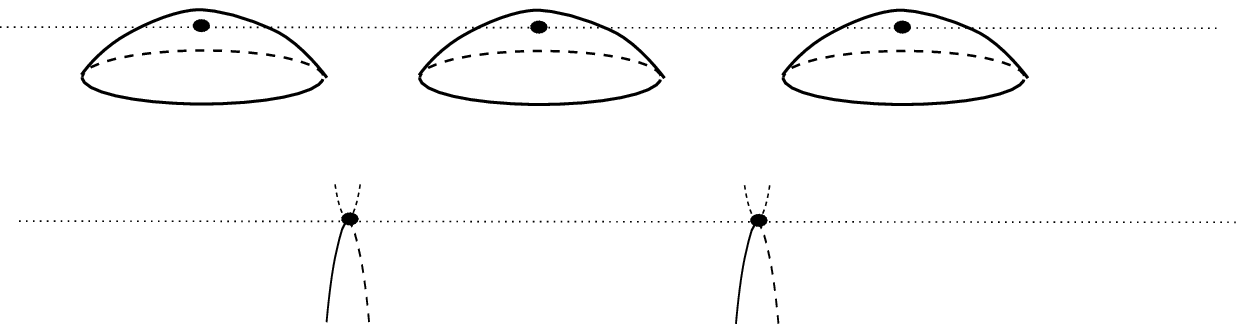}
\captionsetup{labelformat=empty,width=14cm}
\caption{\textbf{Figure 21:} In dimension $2$ we have represented  $3$ critical points
  with index $2$ at the value $c'$ and $2$ critical points with index
  $1$ at the value $c$\,. The unstable (and stable manifold for the
  index $1$) of $-\nabla f$ are considered in the level sets
  $f_{c'-\varepsilon}^{c'+\varepsilon}$ and
  $f_{c-\varepsilon}^{c+\varepsilon}$\,. The homotopy with respect to
  $\delta$ consists simply to move separately up or down, the
  disconnected parts of this picture.
  }  
}
\end{figure}

\begin{prop}
\label{pr:genalMorse}
In the above framework with $\delta$ small enough, 
the singular values $\mu_{h}$ of
  $d_{f_{\delta},(f_{\delta})^{-1}([a,b]),h}^{(p)}$ which satisfy $\lim_{h\to
    0}-h\log\mu_{h}=c-c'+\mathcal{O}(\delta)$ are equal to
  $(1+\mathcal{O}(h))\times$ the non zero
  singular values of the $K\times K'$ matrix
$$
\left(\frac{h}{\pi}\right)^{1/2}
(D^{(p)})^{-1}\bm\kappa
D^{(p+1)}
$$
where $D^{(p)}$  (resp. $D^{(p+1)}$) is the diagonal matrix with entries
\begin{eqnarray*}
&&
\frac{|\lambda_{1}(\underline{x}_{k}^{(p)})\cdots
  \lambda_{p}(\underline{x}_{k}^{(p)})|^{1/4}}{|\lambda_{p+1}(\underline{x}_{k}^{(p)})\cdots
  \lambda_{d}(\underline{x}_{k}^{(p)})|^{1/4}}
e^{-\frac{f_{\delta}(\underline{x}_{k}^{(p)})}{h}}\quad,\quad k=1,\ldots,K\,,\\
\text{resp.}  &&\frac{|\lambda_{1}(\underline{y}_{k'}^{(p+1)})\cdots
  \lambda_{p+1}(\underline{y}_{k'}^{(p+1)})|^{1/4}}{|\lambda_{p+2}(\underline{y}_{k'}^{(p+1)})\cdots
  \lambda_{d}(\underline{y}_{k'}^{(p+1)})|^{1/4}}e^{-\frac{f_{\delta}(\underline{y}_{k'}^{(p+1)})}{h}}\quad,\quad k'=1,\ldots,K'\,.
\end{eqnarray*}
\end{prop}
\begin{proof}
 Set $x_{k,\delta}^{(p)}=f_{\delta}(\underline{x}_{k}^{(p)})=c+\delta
 t_{k}^{(p)}$  and
 $y_{k',\delta}^{(p+1)}=f_{\delta}(\underline{y}_{k'}^{
     (p+1)})=c'+\delta t_{k'}^{(p+1)}$\,, for $k=1,\ldots, K$ and
 $k'=1,\ldots,K'$\,. An orthonormal basis of
 $\ker(\Delta_{f_{\delta},f_{\delta}[c'-\varepsilon,c'+\varepsilon],h}^{(p+1)})$
 is well approximated by the local quasimodes $\psi_{k'}^{(p+1)}$
 which is the ground state of the full Dirichlet
 realization of $\Delta_{f,h}^{(p+1)}$ in
 $B(\underline{y}_{k'}^{(p)},R\sqrt{\varepsilon})$ which do not depend on $\delta$\,. The same holds for
 $\ker(\Delta_{f_{\delta},f_{\delta}([c-\varepsilon,c+\varepsilon],h)}^{(p)})$
 with the notation $\psi_{k}^{(p)}$\,, $k=1,\ldots, p$\,. Hence
 $\mathop{\oplus}_{k=1\ldots,K}^{\perp}\cz \psi_{k}^{(p)}$ provides a
 good approximation in the energy interval
 $[c-\varepsilon,c+\varepsilon]$ for $f_{\delta}$ for the vector space
 of global quasimodes $\varphi_{k,\delta}^{(p)}$ for $f_{\delta}$ associated
 with the bars $[x_{k,\delta}^{(p)},y_{k,\delta}^{(p+1)}[$ for
 $k=1,\ldots, K_{0}$ and $[x_{k,\delta}^{(p)},b[$ for
 $k=K_{0}+1\,,\ldots, K$\,. Let us chose the basis
 $(\varphi_{k,\delta}^{(p)})_{k=1,\ldots,K}$ as an orthonormal basis
 such that
 $\|\varphi_{k,\delta}^{(p)}-\psi_{k}^{(p)}\|_{L^{2}}=\tilde{o}(1)$\,,
 while such a $\tilde{o}(1)$ estimate also holds in any Sobolev norm
 in $f_{c-\frac\varepsilon2}^{c+\frac\varepsilon 2}\cap B(\underline{x}_{k,\delta}^{(p)},\frac R2\sqrt{\varepsilon})$\,. Those global quasimodes are assumed to solve
 $d_{f_{\delta},h}\varphi_{k,\delta}^{(p)}=0$ in $f_{\delta}^{-1}([a,c'-M\delta])$
 for some $M>0$ large enough and we assume $M\delta<<\delta_{2}<<\varepsilon$\,. 
 We now compute the interaction $K'\times K$ matrix $\langle
 \psi_{k'}^{(p)}\,,\,
 d_{f,h}\chi_{\delta_{2}}(f_{\delta})\varphi_{k,\delta}^{(p)}\rangle$ where
 $\chi_{\delta_{2}}$ smoothly vanishes in $[c'-\delta_{2},b]$ and
 equals $1$ in $[a,c'-2\delta_{2}]$ for all $k=1,\ldots,K$\,. Because
 $d_{f,h}\varphi_{k,\delta}^{(p)}=0$ in
 $f_{\delta}^{-1}([b,c'-M\delta])$\,, the local computation around
 $y_{k'}^{(p+1)}$ done in \cite{LNV}-Section~4 are the same  and
 they say:
\begin{multline*}
\langle \psi_{k'}^{(p+1)}\,,\,
d_{f,h}T_{\delta_{2}}\varphi_{k,\delta}^{(p)}\rangle
=\pm \left(\frac{h}{\pi}\right)^{1/2}
\times 
\frac{|\lambda_{1}(\underline{y}_{k'}^{(p+1)})\cdots
  \lambda_{p+1}(\underline{y}_{k'}^{(p+1)})|^{1/4}}{|\lambda_{p+2}(\underline{y}_{k'}^{(p+1)})\cdots
  \lambda_{d}(\underline{y}_{k'}^{(p+1)})|^{1/4}}\times(\pi h)^{\frac{d}{4}-\frac{p}{2}}\\
\times \int_{\partial e_{k'}^{(p+1)}}e^{\frac{f_{\delta}}{h}}\varphi_{k,\delta}^{(p)}\times
e^{-\frac{y_{k',\delta}^{(p+1)}}{h}}\times (1+\mathcal{O}(h))\,.
\end{multline*}
By Stokes's formula applied with
$d[e^{\frac{f_{\delta}}{h}}\varphi_{k,\delta}^{(p)}]=0$ in
$f_{\delta}^{-1}([b,c'-M\delta])$ we obtain
\begin{multline*}
\langle \psi_{k'}^{(p+1)}\,,\,
d_{f,h}T_{\delta_{2}}\varphi_{k,\delta}^{(p)}\rangle
=\pm \left(\frac{h}{\pi}\right)^{1/2}
\times 
\frac{|\lambda_{1}(\underline{y}_{k'}^{(p+1)})\cdots
  \lambda_{p+1}(\underline{y}_{k'}^{(p+1)})|^{1/4}}{|\lambda_{p+2}(\underline{y}_{k'}^{(p+1)})\cdots
  \lambda_{d}(\underline{y}_{k'}^{(p+1)})|^{1/4}}\times(\pi
h)^{\frac{d}{4}-\frac{p}{2}}\\
\times\left[\sum_{j=1}^{K}\bm\kappa_{j,k'} \int_{e_{j}^{(p)}}e^{\frac{f_{\delta}}{h}}\varphi_{k,\delta}^{(p)}\right]\times
e^{-\frac{y_{k',\delta}^{(p+1)}}{h}}\times (1+\mathcal{O}(h))\,.
\end{multline*}
By approximating  $\varphi_{k,\delta}^{(p)}$ by
$\psi_{k,\delta}^{(p)}$ and its WKB approximation in
$B(\underline{x}_{k}^{(p)},\frac R2 \sqrt{\varepsilon})$ we have
\begin{eqnarray*}
(\pi h)^{\frac{d}{4}-\frac{p}{2}}\int_{e_{j}^{(p)}}e^{\frac{f_{\delta}}{h}}\varphi_{k,\delta}^{(p)}
&=&(\pi
h)^{\frac{d}{4}-\frac{p}{2}}\int_{e_{j}^{(p)}}e^{\frac{f_{\delta}}{h}}\psi_{k}^{(p)}\times
(1+\tilde{o}(1))\\
&=&\pm 1\frac{|\lambda_{p+1}(\underline{x}_{k}^{(p)})\ldots
  \lambda_{d}(\underline{x}_{k}^{(p)})|^{1/4}}{|\lambda_{1}(\underline{x}_{k}^{(p)})\ldots
  \lambda_{p}(\underline{x}_{k}^{(p)})|^{1/4}}e^{\frac{x_{k,\delta}^{(p)}}{h}}\times(1+\mathcal{O}(h))\,.
\end{eqnarray*}
The error terms  actually occur as matricial products on the left-hand side
for the approximation of $\psi_{k'}^{(p+1)}$ and on the right-hand
side for $\varphi_{k,\delta}^{(p)}$\,.\\
 The interaction matrix $\langle
(\psi_{k'}^{(p+1)}\,,\,
d_{f_{\delta},h}\chi_{\delta_{2}}\varphi_{k,\delta}^{(p)}\rangle)_{1\leq
k'\leq K', 1\leq k \leq K}$ is thus equal to 
$$
\textrm{diag}\big(\pm
1+\mathcal{O}(h)\big)\left(\frac{h}{\pi}\right)^{1/2}D^{(p+1)}({}^{t}\bm\kappa)
(D^{(p)})^{-1}\textrm{diag}\big(\pm 1+\mathcal{O}(h)\big)\,.
$$
Its singular values are thus equal up to a $\mathcal{O}(h)$-relative
error to the singular values of 
$$\left(\frac{h}{\pi}\right)^{1/2}
D^{(p+1)}({}^{t}\bm\kappa) (D^{(p)})^{-1}
$$
 or equivalently
of 
$$\left(\frac{h}{\pi}\right)^{1/2}(D^{(p)})^{-1}\bm\kappa D^{(p+1)}\,.
$$
\end{proof}
\begin{remark}
  The result of Propostion~\ref{pr:genalMorse}, in a specific case,
  show that it is possible to get a matricial robust accurate formula
  for the exponentially small eigenvalues of Witten Laplacians for
  general Morse potentials. This provides another stability property
  valid for the first term asymptotics of the subexponential factor, which allows
  to study the transition from the generic Morse case with simple
  critical values to the general case. Note that here the power of $h$
  in this subexponential factor is not changed. But discontinuities appear on the
  constants as it is shown in the next simple examples. Actually we
  have considered a simple case where only one multiple bar $[c,c'[$
  has to be taken into account.  A more general form would consist in
  following the induction scheme of Theorem~\ref{th:induc} and would
  lead to some complicated linear matricial structure for which we do
  not have an elegant presentation at the moment. In the degree
  $p=0$\,, L.~Michel in \cite{Mic} proposed an interpretration in terms
  of the spectrum of a discrete Laplacian on a finite graph with
  vertices given by the local minima and edges given by saddle
  points. This formulation is written for a fixed Morse function with
  possible multiple local minima and saddle points, the perturbative
  issue is not really clarified there. In our specific example, the discrete Laplacian proposed by
  L.~Michel is actually the square 
$$
\frac{h}{\pi} (D^{(0)})^{-1}\bm\kappa D^{(1)}D^{(1),*}\bm\kappa^{*} (D^{(0)})^{-1,*}\,.
$$
 It would be interesting to find such a general robust formulation,
 with several multiple critical values,
 in degree $p>0$\,.
\end{remark}

\textbf{Examples:}
\begin{enumerate}
\item Consider a Morse function $f$
  on $[s,t]$ such that $\min_{x\in[s,t]}f(x)=f(s)=a$\,,
  $\max_{x\in[s,t]}f(x)=f(t)=b$\,, with $f'(s)>0$ and $f'(t)>0$\,,
  with two local  maxima and two local minima $s<\underline{y}_{1}^{(1)}<\underline{x}_{1}^{(0)}<\underline{y}_{2}^{(1)}<\underline{x}_{2}^{(0)}<t$\,,
  $f(\underline{y}_{1}^{(1)})=f(\underline{y}_{2}^{(1)})=c'$ and 
$f(\underline{x}_{1}^{(0)})=f(\underline{x}_{2}^{(0)})=c$\,. For the
perturbation of $f$ we will consider the cases when
$(t_{1}^{(1)},t_{2}^{(1)})=(0,0)$\,,
$(t_{1}^{(0)},t_{2}^{(0)})\in\left\{(0,0),(0,-1),(-1,0)\right\}$\,.
The matrix $\bm\kappa$ equals
$$
\begin{pmatrix}
  1&-1\\0&1
\end{pmatrix}
$$ 
while the matrices $D^{(0)}$ and $D^{(1)}$ are given by
\begin{eqnarray*}
D^{(0)}&=&
\begin{pmatrix}
  |\lambda_{1}(\underline{x}_{1}^{(0)})|^{-1/4}e^{-\frac{c+\delta
      t_{1}^{(0)}}{h}}&0\\
0& |\lambda_{1}(\underline{x}_{2}^{(0)})|^{-1/4}e^{-\frac{c+\delta
      t_{2}^{(0)}}{h}}
\end{pmatrix}
=
\begin{pmatrix}
  \alpha_{1}^{-1}&0\\
0&\alpha_{2}^{-1}
\end{pmatrix}
e^{-\frac{c}{h}}\,,
\\
D^{(1)}&=&
\begin{pmatrix}
  |\lambda_{1}(\underline{y}_{1}^{(0)})|^{1/4}e^{-\frac{c'}{h}}&0\\
0& |\lambda_{2}(\underline{y}_{2}^{(0)})|^{1/4}e^{-\frac{c'}{h}}
\end{pmatrix}
=
\begin{pmatrix}
  \beta_{1}&0\\
0&\beta_{2}
\end{pmatrix}
e^{-\frac{c'}{h}}\,.
\end{eqnarray*}
The singular values of the matrix $(D^{(0)})^{-1}\bm\kappa D^{(1)}$ are
the square roots of the eigenvalues of the symmetric square matrix
$$
\begin{pmatrix}
  \alpha_{1}^{2}(\beta_{1}^{2}+\beta_{2}^{2})&
  -\alpha_{1}\alpha_{2}\beta_{2}^{2}\\
  -\alpha_{1}\alpha_{2}\beta_{2}^{2}& \alpha_{2}^{2}\beta_{2}^{2}
\end{pmatrix}e^{-2\frac{c'-c}{h}}\,. 
$$
Those eigenvalues equal
$$
\frac{[(\alpha_{1}^{2}(\beta_{1}^{2}+\beta_{2}^{2})+\alpha_{2}^{2}\beta_{2}^{2})]\pm
\sqrt{[\alpha_{1}^{2}(\beta_{1}^{2}+\beta_{2}^{2})-\alpha_{2}^{2}\beta_{2}^{2})]^{2}+4\alpha_{1}^{2}\alpha_{2}^{2}\beta_{2}^{4}}}{2}\times e^{-2\frac{c'-c}{h}}\,.
$$

Depending on the three considered cases, we obtain:

\begin{description}
\item[$t_{1}^{(0)}=t_{2}^{(0)}=0$:] The $2$ exponentially small
  eigenvalues of $\Delta_{f,[s,t],h}^{(0)\text{~or~}(1)}$ have the
  form $C_{k}\frac{h}{\pi}e^{-2\frac{c'-c}{h}}(1+\mathcal{O}(h))$\,, $k=1,2$\,, where the constants
  $C_{1}$ and $C_{2}$ clearly depend on the four hessians at the
  critical points.
\item[$t_{1}^{(0)}=-1, t_{2}^{(1)}=0$:] The $2$ exponentially
  small eigenvalues of    $\Delta_{f,[s,t],h}^{(0)\text{~or~}(1)}$
  are equal to: 
  \begin{eqnarray*}
    && \frac{h}{\pi}|\lambda_{1}(\underline{x}_{2}^{(0)})|^{1/2}|\lambda_{1}(\underline{y}_{2}^{(1)})|^{1/2}
e^{-2\frac{c'-c}{h}} (1+\mathcal{O}(h))\\
&& \frac{h}{\pi}|\lambda_{1}(\underline{x}_{1}^{(0)})|^{1/2}|\lambda_{1}(\underline{y}_{1}^{(1)})|^{1/2}
e^{-2\frac{c'-c+\delta}{h}}(1+\mathcal{O}(h))\,.
  \end{eqnarray*}
In particular the smallest one depends on the hessians of
$f_{\delta}$ at the only points $\underline{x}_{1}^{(0)}$ and
$\underline{y}_{1}^{(1)}$\,.
\item[$t_{1}^{(0)}=0\,,\, t_{2}^{(0)}=-1$:] 
The $2$ exponentially
  small eigenvalues of    $\Delta_{f,[s,t],h}^{(0)\text{~or~}(1)}$
  are equal to:
  \begin{eqnarray*}
    && 2 \frac{h}{\pi}|\lambda_{1}(\underline{x}_{1}^{(0)})|^{1/2}\frac{|\lambda_{1}(\underline{y}_{1}^{(1)})|^{1/2}+|\lambda_{1}(\underline{y}_{1}^{(1)})|^{1/2}}{2}e^{-2\frac{c'-c}{h}} (1+\mathcal{O}(h))\\
&&\frac12  \frac{h}{\pi}|\lambda_{1}(\underline{x}_{2}^{(0)})|^{1/2}H(|\lambda_{1}(\underline{y}_{1}^{(1)})|^{1/2},|\lambda_{1}(\underline{y}_{2}^{(1)})|^{1/2})e^{-2\frac{c'-c+\delta}{h}}(1+\mathcal{O}(h))\,,
  \end{eqnarray*}
where $H(u,v)=\frac{2uv}{u+v}$ denotes the harmonic mean.\\
In this case the smallest eigenvalue depends on the Hessians of
$f_{\delta}$ at the points $\underline{x}_{2}^{(0)}$\,,
$\underline{y}_{1}^{(1)}$ and $\underline{y}_{2}^{(1)}$\,.
\end{description}

\begin{figure}[h]
\centering{
\includegraphics[width=10cm]{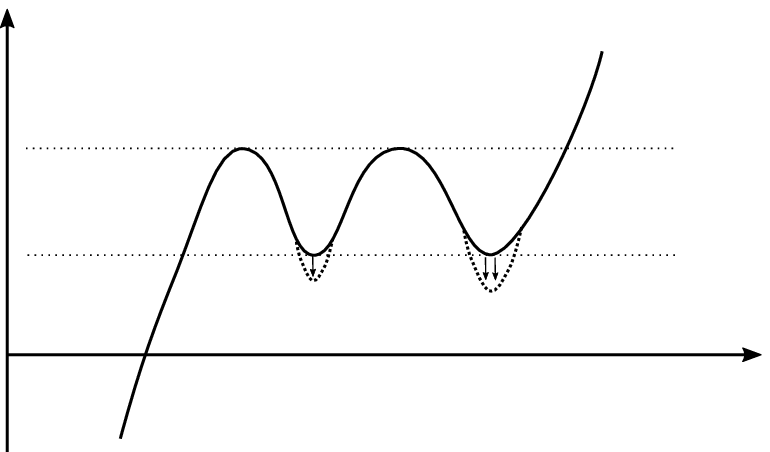}
\captionsetup{labelformat=empty,width=14cm}
\caption{\textbf{Figure 22:} The three considered cases: $(t_{1},t_{2})=(0,0)$ plain line;
  $(t_{1},t_{2})=(-1,0)$ move the curve downward with $(\downarrow)$\,,
  $(t_{1},t_{2})=(0,-1)$ move the curve downward with $(\downarrow~\downarrow)$\,.
  }  
}
\end{figure}

The general formula is again a robust formula which allow to follow
the dependence on the parameter $\delta$ of the asymptotic expressions
although those at the end  are not continuous with respect to
$\delta$\,.
\item Consider in $\rz^{d}$ a function $f$ with a unique
  local minimum at $x_{1}^{(0)}=0$ with $f(0)=c$\,, such that
  $\lim_{x\to\infty}f(x)=-\infty$ and surrounded by $K'$ saddle
  points, critical points with index $1$\,, such that
  $f(\underline{y}_{k'}^{(1)})=c'$\,, while all the other crtical
  values are larger than $c'$ with an index  $p\geq 2$\,.
For the perturbation we will consider the two cases when $t_{1}^{(0)}=0$
and $t_{1}^{(1)}\in \left\{0,-1\right\}$\,. The matrix $\bm\kappa$ is the
$1\times K'$ matrix
$(1\,,\,1\,,\,\ldots\,,\,1)\,.$
Thus the smallest eigenvalue of
$\Delta_{f_{\delta},(f_{\delta})^{-1}([a,b]),h}^{(0)}$\,, which is the
unique exponentially small eigenvalue, equals

\begin{eqnarray*}
  &\frac{h}{\pi}|\det(\mathrm{Hess}f(\underline{x}_{1}^{(0)}))|^{1/2}\sum_{k'=1}^{K'}\frac{|\lambda_{1}(\underline{y}_{k'}^{(1)})|^{1/2}}{|\lambda_{2}(\underline{y}_{k'}^{(1)})\ldots\lambda_{d}(y_{k'}^{(1)})|^{1/2}}e^{-2\frac{c'-c}{h}}(1+\mathcal{O}(h))&\text{if}~\delta=0\,,
\\
  &\frac{h}{\pi}|\det(\mathrm{Hess}f(\underline{x}_{1}^{(0)}))|^{1/2}\frac{|\lambda_{1}(\underline{y}_{1}^{(1)})|^{1/2}}{|\lambda_{2}(\underline{y}_{1}^{(1)})\ldots\lambda_{d}(y_{1}^{(1)})|^{1/2}}e^{-2\frac{c'-c-\delta}{h}}(1+\mathcal{O}(h))&\text{if}~\delta>0\,.
\end{eqnarray*}

Similar formulas are  obtained for various configurations in \cite{DLLN1,DLLN2,LeNe2,LeMi}.

\begin{figure}[h]
\centering{
\includegraphics[width=6cm]{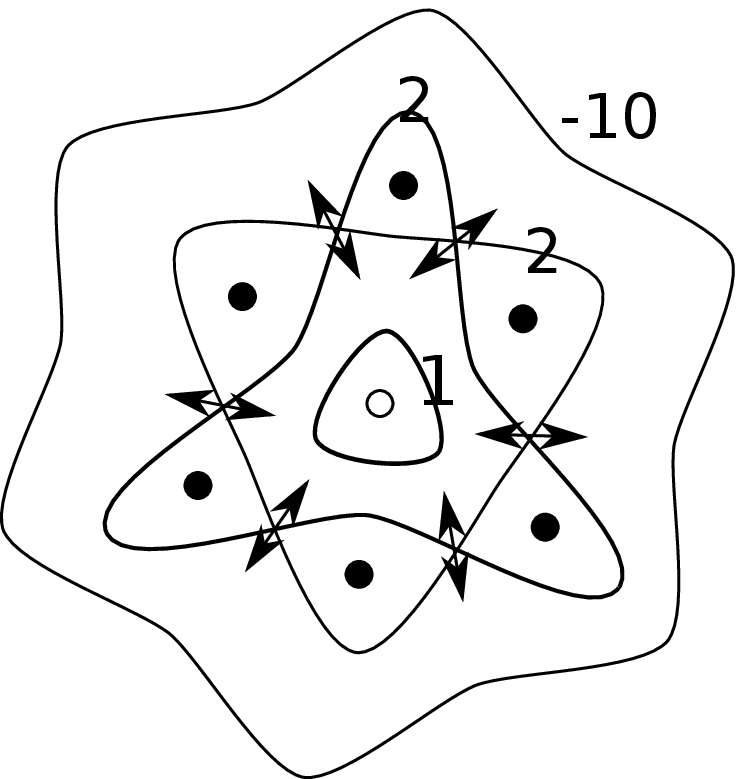}
\captionsetup{labelformat=empty,width=14cm}
\caption{\textbf{Figure 23:} An example with $K=6$\,. Level curves at
  the level $1,2,-10$ are represented, the global minimum is denoted
  by $\circ$\,, the saddle points by $\leftrightarrow$ and the local
  maximum by $\bullet$\,.
  }  
}
\end{figure}

\item \textbf{A case with symmetries:} Consider in $\rz^{2}$ a Morse
  function $f$ with a local maximum at $\underline{y}_{1}^{(2)}=0$\,, $f(\underline{y}_{1}^{(2)})=c_{2}$
  surounded by $K$ saddle points at
  $\underline{x}_{k}^{(1)}=\underline{y}_{k}^{(1)}$\,, $k=1\ldots
  K$\,, $f(\underline{x}_{k}^{(1)})=c_{1}$\,, and $K$ local minima at
  $\underline{x}_{k}^{(0)}$\,, $k=1\ldots, K$\,,
  $f(\underline{x}_{k}^{(0)})=c_{0}$\,, $c_{0}<c_{1}<c_{2}$\,. We also
  assume that $\lim_{x\to\infty}f(x)=+\infty$ and that $f$ has no
  other critical points. 
 When $j\in
\left\{1,2\right\}$ or $p\in \left\{0,1\right\}$ are fixed
$\lambda_{j}(\underline{x}_{k}^{(p)})=\lambda_{j}^{(p)}$ do not depend
on $k=1,\ldots,K$\,. We study the eigenvalues of
$\Delta_{f,\rz^{2},h}^{(p)}$\,, $p=0,1,2$ by restricting to the case
$c_{0}<a<c=c_{1}<c'=c_{2}<b$ for $p=2$ and to the case
$a<c=c_{0}<c'=c_{1}<b< c_{2}$ for $p=0$\,. By supersymmetry, the non
zero eigenvalues of $\Delta_{f,\rz^{2},h}^{(1)}$ are obtained by
gathering the ones of $\Delta_{f,\rz^{2},h}^{(0)}$ and
$\Delta_{f,\rz^{2},h}^{(2)}$\,. 

\begin{figure}[h]
\centering{
\includegraphics[width=6cm]{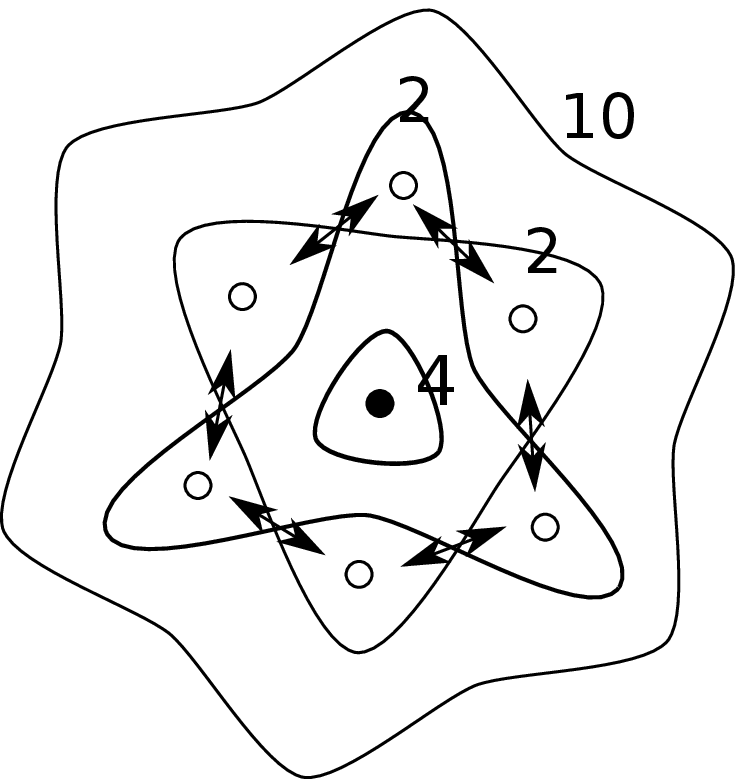}
\captionsetup{labelformat=empty, width=14cm}
\caption{\textbf{Figure 24:} An example with $K=6$\,. Level curves at
  the level $4,2,10$ are represented, the local minima are denoted
  by $\circ$\,, the saddle points by $\leftrightarrow$ and the global
  maxima by $\bullet$\,.
  }  
}
\end{figure}
\begin{description}
\item[For $p=2$\,, $c_{0}<a<c=c_{1}<c'=c_{2}<b$:] The matrix $\bm\kappa$ equals
  the $K\times 1$ matrix
$$
\bm\kappa=
\begin{pmatrix}
  1\\\vdots\\1
\end{pmatrix}\,.
$$
The smallest  eigenvalue of $\Delta_{f,f^{-1}([a,b]),h}^{(2)}$\,,
which is the only exponentially small one, then equals:
$$
\frac{h}{\pi}|\det(\mathrm{Hess}f(\underline{y}_{1}^{(2)}))|^{1/2}K\frac{|\lambda_{2}^{(1)}|^{1/2}}{|\lambda_{1}^{(1)}|^{1/2}}e^{-2\frac{c_{2}-c_{1}}{h}}(1+\mathcal{O}(h))\,.
$$
\item[For $p=0$\,, $a<c=c_{0}<c'=c_{1}<b<c_{2}$:] The matrix $\bm\kappa$
  is the $K\times K$ matrix
$$
\bm\kappa=
\begin{pmatrix}
  1&-1&0&\ldots&0\\
0&1&-1&\ddots&\vdots\\
\vdots&\ddots&\ddots&\ddots&\vdots\\
\vdots&\vdots&0&\ddots&-1\\
-1&0&\ldots&0&1
\end{pmatrix}
$$
 of which the singular values equal $|1-\omega^{k}|$\,, $k=1,\ldots,
 K$\,, where $\omega^{k}=e^{2i\pi\frac{k}{K}}$ for $k=1,\dots,K$\,. 
Owing to 
$$
(D^{(0)})^{-1}\bm\kappa
D^{(1)}=|\lambda_{1}^{(0)}\lambda_{2}^{(0)}|^{1/4}\frac{|\lambda_{1}^{(1)}|^{1/4}}{|\lambda_{2}^{(1)}|^{1/4}}e^{-\frac{c_{1}-c_{0}}{h}}\bm\kappa
\,,
$$
we deduce that the $K$ exponentially small eigenvalues of
$\Delta_{f,f^{-1}([a,b]),h}^{(0)}$ are equal to 
$$
|\lambda_{1}^{(0)}\lambda_{2}^{(0)}|^{1/2}\frac{|\lambda_{1}^{(1)}|^{1/2}}{|\lambda_{2}^{(1)}|^{1/2}}|1-\omega^{k}|^{2}e^{-2\frac{c_{1}-c_{0}}{h}}(1+\mathcal{O}(h))\,,\quad
k=1,\ldots, K\,.
$$
\end{description}
\end{enumerate}
This case with $p=0$ was considered by Michel in \cite{Mic} for the Witten Laplacian and by H\'erau-Hitrik-Sj\"ostrand
in \cite{HHS} for the non-self-adjoint Kramers-Fokker-Planck operator.

\section{Broadening the scope}
\label{sec:persp}
Our work provides a general method for analyzing the
exponentially small eigenvalues of Witten Laplacians with a general
potential function. However, it does not answer all the
questions that arose along this analysis.  Here is a short list of
still open questions and of connections with  closely related fields.
\\
\noindent\textbf{a) General $\mathcal{C^{\infty}}$ potential:} A general $\mathcal{C}^{\infty}$-function on a compact manifold
  $M$
may have an infinite
  number of critical values and bars in its bar code. Nevertheless, 
for any $\varepsilon>0$, the set of bars of length larger than
$\varepsilon$ is finite. In order to realize this, take a covering
 $[\min f,\max f]\subset \cup_{i=1}^{N_{\varepsilon}}[a_{i},a_{i+1}]$, where
 the $a_{i}$'s are not critical values and such that
 $0<a_{i+1}-a_{i}\leq \varepsilon$ for all $i\in 1,\ldots, N-1$\,. 
Any bar $\alpha^{(p)}$ of degree $p$ and
length larger than $\varepsilon$ has at most two endpoints lying in
different intervals $[a_{i},a_{i+1}]$ and appearing as an element of
$\mathcal{Z}^{(p)}(a_{i},a_{i+1})$ for the possible lower endpoint and
an element of $\mathcal{Z}^{(p+1)}(a_{i'},a_{i'+1})$ for the possible
upper endpoint with $i\neq i'$\,. Therefore, the set of bars of  degree
$p$ and length larger than $\varepsilon$ is bounded by
\begin{multline*}
\sum_{i=1}^{N_{\varepsilon}-1}\sharp
\mathcal{Z}^{(p)}(a_{i},a_{i+1})+\sharp
\mathcal{Z}^{(p+1)}(a_{i},a_{i+1})\\
=\sum_{i=1}^{N_{\varepsilon}-1}\beta^{(p)}(f^{a_{i+1}},f^{a_{i}})+\beta^{(p+1)}(f^{a_{i+1}},f^{a_{i}})<+\infty\,.
\end{multline*}
The conjecture stated in the introduction
for a general
$\mathcal{C}^{\infty}$ function $f$ has now the following more precise
version:
For $\varepsilon>0$, the $\tilde{o}(e^{-\frac{2\varepsilon}{h}})$
eigenvalues of
$\Delta^{(p)}_{f,M,h}$ are given by the $\lambda_{\alpha}^{(p)}(h)$ such that
$\alpha$ is of length
larger than $\varepsilon$,
$\alpha\in A^{(p)}$ or $(\alpha\in
A^{(p-1)}~\text{and}~b_{\alpha}^{(p)}<+\infty)$, and
$$
\lim_{h\to 0}-h\log(\lambda_{\alpha}(h))=2(b_{\alpha}-a_{\alpha})\,.
$$
Our proof, relying on a recurrence on the number of critical values by
following  their increasing (and decreasing)
order, is not  adapted to the more general case with an
infinite number of critical values. One may think of a different type
of induction: Starting from our result for  finitely many critical values, one may
increase the number of critical values by perturbing the function such
that it creates small bars in a given interval $[a,b]$, and then try to
obtain  spectral and resolvent 
estimates for the spectral parameter $\lambda
\in [0,\tilde{o}(e^{-\frac{2\varepsilon}{h}})]$, which are uniform with
respect to the additional small bars.
\\
\noindent \textbf{b) What about $\mathcal{C}^{0}$-potentials ?} The
stability of bar codes makes sense in the $\mathcal{C}^{0}$ topology
while a finite bar code can be associated with a continuous function
which satisfies Hypothesis~\ref{hyp:weakreg}. The relation between
the exponentially small eigenvalues of $\Delta_{f,h}$ and the bar code of
$f$ suggests that the bottom of the spectrum of $\Delta_{f,h}$ makes
sense only under Hypothesis~\ref{hyp:weakreg}.  Is there a natural
self-adjoint operator ``$\Delta_{f,h}$'' on $M$ when $f$ is only
continuous and for which Theorem~\ref{th:mainsimple} could be
extended~?
\\
\noindent
\textbf{c) Applications of the result on $p$-forms:}
Over decades, the case of functions has received a lot of attention
  with an easy interpretation in terms of Fokker-Planck equation
  associated with reversible processes at low temperature and within
  the modelling e.g. in chemistry as points the title of this
  text. Here is an attempt to interpret our spectral results for
  $p$-forms. This deserves more precise studies and we hope that
  relevant applications will be found.
Within the stochastic approach, the Witten Laplacian is better written
as
$$
\mathcal{L}_{f,h}=
e^{\frac{f}{h}}\Delta_{f,h}e^{-\frac{f}{h}}=h^{2}\Delta_{0,1}+2h\mathcal{L}_{\nabla
f}=d_{0,h}d_{2f,h}^{*}+d_{2f,h}^{*}d_{0,h}\,,
$$ considered in the $L^{2}$-space associated with the invariant
measure $e^{-\frac{2f}{h}}~dx$\,,
 $L^{2}(M,e^{-\frac{2 f}{h}}~dx;
\Lambda T^{*}M)$\,, and where $\Delta_{0,1}=dd^{*}+d^{*}d$ is the
usual Hodge Laplacian.
  There are formulas to express the semigroups
  associated with Hodge and Witten Laplacians, in terms of
  expectations values along brownian motion:
  $e^{-t\mathcal{L}_{f,h}}v=\mathbb{E}(\xi_{t}^{*}v)$ for $v\in
  \mathcal{C}^{\infty}(M;\Lambda T^{*}M)$\,, where $\xi_{t}$ is the
  flow associated with a stochastic differential equation of the type
 $dx=X(x_{t})\circ
  dB_{t}-2\nabla f(x_{t})dt$ where $B$ is an $m$-dimensional brownian
  motion in $\rz^{m}$ and $X:M\times \rz^{m}\to TM$ is a submersion specified by
  the metric on $M$ (see in particular \cite[Theorem~1.1.2,  formula 1.2.5, and Section~2.4]{ELJL}). Due to the
  supesymmetric argument, eigenforms of
  $\Delta_{f,h}$  (resp. $\mathcal{L}_{f,h}$) can be assumed to solve
  $d_{f,h}^{*}\omega=0$
  (resp. $d_{2f,h}^{*}\tilde{\omega}=0$ with
  $\tilde{\omega}=e^{\frac{f}{h}}\omega$), because when
  $d_{f,h}^{*}\omega\neq 0$ (resp. $d_{2f,h}^{*}{\omega})\neq 0$) then
  $d_{f,h}^{*}\omega$ (resp $d_{2f,h}^{*}\tilde{\omega}$) is 
 another eigenform of
  $\Delta_{f,h}$ (resp. of $\mathcal{L}_{f,h}$) with degree decreased
  by $1$ and associated with the
  same eigenvalue. Let $\tilde{\omega}$ be such an eigenform with
  $d^{*}(e^{\frac{2f}{h}}\tilde{\omega})=0$ and
  $\mathcal{L}_{f,h}\tilde{\omega}_{h}=\lambda_{h}\tilde{\omega}_{h}$\,. By
  assuming that $\tilde{\omega}$ is a $p$-form and after 
  normalization,
  $A_{h}e^{\frac{2f}{h}}\tilde{\omega}$ may be identified with
  a $p$-cycle via 
$$
\int_{M}\eta\wedge (\star
e^{\frac{2f}{h}}A_{h}\tilde{\omega})=\int_{C_{\tilde{\omega},h}}\eta\,,
$$
where $\partial C_{\tilde{\omega},h}=0$ is a consequence of
$d^{*}(e^{\frac{2f}{h}}\tilde{\omega})=0$\,. It would be better to
think of $C_{\tilde{\omega},h}$ as a courant but let us forget the
regularity issues. When $f$ is a Morse
function with
$$f(x_{1},\ldots,x_{p},x_{p+1},\ldots,x_{d})=-\varphi_{-}(x_{1},\ldots,x_{p})+\varphi_{+}(x_{p+1},\ldots,x_{d})$$
around a critical point of index $p$ with critical value $0$
which is a lower endpoint of a bar of degree  $p$\,,
the leading term of the WKB-approximation says
$e^{\frac{2f}{h}}\tilde{\omega}=e^{-\frac{2\varphi_{+}(x_{p+1},\ldots,x_{d})}{h}}dx_{1}\wedge\ldots\wedge
dx_{p}$
and $C_{\tilde{\omega},h}$ is assymptotically equal to some
fixed cycle $C_{\tilde{\omega},0}$ supported by the unstable manifold
of $-\nabla f$\,. We may expect such a behaviour in general.
The evolution
$\tilde{\omega}_{h}(t)=e^{-t\mathcal{L}_{f,h}}\tilde{\omega}_{h}=e^{-t\lambda_{h}}\tilde{\omega}_{h}$\,
says that this cycle is not changed by the dynamics when
$t<<\frac{1}{\lambda_{h}}$ and disappears when
$t>>\frac{1}{\lambda_{h}}$\,. The reverse eigenvalue
$\frac{1}{\lambda_{h}}$ appear as the lifetime of the cycle
$C_{\tilde{\omega},h}$ of which an asymptotic form $C_{\tilde{\omega},0}$ is expected when
the normalization factor $A_{h}$ is well chosen. Below is a picture
for the brownian dynamics of a $1$-cycle, which shows the
generalization of the metastability picture that we expect.
\begin{figure}[h]
\centering{
\includegraphics[width=8cm]{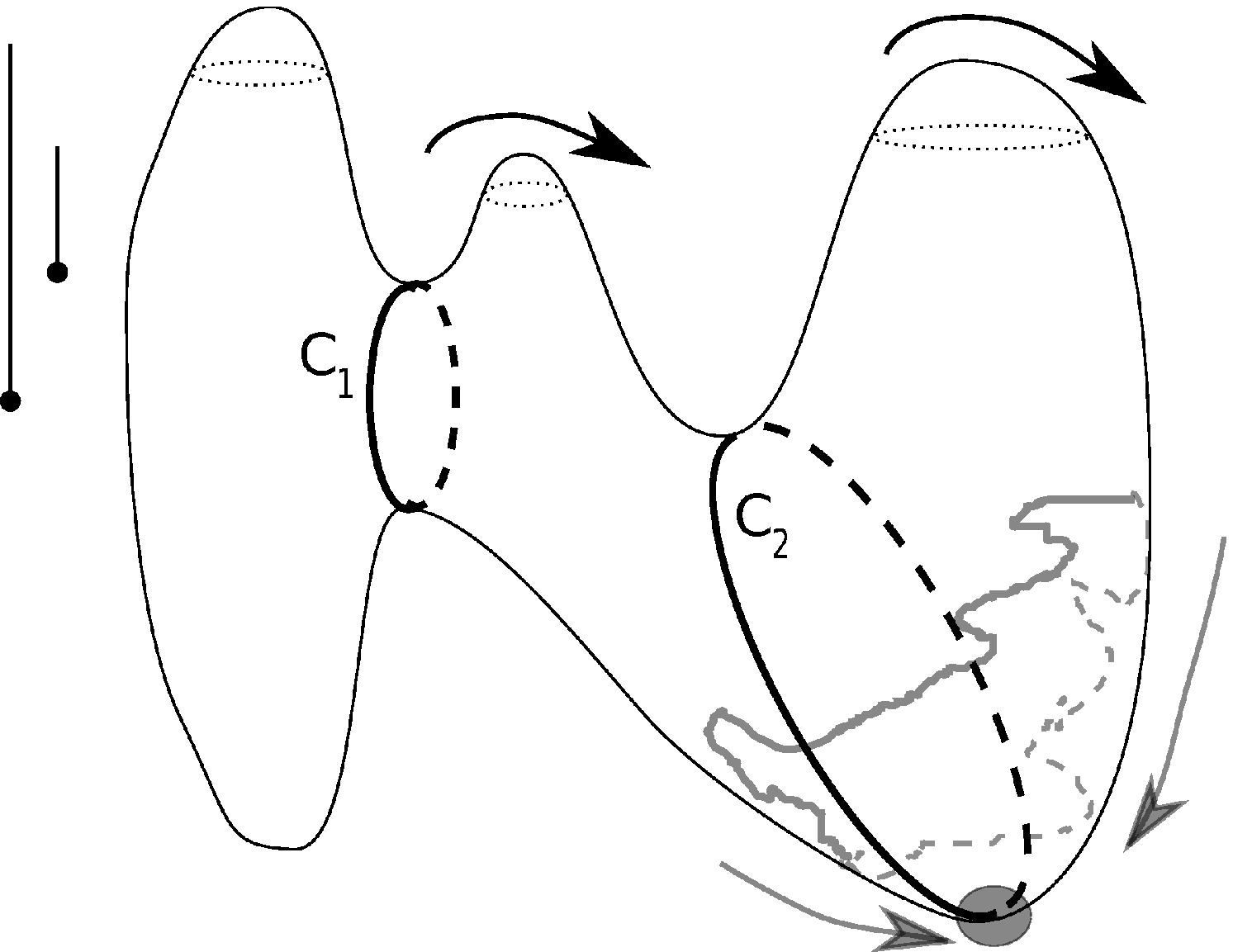}}
\captionsetup{labelformat=empty,width=14cm}
\caption{\textbf{Figure 25:} Metastability of cycles: The bars of degree $1$ represented
  on the left-hand side, with lengths
  $\ell_{1}<\ell_{2}$, provide the lifetime $e^{\frac{2
      \ell_{1}}{h}}$ (resp $e^{\frac{2\ell_{2}}{h}}$) of the cycle
  $C_{1}$  (resp. $C_{2}$).  After a time larger than the lifetime,
  $C_{1}$ is first deformed  into $C_{2}$ and $C_{2}$ is  then 
  deformed
   into the
  grey cycle which is rapidly retracted to the global minimum.}
\end{figure}
\\
\noindent
\textbf{d) General statement for subexponential factors:}
Specifying the exponential scales  of the exponentially small eigenvalues of
  $\Delta_{f,M,h}^{(p)}$ associated with the lengths of  the
 bar code of $f$ was
  done in
  Theorem~\ref{th:mainsimple} and Theorem~\ref{th:specres}, while the
  spectral version of the stability was given
 in Corollary~\ref{cor:mainsimple} and
  Theorem~\ref{th:stab1}. Those results 
 are general statements which
  hold under simple general assumptions like
 Hypothesis~\ref{hyp:mainf} or
  Hypothesis~\ref{hyp:realana}. The situation is different when we
  want to specify the subexponential factors.
In Section~\ref{sec:applications},
  the general construction was used in order to specify the
  subexponential factors and to show that they were keeping some kind
  of stability property, possibly within a finite dimensional matricial
  writing (see Proposition~\ref{pr:genalMorse}). Although the method is clear and heavily relies on
  Theorem~\ref{th:induc} and the use of Stokes' theorem like in
  \cite{LNV}, we were not able to take into account
 all the possible
  configurations in a uniform and satisfactory presentation. Although
  the stability of individual subexponential factors cannot hold,
a general robust
  statement or formula for the determination of
  the subexponential factors would be valuable.
\\
\noindent
\textbf{e) Piecewise affine functions and discretization via triangulation:} In the one dimensional case, a schematic Witten Laplacian for
  which everything relies on simple linear algebra is provided by a
  piecewise affine function $f$\,. Eigenforms of degree $0$ or $1$ are
  computed by matching exponentials at the discontinuities of the
  slope of $f$\,. It becomes a fully discrete model, in its coding and
  in the computation of the eigenforms.
 The generalization of a piecewise affine function 
 after a triangulation of 
 $\rz^{d}$ or $\tz^{d}$  (and for further generalizations, one should consider  a
 Lipschitz triangulated riemannian manifold like in \cite{GMM})
 enters in our general
 assumption  Hypothesis~\ref{hyp:realana}.  Away from the
 singularities of $f$\,, the Witten Laplacian is nothing but a scalar
 operator $-\Delta+V(x)$, where $V$ is a piecewise constant function,  while the Hessian of $f$
brings a measure potential carried by the singularities of $f$\,. We
are led to consider a specific self-adjoint extension of
$-\Delta+V(x)$ on $\mathcal{C}^{\infty}_{0}(\Omega_{reg};\Lambda T^{*}M)$,
where $\Omega_{reg}$ is the open domain where $f$ is differentiable
with a locally constant gradient. 
Many things have been done on the scalar
 Laplacian plus simple or double layer potentials, or more general
 interface conditions (see \cite{AGHKH,BGP}). Here we work with
 Hodge-type Laplacians and
 discriminating with respect to the degree will lead to different types
 of interface conditions and we wonder whether cohomology brings
additional restrictions along strata of codimension $>1$\,.
It would be interesting to see 
if such a finitely coded potential $f$ leads
to a completely solvable linear algebra 
problem like in dimension $1$\,. It could be an alternative model
problem as compared to the case of Morse functions, which could be
useful to understand some non trivial boundary or corner problems.
\\
\noindent
\textbf{f) Infinite or large dimensional problems:} 
After specifying the geometrical problems, especially concerned
  with the domain issues for the differential, codifferential, and
  Witten Laplacian, all the analysis is carried out along the real
  axis of values of $f$\,, $\rz\supset f(M)$\,. In this projective
  perspective, the dimension of $M$ does not count until the
  computation of the subexponential factors, which involve the
  asymptotics of Laplace integrals. This raises the question of the
  validity of such an approach for infinite dimensional --
   or large
  dimensional --  problems,
  which have applications in
  statistical physics, and
  where the asymptotic behaviour of the dimension
  is related with the small parameter $h\to 0^{+}$ (see e.g. \cite{HelW,DGLP}, or the recent
  \cite{DGMo} where the estimates when $h\to 0^{+}$ are even shown to be uniform in the dimension,  and references
  therein).
\\
\noindent
\textbf{g) Other boundary conditions for Witten Laplacians:}
Our results include the case of Witten Laplacians on bounded
  domains like $f_{a}^{b}$\,, provided that one considers Neumann boundary
  conditions on the upper boundary $f^{-1}(\left\{b\right\})$ and
  Dirichlet boundary conditions on $f^{-1}(\left\{a\right\})$\,. In
  some applications like in the analysis of quasi-stationary
  distributions, it is relevant to put Dirichlet boundary conditions
  everywhere on $\partial \Omega$ when the manifold $M$ is replaced at
  the beginning by some regular domain
  (see \cite{LeNi, DLLN1,DLLN2, LeNe1, LeNe2}). The cohomology groups $H^{*}(f^{b};f^{a})$
  have to be replaced by $H^{*}(f^{b};f^{a}\cup\partial \Omega)$, but
  additional corner problems at the intersection $\partial
  \Omega\cap f^{-1}(\left\{a,b \right\})$ have to 
be analyzed carefully. 
\\
\noindent
\textbf{h) Non reversible dynamics and spectral analysis of non
  self-adjoint related problems:} }The analysis of Witten Laplacians enters in the general scope of
  the semiclassical analysis of self-adjoint Schr{\"o}dinger-type
  operators. Within the stochastic analysis, several models, motivated
  by applications where a non reversible drift is considered, 
lead to non self-adjoint operators for which a
  similar analysis can be carried out in the case of functions, $p=0$
  (see e.g \cite{LeMi}). An interesting non self-adjoint (and non
  elliptic) operator which has many connections with Witten Laplacian
  is Bismut's hypoelliptic Laplacian, which is defined in any degree
  $0\leq p\leq 2d$ when we work on $\mathcal{X}=T^{*}Q$ with $\dim
  Q=d$\,. The asymptotic behaviour of exponentially small eigenvalues
  has been studied so far only when $p=0$ and $Q=\rz^{d}$ in
  \cite{HHS}, where Bismut's hypoelliptic Laplacian is nothing but the
  Kramers-Fokker-Planck operator of kinetic theory. For studying the
  case of general $p$-forms on a manifold, a better understanding of
  boundary conditions for Bismut's hypoelliptic Laplacians (defined in
  \cite{Nie}) is necessary. When $f:Q\to \rz$ is the potential,
  adapting the analysis of this text would lead to  ``Dirichlet
  boundary conditions'' on $T^{*}_{f^{-1}(\left\{a\right\})}Q$ and ``Neumann
  boundary conditions'' on $T^{*}_{f^{-1}(\left\{b\right\})}Q$ for the
  hypoelliptic Laplacian acting in $\pi^{-1}(f_{a}^{b})$\,, where
  $\pi:T^{*}Q\to Q$ is the fiber projection.  Additionally, the non
  self-adjoint nature of the problem requires different techniques
  relying on complex deformations in order to handle the exponential
  decay of resolvents and eigenfunctions.
\\
\noindent
\textbf{i) Remarks about the subanalytic case:} In the subanalytic case and for at least the
 second time (a previous time was in \cite{GeNi} for the
analysis of Mourre estimates for analytically fibered operators), the differentiation along regular strata has been
used in order to prove estimates.
Instead of considering a non regular solution $\phi$ to the
Hamilton-Jacobi equation $|\nabla f|^{2}=|\nabla \phi|^{2}$\,, we
constructed a finite family of regular functions $\phi_{k}$\,,
$k=1,\ldots, K$\,, $|\nabla
f|^{2}\geq |\nabla \phi_{k}|^{2}$, finally leading to a good enough
exponential decay estimate. We were not able to make a direct use of
viscosity solutions, which did not allow to absorb all the singular
terms in Agmon's type estimates.  In a different context,  global 
subanalytic viscosity 
solutions to Hamilton-Jacobi  with analytic coefficients (which is not
the case here) 
were studied in \cite{Tre}. Is there a better way to introduce
 viscosity solutions in our problem~? In the other way,
 differentiating along the regular strata could it be used for constructing
 subsolutions to Hamilton-Jacobi type equations~?
\\
\noindent
\textbf{j) Fukaya conjecture and multidimensional persistence:}
Determining the homotopy type of a compact manifold $M$ such that
$\pi_{1}(M)=0$ and the $A_{\infty}$ structure on
harmonic forms induced by the pullback of the wedge product, can be
attacked via Witten's deformation. This was proposed by Fukaya in
\cite{Fuk} and more precisely studied via WKB methods a la
Helffer-Sj{\"o}strand in 
\cite{CLM}. It consists in considering several Witten's deformations
of the differential and the Hodge Laplacian,
$d_{f_{ij},h}=e^{-\frac{f_{ij}}{h}}(hd)e^{\frac{f_{ij}}{h}}$\,,
associated with a sequence $(f_{0},f_{1},\ldots, f_{k})$ such that
$f_{ij}=f_{j}-f_{i}$\,, $0\leq i<j\leq k$\,, are Morse
functions. Although it may not bring an additional topological
information, replacing Morse functions by more general
$\mathcal{C}^{\infty}$ functions means the understanding of the
$\frac{k(k+1)}{2}$-dimensional version of persistence diagrams, bars
being replaced by multidimensional objects. The multidimensional
version of persistence homology, partly motivated by applications in
statistical data analysis, is only emerging. We refer again to
\cite{KaSc} for a theoretical presentation of multidimensional persistence.
\\
\noindent\textbf{k) Comparison with the instantonic picture:}
The instantonic picture makes sense within Thom-Smale transversality
condition, which ensures that any critical point
of index $p+1$ is
related to some critical points of index $p$ by a finite number of
regular integral curves of $-\nabla f$\,. This gives rise to the standard 
Thom-Smale complex
structure. More recently, it has received an accurate analysis in terms
of the analysis of the dynamical system of $-2\nabla f$ perturbed by a
brownian motion in \cite{DaRi} by applying Faure-Sj{\"o}strand theory of
weighted Sobolev spaces. We already mentioned that our approach is
orthogonal to the instantonic point of view: Instead of exploring the
geometry of the potential landscape $M\ni x\mapsto f(x)\in \rz$\,, we
considered globally the sublevel sets $f^{\lambda}$ and their
homological properties.  We can parallel this with the comparison
between Riemann's and Lebesgue's integration theory. This global
approach avoids to consider possibly complicated cancellation
phenomena in the general method of tunnel effect computations
described in \cite{HeSj2,HeSj3}. It is a question whether such a
global and topological approach makes sense for other spectral
problems related with dynamical systems.

\appendix

\section{Abstract Hodge theory}
\label{sec:abstHodge}

The abstract version of Hodge theory provides spectral results,
like \eqref{eq.Res-commut} or  Corollary~\ref{co.abs-hodge} below, which hold in general with
weak regularity assumptions. For a proof, we refer for example to
\cite[Section~2]{GMM} (see in particular Propositions 2.3 and 2.4, Corollary 2.5, and Theorem 2.8 there).
\begin{prop}
\label{pr.abs-hodge}
Let $(H, \|\cdot\|_{H})$ be a Hilbert space and let $T:D(T)\subset H \to H$
be a closed  densely defined unbounded linear operator such that  
$$\Ran\, T\subset \ker T\ \ \text{and}\ \ D(T)\cap D(T^{*})\ \ \text{embeds compactly into $H$}\,,$$ 
where $D(T)\cap D(T^{*})$ is equipped with the graph norm
$$
\|u\|_{D(T)\cap D(T^{*})}\ :=\ \sqrt{\|u\|^{2}_{H}+\|Tu\|^{2}_{H}+\|T^{*}u\|^{2}_{H}}\,.
$$
We then have the following properties:
\begin{enumerate}
\item[i)] The operator  $(T+T^{*}, D(T)\cap D(T^{*}))$ is self-adjoint with a compact resolvent
and satisfies 
$$\ker(T+T^{*})=\ker T \cap \ker T^{*}\,.$$
In particular, the linear space $D(T)\cap D(T^{*})$ is dense in $H$ and  $T+T^{*}$
is a self-adjoint Fredholm operator with index $0$\,, that is more precisely
$$
\ker T \cap \ker T^{*}\  \ \text{has finite dimension}
\quad\text{and}
\quad 
\Ran(T+T^{*})\ =\ \big(\ker T \cap \ker T^{*}\big)^{\perp}\,.
$$
\item[ii)] The operator $\Delta:=TT^{*}+T^{*}T$ with domain 
$$D(\Delta)\ :=\ \{u\in D(T)\cap D(T^{*})\ \text{s.t.}\ Tu \in D(T^{*})\ \text{and}\ T^{*}u\in D(T)\}$$
is a nonnegative self-adjoint operator with kernel
$$
\ker \Delta\ =\ \ker T \cap \ker T^{*}\ =\ \ker(T+T^{*})\,.
$$
In particular,  $\Delta$ has a compact resolvent (since $D(\Delta)$
with its graph norm  embeds continuously into
$D(T)\cap D(T^{*})$)
and is the Friedrichs extension
associated with the closed nonnegative quadratic form $Q$ on $D(T)\cap D(T^{*})$
defined by
$$
Q(u,v)\ :=\ \langle Tu,Tv\rangle_{H}+\langle T^{*}u,T^{*}v\rangle_{H}\,.
$$ 
\end{enumerate}
\end{prop}

Let  us also note the following  consequences  of Proposition~\ref{pr.abs-hodge}
underlining the supersymmetric structure of the operator $\Delta$ defined there :
when $T$ is as in the statement of Proposition~\ref{pr.abs-hodge}, the
resolvent satisfies 
for every $z\in\cz\setminus \sigma(\Delta)$\,,  $u\in D(T)$\,, and $u'\in D(T^{*})$\,,
\begin{equation}
\label{eq.Res-commut} 
(z-\Delta)^{-1}\, T\,u\ =\ T\,(z-\Delta)^{-1}\, u\quad\text{and}
\quad 
(z-\Delta)^{-1}\, T^{*}\,u'\ =\ T^{*}\,(z-\Delta)^{-1}\, u'\,.
\end{equation}
Let us prove the first relation, the proof of the second one being similar.
Let us then consider 
 $u\in D(T)$ and let us define
$v=(z-\Delta)^{-1} u$ for some $z\in\cz\setminus \sigma(\Delta)$\,.
Then $v\in D(\Delta)$ and $(z-\Delta)v=u\in D(T)$\,,
which implies $\Delta v=T^{*}Tv+TT^{*}v\in D(T)$ and hence,
 since $\Ran T\subset \ker T$\,,
$T^{*}Tv\in D(T)$\,.
In particular, one has $Tv\in D(TT^{*})$\,, and hence
$Tv\in D(\Delta)$\,,
and 
$$
(z-\Delta)Tv\ =\ zTv-TT^{*}Tv\ =\ T(z-\Delta)v\ =\ Tu
\quad\text{and then}\quad  Tv\ =\ (z-\Delta)^{-1}Tu\,,
$$
that is precisely the first relation in \eqref{eq.Res-commut}.\medskip

\noindent
An easy consequence of \eqref{eq.Res-commut} is the following: 
for any eigenvalue $\lambda$ of $\Delta$
and associated eigenvector $u\in D(\Delta)$\,, we have
$
Tu\in D(\Delta)
$ and
$
T^{*}u\in D(\Delta)
$\,,
with
\begin{equation}
\label{eq.Res-commut-2} 
T\,\Delta\,u\ =\ 
\Delta\,T\,u
\ =\ \lambda\,T\,u
\quad\text{and}\quad
T^{*}\,\Delta\,u\ =\ 
\Delta\,T^{*}\,u
\ =\ \lambda\,T^{*}\,u
\end{equation}
Note that if in addition $\lambda\neq 0$\,, one element among $Tu,T^*u$ is nonzero
(since in this case $u\notin\ker \Delta = \ker T \cap \ker T^{*}$).
\begin{cor}
\label{co.abs-hodge} Assume the hypotheses of Proposition~\ref{pr.abs-hodge}
and define $\Delta:=TT^{*}+T^{*}T$ as there.
The following orthogonal decompositions then hold:
$$
H\ =\ \Ran T\,\mathop{ \oplus}^{\perp}\, \Ran T^{*}\,\mathop{ \oplus}^{\perp}\,\ker \Delta\quad \text{and, for
$\mathbf T=T$ or $\mathbf T=T^{*}$\,,}\quad
\ker \mathbf T\ =\ \Ran \mathbf T\,\mathop{ \oplus}^{\perp}\,\ker \Delta\,.
$$
In particular, the operators $T$ and $T^{*}$ have  closed ranges
and 
$$
\ker T/\Ran T \ \simeq\ \ker T^{*}/\Ran T^{*}\ \simeq\ \ker \Delta\,.
$$
\end{cor}

\begin{proof}
This result is the statement of \cite[Proposition~2.9]{GMM} 
and is an easy consequence of Proposition~\ref{pr.abs-hodge}.
First, since 
$\Ran(T+T^{*})= \big(\ker T \cap \ker T^{*}\big)^{\perp}=(\ker \Delta)^{\perp}$
according to Proposition~\ref{pr.abs-hodge},
we deduce the inclusions (since $T$ and $T^{*}$ are closed),
$$
\Ran T + \Ran T^{*}\ \supset\ (\ker \Delta)^{\perp}\ =\ \overline{\Ran T + \Ran T^{*}}\ \supset\ 
\Ran T + \Ran T^{*}\,.
$$
The linear space $\Ran T + \Ran T^{*}$ is then closed in $H$
and, owing to $T^{2}=0$\,,  this sum 
is moreover orthogonal. The spaces $\Ran T$ and $\Ran T^{*}$
are consequently closed and 
$$
H\ =\ (\ker \Delta)^{\perp}\,\mathop{ \oplus}^{\perp}\,\ker \Delta
\ =\ \Ran T\, \mathop{ \oplus}^{\perp}\, \Ran T^{*}\,\mathop{ \oplus}^{\perp}\,\ker \Delta\,.
$$
Furthermore, the inclusion $\ker T\supset \Ran T\mathop{ \oplus}^{\perp}\ker \Delta$
is clear, owing again to $T^{2}=0$\,. To prove the reverse inclusion,
just notice that any $v\in \ker T$ writes as the sum 
$v=u_{0}+Tu_{1}+ T^{*}u_{2}$\,, where $u_{0}\in \ker \Delta$\,,
$u_{1}\in D(T)$\,, and $u_{2}\in D(T^{*})$\,. It follows that
$T^{*}u_{2}\in D(T)$ and $TT^{*}u_{2}=0$\,, which implies $T^{*}u_{2}=0$ (by taking the scalar product with $u_{2}$)
and then $v=u_{0}+Tu_{1}\in \Ran T\,\mathop{ \oplus}^{\perp}\,\ker \Delta$\,.\\
Lastly, the relation $\ker T^{*} = \Ran T^{*}\mathop{ \oplus}^{\perp}\ker \Delta$
follows by applying the relation $\ker T = \Ran T\mathop{ \oplus}^{\perp}\ker \Delta$
 with  $T$ replaced by $T^{*}$\,, which  satisfies $\Ran T^{*}\subset \ker T^{*} $ and $T^{**}=T$\,.
\end{proof}

\section{Persistent cohomology and bar codes}
\label{app:perscohom}

\subsection{A sheaf theoretic presentation}
\label{sec:sheaf}
Let $f$ be  a $\mathcal{C}^{\infty}$ function on the compact manifold $M$ having finitely many critical values (but we do not assume $f$ is a Morse function). 
We shall define its bar code  following the sheaf theoretic presentation of \cite{KaSc}.

The following assumption on  $f$
which is weaker than Hypothesis~\ref{hyp:mainf} allows us to use this
construction in a low regularity setting. We keep the notation of Definition~\ref{de:fab}
$$
f^{t}=\left\{x\in M\,,\quad f(x)<t\right\}\quad\text{and}\quad
f^{\leq t}=\left\{x\in M\,,\quad f(x)\leq t\right\}\,.
$$
\begin{hyp}
\label{hyp:weakreg}
The function $f:M\to \rz$ is continuous and there exist finitely many values
$\min f=c_{1}<\ldots<c_{N_{f}}=\max f$ with the following
property:
For any $n\in \left\{1,\ldots, N_{f}-1\right\}$ 
and all $a < b \in  ]c_{n},c_{n+1}[$\,,  $f^{\leq a}$ is a deformation retract of $f^{\leq b}$\,.
The values $c_{1},\ldots, c_{N_{f}}$ are called ``critical values'' of $f$\,.
\end{hyp}
 We shall need the following 
 
 \begin{lem} 
 With the assumptions of Hypothesis \ref{hyp:weakreg}, we the space  $H^*(f^b,f^a)$ is finite dimensional. 
 \end{lem} 
 \begin{proof} 
 It is enough to prove that if $t$ is  in some $]c_j, c_{j+1}[$\,,
 then $H^*(f^{\leq t})$ is finite dimensional.  The general case
 follows by applying the long exact sequence of the pair $(f^{\leq b},
 f^{\leq a})$\,. Now let  $\varepsilon $ be small enough, $g$ a smooth
 function such that $ \Vert g-f \Vert \leq \varepsilon $\,. Then the inclusions
 \begin{displaymath} 
f^{\leq t} \subset g^{\leq t + \varepsilon } \subset f^{\leq t + 2 \varepsilon }
 \end{displaymath} 
hold true
 and for  $ \varepsilon $ small enough, 
 \begin{displaymath} 
 f^{\leq t- 2 \varepsilon } \subset f^{\leq t} \subset f^{\leq t + 2 \varepsilon }
 \end{displaymath} 
 are homotopy equivalences. Notice that $g$ being smooth and $t + \varepsilon $ being a regular value for $ \varepsilon $ generic , the cohomologies $H^*(g^{t+ \varepsilon })$ are finite dimensional, and we have maps 
 \begin{displaymath} 
 H^*(f^{\leq t+ 2 \varepsilon })\longrightarrow H^*(g^{\leq t+ \varepsilon })\longrightarrow H^*(f^{\leq t })  \end{displaymath} 
 but the composition of the above two arrows must be an isomorphism, and it factors through a finite dimensional space, therefore $H^*(f^{\leq t })$ is finite  dimensional and we have a uniform bound for $t$ in $]c_j, c_{j+1}[$\,. 
 \end{proof} 
By using the deformation along the gradient flow away from the ``critical
values'' $c_{1},\ldots, c_{N_{f}}$\,, Hypothesis \ref{hyp:weakreg} is obviously
true when $f$ satisfies Hypothesis~\ref{hyp:mainf}. It is also true
 for a general Lipschitz function satisfying Hypothesis \ref{hyp:Lipbar} as mentioned in Subsection \ref{sec:moregenLipgen}. 
It implies that for any $a,b\not\in \left\{c_{1},\ldots,
  c_{N_{f}}\right\}$\,, $a<b$\,, the relative homology groups
($\kz$-vector spaces) $H^{*}(f^{\leq b},f^{\leq a};\kz)$ are finite
dimensional and change only when $a$ or $b$ passes a ``critical
value'', $c_{1},\ldots, c_{N_{f}}$\,.\\
For the introduction of a persistent sheaf on $\rz$\,, we need to
consider all the sublevel sets, and only at the end, do we restrict our
attention 
to the relative cohomology groups $H^{*}(f^{b},f^{a};\kz)$ with
$a<b$\,, $a,b\not\in \left\{c_{1},\ldots,c_{N_{f}}\right\}$\,. 

In order to use standard results of sheaf theory it is better to work
with the closed sublevel set $f^{\leq t}$ for a general $t\in\rz$
which may be a ``critical value''.\\
For a field $\kz$\,, $\kz_{M}$ denotes the locally constant sheaf on
$M$ and we consider a $c$-soft injective   resolution
$$
\xymatrix{
0\ar[r]& \kz_{M}\ar[r] &{\cal L}^{0}\ar[r]& {\cal L}^{1}\ar[r] &\ldots
}\,,
$$
$c$-soft meaning that the restriction morphism
 $\Gamma(M;{\cal L}^{q})\to
\Gamma(K;{\cal L}^{q})$ is surjective for any compact subset $K\subset
M$ and any $q\in\nz$\,. A bounded c-soft resolution ending with ${\cal
  L}^{\dim M}\to 0$  exists
because $M$ is a compact manifold.\\
Such a resolution can be obtained by introducing the canonical
injective resolution or the sheaf of $\kz$-valued
 Alexander-Spanier cochains on $M$\,. When $\kz=\rz$ or $\cz$ we can use the de Rham
 resolution
$$
\xymatrix{
0\ar[r]& \kz_{M}\ar[r] &{\cal C}^{\infty}(M;\kz)\ar[r]^-{d}& {\cal C}^{\infty}(M;\Lambda^{1}T^{*}M)\ar[r] ^-{d}&\ldots
}\,.
$$
 showing that $\kz_{M}$ is quasi-isomorphic to the de~Rham complex
$$
\xymatrix{
0\ar[r] &{\cal C}^{\infty}(M;\kz)\ar[r]^-{d}& {\cal C}^{\infty}(M;\Lambda^{1}T^{*}M)\ar[r] ^-{d}&\ldots
}\,.
$$
and the homology groups of $\kz_{M}$\,, denoted $H^{i}(M;\kz)$\,, are obtained by computing the homology of the complex ${\cal L}^\bullet$\,.
\\
For any locally closed subset (i.e. the intersection of a closed and
an open set) $A$\,,  ${\cal L}_{A}$ is $c$-soft. When $A$ and $B$  are closed,
$A\subset B$\,, the
short exact sequence
$$
\xymatrix{
0\ar[r]&{\cal L}^{\bullet}_{B\setminus A}\ar[r]&{\cal L}^{\bullet}_{B}\ar[r]&{\cal L}^{\bullet}_{A}\ar[r]&0
}
$$
leads to the long exact sequence
$$
\xymatrix{
\cdots\ar[r]&
H_{c}^{*}(B\setminus A, {\cal L}^{\bullet})\ar[r]&
H^{*}(B,{\cal L}^{\bullet})\ar[r]&
H^{*}(A,{\cal L}^{\bullet})\ar[r]&
H^{*+1}_{c}(B\setminus A, {\cal L}^{\bullet})\ar[r]&
\cdots
}
$$
With our choice of ${\cal L}^{\bullet}$\,, this says
\begin{equation}
  \label{eq:longexact}
\xymatrix{
\cdots\ar[r]&
H_{c}^{*}(B\setminus A, \kz)\ar[r]&
H^{*}(B,\kz)\ar[r]&
H^{*}(A,\kz)\ar[r]&
H^{*+1}_{c}(B\setminus A, \kz)\ar[r]&
\cdots
}
\end{equation}
when $A$ is a closed subspace of $M$\,. We have just summarized
Godement's arguments for Theorem~4.10.1 of \cite{God} defining the long exact sequence associated to a closed subset. For general values $a<b$ in $\rz$\,, 
 the
relative cohomology groups
$H^{*}(f^{\leq b},f^{\leq a};\kz)$ can be understood
in terms of  the cohomology 
with compact support in $\left\{x\in M,
  a<f(x)\leq b\right\}$\,. 
Under Hypothesis~\ref{hyp:weakreg},
$H^{*}(f^{\leq a-\varepsilon'},\kz)\sim
H^{*}(f^{\leq a-\varepsilon},\kz)$ for any
$\varepsilon,\varepsilon'>0$ small enough,
the Mittag-Leffler condition (see \cite{KaScBook}-chap~I) is satisfied and the
cohomology groups 
of open sublevel sets  are given by 
 the projective limits
$H^{*}(f^{a};\kz)=\ds\varprojlim_{\varepsilon\to
  0^{+}}H^{*}(f^{\leq a-\varepsilon};\kz)\sim
H^{*}(f^{\leq a-\varepsilon_{0}},\kz)$ for $\varepsilon_{0}>0$
small enough\,.\\
Persistent cohomology is introduced in this way in \cite{KaSc} (we
refer the reader to \cite{CaZo}\cite{EdHa}\cite{LSV} for other presentations) via the
direct image functor $Rp_{*}$\,, in the derived category,
 applied to
 the locally constant sheaf
$\kz_{\Gamma_{f}^{+}}$ on
$\Gamma_{f}^{+}=\left\{(x,t)\in M\times \rz, f(x)\leq t\right\}$ where
$p:M\times \rz\to (\rz,\gamma)$ is given by $p(x,t)=t$\,. The notation 
$(\rz,\gamma)$ means that $\rz$ is endowed with the non-Hausdorff  $\gamma$-topology
for which open (resp. closed) sets are $]-\infty,t[$
(resp. $[t,+\infty[$), $t\in\rz$\,. Note that here we do not need to
consider the values $\pm \infty$ because $M$ is compact.\\ 
So we set
${\cal P}=Rp_{*}\kz_{\Gamma_{f }^{+}}$\,. For a $\gamma$-open set
$]-\infty,t[$ the set of sections $\Gamma(]-\infty,t[;\mathcal{P})$
 is quasi-isomorphic to the de~Rham complex
$$
\xymatrix{
0\ar[r]&{\cal C}^{\infty}(f^{t};\kz)\ar[r]^-{d}& {\cal C}^{\infty}(f^{t};\Lambda^{1}T^{*}M)\ar[r] ^-{d}&\ldots
}\,,\text{when}~\kz=\rz~\text{or}~\cz\,,
$$
while the stalk  at $t\in \rz$\,, ${\cal
  P}_{t}=\ds\varinjlim_{t'>t}\Gamma(]-\infty,t'[;{\cal P})$ is
quasi-isomorphic to the de~Rham complex of $f^{\leq{t}}$\,.
With the $\gamma$-topology on $\rz$ an example of a locally constant
sheaf is $\kz_{[a,b[}$\,, $-\infty< a<b\leq +\infty$ with
$$
\mathrm{Hom}(\kz_{[a,b[};\kz_{[c,d[})=\left\{
  \begin{array}[c]{ll}
    \kz&\text{if}~a\leq c<b\leq d\\
0&\text{else}
  \end{array}
\right.
$$
Under Hypothesis~\ref{hyp:weakreg},
the cohomology $H^{*}(f^{<t};\kz)$ is finite dimensional and locally
constant on $\rz\setminus\left\{c_{1},\ldots,c_{N_{f}}\right\}$\,. 
Therefore the sheaf $\mathcal{P}$ is an $\rz$-constructible
sheaf of $\kz$-vector spaces. By applying results of 
Crawley-Boevey in \cite{Cra} (see also Guillermou in \cite{Gui}), Kashiwara and Schapira show in
\cite{KaSc}
that 
$$
{\cal P}\sim \mathop{\oplus}_{p=0}^{\dim M}
\mathop{\oplus}_{\alpha\in A^{(p)}}\kz_{[a_{\alpha}^{(p)},b_{\alpha}^{(p+1)}[}[p]\,,\quad
-\infty<a_{\alpha}^{(p)}< b_{\alpha}^{(p+1)}\leq+\infty\,.
$$
As pointed out in \cite{KaSc} this equivalence has to be
  understood as an equivalence of the \underline{objects} in the
  bounded derived category, for
  $\mathop{Ext}^{1}(\rz_{[0,+\infty[}, \rz_{]-\infty,0]})=\rz$\,. This
  subtlety has no consequence as long as we focus on those objects
  which are the $H^{j}(f^{<t},\rz)$\,.\\
Because the sheaf is locally constant in
$\rz\setminus\left\{c_{1},\ldots, c_{N_{f}}\right\}$\,, the endpoints
$a_{\alpha}$ belong to $\left\{c_{1},\ldots,c_{N_{f}}\right\}$ and the
endpoints $b_{\alpha}$ to
$\left\{c_{2},\ldots,c_{N_{f}},+\infty\right\}$\,.
The reason why we put the exponent ${}^{(p+1)}$ for $b_{\alpha}$ will
become clear below. When we allow
$a_{\alpha}= b_{\alpha}$ the finite
cardinal of $A$ can be augmented arbitrarily by adding
$[a_{\alpha},b_{\alpha}[=\emptyset$\,, $\kz_{\emptyset}=0$\,,
 with $b_{\alpha}=a_{\alpha}$\,.\\

Remember that when $F$  is a sheaf on the topological space $X$ and $Z$ is locally
closed, $F_{Z}$ is the sheaf on $X$ characterized by
$$
\left\{
  \begin{array}[c]{l}
F_{Z}\big|_{Z}=F\big|_{Z}\\
F_{Z}\big|_{X\setminus Z}=0
  \end{array}
\right.
$$
and when $Z$ is closed one has the exact sequence 
$$
\xymatrix{
0\ar[r]&F_{X\setminus Z}\ar[r]&F\ar[r]&F_{Z}\ar[r]&0}\,.
$$
Applied to $X=(\rz,\gamma)$ and $F={\cal P}\sim \ds\mathop{\oplus}_{\alpha\in
A}\kz_{[a_{\alpha},b_{\alpha}[}$ we obtain
\begin{eqnarray*}
  {\cal P}_{[t_{0},+\infty[}&\sim& \mathop{\oplus}_{\alpha\in A,
  t_{0}<b_{\alpha}}\kz_{[\max(a_{\alpha},t_{0}),b_{\alpha}[}\,,\\
{\cal P}_{]-\infty,t_{0}[}&\sim&\mathop{\oplus}_{\alpha\in A,
  b_{\alpha}\leq t_{0}}\kz_{[a_{\alpha},b_{\alpha}[}\,,\\
{\cal P}_{[a,b[}&\sim& \mathop{\oplus}_{\alpha\in A, a<b_{\alpha}\leq
  b}\kz_{[\max(a,a_{\alpha}), b_{\alpha}[}\,.
\end{eqnarray*}
and the obvious graded analogous result holds.
From the long exact sequence \eqref{eq:longexact} written
$$
\xymatrix{\scriptstyle
\cdots\ar[r]&\scriptstyle
H^{(p-1)}(f^{\leq t})\ar[r]&
\scriptstyle
H^{(p-1)}(f^{\leq a})\ar[r]&
\scriptstyle
H_{c}^{(p)}(f^{\leq t}\setminus f^{\leq a})\ar[r]&
\scriptstyle
H^{(p)}(f^{\leq t})\ar[r]&
\scriptstyle
H^{(p)}(f^{\leq a})\ar[r]&
\scriptstyle
\cdots
}
$$
and because we are working with
$\kz$-vector spaces we obtain
$$
{\cal P}(a)^{(p)}\big|_{t}
\sim \ker[H^{(p)}(f^{\leq t};\kz)\to
H^{(p)}(f^{\leq a};\kz)]\oplus
\mathrm{coker}[H^{(p-1)}(f^{\leq t};\kz)\to H^{(p-1)}(f^{\leq a};\kz)]\,,
$$
or
$$
{\cal P}(a)^{(p)}\sim \ker({\cal P}^{(p)}_{[a,+\infty[}\to {\cal P}^{(p)}_{a})\oplus
\mathrm{coker}({\cal P}^{(p-1)}_{[a,+\infty[}\to {\cal
  P}^{(p-1)}_{a})\,.
$$
Using ${\cal P}_{[a,+\infty[}\sim \mathop{\oplus}_{\alpha\in A,
  a<b_{\alpha}}\kz_{[\max(a_{\alpha},a),b_{\alpha}[}$\,, we deduce
\begin{eqnarray*}
\ker({\cal P}^{(p)}_{[a,+\infty[}\to {\cal
  P}^{(p)}_{a})&\sim&\mathop{\oplus}_{\alpha\in A^{(p)}, a<a_{\alpha}^{(p)}}\kz_{[a_{\alpha}^{(p)},b_{\alpha}^{(p+1)}[}
\\
\mathrm{coker}({\cal P}^{(p-1)}_{[a,+\infty[}\to {\cal
  P}^{(p-1)}_{a})&\sim&
\mathop{\oplus}_{\alpha\in A^{(p-1)}, a_{\alpha}^{(p-1)}\leq
                     a<b_{\alpha}^{(p)}<+\infty}\kz_{[b_{\alpha}^{(p)},+\infty[}\,.
\end{eqnarray*}
We obtain 
$$
H_{c}^{(p)}(f^{\leq b}\setminus f^{\leq a};\kz)\sim
\left(
\mathop{\oplus}_{\alpha\in
A^{(p)}, a<a_{\alpha}^{(p)}\leq b<b^{(p+1)}_{\alpha}}\kz\right)\bigoplus\left(
\mathop{\oplus}_{\alpha\in A^{(p-1)}, a_{\alpha}^{(p-1)}\leq a<
  b_{\alpha}^{(p)}\leq b}\kz\right)\,.
$$
When $a,b$ do not belong to $\left\{c_{1},\ldots,c_{N_{f}}\right\}$\,, the inequalities in the
sums can be replaced by strict inequalities and
$$
H^{(p)}(f^{b},f^{a};\kz)\sim
\left(
\mathop{\oplus}_{\alpha\in
A^{(p)}, a<a_{\alpha}^{(p)}<b<b^{(p+1)}_{\alpha}}\kz\right)\bigoplus\left(
\mathop{\oplus}_{\alpha\in A^{(p-1)}, a_{\alpha}^{(p-1)}< a<
  b_{\alpha}^{(p)}<b}\kz\right)\,.
$$
\subsection{Trivialized complex}
\label{sec:triv}
We now establish the
relationship with the bar codes used in  \cite{LNV} which was inspired by
Barannikov's presentation of Morse theory in \cite{Bar} (see also \cite{LSV}). 

With the definitions of \cite{LNV}, the equality $\partial_Bb=a$ holds
true for two critical values $a, b$  if and only
the map 
\begin{displaymath} 
H^{p}(f^{b+ \varepsilon }, f^{a - \varepsilon }) \longrightarrow H^{p}(f^{\leq a+ \varepsilon }, f^{a - \varepsilon })
\end{displaymath} 
vanishes, while
\begin{displaymath} 
H^{p}(f^{b- \varepsilon }, f^{a - \varepsilon }) \longrightarrow H^{p}(f^{\leq a+ \varepsilon }, f^{a - \varepsilon })
\end{displaymath}
is non-zero. 
But we have $H^*(f^{y}, f^x) = H^*([x,y[, {\cal P})$\,, where ${\cal P}=Rp_{*}\kz_{\Gamma_{f }^{+}}$ and by assumption 
\begin{displaymath} {\cal P}=\mathop{\oplus}_{p=0}^{\dim M}
\mathop{\oplus}_{\alpha\in A^{(p)}}\kz_{[a_{\alpha}^{(p)},b_{\alpha}^{(p+1)}[}[p]\,,\quad
-\infty<a_{\alpha}^{(p)}< b_{\alpha}^{(p+1)}\leq+\infty\,.
\end{displaymath}  
so that 
\begin{displaymath} H^*([x,y[,{\cal P})= \mathop{\oplus}_{p=0}^{\dim M}
\mathop{\oplus}_{\alpha\in A^{(p)}}H^*([x,y[,\kz_{[a_{\alpha}^{(p)},b_{\alpha}^{(p+1)}[}[p])\,,\quad
-\infty<a_{\alpha}^{(p)}< b_{\alpha}^{(p+1)}\leq+\infty\,.
\end{displaymath} 
so it is enough to consider the case of ${\cal P}= \kz_{[a_{\alpha}^{(p)},b_{\alpha}^{(p+1)}[}[p]$ and then it is obvious that 
$\partial_B b_\alpha^{(p+1)}=a_\alpha^{(p)}$\,. 
We thus proved
\begin{prop}
If \begin{displaymath} {\cal P}=Rp_{*}\kz_{\Gamma_{f }^{+}} \end{displaymath}  and 
\begin{displaymath} \partial_B b_\alpha^{(p+1)}=a_\alpha^p ,\; \partial_B a_\alpha^{(p)}=0
\end{displaymath} 
\end{prop} 
 With the
persistent cohomology described above,
 we are now able to extend it
under the general Hypothesis~\ref{hyp:weakreg} and
 we fix the corresponding  notations.\\ 
The bar code $\mathcal{B}(f)=([a_{\alpha},b_{\alpha}[)_{\alpha\in
  A}$ associated with $f$  with $a_{\alpha}<b_{\alpha}$\,,
$a_{\alpha}\in \left\{c_{1},\ldots,c_{N_{f}}\right\}$\,, $b_{\alpha}\in
\left\{c_{2},\ldots, c_{N_{f}},+\infty\right\}$ and graded according to
$\mathcal{B}^{(p)}(f)=([a_{\alpha}^{(p)},b_{\alpha}^{(p+1)}[)_{\alpha\in
  A^{(p)}}$ is the one introduced in the previous paragraph\,. We use the
superscript ${}^{*}$ when we do not want to specify $(p)$\,.
When $a<b$ are not ``critical values'' we write
\begin{eqnarray*}
 A^{*}(a,b)&=&\left\{\alpha\in A^{*}, [a_{\alpha}^{*},b_{\alpha}^{*+1}[\cap ]a,b[ \not\in\left\{\emptyset,]a,b[\right\} \right\}\,,\\
 A_{c}^{*}(a,b)&=&\left\{\alpha\in A^{*}(a,b),
   [a_{\alpha}^{*},b_{\alpha}^{*+1}[\cap]a,b[~\text{relatively~compact~in
   }~]a,b[\right\}\,,\\
&&
\alpha\in A^{*}(a,b) \Leftrightarrow
   a<a_{\alpha}^{*}<b~\text{or}~a<b_{\alpha}^{*+1}<b\,,\\
&&
\alpha\in A_{c}^{*}(a,b)\Leftrightarrow a<a_{\alpha}^{*}< b_{\alpha}^{*+1}<b\,.\\
\end{eqnarray*}
In order to keep track of the possible multiplicities of the values
$a_{\alpha}$ and $b_{\alpha}$\,, we set
\begin{eqnarray*}
  {\cal X}^{*}(a,b)&=&\left\{(\alpha,a_{\alpha}^{*})\,, \alpha\in
                       A_{c}^{*}(a,b)\right\}\\
{\cal Y}^{*}(a,b)&=&\left\{(\alpha,b_{\alpha}^{*}), \alpha\in
                       A_{c}^{*-1}(a,b)\right\}\\
{\cal Z}^{*}(a,b)&=&\left\{(\alpha,a_{\alpha}^{*})\,, \alpha\in A^{*}(a,b)\setminus
                     A_{c}^{*}(a,b)\,, a<a_{\alpha}<b\right\}\\
&&\hspace{2cm}\sqcup
\left\{(\alpha,b_{\alpha}^{*})\,, \alpha\in A^{*-1}(a,b)\setminus
                     A_{c}^{*-1}(a,b), a<b_{\alpha}^{*}<b\right\}\,,\\
{\cal J}^{*}(a,b)&=&{\cal X}^{*}(a,b)\sqcup {\cal Y}^{*}(a,b)\sqcup {\cal Z}^{*}(a,b)\,.
\end{eqnarray*}
We now consider the complex defined on 
$$
\kz^{{\cal J}(a,b)}=\mathop{\oplus}_{p=0}^{\dim M}\kz^{{\cal J}^{(p)}(a,b)}\sim \kz^{2\sharp A_{c}(a,b)+\sharp
  (A(a,b)\setminus A_{c}(a,b)) }
$$ 
with natural basis $(x\in {\cal X}(a,b),y\in {\cal Y}(a,b),z\in {\cal
  Z}(a,b))$ and
with the differential $\mathbf{d}_{{\cal B}}$ defined by
\begin{eqnarray*}
  &&
\mathbf{d}_{{\cal B}}x^{(p)}=y^{(p+1)}\quad\text{if}~x^{(p)}\in {\cal
     X}^{(p)}(a,b)\,, y^{(p+1)}\in {\cal Y}^{(p+1)}(a,b)\,,
     p_{1}(x)=\alpha=p_{1}(y)\,,\\
&&
\mathbf{d}_{{\cal B}}y^{(p)}=0\quad\text{if}~y^{(p)}\in {\cal
   Y}^{(p)}(a,b)\,,\\
&& \mathbf{d}_{{\cal B}}z^{(p)}=0 \quad\text{if}~z^{(p)}\in {\cal
   Z}^{(p)}(a,b)\,.
\end{eqnarray*}
By construction, when $-\infty<a<b<+\infty$ are not ``critical values'' of $f$\,,
$$
H^{p}(\kz^{{\cal J}(a,b)},\mathbf{d}_{{\cal B}})=\mathop{\oplus}_{z\in
  {\cal Z}^{(p)}(a,b)}\kz z
\sim
\left(
\mathop{\oplus}_{\alpha\in
A^{(p)}, a<a_{\alpha}^{(p)}<b<b^{(p+1)}_{\alpha}}\kz\right)\bigoplus\left(
\mathop{\oplus}_{\alpha\in A^{(p-1)}, a_{\alpha}^{(p-1)}< a<
  b_{\alpha}^{(p)}<b}\kz\right)\,,
$$
and the complex $(\kz^{{\cal J}(a,b)},\mathbf{d}_{{\cal B}})$ computes all the relative
cohomology groups $H^{*}(f^{b},f^{a};\kz)$\,.\\
The sets ${\cal X}(a,b)$\,, ${\cal Y}(a,b)$\,, $A_{c}(a,b)$\,, play
a role when we compute the positive 
exponentially small 
eigenvalues of Witten
Laplacians with Dirichlet boundary conditions on $f^{-1}(\left\{a\right\})$ and
Neumann boundary conditions on $f^{-1}(\left\{b\right\})$\,.
\subsection{Stability theorem}
\label{sec:stab}

The bar code associated with $f$ is given by a family
${\cal B}(f)=([a_{\alpha},b_{\alpha}[)_{\alpha\in A}$\,, now containing possibly
empty sets when $a_{\alpha}=b_{\alpha}$\,, with the equivalence
$([a_{\alpha},b_{\alpha}[)_{A}\sim ([c_{\beta},d_{\beta}[)_{\beta\in
  B}$ if there is a bijection between $j:\left\{\alpha\in A,
  a_{\alpha}<b_{\alpha}\right\}\to \left\{\beta\in B,
  c_{\beta}<d_{\beta}\right\}$ such that $c_{j(\alpha)}=a_{\alpha}$ and
$d_{j(\alpha)}=b_{\alpha}$\,. 
Following \cite{CEH} they can be
represented as a family of  points
$((a_{\alpha},b_{\alpha}))_{\alpha\in A}$ 
in $\left\{(x,y)\in \rz\times (\rz\cup
  \left\{+\infty\right\}),~x\leq y\right\}$\,, appearing with
multiplicities, and the bottleneck distance between two general 
bar codes
${\cal B}_{A}=([a_{\alpha},b_{\alpha}[)_{\alpha\in A}$ and
${\cal B}_{B}=([c_{\beta},d_{\beta}[)_{\beta\in B}$\,, where  $A$ and $B$ can be
assumed with the same cardinal when we authorize $a_{\alpha}=
b_{\alpha}$\,, $c_{\beta}=d_{\beta}$\,, is given by 
$$
d_{bot}({\cal B}_{A},{\cal B}_{B})=\inf_{j:A\stackrel{\mathrm{bij}}{\to}B}\max_{\alpha\in
A}
\max(|a_{\alpha}-c_{j(\alpha)}|, |b_{\alpha}-d_{j(\alpha)}|)\,,
$$
with the convention $|(+\infty)-(+\infty)|=0$\,.
The stability theorem says that for two different functions $f,g$ on
$M$ which satisfy Hypothesis~\ref{hyp:weakreg}, the bottleneck distance
between the bar codes ${\cal B}(f)$  and ${\cal B}(g)$ associated with $f$ and $g$
satisfies
$$
d_{bot}({\cal B}(f),{\cal B}(g))\leq \|f-g\|_{{\cal C}^{0}}\,.
$$
It is proved in \cite{KaSc} by using the convolution
of sheaves. In the one-dimensional case and for $\varepsilon\geq 0$ we have
$\kz_{[-\varepsilon,\varepsilon]}*\kz_{[a,b[}=\kz_{[a-\varepsilon,b-\varepsilon[}$
(in terms of constructible functions  according to
\cite{Sch}, simply use $1_{[a,b[}=1_{[a,+\infty]}-1_{[b,+\infty[}$ and
$1_{[-\varepsilon,\varepsilon]}*1_{[a,+\infty[}=1_{[a-\varepsilon,+\infty[}$)
and this convolution is nothing but a translation by $-\varepsilon$ on
the real axis. Two $\rz$-constructible sheaves on $(\rz,\gamma)$\,, 
 $F,G$ are said
$\varepsilon$-isomorphic, $F\stackrel{\varepsilon}{\sim}G$\,, if there are morphisms $i:\kz_{[-\varepsilon,\varepsilon]}*F\to G$ and $j:\kz_{[-\varepsilon,\varepsilon]}*G\to F$ such that   natural morphisms
$\kz_{[-2\varepsilon,2\varepsilon]}*F\to F$ and
$\kz_{[-2\varepsilon,2\varepsilon]}*G\to G$
are  factored via
\begin{eqnarray*}
  &&
\kz_{[-2\varepsilon,2\varepsilon]}*F\stackrel{\kz_{[-\varepsilon,\varepsilon]}*i}{\to}\kz_{[-\varepsilon,\varepsilon]}*G\stackrel{j}{\to}F
\\
&&
\kz_{[-2\varepsilon,2\varepsilon]}*G\stackrel{\kz_{[-\varepsilon,\varepsilon]}*j}{\to}\kz_{[-\varepsilon,\varepsilon]}*F\stackrel{i}{\to}G\,.
\end{eqnarray*}
The bottleneck distance is then equal to 
$$
d_{bot}(F,G)=\inf\left\{\varepsilon\geq 0,  F\stackrel{\varepsilon}{\sim}G\right\}\,,
$$
and coincides with $d_{bot}({\cal B}_{A},{\cal B}_{B})$ after writing $F\sim
\mathop{\oplus}_{\alpha\in A}\kz_{[a_{\alpha},b_{\alpha}[}$ and
$G\sim \mathop{\oplus}_{\beta\in B}\kz_{[c_{\beta},d_{\beta}[}$\,.

\textbf{Acknowledgement:} The second author thanks C. Ausoni, G.~Ginot
and P. Schapira for discussions about persistent homology
and J.M.~Delort for discussions about Lipschitz subanalytic functions. The third
author acknowledges the support of the french ANR-project, ANR
MIcrolocal ANR-15-CE40-0007, and the two first authors the support of
ANR QuAMProcs ANR-19-CE40-0010-01.

\bibliographystyle{amsalpha}

\end{document}